\documentclass{cmfnre}
\usepackage[cp1251]{inputenc}
\usepackage[russian]{babel}
\usepackage{cite}
\usepackage{graphicx,float}
\usepackage{epigraph}
\usepackage{caption}
\usepackage{subcaption}
\allowdisplaybreaks

\def\currentyear{2018}
\def\currentvolume{64}
\def\currentissue{2}
\firstpage{211}  \doi{\doipaper}

\title[\sМ\sе\sт\sо\sд\s{ }\sо\sп\sе\sр\sа\sт\sо\sр\sо\sв\s{ }\sп\sр\sе\sо\sб\sр\sа\sз\sо\sв\sа\sн\sи\sя\s{
}\sи\s{ }\sк\sр\sа\sе\sв\sы\sе\s{ }\sз\sа\sд\sа\sч\sи\s{
}\sд\sл\sя\s{ }\sс\sи\sн\sг\sу\sл\sя\sр\sн\sы\sх\s{
}\sэ\sл\sл\sи\sп\sт\sи\sч\sе\sс\sк\sи\sх\s{
}\sу\sр\sа\sв\sн\sе\sн\sи\sй]
{Метод операторов преобразования и краевые задачи для сингулярных эллиптических уравнений} 
\author[Валерий Вячеславович Катрахов (1949--2010)]{В.\,В.~Катрахов} 
\address{Воронеж, Владивосток.\\~}
\author[Сергей Михайлович Ситник]{С.\,М.~Ситник}
\address{Белгородский государственный национальный исследовательский университет (<<БелГУ>>),\\
Институт инженерных технологий и естественных наук,\\
кафедра дифференциальных уравнений,\\
308015, г.~Белгород, ул. Победы, д. 85.} 
\email{sitnik@bsu.edu.ru} 

\engtitle{The Transmutation Method and Boundary-Value Problems\\for Singular Elliptic Equations} 
\engauthor[Valeriy V. Katrakhov (1949--2010)]{V.\,V.~Katrakhov} 
\engaddress{Voronezh and Vladivostok, Russia.\\~}
\engauthor[Sergey M. Sitnik]{S.\,M.~Sitnik}
\engaddress{Belgorod National Research University ``Belgorod State University,''\\
Institute of Engineering Technology and Natural Science, Belgorod, Russia.} 
\engemail{sitnik@bsu.edu.ru} 

\def\s{\hspace{-0.17ex}}

\theoremstyle{plain} 
\newtheorem{theorem}{Теорема}[section]
\newtheorem*{theorem*}{Теорема}
\newtheorem{lemma}{Лемма}[section]
\newtheorem*{lemma*}{Лемма}
\newtheorem{corollary}{Следствие}[section]
\newtheorem*{corollary*}{Следствие}
\newtheorem{property}{Свойство}[section]
\theoremstyle{definition} 
\newtheorem{definition}{Определение}[section]
\newtheorem{remark}{Замечание}[section]

\numberwithin{equation}{section}
\numberwithin{section}{chapter}

\renewcommand{\chaptermark}[1]{\markboth{\sc #1}{}}
\newcommand{\chaptermarknum}[1]{\markboth{\sc Глава \thechapter. #1}{}}
\renewcommand{\sectionmark}[1]{\markright{\sc #1}}
\newcommand{\sectionmarknum}[1]{\markright{\sc\thesection. #1}}

\def\le{\leqslant}
\def\ge{\geqslant}
\def\leq{\leqslant}
\def\geq{\geqslant}

\DeclareCaptionLabelFormat{mySubCaptionLabel}{}
\makeatletter
\renewcommand*\l@subsection{\@dottedtocline{2}{2em}{3.3em}}
\makeatother

\def\p{{,}}
\def\Re{\mathop{\rm Re}\nolimits}
\def\supp{\mathop{\rm supp}\nolimits}
\def\arg{\mathop{\rm arg}\nolimits}
\def\const{\mathop{\rm const}\nolimits}
\def\ds{\displaystyle}

\newcommand{\lr}[1]{\left(#1\right)}
\newcommand{\sq}{\sqrt}
\newcommand{\ov}{\overline}
\newcommand{\pr}{\partial}
\newcommand{\pd}{\partial}
\newcommand{\lrs}[1]{\left[#1\right]}
\newcommand{\N}{\mathbb{N}}

\newcommand{\R}{\mathbb{R}}
\newcommand{\Cbb}{\mathbb{C}}

\begin{document}

\begin{abstract}
Основное содержание книги составлено из материалов двух докторских
диссертаций: В.\,В.~Катрахова 1989~г. и С.\,М.~Ситника 2016~г. В
книге впервые  в формате монографии систематически изложена теория
операторов преобразования и их приложений для дифференциальных
уравнений с особенностями в коэффициентах, в том числе содержащих
операторы Бесселя. Наряду с детальной обзорной информацией и
библиографией по указанной тематике книга содержит оригинальные
результаты авторов, существенная часть которых с подробными
доказательствами публикуется впервые. В первой главе излагаются
исторические сведения, необходимые обозначения, определения и
вспомогательные факты. Во второй главе изложена подробная теория
операторов преобразования типа Сонина и Пуассона. В третьей главе
изложена теория специального важного класса операторов
преобразования Бушмана---Эрдейи и их приложения. В четвёртой главе
рассматриваются новые весовые краевые задачи с операторами
преобразования типа Сонина и Пуассона. В пятой главе рассмотрены
приложения операторов преобразования типа Бушмана---Эрдейи
специального вида к новым краевым задачам для эллиптических
уравнений с существенными особенностями в решениях. В шестой главе
излагается универсальный композиционный метод построения
операторов преобразования и приведены его приложения. В
заключительной седьмой главе изложены приложения теории операторов
преобразования к  задачам для дифференциальных уравнений с
переменными коэффициентами: к задаче о построении нового класса
операторов преобразования с точными оценками ядер для возмущённого
дифференциального уравнения с оператором Бесселя, а также к
специальным случаям известной задачи Е.\,М.~Ландиса об
экспоненциальных оценках роста решений для стационарного уравнения
Шрёдингера. В заключение книги приведён краткий биографический
очерк о Валерии Вячеславовиче Катрахове, а также подробная
библиография, содержащая~\ref{BibEnd} ссылок.
\end{abstract}

\begin{engabstract}
The main content of this book is composed from two doctoral
theses: by V.\,V.~Katrakhov (1989) and by S.\,M.~Sitnik (2016). In
our work, for the first time in the format of a monograph, we
systematically expound the theory of transmutation operators and
their applications to differential equations with singularities in
coefficients, in particular, with Bessel operators. Along with
detailed survey and bibliography on this theory, the book contains
original results of the authors. Significant part of these results
is published with detailed proofs for the first time. In the first
chapter, we give historical background, necessary notation,
definitions, and auxiliary facts. In the second chapter, we give
the detailed theory of Sonin and Poisson transmutations. In the
third chapter, we describe an important special class of the
Buschman---Erd\'elyi transmutations and their applications. In the
fourth chapter, we consider new weighted boundary-value problems
with Sonin and Poisson transmutations. In the fifth chapter, we
consider applications of the Buschman---Erd\'elyi transmutations
of special form to new boundary-value problems for elliptic
equations with significant singularities of solutions. In the
sixth chapter, we describe a universal compositional method for
construction of transmutations and its applications. In the
concluding seventh chapter, we consider applications of the theory
of transmutations to differential equations with variable
coefficients: namely, to the problem of construction of a new
class of transmutations with sharp estimates of kernels for
perturbed differential equations with the Bessel operator, and to
special cases of the well-known Landis problem on exponential
estimates of the rate of growth for solutions of the stationary
Schr\"odinger equation. The book is concluded with a brief
biographic essay about Valeriy V.~Katrakhov, as well as detailed
bibliography containing~\ref{BibEnd} references.
\end{engabstract}

\maketitle

\vspace{-3ex}
\tableofcontents

\chapter*{Предисловие редактора}
\addcontentsline{toc}{chapter}{Предисловие редактора}

Теория    сингулярных дифференциальных уравнений, содержащих
оператор Бесселя $$\displaystyle{B_\nu = {1\over x^\nu}{d\over d
x} \left(x^\nu{d\over d x}\right)},$$ и неразрывно связанная с ней
теория соответствующих весовых функциональных пространств,
бесспорно, относятся к тем математическим направлениям,
теоретическое и прикладное значение которых трудно переоценить. И,
как это и характерно для активно развивающихся направлений
переднего края науки, она полностью подтверждает известный афоризм
<<основоположники учебников не пишут>>. Действительно, за
последние четыре десятилетия прошлого века (в течение которых,
собственно, эта теория и создана в современном виде) опубликовано
значительное количество статей с абсолютно прорывными
результатами, выращена достойная плеяда высококвалифицированных
исследователей (достаточно отметить, что только докторских
диссертаций защищено как минимум пять), а вот монография издана
только одна "--- это фундаментальный труд основоположника этой
науки Ивана Александровича Киприянова <<Сингулярные эллиптические
краевые задачи>> (М., Наука, 1997).

Это совершенно закономерно для этапа бурного развития новой
научной области, однако в первые полтора десятилетия нынешнего
века эта наука понесла тяжелейшие утраты "--- ушли из жизни и
И.\,А.~Киприянов, и такие её классики, как В.\,В.~Катрахов и
В.\,З.~Мешков, и другие математики, упорно и талантливо
развивавшие это направление. Настал новый этап, на котором крайне
важно (разумеется, наряду с дальнейшим развитием науки) разбирать
архивы учителей, сберегать их научное наследие и делать его
доступным широкой математической общественности. Достойным
примером такой деятельности является предлагаемая вниманию
читателей монография "--- итог без преувеличения титанического
труда С.\,М.~Ситника по систематизации наследия своего учителя
В.\,В.~Катрахова.

Этот труд начал еще сам Катрахов (более десяти лет назад), и вот
теперь он, наконец, завершен. Очень бережная и кропотливая работа
второго автора по сохранению наследия первого заслуживает самой
высокой оценки (во многих местах текста даже сохранен неповторимый
стиль изложения Катрахова, столь привлекавший его читателей и
студентов), однако значение этого труда выходит за рамки простого
разбора архива учителя: Ситник включил в монографию не только
результаты Катрахова (в том числе "--- те, которые не были
опубликованы при его жизни) и не только результаты, полученные в
соавторстве с ним, но и собственные результаты, совершенно
самостоятельно полученные за те почти десять лет, что Катрахова
нет с нами.

Указанные результаты образуют стройную и обоснованно
структурированную книгу, посвященную одному из самых эффективных и
(по-прежнему) перспективных методов изучения сингулярных задач
киприяновского направления "--- методу операторов преобразования
(сплетающих операторов). Безусловно, сфера применения этого метода
отнюдь не исчерпывается сингулярными дифференциальными
уравнениями. Однако, как показано в последние три десятилетия (в
первую очередь "--- усилиями авторов монографии,  предлагаемой
вниманию читателей), для уравнений с  операторами Бесселя он дает
так много, что есть все основания считать операторы преобразования
такой же органичной частью сингулярной теории, что и, скажем,
весовые функциональные пространства Киприянова.

Выход в свет этой книги, несомненно, станет событием
математической жизни.

\vspace{5mm}

\begin{flushright}
{\it А.\,Б.~Муравник}
\end{flushright}

\newpage

\chapter*{Предисловие автора}
\addcontentsline{toc}{chapter}{Предисловие автора}

   Теория операторов преобразования является хорошо разработанным самостоятельным разделом математики.
Значительный вклад в эту теорию и её приложения к дифференциальным
уравнениям с частными производными внесли работы воронежского
математика Валерия Вячеславовича Катрахова (1949--2010), ученика
Ивана Александровича Киприянова.

К числу важных результатов В.\,В.~Катрахова следует отнести
исследование весовых и спектральных задач для дифференциальных
уравнений и систем с операторами Бесселя с использованием техники
операторов преобразования. Им также совместно с И.\,А.~Киприяновым
были введены и изучены уравнения с псевдодифференциальными
операторами, которые определялись через преобразование Ханкеля при
помощи операторов преобразования Сонина и Пуассона.

Особо следует выделить введённый В.\,В.~Катраховым  новый класс
краевых задач для уравнения Пуассона, решения которого могут иметь
существенные особенности. На основе введённого им нового класса
операторов преобразования, получаемых из известных операторов
Сонина и Пуассона композициями с дробными интегралами
Римана---Лиувилля,   В.\,В.~Катраховым были введены специальные
функциональные пространства, содержащие функции с существенными
особенностями, доказаны для них теоремы вложения, прямые и
обратные теоремы о следах. Для функций без особенностей указанные
пространства сводятся к пространствам С.\,Л.~Соболева, таким
образом являясь их прямыми обобщениями. Для корректности задач с
существенными особенностями В.\,В.~Катраховым было предложено
новое естественное краевое условие во внутренней точке области,
которое заключается в задании предела свёртки решения с некоторым
сглаживающим ядром типа ядра Пуассона. Мы предлагаем называть это
новое краевое условие <<$K$-следом>> в честь В.\,В.~Катрахова,
который ввёл это  условие и подробно изучил краевые задачи с ним.
В терминах <<$K$-следа>> получается полная характеризация решений
уравнения Лапласа с внутренней особой точкой, в том числе для
решений с существенными особенностями в этой точке. Для данной
задачи в указанных функциональных пространствах В.\,В.~Катраховым
была доказана корректность постановки, включая существование и
единственность решения, априорные оценки. Этот результат обобщает
теоремы о разрешимости эллиптических уравнений в классах
С.\,Л.~Соболева для гладких решений без особенностей. Кроме того,
в последующих работах В.\,В.~Катрахова с соавторами  были
рассмотрены обобщения новых краевых задач для уравнений с
операторами Бесселя и сингулярным потенциалом, для областей в
пространствах Лобачевского и случая угловых точек на границе
области. Краткое перечисление основных результатов
В.\,В.~Катрахова также приведено в~\cite{S92, S95}.

Данная монография составлена из результатов, вошедших в докторские
диссертации  В.\,В.~Катрахова (1989~г.) и С.\,М.~Ситника (2016~г.)
(см., соответственно,~\cite{KatDis} и~\cite{SitDis}). Результаты
второго автора (ученика) развивают результаты первого автора
(учителя). Надеюсь, что публикуемая книга будет способствовать
более широкой известности результатов В.\,В.~Катрахова, которые
представляют существенный интерес для теории вырождающихся и
сингулярных дифференциальных уравнений, а также их разработке в
русле идей и методов теории операторов преобразования. Кроме того,
в книге отражён вклад Ивана Александровича Киприянова и созданной
им Воронежской математической школы по сингулярным и вырождающимся
дифференциальным уравнениям в развитие теории дифференциальных
уравнений и теории функций.

\vspace{5mm}

\begin{flushright}
\it С.\,М.~Ситник. Белгород---Воронеж, 2018.
\end{flushright}

\newpage

\chapter{Введение}\label{ch1}
\section{Исторические сведения и краткое содержание книги}\label{sec1}

В последние десятилетия в связи с потребностями приложений возрос
интерес к сингулярным и вырождающимся  эллиптическим краевым
задачам для дифференциальных уравнений в частных производных.
Основополагающей здесь была известная работа
М.\,В.~Келдыша~\cite{Kel}, в которой на примере уравнения второго
порядка со степенным вырождением были выявлены основные
особенности постановки краевых условий для таких уравнений. Было
показано, что при одних соотношениях между коэффициентами на
характеристической части границы нужно ставить условие Дирихле
(задача $D$), при других соотношениях это условие должно быть
заменено на условие ограниченности решения (задача $E$). В
последнем случае по терминологии С.\,М.~Никольского и
П.\,И.~Лизоркина принято говорить о сильном вырождении
соответствующей задачи. Аналог краевых задач $D$ и $E$ был изучен
многими авторами для весьма общих эллиптических уравнений. Мы
сошлемся в связи с этим на работы С.\,М.~Никольского~\cite{67},
С.\,М.~Никольского, П.\,И.~Лизоркина~\cite{57, 58, 59, 68},
Л.\,Д.~Кудрявцева~\cite{48}, И.\,А.~Киприянова~\cite{39, 40, 41},
М.\,И.~Вишика, В.\,В.~Грушина~\cite{15}, Л.\,А.~Ройтберга,
З.\,Г.~Шефтеля~\cite{72} и на монографии А.\,В.~Бицадзе~\cite{9,
Bitz2}, М.\,М.~Смирнова~\cite{76}, Х.~Трибеля~\cite{Trib1}, где
имеются подробные библиографии по данному вопросу. В случае
сильного вырождения уравнение кроме ограниченных решений имеет и
неограниченные (сингулярные) вблизи характеристической части
границы решения. А.\,В.~Бицадзе (см.~\cite{9, Bitz2}) предложил в
этом случае задавать на характеристической части не само решение
или его нормальную производную, а их произведения с заранее
подобранной весовой функцией. Такие задачи мы называем весовыми.
Весовая задача Коши для гиперболических уравнений изучалась
Ж.\,Л.~Лионсом~\cite{Lio1}, Р.~Кэрролом, Р.~Шоуолтером~\cite{CSh}
(см. также библиографию в этих книгах), а весовая краевая задача
для эллиптических уравнений "--- А.\,А.~Вашариным,
П.\,И.~Лизоркиным~\cite{12}, Г.\,Н.~Яковлевым~\cite{87},
А.\,И.~Янушаускасом~\cite{Yanu} и многими другими авторами (см.
монографию С.\,Г.~Самко, А.\,А.~Килбаса,
О.\,М.~Маричева~\cite{SKM} и библиографию в ней).

В другой ситуации краевые задачи с сильным вырождением возникают в
теории особых точек решений (даже регулярных) эллиптических
уравнений. Классические результаты в этой теории получены
Н.~Винером~\cite{14} и М.\,В.~Келдышем~\cite{37}. Более
современное состояние вопроса отражено в монографиях
Н.\,С.~Ландкофа~\cite{54}, Е.\,М.~Ландиса~\cite{53} и
И.\,А.~Шишмарева~\cite{85}, а также в работах
А.\,А.~Новрузова~\cite{69}, А.\,И.~Ибрагимова~\cite{18}. Во всех
указанных работах рассматривались в основном вопросы, связанные с
нахождением условий, обеспечивающих устранимость особенностей. Для
уравнений математической физики соответствующие факты приведены в
книге А.\,Н.~Тихонова, А.\,А.~Самарского~\cite{79}. При этом
техника вырождающихся уравнений здесь не использовалась, так как
постановка задачи приводит по сути дела к устранению особенности
решения, а постановка весовых краевых условий без дополнительных
ограничений на рост решения невозможна.

Ещё раз подчеркнём, что первой фундаментальной работой, с которой
начался отсчёт  изучения вырождающихся и сингулярных
дифференциальных уравнений в частных производных с переменными
коэффициентами, является статья М.\,В.~Келдыша~\cite{Kel} (задачи
E и D). Мы понимаем сингулярные и вырождающиеся уравнения в том
смысле, как они определены во введении в классической
монографии~\cite{CSh}. Сингулярные, вырождающиеся и тесно
связанные с ними уравнения смешанного типа можно объединить в один
класс "--- \textit{неклассические уравнения математической
физики}, этот термин был предложен В.\,Н.~Враговым~\cite{Vrag}.
Теория для уравнений указанных типов разрабатывалась многими
математиками, в том числе С.~Геллерстедом, М.~Проттером,
В.~Радулеску, Ф.~Трикоми, Г.~Фикерой, Е.~Хольмгреном, М.~Чибрарио,
С.\,А.~Алдашевым, А.\,А.~Андреевым, Ф.\,Т.~Барановским,
А.\,В.~Бицадзе, Б.\,А.~Бубновым, А.\,А.~Вашариным, И.\,Н.~Векуа,
М.\,И.~Вишиком, В.\,Ф.~Волкодавовым, В.\,Н.~Враговым,
В.\,П.~Глушко, Г.\,В.~Джаяни, И.\,Е.~Егоровым, В.\,И.~Жегаловым,
А.\,Н.~Зарубиным, В.\,А.~Ильиным, А.\,М.~Ильиным,
Т.\,Ш.~Кальменовым, М.\,Б.~Капилевичем, А.\,А.~Килбасом,
С.\,Б.~Климентовым, А.\,И.~Кожановым, Л.\,Д.~Кудрявцевым,
П.\,И.~Лизоркиным, О.\,И.~Маричевым, Л.\,Г.~Михайловым,
С.\,Г.~Михлиным, Е.\,И.~Моисеевым, А.\,М.~Нахушевым,
Н.\,Я.~Николаевым, С.\,М.~Никольским, О.\,А.~Олейник,
Л.\,С.~Парасюком, С.\,В.~Поповым, С.\,П.~Пулькиным,
Л.\,С.~Пулькиной, С.\,Г.~Пятковым, А.\,Б.~Расуловым,
О.\,А.~Репиным, К.\,Б.~Сабитовым, М.\,С.~Салахитдиновым,
А.\,Л.~Скубачевским, М.\,М.~Смирновым, А.\,П.~Солдатовым,
С.~Руткаускасом, С.\,А.~Терсеновым, З.\,Д.~Усмановым,
В.\,Евг.~Фёдоровым, В.\,Евс.~Фёдоровым, Ф.\,И.~Франклем,
Л.\,И.~Чибриковой, А.\,И.~Янушаускасом
 и многими другими. Из многочисленной литературы укажем лишь ссылки~\cite{Bitz2,EgFed,EPP,Glushko,Kip1,MaKiRe,Moi,Nah4,Pul2,Pyat,Rep,Smi,Sku1,Sku2,Sku3,Sku4,Ter1,
 Tricomi1,Yanu,CSh,Rad1}.

Специальный класс составляют функционально-дифференциальные
уравнения, или, как их называют в частных случаях,
дифференциально-разностные уравнения, дифференциальные уравнения с
запаздывающим или отклоняющимся аргументом. Их теория была начата
в работах А.\,Д.~Мышкиса, Дж.~Хейла, Г.\,А.~Каменского,
Л.\,Э.~Эльсгольца и изложена в
монографиях~\cite{Mys,Hale,Els,ElNor,GKE,ZKNE}, существенные
результаты в этом направлении получены Н.\,В.~Азбелевым,
А.\,Л.~Скубачевским, В.\,П.~Максимовым,  Л.\,Ф.~Рахматуллиной,
Л.\,Е.~Россовским, А.\,Б.~Муравником,  укажем лишь
ссылки~\cite{AMR,Sku1,Sku2,KaSk,Skub,Ross1,Ros3,VaRo,Mur}. Отметим
использование функционально-дифференциальных уравнений в работах
В.\,А.~Рвачёва и В.\,Л.~Рвачёва для определения специального
класса атомарных функций, которые нашли важные приложения в теории
приближений~\cite{Rv1, Rv2, Rv3}, а также в многочисленных
прикладных задачах~\cite{Kra1, Kra2}. Принципиальной отличительной
чертой функционально-дифференциальных уравнений является наличие
даже у модельных задач финитных решений, например, таких, как
атомарные функции Рвачёвых, что невозможно по отдельности ни для
чисто дифференциальных, ни для чисто разностных уравнений. К
функционально-дифференциальным относится также важный класс
уравнений с операторами Дункла, которые родственны
дифференциальным операторам Бесселя и возникают на стыке теории
групп и симметрий в них (группы Коксетера), дифференциальных
уравнений и интегральных преобразований,  а также квантовой физики
и кристаллографии, см.~\cite{Me1,DHS, Dun1, Dun2, Dun3, Gal1,Gal2,
Ros1,Ros2, Tri6}. К классу функционально-дифференциальных
уравнений также относятся задачи с инволютивным или Карлемановским
сдвигом, см.~\cite{KaSa,Lit1}.

Особо отметим один класс уравнений с частными производными с
особенностями, типичным представителем которого является
$B$-эллипти\-ческое уравнение с операторами Бесселя по каждой
переменной вида
\begin{equation}
\label{Bes2} \sum\limits_{k=1}^{n}B_{\nu,x_k}u(x_1,\dots, x_n)=f,
\end{equation}
аналогично рассматриваются $B$-гиперболические и
$B$-параболические уравнения, эта удобная терминология была
введена И.\,А.~Киприяновым~\cite{Kip1}. Изучение этого класса
уравнений было начато в работах Эйлера, Пуассона, Дарбу,
продолжено в теории обобщённого осесимметрического потенциала
А.~Вайнштейна (теория GASPT "--- Generalized Axially Symmetric
Potential Theory,~\cite{Wei1, Wei2, Wei3}), Л.~Берса~\cite{Bers1,
Bers2, Bers3} и в трудах математиков С.\,А.~Алдашева,
И.\,Е.~Егорова, Я.\,И.~Житомирского, М.\,Б.~Капилевича,
Ш.\,Т.~Каримова, Э.\,Т.~Каримова, А.\,А.~Килбаса,
Л.\,Д.~Кудрявцева, П.\,И.~Лизоркина, О.\,И.~Маричева,
М.\,И.~Матийчука, Л.\,Г.~Михайлова, М.\,Н.~Олевского,
С.\,П.~Пулькина, М.\,М.~Смирнова, С.\,А.~Терсенова, А.~Хасанова,
Хе Кан Чера, А.\,И.~Янушаускаса   и других. Важность уравнений из
этих классов определяется также их использованием в приложениях к
задачам теории осесимметрического потенциала~\cite{Wei1, Wei2,
Wei3}, уравнениям Эйлера---Пуассона---Дарбу
(ЭПД)~\cite{VoNi,Dza1}, преобразованию Радона и
томографии~\cite{Nat1, Rub3, Lud, Rub1, Rub2, Rub4}, газодинамики
и акустики~\cite{Bers1, Bers2, Bers3}, теории струй в
гидродинамике (М.\,И.~Гуревич~\cite{Gur}), линеаризованным
уравнениям Максвелла---Эйнштейна (А.\,В.~Бицадзе,
В.\,И.~Пашковский~\cite{Bitz1,Bitz12}), механике, теории упругости
и пластичности~\cite{Dza2} и многим другим.

В определённом приближении можно сказать, что указанные три класса
дифференциальных уравнений по терминологии И.\,А.~Киприянова в
своё время были рассмотрены в трёх известных монографиях:
$B$-эллиптические уравнения в монографии
И.\,А.~Киприянова~\cite{Kip1}, $B$-гиперболические уравнения в
монографии Р.~Кэрролла и Р.~Шоуолтера~\cite{CSh},
$B$-параболические уравнения в монографии
М.\,И.~Матийчука~\cite{Mat1}. Разумеется, в указанных книгах
рассматриваются и многие другие вопросы.

Наиболее полно весь круг вопросов для  уравнений с операторами
Бесселя был изучен воронежским математиком Иваном Александровичем
Киприяновым и его учениками Л.\,А.~Ивановым, В.\,В.~Катраховым,
А.\,В.~Рыжковым, А.\,А.~Азиевым, В.\,П.~Архиповым,
А.\,Н.~Байдаковым, Б.\,М.~Богачёвым, А.\,Л.~Бродским,
Г.\,А.~Виноградовой, В.\,А.~Зайцевым, Ю.\,В.~Засориным,
Г.\,М.~Каганом, А.\,А.~Катраховой, Н.\,И.~Киприяновой,
В.\,И.~Кононенко, М.\,И.~Ключанцевым, А.\,А.~Куликовым,
А.\,А.~Лариным, М.\,А.~Лейзиным, В.\,М.~Ляпиным, Л.\,Н.~Ляховым,
А.\,Б.~Муравником, И.\,П.~Половинкиным, А.\,Ю.~Сазоновым,
С.\,М.~Ситником, П.\,С.~Украинским, В.\,П.~Шацким,
Э.\,Л.~Шишкиной, В.\,Я.~Ярославцевой; основные результаты этого
направления представлены в~\cite{Kip1}.

\begin{figure}[H]
\centering
\begin{subfigure}[t]{0.45\textwidth}
\centering
\includegraphics[width=0.7\textwidth]{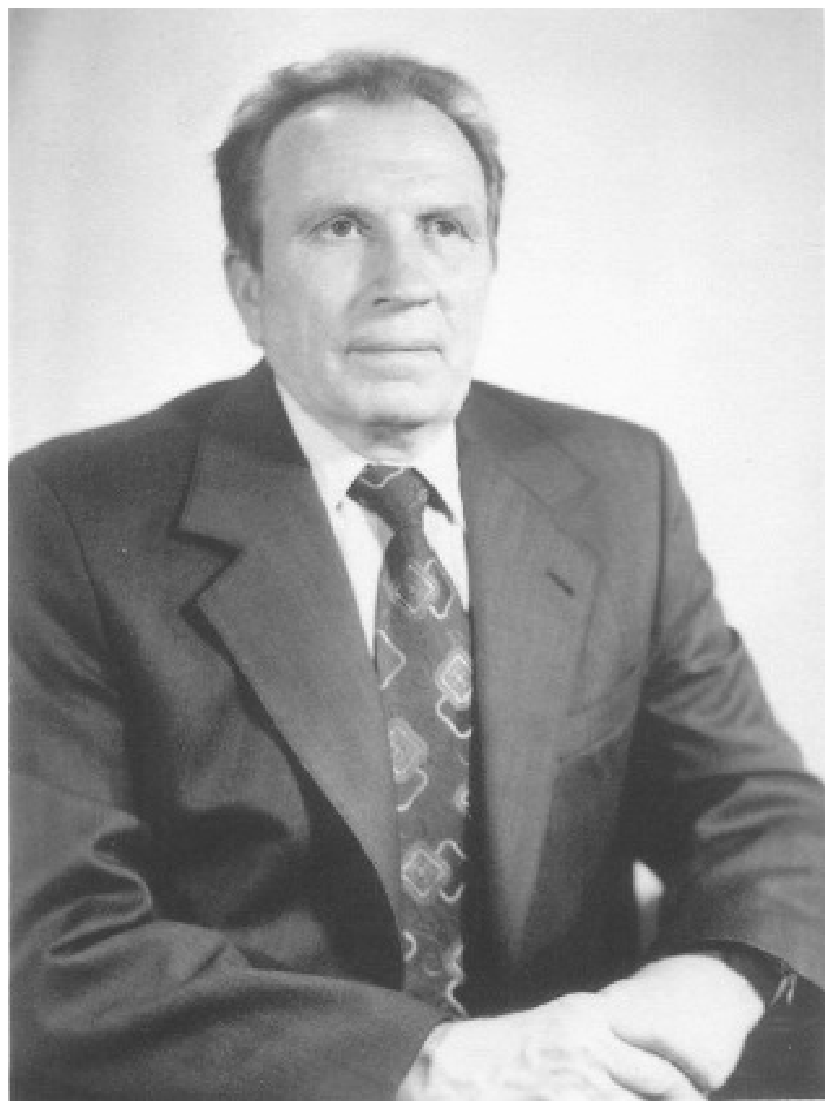}
\subcaption{\hfil Иван Александрович Киприянов\hfil\\
\hphantom{\hspace{2.9cm}}(1923--2001)} \label{pic1}
\end{subfigure}
\qquad
\begin{subfigure}[t]{0.45\textwidth}
\centering
\includegraphics[height=0.936\textwidth]{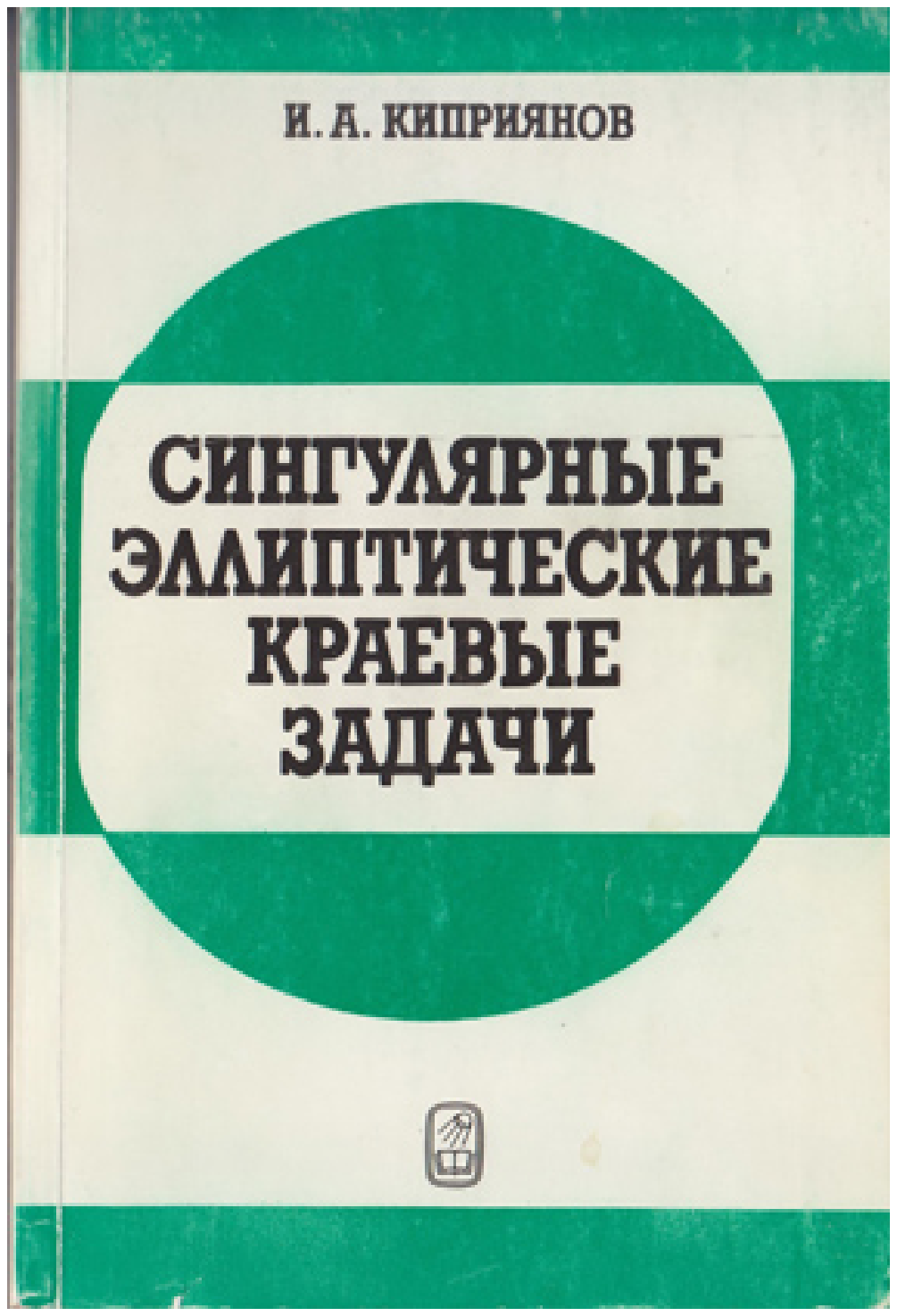}
\subcaption{}\label{pic2}
\end{subfigure}
\end{figure}

Для описания классов решений соответствующих уравнений
И.\,А.~Киприяновым были введены и изучены функциональные
пространства~\cite{Kip2},  позднее названные его именем (см.
монографии Х.~Трибеля~\cite{Trib1}, Л.\,Д.~Кудрявцева и
С.\,М.~Никольского~\cite{KN}, в которых пространствам Киприянова
посвящены отдельные параграфы). Интересные результаты по изучению
псевдодифференциальных операторов на основе операторов
преобразования Сонина---Пуассона---Дельсарта были получены
В.\,В.~Катраховым~\cite{KiKa1}, они также изложены в специально
переработанном Р.~Кэрролом виде в отдельной главе в
монографии~\cite{Car2}. В работах И.\,А.~Киприянова и
Л.\,А.~Иванова  было начато исследование $B$-гиперболических
уравнений, распадающихся на множители~\cite{Iva1,Iva2,KipIv1}, а
также изучены $B$-гиперболические уравнения в римановых
пространствах и пространстве Лобачевского~\cite{42, KipIv2,
KipIv3, KipIv4}. Фундаментальные решения для $B$-эллиптических
уравнений были построены и изучены в работах И.\,А.~Киприянова,
В.\,И.~Кононенко и А.\,А.~Куликова~\cite{KipKo1, KipKo2, KipKo3,
KipKu1}. Л.\,Н.~Ляховым были рассмотрены весовые сферические
функции, обобщённые потенциалы и гиперсингулярные интегралы,
задачи для псевдодифференциальных уравнений с операторами Бесселя
и некоторые обобщения преобразования Радона~\cite{Lyah1, Lyah2,
Lyah3}. В работах А.\,Н.~Байдакова были рассмотрены априорные
оценки гёльдеровых норм решений квазилинейных $B$-эллиптических
уравнений~\cite{Baid}, в работах В.\,А.~Зайцева были изучены
распадающиеся на множители уравнения с особенностями в
коэффициентах, а также предложен оригинальный способ исследования
для них принципа Гюйгенса, что позволило по-новому осмыслить
понятия <<слабая лакуна>>, <<принцип Гюйгенса>> и <<диффузия
волн>>, введя понятие <<принципа Гюйгенса порядка $(p,q)$>>, а
также  построить примеры уравнений, для которых выполняется
принцип Гюйгенса в чётномерных пространствах, что до этих
результатов считалось невозможным~\cite{Zai1,Zai2}. В работах
В.\,П.~Шацкого были рассмотрены сингулярные и вырождающиеся
системы первого порядка с особенностями~\cite{Sha1, Sha2, Sha3}; в
работах А.\,А.~Ларина изучены различные задачи для преобразования
Фурье---Бесселя~\cite{Lar1} и сингулярных дифференциальных
уравнений~\cite{Lar2, Lar3, Lar4, Lar5, Lar6}; в работах
А.\,Б.~Муравника рассмотрены свойства преобразования
Фурье---Бесселя~\cite{Mur3, Mur4, Mur5, Mur6, Mur7}, сингулярные
дифференциальные уравнения~\cite{Mur2} и вопросы их
стабилизации~\cite{Mur8,Mur9}. В.\,В.~Катраховым  и
С.\,М.~Ситником рассмотрены краевые задачи с $K$-следом  для
стационарного уравнения Шрёдингера с сингулярным
потенциалом~\cite{S1}; С.\,М.~Ситником разработан общий
композиционный метод для построения различных классов операторов
преобразования~\cite{S46,14,SitDis,S42,FJSS},  введены новые
классы операторов преобразования типа
Бушмана---Эрдейи~\cite{S3,S4,S66,S6,S46,S103,S400,SitDis,S42,S38,S401,S402}
и Векуа---Эрдейи---Лаундеса~\cite{S125,S46,S14,S42,S46}, с
использованием полученных оценок для норм операторов
Бушмана---Эрдейи установлена в простейшем одномерном случае
эквивалентность норм пространств И.\,А.~Киприянова и весовых
пространств С.\,Л.~Соболева, а также получены некоторые обобщения
неравенства Харди, методом операторов преобразования изучены
специальные случаи известной задачи Е.\,М.~Ландиса о предельной
скорости убывания решений стационарного уравнения
Шрёдингера~\cite{S75,S71,S3,S42}, модифицирован известный метод
построения операторов преобразования для возмущённого уравнения
Бесселя~\cite{S63,S19,S43,S42,SitDis}, построены в явном
интегральном виде дробные степени оператора Бесселя и изучены их
основные свойства~\cite{S135,S133,S127,S123,S18,S700,SS,S42},
совместно с Э.\,Л.~Шишкиной рассмотрены приложения метода
операторов преобразования к уравнениям
Эйлера---Пуассона---Дарбу~\cite{SSfiz,ShiE3}.

Кроме того, в работах И.\,А.~Киприянова и математиков его школы
были получены наиболее глубокие результаты для важного класса
сингулярных дифференциальных уравнений с параметром. Как уже было
отмечено выше, в~\cite{Kip1, 39} И.\,А.~Киприяновым были введены
функциональные весовые пространства путём замыкания гладких и
чётных по нормальному к границе направлению функций. Поскольку в
процессе замыкания сохраняются не все следы (это характерный
момент теории весовых пространств), то рассматриваемая в этих
пространствах краевая задача~\cite{39, 40, 41} аналогична задаче
типа $E,$ причём с однородными краевыми условиями чётности. Кроме
того, применение в качестве основного метода преобразования
Фурье---Бесселя обусловило вещественность упомянутого параметра.
Использованный нами другой подход позволил освободиться от всех
этих ограничений и рассмотреть случай произвольного комплексного
параметра, что не является формальным обобщением ни с точки зрения
постановки краевых задач, ни с точки зрения приложений к физике. В
отличие от~\cite{40,41} мы изучаем общие неоднородные весовые
краевые условия и, следовательно, рассматриваем сингулярные
решения в новых, более широких функциональных пространствах. Если
же применить наш метод при указанных выше ограничениях, то
получающиеся при этом результаты совпадают, по-существу, с
результатами И.\,А.~Киприянова. Укажем в этом плане также на
работу И.\,А.~Киприянова, М.\,И.~Ключанцева~\cite{45}, в которой
рассмотрены в случае вещественного параметра общие невесовые
краевые условия, которые, однако, ставятся на нехарактеристической
части границы, а в противном случае ставятся условия чётности.

Задачи для операторно-дифференциальных (абстрактных)
дифференциальных уравнений  с операторами Бесселя, берущие начало
в известной монографии~\cite{CSh}, начали рассматриваться
воронежскими математиками в работе А.\,В.~Глушака,
В.\,И.~Кононенко, С.\,Д.~Шмулевича~\cite{Glu1}, и затем
последовательно изучались в работах А.\,В.~Глушака~\cite{Glu2,
Glu3, Glu4, Glu55, Glu5, Glu6, Glu7, Glu8, Glu9, Glu10, Glu11,
Glu12, Glu13}. В указанных работах А.\,В.~Глушака для абстрактных
дифференциальных уравнений рассматривались вопросы решения
уравнений Эйлера---Пуассона---Дарбу, Бесселя---Струве,
стабилизации решений, операторные косинус функции, функции Бесселя
и Лежандра.

Таким образом, к настоящему времени сложилась такая ситуация, что
краевые задачи типа $E$ глубоко и разносторонне изучены, весовые
же краевые задачи рассмотрены лишь для довольно простых уравнений,
причём только второго порядка. Другие постановки краевых условий
для уравнений с сильным вырождением, которые заменяли бы весовые в
случае невозможности постановки последних, с должной глубиной до
сих пор не изучались. Такое положение, на наш взгляд, можно
объяснить двумя причинами. Во-первых, обычные методы теории
уравнений с частными производными мало приспособлены для решения
указанных задач, а некоторые из них (например, вариационный метод
для уравнений, не являющихся сильно эллиптическими) в принципе не
применимы. Во-вторых, отсутствовала теория функциональных
пространств, для которых имели бы место соответствующие теоремы о
следах. Известные весовые функциональные пространства
И.\,А.~Киприянова (см. монографии С.\,М.~Никольского~\cite{66},
Л.\,Д.~Кудрявцева~\cite{48}), обычно применяемые для вырождающихся
уравнений, обладают тем свойством, что порядок особенности функций
из этих пространств зависит от показателя гладкости и убывает с
возрастанием гладкости. Последнее же обстоятельство неприемлемо в
рассматриваемой нами ситуации. Поэтому, постановка и изучение
новых краевых задач для уравнений с сильным вырождением в
соответствующих им функциональных пространствах, а также создание
эффективных методов их решения являются актуальными.

Кратко изложим содержание книги по главам с некоторыми комментариями и пояснениями.

В первой главе, имеющей вводный характер,  приводится план работы, краткий исторический очерк развития теории операторов преобразования, собраны основные сведения о специальных функциях, интегральных преобразованиях и функциональных пространствах.

Во второй главе разработан аппарат операторов преобразования для
одномерного и многомерного случаев, который используется в
дальнейших главах (подробнее см. следующий параграф этой главы).
Своим появлением метод операторов преобразования обязан работам
Ж.~Дельсарта, Ж.-Л.~Лионса~\cite{Lio1,Lio2},
Б.\,М.~Левитана~\cite{Lev7}, который дал название операторов
Сонина и Пуассона известному классу операторов преобразования.
Операторы Сонина и Пуассона нашли широкое применение в теории
сингулярных гиперболических уравнений типа
Эйлера---Пуассона---Дарбу (см. цитированные монографии
Ж.-Л.~Лионса, Р.~Кэррола и Р.~Шоуолтера). Ряд важных результатов в
этом направлении для гиперболических уравнений получен
И.\,А.~Киприяновым и Л.\,А.~Ивановым~\cite{42}. Заметим, что
отдельные элементы этого метода в связи с применением метода
сферических средних использовались еще при изучении задачи Коши
для классического волнового уравнения (см., например, книгу
Р.~Куранта~\cite{49}). Соответствующие операторы преобразования
применялись Г.\,Н.~Положим~\cite{76} в теории обобщённых
аналитических функций. Для других видов сингулярных уравнений
операторы преобразования были построены
М.\,И.~Ключанцевым~\cite{Kly1} (эти результаты являются частными
случаями полученных ранее И.~Димовски операторов преобразования
для общего случая гипербесселевых операторов,
см.~\cite{DK2,Dim,Kir1}) и В.\,Я.~Ярославцевой~\cite{Yar1} (для
гиперболических аналогов операторов Сонина и Пуассона в
тригонометрической форме). В другой ситуации метод операторов
преобразования был использован в спектральной теории уравнений
Шрёдингера (см. монографии В.\,А.~Марченко~\cite{Mar2}, К.~Шадана,
П.~Сабатье~\cite{ShSa} и библиографию в них).

Подробные результаты для специального важного класса операторов
преобразования типа Сонина и Пуассона изложены во второй главе. Мы
используем исторически более точное название "--- операторы
Сонина---Пуассона---Дельсарта (СПД). Такие операторы
преобразования как интегральные  относятся к классу дробных
интегралов Эрдейи---Кобера.

Почти одновременно с первыми работами Ж.~Дельсарта,
Х.~Кобером~\cite{47} и А.~Эрдейи~\cite{86} были введены другие
операторы преобразования. История развития теории операторов
Эрдейи---Кобера с достаточной полнотой изложена в
монографии~\cite{SKM}. Круг приложения этих и связанных с ними
операторов ограничивался в основном получением явных представлений
классических решений краевых задач для уравнений второго порядка
типа Эйлера---Пуассона---Дарбу.

В этой книге метод операторов преобразования получает дальнейшее
развитие. В первую очередь это связано с построением на их основе
новых функциональных пространств типа пространств С.\,Л.~Соболева.
Ранее подобные конструкции с операторами преобразования не
приводились, причём применение известных операторов преобразования
для этих целей оказалось невозможным. В связи с этим в книге
введены новые операторы преобразования, отличающиеся от изученных
интегральных операторов на конечных промежутках  Р.~Бушманом,
А.~Эрдейи и другими. Некоторые новые результаты получены и для
операторов Эрдейи---Кобера, например, соответствующие утверждения
о весовых следах.

Аналогичные определённые пробелы оставались в вопросах корректной
постановки краевых задач, построении операторов преобразования и
подборе подходящих функциональных пространств  и в многомерном
случае. Оператор Радона (подробно изученный в монографиях
И.\,М.~Гельфанда, М.\,И.~Граева, Н.\,Я.~Виленкина~\cite{17} и
С.~Хелгасона~\cite{82}) преобразует многомерный оператор Лапласа в
одномерный оператор двукратного дифференцирования по радиальной
переменной и является, к сожалению, в определенном смысле
сглаживающим. В результате <<достоинства>> прямого преобразования
перечёркиваются принципиальными <<недостатками>> обратного, в чём,
по существу и заключаются основные трудности в теории обращения
преобразования Радона и его численных реализациях. Для
определённого решения этой проблемы в книге подобран погашающий
эту сглаживаемость множитель (дифференциальный оператор в
нечётномерном пространстве и лиувиллевский "--- в чётномерном), и
так сконструированный оператор преобразования применен для
построения новых функциональных пространств. Заметим попутно, что
используемая  техника позволяет обойтись без использования
проекционных пространств, обычно возникающих в теории
преобразования Радона.

Во второй главе изучаются в одномерном случае отображения
введёнными операторами преобразования известных функциональных
пространств, а также строятся новые пространства $H_{\nu}^s
\lr{E_{+}^1},$ даны обобщения этих пространств на многомерный
случай. Здесь же, пользуясь одномерной спецификой, мы в ряде
случаев получаем точные постоянные в оценках. В результатах такого
рода мы будем нуждаться в последующих главах.

В книгу включены лишь результаты по операторам преобразования и
краевым задачам для дифференциальных уравнений. Тем же методом
операторов преобразования были рассмотрены и другие вопросы, в
частности, построена~\cite{25, 26, 33, KiKa1} алгебра сингулярных
псевдодифференциальных операторов и теория комплексных степеней
сингулярных эллиптических операторов, обобщающая известные
результаты В.\,А.~Ильина~\cite{20}, Ш.\,А.~Алимова~\cite{2},
Р.\,Т.~Сили~\cite{75} на сингулярный случай.

Касаясь теории введённых в книге двух классов функциональных
пространств (гильбертовых и счётно-нормируемых), заметим, что
ранее такого типа пространства, для которых справедливы теоремы о
весовых или нелокальных следах, не изучались. При этом наши
пространства нельзя отнести к классу весовых, поскольку при
определении их норм использованы интегро-дифференциальные
операторы. В первый из указанных классов вложены пространства,
введённые И.\,А.~Киприяновым~\cite{Kip1}, а во второй "---
пространства С.\,Л.~Соболева~\cite{77}.

В третьей главе подробно рассматриваются операторы преобразования
из другого важного класса "--- Бушмана---Эрдейи и их модификации.
Они содержат все ранее известные классы операторов преобразования
и имеют многочисленные приложения. Введены операторы
преобразования Бушмана---Эрдейи различных классов: первого и
второго родов, нулевого порядка гладкости, унитарные операторы
преобразования Сонина---Катрахова и Пуассона---Катрахова. Среди
приложений этих классов операторов преобразования рассмотрены
свойства унитарности и оценки норм в лебеговых пространствах на
полуоси, приложения к лемме Копсона и постановке задач Коши с
нелокальными начальными условиями, связи с классическими
операторами Харди, задача А.\,В.~Бицадзе и В.\,И.~Пашковского для
уравнений Максвелла---Эйнштейна, доказательство в одномерном
случае эквивалентности норм в пространствах И.\,А.~Киприянова и
С.\,Л.~Соболева со степенным весом.

В четвёртой главе излагаются основные результаты  по весовым
краевым задачам. Для уравнений высших порядков весовые краевые
условия поставлены в данной книге впервые, эти результаты
приведены в четвёртой главе. Основное их отличие от результатов
М.\,В.~Келдыша, С.\,М.~Никольского, А.\,В.~Бицадзе,
Л.\,Д.~Кудрявцева, П.\,И.~Лизоркина, И.\,А.~Киприянова,
М.\,И.~Вишика, В.\,В.~Грушина и других авторов, изучавших задачу
$E$ и её аналоги, заключается как в виде краевых условий, так и в
их числе, поскольку у нас рассматривается число краевых условий,
равное половине порядка уравнения. Всё это позволяет рассмотреть в
отличие от указанных авторов вместе с регулярными и сингулярные
решения. В случае отсутствия особенностей у коэффициентов
уравнения весовая функция в краевых условиях исчезает, и  наши
результаты совпадают с классическими результатами общей теории
эллиптических краевых задач (см. С.~Агмон, А.~Дуглис,
Л.~Ниренберг~\cite{1}, Л.~Хермандер~\cite{83}, Ж.-Л.~Лионс,
Э.~Мадженес~\cite{62}, О.\,А.~Ладыженская,
Н.\,Н.~Уральцева~\cite{51}).

В первом пункте дается постановка весовой краевой задачи для
однородных операторов с постоянными коэффициентами, которая
методом операторов преобразования сводится к краевой задаче,
изученной в предыдущем параграфе. Здесь следует отметить, что
применяя операторы преобразования для уравнений с переменными
коэффициентами, мы приходим, как это показано в работах~\cite{26,
KiKa1}, вообще говоря, к псевдодифференциальным операторам, у
которых символ по двойственным переменным имеет особенность
высокого порядка на последней координатной гиперплоскости. Теория
таких псевдодифференциальных операторов до сих пор не разработана.
Поэтому в случае переменных коэффициентов мы пользуемся шаудеровой
техникой, связанной с <<замораживанием>> коэффициентов и
последующей склейкой с применением разбиения единиц. В
соответствии с этой методикой рассматривается краевая задача с
маломеняющимися коэффициентами. Здесь пришлось преодолеть ряд
технических трудностей как общего характера, так и связанных с
весовым характером краевых условий и  сингулярностью
коэффициентов.

Сравнивая представленные результаты  с соответствующими
результатами А.\,В.~Бицадзе, А.\,А.~Вашарина и П.\,И.~Лизоркина,
Г.\,Н.~Яковлева,  А.\,И.~Янушаускаса и другими результатами по
весовым краевым задачам для эллиптических уравнений второго
порядка, нетрудно увидеть, что и в этом случае многие из наших
результатов являются новыми. К ним, в частности, относятся
основной результат о нётеровости, сам вид краевых условий (и их
порядок), теорема о повышении гладкости и другие. Также замечено,
что в некоторой ситуации, которая, по-видимому, ранее не
рассматривалась, невозможна постановка весовой задачи типа
Дирихле, а корректной при этом будет весовая задача типа Неймана,
которая является в приводимых ниже физических примерах аналогом
условий излучения Зоммерфельда.

В пятой главе в книге изучается новая краевая задача для уравнения
Пуассона с особенностью в изолированных граничных точках. В этих
точках решению и правой части разрешается иметь произвольный рост.
В случае однородного уравнения (Лапласа) никаких априорных
ограничений на поведение решения вообще не накладывается, и
поэтому его особенность может иметь тип существенной особенности
аналитических функций. В такой общей ситуации все известные
постановки краевых условий оказываются неприменимыми.  В самом
деле, в цитированных работах Н.~Винера, М.\,В.~Келдыша и в других
работах этого направления изучались только устранимые особенности.
В работе Н.\,С.~Казаряна~\cite{21} у решения допускается наличие
изолированных особенностей, но эти особые точки не являются
изолированными граничными точками, кроме того порядок роста в них
решения не может превышать степенного. Некоторые  результаты в
рассматриваемом нами направлении получены в работах
В.\,М.~Ивакина~\cite{19}, Б.\,В.~Квядараса~\cite{35},
С.~Руткаускаса~\cite{Rut1, Rut2, Rut3} и ряде других.

В связи с только что указанной задачей в книге вводится
принципиально новое понятие (в определенном смысле нелокального)
следа для функции, имеющей, вообще говоря, особенность в
изолированных граничных точках, и строятся новые функциональные
пространства типа Фреше, для которых оказываются справедливы
соответствующие теоремы о следах. Основной результат здесь состоит
в доказательстве однозначной разрешимости поставленной краевой
задачи и устойчивости по Адамару в с смысле указанных пространств.
В терминах введённого следа дана также классификация изолированных
особенностей гармонических функций. Отметим, что тем же методом
аналогичные результаты для других уравнений получены в
работах~\cite{S1, 44}.

Свёрточный характер рассматриваемого краевого условия позволяет отнести данный класс задач как к нелокальным, так и к нагруженным с заданием специального следа решений на части границы.

Более конкретно, пусть $\Omega \subset E^n$ "---  ограниченная
область с гладкой границе $\pr \Omega.$ Не ограничивая общности,
предположим, что начало координат принадлежит области. Пусть
$\Omega_0 = \Omega \setminus 0.$ Рассмотрим в $\Omega_0$ уравнение
Пуассона
\begin{equation}
\Delta u = f (x), \  x \in \Omega_0,
\label{50}
\end{equation}
с краевым условием Дирихле на $\pr \Omega$
\begin{equation}
 u|_{\pr \Omega} = g (x), \  x \in \pr \Omega.
\label{51}
\end{equation}
Проблема, решаемая в пятой главе, состоит в постановке такого
краевого условия в точке $0$ и в выборе таких функциональных
пространств, в которых краевая задача была бы однозначно
разрешима. При этом на порядок особенности решения в точке $0$ не
накладывается никаких ограничений. Такая ситуация встречается в
классических задачах электростатики. Пусть, например,
пространственная область $\Omega_0$ свободна от зарядов ($f=0$) и
окружена заземленной поверхностью ($g=0$). Тогда функция $u$ будет
потенциалом электростатического поля, созданного помещенным в
точку $0$ заряженным объектом. Если таким объектом будет просто
единичный заряд, то  функция $u$ называется  в   этом случае
функцией Грина. Она имеет особенность вида $|x|^{2-n}$ при $n \geq
3$   и $\ln |x|$ при $n=2.$  При помещении   в точку $0$ диполя
порядок особенности функции $u$ увеличится на единицу   и   т.~д.
В общем случае, помещая в точку $0$ бесконечную комбинацию
мультиполей различных порядков, получим у функции $u$ особенность
бесконечного порядка. Более того, поведение функции $u$ в
окрестности нуля аналогично поведению аналитической функции в
окрестности существенно особой точки.

Известные функциональные пространства и известные постановки
краевых условий в исследуемой ситуации непригодны. Например,
функции с существенной в точке особенностью не входят даже в
классы распределений Л.~Шварца~\cite{16}. Уравнение~\eqref{50} в
сферических координатах $r=|x|,$ $\vartheta = \dfrac{x}{|x|}$
принимает вид
\begin{equation}
\frac{\pr^2 u}{\pr r^2} +  \frac{n-1}{r} \frac{\pr u}{\pr r} + \frac{1}{r^2} \Delta_{\Theta} u = f \lr{r, \vartheta},
\label{52}
\end{equation}
где $\Delta_{\Theta}$ "--- оператор Лапласа---Бельтрами на сфере.
Постановка краевой задачи $E$ для такого уравнения хотя и
возможна, однако приводит (при гладкой функции $f$) к устранимой
особенности решения. Это нетрудно установить классическими
приёмами (см.~\cite{79}). Постановка весовых краевых условий
невозможна, поскольку уравнение~\eqref{52}  хотя и относится к
классу известных вырождающихся уравнений, рассмотренных
М.\,В.~Келдышем, но в нашем случае параметр $p=2,$ а для
возможности постановки весовых краевых условий необходимо $p<2.$
Отображение точки $0$ в бесконечность преобразованием
Кельвина~\cite{79} также не меняет сути дела, поскольку внешние
краевые задачи такого рода не изучались. Наоборот, это
преобразование позволяет перенести полученные ниже результаты на
внешние краевые задачи, причём на порядок роста решения в
бесконечности не накладывается никаких ограничений.

В настоящей работе в точке $0 \in \pr \Omega_0$ для решения $u$
ставится  нелокальное краевое условие следующего вида:
\begin{equation}
\lim\limits_{r \to +0} r^{n-2} \int\limits_{\Theta} u \lr{r, \vartheta'} K_n \lr{r \vartheta, \vartheta'}  d \vartheta' = \Psi \lr{\vartheta'},
\label{53}
\end{equation}
где $n \geq 3$ и $K_n (x, y)$   является    ядром   Пуассона, ассоциированным с единичной сферой, то есть
\begin{equation*}
K_n (x, y) = \frac{\Gamma \lr{\frac{n}{2}}}{2 \pi^{\frac{n}{2}}} \frac{1- |x|^2}{ |x-y|^n}, \  x,y \in E^n.
\end{equation*}
При $n=2$ в полярных координатах $x_1 = r \cos \varphi,$ $x_2 = r
\sin \varphi$ соответствующее краевое условие имеет вид
\begin{equation}
\lim\limits_{r \to +0} \frac{1}{2 \pi} \int\limits_{- \pi}^{\pi} u
\lr{r, \varphi'}  \lr{ \frac{2r \cos \lr{\varphi-\varphi'} -2
r^2}{1- 2r \cos \lr{\varphi-\varphi'} + r^2} +\frac{1}{\ln r}}  d
\varphi' = \Psi \lr{\varphi}. \label{54}
\end{equation}

Суть условий~\eqref{53},~\eqref{54} заключается в предварительном
усреднении с подходящим ядром функции $u \lr{r, \vartheta}$ при
фиксированном $r>0$ по угловым переменным $\vartheta.$ После этого
возможен   предельный переход по $r \to 0.$ Выражение, стоящее
слева в~\eqref{53} или~\eqref{54}, будем называть \mbox{\it
$\sigma$-следом} (или {\it $K$-следом} в честь В.\,В.~Катрахова,
который ввёл это условие) и обозначать символом $\sigma u|_0$ (или
$K u|_0$). $K$-след отличен от нуля лишь для функций, имеющих в
точке $0$ особенность, порядок которой не меньше порядка
особенности фундаментального решения уравнения Лапласа.

Оказывается, у  любой гармонической в   окрестности     точки $0,$
исключая саму эту точку, функции $\sigma$-след существует   и он
однозначно определяет её сингулярную часть. Для строгой
формулировки этого утверждения вводится пространство $A
\lr{\Theta},$ состоящее    из вещественно-аналитических на сфере
$\Theta$ функций $\Psi,$ для    которых конечны при любом $h>0$
 нормы
\begin{equation}
\| \Psi \|_h = \lr{\sum\limits_{k=0}^{\infty}
\sum\limits_{l=1}^{d_k} |\Psi_{k, l}|^2 h^{-2k}}^{\frac{1}{2}},
\label{55}
\end{equation}
где
\begin{equation*}
\Psi_{k, l} = \int\limits_{\Theta} \Psi \lr{\vartheta} Y_{k, l} \lr{\vartheta}  d \vartheta,
\end{equation*}
через $Y_{k,l}$ обозначена полная ортонормированная система сферических гармоник.

Таким образом, к основным результатам пятой главы относятся:
постановка краевых задач для уравнения Пуассона с изолированными
особыми точками, решения которого могут иметь в этих точках
существенные особенности, введение новых нелокальных краевых
условий в особых точках в форме $\sigma$-следа (или $K$-следа),
определение соответствующих функциональных пространств с
использованием операторов преобразования, доказательство прямых и
обратных теорем о следах, доказательство корректности
рассматриваемых краевых задач во введённых функциональных
пространствах.

Далее, в  шестой главе   излагается универсальный композиционный
метод построения операторов преобразования и приведены его
приложения. В этом методе операторы преобразования строятся как из
<<кирпичиков>> из классических интегральных преобразований. На
этом пути удаётся единым способом построить как все уже известные
явные представления для операторов преобразования, так и получить
их многочисленные новые классы.

Наконец, в заключительной седьмой главе рассмотрены некоторые
приложения метода операторов преобразования к дифференциальным
уравнениям с переменными коэффициентами. Сначала мы рассматриваем
задачу о построении нового класса операторов преобразования для
возмущённого уравнения Бесселя с использованием модифицированного
интегрального уравнения для ядра. При этом получены точные оценки
ядер через специальные функции, расширен возможный класс
допустимых потенциалов, которые включают потенциалы Баргмана,
Юкавы, а также сильно сингулярные потенциалы. Для частных случаев
потенциалов установленные общие оценки уточнены. Затем
рассматривается известная задача Е.\,М.~Ландиса об оценке
экспоненциальной скорости убывания решений стационарного уравнения
Шрёдингера, см.~\cite{Lan}. Несмотря на общий отрицательный ответ
в этой задаче, полученный В.\,З.~Мешковым~\cite{Mesh1,Mesh2},
оказывается, что для некоторых потенциалов специального вида
гипотеза Е.\,М.~Ландиса выполняется. Интересно, что этот результат
получен также с использованием техники операторов преобразования с
ядрами специального вида.

\section{Краткий очерк истории и современного состояния теории операторов
преобразования}\label{sec2}

\vskip\baselineskip
Начнём с основного определения.

\begin{definition}
 Пусть дана пара операторов $(A,B).$ Ненулевой оператор
$T$ называется \textit{оператором преобразования} (ОП, {\it transmutation}), если
выполняется соотношение
\begin{equation}
\label{1.1} {T\,A=B\,T.}
\end{equation}
\end{definition}

Соотношение~\eqref{1.1} называется иначе \textit{сплетающим
свойством}, тогда говорят, что ОП $T$ \textit{сплетает} операторы
$A$ и $B$ ({\it intertwining operator}). Для
превращения~\eqref{1.1} в строгое определение необходимо задать
пространства или множества функций, на которых действуют операторы
$A,$ $B,$ $T.$ Обычно в определение ОП закладывают также
требования обратимости и непрерывности, которые являются
желательными, но не обязательными условиями.
 В конкретных реализациях операторы $A$ и $B$ чаще всего (но не обязательно!)
являются дифференциальными, $T$ "--- интегральный линейный
непрерывный оператор на стандартных пространствах.

Ясно, что понятие ОП является прямым и далеко идущим обобщением
понятия подобия матриц из линейной алгебры~\cite{Gan, Hor, Tyr}.
Следует отметить, что не существует эффективных методов проверки
подобия двух конечных матриц, так как не существует эффективных
способов проверки совпадения их Жордановых форм, кроме прямого
вычисления.

 Вместе с тем ОП \textit{не сводятся к подобным {\rm (}или
эквивалентным{\rm )} операторам},  так как сплетаемые операторы
как правило являются неограниченными в естественных пространствах,
к тому же обратный к ОП не обязан существовать, действовать в том
же пространстве или быть ограниченным. Так что спектры операторов,
сплетаемых ОП, как правило не совпадают. Кроме того, сами ОП могут
быть неограниченными. Так бывает, например, в теории
преобразований Дарбу, предметом которой является нахождение
дифференциальных операторов преобразования (подстановок или замен)
между парой дифференциальных же операторов, таким образом в этом
случае все три рассматриваемых оператора являются неограниченными
в естественных пространствах. При этом теория преобразований Дарбу
как соответствующий раздел теории дифференциальных уравнений также
вписывается в общую схему теории операторов преобразования при её
расширенном понимании. Операторы в паре, для которой ищется ОП, не
обязаны быть только дифференциальными. В теории ОП встречаются
задачи для  следующих разнообразных типов операторов:
интегральных,  интегро-дифференциальных,
дифференциально-разностных (например,  типа Дункла),
дифференциальных или интегро-дифференциальных   бесконечного
порядка (например, в вопросах, связанных с леммой Шура о
дополняемости), общих линейных  в фиксированных функциональных
пространствах, псевдодифференциальных и
операторно-дифференциальных (абстрактных дифференциальных).

Для примера кратко изложим модельную схему, иллюстрирующую
использование операторов преобразования для получения формул связи
между решениями возмущённого и невозмущённого уравнений, для
которых доказано сплетающее свойство. Пусть, например,  мы изучаем
некоторый достаточно сложно устроенный оператор $A.$ При этом
нужные свойства уже известны для модельного более простого
оператора $B.$ Тогда, если существует ОП~\eqref{1.1}, то часто
удаётся перенести свойства модельного оператора  $B$ и на $A.$
Такова в нескольких словах примерная схема типичного использования
ОП в конкретных задачах.

В частности, если рассматривается уравнение $Au=f$ с оператором
$A,$ то применяя к нему ОП $T$ со сплетающим
свойством~\eqref{1.1}, получаем уравнение с оператором $B$ вида
$Bv=g,$ где обозначено $v=Tu,$ $g=Tf.$ Поэтому, если второе
уравнение с оператором $B$ является более простым, и для него уже
известны формулы для решений, то мы получаем и представления для
решений первого уравнения $u=T^{-1}v.$ Разумеется, при этом
обратный оператор преобразования должен существовать и действовать
в рассматриваемых пространствах, а для получения явных
представлений решений должно быть получено и явное представление
этого обратного оператора. Таково одно из простейших применений
техники ОП в теории дифференциальных уравнений, как обыкновенных,
так и  с частными производными. При этом следует отметить, что при
наличии пары линейных ОП они осуществляют одновременно связь и
представление решений как для однородных, так и для неоднородных
уравнений, а также для уравнений со спектральным параметром.

Сделаем  терминологические замечания. В западной литературе принят
для ОП термин <<{\it transmutation}>>, восходящий к Ж.~Дельсарту.
Как отмечает Р.~Кэрролл, похожий термин <<{\it transformation}>>
при этом закрепляется за классическими интегральными
преобразованиями Фурье, Лапласа, Меллина, Ханкеля и другими
подобными им. Кроме того, термин <<{\it transmutation}>> имеет в
романских языках дополнительный оттенок <<волшебного
превращения>>, что довольно образно характеризует действие ОП.
Приведём дословную цитату из~\cite{Car3}: <<{\it Such operators
are often called transformation operators by the Russian school
{\rm (}Levitan, Naimark, Marchenko et. al.{\rm )}, but
transformation seems too broad a term, and, since some of the
machinery seems {\rm ``}magical{\rm ''} at times, we have followed
Lions and Delsarte in using the word transmutation}>>. Лучше и
точнее не скажешь. Название <<операторы преобразования>> на
русском языке было предложено в 1940-е годы
В.\,А.~Марченко~\cite{Mar9}.

Необходимость теории операторов преобразования  доказана большим числом её приложений.
Особую роль методы операторов преобразования играют в теории дифференциальных уравнений различных типов.
С их помощью были доказаны многие фундаментальные результаты для различных классов дифференциальных уравнений.

В настоящее время теория операторов преобразования представляет собой полностью оформившийся самостоятельный раздел математики, находящийся на стыке дифференциальных, интегральных и интегродифференциальных уравнений, функционального анализа, теории функций, комплексного анализа, теории специальных функций и дробного интегродифференцирования, гармонического анализа, теории оптимального управления, обратных задач и задач рассеяния.

В развитии теории операторов преобразования можно условно выделить
три основных периода. В первом начальном периоде становления
теории операторов преобразования закладывались базовые идеи и
определения, их источником была теория подобия конечных
матриц~\cite{Gan, Hor, Tyr}, отдельные результаты по подобию
операторов, а также некоторые результаты для простейших
дифференциальных уравнений. Считается, что идея операторов
преобразования в операторной формулировке была высказана
Фридрихсом~\cite{Fri}. На самом деле метод операторов
преобразования для получения представлений решений
дифференциальных уравнений был разработан и впервые применён
намного раньше в 19 веке в работах А.\,В.~Летникова, кроме того
это было и по существу первое реальное применение дробного
интегродифференцирования как ОП к задачам дифференциальных
уравнений~\cite{Shos, Koo1}.

Второй период условно продолжался в течение 1940--1980~гг., его
можно назвать классическим. В этот период были получены
многочисленные результаты в теории операторов преобразования и их
приложениям. Перечислим основные направления и результаты этого
периода.

Методы операторов преобразования были с успехом применены в теории
обратных задач, определяя обобщённое преобразование Фурье,
спектральную функцию и решения знаменитого уравнения
Гельфанда---Левитана (З.\,С.~Агранович, В.\,А.~Марченко~\cite{AM,
Mar1, Mar2, Mar3, Mar4, Mar5, Mar6, Mar7},
Б.\,М.~Левитан~\cite{Lev1, Lev2, Lev3, Lev4, Lev5, Lev6, Lev7,
Lev8}); в теории рассеяния через операторы преобразования
выписывается не менее знаменитое уравнение Марченко
(Б.\,М.~Левитан~\cite{Lev1, Lev2, Lev3, Lev4, Lev5, Lev6, Lev7,
Lev8}, В.\,А.~Марченко~\cite{AM, Mar1,Mar2},
Л.\,Д.~Фаддеев~\cite{Fad1,Fad2}); для обоих классов обратных задач
операторы преобразования являются основным инструментом, так как
перечисленные классические уравнения выписываются для ядер
операторов преобразования, а значения ядер на диагонали
восстанавливают неизвестные потенциалы в обратных задачах по
спектральной функции или данным рассеяния~\cite{Laks1,Laks2, ShSa,
Kol1, Nizh1,Nizh2, New, Bloh}. Для операторов Штурма---Лиувилля
были построены ставшие классическими ОП на отрезке
(А.\,Я.~Повзнер~\cite{Povz}) и полуоси
(Б.\,Я.~Левин~\cite{Levin2}). В спектральной теории были получены
известные формулы следов и асимптотика спектральной функции
(В.\,А.~Марченко~\cite{Mar1,Mar2}, Б.\,М.~Левитан~\cite{Lev1,
Lev2, Lev3, Lev4, Lev5, Lev6, Lev7, Lev8}); оценки ядер операторов
преобразования, отвечающие за устойчивость обратных задач и задач
рассеяния (В.\,А.~Марченко~\cite{AM, Mar1,Mar2}); оценки решений
Йоста в квантовой теории рассеяния (З.\,С.~Агранович,
В.\,А.~Марченко~\cite{AM, Mar1,Mar2}, Б.\,М.~Левитан~\cite{Lev1,
Lev2, Lev3, Lev4, Lev5, Lev6, Lev7, Lev8},
В.\,В.~Сташевская~\cite{Sta1,Sta2}, А.\,С.~Сохин~\cite{Soh1, Soh2,
Soh3, Soh4}). В результате применения ОП можно сказать, что теория
операторов Штурма---Лиувилля с переменным коэффициентом была
тривиализирована до уровня простейшего уравнения с
тригонометрическими или экспоненциальными решениями. Также была
исследована система Дирака и другие матричные системы
дифференциальных уравнений (Б.\,М.~Левитан,
И.\,С.~Саргсян~\cite{Lev6}).

Была развита теория обобщённых аналитических функций, которую
можно трактовать как раздел теории операторов преобразования,
сплетающих невозмущённые и возмущённые уравнения Коши---Римана
(Л.~Берс~\cite{Bers1,Bers2}, С.~Бергман~\cite{Berg},
И.\,Н.~Векуа~\cite{Vek3, Vek4}, Б.~Боярский~\cite{Boyar},
Г.\,Н.~Положий~\cite{Pol1, Pol2, Pol3}) с приложениями в задачах
механики, теории упругости и газодинамики. На основе методов
операторов преобразований был создан новый раздел гармонического
анализа, изучающий различные модификации операторов обобщённого
сдвига и обобщённых операторных свёрток
(Ж.~Дельсарт~\cite{Del1,Del2, Del7},
Я.\,И.~Житомирский~\cite{Zhit}, Б.\,М.~Левитан~\cite{Lev2,Lev3}).
Была установлена глубокая связь операторов преобразования с
теоремами типа Пэли---Винера (В.\,В.~Сташевская~\cite{Sta1,Sta2},
А.\,И.~Ахиезер~\cite{Ahi1}, Х.~Шабли~\cite{Che1, Che2, Che3,
Che4}, Х.~Тримеш~\cite{Tri1,Tri2}). Теория операторов
преобразования позволила дать новую классификацию специальных
функций и интегральных операторов со специальными функциями в
ядрах (Р.~Кэрролл~\cite{Car1, Car2, Car3},
Т.~Корвиндер~\cite{Koo1}). При этом нахождение ядер ОП использует
существование и явный вид функций Грина или Римана для различных
классов дифференциальных уравнений~\cite{Sob, VNN, VZ}, стимулируя
нахождение этих функций для различных задач.

В теории нелинейных дифференциальных уравнений был разработан
метод Лакса, который  использует операторы преобразования для
доказательства существования решений и построения
солитонов~\cite{Zhu, ZMNP, AbSi, Car4}, также широкие применения
нашли преобразования Дарбу~\cite{MaSa}, которые можно
рассматривать как операторы преобразования, в которых и сплетаемые
и сплетающий операторы являются дифференциальными; о связи теорий
преобразования Дарбу и ОП см. обзор~\cite{BaSa}.  В квантовой
физике при рассмотрении уравнения Шрёдингера и задач теории
рассеяния был изучен специальный класс операторов преобразования
"--- волновые операторы. Рассмотрение общих задач рассеяния и
обратных задач с точки зрения ОП изложено в
обзорах~\cite{Fad1,Fad2, Mar9}. В работе~\cite{Kach1} волновые
операторы построены для задач теории рассеяния с потенциалом
Штарка, к сожалению эта работа В.\,П.~Качалова и Я.\,В.~Курылёва
1989~г. практически забыта, например,  в сборнике~\cite{AppHyp}
1995~г. в статье~\cite{Lev1995} Б.\,М.~Левитан формулирует задачу
о построении соответствующего оператора преобразования как
нерешённую.

В самой теории ОП были обнаружены ограничения, связанные с
порядком дифференциального оператора. Так, для дифференциальных
операторов порядков выше третьего было показано, что классические
ОП в форме Вольтерра существуют только для аналитических
коэффициентов (В.\,И.~Мацаев~\cite{Mats},
Л.\,А.~Сахнович~\cite{Sah1, Sah2, Sah3},
М.\,М.~Маламуд~\cite{Mal1, Mal2, Mal3, Mal4, Mal5}), а в общем
случае ОП имеют более сложную структуру, требующую выхода в
комплексную плоскость даже для построения действительных решений
(А.\,Ф.~Леон\-тьев~\cite{Leo}, Ю.\,Н.~Валицкий~\cite{Val},
И.\,Г.~Хачатрян~\cite{Hach1,Hach2}, М.\,М.~Маламуд~\cite{Mal1,
Mal2, Mal3, Mal4, Mal5}, А.\,П.~Хромов~\cite{Hrom1}). Вместе с тем
в пространствах аналитических функций была доказана
эквивалентность дифференциальных операторов одного порядка и
изучен целый ряд задач (Д.\,К.~Фаге~\cite{FN, Fage1, Fage2, Fage3,
Fage4, Fage5, Fage6}, В.\,А.~Марченко~\cite{Mar5, Mar6, Mar7},
Ю.\,Ф.~Коробейник~\cite{Kor1, Kor2}, М.\,К.~Фишман~\cite{Fish}). С
целью применения к теории ОП была построена теория разрешимости
для известного уравнения Бианки (Д.\,К.~Фаге~\cite{FN}).

Отдельной областью применения ОП стала теория дифференциальных
уравнений с особенностями в коэффициентах, прежде всего с
операторами Бесселя
\begin{equation}\label{Bes1}{B_{\nu}u(x)=\frac{d^2
u}{dx^2}+\frac{2\nu+1}{x} \frac{du}{dx}.}\end{equation} На
первоначальном этапе исследований уравнений этого класса
применялась пара известных ОП Сонина и Пуассона (см. определения в
главе~\ref{ch1}, подробное рассмотрение в главе~\ref{ch2}). Как ОП
эти операторы впервые были введены в работах Жана
Дельсарта~\cite{Del1, Del2, Del3, Del4}, а затем на основе идей
Дельсарта их изучение продолжилось в работах Дельсарта и
Лионса~\cite{Del5,Del6, Del7, Lio1, Lio2, Lio3}. Поэтому мы будем
использовать термин ОП Сонина---Пуассона---Дельсарта (СПД). Об
операторах СПД см. также известную статью
Б.\,М.~Левитана~\cite{Lev7}, во многом основанную на классических
результатах Дарбу.

Ж.~Дельсартом на базе ОП СПД было введено фундаментальное понятие
операторов обобщённого сдвига (ООС, см. определение). Были
разработаны многочисленные конструкции обобщённого гармонического
анализа, основанные на определениях обобщённого сдвига и вводимых
с его помощью групповых структурах. Направление обобщённых
почти-периоди\-ческих функций с использованием ОП типа СПД и ООС
было заложено в работах Ж.~Дельсарта 1938~г.~\cite{Del1,Del2} и
продолжено Дельсартом и Лионсом в~\cite{Del5,Del6, Lio1, Lio2,
Lio3}. Этот круг вопросов параллельно был исчерпывающе изучен в
работах Б.\,М.~Левитана 1940~г., и особенно 1947--1949~гг.,
результаты вошли в его классические монографии~\cite{Lev2, Lev3,
Lev4}. (Отметим неточности в статье Большой Советской Энциклопедии
<<Почти периодическая функция>>. Там имя  Дельсарта не
упоминается, а работы Б.\,М.~Левитана датированы 1938~г. Насколько
известно автору, первая печатная работа Б.\,М.~Левитана на
немецком языке вышла в 1940~г.). Также в работах Жана Дельсарта
были впервые построены разложения по обобщённым рядам Тэйлора,
справедливо названных рядами Тэйлора---Дельсарта~\cite{Del2, Del7,
Lev1, Lev2, Lev3}, такие ряды изучаются до сих пор во многих
работах, см., например,~\cite{Kam}. Следует отметить, что
первоначальным источником и прототипом большинства вариантов
обобщённого гармонического анализа были  операторы Бесселя и
связанные с ними дифференциальные уравнения.

Огромную роль конструкции ОП и ООС сыграли в теории уравнений с
частными производными, см. краткий исторический обзор
В.\,А.~Марченко~\cite{Mar9}. ОП позволяют преобразовывать более
сложные уравнения в более простые, ООС помогают в сингулярных
уравнениях переносить особенность из начала координат в
произвольную точку, а также с их использованием можно строить
фундаментальные решения, а затем с помощью обобщённой свёртки
интегральные представления для решений соответствующих
дифференциальных уравнений.

Рассматривались также ОП для многочисленных обобщений оператора
Бесселя. Важным обобщением операторов
Сонина---Пуассона---Дель\-сарта являются ОП для гипербесселевых
функций. Теория таких функций была первоначально заложена в
работах Куммера и Делерю. Полное исследование гипербесселевых
функций, дифференциальных уравнений для них и соответствующих
операторов преобразования  было исчерпывающе проведено в работах
И.~Димовски и его учеников~\cite{Dim, DK1,DK2}. Соответствующие ОП
заслуженно получили в литературе названия ОП Сонина---Димовски и
Пуассона---Димовски, они также изучались в работах ученицы
И.~Димовски "--- В.~Киряковой~\cite{DK1,DK2, Kir1, Kir2, Kir3,
Kir4, Kir5}. В теории гипербесселевых функций, дифференциальных
уравнений и операторов преобразования для них центральную роль
играет знаменитое интегральное преобразование Обрешкова, введённое
болгарским математиком Н.~Обрешковым. Это преобразование, ядро
которого выражается в общем случае через $G$-функцию Майера,
является одновременным обобщением преобразований Лапласа, Меллина,
синус- и косинус-преобразований Фурье, Ханкеля, Майера и других
классических интегральных преобразований. Различные формы
гипербесселевых функций, дифференциальных уравнений и  операторов
преобразований для них, а также частные случаи преобразования
Обрешкова многократно впоследствии переоткрывались. По мнению
автора, преобразование Обрешкова, наряду с преобразованиями Фурье,
Меллина и Лапласа относится к небольшому числу фундаментальных, из
которых, как из кирпичиков, складываются многие другие
преобразования, а также основанные на них конструкции и
приложения. Преобразование Обрешкова было исторически первым
интегральным преобразованием, ядро которого выражается через
$G$-функцию Майера, но не выражается через одну (!) обобщённую
гипергеометрическую функцию. Аналогично было введено сербским
математиком Б.~Станковичем другое важное интегральное
преобразование "--- преобразование Станковича "--- ядро которого
выражается через $H$-функцию Райта---Фокса, но не выражается через
более простую $G$-функцию Майера, преобразование Станковича
находит важные применения при изучении дифференциальных уравнений
дробного порядка типа дробной диффузии~\cite{Koch1,Koch2,
EiIvKoch, Pshu1,Pshu2}.

Вместе с тем были построены аналогичные ОП и для некоторых других
модельных операторов, например, таких~\cite{Car1, Car2, Car3,
Yar1,Yar2}:
\begin{gather}
\label{510}
A=\frac{1}{v(x)}\frac{d}{dx}v(x)\frac{d}{dx},\\
v(x)=\sin^{2\nu+1}x,\sh^{2\nu+1}x,
(e^x-e^{-x})^{2\nu+1}(e^x+e^{-x})^{2\mu+1}.
\end{gather}
Важность операторов $A$ вида~\eqref{510} для теории заключается в
том, что по знаменитой формуле Гельфанда они представляют
радиальную часть оператора Лапласа на симметрических
пространствах~\cite{Hel1}. При этом оператор Бесселя получается
при выборе в~\eqref{510} $v(x)=x^{2\nu+1}.$ Другим модельным
оператором, для которого построены ОП, является оператор Эйри
$D^2+x,$ рассматривался также его возмущённый вариант, связанный с
эффектом Штарка из квантовой механики~\cite{Kach1}, и случай
возмущённого оператора Хилла с периодическим
потенциалом~\cite{Lev1995}. Были изучены операторы сдвига по
спектральному параметру Векуа---Эрдейи---Лаундеса~\cite{Low1,
Low2, Low3}.

К третьему современному периоду развития теории ОП можно отнести
работы с 1990-х годов и до настоящего времени, за этот период были
получены и продолжают появляться многие важные исследования, см.,
например, обзоры~\cite{Mar9, CB, Lev1995, S42, S46, S38, S400,
S401} и диссертацию~\cite{SitDis}. Перечислим некоторые из
направлений исследований, связанных с применением методов
операторов преобразования. Продолжено развитие теории обобщённых
аналитических функций (А.\,П.~Солдатов~\cite{Sol},
С.\,Б.~Климентов~\cite{Kli1, Kli2, Kli3, Kli5},
В.\,В.~Кравченко~\cite{Krav1}). Были найдены приложения операторов
преобразования к вложениям функциональных пространств и обобщению
операторов Харди~\cite{S66, S6, S42, S46, S38, S400, S401},
построению различных конструкций обобщённого сдвига и основанным
на них обобщённых вариантов гармонического анализа, см.
обзор~\cite{Mar9}, а также А.\,Д.~Гаджиев~\cite{Gad1},
В.~Гулиев~\cite{Gul1,Gul2}, С.\,С.~Платонов~\cite{Plat1, Plat2,
Plat3}, Л.\,Н.~Ляхов~\cite{Lyah1, Lyah2, Lyah3, LPSh1, LShFrac,
LPSh2}. В работах Ф.\,Г.~Мухлисова и его учеников рассматривались
задачи для $B$-потенциалов. В работах Э.\,Л.~Шишкиной исследованы
теоремы о среднем типа Айсгерссона, новые задачи для уравнения
Эйлера---Пуассона---Дарбу, построены новые классы потенциалов типа
Рисса с $B$-гиперболическими и ультрагиперболическими операторами,
рассмотрены применения этих результатов для соответствующих типов
дифференциальных уравнений~\cite{ShiR1, ShiR2, ShiR3, ShiR4,
ShiE1, ShiE2, SS, ShiE3, ShiE4, ShiE5, ShiE6}. Продолжилось
применение ОП и родственных методов в теории обратных задач и
теории рассеяния~\cite{But, Ram1, Mar9, Yurko, ChCPR, PiSa}. Для
дифференциальных уравнений продолжается развитие метода Дарбу и
его модификаций~\cite{Kap1, MaSa}. Были рассмотрены новые классы
задач для решений с существенными особенностями на части границы
во внутренних или угловых точках (В.\,В.~Катрахов~\cite{Kat1, 30,
Kat2, 32, KatDis, Kat3, Kat4},
И.\,А.~Киприянов~\cite{KiKa3,KiKa4}), получены точные оценки
скорости убывания решений некоторых эллиптических и
ультраэллиптических уравнений (В.\,З.~Мешков~\cite{Mesh1,Mesh2,
S3}). Отдельной тематикой стало использование ОП при исследовании
различных операторов дробного интегродифференцирования
(И.~Димовски, В.~Кирякова~\cite{AKK, Dim, DK1,DK2, Kir1, Kir2,
Kir3, Kir4, Kir5}, Н.\,А.~Вирченко~\cite{ViRy, Vir1,Vir2}, а
также~\cite{Bac,SS,ShiE3,FJSS}). Было продолжено  с использованием
методов ОП изучение сингулярных и вырождающихся краевых задач,
псевдодифференциальных операторов (В.\,В.~Катрахов~\cite{Kat1, 30,
Kat2, 32, KatDis, Kat3, Kat4},
И.\,А.~Киприянов~\cite{KiKa3,KiKa4}, Хе Кан Чер~\cite{Kan},
Л.\,Н.~Ляхов~\cite{Lyah1,Lyah2}, О.\,А.~Репин~\cite{Rep, Lar1,
Lar2, Lar3, Lar4, Lar5, Lar6}), операторных уравнений
(А.\,В.~Глушак~\cite{Glu1, Glu2, Glu3, Glu4, Glu55, Glu5, Glu6,
Glu7, Glu8, Glu9, Glu10, Glu11, Glu12},
В.\,Е.~Фёдоров~\cite{FGP,FeIv}). Уравнения с оператором Бесселя и
связанные с ними вопросы изучают А.\,В.~Глушак~\cite{Glu1, Glu2,
Glu3, Glu4, Glu55, Glu5, Glu6, Glu7, Glu8, Glu9, Glu10, Glu11,
Glu12, Glu13}, В.\,С.~Гулиев~\cite{Gul1,Gul2},
Л.\,Н.~Ляхов~\cite{Lyah1, Lyah2, Lyah3, LPSh1, LShFrac, LPSh2},
Л.\,С.~Пулькина~\cite{Pulkina}, К.\,Б.~Сабитов~\cite{SaIl},
В.\,В.~Кравченко~\cite{CKT1, CKT2, Krav1, Krav2, Krav3, Krav4,
Krav5, Krav6, Krav7, Krav8, Krav9, Krav10} со своими коллегами и
учениками, а также А.\,Б.~Муравник
(дифференциально-функцио\-нальные уравнения, стабилизация
решений~\cite{Mur,Mur8,Mur9}, свойства преобразований
Фурье---Бесселя~\cite{Mur3,Mur4}), В.\,В.~Волчков~\cite{Volch},
И.\,П.~Половинкин (теоремы о среднем для уравнений с операторами
Бесселя~\cite{LPSh1, LShFrac, LPSh2}), Э.\,Л.~Шишкина
($B$-гиперболические потенциалы и обобщённые средние~\cite{ShiR1,
ShiR2, ShiR3, ShiR4, ShiE1, ShiE2, SS, ShiE3, ShiE4, ShiE5, ShiE6,
LPSh1, LShFrac, LPSh2}), Ш.\,Т.~Каримов~\cite{KarST, KarST1,
KarST2, KarST3, KarST4}, А.~Хасанов~\cite{SaHa1},
Э.\,Т.~Каримов~\cite{HaKa1}, Т.\,Г.~Эргашев~\cite{Erg1},
И.\,Б.~Гарипов, Н.\,В.~Зайцева, Р.\,М.~Мавлявиев, Ф.\,Г.~Хуштова и
другие.

В последнее время также разработаны эффективные численные методы
для применения операторов преобразования при расчётах решений
дифференциальных уравнений, их собственных функций и спектральных
характеристик (В.\,В.~Кравченко, С.~Торба~\cite{CKT1, CKT2, Krav1,
Krav2, Krav3, Krav4, Krav5, Krav6, Krav7, Krav8, Krav9, Krav10}),
эти методы основаны на разработанном В.\,В.~Кравченко методе
представления решений уравнений Штурма---Лиувилля и их обобщений с
операторами Бесселя или Дирака (метод SPPS "--- {\it Spectral
Parameter Power Series}). Отметим, что идея представления ядер
операторов преобразований рядами как альтернатива интегральному
представлению ОП является естественной, так как ядра получаются
методом последовательных приближений, то есть в виде рядов
Неймана. Другое представление ОП в виде рядов с использованием
обобщённых базисов изучено в~\cite{FN}. Для возмущённых уравнений
Бесселя с переменным потенциалом подробное изучение представлений
ядер ОП в виде рядов получено в~\cite{CFH,FH}.

Отдельный класс задач составляют задачи типа Дирихле---Нейман,
Нейман---Дирихле, при которых оператор преобразования действует на
краевые или начальные условия, сохраняя дифференциальное
выражение; такие задачи нашли важные приложения в механике
(О.\,Э.~Яремко~\cite{BMYa, Yarem}), кроме того этот класс задач
тесно связан со спектральной теорией и теорией вероятностей.
Достаточно законченные модификации гармонического анализа для
операторов Бесселя построены в работах
С.\,С.~Платонова~\cite{Plat1, Plat2, Plat3},  для возмущённого
оператора типа Бесселя с переменными коэффициентами "--- в работах
Х.~Тримеша~\cite{Tri3, Tri4, Tri5}, в последнее время активно
создаётся гармонический анализ для дифференциально-разностных
операторов типа Дункла~\cite{Me1, DHS, Dun1, Dun2, Dun3,
Gal1,Gal2, Rod, Ros1, Ros2, Tri6} на основе соответствующих
обобщений операторов СПД. Наличие благодаря ОП соответствующих ООС
позволяет также определить обобщённую свёртку, новые
алгебраические и групповые структуры, рассматривать различные
задачи аппроксимации функций~\cite{AppHyp}. Идеи М.\,К.~Фаге,
развитые  для уравнения Бианки в связи с построением ОП для
дифференциальных уравнений высоких порядков, нашли своё
продолжение в исследовании более общих уравнений в работах
В.\,И.~Жегалова, А.\,Н.~Миронова, Е.\,А.~Уткиной~\cite{ZhM1,ZhM2}.
В теории уравнений дробного порядка появились работы, которые
можно трактовать как метод ОП для представления решений уравнений
дробного порядка через решения уравнений целого порядка (так
называемый принцип субординации "--- Я.~Прусс, А.\,Н.~Кочубей,
Э.~Бажлекова, А.~Псху). ОП находят применения в теории
преобразования Радона и математической томографии~\cite{Nat1,
Hel1, Rub3, Rub1, Rub2, Rub4}, а также при разложении функций в
различные ряды по специальным функциям~\cite{Kam}. В работах
В.\,А.~Марченко продолжилось применение ОП к квантовой
теории~\cite{Mar8,Mar9}.

Теория ОП также связана с вопросами факторизации дифференциальных
операторов (Л.\,М.~Беркович~\cite{Berk}, А.\,Б.~Шабат). При
изучении групповых свойств дифференциальных уравнений важное
значение имеет лемма Шура о дополняемости, которую можно
трактовать как существование формального ОП между оператором
дробного интегрирования и некоторым дифференциальным оператором
бесконечного порядка, подобная задача рассматривалась ещё в
учебнике Н.~Бурбаки~\cite{Burb}. В  последних работах Р.~Кэрролла
была сделана попытка построения <<квантовых>>   ОП для $q$-дифференциальных операторов~\cite{Car10}.
 Разные задачи,  использующие идеи ОП или родственные методы, также рассматривались в~\cite{Arsh, Bas, Han, GSPP, Kam, Matv, SPP}.

Возможность, чтобы исходная и преобразованная функции принадлежали
различным пространствам, что принято подчёркивать использованием
различных обозначений для переменных, позволяет включить в общую
схему ОП все классические интегральные преобразования: Фурье,
Лапласа (на самом деле Петцваля), Меллина, Ханкеля, Вейерштрасса,
Конторовича---Лебедева, Фока, Обрешкова, Станковича и
другие~\cite{BE1,BE2, Ome}. При вычислении интегралов, необходимых
для реализации метода ОП, фундаментальные приложения нашла теорема
Слейтер---Мари\-чева, соединившая методы преобразования Меллина с
теорией гипергеометрических функций~\cite{Marich1, PBM}. В общую
схему ОП также включаются конечные интегральные преобразования
Г.\,А.~Гринберга~\cite{Grin}. Продолжилось исследование функций
Грина и Римана, которые используются в методе ОП~\cite{Ler,
MaKiRe, Sob}. Существуют связи ОП с теорией дробных (квадратичных)
интегральных преобразований Фурье~\cite{OZK} и Ханкеля. Решение
В.\,А.~Чернятиным знаменитой задачи о нахождении необходимых и
достаточных условий на потенциалы для обоснования метода Фурье для
волнового уравнения с переменным коэффициентом открывает
перспективы использования этого метода для оценки ядер в теории
ОП~\cite{Cher}.

Коммутирующие операторы любой природы также подходят под
определение ОП. Наиболее близко к духу и задачам теории ОП
относится изучение операторов, коммутирующих с производными. Сами
ОП в этом случае зачастую представляются формальными рядами,
псевдо-дифференциальными операторами или дифференциальными
операторами бесконечного порядка. Описание коммутантов напрямую
связано с описанием всего семейства ОП для заданной пары по его
единственному представителю. В этом классе задач фундаментальные
приложения нашла теория операторных свёрток, особенно свёртки
Берга---Димовски (Л.~Берг, И.~Димовски~\cite{Dim, Bozh}). Начинают
находить приложения в теории ОП и результаты для коммутирующих
дифференциальных операторов, восходящие к классическим работам
Бёчнела и Чонди (J.\,L.~Burchnall, T.\,W.~Chaundy).

Важным разделом теории ОП стал специальный класс "--- ОП
Бушмана---Эрдейи (см. главу~\ref{ch3}). Это класс ОП, который при
определённом выборе параметров является одновременным обобщением
ОП Сонина---Пуассона---Дельсарта и их сопряжённых, операторов
дробного интегродифференцирования Римана---Лиувилля и
Эрдейи---Кобера, а также интегральных преобразований
Мелера---Фока. Термин <<операторы Бушмана---Эрдейи>>  как наиболее
исторически оправданный был введён одним из авторов в~\cite{S66,
S6}, впоследствии он использовался и другими авторами. Важность
операторов Бушмана---Эрдейи во многом обусловлена их
многочисленными приложениями. Например, они встречаются в
следующих вопросах теории уравнений с частными
производными~\cite{SKM}: при решении задачи Дирихле для уравнения
Эйлера---Пуассона---Дарбу в четверти плоскости и установлении
соотношений между значениями решений уравнения
Эйлера---Пуассона---Дарбу на многообразии начальных данных и
характеристике (см. лемму Копсона выше), теории преобразования
Радона, так как в силу результатов Людвига~\cite{Lud, Nat1, Deans,
Rub3, Rub1, Rub2, Rub4} действие преобразования Радона при
разложении по сферическим гармоникам сводится как раз к операторам
Бушмана---Эрдейи по радиальной переменной, при исследовании
краевых задач для различных уравнений с существенными
особенностями внутри области, доказательству вложения пространств
И.\,А.~Киприянова  в весовые пространства С.\,Л.~Соболева,
установлению связей между операторами преобразования и волновыми
операторами теории рассеяния, обобщению классических интегральных
представлений Сонина и Пуассона и операторов преобразования
Сонина---Пуассона---Дельсарта.

Наиболее полное изучение операторов Бушмана---Эрдейи  было
проведено в работах С.\,М.~Ситника в 1980--1990-е
годы~\cite{S1,S70,S72,S2,S73,S4,S66,S65,S6,S5,S7}, и затем
продолжено в последующие годы в~\cite{S46,S14,S103,S400,SitDis,
SSfiz,S42,S94,S38,S401,S402}. При этом необходимо отметить, что
роль операторов Бушмана---Эрдейи как ОП до указанных работ вообще
ранее нигде не отмечалась и не рассматривалась.

В литературе изложению теории ОП и их приложениям посвящены
существенные части монографий~\cite{CSh, Car1, Car2, Car3, Car4,
Berg, GB, Lio1, BMYa, Yarem}, кроме того различные вопросы ОП
рассматриваются также в~\cite{SKM, Col1,Col2, Gil1,Gil2} и целом
ряде других известных монографий, а также в
обзорах~\cite{CB,Lev1995,Mar9,S46,S42}. К сожалению, на русском
языке пока нет книг, полностью посвящённых  ОП, таких, как
содержательные книги Роберта Кэрролла на английском~\cite{CSh,
Car1, Car2, Car3, Car4}, поэтому данная монография призвана в
определённой степени заполнить этот пробел. При этом следует
безусловно отметить монографию Д.\,К.~Фаге и
Н.\,И.~Нагнибида~\cite{FN}. В этой монографии практически никак не
отражены известные к тому времени результаты теории ОП, что
полностью компенсируется изложением в основном собственных
результатов авторов по одной из самых трудных задач теории ОП "---
их построении  для дифференциальных операторов высоких порядков с
переменными коэффициентами. Кроме того, в эту монографию вошли и
многие другие вопросы: решение задачи об операторах, коммутирующих
с производными в пространствах аналитических функций (включая
исправление ошибочных результатов Дельсарта и Лионса), создание
законченной теории разрешимости для уравнения Бианки, теория
операторно-аналитических функций (первоначально возникшая в
работах В.\,А.~Марченко~\cite{Mar5, Mar6, Mar7}), исследование
операторов дифференцирования, интегрирования и корней из них в
пространствах аналитических функций.

Отметим также, что С.\,М.~Ситником был разработан метод уточнения
неравенств Коши---Буняковского методом средних~\cite{S41, S45,
S140, S53}, интегральные уточнения по этому методу могут быть
применены к оценкам норм различных операторов преобразования.
Кроме того, полезными являются различные оценки специальных
функций, которые позволяют оценивать ядра операторов
преобразования, см., например,~\cite{S61, S9, S12, S13, S15, S16,
S24, MS1, MS2, MS3}.

Таким образом, методы теории ОП и связанные с ними задачи в той
или иной степени применялись в работах многих математиков.
Перечислим некоторых из них: A.\,I.~Aliev,  H.~Begehr, S.~Bergman,
J.~Betancor,  A.~Boumenir, B.~Braaksma, L.~Bragg, R.~Carroll,
R.~Castillo-P\'erez, H.~Chebli, I.~Dimovski, C.~Dunkl,
J.~Delsarte, A.~Fitouhi, A.~Gasmi, R.~Gilbert, M.\,M.~Hamza,
M.~Holzleitner, R.\,O.~Hriniv, V.~Hristov, V.~Hutson,
G.\,К.~Каlisсh, S.\,L.~Kalla, T.\,H.~Koornwinder, V.~Kiryakova,
J.~L\H{o}ffstr\H{o}m, J.~Lions, Y.~Luchko, M.~Moro,
Y.\,V.~Mykytyuk, J.\,S.~Pym, B.~Rubin, J.\,Yu.~Santana-Bejarano,
F.~Santosa, J.~Siersma, M.~Sifi, A.\,M.~Sinclair, S.\,V.~de Snoo,
M.~Spiridonova, K.~Stempak, G.~Teschl, V.~Thyssen,
K.~Trim\`{e}che, Y.~Tsankov, M.~Voit, Vu Kim Tuan,
З.\,С.~Агранович, А.\,А.~Андрощук, И.\,И.~Баврин, А.\,Г.~Баскаков,
Л.\,Е.~Бритвина, С.\,А.~Бутерин, Ю.\,Н.~Валицкий, В.\,Я.~Волк,
В.\,В.~Волчков, А.\,Д.~Гаджиев, А.\,В.~Глушак, М.\,Л.~Горбачук,
И.\,Ц.~Гохберг, В.\,С.~Гулиев, И.\,М.~Гусейнов,
Я.\,И.~Житомирский, Л.\,А.~Иванов, Х.\,К.~Ишкин, Т.\,В.~Елисеева,
М.\,С.~Ерёмин, М.\,Б.~Капилевич, Ш.\,Т.~Каримов, Д.\,Б.~Карп,
В.\,В.~Катрахов, А.\,П.~Качалов, А.~КилбасА, И.\,А.~Киприянов,
М.\,И.~Ключанцев, В.\,И.~Кононенко, Ю.\,Ф.~Коробейник,
А.\,С.~Костенко, В.\,В.~Кравченко, М.\,Г.~Крейн, П.\,П.~Кулиш,
И.\,Ф.~Кушнирчук, Г.\,И.~Лаптев, Б.\,Я.~Левин, Б.\,М.~Левитан,
А.\,Ф.~Леонтьев, Н.\,Е.~Линчук, С.\,С.~Линчук, А.\,Р.~Лятифова,
Л.\,Н.~Ляхов, Г.\,В.~Ляховецкий, М.\,М.~Маламуд, В.\,А.~Марченко,
В.\,Л.~Матросов, В.\,И.~Мацаев, А.\,Б.~Муравник, Н.\,И.~Нагнибида,
Л.\,П.~Нижник, М.\,Н.~Олевский, Ю.\,А.~Парфенова,
А.\,К.~Прикарпатский, С.\,С.~Платонов, А.\,Я.~Повзнер,   Б.~Рубин,
Ф.\,С.~Рофе-Бекетов, К.\,Б.~Сабитов, А.\,М.~Самойленко,
Л.\,А.~Сахнович,  С.\,М.~Ситник, А.\,С.~Сохин, В.\,В.~Сташевская,
С.\,М.~Торба, Л.\,Д.~Фаддеев, Д.\,К.~Фаге, К.\,М.~Фишман,
И.\,Г.~Хачатрян, А.\,П.~Хромов, Л.\,А.~Чудов, Э.\,Л.~Шишкина,
С.\,Д.~Шмулевич, В.\,А.~Юрко, О.\,Э.~Яремко, В.\,Я.~Ярославцева.
Разумеется, этот список не полон и может быть существенно
расширен.

Из приведённого анализа следует, что метод операторов
преобразования является эффективным методом в теории
дифференциальных уравнений и смежных разделах математики, ему
посвящено большое число работ и на его основе решаются
многочисленные классы задач. Однако, несмотря на вышеизложенное, в
теории операторов преобразования остаются существенные пробелы и
многие нерешённые задачи. Так, для операторов преобразования,
сплетающих дифференциальные операторы или решения дифференциальных
уравнений с особенностями в коэффициентах, включая
дифференциальные операторы Бесселя, отсутствует подробная
классификация с описанием основных свойств. Подробно изучены и
описаны в литературе свойства и приложения простейшего класса ОП
Сонина и Пуассона, но отсутствует систематическое изложение и
доказательства многих свойств для их важных обобщений "---
операторов Бушмана---Эрдейи. До работ С.\,М.~Ситника не отмечалось
вообще, что интегральные операторы  Бушмана---Эрдейи являются ОП
для дифференциальных уравнений с операторами Бесселя. Не
существует общих схем для построения ОП нужных классов, доведённых
до возможности построения на их основе явных формул, сплетающих
решения различных дифференциальных уравнений. Практически
отсутствуют работы, вскрывающие связь ОП с основными конструкциями
дробного исчисления. Не рассматривались возможные применения ОП к
доказательствам вложений функциональных пространств, таких, как
пространства С.\,Л.~Соболева и И.\,А.~Киприянова, в том числе
энергетических пространств для сингулярных дифференциальных
уравнений в частных производных с операторами Бесселя по одной или
нескольким переменным. Для операторов второго порядка с
переменными коэффициентами при построении ОП использовались
методы, дающими грубые оценки ядер ОП с неопределёнными
постоянными, неточные требования на коэффициенты дифференциальных
уравнений приводили к сужению их классов, например классов
допустимых потенциалов для задач Штурма---Лиувилля и их обобщений
для дифференциальных уравнений с особенностями в коэффициентах.
Методы ОП практически не применялись к получению точных оценок для
решений дифференциальных уравнений, например, в таких задачах, как
известная задача Е.\,М.~Ландиса. Также сложилась парадоксальная
ситуация, когда дробные степени оператора Бесселя, которые
используются во многих работах, определяются исключительно неявно
в терминах преобразования Фурье---Бесселя или Ханкеля, а при этом
отсутствуют формулы для их явного определения в интегрального
виде, хотя именно с таких представлений начиналась теория
классических дробных интегралов Римана---Лиувилля. Многие простые
и естественные конструкции ОП для стандартных пар дифференциальных
операторов не построены в явном виде. Также не вводились и не
рассматривались общие схемы для оценок ядер ОП в широко
используемых функциональных пространствах, требующие уточнений и
обобщений классических неравенств. Полученные в последнее время
точные неравенства для многих специальных функций не находили
применения для оценки ядер ОП.

Решения ряда перечисленных выше задач и приведено в данной     монографии.

\section{Основные определения, обозначения и свойства: специальные функции, функциональные пространства, интегральные
преобразования}\label{sec3}
\sectionmarknum{Основные определения, обозначения и свойства}

\subsection{Специальные функции}\label{sec3.1}

Здесь даны краткие определения и пояснения для специальных
функций, которые использованы в работе, следуя монографиям и
справочникам~\cite{AS, BE1, BE2, BE3, Luke3, Luke1,Luke2, Wit1,
AAR, PBM123, NIST}. Также приведён ряд кратких исторических
комментариев.

\subsubsection{Гамма-функция, бета-функция, пси-функция, символ
    Похгаммера и биномиальные коэффициенты}\label{sec3.1.1}

Гамма-функция  является обобщением понятия факториала на случай
чисел, не являющихся натуральными. Бета-функция в общем случае
определяется через гамма-функции. Пси-функция является
логарифмической производной гамма-функции.

Пусть $z\in\mathbb{C}.$ \textit{Гамма-функция} $\Gamma(z)$
определялась Эйлером как предел
$$
\Gamma(z)=\lim\limits_{N\rightarrow\infty}\frac{N!N^z}{z(z+1)(z+2)\ldots(z+N)},\,\,
z\in\mathbb{C},
$$
но чаще используется определение в виде интеграла Эйлера второго
рода
\begin{equation}\label{Gamma}
\Gamma(z)=\int\limits_{0}^\infty y^{z-1}e^{-y}dy,\qquad \Re z>0,
\end{equation}
который сходится при всех $z\in\mathbb{C},$ для которых $\Re x>0.$

Интегрирование по частям выражения~\eqref{Gamma} приводит к
рекуррентной формуле
\begin{equation}\label{Rec}
\Gamma(z+1)=z\Gamma(z).
\end{equation}
Поскольку $\Gamma(1)=1,$ то рекуррентная формула~\eqref{Rec} для
положительных целых $n$ приводит к равенству
$$
\Gamma(n+1)=n\Gamma(n)=n(n-1)\Gamma(n-1)=\ldots=n(n-1)\cdot\ldots\cdot2\cdot
1\cdot\Gamma(1)
$$
или
$$
\Gamma(n+1)=n!,
$$
которое и позволяет рассматривать гамма-функцию как обобщение
понятия факториала. Переписав формулу~\eqref{Rec} в виде
\begin{equation}\label{ReGa}
\Gamma(z-1)=\frac{\Gamma(z)}{z-1},
\end{equation}
мы получим выражение, позволяющее определить гамма-функцию от
отрицательных аргументов, для которых определение~\eqref{Gamma}
неприемлемо. Формула~\eqref{ReGa} показывает, что $\Gamma(z)$
имеет в точках $z=0,-1,-2,-3,\ldots$ разрывы второго рода.

После многократного применения равенства~\eqref{ReGa} получим
формулы \emph{понижения и повышения}, которые, соответственно,
имеют вид
\begin{equation}\label{Povysh}
\Gamma(z+n)=z(z+1)\ldots(z+n-1)\Gamma(z),\qquad n=1,2,\ldots
\end{equation}
и
\begin{equation}\label{Ponizh}
\Gamma(z-n)=\frac{\Gamma(z)}{(z-n)(z-n+1)\ldots(z-1)},\qquad
n=1,2,\ldots.
\end{equation}

Отметим, что
$$\Gamma\left(\frac{1}{2}\right)=\sqrt{\pi}, $$
$$
\Gamma\left(\frac{1}{2}+n\right)=\frac{(2n)!\sqrt{\pi}}{4^nn!},\qquad
\Gamma\left(\frac{1}{2}-n\right)=\frac{(-4)^nn!\sqrt{\pi}}{(2n)!}.
$$

Имеют место следующие соотношения: формула дополнения
\begin{equation}\label{Dopol}
\Gamma(z)\Gamma(1-z)=\frac{\pi}{\sin{z\pi}},
\end{equation}
формула удвоения (формула Лежандра)
\begin{equation}\label{Lezh}
\Gamma(2z)=\frac{2^{2z-1}}{\sqrt{\pi}}\Gamma(z)\Gamma\left(z+\frac{1}{2}\right).
\end{equation}

\emph{Бета-функция} $B(z,w)$ тесно связана с гамма-функцией.  Для
двух  параметров $z$ и $w,$ удовлетворяющих условиям $\Re z>0$ и $\Re w>0,$ бета-функция Эйлера определяется
интегралом Эйлера первого рода
\begin{equation}\label{Beta}
B(z,w)=\int\limits_{0}^1t^{z-1}(1-t)^{w-1}dt.
\end{equation}
Если $\Re z\leq 0$ и $\Re w\leq 0$ не положительны,
то бета-функция определяется формулой
\begin{equation}\label{BetaGamma}
B(z,w)=\frac{\Gamma(z)\Gamma(w)}{\Gamma(z+w)}.
\end{equation}

\emph{Пси-функция} $\psi(z)$ определяется как логарифмическая
производная гамма-функции
$$
\psi(z)=\frac{d\,\ln\Gamma(z)}{dz}=\frac{\Gamma'(z)}{\Gamma(z)}.
$$
Функция $\psi(z)$  терпит разрывы второго рода в точках
$z=0,-1,-2,\ldots.$ Для пси-функции справедливо представление
$$
\psi(z)=-\gamma+(z-1)\sum\limits_{n=0}^\infty\frac{1}{(n+1)(z+n)},
$$
где
$$\gamma=\lim\limits_{m\rightarrow\infty}\left(\sum\limits_{n=1}^m\frac{1}{n}-\ln{m}\right)=0\p5772156649\ldots$$
обозначает постоянную Эйлера---Маскерони~\cite{AS}. Очевидно,
$\psi(1)=-\gamma.$ Отметим справедливость формулы
\begin{equation}\label{Psi}
\int\limits_0^1\frac{t^x-t^y}{1-t}dt=\psi(y+1)-\psi(x+1).
\end{equation}

\emph{Символ Похгаммера} $(z)_n$ при целых $n$ определяется
равенством
$$
(z)_n=z(z+1)\ldots(z+n-1),\,\, n=1,2,\ldots,\,\, (z)_0\equiv 1.
$$
Справедливы равенства $ (z)_n=(-1)^n(1-n-z)_n,$ $ (1)_n=n! $ и
\begin{equation}\label{Poh}
(z)_n=\frac{\Gamma(z+n)}{\Gamma(z)}.
\end{equation}
Равенство~\eqref{Poh} можно использовать для введения символа
$(z)_n$  при действительных (комплексных) $n.$

\emph{{Биномиальные коэффициенты}} определяются по формуле
$$
\biggl(\begin{array}{c}
$$\alpha$$ \\
$$n$$ \\
\end{array}\biggr)=\frac{(-1)^{n-1}\alpha\Gamma(n-\alpha)}{\Gamma(1-\alpha)\Gamma(n+1)}.
$$
В частности, при целых $\alpha=m,$ $m=1,2,\ldots,$ имеем равенства
$$
\biggl(\begin{array}{c}
$$m$$ \\
$$n$$ \\
\end{array}\biggr)=\frac{m!}{n!(m-n)!},\qquad m\geq n.
$$
В случае произвольных (комплексных) $\beta$ и $\alpha,$
$\alpha\neq -1,-2,\ldots,$ полагают
$$
\biggl(\begin{array}{c}
$$\alpha$$ \\
$$\beta$$ \\
\end{array}\biggr)=\frac{\Gamma(\alpha+1)}{\Gamma(\beta+1)\Gamma(\alpha-\beta+1)}.
$$
Для натуральных $k$ справедлива формула
\begin{equation}\label{bin}
(-1)^k\biggl(\begin{array}{c}
$$\alpha$$ \\
$$k$$ \\
\end{array}\biggr)=\biggl(\begin{array}{c}
$$k-\alpha-1$$ \\
$$k$$ \\
\end{array}\biggr)=\frac{\Gamma(k-\beta)}{\Gamma(-\alpha)\Gamma(k+1)}.
\end{equation}


\subsubsection{Функции Бесселя}\label{sec3.1.2}

Функции Бесселя  определяются как решения дифференциального
уравнения Бесселя:
$$
x^2 \frac{d^2 y}{dx^2} + x \frac{dy}{dx} + (x^2 - \alpha^2)y = 0,
$$
где порядок $\alpha$ "--- произвольное  комплексное число.

Функциями Бесселя первого рода, обозначаемыми $J_\alpha(x),$
являются решения, конечные в точке $x=0$ при целых или
неотрицательных $\alpha.$ Можно определить эти функции с помощью
разложения в ряд Тейлора около нуля или в более общий степенной
ряд при нецелых $\alpha$:
$$
J_\alpha(x) = \sum\limits_{m=0}^\infty \frac{(-1)^m}{m!\, \Gamma(m+\alpha+1)} {\left({\frac{x}{2}}\right)}^{2m+\alpha}.
$$
Если $\alpha$ не является целым числом, функции $J_\alpha (x)$ и
$J_{-\alpha} (x)$ линейно независимы и, следовательно, являются
решениями уравнения. Но если $\alpha$ целое, то верно следующее
соотношение:
$$
J_{-\alpha}(x) = (-1)^{\alpha} J_{\alpha}(x).
$$
Оно означает, что в этом случае функции линейно зависимы. Тогда
вторым решением уравнения станет функция Бесселя второго рода "---
\textit{функция Неймана}, то есть решение $N_\alpha(x)$ уравнения
Бесселя, бесконечное в точке $x=0.$ Эта функция связана с
$J_\alpha(x)$ следующим соотношением:
$$
N_\alpha(x) = \frac{J_\alpha(x) \cos(\alpha\pi) - J_{-\alpha}(x)}{\sin(\alpha\pi)},
$$
где в случае целого $\alpha$ берётся предел по $\alpha,$
вычисляемый, например, с помощью правила Лопиталя. Функции Неймана
также называются функциями Бесселя второго рода. Линейная
комбинация функций Бесселя первого и второго родов являет собой
полное решение уравнения Бесселя:
$$
y(x) = C_1 J_\alpha(x) + C_2 N_\alpha(x).
$$
Часто также используется обозначение $N_\alpha(x)=Y_\alpha(x).$

\textit{Модифицированными функциями Бесселя, или функциями Бесселя мнимого аргумента} называются
 функция
$$
I_\nu(x)=i^{-\nu} J_\nu(ix)
$$
и \textit{функция Макдональда }
$$
K_\nu(x)= \frac{\pi}{2\sin(\pi\nu)} \left[I_{-\nu}(x) - I_\nu(x)\right], \nu\notin\mathbb{Z}.
$$
В случае целого $ \nu\in\mathbb{Z}$ функция Макдональда вычисляется предельным переходом по индексу с помощью правила Лопиталя.

\textit{Функции Ханкеля} (Ханкеля) или функции Бесселя третьего
рода "--- это линейные комбинации функций Бесселя первого и
второго родов и, следовательно, также решения уравнения Бесселя.
Названы в честь немецкого математика Хермана Ханкеля.
$$H_{\nu}^{(1)}(z)=J_{\nu}(z)+iN_{\nu}(z)$$ "--- функция Ханкеля первого рода;
$$H_{\nu}^{(2)}(z)=J_{\nu}(z)-iN_{\nu}(z)$$ "--- функция Ханкеля второго рода.
Функции Ханкеля с индексом $0$ являются фундаментальными решениями уравнения Гельмгольца.
Представление функциями Бесселя первого рода:
$$
H_{\nu}^{(1)} (z) = \frac{J_{-\nu} (z) - e^{-\nu\pi i} J_{\nu} (z)}{i \sin (\nu\pi)},
$$
$$
H_{\nu}^{(2)} (z) = \frac{J_{-\nu} (z) - e^{\nu\pi i} J_{\nu} (z)}{- i \sin (\nu\pi)}.
$$

\textit{Нормированная функция Бесселя} ($j$-малая функция Бесселя)
$j_\nu,$ определяется  формулой
(см.~\cite[с.~10]{Kip1},~\cite{Lev1})
\begin{equation}\label{FBess1}
j_\nu(x) ={2^\nu\Gamma(\nu+1)\over x^\nu}\,\,J_\nu(x),
\end{equation}
где  $J_\nu$ "--- функция Бесселя первого рода. Для~\eqref{FBess1}
справедливо равенство (см., например,~\cite{Lev1})
\begin{equation}\label{RavDBes1}
T^y_x j_{\frac{\gamma-1}{2}}(x)=j_{\frac{\gamma-1}{2}}(x)\,j_{\frac{\gamma-1}{2}}(y).
\end{equation}

Отметим, что полезные свойства функций Бесселя с приложениями в
аналитической теории чисел  получены в работах
Н.\,В.~Кузнецова~\cite{Kuz1,Kuz2}.

\subsubsection{Гипергеометрическая функция Гаусса}\label{sec3.1.3}

Гипергеометрическая функция Гаусса определяется внутри круга
$|z|{<}1$ как сумма гипергеометрического ряда (см.~\cite[с.~373,
формула~15.3.1]{AS})
\begin{multline}\label{FG}
\,_2F_1(a,b;c;z)=F(a,b,c;z)=\sum\limits_{k=0}^\infty\frac{(a)_k(b)_k}{(c)_k}\frac{z^k}{k!}=
\\
=\frac{\Gamma(c)}{\Gamma(b)\Gamma(c-b)}\int\limits_0^{1}t^{b-1}(1-t)^{c-b-1}(1-zt)^{-a}dt,
\,\,\Re c>\Re b>0,
\end{multline}
а при $|z|\geq1$ получается аналитическим продолжением этого ряда
(см.~\cite{AKdF, Bai}). В формуле~\eqref{FG} параметры $a,b,c$ и
переменная $z$ могут быть комплексными, причём $c \not = 0, -1,
-2, \dots,$ а $(a)_k$ есть символ Похгаммера~\eqref{Poh}.

Гипергеометрический ряд~\eqref{FG} сходится только в единичном
круге комплексной плоскости, поэтому возникает необходимость
построения аналитического продолжения гипергеометрической функции
за границу этого круга, на всю комплексную плоскость. Один из
способов аналитического продолжения "--- использование
интегрального представления Эйлера
$$_2F_1(\alpha,\beta,\gamma;z) = { \Gamma(\gamma) \over
    \Gamma(\beta)\Gamma(\beta-\gamma) } \int\limits_{0}^{1}
t^{\beta-1} (1-t)^{\gamma-\beta-1} (1-tz)^{-\alpha} \,dt,
$$
$$
0<\Re \beta<\Re \gamma,\,\,\,\,|\arg(1-z)|<\pi,
$$
в котором правая часть определена при указанных условиях,
обеспечивающих сходимость интеграла.

Важным свойством гипергеометрической функции является то, что
многие специальные и элементарные функции могут быть получены из
неё при определённых значениях параметров и преобразовании
независимого аргумента.

Примеры для элементарных функций:
$$
(1+x)^n = F(-n,\beta,\beta;-x),\qquad {1 \over x} \ln(1+x) =
F(1,1,2;-x),\qquad e^x = \lim\limits_{n \to \infty} F(1,n,1;{x
\over n}),
$$
$$ \cos x = \lim\limits_{\alpha,\;\beta \to \infty}
F\left(\alpha,\beta,\frac{1}{2}; -\frac{x^2}{4 \alpha
    \beta}\right),\qquad \cosh x = \lim\limits_{\alpha,\;\beta \to \infty}
F\left(\alpha,\beta,\frac{1}{2};{ x^2 \over 4 \alpha
    \beta}\right).
$$

Функция Бесселя первого рода и гипергеометрическая функция Гаусса
связаны формулой
$$ J_\nu(z)= \lim\limits_{\alpha,\;\beta \to \infty} \left[
\frac{\left(\dfrac{z}{2}\right)^\nu}{\Gamma(\nu+1)}
F\left(\alpha,\beta,\nu+1; -\frac{z^2}{4 \alpha \beta}\right)
\right].
$$

\subsubsection{Функции Лежандра}\label{sec3.1.4}

Функции Лежандра  $P_{\nu}^\mu (x), Q_{\nu}^\mu (x)$ "--- это
обобщения полиномов Лежандра на нецелую степень (см.~\cite{BE1}).

 Функции Лежандра $P_{\nu}^\mu (x), Q_{\nu}^\mu (x)$  являются решениями общего уравнения Лежандра
 $$
 (1-x^{2})\,y''-2xy'+\left[\lambda (\lambda +1)-{\frac {\mu ^{2}}{1-x^{2}}}\right]\,y=0,\,
$$
где комплексные числа $\lambda$ и $\mu$ называются степенью и порядком  функций Лежандра, соответственно. Эти функции могут быть определены для комплексных значений параметров и аргумента:
$$
 P_{\lambda }^{\mu }(z)={\frac {1}{\Gamma (1-\mu )}}\left[{\frac {1+z}{1-z}}\right]^{\mu /2}\,_{2}F_{1}(-\lambda,\lambda +1;1-\mu;{\frac {1-z}{2}}),\qquad |1-z|<2,
$$
$$ Q_{\lambda }^{\mu }(z){=}{\frac {{\sqrt {\pi }}\ \Gamma (\lambda {+}\mu {+}1)}{2^{\lambda +1}\Gamma (\lambda {+}\frac{3}{2})}}{\frac {e^{i\mu \pi }(z^{2}-1)^{\mu /2}}{z^{\lambda +\mu +1}}}\,_{2}F_{1}\left({\frac {\lambda +\mu +1}{2}},{\frac {\lambda +\mu +2}{2}};\lambda +{\frac {3}{2}};{\frac {1}{z^{2}}}\right),\quad |z|>1,
$$
где $\Gamma$ "--- гамма-функция и $_{2}F_{1}$ "---
гипергеометрическая функция. В книге используются также прямые
значения  функций $P_{\nu}^\mu (x), Q_{\nu}^\mu (x)$ на разрезе
$x\in [-1;1],$ которые обозначаются   $\mathbb{P}_{\nu}^\mu (x),
\mathbb{Q}_{\nu}^\mu (x).$

Функции Лежандра являются частными случаями гипергеометрической функции Гаусса, поэтому для них известны многочисленные разложения в ряды, интегральные представления, формулы продолжения, а также выражения через элементарные функции при специальных значениях параметров.

Многочлены Лежандра являются функциями Лежандра порядка $\mu=0$
при целых неотрицательных $\lambda=n.$

\subsubsection{Функция Миттаг-Лефлера}\label{sec3.1.5}

Функция Миттаг-Лефлера $E_{\alpha,\beta}(z)$ "--- это целая
функция (в $z\in\mathbb{C}$) порядка $1/\alpha,$ определяемая
степенным рядом (см.~\cite{Dzh1, Djr2, BE3, GKMR, Kir5,Kir6, Jor,
ViGa})
\begin{equation}\label{ML}
E_{\alpha,
\beta}(z)=\sum\limits_{n=0}^{\infty}\frac{z^n}{\Gamma(\alpha
n+\beta)},\;\; z\in\mathbb{C},\;\alpha, \beta\in\mathbb{C},\;
\Re\alpha>0,\; \Re\beta>0.
\end{equation}

Функция~\eqref{ML} была введена Гёстой Миттаг-Лефлером в 1903~г.
для  $\alpha=1$ и А.~Виманом в 1905~г. в общем случае. Первыми
приложениями этих функций у Миттаг-Лефлера и Вимана были
приложения в ТФКП (нетривиальные примеры целых функций с нецелыми
порядками роста и обобщённые методы суммирования). В СССР эти
функции стали в основном известны после опубликования знаменитой
монографии М.\,М.~Джрбашяна~\cite{Dzh1}, см. также его более
позднюю монографию~\cite{Djr2}. Самым известным применением
функций Миттаг-Лефлера в теории интегродифференциальных уравнений
и дробном исчислении является тот факт, что через них в явном виде
выражается резольвента дробного интеграла Римана---Лиувилля по
знаменитой формуле Хилле---Тамаркина---Джрбашяна~\cite{HT, SKM,
Pshu1}. Ввиду многочисленных приложений к решению дифференциальных
уравнений дробного порядка эта функция была заслуженно названа
в~\cite{GoMa} <<\textit{королевской функцией дробного
исчисления}>>.

Производная функции Миттаг-Леффлера вычисляется по формуле
$$
E_{\alpha,\beta}'(z)=\frac{d\,E_{\alpha,\beta}(z)}{dz}=\sum\limits_{k=0}^\infty\frac{(1+k)z^k}{\Gamma(\beta+\alpha
    (1+k))}.
$$
Отметим, что $E_{1,1}(z)=e^z.$

\subsubsection{Обобщённые функции гипергеометрического
типа}\label{sec3.1.6}

Функции гипергеометрического типа определяются комплексным
интегралом типа Меллина---Барнса~\cite{Marich1,PBM}.
 Основные обобщённые функции гипергеометрического типа  рассматриваются в монографиях~\cite{Marich1, Sla, VN, Dwo1, Ext, KaSr, KiSa} и ряде других.

Фундаментальное значение имеет метод вычисления определённых
интегралов с использованием обобщённых гипергеометрических
функций, основанный на теореме Слейтер---Маричева,
см.~\cite{Marich1,PBM}. Для суммирования достаточно произвольных
рядов в терминах обобщённых гипергеометрических и дзета функций
оригинальный метод был разработан В.\,С.~Рыко,
см.~\cite{Ryko1,Ryko2}, сейчас он практически забыт и не
используется. Современные методы теоретического и компьютерного
суммирования на основе гипергеометрических функций рассматриваются
в известных монографиях~\cite{PWZ,Koepf}.

Известные многочисленные обобщения гипергеометрических функций,
например, $q$-ги\-пер\-гео\-мет\-ри\-чес\-кие функции, теория которых
восходит к Э.~Гейне, $A$ "--- гипергеометрические системы
Гельфанда---Граева, эллиптические гипергеометрические функции
В.\,П.~Спиридонова и другие.

Следует отметить, что в последнее время функции
гипергеометрического типа, такие как функции Миттаг-Лефлера, Райта
и Фокса, получили многочисленные приложения в задачах теории
вероятностей и математической статистики. В частности, через них
выражаются некоторые важные функции распределения, плотности
вероятностей и их характеристики.

\textit{Функции Райта} были введены Эдвардом Мэйтлендом Райтом в
1935~г., см.~\cite{BE1, ViGa, Jor, Kir1, Kir5,Kir6}. Являются
непосредственными обобщениями функций Миттаг-Лефлера и
определяются в виде рядов, аналогичным~\eqref{ML}, но с
произвольным конечным числом гамма-функций в числителе и
знаменателе дроби для общего члена ряда.

Пусть $p,q\in\mathbb{N}_0=\{0,1,2,\ldots\},$ $p^2+q^2\neq 0,$
$a_i,b_j\in\mathbb{C},$ $\alpha_i,\beta_j\in\mathbb{R}$
($\alpha_i,\beta_j\neq0;$ $i=1,2,\ldots,p;$ $j=1,2,\ldots,q$).
Тогда функция Райта определяется степенным рядом
$$
\,_p\Psi_q(z)=\sum\limits_{m=0}^\infty
\frac{\prod\limits_{i=1}^p\Gamma(a_i+\alpha_ik)}{\prod\limits_{j=1}^q\Gamma(b_j+\beta_jk)}
\frac{z^k}{k!},\qquad z\in\mathbb{C}.
$$

Иногда эти функции используются под другими названиями, например,
обобщённая функция Бесселя~\cite{Jor},  обобщённая функция
Миттаг-Лефлера~\cite{GKMR}, многопараметрическая функция
Миттаг-Лефлера~\cite{Kir5,Kir6, ViGa}, и некоторыми другими.
Функции Райта, представленные в образах преобразования Меллина или
интегралами Меллина---Барнса, также являются функциями Фокса и
наоборот, хотя не при всех значениях параметров такие
преобразования очевидны или даже возможны. Для этого класса
функций используются также названия многопараметрические
гипергеометрические функции, обобщённые или многопараметрические
функции Миттаг-Лефлера, многопараметрические функции Бесселя и ряд
других. Иногда используется неграмотное название <<функции
Бесселя---Мэйтленда>> для одного подобного класса функций, при
этом за имя никогда не существовавшего математика выдаётся средняя
часть имени Райта, к сожалению, эта неточность воспроизведена в
такой уважаемой книге, как~\cite{BE2}.

В качестве примера приложений отметим, например, что при
представлении резольвенты для дробных степеней оператора Бесселя
используется такой вариант функции Райта~\cite{S700,SS}:
\begin{equation}\label{GML}
J_{\gamma,\lambda}^\mu(z)=\sum\limits_{m=0}^\infty\frac{(-1)^m}
{\Gamma(\gamma+m\mu+\lambda+1)\Gamma(\lambda+m+1)}\left(\frac{z}{2}\right)^{2m+\gamma+2\lambda}.
\end{equation}

\textit{$G$-функции Майера}~\cite{Dwo2, Sie}. Эта функция была
введена в 1946~г. голландским математиком Корнелиусом Симоном
Майером, он был учеником Ван дер Корпута, его учеником был Боэль
Брааксма.

Общее определение $G$-функции Майера дается следующим интегралом
в комплексной плоскости:
$$
G_{p,q}^{\,m,n} \!\left( \left. \begin{matrix} a_1, \dots, a_p \\
b_1, \dots, b_q \end{matrix} \; \right| \, z \right) =\frac{1}{2
\pi i} \int\limits_L \frac{\prod\limits_{j=1}^m \Gamma(b_j - s)
\prod\limits_{j=1}^n \Gamma(1 - a_j +s)} {\prod\limits_{j=m+1}^q
\Gamma(1 - b_j + s) \prod\limits_{j=n+1}^p \Gamma(a_j - s)} \,z^s
\,ds.
$$
Это интеграл так называемого типа Меллина---Барнса. Определение
корректно при следующих предположениях: $ 0 \leq m' \leq q,$
$0\leq n\leq p,$ $m, n, p, q\in\mathbb{Z},$ $a_k- b_j\neq 1, 2,
3,\ldots$ для $k = 1, 2,\ldots, n$  и $j = 1, 2,\ldots, m,$ что
означает, что никакие полюса $\Gamma(b_j-s),$ $j = 1, 2,\ldots,
m,$ не совпадают ни с какими  полюсами $\Gamma(1-a_k+s),$ $k = 1,
2,\ldots, n,$ $z\neq0.$ Кроме того, контур интегрирования после
обхода области с полюсами должен продолжаться по одной
вертикальной прямой вверх и вниз.

В случае другого выбора контура с уходом на бесконечность по
горизонтальным прямым, получаются новые варианты функций Майера,
подробный анализ этих модификаций, а также их отличий от
классического случая см. в недавних работах~\cite{KaPr3, KaLo1,
KaLo2}. В этих случаях, в частности, функция Майера может быть
определена, но для неё может не существовать преобразование
Меллина, или даже преобразование Меллина может существовать, но
сама функция не восстанавливаться обратным преобразованием Меллина
через интеграл Меллина---Барнса.

Важным свойством  $G$-функции Майера является тот факт, что во
всех случаях она представима в виде конечной суммы обобщённых
гипергеометрических функций. Кроме того, $G$-функция Майера
удовлетворяет явно выписываемому дифференциальному уравнению с
полиномиальными коэффициентами. Отметим, что  функции Майера
используются, например, в~\cite{Dim, Kir1} при представлениях
преобразования Обрешкова.

В приложениях выделяется специальный случай функции Майера
$G_{p,0}^{\,p,p},$ он был подробно изучен в работах Нёрлунда,
поэтому в ряде недавних работ предлагается для этого случая
название: функция Майера---Нёрлунда, см.~\cite{KaPr3, KaLo1,
KaLo2}.

\textit{Функции Фокса}~\cite{KiSa}. $H$-функцию ввел английский
математик Чарльз Фокс  в 1961~г., см.~\cite{KiSa, Bra}.

$H$-функция также определяется с помощью  интеграла
Меллина---Барнса:
$$
H(x)=H^{m,n}_{p,q}(z)=H^{m,n}_{p,q}\left[z\biggl|\left(
\begin{array}{c}
$$(a_p,A_p)$$ \\
$$(b_q,B_q)$$ \\
\end{array}
\right) \right]= H_{p,q}^{\,m,n} \!\left[ z \left| \begin{matrix}
(a_1, A_1) & (a_2, A_2) & \ldots & (a_p, A_p) \\
(b_1, B_1) & (b_2, B_2) & \ldots & (b_q, B_q) \end{matrix} \right.
\right] =
$$
$$
=\frac{1}{2\pi i}\int\limits_L
\frac
{\left(\prod\limits_{j=1}^m\Gamma(b_j+B_js)\right)\left(\prod\limits_{j=1}^n\Gamma(1-a_j-A_js)\right)}
{\left(\prod\limits_{j=m+1}^q\Gamma(1-b_j-B_js)\right)\left(\prod\limits_{j=n+1}^p\Gamma(a_j+A_js)\right)}
z^{-s} \, ds.
$$
Здесь пустое произведение всегда интерпретируется как единица;
\\$m,n,p,q\in \mathbb{N}_0,$ $0\leq n\leq p,$ $1\leq m\leq q,$
$A_i, B_j\in \mathbb{R}_+,$ $a_i, b_j\in\mathbb{R}$ или
$\mathbb{C},$ $i=1,\ldots,p;$ $j=1,\ldots,q,$ $L$ "--- некоторый
подходящий контур, отделяющий полюса
$$
\zeta_{j\nu}=-\left(\frac{b_j+\nu}{B_j}\right),\qquad
j=1,\ldots,m;\qquad \nu=0,1,2,\ldots
$$
гамма-функции $\Gamma(b_j+sB_j)$ от полюсов
$$
\omega_{\lambda
k}=\left(\frac{1-a_\lambda+k}{A_\lambda}\right),\qquad
\lambda=1,\ldots,n;\qquad k=0,1,2,\ldots
$$
гамма-функции $\Gamma(1-a_\lambda-sA_\lambda)$ так чтобы
$$
A_\lambda(b_j+\nu)\neq B_j(a_\lambda-k-1),  j=1,\ldots,m;
\lambda=1,\ldots,n; \nu,k=0,1,2,\ldots
$$
Кроме того, контур интегрирования после обхода области с полюсами должен продолжаться по одной вертикальной прямой вверх и вниз.

В случае другого выбора контура с уходом на бесконечность по
горизонтальным прямым, получаются новые варианты функций Фокса,
подробный анализ этих модификаций, а также их отличий от
классического случая см. в недавних работах~\cite{KaPr1,KaPr2}. В
этих случаях, в частности, функция Фокса может быть определена, но
для неё может не существовать преобразование Меллина, или даже
преобразование Меллина может существовать, но сама функция не
восстанавливаться обратным преобразованием Меллина через интеграл
Меллина---Барнса.

Особым случаем, для которого $H$-функция Фокса сводится к
$G$-функции Майера, является набор параметров $A_j=B_k=C,$ $C>0$
для $j =1,\ldots,p,$ $k=1,\ldots,q$:
$$
H_{p,q}^{\,m,n} \!\left[ z \left| \begin{matrix}
(a_1, C) & (a_2, C) & \ldots & (a_p, C) \\
(b_1, C) & (b_2, C) & \ldots & (b_q, C) \end{matrix} \right.
\right] = \frac{1}{C} G_{p,q}^{\,m,n} \!\left( \left.
\begin{matrix} a_1, \dots, a_p \\
b_1, \dots, b_q \end{matrix} \; \right| \, z^{1/C} \right).
$$

В отличии от  $G$-функции Майера $H$-функция Фокса в общем случае
не представляется в виде конечной суммы обобщённых
гипергеометрических функций или $G$-функций. В том случае, когда
она представляется в виде ряда, используется название функция
Райта---Фокса. Принципиальным отличием функции Фокса от функции
Райта является тот факт, что она не является решением никакого
дифференциального уравнения конечного порядка с полиномиальными
коэффициентами (известная недоказанная гипотеза, в которую верят
большинство специалистов).

Самым известным применением функций Фокса в теории
дифференциальных уравнений является их использование
А.\,Н.~Кочубеем при построения функции Грина для уравнения дробной
диффузии~\cite{Koch1,Koch2, EiIvKoch}. При этом по существу
использовалось мало известное большинству математиков, но
чрезвычайно полезное преобразование Станковича.

Отметим,  что используемые при определениях функций Майера и Фокса
интегралы типа Меллина---Барнса были впервые введены, по-видимому,
итальянским математиком Сальваторе Пинкерле в 1888~г. (работы
Эрнеста Уильяма Барнса были выполнены в 1908--1910~гг., работы
Ялмара Меллина "--- в 1895 и 1909~гг.).

Представляется, что для введённых выше функций исторически
наиболее правильным является общее название <<функции
Райта---Фокса>>  с указанием их формы представления: или как
степенной ряд в форме Райта, или представление через комплексный
интеграл в форме Фокса. Тем не менее, приняты такие названия:
функция Райта для определения в виде ряда, функция Майера или
функция Фокса для определений в виде комплексного интеграла,
функция Райта---Фокса для тех случаев, когда функция Фокса
представляется рядом типа Райта.

Кроме функций гипергеометрического типа от одной переменной
используются гипергеометрические функции от нескольких переменных:
Горна, Аппеля, Лауричелы, Кампе де Ферье и ряд других. Эти функции
играют существенную роль в теории вырождающихся и сингулярных
дифференциальных уравнения, см., например,~\cite{VN,
KarST2,KarST3, Erg1, CSh, HaKa1, KarST, SaHa1, ShiE5}.

\subsection{Функциональные пространства}\label{sec3.2}

Приведем некоторые понятия и утверждения, которые далее
потребуются для введения операторов преобразования и изучения их
свойств. Эти сведения взяты  из известного учебника~\cite{KF},
монографии~\cite{SKM}, а также из~\cite{Bur1, Bur2}.

Как обычно, через $\mathbb{R}$ будем обозначать множество
вещественных чисел, а через $\mathbb{C}$ "--- множество
комплексных чисел. Рассматриваемые далее множества и функции мы,
не оговаривая это отдельно, считаем измеримыми.

\subsubsection{Гёльдеровы функции, абсолютно непрерывные функции,
    класс $AC^n$}\label{sec3.2.1}

В этом пункте определим (см.~\cite{KF} и~\cite{SKM}) локальное и
глобальное условия Гёльдера, а также классы функций $AC$ и $AC^n.$

Пусть $\Omega=[a,b],$ $-\infty{<}a{<}b{<}\infty$ обозначает
отрезок вещественной оси.

\begin{definition}
    Говорят, что
    функция $f(x)$ удовлетворяет на $\Omega$
    \textit{условию Гёльдера} порядка $\lambda,$ если
    \begin{equation}\label{Ge}
    |f(x_1)-f(x_2)|\leq A|x_1-x_2|^\lambda
    \end{equation}
    для всех $x_1,\,x_2\in\Omega,$ где $A$ "--- постоянная, а $\lambda$ "--- показатель Гёльдера.
\end{definition}

\begin{definition}
    Через $H^\lambda=H^\lambda(\Omega),$ обозначается класс всех
    (вообще говоря, комплекснозначных) функций, удовлетворяющих на
    $\Omega$ условию Гёльдера фиксированного порядка $\lambda.$
\end{definition}

При $\lambda> 1$ класс $H^\lambda$
содержит только постоянные $f(x)\equiv \const$:
$$
|f'(x)|=\left|\lim\limits_{\Delta
x\to0}\frac{f(x_1)-f(x_2)}{x_1-x_2} \right| \le
A\lim\limits_{\Delta
    x\to0}\frac{|x_1-x_2|^\lambda}{|x_1-x_2|} | \le
A|x_1-x_2|^{\lambda-1}=0.
$$
Поэтому класс $H^\lambda$ интересен лишь в случае $0<\lambda\leq
1.$

Класс $H^1(\Omega)$ называют липшицевым классом. Приведем
определение более широкого, чем $H^1(\Omega),$ класса $AC(\Omega)$
абсолютно непрерывных функций.

\begin{definition}
    Функция $f(x)$  называется \textit{абсолютно непрерывной} на
    отрезке $\Omega,$ если для любого $\varepsilon>0$ можно найти
    такое $\delta>0,$ что для любой конечной системы попарно
    непересекающихся отрезков $[a_k,b_k]\subset\Omega,$ $k=1,2,\ldots,n,$
    такой, что $$\sum\limits_{k=1}^n(b_k-a_k)<\delta,$$ справедливо
    неравенство $$\sum\limits_{k=1}^n|f(b_k)-f(a_k)|<\varepsilon.$$
    Класс всех таких функций обозначается $AC(\Omega).$
\end{definition}

Известно (см.~\cite{KF}), что класс $AC(\Omega)$ совпадает с
классом первообразных от суммируемых по Лебегу функций, т.~е.
\begin{equation}\label{AC}
f(x)\in AC(\Omega)\Leftrightarrow
f(x)=c+\int\limits_a^x\varphi(t)dt,\,\,\,
\int\limits_a^b|\varphi(t)|dt<\infty,\,\,\varphi(t)=f'(t).
\end{equation}
Поэтому абсолютно непрерывные функции имеют почти всюду
суммируемую производную $f'(x)$ (обратно, из существования почти
всюду суммируемой производной еще не вытекает абсолютная
непрерывность).

\begin{definition}
    Через $AC^n(\Omega),$ где $n=1,2,\ldots,$ и $\Omega$ "--- отрезок,
    обозначим класс функций $f(x),$ непрерывно дифференцируемых на
    $\Omega$ до порядка $n-1,$ причём $f^{(n-1)}(x)\in AC(\Omega).$
\end{definition}

Очевидно, что $AC^1(\Omega)=AC(\Omega)$ и класс $AC^n(\Omega)$
состоит из функций, представимых $n$-кратным интегралом Лебега с
переменным верхним пределом от суммируемой функции с заменой
постоянной в~\eqref{AC} на многочлен порядка $n-1$:
\begin{equation}\label{ACn}
f(x)\in AC^n(\Omega)\Leftrightarrow
f(x)=\sum\limits_{k=0}^{n-1}c_k(x-a)^k+\underbrace{\int\limits_a^xdt\ldots\int\limits_a^{x}dt\int\limits_a^{x}}_n\varphi(t)dt,
\end{equation}
$$
\int\limits_a^b|\varphi(t)|dt<\infty,\,\,
c_k=\frac{f^{(k)}(a)}{k!},\,\,\varphi(t)=f^{(n)}(t).
$$

Пусть теперь $\Omega$ "--- ось или полуось. В этом случае при
определении класса $H^\lambda(\Omega)$ дополнительно оговаривается
<<гёльдеровское>>, поведение на бесконечности. Именно, говорят,
что функция $f(x)$ удовлетворяет условию Гёльдера в окрестности
бесконечно удаленной точки, если
\begin{equation}\label{Geld}
\biggl|f\left(\frac{1}{x_1}\right)-f\left(\frac{1}{x_2}\right)\biggr|\leq
A\biggl|\frac{1}{x_1}-\frac{1}{x_2}\biggr|^\lambda
\end{equation}
для всех $x_1,\,x_2,$ достаточно больших по абсолютной величине.

\begin{definition}
    Пусть $\Omega$ "--- ось или полуось. Через
    $H^\lambda=H^\lambda(\Omega)$ обозначается класс функций,
    удовлетворяющих условию Гёльдера~\eqref{Ge} на любом конечном
    отрезке в $\Omega$ и условию~\eqref{Geld} в окрестности бесконечно
    удаленной точки.
\end{definition}

Совокупность двух условий~\eqref{Ge}  и~\eqref{Geld}, определяющих
класс $H^\lambda(\Omega)$ для бесконечного интервала $\Omega,$
равносильна одному условию
\begin{equation}\label{GlGe}
|f(x_1)-f(x_2)|\leq
A\frac{|x_1-x_2|^\lambda}{(1+|x_1|)^\lambda(1+|x_2|)^\lambda}.
\end{equation}

Условие~\eqref{GlGe} называется <<{\it глобальным условием
    Гёльдера}>>.

\subsubsection{Класс $L_p$ и его свойства}\label{sec3.2.2}

Введем класс суммируемых в $p$-й степени функций  и приведем
некоторые неравенства и теоремы, справедливые для функций из этого
класса (см.~\cite{KF} и~\cite{SKM}). Эти пространства были введены
Ф.~Риссом, в последнее время важные применения нашли их обобщения
на случай переменного показателя $p=p(x),$ см.,
например,~\cite{KMRS,Rad2}.

Пусть теперь $\Omega=[a,b],$ где $-\infty\leq a<b\leq\infty.$

\begin{definition}
    Через $L_p=L_p(\Omega)$ обозначается множество всех измеримых на
    $\Omega$ функций $f(x),$ вообще   говоря, комплекснозначных, для
    которых выполнено неравенство $$\int\limits_{\Omega}|f(x)|^pdx <
    \infty,\qquad1 \leq p <\infty.$$
\end{definition}

Норма в $L_p(\Omega)$ определяется формулой
$$
||f||_{L_p(\Omega)}=\left(\int\limits_{\Omega}|f(x)|^pdx\right)^{1/p}.
$$

Две  отличающиеся на множестве меры нуль функции не различаются
как элементы пространства $L_p(\Omega).$

Для функций $f\in L_p(\Omega)$ и $g\in L_p(\Omega)$ справедливо
\emph{неравенство Минковского }
$$
||f+g||_{L_p(\Omega)}\leq ||f||_{L_p(\Omega)}+||g||_{L_p(\Omega)},
$$
с учётом которого $L_p(\Omega)$  является нормированным
пространством.

Если же  $f(x)\in L_p(\Omega),$ $g(x)\in L_{p'}(\Omega),$ где
$p'=\dfrac{p}{p-1},$ то выполняется \emph{неравенство Гёльдера }
(или точнее "--- неравенство Роджерса---Гёльдера---Рисса)
$$
\int\limits_{\Omega}|f(x)g(x)|dx\leq||f||_{L_p(\Omega)}||g||_{L_{p'}(\Omega)},\qquad
p'=\frac{p}{p-1}.
$$
Известно, что $L_p(\Omega)$ "--- полное пространство.

Нам также потребуется теорема, позволяющая менять
порядок   интегрирования в повторных
интегралах.

\begin{theorem*}[Фубини] Пусть $\Omega_1 =[a, b],$ $\Omega_2 = [c, d],$ $-\infty\leq a<b\leq\infty,$ $-\infty\leq c<d\leq\infty,$
    и пусть $f(x, y)$ "--- определенная на $\Omega_1\times\Omega_2$
    измеримая функция. Если сходится {\rm (}абсолютно{\rm )} хотя бы один из
    интегралов
    $$
    \int\limits_{\Omega_1}dx\int\limits_{\Omega_2}f(x,y)dy,\qquad\int\limits_{\Omega_2}dy\int\limits_{\Omega_1}f(x,y)dx,\qquad\iint\limits_{\Omega_1\times\Omega_2}f(x,y)dxdy,
    $$ то они совпадают.
\end{theorem*}

Имеет место следующий частный случай теоремы Фубини:
\begin{equation}\label{Dirihle}
\int\limits_{a}^bdx\int\limits_{a}^xf(x,y)dy=\int\limits_{a}^bdy\int\limits_{y}^bf(x,y)dx
\end{equation}
в предположении, что абсолютно сходится один    из   этих
интегралов. Равенство~\eqref{Dirihle} называется \emph{формулой
Дирихле}.

Справедливо также обобщённое неравенство Минковского:
$$
\left(\int\limits_{\Omega_1}dx\left|\int\limits_{\Omega_2}f(x,y)dy\right|^p\right)^{1/p}\leq
\int\limits_{\Omega_2}dy\left(\int\limits_{\Omega_1}|f(x,y)|^pdx\right)^{1/p}.
$$

\subsubsection{Пространства $C^m_{ev},$ $C^\infty_{ev}$ и
$L_p^\gamma$}\label{sec3.2.3}

Функция $\varphi$ называется основной функцией, если она
бесконечно дифференцируема, чётна и удовлетворяет оценкам
$$
|D^q \varphi(x)|\leq\frac{C_{qr}}{(1+x^2)^r}.
$$
Здесь $q\geq 0$ и $r\geq 0$ "--- произвольные целые числа,
$C_{qr}$ "--- постоянные, не зависящие от $x.$ Ясно, что к
основным функциям можно применять оператор Бесселя сколько угодно
раз (так как если $\varphi$ "--- основная функция, то и
$\dfrac{\varphi'(x)}{x}$ тоже основная функция), причём
выполняются оценки
$$
|B^q_\gamma \varphi(x)|\leq\frac{A_{qr}}{(1+x^2)^r}.
$$
при любых целых $r\geq 0$ и $q\geq 0.$ Обозначим через $S$
пространство всех основных функций.

Обозначим через $L_{2} (\mathbb{R}_{+}^1)$ гильбертово
пространство функций $f(y),$ $y>0,$ для которых конечна норма
\begin{equation}
\| f \|_{L_{2} \lr{\mathbb{R}^1_+}} = \lr{\int\limits_0^{\infty} |f(y)|^2 \, dy}^{\frac{1}{2}}.
\label{141}
\end{equation}

Обозначим через $L_{2, \nu} (\mathbb{R}_{+}^1),$ $\nu \geq -
\dfrac{1}{2},$ весовое гильбертово   пространство функций $f(y),$
$y>0,$ для которых конечна норма
\begin{equation}
\| f \|_{L_{2, \nu} \lr{\mathbb{R}^1_+}} = \lr{\int\limits_0^{\infty} |f(y)|^2 y^{2 \nu +1} \, dy}^{\frac{1}{2}}.
\label{1.4.1}
\end{equation}

Хорошо известно, что преобразование Фурье---Бесселя унитарно в\\
$L_{2, \nu} (\mathbb{R}_{+}^1)$ и справедливо равенство
Парсеваля
\begin{equation}
\|F_{\nu} f \|_{L_{2, \nu} \lr{\mathbb{R}^1_+}} =  \| f \|_{L_{2, \nu} \lr{\mathbb{R}^1_+}}.
\label{1.4.2}
\end{equation}
Функциональное пространство $H^s_{\nu,+} (\mathbb{R}_{+}^1)$
(пространство И.\,А.~Киприянова), $s \geq 0,$ $\nu \geq -
\dfrac{1}{2},$    введенное в работе~\cite{Kip2}, определяется как
замыкание по норме
\begin{equation}\label{KipSp}
\| f \|_{H_{ \nu, +}^s \lr{\mathbb{R}^1_+}} = \frac{1}{2^{\nu}\, \Gamma (\nu +1)} \| (1+\eta^2)^{\frac{s}{2}} F_{\nu} f \|_{L_{2, \nu} \lr{\mathbb{R}^1_+}}
\end{equation}
множества чётных по И.\,А.~Киприянову функций
$\mathring{C}^{\infty}_{+} (\ov{\mathbb{R}_{+}^{1}})$
(см.~\cite{Kip2}).   Предположение о    чётности    здесь
существенно, поскольку в противном случае норма~\eqref{KipSp}
может быть равной бесконечности. Для пространств Киприянова
возможны и другие эквивалентные определения нормы.

Определим пространство С.\,Л.~Соболева $\mathring{H}^s (0, R),$ $s
\geq 0,$ $0<R<\infty,$ как замыкание множества
$\mathring{C}^{\infty} [0, R)$ по норме
$$
\| f \|_{\mathring{H}^s (0, R)} =  \|D^s  f \|_{L_{2, \nu} (0, R)}.
$$

В данной книге вводится целый ряд новых функциональных пространств, как Банаха, так и Фреше, приспособленных для решения краевых задач для дифференциальных уравнений с частными производными при наличии особенностей в коэффициентах уравнений и в допускаемых решениях.

\subsection{Основные интегральные преобразования}\label{sec3.3}

Основные определения интегральных преобразований и их свойства
можно найти в монографиях и справочниках~\cite{BE1, BE2, BE3,
PBM123, Ome, Ahi2, KGS, Pou, Dzh1, KiSa, Marich1, KiSa, GSS}.

Досадные ошибки в обращении некоторых интегральных преобразований,
которые допущены в~\cite{BE1, BE2, BE3}, исправлены
в~\cite{Kuz2,Kuz3}.

\subsubsection{Преобразование Фурье, синус и косинус преобразования, преобразование
Ханкеля}\label{sec3.3.1}

Подробнее о рассматриваемых преобразованиях см.~\cite{Tit1, Dzh1,
Pou}. Важные свойства преобразований Ханкеля (Фурье---Бесселя)
установлены в~\cite{Lar1, Mur3, Mur4, Mur5, Mur6, Mur7}.

\textit{Преобразование Фурье, синус и косинус преобразования, преобразование Ханкеля}  имеют, соответственно, вид
\begin{eqnarray*}
    (Ff)(t)=\frac{1}{\sqrt{2\pi}}\int\limits_0^\infty \exp(-ity)f(y)\,dy,\\
    (F_c f)(t)=\sqrt{\frac{2}{\pi}}\int\limits_0^\infty \cos(ty)f(y)\,dy,\\
    (F_s f)(t)=\sqrt{\frac{2}{\pi}}\int\limits_0^\infty \sin(ty)f(y)\,dy.
\end{eqnarray*}
При введённых нормировках все эти преобразования унитарны в $L_2(0,\infty)$ и совпадают с обратными.

Преобразование Ханкеля (Ханкеля, Фурье---Бесселя) имеет вид
$$
(H_\nu f)(t)=\frac{1}{t^\nu}\int\limits_{0}^{\infty}  J_\nu(ty)f(y)\,dy
$$
или
$$(H_\nu f )(\xi)=\int\limits_{0}^{\infty} {j}_{\frac{\nu-1}{2}} (x\xi)\,
f(x)x^\nu\,dx.
$$
В связи со свойством ${j}_{\nu}(0)=1$ удобней использовать преобразование Ханкеля с функцией ${j}_{\frac{\nu-1}{2}} (x\xi)$ в ядре.

\subsubsection{Преобразование Меллина. Теорема
Слейтер}\label{sec3.3.2}

При вычислении интеграла от произведения гипергеометрических
функций будем пользоваться методом из~\cite{Marich1, PBM},
основанном на применении преобразования Меллина. Преобразование
Меллина и теорема Слейтер  изучены в~\cite{Marich1, PBM, PBM123}.

Преобразованием Меллина функции $f(x)$ называется функция $g(s),$
которая определяется по формуле
\begin{equation}
\label{1710}
g(s)=M{f}(s)=\int\limits_0^\infty x^{s-1} f(x)\,dx.
\end{equation}
Определим также свёртку Меллина
\begin{equation}
\label{1711}
(f_1*f_2)(x)=\int\limits_0^\infty  f_1\left(\frac{x}{y}\right) f_2(y)\,\frac{dy}{y},
\end{equation}
при этом оператор свёртки с ядром $K$ действует в образах преобразования Меллина как умножение на мультипликатор
\begin{eqnarray}
\label{1712}
M{Af}(s)=\int\limits_0^\infty  K\left(\frac{x}{y}\right) f(y)\,\frac{dy}{y}=M{K*f}(s)
=m_A(s)M{f}(s),\\\nonumber m_A(s)=M{K}(s).\phantom{1111111111111111111}
\end{eqnarray}

Для изучения операторов типа свёртки Меллина~\eqref{1712} одним из
авторов в~\cite{S6, S66} был предложен удобный алгебраический
подход, который не содержит ничего нового, но в удобной форме
позволяет быстро получать нужные оценки. Полезные факты будут
собраны вместе как

\begin{theorem}  \label{1tMel}
Пусть оператор свёртки $A$ действует по формуле~\eqref{1712} в
образах преобразования Меллина как умножение на мультипликатор.
Тогда справедливы следующие условия ограниченности для прямого и
обратного операторов и формулы для их норм.
    \begin{enumerate}
    \item[а)] Для того, чтобы он допускал расширение до ограниченного оператора в $L_2(0,\infty)$ необходимо и достаточно, чтобы
    \begin{equation}
    \label{1713}
    \sup\limits_{\xi\in\mathbb{R}} \left|m_A\left( i\xi+\frac{1}{2}\right) \right|=M_2<\infty,
    \end{equation}
    при этом $\|A\|_{L_2}=M_2.$

    \item[б)] Для того, чтобы он допускал расширение до ограниченного оператора в $L_p(0,\infty),$ ${p>1}$ при дополнительном условии неотрицательности ядра $K$ необходимо и достаточно, чтобы
    \begin{equation}
    \label{1714}
    \sup\limits_{\xi\in\mathbb{R}}\left|m_A\left( i\xi+\frac{1}{p}\right) \right|=M_p<\infty,
    \end{equation}
    при этом $\|A\|_{L_p}=M_p.$

    \item[в)] Обратный оператор $A^{-1}$ действует также по формуле~\eqref{1712} с мультипликатором~$\dfrac{1}{m_A}.$ Для того, чтобы он допускал расширение до ограниченного оператора в $L_2(0,\infty),$ необходимо и достаточно, чтобы
    \begin{equation}
    \label{1715}
    \inf\limits_{\xi\in\mathbb{R}}\left|m_A\left( i\xi+\frac{1}{2}\right) \right|=m_2>0,
    \end{equation}
    при этом $\|A^{-1}\|_{L_2}=\dfrac{1}{m_2}.$

    \item[г)] Пусть операторы $A,A^{-1}$ определены и ограничены в $L_2(0,\infty).$ Они унитарны тогда и только тогда, когда выполнено равенство
    \begin{equation}
    \label{1716}
    \left|m_A\left( i\xi+\frac{1}{2}\right) \right|=1
    \end{equation}
    для почти всех $\xi.$
    \end{enumerate}
\end{theorem}

Последняя теорема суммирует результаты многих математиков: Шура,
Харди, Литтлвуда, Пойа, Кобера, Михлина и Хермандера. Неизвестно,
можно ли в формулировке пункта б) при дополнительных условиях
опустить требование неотрицательности ядра. В  диапазоне $0<p<1$ в
общем случае  оценок нет, на что было  указано В.\,И.~Буренковым
на примере операторов Харди. Впоследствии мы увидим, что операторы
Харди тесно связаны с ОП Бушмана---Эрдейи. Первым математиком,
использовавшим технику преобразования Меллина для оценки норм
операторов Римана---Лиувилля для случая чисто мнимых степеней,
был, насколько нам известно, Кобер~\cite{Kob1}. Поэтому иногда
часть б) приведённой теоремы называется леммой Кобера, что не
совсем точно, так как он на самом деле доказал формулу для нормы
из части а) для случая знакопеременных функций.

Заметим, что преобразование Меллина является обобщённым
преобразованием Фурье на полуоси по мере Хаара
$\dfrac{dy}{y}$\cite{Hel1}. Его роль велика в теории специальных
функций, например, гамма-функция является преобразованием Меллина
экспоненты. С преобразованием Меллина связан важный прорыв в
1970-х годах, когда в основном усилиями О.\,И.~Маричева была
полностью доказана и приспособлена для нужд вычисления интегралов
известная теорема Джоан Люси Слейтер, позволяющая для большинства
образов преобразований Меллина восстановить оригинал в явном виде
по простому алгоритму через гипергеометрические
функции~\cite{Marich1, Sla, PBM}. Этот результат, который несложно
получить формально по общей формуле обращения Меллина---Барнса
через вычеты, для своего строгого обоснования потребовал
достаточно сложных и тщательных выкладок, связанных с обработкой
асимптотик гипергеометрических функций вблизи полюсов и на
бесконечности, а такие асимптотики весьма разнообразны и
разнородны. Эта работа была только начата Люси Джоан Слейтер, а в
основном проведена до конца Олегом Игоревичем Маричевым, данное
обстоятельство часто недооценивается. Теорема Слейтер---Маричева
позволила создать универсальный мощный метод вычисления
интегралов, который впоследствии позволил решить многие задачи в
теории уравнений с частными производными, а также воплотился в
передовые технологии символьного интегрирования пакета MATHEMATICA
(О.\,И.~Маричев работает в Wolfram Research).

Приведем теорему Слейтер (см.~\cite{Sla, Marich1, PBM}). Пусть
\begin{equation}\label{4.13}
\Gamma\left[\begin{array}{cccc}
a_1, & a_2, &\ldots, & a_A \\
b_1, & b_2, &\ldots & b_B
\end{array}
\right]=\Gamma[(a),(b)]=\frac{\Gamma(a_1)\Gamma(a_2)\ldots\Gamma(a_A)}{\Gamma(b_1)\Gamma(b_2)\ldots\Gamma(b_B)},
\end{equation}
пустое произведение заменяется единицей,
\begin{equation*}
(a)+s=a_1+s,a_2+s,\ldots,a_A+s,
\end{equation*}
\begin{equation}\label{4.14}
(b)'-b_k=b_1-b_k,\ldots,b_{k-1}-b_k,b_{k+1}-b_k,\ldots,b_B-b_k,
\end{equation}
\begin{equation}\label{4.15}
\arraycolsep=2.5pt \Sigma_A(z)=\sum\limits_{j=1}^A
z^{a_j}\Gamma\left[\begin{array}{cc}
(a)'-a_j, & (b)+a_j \\
(c)-a_j & (d)+a_j
\end{array}
\right] \,_{B+C}F_{A+D-1}\!\left(\begin{array}{ccc}
(b)+a_j, & 1+a_j-(c); & (-1)^{C-A}z \\
1+a_j-(a)', & (d)+a_j &
\end{array}
\right),
\end{equation}
\begin{equation}\label{4.16}
\arraycolsep=2pt \Sigma_B(1/z)=\sum\limits_{k=1}^B
z^{-b_k}\Gamma\left[\begin{array}{cc}
(b)'-b_k, & (a)+b_k \\
(d)-b_k & (c)+b_k
\end{array}
\right] \,_{A+D}F_{B+C-1}\!\left(\begin{array}{ccc}
(a)+b_k, & 1+a_k-(d); & \dfrac{(-1)^{D-B}}{z} \\
1+b_k-(b)', & (c)+b_k &
\end{array}
\right),
\end{equation}
$$
|\arg{z}|<\pi.
$$

Если ряды сходятся, то функции $\Sigma_A(z),$ $\Sigma_B(1/z)$
являются функциями гипергеометрического типа, причём переходят
друг в друга, если поменять местами $A$-мерный комплексный вектор
$(a)=a_1,a_2,\ldots,a_A$ с аналогичным $B$-мерным вектором $(b),$
$C$-мерный вектор $(c)$ с $D$-мерным $(d),$ а $z$ заменить на
$1/z.$ Эти функции аналитически зависят от комплексных параметров
$(a),$ $(b),$ $(c),$ $(d)$ и переменной $z.$ Если некоторые
параметры векторов $(a)$ (или $(b)$) совпадают между собой или
отличаются на целое число, то векторы $(a)'-a_j$  ($(b)'-b_k$)
содержат ненулевые или отрицательные целые компоненты и в силу
свойства $\Gamma(-n)=\infty,$ $n=0,1,2,\ldots,$ у функции
$\Sigma_A(z),$ $\Sigma_B(1/z),$ вообще говоря, могут возникнуть
неопределенности типа $\infty-\infty.$ В таких логарифмических
случаях под значениями $\Sigma_A(z),$ $\Sigma_B(1/z)$ будем
понимать соответствующие пределы <<регулярных>>, функций
$\Sigma_A(z),$ $\Sigma_B(1/z),$ когда их параметры непрерывно
стремятся к рассматриваемым особым значениям.

Следует отметить, что без ограничения $|\arg z|<\pi$ функции
$\Sigma_A(z),$ $\Sigma_B(1/z)$ в общем случае являются
многозначными.

\begin{theorem*}[Слейтер] Пусть функция $K^*(s)$ имеет вид
\begin{equation}\label{4.17}
K^*(s)=\Gamma\left[\begin{array}{cccc}
(a)+s & (b)-s \\
(c)+s, & (d)-s
\end{array}\right],
\end{equation}
где векторы $(a),$ $(b),$ $(c),$ $(d)$ имеют соответственно $A,$
$B,$ $C,$ $D$ компонент $a_j,$ $b_k,$ $c_l,$ $d_m.$ Тогда если
выполняются следующие две группы условий{\rm :}
\begin{equation}\label{4.181}
-\Re a_j<\Re s< \Re b_k,\,\, j=1,2,\ldots,A,\,\, k=1,2,\ldots,B,
\end{equation}
\begin{equation}\label{4.182}
\left\{
\begin{array}{ll}
$$ A+B>C+D,$$ &  \\
$$ A+B=C+D,$$ & \hbox{$\Re s(A+D-B-C)<-\Re \nu,$} \\
$$ A=C,\, B=D,$$ & \hbox{$\Re \nu<0,$}
\end{array}
\right.
\end{equation}
где
$$\nu=\sum\limits_{j=1}^Aa_j+\sum\limits_{k=1}^Bb_k-\sum\limits_{l=1}^Cc_l-\sum\limits_{m=1}^Dd_m,
$$ то для таких $s$ справедливы равенства
\begin{equation}\label{4.19}
K^*(s)=\left\{%
\begin{array}{ll}
$$\int\limits_0^\infty x^{s-1}\Sigma_A(x)dx,$$ & \hbox{$A+D>B+C;$} \\
$$\int\limits_0^1 x^{s-1} \Sigma_A(x)dx+\int\limits_1^\infty x^{s-1}\Sigma_B(1/x) dx,$$ & \hbox{$A+D=B+C;$} \\
$$\int\limits_0^\infty x^{s-1} \Sigma_B(1/x)dx,$$ & \hbox{$A+D<B+C,$} \\
\end{array}%
\right.
\end{equation}
$\Sigma_A(1)=\Sigma_B(1),$ если $A+D=B+C,$ $\Re \nu+C-A+1<0,$ $A\geq C.$
\end{theorem*}

\begin{corollary*}
При условиях~\eqref{4.181},~\eqref{4.182} прообразом функции
$$
K^*(s)=\Gamma\left[\begin{array}{cccc}
(a)+s & (b)-s \\
(c)+s, & (d)-s
\end{array}\right],
$$
является функция $K(x)$ гипергеометрического типа, равная одному
из следующих выражений{\rm :}
\begin{equation}\label{4.20}
K(x)=\left\{%
\begin{array}{ll}
$$\Sigma_A(x),$$ & \hbox{$x>0,A+D>B+C;$} \\
$$\Sigma_A(x),$$ & \hbox{$0<x<1,A+D=B+C;$} \\
$$\Sigma_B(1/x),$$ & \hbox{$x>1,A+D=B+C;$} \\
$$\Sigma_B(1/x),$$ & \hbox{$x>0,A+D<B+C,$} \\
\end{array}%
\right.
\end{equation}
$K(1)=\Sigma_A(1)-\Sigma_B(1),$ если $A+D=B+C,$ $\Re \nu+C-A+1<0,$ $A\geq C.$
\end{corollary*}

\pagebreak

\begin{remark}
\label{r1} Соответствующая теорема~1 из~\cite[\S~4.8]{Sla}
содержит ряд неточностей:
\begin{enumerate}
\item[1)] вместо условий
$$
\left\{
\begin{array}{ll}
$$A+B>C+D,$$ &  \\
$$A+B=C+D,$$ & \hbox{$\Re s(A+D-B-C)<-\Re \nu,$} \\
$$A=C,\, B=D,$$ & \hbox{$\Re \nu<0,$}
\end{array}
\right.
$$
указано лишь условие $A+B\geq C+D;$

\item[2)] вместо $\Re \nu+C-A+1<0$ приводится условие $\Re \nu<0;$

\item[3)] говорится, что при $A+D=B+C$ функции $\Sigma_A(z),$
$\Sigma_B(1/z)$ аналитически продолжают друг друга, а это верно
лишь в случае $A+B>C+D.$
\end{enumerate}
\end{remark}

\begin{remark}
\label{r2} В случае $|A+D-B-C|>1,$ $A+B=C+D$ ограничение на $\Re s,$ указанное в~\eqref{4.182}, может быть несколько ослаблено
до условия
\begin{equation}\label{4.183}
\Re s(A+D-B-C)<\frac{1}{2}-\Re \nu.
\end{equation}
\end{remark}

\begin{remark}
\label{r3} Если для некоторых параметров выполняется одно или
несколько условий $a_j=c_l+n$  (или $a_j=-d_m-n$) [$b_k=d_m+n$
(или $b_k=-c_l-n$)], где $n=0,1, 2,\ldots,$ и при этом векторы
$(a)'-a_j,$ $(b)'-b_k$ не содержат целые компоненты, то из группы
условий~\eqref{4.181} можно устранить (или ослабить) условия,
относящиеся к этим параметрам. При $a_j=c_l+n$  [$b_k=d_m+n$]
соответствующие условия $\Re (s+a_j)>0$  [$\Re (b_k-s)>0$]
устраняются, а при $(a)'-a_j,$ $(b)'-b_k$ они заменяются на
ослабленные требования $\Re (s+a_j)>-n-1$  [$\Re (b_k-s)>-n-1$].
Если же векторы $(a)'-a_j,$ $(b)'-b_k$ содержат целые компоненты,
то вопрос об ослаблении ограничений~\eqref{4.181} требует
специальных исследований.
\end{remark}

\begin{remark}
\label{r4} Ограничения~\eqref{4.181},~\eqref{4.182} обеспечивают,
по крайней мере, условную сходимость интегралов~\eqref{4.19} в
$0,$ $\infty$ (и $1$). Если эти ограничения нарушаются, то
интегралы~\eqref{4.19} как несобственные, вообще говоря,
расходятся, однако в отдельных случаях они существуют в смысле
главного значения.
\end{remark}

\subsubsection{Различные формы операторов дробного
интегродифференцирования}\label{sec3.3.3}

Операторы дробного интегродифференцирования изучены, например,
в~\cite{SKM, Nah1, Nah2, Nah3, KST, ViRy, GoMa, KT1,KT2}. Полезные
свойства дробных интегралов и их обращений для мер приведены
в~\cite{KaPr4}, см. также удачное учебное пособие~\cite{LShFrac}.

Операторы дробного интегродифференцирования играют важную роль во
многих современных разделах математики. Для теории специальных
функций важность дробного интегродифференцирования отражена в
названии известной статьи~\cite{Kir4}: <<\textit{Все специальные
функции получаются дробным интегродифференцированием элементарных
функций}!>> (Замечание проф. А.\,А.~Килбаса: <<{\it кроме функций
Фокса}!>>)

Рассмотрим дробные интегралы и производные Римана---Лиувилля

\begin{definition}
    Пусть $\varphi (x)\in L_{1} (a,b),$ тогда интегралы
    \begin{equation} \label{RLI1}
    (I_{a+}^{\alpha } \varphi )(x)=\frac{1}{\Gamma (\alpha )} \int\limits
    _{a}^{x}\frac{\varphi(t)}{(x-t)^{1-\alpha } } dt,\qquad x>a,
    \end{equation}
    \begin{equation} \label{RLI2}
    (I_{b-}^{\alpha } \varphi )(x)=\frac{1}{\Gamma (\alpha )} \int\limits
    _{x}^{b}\frac{\varphi(t)}{(t-x)^{1-\alpha } } dt,\qquad x<b,
    \end{equation}
    где $\alpha >0,$ называются соответственно {\it левосторонним}~\eqref{RLI1}  и {\it правосторонним}~\eqref{RLI2}
        {\it дробными
        интегралами Римана---Лиувилля} порядка $\alpha$ $(0<\alpha).$

    Для функции $f(x),$ $x\in[a,b]$ каждое из выражений
    \begin{equation} \label{RLD3} (D_{a+}^{\alpha }f)(x)=\frac{1}{\Gamma (n-\alpha )} \left(\frac{d}{dx} \right)^{n} \int\limits _{a}^{x}\frac{f(t)dt}{(x-t)^{\alpha -n+1} },   \end{equation}
    \begin{equation} \label{RLD4} (D_{b-}^{\alpha }f)(x)=\frac{1}{\Gamma (n-\alpha )} \left(\frac{d}{dx} \right)^{n} \int\limits _{x}^{b}\frac{f(t)dt}{(t-x)^{\alpha -n+1} },  \end{equation}
    где $n=[\alpha ]+1,$ ${\alpha >0},$ называется {\it дробной производной
    Римана---Лиувилля} порядка $\alpha,$ соответственно, {\it левосторонней} и
    {\it правосторонней}.
\end{definition}

В частных случаях наиболее важные в приложениях операторы
Римана---Лиувилля определяются при $\alpha > 0$ по формулам:
\begin{eqnarray}
\label{161}
I_{0+,x}^{\alpha}f=\frac{1}{\Gamma(\alpha)}\int\limits_0^x \left( x-t\right)^{\alpha-1}f(t)d\,t,\\
\nonumber
I_{-,x}^{\alpha}f=\frac{1}{\Gamma(\alpha)}\int\limits_x^\infty
\left( t-x\right)^{\alpha-1}f(t)d\,t.
\end{eqnarray}
Для остальных значений  $\alpha$ они определяются при помощи аналитического продолжения (регуляризации).

Отметим, что существуют многочисленные варианты операторов
дробного интегродифференцирования: Герасимова---Капуто, введённые
советским механиком Герасимовым в 1948~г.~\cite{Nov} и Капуто в
1968~г., Маршо, Вейля, Рисса, Эрдейи---Кобера, Адамара,
Гельфонда---Леонтьева (двойное название предложено
Ю.\,Ф.~Коробейником), Джрбашяна---Нерсесяна и т.~д.
Дифференциальные уравнения типа Эйлера с основными типами
операторов дробного интегродифференцирования рассмотрены,
например, в~\cite{KiZhu,ZhuSi1}.

Теперь рассмотрим дробные интегралы и производные Эрдейи---Кобера.

\textit{Операторы Эрдейи---Кобера} определяются при $\alpha > 0$
по формулам:
\begin{eqnarray}
\label{162}
& & I_{0+;\, 2,\, y}^{\alpha} f = \frac{2}{\Gamma(\alpha)}x^{-2(\alpha+y)}
\int\limits_0^x (x^2-t^2)^{\alpha-1}t^{2y+1}f(t)\,dt, \\
& & I_{-;\, 2,\, y}^{\alpha} f = \frac{2}{\Gamma(\alpha)}x^{2y}
\int\limits_x^{\infty}
(t^2-x^2)^{\alpha-1}t^{2(1-\alpha-y)-1}f(t)\,dt,
\end{eqnarray}
а при значениях $\alpha > -n,$ $n \in \mathbb{N}$ по формулам
\begin{eqnarray}
& & I_{0+;\, 2, y}^{\alpha} f =x^{-2(\alpha+y)} {\lr{\frac{d}{d x^2}}}^n  x^{2(\alpha
    +y+n)}I^{\alpha+n}_{0+; \, 2,\, y}f, \label{1.15} \\
& & I_{-;\, 2, y}^{\alpha} f = x^{2y} {\lr{-\frac{d}{d x^2}}}^n
x^{2(\alpha
    -y)}I^{\alpha+n}_{-; \, 2,\, y-n}f \label{1.16}.
\end{eqnarray}
Для остальных значений  $\alpha$ они определяются при помощи аналитического продолжения, аналогично операторам дробного интегродифференцирования Лиувилля.

Отметим, что в классической монографии~\cite{SKM} случаи выбранных
нами пределов интегрирования $0,\infty$ не рассматриваются. В
последующей английской версии~\cite{SaKiMar} эти особые случаи
пределов допускаются, но определения содержат неточности, в
частности, приводящие к комплексным величинам под знаком
интеграла.

Теперь рассмотрим наиболее общую форму подобных операторов "---
дробный интеграл по произвольной функции.

\textit{Дробный интеграл по произвольной функции $g(x)$}:
\begin{eqnarray}
\label{163}
I_{0+,g}^{\alpha}f=\frac{1}{\Gamma(\alpha)}\int\limits_0^x \left( g(x)-g(t)\right)^{\alpha-1}g'(t)f(t)d\,t,\\
\nonumber
I_{-,g}^{\alpha}f=\frac{1}{\Gamma(\alpha)}\int\limits_x^\infty
\left( g(t)-g(x)\right)^{\alpha-1}g'(t)f(t)d\,t,
\end{eqnarray}
во всех случаях предполагается, что $\Re\alpha>0,$ на оставшиеся
значения $\alpha$ формулы также без труда продолжаются~\cite{SKM}.
При этом обычные дробные интегралы~\eqref{161} получаются при
выборе в~\eqref{163} $g(x)=x,$ Эрдейи---Кобера~\eqref{162} при
выборе $g(x)=x^2,$ Адамара при $g(x)=\ln x.$

Связь с ОП проявляется в том, что операторы преобразования
Сонина---Пуассона---Дельсарта с точностью до множителей как раз и
являются операторами Эрдейи---Кобера, то есть дробными степенями
$$\left(\frac{d}{dg(x)}\right)^{-\alpha}=\left(\frac{d}{2xdx}\right)^{-\alpha},\quad g(x)=x^2.$$
Поэтому основные свойства этих ОП можно получить из теории
операторов дробного интегродифференцирования, а не изобретать
заново, что нередко и делалось. А.\,М.~Джрбашян обратил наше
внимание на тот факт, что операторы дробного интегрирования по
функции~\eqref{163} являются частными случаями несколько более
общих операторов, которые были введены и изучались его отцом
М.\,М.~Джрбашяном~\cite{SKM}.

\subsubsection{Квадратичные {\rm (}дробные{\rm )} преобразования Фурье и
Ханкеля}\label{sec3.3.4}

Кратко изложим основные факты, относящиеся к теории квадратичного
или дробного преобразования Фурье---Френеля, следуя~\cite{OZK}
и~\cite{Kar1}, где можно найти более подробную информацию и
ссылки.

Целые степени  (орбита) классического преобразования Фурье
образуют циклическую группу порядка 4, при этом четвёртая степень
этого преобразования даёт тождественный оператор. Поэтому, в
частности,  спектр классического преобразования Фурье в
$L_2(-\infty,\infty)$ состоит из четырех точек, расположенных на
единичной окружности:  $1,$ $i,$ $-1$ и $-i.$ Идея включить эту
дискретную группу в непрерывную со спектром, целиком заполняющим
единичную окружность, принадлежит Винеру  и была реализована им в
работе 1929~г.  Несколько позже Кондон  (в 1937~г.),  а затем
Кобер (в 1939~г.) независимо переоткрыли эту группу, которая стала
называться дробным преобразованием Фурье (ДПФ). Баргманну
принадлежит обобщение на многомерный случай.

Впоследствии дробное преобразование Фурье неоднократно
переоткрывалось целым рядом авторов.  В книге Антосика,
Микусинского и Сикорского  под названием <<преобразование
Фурье---Мелера>> упоминается циклическая группа произвольного
порядка, в которую можно включить преобразование Фурье. ДПФ
изучалось Гинандом  и Вольфом. Вавржинчик в  приходит к ДПФ
рассматривая классическое преобразование Фурье в виде экспоненты
от производящего оператора. В.\,Ф.~Осиповым  независимо от
предыдущих авторов построена теория ДПФ на группах и введены
соответствующие почти-периодические функции Бора---Френеля,
изучены асимптотические свойства этого преобразования, рассмотрены
приложения в гармоническом анализе и теории чисел~\cite{Os, AbOs}.
Намиас  переоткрыл дробное преобразования Фурье и использовал его
для решения некоторых задач для уравнения Шр\"{e}\-дин\-гера. Керр
исследовала дробное преобразование Фурье в пространстве $L_2$ и
пространстве Шварца $S.$

Отдельное направление исследований связано с дробным
преобразованием Фурье при чисто мнимых значениях группового
параметра. В этом случае чаще всего встречается название
<<полугруппа Эрмита>>. Начало этому направлению положила статья
Хилле 1926~г., в которой оператор ДПФ при мнимых значениях
параметра возникает в связи с суммированием методом Абеля
разложений по полиномам Эрмита. Позднее, <<полугруппа Эрмита>>
была использована Бабенко, Бекнером и Вайслером
 для получения неравенств в теории классического преобразования Фурье.

Отметим, что  квадратичное преобразование Фурье является одной из
двух основных составляющих (наряду с неравенствами для средних
значений в комплексной плоскости специального вида), которые были
использованы сначала К.\,И.~Бабенко для частных случаев, а затем
Бекнером для общего случая при доказательстве знаменитых условий
ограниченности обычного   преобразования Фурье в пространствах
$L_p$ с точными постоянными. Другим интересным применением
квадратичного преобразования Фурье является круг вопросов,
связанных со знаменитой задачей Паули по определению функции по
некоторому набору спектральных данных. Подобные задачи, обычно
неразрешимые для классического преобразования, находят элегантные
решения в рамках КПФФ.

Применения дробного преобразования Фурье столь же разнообразны и
обширны как и классического. Перечислим лишь некоторые.
Приложениям в квантовой механике посвящена уже упоминавшаяся
работа Намиаса. Использование ДПФ в оптике и анализе сигналов
разрабатывается  группой исследователей под руководством Оцактаса.
Участниками этой группы написана книга~\cite{OZK}, целиком
по\-свя\-щ\"{e}н\-ная теории и приложениям дробного преобразования
Фурье,  в которой процитировано несколько сотен работ по данной
теме (и в том числе выражена благодарность одному из авторов этой
книги за консультации и обсуждения в период подготовки монографии
к печати). Мастардом найдены аналоги неравенства Гейзенберга,
инвариантные относительно дробного преобразования Фурье.  Было
показано, что многомерное преобразование Вигнера равно корню
шестой степени из обратного преобразования Фурье, и,
следовательно, также является частным случаем дробного
преобразования Фурье. Бьюн  получил некоторые новые формулы
обращения для дробного преобразования Фурье.  Аналог дробного
преобразования Фурье для $q$-полиномов Эрмита был введён Аски,
Н.\,М.~Атакишиевым и С.\,К.~Сусловым.  Дальнейшие обобщения на
операторы с ядрами в форме билинейных производящих функций
полиномов Аски---Вилсона рассматривались в работах Аски и Рахмана.

Дробное преобразование Ханкеля (ДПХ) гораздо менее изучено.  Оно было введено Кобером в
 и изучалось Гинандом.  Затем было несколько раз переоткрыто, например,  Намиасом.  ДПХ было рассмотрено
Керр при действительных значениях группового параметра в
пространстве $L_2(0,\infty)$  и в пространствах Земаняна, см.
также~\cite{S60}.

В диссертации Д.\,Б.~Карпа~\cite{Kar1} рассмотрен достаточно общий
подход к построению подобных преобразований при помощи разложения
в ряды по известным системам  ортогональных функций. В частности,
при выборе системы функций Эрмита получается классическое
преобразование Фурье, при выборе системы функций Лагерра "---
квадратичное преобразование Фурье---Френеля, а при выборе систем
функций Лежандра, Чебышёва или Гегенбауэра построены новые
полугруппы интегральных преобразований.

Выпишем явно интегральные формулы для квадратичных преобразований
Фурье---Френеля (КПФФ) и Ханкеля (КПХ), следуя~\cite{Kar1}:
\begin{equation}\label{FrF}
(F^{\alpha}f)(y)=\frac{1}{\sqrt{\pi(1-e^{2i\alpha})}}\int\limits_{-\infty}^{\infty}\!\!
e^{-\frac{1}{2}i(x^2+y^2)\ctg\alpha}e^{ixy\cosec\alpha}f(x)dx,
\end{equation}

\begin{equation}\label{FrH}
(H_\nu^{\alpha}f)(y)=\frac{2\left(-e^{i\alpha}\right)^{-\frac{\nu}{2}}}{1-e^{i\alpha}}
\int\limits_0^{\infty}\!\!e^{-\frac{1}{2}i\ctg\frac{\alpha}{2}\left(x^{2}+y^{2}\right)}
(xy)^{\frac{1}{2}}J_{\nu}\left(\frac{2xy\sqrt{-e^{i\alpha}}}{1-e^{i\alpha}}\right)f(x)dx.
\end{equation}

Рассмотрим кратко формулы преобразований, которые образуют операционное исчисление для КПХ.

Введ\"{e}м следующие дифференциальные операторы:
\begin{gather*}
A^{-}_\nu=x^{\nu +\frac{1}{2}}e^{-\frac{x^2}{2}}\frac{d}{dx}
x^{-\nu -\frac{1}{2}}e^{\frac{x^2}{2}}=-\frac{\nu +
\frac{1}{2}}{x}+x+\frac{d}{dx},
\\
A^{+}_\nu=x^{-\nu -\frac{1}{2}}e^{\frac{x^2}{2}}
\frac{d}{dx}x^{\nu
+\frac{1}{2}}e^{-\frac{x^2}{2}}=\frac{\nu+\frac{1}{2}}{x}-x+\frac{d}{dx},
\\
N_\nu=x^{\nu+\frac{1}{2}}\frac{d}{dx}x^{-\nu-\frac{1}{2}}=-\frac{\nu+\frac{1}{2}}{x}+\frac{d}{dx},
\\
M_\nu=x^{-\nu-\frac{1}{2}}\frac{d}{dx}x^{\nu+\frac{1}{2}}=\frac{\nu+\frac{1}{2}}{x}+\frac{d}{dx}.
\end{gather*}
Эти операторы связаны соотношениями
\begin{eqnarray}\label{AMN}
L_{\nu}= -\frac{1}{4}D^2-\frac{\nu^2-1/4}{x^2} + \frac{1}{4}x^2 - \frac{\nu+1}{2}=-\frac{1}{4}A^{+}_{\nu}A^{-}_{\nu},\nonumber\\~~~~A^{-}_{\nu}=N_{\nu}+x,
~~~~A^{+}_{\nu}=M_{\nu}-x,
\end{eqnarray}
\begin{equation}\label{LMN}
L_{\nu}=M_{\nu} N_{\nu},~~~~
x\frac{d}{dx}+\frac{d}{dx}x=N_{\nu}x+xM_{\nu}=M_{\nu}x+xN_{\nu}.
\end{equation}

Обозначим через $X$ оператор умножения на независимую переменную.
Теперь мы можем, следуя~\cite{Kar1},  выписать набор формул
преобразования операций для КПХ.

\begin{align}\label{H+tD}
H_{\nu+1}^{\alpha}Xf&=\frac{1}{2}\left[(e^{-i\alpha}-1)N_{\nu}+
(e^{-i\alpha}+1)X\right]H_{\nu}^{\alpha}f,
\\
\label{H+ND}
H_{\nu+1}^{\alpha}N_{\nu}f&=\frac{1}{2}\left[(e^{-i\alpha}+1)N_{\nu}+
(e^{-i\alpha}-1)X\right]H_{\nu}^{\alpha}f,
\\
\label{NHD} N_{\nu}H_{\nu}^{\alpha}f&=
\frac{1}{2}H_{\nu+1}^{\alpha}\left[(e^{i\alpha}+1)N_{\nu}+(e^{i\alpha}-1)X\right]f,
\\
\label{xHD}
XH_{\nu}^{\alpha}f&=\frac{1}{2}H_{\nu+1}^{\alpha}\left[(e^{i\alpha}-1)N_{\nu}+
(e^{i\alpha}+1)X\right]f,
\\
\label{HtD} H_{\nu}^{\alpha}Xf&=
\frac{1}{2}\left[(1-e^{i\alpha})M_{\nu}+(1+e^{i\alpha})X\right]H_{\nu+1}^{\alpha}f,
\\
\label{HMD} H_{\nu}^{\alpha}M_{\nu}f&=
\frac{1}{2}\left[(1+e^{i\alpha})M_{\nu}+(1-e^{i\alpha})X\right]H_{\nu+1}^{\alpha}f,
\\
\label{MH+D} M_{\nu}H_{\nu+1}^{\alpha}f&=
\frac{1}{2}H_{\nu}^{\alpha}\left[(1+e^{-i\alpha})M_{\nu}+(1-e^{-i\alpha})X\right]f,
\\
\label{xH+D} XH_{\nu+1}^{\alpha}f&=
\frac{1}{2}H_{\nu}^{\alpha}\left[(1-e^{-i\alpha})M_{\nu}+(1+e^{-i\alpha})X\right]f,
\\
\label{A-H}
H_{\nu+1}^{\alpha}A^{-}_{\nu}f&=e^{-i\alpha}A^{-}_{\nu}H_{\nu}^{\alpha}f,
\qquad
A^{-}_{\nu}H_{\nu}^{\alpha}f=H_{\nu+1}^{\alpha}e^{i\alpha}A^{-}_{\nu}f,
\\
\label{A+H}
H_{\nu}^{\alpha}A^{+}_{\nu}f&=e^{i\alpha}A^{+}_{\nu}H_{\nu+1}^{\alpha}f,\qquad
A^{+}_{\nu}H_{\nu+1}^{\alpha}f=H_{\nu}^{\alpha}e^{-i\alpha}A^{+}_{\nu}f,
\\
\label{Ht^2D} H_{\nu}^{\alpha}X^{2}f&=
\left[X^{2}\cos^{2}{\frac{\alpha}{2}}-\frac{1}{2}i\sin{\alpha}\left[XD+DX\right]-
\sin^{2}{\frac{\alpha}{2}}L_{\nu}\right]H_{\nu}^{\alpha}f,
\\
\label{x^2HD} X^{2}H_{\nu}^{\alpha}f&=
H_{\nu}^{\alpha}\left[X^{2}\cos^{2}{\frac{\alpha}{2}}+\frac{1}{2}i\sin{\alpha}\left[XD+DX\right]-
\sin^{2}{\frac{\alpha}{2}}L_{\nu}\right]f,
\\
\label{HLD} H_{\nu}^{\alpha}L_{\nu}f&=
\left[-X^{2}\sin^{2}{\frac{\alpha}{2}}-\frac{1}{2}i\sin{\alpha}\left[XD+DX\right]+
\cos^{2}{\frac{\alpha}{2}}L_{\nu}\right]H_{\nu}^{\alpha}f,
\\
\label{LHD} L_{\nu}H_{\nu}^{\alpha}f&=
H_{\nu}^{\alpha}\left[-X^{2}\sin^{2}{\frac{\alpha}{2}}+\frac{1}{2}i\sin{\alpha}\left[XD+DX\right]+
\cos^{2}{\frac{\alpha}{2}}L_{\nu}\right]f,
\\
\label{HtD+DtD} H^\alpha_\nu\left[XD+DX\right]f&=
\left[-i\sin{\alpha}(X^2+L_\nu)+\cos\alpha\left[XD+DX\right]\right]H_{\nu}^{\alpha}f,
\\
\label{xD+DxHD}
\left[XD+DX\right]H_{\nu}^{\alpha}f&=H_{\nu}^{\alpha}\left[i\sin{\alpha}(X^2+L_\nu)+
\cos\alpha\left[XD+DX\right]\right]f.
\end{align}

Приведённые выше формулы операционного исчисления для КПХ
позволяют использовать это интегральное преобразование в
композиционном методе построения операторов преобразования,
излагаемом в главе~\ref{ch6}.

Итак, мы рассмотрели основные интегральные преобразования, которые
используются в данной работе. По поводу общих вопросов теории для
различных классов  операторов и функциональных пространств
см.~\cite{Nai, SvFe, KF, KPS, Bas, Kus, Pas1, Sol, Fet}.

\subsubsection{Основные классы дифференциальных уравнений с операторами Бесселя и связанных с ними операторов
преобразования}\label{sec3.3.5}

В книге принята следующая терминология для названий дифференциальных уравнений с операторами Бесселя
\begin{equation}
\label{c1Bessel}
B_{\nu}u(x)=\frac{d^2 u}{dx^2}+\frac{2\nu+1}{x} \frac{du}{dx},
\end{equation}
в основном введённая И.\,А.~Киприяновым~\cite{Kip1}.

{\it $B$-эллиптическим} называется уравнение с операторами Бесселя
вида
\begin{equation}
\label{c1Bes1} \sum\limits_{k=1}^{n}B_{\nu,x_k}u(x_1,\dots,
x_n)=f(t,x),
\end{equation}
иногда такое уравнение называется {\it уравнением
Лапласа---Бесселя}.

{\it $B$-гиперболическим} называется уравнение с операторами
Бесселя вида
\begin{equation}
\label{c1Bes2} B_{\nu,t}u(t,x_1,\dots,
x_n)=\sum\limits_{k=1}^{n}B_{\nu,x_k}u(t,x_1,\dots, x_n) + f(t,x),
\end{equation}
в случае одной пространственной переменной получаем {\it уравнение
Эйлера---Пуассона---Дарбу}.

{\it $B$-параболическим} называется уравнение с операторами
Бесселя вида
\begin{equation}
\label{c1Bes3} \frac{\pr}{\pr t}u(t,x_1,\dots,
x_n)=\sum\limits_{k=1}^{n}B_{\nu,x_k}u(t,x_1,\dots, x_n) + f(t,x).
\end{equation}

Те же названиях сохраняем для неполных уравнений, в которых один
или несколько операторов Бесселя сводятся ко вторым производным, а
также к уравнениям добавляются спектральные параметры. Указанные
три класса дифференциальных уравнений по терминологии
И.\,А.~Киприянова в своё время были рассмотрены в трёх известных
монографиях: $B$-эллиптические уравнения в монографии
И.\,А.~Киприянова~\cite{Kip1}, $B$-гиперболические уравнения в
монографии Р.~Кэрролла и Р.~Шоуолтера~\cite{CSh},
$B$-параболические уравнения в монографии
М.\,И.~Матийчука~\cite{Mat1}. Перечислим некоторые известные
классы дифференциальных уравнений с операторами Бесселя, не
претендуя на полноту этого списка, также опустив ссылки на работы,
кроме нескольких основных. Более полная информация приведена в
монографии~\cite{SSfiz}, которая как раз и посвящена рассмотрению
дифференциальных уравнений с операторами Бесселя.

$B$-эллиптическое уравнение вида
\begin{equation}\label{BEll}
\sum\limits_{i=1}^n\left( \frac{\partial^2 u}{\partial x_i^2}+\frac{k_i}{x_i}\frac{\partial u}{\partial x_i}\right)=f(x).
\end{equation}
К этому типу относится  уравнение обобщённой теории
осесимметрического потенциала А.~Вайнштейна (GASPT "--- {\it
Generalized Axially Symmetric Potential Theory},~\cite{Wei1, Wei2,
Wei3}), уравнения этого класса были достаточно полно изучены в
работах И.\,А.~Киприянова и его школы~\cite{Kip1}.

В случае, когда оператор Бесселя действует по пространственной переменной мы получаем сингулярный вариант волнового уравнения с осевой или центральной симметрией
\begin{equation}\label{EPDWe}
\frac{\partial^2 u}{\partial t^2}=
\frac{\partial^2 u}{\partial x^2}+\frac{\nu}{x}\frac{\partial u}{\partial x},\qquad u=u(x,t),\qquad x>0,\qquad t\in\mathbb{R},\qquad
\nu=\const.
\end{equation}
Представления решений уравнения~\eqref{EPDWe} было получено
Пуассоном~\cite{Poisson}. К этому классу относится обобщённая
радиальная задача излучения Вайнштейна~\cite{CSh}.

Для случая, когда оператор Бесселя действует по временной
переменной $t,$ мы получаем знаменитое уравнение
Эйлера---Пуассона---Дарбу (ЭПД)
\begin{equation}\label{EPD0}
\frac{\partial^2 u}{\partial t^2}+\frac{\nu}{t}\frac{\partial u}{\partial t}=
a^2\frac{\partial^2 u}{\partial x^2},\qquad u=u(x,t),\qquad t>0,\qquad x\in\mathbb{R},\qquad
a, \ \nu=\const.
\end{equation}
Это уравнение впервые появилось в работе
Л.~Эйлера~\cite[с.~227]{Euler}, и затем изучалось
С.~Пуассоном~\cite{Poisson}, Б.~Риманом~\cite{Riman} и
Г.~Дарбу~\cite{Darboux}.

В случае нескольких пространственных переменных при условии
$\nu>n-1$ Диаз и Вайнбергер~\cite{Diaz} рассмотрели задачу Коши
\begin{equation}\label{EPD1}
\frac{\partial^2 u}{\partial t^2}+\frac{\nu}{t}\frac{\partial u}{\partial t}=
\sum\limits_{i=1}^n\frac{\partial^2 u}{\partial x_i^2},\qquad u=u(x,t),\qquad t>0,\qquad x\in\mathbb{R}^n,\qquad
\nu=\const,
\end{equation}
\begin{equation}\label{UslEPD1}
u(x,0)=f(x),\qquad u_t(x,0)=0.
\end{equation}
Для произвольных $\nu$ задача~\eqref{EPD1}-\eqref{UslEPD1} была
решена А.~Вайнштейном, см.~\cite{Wei1, Wei2, Wei3}.

Обобщённое уравнение ЭПД
\begin{equation}\label{EPDWG}
\frac{\partial^2 u}{\partial t^2}+\frac{\nu}{t}\frac{\partial u}{\partial t}=
\frac{\partial^2 u}{\partial x^2}+\frac{k}{x}\frac{\partial u}{\partial x},\qquad u=u(x,t),\qquad t>0,\qquad x>0,\qquad
\nu,k=\const,
\end{equation}
и его многомерный вариант
\begin{equation}\label{EPDM}
\frac{\partial^2 u}{\partial t^2}+\frac{\nu}{t}\frac{\partial u}{\partial t}=
\sum\limits_{i=1}^n\left( \frac{\partial^2 u}{\partial x_i^2}+\frac{k_i}{x_i}\frac{\partial u}{\partial x_i}\right),
\end{equation}
$$
 u=u(x_1,..,x_n,t),\qquad t>0,\qquad x_i>0,\qquad
 \nu,k_i=\const, \qquad i=1,..,n,
$$
были рассмотрены в~\cite{CSh,Fox,LPSh1,LPSh2,ShiE3,Smi,76,Ter1}.
Уравнение ЭПД со спектральным параметром
\begin{equation}\label{EPDSP}
\frac{\partial^2 u}{\partial t^2}+\frac{\nu}{t}\frac{\partial u}{\partial t}=
\frac{\partial^2 u}{\partial x^2}\pm\lambda^2 u,\qquad \lambda\in\mathbb{R}
\end{equation}
было рассмотрено в~\cite{Bresters2,Smi}, а обобщённое уравнение
ЭПД со спектральным параметром
\begin{equation}\label{EPDSP1}
\frac{\partial^2 u}{\partial t^2}+\frac{\nu}{t}\frac{\partial u}{\partial t}=
\sum\limits_{i=1}^n\left( \frac{\partial^2 u}{\partial x_i^2}+\frac{k_i}{x_i}\frac{\partial u}{\partial x_i}\right) -\lambda^2 u,\qquad \lambda\in\mathbb{R}
\end{equation}
было изучено в~\cite{ShiE5}.

Обобщённые волновые уравнения с переменными потенциалами вида
\begin{equation}\label{wpot}
\frac{\partial^2 u}{\partial t^2} = \frac{\partial^2 u}{\partial x^2}+p(x)u,\ \ \
\frac{\partial^2 u}{\partial t^2} +q(t)u= \frac{\partial^2 u}{\partial x^2}+p(x)u,
\end{equation}
далее обобщаются на $B$-гиперболические уравнения с переменными
потенциалами
\begin{equation}\label{Bpot1}
B_{\nu,t}u = B_{\nu,x}u + p(x)u,\ \ \
B_{\nu,t}u +q(t)u= B_{\nu,x}u + p(x)u,
\end{equation}
их многомерные аналоги
\begin{equation}\label{Bpot2}
B_{\nu,t}u +q(t)u= \sum\limits_{k=1}^n B_{\nu,x_k}u + p(x_k)u
\end{equation}
и $B$-ультрагиперболические уравнения с переменными потенциалами
\begin{equation}\label{Bpot3}
\sum\limits_{j=1}^m \left(B_{\nu,t_j}u +q(t_j)u\right)=
\sum\limits_{k=1}^n \left(B_{\nu,x_k}u + p(x_k)u\right).
\end{equation}
Эти уравнения, явно или неявно, обычно решаются с использованием
ОП для операторов Штурма---Лиувилля или возмущённых операторов
Бесселя.

Другим важным классом, который начал разрабатываться только в
последнее время, являются интегро-дифференциальные уравнения
дробного порядка с дробными степенями операторов Бесселя. Основы
этой теории были заложены в~\cite{McB,Spr,Dim,Kir1, S140p, S135,
S133} и развивались в~\cite{S127, S123, 18, S700, SS, FJSS}.

В цитированной литературе приведены также многочисленные приложения дифференциальных уравнений с операторами Бесселя. Таким образом, этот класс уравнений имеет важное значение как для теории дифференциальных уравнений в частных производных, так и для практических приложений.

Отметим, что мы ограничиваемся  в этой книге рассмотрением
линейных дифференциальных уравнений, нелинейные уравнения требуют
других подходов и являются отдельным направлением исследований,
см., например, обзорные монографии~\cite{Kap1,ZaiPol,Rad3,Rad4}.

Приведём полезную непосредственно проверяемую формулу
Дарбу---Вайнштейна (это название дано Ж.-Л.~Лионсом)
\begin{equation}
B_{\nu} \lr{ y^{-2 \nu} f(y)} = y^{-2 \nu} B_{\nu} f(y),
\label{9}
\end{equation}
с помощью которой случай $\Re \nu < 0$ элементарно сводится к
случаю $\Re \nu > 0.$

Теперь определим самый известный класс ОП, сплетающих дифференциальный оператор Бесселя со второй производной:
\begin{equation}
\label{151} T\lr{B_\nu} f=\lr{D^2} Tf,\quad B_{\nu}=D^2+\frac{2\nu
+1}{x}D,\quad D^2=\frac{d^2}{dx^2},\quad \nu \in \mathbb{C}.
\end{equation}

 Одним из способов построения ОП является установление соответствий между решениями  дифференциальных уравнений. Решениями уравнения вида $B_\nu f=\lambda f$ являются функции Бесселя, а уравнения $D^2f=\lambda f$ "--- тригонометрические функции или экспонента. Поэтому прообразами ОП вида~\eqref{151} были формулы Пуассона и Сонина:
\begin{eqnarray}
\label{152}
J_{\nu} (x)=\frac{1}{\sqrt{\pi}\Gamma({\nu+\frac{1}{2}})2^{\nu-1}x^\nu}
\int\limits_0^x \left( x^2-t^2\right)^{\nu-\frac{1}{2}}\cos(t)\,dt,\quad \Re \nu> \frac{1}{2},\\
\label{153} J_{\nu}
(x)=\frac{2^{\nu+1}x^{\nu}}{\sqrt{\pi}\Gamma({\frac{1}{2}-\nu})}
\int\limits_x^\infty \left(
t^2-x^2\right)^{-\nu-\frac{1}{2}}\sin(t)\,dt,\quad -\frac{1}{2}<
\Re \nu<\frac{1}{2}.
\end{eqnarray}

Интеграл~\eqref{152}  начал изучать Эйлер в 1769~г. Затем
Парсеваль посчитал интеграл при $\nu=0$ в 1805~г.,  для целых
$\nu$ формулу~\eqref{152} получил Плана в 1821~г., Пуассон вывел
её для полуцелых $\nu$ в 1823~г., его метод применим и для целых
$\nu,$ но он этого не заметил.  Далее этот интеграл встречался в
работах Куммера, Лоббато и Дюамеля. Окончательно
формулу~\eqref{152}, которую мы приписываем Пуассону, установил в
общем случае Ломмель в 1868~г., а Сонин вывел формулу~\eqref{153}
в 1880~г.

\begin{definition}
 {\it ОП Пуассона} называется выражение
\begin{equation}
\label{154} P_{\nu}f=\frac{1}{\Gamma(\nu+1)2^{\nu}x^{2\nu}}
\int\limits_0^x \left(
x^2-t^2\right)^{\nu-\frac{1}{2}}f(t)\,dt,\quad \Re \nu>
-\frac{1}{2}.
\end{equation}
{\it ОП Сонина} называется выражение
\begin{equation}
\label{155}
S_{\nu}f=\frac{2^{\nu+\frac{1}{2}}}{\Gamma(\frac{1}{2}-\nu)}\frac{d}{dx}
\int\limits_0^x \left(
x^2-t^2\right)^{-\nu-\frac{1}{2}}t^{2\nu+1}f(t)\,dt,\quad \Re \nu<
\frac{1}{2}.
\end{equation}

Операторы~\eqref{154}--\eqref{155} действуют как ОП по формулам
\begin{equation}
\label{156} S_\nu B_\nu=D^2 S_\nu,\quad  P_\nu D^2=B_\nu P_\nu.
\end{equation}
\end{definition}
Их можно доопределить на все значения $\nu\in\Cbb.$

Идею изучения операторов подобных~\eqref{154}--\eqref{155}
высказывал ещё Лиувилль, их реальное использование в контексте
теории функций Бесселя начал Н.\,Я.~Сонин. Как ОП эти операторы
впервые были введены в работах Жана Дельсарта, а затем на основе
идей Дельсарта их изучение продолжилось в работах Дельсарта и в
совместных работах Дельсарта и Лионса. Поэтому мы будем
называть~\eqref{154}--\eqref{155} ОП Сонина---Пуассона---Дельсарта
(СПД). Об операторах СПД см. также статью
Б.\,М.~Левитана~\cite{Lev7}.

Не будет преувеличением сказать, что  операторы
СПД~\eqref{154}--\eqref{155} являются самыми знаменитыми объектами
всей теории ОП, их изучению, приложениям и обобщениям посвящены
сотни работ.

Кроме операторов Сонина и Пуассона  нам потребуются и другие
подобные, которые получаются заменой в  операторов
Римана---Лиувилля $I^{\mu}$ на операторы $I^{\mu}_e,$ определяемые
по формуле
\begin{equation}
I^{\mu}_e =  \mathcal{E} I^{\mu}  \mathcal{E}^{-1},
\label{10}
\end{equation}
где $\mathcal{E}$ "--- оператор умножения на функцию $e^x.$ Эти
операторы введены в главе~\ref{ch2}, где установлена, в частности,
формула
\begin{equation}
\mathcal{J}_{\mu, e} f(y) \equiv  I^{\mu}_e I^{-\mu} f(y) = f(y) - \mu \int\limits_y^{\infty} f(t) \, \Phi (\mu+1, 2; y-t) \, dt,
\label{11}
\end{equation}
в которой $\Phi (a, c; z)$ "--- вырожденная гипергеометрическая
функция (относительно других свойств операторов $I^{\mu}_e$ см.
книгу~\cite{SKM}). Тогда справедливы формулы
\begin{equation}
P_{\nu, e} \equiv P_{\nu}^{\frac{1}{2}-\nu} I_e^{\nu-\frac{1}{2}}=P_{\nu} \mathcal{J}_{\nu-\frac{1}{2}, e},
\label{12}
\end{equation}
\begin{equation}
S_{\nu, e} \equiv I_e^{\frac{1}{2}-\nu} S_{\nu}^{\nu-\frac{1}{2}}= \mathcal{J}_{\frac{1}{2}-\nu, e} S_{\nu}.
\label{13}
\end{equation}

В главе~\ref{ch2} также изучаются связи операторов типа
Римана---Лиувилля с преобразованиями Фурье и Ханкеля, имеющими
следующий вид:
\begin{equation}
F f \lr{\eta} = \int\limits_{-\infty}^{\infty} f(y) e^{- i y \eta} dy, \  F^{-1} g(y) = \frac{1}{2 \pi}
\int\limits_{-\infty}^{\infty} g(\eta) e^{ i y \eta} d \eta,
\label{14}
\end{equation}
\begin{equation}
F_{-} f \lr{\eta} = \int\limits_{0}^{\infty} f(y) \sin \lr{y \eta} dy, \  F^{-1}_{-}  g(y) = \frac{2}{\pi}
\int\limits_{0}^{\infty} g(\eta) \sin \lr{y \eta} d \eta,
\label{15}
\end{equation}
\begin{equation}
F_{\nu} f \lr{\eta} = \int\limits_{0}^{\infty} f(y) j_{\nu} \lr{y \eta} y^{2 \nu+1} dy, \  F^{-1}_{\nu}  g(y) = \frac{2^{-2 \nu}}{\Gamma^2 (\nu+1)}
\int\limits_{0}^{\infty} g(y) j_{\nu} \lr{y \eta} \eta^{2 \nu+1} d \eta,
\label{16}
\end{equation}
где $j_{\nu} (t) = \dfrac{2^{\nu} \Gamma (\nu+1) J_{\nu}
(t)}{t^{\nu}}$ "--- нормированная или малая функция Бесселя,
$J_{\nu}$ "--- обычная функция Бесселя первого рода.

По аналогии  можно ввести целое семейство других операторов
преобразования, которые являются обобщением бесселевых и
риссовских одномерных потенциалов. Некоторые приложения таких
операторов даны в работе~\cite{Kat1}.

Ж.~Дельсартом на базе ОП СПД было введено фундаментальное понятие
обобщённого сдвига.

\begin{definition}
{\it Оператором обобщённого сдвига} (ООС) называется решение
$u(x,y)=T_x^yf(x)$ задачи
\begin{equation}
\label{157}
(B_\nu)_y u(x,y)=(\frac{\partial^2}{\partial y^2}+\frac{2\nu +1}{y}\frac{\partial}{\partial y})\  u(x,y)= \frac{\partial^2}{\partial x^2} \ u(x,y),
\end{equation}
$$u(x,0)=f(x),\quad  u_y (x,0)=0.$$
\end{definition}

Название объясняется тем, что ООС в частном случае
$\nu=-\dfrac{1}{2}$  сводится к почти обычному сдвигу
$$T_x^yf(x)=\frac{1}{2}\left( f(x+y)+f(x-y)\right).$$
 Для ООС~\eqref{157} Дельсартом была получена явная формула
\begin{equation}
\label{158} T_x^y
f(x)=\frac{\Gamma(\nu+1)}{\sqrt{\pi}\Gamma(\nu+\frac{1}{2})}
\int\limits_0^{\pi}f(\sqrt{x^2+y^2-2xy \cos(t)}) \sin^{2\nu}t\,dt.
\end{equation}
Можно рассматривать в определении~\eqref{157} и произвольные пары
дифференциальных (или даже любых) операторов. Например, при таком
определении получаем привычный сдвиг:
$$\frac{\partial u}{\partial x}=\frac{\partial u}{\partial y},\quad u(x,0)=f(x),\quad  T_x^y f(x)=f(x+y).$$
Отметим, что ООС~\eqref{157}--\eqref{158} явно выражаются через ОП
СПД~\eqref{154}--\eqref{155} (см.~\cite{Lev2,Lev3, Mar9}).

Сделаем важное для дальнейшего замечание. С точки зрения
приложений к исследуемым в данной работе решениям дифференциальных
уравнений в частных производных с особенностями в коэффициентах
указанные операторы СПД обладают рядом недостатков, которые не
позволяют применять во многих важных случаях. К этим недостаткам
относится следующее: во-первых,  введённые выше операторы СПД
являются ОП лишь на множестве чётных функций, что исключает
возможность рассмотрения функций с особенностями в нуле;
во-вторых, они не сохраняют финитность и быстроубываемость на
бесконечности функций; в-третьих, они изменяют гладкость
преобразуемых функций. На этот факт впервые обратили внимание
Ж.-Л.~Лионс~\cite{Lio1}.

Таким образом, возникает необходимость введения и изучения других классов ОП для дифференциальных уравнений, содержащих операторы Бесселя.

\subsubsection{Дробные степени операторов Бесселя}\label{sec3.3.6}

Отметим, что существует теория дробных степеней операторов Бесселя
и их приложениям к дифференциальным уравнениям дробного порядка.
Этот подход позволяет определить дробные степени операторов
Бесселя не в образах интегрального преобразования Ханкеля неявно,
а в явном  виде как конкретный интегральный оператор со
специальными функциями в ядрах. Из работ этого направления
укажем~\cite{McB, Spr, Dim, Kir1, S140p, S135, S133, S127, S123,
18, S700, SS, FJSS}.

Безусловно, по аналогии с обычными производными, можно определить
такие дробные степени на полуоси при помощи естественного
свойства, что они действуют как умножение на степень в образах
преобразования Ханкеля. Такой подход оправдан за неимением лучшего
и позволяет получить ряд интересных результатов, хотя на этом пути
невозможно, по-видимому, получить явные представления дробных
степеней. Но представим на минуту, что мы умеем определять обычные
операторы дробного интегрирования Римана---Лиувилля только через
их действие в образах преобразований Лапласа или Меллина, а
интегральные формулы для этих операторов нам неизвестны. Тогда
сразу становится ясным, что подобная теория будет достаточно
бедной и лишится большинства своих наиболее полезных и красивых
результатов. Примерно в таком состоянии до недавнего времени
находилась теория дробных степеней оператора Бесселя, поэтому их
построение в явном интегральном виде является актуальной и
интересной задачей.

Приведём основные определения и свойства дробных степеней оператора Бесселя.

Мы рассматриваем вещественные степени сингулярного дифференциального оператора Бесселя
\begin{equation}\label{Bess}
B_\nu= D^2+\frac{\nu}{x}D,\qquad \nu\geq 0
\end{equation}
на вещественной полуоси $(0,\infty).$

\begin{definition}
Пусть $f(x)\in C^{2k}(0,b].$ Определим {\it правосторонний
оператор дробного интегрирования Бесселя} при условии
$f^{(i)}(b)=0, 0\leq i \leq 2k-1,$ $k \in \N$  по формуле
\begin{multline*}
(B_{b-}^{{\nu},k}f)(x)=\frac{1}{{\Gamma (2k)}}\int\limits_{x}^{b}
\left(\frac{y^{2}-x^{2}}{2y}\right)^{2k-1}
{_2}F_1(k+\frac{{\nu}-1}{2},k;2k;1-\frac{x^{2}}{y^{2}})
f(y)\,dy=\\
=\frac{\sqrt{{\pi}}}{2^{2k-1}{\Gamma (k)}}
\int\limits_{x}^{b}(y^{2}-x^{2})^{k-\frac{1}{2}}
\left(\frac{y}{x}\right)^{\frac{{\nu}}{2}}
P_{\frac{{\nu}}{2}-1}^{\frac{1}{2}-k}
\left(\frac{1}{2}\left(\frac{x}{y}+\frac{y}{x}\right)\right)
f(y)\,dy.
\end{multline*}
\end{definition}

\begin{definition} Определим {\it левосторонний оператор дробного интегрирования Бесселя} при условии $f^{(i)}(a)=0, 0\leq i \leq 2k-1, k \in N$  по формуле
\begin{multline*}
(B_{a+}^{{\nu},k}f)(x)=\frac{1}{{\Gamma (2k)}}\int\limits_{a}^{x}
\Bigl(\frac{x^{2}-y^{2}}{2x}\Bigr)^{2k-1}
{_2}F_1(k+\frac{{\nu}-1}{2},k;2k;;1-\frac{y^{2}}{x^{2}}) f(y)\,dy=
\\
=\frac{\sqrt{{\pi}}}{2^{2k-1}{\Gamma (k)}}
\int\limits_{a}^{x}(x^{2}-y^{2})^{\left(k-\frac{1}{2}\right)}\left(\frac{x}{y}\right)^{\frac{{\nu}}{2}}
P_{\frac{{\nu}}{2}-1}^{\frac{1}{2}-k}\left(\frac{1}{2}\left(\frac{x}{y}+\frac{y}{x}\right)\right)
 f(y)\,dy,
\end{multline*}
где ${_2}F_{1}$ "--- гипергеометрическая функция Гаусса,
$P_{\nu}^{\mu}(z)$ "--- функция Лежандра.
\end{definition}

Выражение дробных интегралов Бесселя через функции Лежандра
является полезным и является упрощением первоначального
определения, так как гипергеометрическая функция Гаусса зависит от
трёх параметров, а функция Лежандра "--- от двух.

Существует также версии дробных интегралов Бесселя с произвольными
пределами интегрирования, а также их дальнейшие модификации,
см.~\cite{S140p, S135, S133, S127,S123, 18, S700, SS}. Чаще всего
используются такие операторы:
$$
B_{0+}^{{\nu},k}, \quad  B_{\infty-}^{{\nu},k}.$$

\begin{property}\label{p1}
При $\nu=0$ дробный интеграл Бесселя на полуоси
$B_{0+}^{0,-\alpha}$ сводится к  дробному интегралу
Римана---Лиувилля, а именно, справедлива формула
$$
(B_{\infty-}^{0,
\alpha}f)(x)=\frac{1}{\Gamma(2\alpha)}\int\limits_x^{\infty}(y-x)^{2\alpha-1}f(y)dy=
(I_{-}^{2\alpha}f)(x).
$$
\end{property}

\begin{property}\label{p2}
Имеет место равенство
\begin{equation}\label{Prop2}
(B_{\infty-}^{{\nu},\alpha}f)(x)=\frac{1}{2^{2\alpha}}J_{x^2}^{2\alpha,\frac{\nu{-}1}{2}-\alpha,-\alpha}\left(x^{\frac{\nu{-}1}{2}}f(\sqrt{x})\right),
\end{equation}
где
\begin{equation}\label{Saigo1}
J_{x}\,^{\gamma,\beta,\eta}f(x)=\frac{1}{\Gamma(\gamma)}\int\limits_x^\infty(t-x)^{\gamma-1}t^{-\gamma-\beta}\,_2F_1\left(\gamma+\beta,-\eta;\gamma;1-\frac{x}{t}\right)f(t)dt
\end{equation} "--- дробный интеграл Сайго   $($см.~\cite{Rep}$)$. В~\eqref{Saigo1} $\gamma>0,\beta,\theta$ "--- вещественные числа.
\end{property}

\begin{property}\label{p3}
При $\lim\limits_{x\rightarrow +\infty}g(x)=0,$
$\lim\limits_{x\rightarrow +\infty}g'(x)=0$ получим, что
$$
(B_{\infty-}^{\nu,-1}B_\nu g)(x)=g(x).
$$
\end{property}

\begin{property}\label{p4}
При $x>0$ и $m+2\alpha+\nu<1$ справедлива формула
\begin{equation}\label{Prop4}
B_{\infty-}^{\nu,\alpha}\,x^m=x^{2\alpha+m}\,2^{-2\alpha}\,\Gamma\left[
                                                           \begin{array}{cc}
                                                              $$-\alpha-\dfrac{m}{2},$$ & $$-\dfrac{\nu-1}{2}-\alpha-\dfrac{m}{2}$$ \\
                                                             $$\dfrac{1-\nu-m}{2},$$ & $$-\dfrac{m}{2}$$  \\
                                                           \end{array}
                                                         \right].
\end{equation}
\end{property}

\begin{property}\label{p5}
Пусть $\alpha>0.$ Преобразования Меллина от дробного интеграла
Бесселя на полуоси имеет вид
\begin{equation}\label{Mellin2}
  M((B_{\infty-}^{\nu,\alpha}f)(x))(s)=\frac{1}{2^{2\alpha}}\,\,\Gamma\left[\begin{array}{cc}
$$\dfrac{s}{2},$$ & $$\dfrac{s}{2}-\dfrac{\nu-1}{2}$$ \\
$$\alpha+\dfrac{s}{2}-\dfrac{\nu-1}{2},$$ & $$\alpha+\dfrac{s}{2}$$  \\
                                                           \end{array}
                                                         \right] f^*(2\alpha+s).
\end{equation}
\end{property}

\begin{property}\label{p6}
Для дробного интеграла Бесселя на полуоси при $\alpha,\beta>0$
справедливо полугрупповое свойство
\begin{equation}\label{SemiGroup3}
    B_{\infty-}^{\nu,\alpha} B_{\infty-}^{\nu,\beta}f= B_{\infty-}^{\nu, \alpha+\beta}f.
\end{equation}
\end{property}

Далее приведём выражение для резольвенты дробных степеней
оператора Бесселя. Оно обобщает знаменитую формулу для дробных
интегралов Римана---Лиувилля, написанную без доказательства
Э.~Хилле и Я.\,Д.~Тамаркиным в работе   1930~г.~\cite{HT}. В этой
работе указывалось, что формула для резольвенты может быть
выведена методом преобразования Лапласа с использованием нового на
тот момент понятия свёртки в духе работ Дёйтча, но этот способ
похоже не был никогда реализован. Формула Тамаркина---Хилле была
на самом деле впервые доказана в монографии
М.\,М.~Джрбашяна~\cite{Dzh1} обычным для теории интегральных
уравнений методом последовательных приближений, хотя в монографии
Мхитара Мкртитевича нет упоминания, что доказательство даётся им
впервые, что характеризует этого замечательного математика.
Поэтому, возможно, исторически правильным было бы называть формулу
для резольвенты операторов дробного интегрирования
Римана---Лиувилля формулой \textit{Тамаркина---Хилле---Джрбашяна}.
Кроме того,  в~\cite{Dzh1} впервые были подробно изучены свойства
функции Миттаг-Лефлера, из этой книги отечественные математики
узнали о существовании подобной функции.

Формула  Тамаркина---Хилле---Джрбашяна  является самым известным
применением функций Миттаг-Лефлера, а также  основой огромного
числа работ по приложениям дробного исчисления к дифференциальным
уравнениям,  итальянские математики R.~Gorenflo и F.~Mainardi
предложили называть функцию Миттаг-Лефлера королевской функцией
теории дробного исчисления (Queen function of the fractional
calculus)~\cite{GoMa}.

\begin{property}\label{p7}
Для резольвенты оператора дробного интегрирования Бесселя
$B_{0+}^{\nu, \alpha},$ ${0\leq \nu <1}$ на подходящих функциях
справедливо интегральное представление
\begin{equation}
R_\lambda f=-\frac{1}{\lambda}f-\frac{1}{\lambda}\int\limits_0^x
K(x,y)f(y)\,dy,
\end{equation}
где ядро $K(x,y)$  выражается по формулам
\begin{gather}\label{R}
K(x,y)=\frac{2y}{x^2-y^2}\int\limits_0^1 S_{\alpha,\nu}(z(t))
\frac{dt} {
\left( t\left(1-t\right) \right) ^{\frac{\nu+1}{2}} },\\
\nonumber
z(t)=\left(\frac{t(1-t)\left(x^2-y^2\right)^2}{\left(1-\left(1-\frac{x^2}{y^2}
\right)t\right)4y^2}\right)^\alpha,\qquad
S_{\alpha,\nu}(z)=\sum\limits_{k=1}^{\infty}\frac{z^k}{\Gamma(\alpha
k+\frac{\nu-1}{2})\Gamma(\alpha k-\frac{\nu-1}{2})}
\end{gather} "--- разновидность гипергеометрической функции Райта---Фокса в
форме Райта $($см.~\cite{Kir4, Kir5, Kir6, KiSa}$)$.
\end{property}

Подобные  функции также применялись в работах
А.\,В.~Псху~\cite{Pshu1}, это специальные случаи более общих
функций Райта---Фокса, которые первоначально вводились как
обобщения функций Бесселя (см. выше о специальных функциях
Райта---Фокса). Интересной задачей является дальнейшее упрощение
представления ядра~\eqref{R}, если оно возможно.

Полученная формула для резольвенты дробных степеней оператора Бесселя позволяет рассматривать  задачи для обыкновенного интегродифференциального уравнения  вида
\begin{equation*}
B_\nu^\alpha u(x)-\lambda u(x)=f(x),
\end{equation*}
где $B_\nu^\alpha$ одна из дробных степеней оператора Бесселя, при
различных краевых условиях. По аналогии с известными результатами
возможно также рассмотрение уравнений в частных производных с
дробными степенями оператора Бесселя и их модификациями по
Герасимову---Капуто, которые можно ввести на основе имеющихся
обобщённых формул Тэйлора. К числу таких уравнений относится
обобщение $B$-эллиптического по терминологии
И.\,А.~Киприянова~\cite{Kip1} дробного уравнения Лапласа---Бесселя
\begin{equation*}
\sum\limits_{k=1}^n B_{\nu_k}^{\alpha_k} u(x_1,x_2,\ldots
x_n)=f(x_1,x_2,\ldots x_n),
\end{equation*}
нестационарное уравнение вида
\begin{equation*}
B_{\nu,t}^\alpha u(x,t)=\Delta_x u(x,t)+f(x,t).
\end{equation*}

Отметим, что  рассмотрение  спектральных свойств подобных
уравнений нуждается в изучении асимптотики функции $K(x,y)$ из
формулы~\eqref{R} в комплексной плоскости, а также распределения
её корней.

Многочисленные приложения операторов дробного
интегродифференцирования Римана---Ли\-у\-вил\-ля основаны на их
вхождении в остаточный член формулы Тэйлора. Поэтому после
определения дробных степеней оператора Бесселя сразу возникает
задача о построении обобщённой формулы Тэйлора, в которой функция
раскладывается по степеням оператора Бесселя. Эта задача возникла
достаточно давно и имеет некоторую историю.

Впервые формулы  разложения по степеням оператора Бесселя были
получены Жаном Дельсартом  (ряды Тэйлора---Дельсарта)~\cite{Del2,
Lev2, Lev3, Lev4}. Общий способ их построения изложен в~\cite{FN}
в терминах операторно-аналитических функций. Но ряды
Тэйлора---Дельсарта позволяют разложить по степеням оператора
Бесселя не обычный, а обобщённый сдвиг. По существу такие
разложения являются операторными вариантами рядов для функции
Бесселя, так же как обычные ряды Тэйлора являются операторными
версиями разложения в ряд экспоненты. Разумеется, ряды
Тэйлора---Дельсарта имеют свою область приложений. Но для
численного решения дифференциальных уравнений в частных
производных нужны обобщённые формулы и ряды Тэйлора несколько
другой природы. При пересчёте решения со слоя на слой, например,
методом сеток формулы для обобщённого сдвига совершенно
бесполезны, а нужны именно формулы для обычного сдвига, выражающие
решение на очередном рассчитываемом слое через его значения на
предыдущих слоях. Оказалось, что строить такие формулы для
обычного сдвига намного труднее, чем для обобщённого, так как они
уже не являются прямыми аналогами известных тождеств для
специальных функций.

Впервые с указанной мотивацией для применения к численному решению
уравнений с оператором Бесселя методом конечных элементов  формула
Тэйлора нужного типа была рассмотрена в работе
В.\,В.~Катрахова~\cite{KaKa}. Но полученный там результат может
рассматриваться только как первое приближение для желаемых формул
в явном виде, так как коэффициенты выражались неопределёнными
постоянными, задаваемыми системой рекуррентных соотношений, а ядро
остаточного члена представлялось некоторым многократным
интегралом. Это не случайно, угадать одновременно явный вид ядер и
остаточных членов невозможно, пока не известны конкретные
выражения для остатка в виде дробных степеней оператора Бесселя.

Тем не менее, поставленная В.\,В.~Катраховым задача была решена и
окончательный вид формулы Тэйлора со степенями операторов Бесселя
и остаточным членом в форме дробной степени оператора Бесселя был
найден Д.\,С.~Коноваловой и С.\,М.~Ситником в~\cite{S140p}, см.
также~\cite{S135,S133, S127,S123, 18, S700, SS}.

\chapter{Операторы преобразования Сонина---Пуассона---Дельсарта и их
модификации}\label{ch2}

В этой главе изучаются свойства операторов преобразования
Сонина---Пуассона---Дельсарта, Эрдейи---Кобера и вводятся новые
операторы преобразования. Устанавливается их связь с
преобразованиями Фурье, Ханкеля и с дробными интегралами
лиувиллевского типа. С этой целью используются как классические
лиувиллевские операторы, так и некоторые другие, занимающие
промежуточное положение между ними и бесселевыми потенциалами и
наследующими положительные свойства тех и других. Необходимость
введения указанных операторов была вызвана тем, что лиувиллевские
операторы не ограничены в соболевских пространствах в случае
неограниченных областей. Рассматривается сведение пространств
С.\,Л.~Соболева с помощью операторов преобразования к весовым
функциональным пространствам, введённым И.\,А.~Киприяновым. Здесь
же вводятся и изучаются новые функциональные пространства (в
одномерном случае).

Далее строится один класс многомерных операторов преобразования, преобразующих многомерный оператор Лапласа в обыкновенный оператор двукратного дифференцирования по радиальной переменной. Эти операторы определены на функциях, у которых допускается наличие особенности в одной точке, в качестве которой выбрано начало координат.

\section{Одномерные операторы преобразования}\label{sec4}

\subsection{Основные конструкции операторов
преобразования}\label{sec4.1}

 Напомним общее определение (см. Ж.-Л.~Лионс~\cite{58})    операторов преобразования. Пусть $A$ и $B$ "--- некоторые линейные операторы. Операторы $P$ и $S$ называются операторами преобразования для $A$ и $B,$   если    имеют место формулы
\begin{equation}
B = P A S, \   A = S B P
\label{1.1.1}
\end{equation}
и операторы $P$ и $S$ взаимно обратны. Данное определение не
является вполне строгим, поскольку не указаны области определения
и значения участвующих операторов. Их описание мы будем давать в
каждом конкретном  случае  отдельно. Простейшим  примером
операторов преобразования может    служить преобразование  Фурье.
В этом случае $P$ и $S,$ соответственно, прямое и обратное
преобразование Фурье, $B$ "--- оператор дифференцирования, $A$
"--- оператор умножения на двойственную переменную. В этом примере
операторы $A$ и $B$ имеют различную природу. Мы же на протяжении
всей книги будем рассматривать случай, когда $A$ и $B$ будут
дифференциальными операторами. В частности, в этом параграфе мы
рассмотрим в качестве $A$ и $B$ следующие   операторы:
\begin{equation*}
A=\frac{\pr^2}{\pr y^2}, \quad    B=B_{\nu}= \frac{\pr^2}{\pr y^2}
+\frac{2 \nu+1}{y}\frac{\pr}{\pr y}.
\end{equation*}
Оператор $B_{\nu}$ принято называть оператором Бесселя с
параметром $\nu.$ В этой ситуации известны~\cite{Lev7, Mar2}
следующие операторы преобразования:
\begin{equation}
P_{\nu, 0}^{\nu+\frac{1}{2}} f(y) = \frac{2\, \Gamma \lr{\nu+1} }{\sqrt{\pi}\, \Gamma \lr{\nu+\frac{1}{2}}} \, y^{-2\nu} \int\limits_0^y \lr{y^2-t^2}^{\nu-\frac{1}{2}} f(t) \, dt,
\label{1.1.2}
\end{equation}
\begin{equation}
S_{\nu, 0}^{-\nu-\frac{1}{2}} f(y) = \frac{\sqrt{\pi}}{ \Gamma \lr{\nu+1} \, \Gamma \lr{\frac{1}{2} - \nu}} \frac{\pr}{\pr y} \int\limits_0^y \lr{y^2-t^2}^{-\nu-\frac{1}{2}} t^{2 \nu+1} f(t) \, dt.
\label{1.1.3}
\end{equation}
Здесь и всюду ниже через $\Gamma(\mu)$ обозначается гамма-функция
Эйлера. Функция $f(y)$ определена на полупрямой $E_{+}^1 = \{y>0
\}$ и является гладкой, финитной и чётной функцией. В тех случаях,
когда интегралы в предыдущих формулах расходятся, нужно перейти к
их регуляризации. Оператор $P_{\nu, 0}^{\nu+\frac{1}{2}}$ принято
называть оператором Пуассона, а оператор $S_{\nu,
0}^{-\nu-\frac{1}{2}}$ "--- оператором  Сонина. Эти названия
заимствованы из теории цилиндрических функций, где так называют
подобные интегралы в формулах~\eqref{1.1.2} и~\eqref{1.1.3} при
некоторой конкретной функции $f.$ Для произвольных функций эти
операторы были введены Ж.~Дельсартом. Поэтому мы также будем
использовать общее название "--- операторы преобразования
Сонина---Пуассона---Дельсарта.

С точки зрения приложений к исследуемым в работе сингулярным
эллиптическим краевым задачам в частных производных указанные
операторы обладают рядом недостатков, которые не позволяют
применять их в данном случае. К этим недостаткам относятся
следующие: во-первых, операторы $P_{\nu, 0}^{\nu+\frac{1}{2}}$ и
$S_{\nu, 0}^{-\nu-\frac{1}{2}}$ удовлетворяют
формулам~\eqref{1.1.1} лишь на множестве чётных функций (нам же
необходимо изучать растущие вблизи начала функции); во-вторых, они
не сохраняют финитность и быстрое убывание функций на
бесконечности; в-третьих, они изменяют гладкость преобразуемых
функций. На последний факт впервые обратил внимание
Ж.-Л.~Лионс~\cite{Lio1, Lio2, Lio3}. Впрочем, это нетрудно
усмотреть и непосредственно.

Введем операторы преобразования, свободные от всех указанных выше
недостатков. Для этого найдем сначала для операторов~\eqref{1.1.2}
и~\eqref{1.1.3}  формально сопряжённые относительно билинейных
форм вида
\begin{equation*}
\lr{f, g}_{\nu} = \int\limits_0^{\infty}  y^{2 \nu+1} f(y) g(y) \, dy.
\end{equation*}
Имеем
\begin{equation*}
\lr{P_{\nu, 0}^{\nu+\frac{1}{2}} g, f}_{\nu} = \int\limits_0^{\infty} t^{2 \nu +1} f(t)  \frac{2 \, \Gamma \lr{\nu+1}}{\sqrt{\pi}  \, \Gamma \lr{ \nu+\frac{1}{2} }} \, t^{- 2 \nu}  \int\limits_0^t \lr{t^2-y^2}^{\nu-\frac{1}{2}} g(y) \, dy dt =
\end{equation*}
\begin{equation*}
 =  \frac{2 \, \Gamma \lr{\nu+1}}{\sqrt{\pi}  \, \Gamma \lr{ \nu+\frac{1}{2} }}  \int\limits_0^{\infty} g(y)      \int\limits_0^t t f(t) \lr{t^2-y^2}^{\nu-\frac{1}{2}}  \, dt dy =  \lr{g,  S_{\nu}^{\nu+\frac{1}{2}} }_{-\frac{1}{2}},
\end{equation*}
где оператор $S_{\nu}^{\nu+\frac{1}{2}}$ определяется по формуле
\begin{equation}
S_{\nu}^{\nu+\frac{1}{2}} f(y)  =  \frac{2 \, \Gamma \lr{\nu+1}}{\sqrt{\pi}  \, \Gamma \lr{ \nu+\frac{1}{2} }}  \int\limits_y^{\infty}  \lr{t^2-y^2}^{\nu-\frac{1}{2}} t f(t)  \, dt.
\label{1.1.4}
\end{equation}
Аналогично получаем
\begin{equation*}
    \lr{S_{\nu, 0}^{-\nu-\frac{1}{2}} g, f}_{-\frac{1}{2}} = \frac{\sqrt{\pi} }{  \Gamma \lr{\nu+1} \, \Gamma \lr{ \frac{1}{2}-\nu }} \int\limits_0^{\infty} \frac{\pr}{\pr t} \int\limits_0^t       \lr{t^2-y^2}^{-\nu-\frac{1}{2}} y^{2 \nu+1} g(y) \, dy  f(t) \, dt =
\end{equation*}
\begin{equation*}
 =\frac{-\sqrt{\pi} }{  \Gamma \lr{\nu+1} \, \Gamma \lr{ \frac{1}{2}-\nu }} \int\limits_0^{\infty} y^{2 \nu+1} g(y)   \int\limits_y^{\infty}       \lr{t^2-y^2}^{-\nu-\frac{1}{2}} \frac{\pr f(t)}{\pr t}   dt dy =  \lr{g,  P_{\nu}^{-\nu-\frac{1}{2}} }_{\nu},
\end{equation*}
где положено
\begin{equation}
    P_{\nu}^{-\nu-\frac{1}{2}} f(y)  =  \frac{-\sqrt{\pi} }{  \Gamma \lr{\nu+1} \, \Gamma \lr{ \frac{1}{2}-\nu }}  \int\limits_y^{\infty}  \lr{t^2-y^2}^{-\nu-\frac{1}{2}} \frac{\pr f(t)}{\pr t} \, dt.
    \label{1.1.5}
\end{equation}
Нам будет удобно выделить другие операторы, которые получаются
из~\eqref{1.1.4} и~\eqref{1.1.5}  перенесением оператора
дифференцирования из второго в первый оператор:
\begin{equation}
S_{\nu}^{\nu-\frac{1}{2}} f(y) = \frac{- 2 \, \Gamma \lr{\nu+1}}{  \sqrt{\pi} \, \Gamma \lr{ \nu+\frac{1}{2} }} \frac{\pr}{\pr y} \int\limits_y^{\infty} \lr{t^2-y^2}^{\nu-\frac{1}{2}} t f(t) \, dt,
\label{1.1.6}
\end{equation}
\begin{equation}
P_{\nu}^{\frac{1}{2}-\nu} f(y) = \frac{\sqrt{\pi} }{\Gamma \lr{\nu+1}\, \Gamma \lr{\frac{1}{2}-\nu}} \int\limits_y^{\infty} \lr{t^2-y^2}^{-\nu-\frac{1}{2}} f(t) \, dt,
\label{1.1.7}
\end{equation}
где, как и ранее, под $a^{\mu},$ $a>0,$ понимается главная ветвь
многозначной функции: $a^{\mu} = \exp (\mu \ln a).$

Предыдущие выкладки носили несколько формальный характер. Они
приведены, чтобы лучше была видна связь операторов Пуассона и
Сонина с операторами Эрдейи---Кобера.

Введем некоторые обозначения. Пусть $R$ обозначает положительное
число или бесконечность. Через $C^{\infty} (0, R)$ обозначается
множество бесконечно дифференцируемых на интервале $(0, R)$
функций. $\mathring{C}^{\infty} (0, R)$ обозначает подмножество
функций из $C^{\infty} (0, R),$ имеющих компактный в $(0, R)$
носитель. Через $C^{\infty} [0, R)$ обозначаем  подмножество
функций из $C^{\infty} (0, R),$ все производные которых непрерывны
вплоть до левого конца. Символ $\mathring{C}^{\infty} [0, R)$
обозначает подмножество функций из $C^{\infty} [0, R),$
обращающихся в нуль в окрестности правого конца. Через
$\mathring{C}^{\infty}_{\{0\}} (0, R)$ обозначим подмножество
функций из $C^{\infty} (0, R)$ обращающихся в нуль в окрестности
правого конца. Отметим, что функции из
$\mathring{C}^{\infty}_{\{0\}} (0, R)$ могут расти произвольным
образом в окрестности точки $y=0.$ Будем называть функцию $f \in
C^{\infty} [0, R)$ чётной (нечётной), если $D^k f(0) = 0$ при всех
нечётных (чётных) неотрицательных значениях $k.$ Здесь    $D=
\dfrac{\pr}{\pr y},$ $D^k = D D^{k-1}.$ Функция $f \in C^{\infty}
[0, R)$ будет чётной тогда и только тогда, когда функция $g(y)=f
(\sqrt{y})$ принадлежит   пространству $C^{\infty} [0, R^2).$ Этот
факт элементарно доказывается с помощью формулы Тейлора. Множество
чётных (нечётных) функций обозначим через $C^{\infty}_{+} [0, R)$
($C^{\infty}_{-} [0, R)$). Пусть также
$\mathring{C}^{\infty}_{\pm} [0, R]=C^{\infty}_{\pm} [0, R) \cap
\mathring{C}^{\infty} [0, R)$    Кроме того, будем использовать
обозначение $(0, \infty) = E^1_+,$ $[0, \infty) = \ov{E^1_+}.$

Определим операторы преобразования  $S_{\nu}^{\nu - \frac{1}{2}}$
при $\Re \nu \geq 0$ на функциях $f \in
\mathring{C}^{\infty}_{\{0\}} \lr{E^1_+}$ по
формуле~\eqref{1.1.6}. Операторы $P_{\nu}^{\nu - \frac{1}{2}}$ на
том же классе функций при $0 \leq \Re \nu < \dfrac{1}{2}$
определим по формуле~\eqref{1.1.7}. Для $\Re \nu \geq
\dfrac{1}{2}$ значение функции $P_{\nu}^{ \frac{1}{2} - \nu }
f(y)$ в точке $y>0$ определим как аналитическое продолжение
интеграла~\eqref{1.1.7} по параметру $\nu.$  Если $\dfrac{1}{2}
\leq \Re \nu <N+ \dfrac{1}{2},$ где $N$ "--- натуральное  число,
то это равносильно заданию оператора $P_{\nu}^{ \frac{1}{2} - \nu
}$ по  формуле
\begin{equation}
P_{\nu}^{\frac{1}{2}-\nu} f(y) = \frac{(-1)^{N} 2^{-N} \sqrt{\pi} } {\Gamma \lr{\nu+1}\, \Gamma \lr{N-\nu+\frac{1}{2}}} \int\limits_y^{\infty} \lr{t^2-y^2}^{N-\nu-\frac{1}{2}} \lr{\frac{\pr}{\pr t} \frac{1}{t}}^N f(t) \, dt,
\label{1.1.8}
\end{equation}
в которой интеграл справа является сходящимся.

Введённые операторы действительно являются операторами
преобразования, поскольку справедлива

\begin{theorem}\label{theorem:1_1_1}
Операторы $P_{\nu}^{ \frac{1}{2} - \nu }$ и  $S_{\nu}^{\nu - \frac{1}{2}}$ при $\Re \nu \geq 0$ взаимно однозначно отображают пространство  $\mathring{C}_{\{0\}}^{\infty} \lr{E_{+}^{1}}$ на себя и являются взаимно обратными. Имеют место формулы
\begin{equation}
B_{\nu} P_{\nu}^{\frac{1}{2}-\nu}  = P_{\nu}^{\frac{1}{2}-\nu} D^2,\
D^2 S_{\nu}^{\nu-\frac{1}{2}}  = S_{\nu}^{\nu-\frac{1}{2}} B_{\nu}.
\label{1.1.9}
\end{equation}
\end{theorem}

\begin{proof}
Пусть $\Re \nu < N + \dfrac{1}{2},$ где $N$ "--- натуральное
число. После замены переменных $t \to ty$ формула~\eqref{1.1.8}
принимает вид
\begin{multline}
P_{\nu}^{\frac{1}{2}-\nu} f(y) = \frac{(-1)^{N} 2^{-N} \sqrt{\pi} \, y^{2(N-\nu)}} {\Gamma \lr{\nu+1}\, \Gamma \lr{N-\nu+\frac{1}{2}}}  \lr{\frac{\pr}{\pr y} \frac{1}{y}}^N  \int\limits_1^{\infty} \lr{t^2-1}^{N-\nu-\frac{1}{2}} t^{-2N} f(t y) \, dt=   \\
 = \frac{(-1)^{N} 2^{-N} \sqrt{\pi}\, y^{2(N-\nu)}} {\Gamma \lr{\nu+1}\, \Gamma \lr{N-\nu+\frac{1}{2}}}  \lr{\frac{\pr}{\pr y} \frac{1}{y}}^N  \int\limits_y^{\infty} y^{2 \nu} \lr{t^2-y^2}^{N-\nu-\frac{1}{2}} t^{-2N} f(t) \, dt.
  \label{1.1.10}
\end{multline}
Отсюда сразу следует, что функция $P_{\nu}^{ \frac{1}{2} - \nu }
f(y)$ бесконечно дифференцируема при $y>0$ и финитна,  если $f \in
\mathring{C}_{\{0\}}^{\infty} \lr{E_{+}^{1}}.$ Более того, точная
верхняя грань носителя при этом не увеличивается. Таким образом,
оператор $P_{\nu}^{ \frac{1}{2} - \nu }$ отображает пространство
$\mathring{C}_{\{0\}}^{\infty} \lr{E_{+}^{1}}$ в себя. Аналогичным
образом это доказывается и для оператора $S_{\nu}^{\nu -
\frac{1}{2}}.$ Покажем, что $P_{\nu}^{ \frac{1}{2} - \nu } =
\lr{S_{\nu}^{\nu - \frac{1}{2}}}^{-1}.$     По теореме Фубини при
$\Re \nu < N + \dfrac{1}{2}$ имеем
$$
P_{\nu}^{\frac{1}{2}-\nu} S_{\nu}^{\nu-\frac{1}{2}} f(y) = \frac{(-1)^{N+1} 2^{1-N} y^{2(N-\nu)}} {\Gamma \lr{\nu+\frac{1}{2}}\, \Gamma \lr{N-\nu+\frac{1}{2}}}  \lr{\frac{\pr}{\pr y} \frac{1}{y}}^N  \int\limits_y^{\infty} y^{2 \nu} \tau^{1-2 \nu} \times
$$
$$
\times \frac{\pr}{\pr \tau} \lr{\tau^{2 \nu+1} f (\tau)}  \int\limits_y^{\tau} \lr{t^2-y^2}^{N-\nu-\frac{1}{2}}  t^{-2N-1} \lr{\tau^2-t^2}^{\nu-\frac{1}{2}} \, dt d \tau.
$$
Во внутреннем интеграле произведем замену переменных по формуле
$\dfrac{1}{t^2} = \dfrac{1}{y^2} + z
\lr{\dfrac{1}{\tau^2}-\dfrac{1}{y^2}}.$ Тогда получим
$$
\int\limits_y^{\tau} \lr{t^2-y^2}^{N-\nu-\frac{1}{2}}  t^{-2N-1}
\lr{\tau^2-t^2}^{\nu-\frac{1}{2}} \, dt = \frac{1}{2} \, y^{2 N -
2 \nu -1 } \tau^{2 \nu -1}  \lr{\frac{1}{y^2}-\frac{1}{\tau^2}}^N
 \int\limits_0^1 z^{N-\nu-\frac{1}{2}} (1-z)^{\nu-\frac{1}{2}} dz
=
$$
$$
=\frac{(\tau^2-y^2)^{2 \nu - 2N-1 }\,\Gamma \lr{N-\nu+\frac{1}{2}}
\Gamma \lr{\nu+\frac{1}{2}}}{2 \, y^{2 \nu+1} \Gamma \lr{\nu+1}}.
$$
Следовательно,
$$
P_{\nu}^{\frac{1}{2}-\nu} S_{\nu}^{\nu-\frac{1}{2}} f(y) = \frac{(-1)^{N}  y^{2N-2\nu}} {2^N N!}   \lr{\frac{\pr}{\pr y} \frac{1}{y}}^N  \int\limits_y^{\infty} \frac{1}{y} (\tau^2-y^2)^N \tau^{-2 N}  \frac{\pr}{\pr \tau} \lr{\tau^{2 \nu+1} f (\tau)} d \tau.
$$
Интегрируя один раз по частям в этом интеграле, а затем дифференцируя по параметру получим
$$
P_{\nu}^{\frac{1}{2}-\nu} S_{\nu}^{\nu-\frac{1}{2}} f(y) = \frac{(-1)^{N} 2^N} {2^N N!} y^{2N-2\nu}  \frac{\pr}{\pr y} \lr{\frac{\pr}{ y \pr y} }^{N-1}  \int\limits_y^{\infty}  (\tau^2-y^2)^{N-1} \tau^{2 \nu-2 N}   f (\tau) d \tau=
$$
$$
=\frac{-N!} { N!} y^{2N-2\nu}  \frac{\pr}{\pr y}
\int\limits_y^{\infty}  \tau^{2 \nu-2 N}   f (\tau) d \tau=f(y).
$$
Аналогичным образом устанавливается формула $S_{\nu}^{\nu -
\frac{1}{2}} P_{\nu}^{\frac{1}{2}-\nu} f(y)=f(y).$ Стало быть
доказано, что операторы $P_{\nu}^{\frac{1}{2}-\nu}$   и
$S_{\nu}^{\nu - \frac{1}{2}}$ осуществляют взаимно однозначное
отображение пространства $\mathring{C}_{\{0\}}^{\infty}
\lr{E_{+}^{1}}$ на себя.

Докажем формулы~\eqref{1.1.9}. Достаточно проверить одну из них,
например, первую, поскольку другая есть её следствие. Ввиду
возможности применения принципа аналитического продолжения
достаточно рассмотреть случай $\Re \nu < \dfrac{1}{2}.$ Для $f \in
\mathring{C}_{\{0\}}^{\infty} \lr{E_{+}^{1}}$ из~\eqref{1.1.10}
получаем
\begin{multline}
B_{\nu} P_{\nu}^{\frac{1}{2}-\nu} f(y) = \frac{\sqrt{\pi}} {\Gamma \lr{\nu+1}\, \Gamma \lr{\frac{1}{2}-\nu}}  \lr{\frac{\pr^2}{\pr y^2}+ \frac{2 \nu+1}{y} \frac{\pr}{\pr y}}  \int\limits_y^{\infty} \lr{t^2-y^2}^{-\nu-\frac{1}{2}}  f(t) \, dt=   \\
= \frac{\sqrt{\pi}} {\Gamma \lr{\nu+1}\, \Gamma
\lr{\frac{1}{2}-\nu}}    \int\limits_1^{\infty}
\lr{t^2-1}^{-\nu-\frac{1}{2}} \lr{\frac{\pr^2}{\pr y^2}+ \frac{2
\nu+1}{y} \frac{\pr}{\pr y}} \lr{y^{- 2 \nu} f(t y)} \, dt.
\label{1.1.11}
\end{multline}
Далее, поскольку
$$
\lr{\frac{\pr^2}{\pr y^2}+ \frac{2 \nu+1}{y} \frac{\pr}{\pr y}} \lr{y^{- 2 \nu} f(t y)} = y^{- 2 \nu} t^2 \left. \frac{\pr^2 f (\tau)}{\pr \tau^2} \right|_{\tau = t y} + (1-2 \nu) y^{- 2 \nu - 1 } t \left. \frac{\pr f (\tau)}{\pr \tau} \right|_{\tau = t y},
$$
$$
 (1-2 \nu) y^{- 2 \nu } \int\limits_1^{\infty}  t
\lr{t^2-1}^{-\nu-\frac{1}{2}} \left. \frac{\pr f (\tau)}{\pr \tau}
\right|_{\tau = t y} dt = -  y^{- 2 \nu - 1 }
\int\limits_1^{\infty} \frac{\pr \lr{t^2-1}^{\frac{1}{2}-\nu}}{\pr
t} \left. \frac{\pr f (\tau)}{\pr \tau} \right|_{\tau = t y} dt =
$$
$$
=- y^{- 2 \nu} \int\limits_1^{\infty}
\lr{t^2-1}^{-\nu-\frac{1}{2}} (t^2-1) \left. \frac{\pr^2 f
(\tau)}{\pr \tau^2} \right|_{\tau = t y} dt,
$$
то правая часть формулы~\eqref{1.1.11} приводится к виду
$$
\frac{\sqrt{\pi}} {\Gamma \lr{\nu+1}\, \Gamma \lr{\frac{1}{2}-\nu}}   y^{- 2 \nu } \int\limits_1^{\infty} \lr{t^2-1}^{-\nu-\frac{1}{2}}  \left. \frac{\pr^2 f (\tau)}{\pr \tau^2} \right|_{\tau = t y} dt.
$$
Этим завершается доказательство формулы~\eqref{1.1.9}, а вместе с
тем и теоремы в целом.
\end{proof}

Теорему~\ref{theorem:1_1_1} следует считать известной, она
доказывалась многими авторами при различных
ограничениях~\cite{SKM}. Мы предпочли дать её полное
доказательство, поскольку сделать конкретную ссылку оказалось
затруднительно, и так как она играет ключевую роль в теории
операторов преобразования.

Лиувиллевский оператор $I^{\mu}$ при $\Re \mu > 0$ на функциях  $f \in \mathring{C}_{\{0\}}^{\infty} \lr{E_{+}^{1}}$ определяется по  формуле
\begin{equation}
I^{\mu} f (y) = \frac{1}{\Gamma (\mu)} \int\limits_y^{\infty} (t-y)^{\mu-1} f(t) \, dt.
\label{1.1.12}
\end{equation}
Для остальных значений комплексного параметра $\mu$ функция
$I^{\mu} f (y)$ определяется с помощью аналитического продолжения
по параметру $\mu.$ Если  $\Re \mu > -M,$ где $M$ "---
натуральное число или нуль, то это эквивалентно заданию $I^{\mu} f
(y)$   по  формуле
\begin{equation}
I^{\mu} f (y) = \frac{(-1)^M}{\Gamma \lr{\mu+M}} \int\limits_y^{\infty} (t-y)^{M+\mu-1} \frac{\pr^M f(t)}{\pr t^M} \, dt.
\label{1.1.13}
\end{equation}
Следовательно, операторы $I^{\mu}$ определены на
$\mathring{C}_{\{0\}}^{\infty} \lr{E_{+}^{1}}$ при всех
комплексных $\mu.$ Они, во-первых, взаимно однозначно отображают
пространство $\mathring{C}_{\{0\}}^{\infty} \lr{E_{+}^{1}}$  на
себя, а во-вторых, справедливо групповое свойство
\begin{equation}
I^{\mu} I^{\nu} = I^{\nu} I^{\mu} = I^{\mu+\nu}.
\label{1.1.14}
\end{equation}

Определим теперь операторы $P_{\nu}$ и $S_{\nu}$ при $\Re \nu \geq 0$ на функциях $f \in \mathring{C}_{\{0\}}^{\infty} \lr{E_{+}^{1}}$    по  формуле
\begin{equation}
P_{\nu} f = P_{\nu}^{\frac{1}{2} - \nu} I^{\nu -\frac{1}{2}} f, \   S_{\nu} f = I^{\frac{1}{2}-\nu} S_{\nu}^{ \nu - \frac{1}{2} }  f.
\label{1.1.15}
\end{equation}

Из перечисленных выше свойств лиувиллевских операторов и
теоремы~\eqref{theorem:1_1_1} следует

\begin{theorem} \label{theorem:1_1_2}
При $\Re \nu \geq 0$ операторы $P_{\nu}$ и $S_{\nu}$ взаимно однозначно отображают пространство   $\mathring{C}_{\{0\}}^{\infty} \lr{E_{+}^{1}}$ на себя  и   являются взаимно обратными. Для $f \in \mathring{C}_{\{0\}}^{\infty} \lr{E_{+}^{1}}$ справедливы формулы
\begin{equation}
B_{\nu} P_{\nu} f  = P_{\nu} D^2 f,\
D^2 S_{\nu} f  = S_{\nu} B_{\nu} f.
\label{1.1.16}
\end{equation}
\end{theorem}

Следовательно, операторы $P_{\nu}$ и $S_{\nu}$ действительно являются операторами преобразования. Они допускают
следующие представления, вывод которых стандартен.

Пусть функция $f \in \mathring{C}_{\{0\}}^{\infty}
\lr{E_{+}^{1}}.$ Тогда   при $0 \leq \Re \nu < \dfrac{1}{2}$ из
формул~\eqref{1.1.15} и~\eqref{1.1.13} получаем
$$
P_{\nu} f (y) =  P_{\nu}^{\frac{1}{2} - \nu} I^{\nu -\frac{1}{2}} f (y)=- P_{\nu}^{\frac{1}{2} - \nu} I^{\nu +\frac{1}{2}} D f (y)=
$$
$$
=\frac{- \sqrt{\pi}} {\Gamma \lr{\nu+1}\, \Gamma \lr{\frac{1}{2}-\nu} \, \Gamma \lr{\nu+\frac{1}{2}} }   \int\limits_y^{\infty} \lr{t^2-y^2}^{-\nu-\frac{1}{2}} \int\limits_t^{\infty} \lr{\tau-t}^{\nu-\frac{1}{2}}  D f(\tau)   d \tau dt =
$$
$$
=\frac{- \sqrt{\pi}} {\Gamma \lr{\nu+1}\, \Gamma \lr{\frac{1}{2}-\nu} \, \Gamma \lr{\nu+\frac{1}{2}} }   \int\limits_y^{\infty} D f(\tau) \int\limits_y^{\tau} \lr{t^2-y^2}^{-\nu-\frac{1}{2}} \lr{\tau-t}^{\nu-\frac{1}{2}}  dt  d \tau.
$$
Для вычисления внутреннего интеграла введем новую переменную по
формуле $z = \dfrac{t-y}{\tau-y}.$   Тогда
\begin{multline}
\int\limits_y^{\tau} \lr{t^2-y^2}^{-\nu-\frac{1}{2}}
\lr{\tau-t}^{\nu-\frac{1}{2}}  dt = (2 y)^{-\nu-\frac{1}{2}}
\int\limits_0^1 z^{-\nu-\frac{1}{2}} (1-z)^{\nu-\frac{1}{2}}
\lr{1-\frac{y - \tau}{2y}z}^{-\nu-\frac{1}{2}} dz = \\
=\frac{\Gamma \lr{\frac{1}{2}-\nu} \, \Gamma \lr{\frac{1}{2}+\nu}
} {(2y)^{\nu+\frac{1}{2}}}\, {_2F_1} \lr{\nu+\frac{1}{2};
\frac{1}{2}-\nu; l; \frac{y - \tau}{2y}}, \label{1.1.17}
\end{multline}
где через ${_2F_1 (a, b; c; \zeta)}$ обозначена
гипергеометрическая функция Гаусса. Выше была использована формула
Эйлера~\cite[с.~72]{BE1}
\begin{equation}
{_2F_1} (a, b; c; \zeta) = \frac{\Gamma (c)}{\Gamma (b) \, \Gamma (c-b)} \int\limits_0^1 t^{b-1} (1-t)^{c-b-1} (1 - \zeta t)^{-a} dt.
\label{1.1.18}
\end{equation}
Функция Лежандра первого рода $P_{\mu}^0 (\zeta)$ при   $\zeta >
0$ может быть определена формулой~\cite[с.~156]{BE1}
\begin{equation}
P_{\mu}^0 (\zeta) = \frac{1}{\pi} \int\limits_0^{\pi} \lr{\zeta+\sqrt{\zeta^2-1} \cos t}^{\mu} dt.
\label{1.1.19}
\end{equation}

Отметим, что функция $P_{\mu}^0 (\zeta)$ является аналитической
функцией комплексного параметра $\mu.$ Гипергеометрическая функция
из формулы~\eqref{1.1.18} выражается через функцию Лежандра по
следующей формуле (см.~\cite[с.~127]{BE1}):
$$
 {_2F_1} \lr{\nu+\frac{1}{2}; \frac{1}{2}-\nu; l; \frac{1}{2} - \frac{1}{2} \frac{\tau}{y}} = P^0_{\nu-\frac{1}{2}} \lr{\frac{\tau}{y}}.
$$
Таким образом, получено следующее представление оператора $P_{\nu}$:
\begin{equation}
P_{\nu} f (y)  = \frac{ -\sqrt{\pi}} {2^{\nu+\frac{1}{2}} \, \Gamma \lr{\nu+1}} \, y^{-\nu-\frac{1}{2}}\, \int\limits_y^{\infty} \frac{\pr f(\tau)}{\pr \tau}\,  P^0_{\nu-\frac{1}{2}}  \lr{\frac{\tau}{y}} d \tau.
\label{1.1.20}
\end{equation}
Эта формула доказана при дополнительном ограничении $0 \leq \Re \nu < \dfrac{1}{2}.$  Однако нетрудно заметить, что при
фиксированной функции $f$ и фиксированном значении $y>0$ и левая и
правая части суть аналитические функции параметра $\nu$ при $\Re \nu \geq 0.$ Поэтому, в силу единственности аналитического
продолжения формула~\eqref{1.1.20} имеет место для всех   $y > 0,$
$f \in \mathring{C}_{\{0\}}^{\infty} \lr{E_{+}^{1}}$ и для всех
$\nu$  с $\Re \nu \geq 0.$

Присоединенная функция Лежандра $P_{\mu}^{-1} (\zeta)$ определяется по формуле
$$
\frac{\pr }{\pr \zeta} \, P_{\mu}^0 (\zeta) = \frac{\mu
(\mu+1)}{\sqrt{\zeta^2 -1}} \, P_{\mu}^{-1} (\zeta), \quad   \zeta
\geq 1.
$$
Тогда, интегрируя по частям в~\eqref{1.1.20}, получаем
\begin{equation}
P_{\nu} f (y) = \frac{ \sqrt{\pi} \, y^{-\nu-\frac{1}{2}}} {2^{\nu+\frac{1}{2}} \, \Gamma \lr{\nu+1}}  \lr{ f(y)+\lr{\nu^2 - \frac{1}{4}} \int\limits_1^{\infty} f (t y)  \frac{P_{\nu-\frac{1}{2}}^{-1} (t)}{\sqrt{t^2 -1}}\, dt}.
\label{1.1.21}
\end{equation}
Из такой записи оператора преобразования $P_{\nu}$ видно, что он
не изменяет гладкость преобразуемой функции, поскольку ядро
интегрального оператора в~\eqref{1.1.21} не имеет несуммируемых
особенностей при   $1 \leq t < \infty.$ Это утверждение, конечно,
справедливо лишь при положительных аргументах. При $y=0$
особенность у функции $P_{\nu} f (y)$ может возникать и в том
случае, если функция $f$ её и не имела. Здесь играет роль
поведение ядра на бесконечности. Например, $P_{\frac{1}{2}} f (y)=
\frac{1}{y} f(y).$

Выведем аналогичное представление для оператора $S_{\nu}.$ При  $0
\leq \Re \nu < \dfrac{1}{2}$ и $f \in
\mathring{C}_{\{0\}}^{\infty} \lr{E_{+}^{1}}$ из
формул~\eqref{1.1.12},~\eqref{1.1.16} имеем
\begin{multline}
S_{\nu} f (y) =  I^{\frac{1}{2} - \nu} S_{\nu}^{\nu -\frac{1}{2}}
f (y) = \frac{- 2 \, \Gamma \lr{\nu+1}} {\sqrt{\pi} \, \Gamma
\lr{\frac{1}{2}+\nu} \, \Gamma \lr{\frac{1}{2}-\nu} } \frac{\pr
}{\pr y} \int\limits_y^{\infty}   \lr{t-y}^{\nu-\frac{1}{2}}
\int\limits_t^{\infty} \lr{\tau^2-t^2}^{-\nu-\frac{1}{2}} \tau f
(\tau) \,     d \tau dt  =  \\
 =\frac{- 2 \, \Gamma \lr{\nu+1}} {\sqrt{\pi}\, \Gamma
\lr{\frac{1}{2}+\nu} \, \Gamma \lr{\frac{1}{2}-\nu} } \frac{\pr
}{\pr y}    \int\limits_y^{\infty}  \tau f (\tau)
\int\limits_y^{\tau} \lr{t-y}^{\nu-\frac{1}{2}}
\lr{\tau^2-t^2}^{-\nu-\frac{1}{2}} dt d \tau. \label{1.1.22}
\end{multline}
Для вычисления внутреннего интеграла введем новую переменную по
формуле $z = \dfrac{t-y}{\tau-y}.$   Тогда   получим
$$
 \int\limits_y^{\tau} \lr{t-y}^{\nu-\frac{1}{2}}
 \lr{\tau^2-t^2}^{-\nu-\frac{1}{2}} dt  = (\tau+
 y)^{\nu-\frac{1}{2}} \int\limits_0^1 z^{-\nu-\frac{1}{2}}
 (1-z)^{\nu-\frac{1}{2}} \lr{1-\frac{y -
 \tau}{y+\tau}z}^{\nu-\frac{1}{2}} dz =
$$
$$
=\frac{\Gamma
 \lr{\frac{1}{2}-\nu} \, \Gamma \lr{\frac{1}{2}+\nu} }
 {(\tau+y)^{\nu-\frac{1}{2}}}\, {_2F_1} \lr{\frac{1}{2}-\nu;
 \frac{1}{2}-\nu; 1; \frac{y - \tau}{y+\tau}}.
$$
Гипергеометрическая функция в последней формуле сводится к функции
Лежандра по формуле~\cite[с.~128]{BE1}
$$
{_2F_1} \lr{\frac{1}{2}-\nu; \frac{1}{2}-\nu; 1; \frac{\zeta - 1}{\zeta+1}} = 2^{\nu -\frac{1}{2}} (\zeta+1)^{\frac{1}{2}-\nu} P^0_{\nu-\frac{1}{2}} \lr{\zeta}.
$$
Учитывая это, из~\eqref{1.1.22} находим
\begin{equation}
S_{\nu} f (y) = \frac {-2^{\nu+\frac{1}{2}} \, \Gamma \lr{\nu+1}}{ \sqrt{\pi} }  \frac{\pr}{\pr y} \int\limits_y^{\infty} \tau^{\nu+\frac{1}{2}} f (\tau)  P^0_{\nu-\frac{1}{2}} \lr{\frac{y}{\tau}} \, d \tau.
\label{1.1.23}
\end{equation}
В силу принципа аналитического продолжения полученное
представление оператора $S_{\nu}$ справедливо при $\Re \nu \geq
0.$ После   дифференцирования интеграла по параметру получаем
$$
S_{\nu} f (y) = \frac {2^{\nu+\frac{1}{2}} \, \Gamma \lr{\nu+1}}{ \sqrt{\pi} }\, y^{\nu+\frac{1}{2}} \lr{ f(y)-\lr{\nu^2 - \frac{1}{4}} \int\limits_1^{\infty} t^{\nu+\frac{1}{2}} f (t y)  \frac{P_{\nu-\frac{1}{2}}^{-1} \lr{\frac{1}{t}}}{\sqrt{t^2 -1}}\, dt}.
$$
Отсюда видно, что оператор $S_{\nu},$ так же как и $P_{\nu},$ не
изменяет гладкость функций при положительных аргументах.

Кроме введённых, мы будем использовать и некоторые другие
операторы преобразования, при построении которых будет использован
другой класс дробных интегралов. Свойства таких интегралов будут
изучены в следующем пункте.

\subsection{Дробные интегралы типа Римана---Лиувилля}\label{sec4.2}

На функциях класса $\mathring{C}_{\{0\}}^{\infty} \lr{E_{+}^{1}}$
определим оператор $I^{\mu}_e$ при $\Re \mu > 0$ по формуле
\begin{equation}
    I^{\mu}_e f(y) = \frac{1}{ \Gamma \lr{\mu}} \int\limits_y^{\infty} \lr{t-y}^{\mu-1}  e^{y-t} f(t) \, dt, \quad   y>0.
    \label{1.2.1}
\end{equation}
Если $\Re \mu > -M,$ где $M$ "--- целое неотрицательное число,
то мы полагаем
\begin{equation}
    I^{\mu}_e f(y) = \frac{(-1)^M}{ \Gamma \lr{\mu+M}}\,  e^{y} \int\limits_y^{\infty} \lr{t-y}^{\mu+M-1} \frac{\pr^M}{\pr t^M} \lr{e^{-t} f(t)} \, dt.
    \label{1.2.2}
\end{equation}
Если через $\mathcal{E}$ обозначим оператор умножения на функцию
$e^y,$ а через $\mathcal{E}^{-1}$ обратный к нему, то справедлива
формула
\begin{equation}
I^{\mu}_e =  \mathcal{E} I^{\mu}  \mathcal{E}^{-1},
\label{1.2.3}
\end{equation}
которая связывает операторы $I^{\mu}_e$ и $I^{\mu}$ при всех
комплексных $\mu.$ Отсюда и из свойства лиувиллевских операторов
вытекает групповое свойство операторов $I^{\mu}_e$:
\begin{equation}
I^{\mu}_e  I^{\nu}_e = I^{\nu}_e I^{\mu}_e =  I^{\mu+\nu}_e.
\label{1.2.4}
\end{equation}
Следовательно, любой оператор $I^{\mu}_e$ отображает пространство
$\mathring{C}_{\{0\}}^{\infty} \lr{E_{+}^{1}}$ на себя, и обратным
к нему будет оператор $I^{-\mu}_e.$

Операторы $I^{\mu}_e$ и $I^{\nu}$ коммутируют. Этот факт сразу
вытекает из их свёрточной природы. В частности, $I^{\mu}_e$
коммутирует с оператором дифференцирования.

Найдем явное выражение операторов вида $I^{\nu} I^{\mu}_e.$ Пусть
сначала $\Re \nu > 0$ и $\Re \mu > 0.$ Тогда для функции  $f
\in \mathring{C}_{\{0\}}^{\infty} \lr{E_{+}^{1}}$   имеем
\begin{multline}
I^{\nu} I^{\mu}_e f(y) = \frac{1} {\Gamma \lr{\nu}\, \Gamma \lr{\mu}} \int\limits_y^{\infty} \lr{t-y}^{\nu-1} \int\limits_t^{\infty} (\tau-t)^{\mu}  e^{t -\tau}  f(\tau) \, d \tau dt=  \\
 =\frac{1} {\Gamma \lr{\nu}\, \Gamma \lr{\mu}} \int\limits_y^{\infty}   f(\tau) \lr{\tau-y}^{\nu+\mu-1} \int\limits_0^1 z^{\mu-1} (1-z)^{\nu-1}  e^{(y -\tau)z}  \, d z d \tau. \label{1.2.5}
\end{multline}
Последний внутренний интеграл выражается через одну из специальных функций, а именно через вырожденную гипергеометрическую функцию $\Phi \lr{a, c; \zeta}$  (другое стандартное обозначение ${_1F_1} \lr{a, c; \zeta}$). Функция $\Phi$ при $\Re c > \Re a > 0$ может быть определена по формуле
\begin{equation}
\Phi \lr{a, c; \zeta} = \frac{\Gamma (c)}{\Gamma (a)  \, \Gamma (c-a)} \int\limits_0^1 z^{a-1} (1-z)^{c-a-1}  e^{\zeta z}  \, d z.
\label{1.2.6}
\end{equation}
Приведем некоторые известные её свойства (см.~\cite[гл.~6]{BE1}).
Функция $\dfrac{1}{\Gamma (c)} \Phi \lr{a, c; \zeta} $
аналитически продолжается до целой функции своих параметров $a,\
c$ и переменной $\zeta.$    При этом в  устранимых особых точках
$c= -m=0,-1,-2, \dots$ полагают
\begin{equation}
\frac{1}{\Gamma (c)} \Phi \lr{a, c; \zeta} = a (a+1) \dots (a+m) \, \frac{\zeta^{m+1}}{(m+1)!}  \, \Phi \lr{a+m+1, m+2; \zeta}.
\label{1.2.7}
\end{equation}
Справедливы следующие формулы:
\begin{equation}
\frac{\pr}{\pr \zeta} \lr{\zeta^{c} \, \Phi \lr{a, c+1; \zeta}} = c \, \zeta^{c-1} \Phi \lr{a, c; \zeta},
\label{1.2.8}
\end{equation}
\begin{equation}
\frac{\pr}{\pr \zeta}\,  \Phi \lr{a, c; \zeta}= \frac{a}{c}\, \Phi \lr{a+1, c+1; \zeta},
\label{1.2.9}
\end{equation}
\begin{equation}
 \Phi \lr{a, c; \zeta}= \frac{\Gamma(c)}{\Gamma(c-a)} (- \zeta)^{-a} \lr{1+O (|\zeta|^{-1})}, \  \zeta \to - \infty.
\label{1.2.10}
\end{equation}

Вернемся к формуле~\eqref{1.2.5}. Учитывая~\eqref{1.2.6}, находим
\begin{equation}
I^{\nu} I^{\mu}_e f(y)  = \frac{1}{\Gamma (\nu+\mu)} \int\limits_y^{\infty} f (\tau) (\tau -y)^{\nu+\mu-1} \Phi \lr{\mu, \mu+\nu; y - \tau} d \tau.
\label{1.2.11}
\end{equation}
Выведем аналогичное представление оператора $I^{\mu}_e I^{-\mu} =
I^{-\mu} I^{\mu}_e$ при комплексных $\mu.$ При $\Re \mu > 0$ и
$m=[\Re \mu]+1$ имеем для $f \in \mathring{C}_{\{0\}}^{\infty}
\lr{E_{+}^{1}}$ по предыдущей формуле
\begin{equation}
I^{-\mu} I^{\mu}_e f(y)  = (-1)^m \frac{1}{\Gamma (m)} \int\limits_y^{\infty} \frac{\pr^m f (\tau)}{\pr \tau^m} (\tau -y)^{m-1} \Phi \lr{\mu, m; y - \tau} d \tau.
\label{1.2.12}
\end{equation}
Интегрируя $m-1$ раз по частям в~\eqref{1.2.12} и используя
формулу~\eqref{1.2.8}, получаем
\begin{equation}
I^{-\mu} I^{\mu}_e f(y)  = - \int\limits_y^{\infty} \frac{\pr f (\tau)}{\pr \tau} \Phi \lr{\mu, 1; y - \tau} d \tau.
\label{1.2.13}
\end{equation}
Еще раз интегрируя по частям в~\eqref{1.2.13}, находим
\begin{equation}
I^{\mu} I^{\mu}_e f(y)  = f(y) \Phi \lr{\mu, 1; 0} + \int\limits_y^{\infty} f (\tau) \frac{\pr }{\pr \tau} \Phi \lr{\mu, 1; y - \tau} d \tau.
\label{1.2.14}
\end{equation}
Так как $\Phi \lr{\mu, 1; 0}=1,$ то с использованием~\eqref{1.2.9}
формула~\eqref{1.2.14} приводится к следующему окончательному
виду:
\begin{equation}
I^{-\mu} I^{\mu}_e f(y)  = f(y) - \mu \int\limits_y^{\infty} f (\tau) \Phi \lr{\mu, 2; y - \tau} d \tau.
\label{1.2.15}
\end{equation}
Эта формула, доказанная нами для $\Re \mu > 0,$ на самом деле
верна для функций $f \in \mathring{C}_{\{0\}}^{\infty}
\lr{E_{+}^{1}}$ при всех комплексных $\mu,$ поскольку здесь
возможно применение аналитического продолжения.

Изучим действие операторов $I^{\mu}_e$ в пространствах
С.\,Л.~Соболева на полупрямой, которые мы обозначим через $H^s
\lr{E_{+}^{1}},$    где целое $s \geq 0.$ Эти пространства
определяются как замыкание множества $\mathring{C}^{\infty}
\lr{\overline{E_{+}^{1}}}$ по норме
\begin{equation}
\| f \|_{H^s \lr{E_{+}^1}} = \lr{\sum\limits_{k=0}^{s} \|D^k f
\|_{L_{2} \lr{E_{+}^1}}^2}^{\frac{1}{2}}, \label{1.2.16}
\end{equation}
где
$$
\| f \|_{L_{2} \lr{E_{+}^1}} = \lr{\int\limits_0^{\infty} |f(y)|^2  \,dy}^{\frac{1}{2}}.
$$

\begin{lemma} \label{lem:1.2.1}
     Оператор $I^s_e$ расширяется до ограниченного оператора, изоморфно отображающего пространство $L_2 \lr{E_{+}^{1}}$ на $H^s \lr{E_{+}^{1}}.$
\end{lemma}

\begin{proof}
Для функции $f \in \mathring{C}^{\infty} (\ov{E_{+}^{1}})$ по
формуле Лейбница при $0 \leq k \leq s,$    получаем
\begin{multline}
D^k I^s_e f = D^k \mathcal{E} I^s \mathcal{E}^{-1} f = \sum\limits_{m=0}^k {k \choose m} \mathcal{E} D^m I^s \mathcal{E}^{-1} f=   \\
=\sum\limits_{m=0}^k (-1)^m {k \choose m} \mathcal{E} I^{s-m}
\mathcal{E}^{-1} f = \sum\limits_{m=0}^k (-1)^m {k \choose m}
I_e^{s-m} f. \label{1.2.17}
\end{multline}
При $\Re \mu > 0$ по обобщённому неравенству Минковского имеем
$$
| I_e^{\mu} f \|_{L_{2} \lr{E_{+}^1}} = \frac{1}{|\Gamma(\mu)|} \left\| \int\limits_0^{\infty} t^{\mu-1} e^{-t} f(y+t)  \,dt \right\|_{L_{2} \lr{E_{+}^1}} \leq
$$
$$
\leq \| f \|_{L_{2} \lr{E_{+}^1}} \, \frac{1}{| \Gamma(\mu)|} \, \int\limits_0^{\infty} t^{\Re \mu-1} e^{-t} dt =
 \| f \|_{L_{2} \lr{E_{+}^1}} \, \frac{\Gamma \lr{\Re \mu}}{| \Gamma(\mu)|}.
$$
Такая же оценка сохраняется и в тривиальном случае $\mu=0.$ Отсюда
и из~\eqref{1.2.17} получаем
\begin{equation}
\| I_e^s f \|_{H^s \lr{E_{+}^1}} \leq c \| f \|_{L_{2} \lr{E_{+}^1}}.
\label{1.2.18}
\end{equation}

Докажем противоположную оценку. По формуле Лейбница имеем для $f
\in \mathring{C}^{\infty} (\ov{E_{+}^{1}})$
$$
I^{-s}_e f = \mathcal{E} I^{-s} \mathcal{E}^{-1} f= (-1)^s
\mathcal{E} D^s \mathcal{E}^{-1} f =(-1)^s \sum\limits_{k=0}^s
(-1)^k {s \choose m} D^{s-k} f.
$$
Отсюда для $f \in \mathring{C}^{\infty} (\ov{E_{+}^{1}})$
получаем
\begin{equation}
 \| f \|_{L_{2} \lr{E_{+}^1}} = \left\|I_e^{-s} I_e^{s} f \right\|_{L_{2} \lr{E_{+}^1}} \leq c \| I_e^s f \|_{H^s \lr{E_{+}^1}}.
\label{1.2.19}
\end{equation}

Так как $I^s_e \lr{ \mathring{C}^{\infty} (\ov{E_{+}^{1}})}=
\mathring{C}^{\infty} (\ov{E_{+}^{1}}),$ то
оценки~\eqref{1.2.18},~\eqref{1.2.19} завершают доказательство
леммы.
\end{proof}

\begin{corollary}\label{cor:1.2.1}
    Норма $\| \tilde{f} \|_{H^s \lr{E_{+}^1}} =  \|I_e^{-s} f \|_{L_{2} \lr{E_{+}^1}}$
    на пространстве $H^s \lr{E_{+}^1}$ при $s \geq 0$ эквивалентна норме~\eqref{1.2.16}.
\end{corollary}

\begin{lemma}\label{lem:1.2.2}
     Пусть $s, s' \geq 0$ и $\mu$ "--- комплексное число. Тогда
     при $s-s'+ \Re \mu > 0$ оператор $I_e^{\mu}$ непрерывно
     отображает пространство $H^s \lr{E_{+}^1}$ в $H^{s'}
     \lr{E_{+}^1}.$   Если    $\mu$   вещественно,    то  это  же
     верно   и   при $s-s'+\mu=0.$
\end{lemma}

\begin{proof}
Используя норму $\widetilde{\|  \|}_{H^s}$ и групповое свойство
операторов $I_e^{\nu},$ получим
$$
\widetilde{| I_e^{\mu} f \|}_{H^{s'} \lr{E_{+}^1}} =  \| I_e^{\mu - s'} f \|_{L_{2} \lr{E_{+}^1}} = \| I_e^{s- s'+\mu } I_e^{-s} f \|_{L_{2} \lr{E_{+}^1}}.
$$
Осталось только заметить, что оператор $I_e^{\nu}$ отображает
непрерывно $L_{2} \lr{E_{+}^1}$ в себя при $\Re \nu > 0,$ а
также при $\nu=0.$ Лемма    доказана.
\end{proof}

Приведем еще один результат, который будет использован во второй главе.

\begin{lemma}\label{lem:1.2.3}
     Пусть функция $a \in C^{\infty} (\ov{E_{+}^1})$ и ограничена вместе со всеми производными. Тогда при $s \geq 0$ и любых комплексных $\mu$ оператор $I_e^{- \mu} a I_e^{\mu}$ непрерывно отображает пространство $H^{s} \lr{E_{+}^1}$ в себя. Справедлива оценка
\begin{equation}
\| I_e^{-\mu} a I_e^{\mu}  f \|_{H^{s} \lr{E_{+}^1}} \leq c  \|  f \|_{H^{s} \lr{E_{+}^1}}  \max\limits_{k \leq |[\Re \mu]|+s+2}\, \sup\limits_{y>0} \left| \frac{\pr^k a(y)}{\pr y^k} \right|,
\label{1.2.20}
\end{equation}
где постоянная $c>0$ не зависит от $f$ и $a.$
\end{lemma}

\begin{proof}
С помощью следствия~\ref{cor:1.2.1} общий случай $s \geq 0$
сводится к $s=0.$ Пусть сначала   $\Re \mu >0.$ По определению
операторов $I_e^{\mu},$ полагая  $m = [\Re \mu]+1,$ получим для
функции $f \in \mathring{C}^{\infty} (\ov{E_{+}^{1}})$
 следующую формулу:
$$
I_e^{-\mu} \lr{a I_e^{\mu} f} (y) = \frac{(-1)^m e^y}{\Gamma
(m-\mu) \, \Gamma(\mu)} \, D^m_y \int\limits_y^{\infty}
(t-y)^{m-\mu-1} a(t)\int\limits_t^{\infty} (\tau-t)^{\mu-1}
e^{-\tau} f (\tau) d \tau dt =
$$
$$
=\frac{(-1)^m e^y}{\Gamma (m-\mu) \,
\Gamma(\mu)} \, D^m_y \int\limits_y^{\infty} e^{- \tau} f (\tau)
 (\tau-t)^{m-1} \int\limits_0^1 z^{m-\mu-1} (1-z)^{\mu-1} a (y + z (\tau-y))\, dz d \tau.
$$
Введем обозначение $a^{(k)} (t) = D^k a(t).$    Тогда
$$
I_e^{-\mu} a I_e^{\mu}  f (y) = a(y) f(y) - \frac{e^y}{ \Gamma(m - \mu) \, \Gamma(\mu)} \int\limits_y^{\infty} e^{- \tau}
 f (\tau)  \sum\limits_{k=1}^m  {m \choose k}   \frac{(m-1)!}{(k-1)!} \times
$$
$$
\times (y - \tau)^{k-1}  \int\limits_0^1 z^{m-\mu-1} (1-z)^{k+\mu-1} a^{(k)} (y + z (\tau-y))\, dz d \tau.
$$
Следовательно,
$$
\left| I_e^{-\mu} a I_e^{\mu}  f (y)  \right| \leq c
\sum\limits_{k=0}^m \sup\limits_{t>0} |a^{(k)} (t)| I^k_e \lr{|f|}
(y).
$$
Применение леммы~\ref{lem:1.2.1} завершает доказательство при $\Re \mu>0.$ Если $\Re \mu < 0$ и $m = [- \Re \mu]+1,$ то
$$
I_e^{-\mu} a I_e^{\mu}  f (y) =  \frac{(-1)^m e^y}{ \Gamma(m +
\mu) \, \Gamma(-\mu)} \int\limits_y^{\infty}  (t - y)^{-\mu-1}
a(t) \int\limits_t^{\infty} (\tau -t)^{m+\mu-1}  D^m \lr{e^{- \tau} f (\tau)} \, d \tau d t=
$$
$$
=  \frac{(-1)^m
e^y}{ \Gamma(m + \mu) \, \Gamma(-\mu)} \int\limits_y^{\infty}  D^m
\lr{e^{- \tau} f (\tau)}  (\tau - y)^{m-1}   \int\limits_0^1 z^{-\mu-1} (1-z)^{m+\mu-1} a (y + z
(\tau-y))\, dz d \tau.
$$
После $m$ "--- кратного интегрирования по частям как и выше мы
убеждаемся в справедливости оценки~\eqref{1.2.20} и при $\Re \mu
< 0.$ Пусть, наконец, $\Re \mu = 0.$ Тогда
$$
I_e^{-\mu} a I_e^{\mu}  f (y)  =  \frac{e^y}{ \Gamma(1 + \mu) \, \Gamma(1-\mu)} \, D_y \int\limits_y^{\infty}   (t - y)^{- \mu} \int\limits_t^{\infty} (\tau -t)^{\mu} D \lr{e^{-\tau} f (\tau)} \, d \tau dt=
$$
$$
=  \frac{e^y}{ \Gamma(1 + \mu) \, \Gamma(1-\mu)} \, D_y \int\limits_y^{\infty} D_{\tau}  \lr{e^{- \tau} f (\tau)}  (\tau - y)  \int\limits_0^1 z^{-\mu} (1-z)^{\mu} a (y + z (\tau-y))\, dz d \tau.
$$
Остальные выкладки аналогичны предыдущим. Лемма доказана.
\end{proof}

\begin{corollary}\label{cor:1.2.2}
    В условиях леммы~\ref{lem:1.2.3} при $\Re (\nu+\mu) > 0$  имеет место оценка
    \begin{equation}
    \| I_e^{\nu} a I_e^{\mu}  f \|_{H^{s} \lr{E_{+}^1}} \leq c  \|  f \|_{H^{s} \lr{E_{+}^1}}
     \max\limits_{k \leq \min\limits \lr{|[\Re \nu]|, |[\Re \mu]|} +s+2} \sup\limits_{y>0} \left|D^k a(y) \right|.
    \label{1.2.21}
    \end{equation}
\end{corollary}

\begin{proof}
Достаточно применить формулу
$$
I_e^{\nu} a I_e^{\mu}  = I_e^{\nu+\mu} I_e^{-\mu} a I_e^{\mu} =I_e^{\nu} a I_e^{-\nu} I_e^{\nu+\mu}
$$
и воспользоваться ограниченностью оператора $I_e^{\lambda}$ с при
$\Re \lambda > 0.$
\end{proof}

\subsection{Связь с преобразованиями Фурье и
Ханкеля}\label{sec4.3}

Определим прямое и обратное преобразование Фурье по формулам
$$
F f (\eta) = \int\limits_{-\infty}^{\infty} f (y) e^{-i y \eta} dy, \  F^{-1} g(y) = \frac{1}{2 \pi} \int\limits_{-\infty}^{\infty} g (\eta) e^{i y \eta} d \eta.
$$
Будем использовать также косинус- и синус-преобразования Фурье:
$$
F_+ f (\eta) = \int\limits_{0}^{\infty} f (y) \cos (y \eta) \, dy, \  F^{-1}_{+} g(y) = \frac{2}{\pi} \int\limits_{0}^{\infty} g (\eta) \cos (y \eta) \, d \eta,
$$
$$
F_- f (\eta) = \int\limits_{0}^{\infty} f (y) \sin (y \eta) \, dy, \  F^{-1}_{-} g(y) = \frac{2}{\pi} \int\limits_{0}^{\infty} g (\eta) \sin (y \eta) \, d \eta.
$$
Прямое и обратное преобразования Ханкеля (Ханкеля) имеют вид
\begin{eqnarray}
F_{\nu} f (\eta) = \int\limits_{0}^{\infty} f (y) j_{\nu} (y \eta)  y^{2 \nu +1}\, dy,
\label{1.3.1}
\end{eqnarray}
$$
F^{-1}_{\nu} g(y) = \frac{1}{2^{2 \nu}\, \Gamma^2 (\nu+1)}   \int\limits_0^{\infty} g(\eta) j_{\nu} (y \eta) \eta^{2 \nu+1}\, d \eta,
$$
где $\nu \geq - \dfrac{1}{2}.$   Нормированная   функция Бесселя
$j_{\nu} (\lambda t)$   есть    решение  следующей задачи:
$$
B_{\nu} f = - \lambda^2 f, \  f(0)=1, \  f'(0)=0.
$$
Она связана с функцией Бесселя первого рода $J_{\nu} (t)$ формулой
\begin{eqnarray}
t^{\nu} j_{\nu} (t) = 2^{\nu} \Gamma (\nu+1) J_{\nu} (t).
\label{1.3.2}
\end{eqnarray}

Обозначим через $\mbox{П}$ какой-либо оператор продолжения
функций, действующий из $\mathring{C}^{\infty} (\ov{E_{+}^{1}})$ в
$\mathring{C}^{\infty} \lr{E^1}.$ Один из таких операторов
построен в работе~\cite{74}. Впрочем, можно предложить и следующий
элементарный оператор продолжения:
$$
\mbox{П} f(y) = \left\{
\begin{array}{ll}
 \chi (y) \int\limits_0^{\infty} f(-\lambda y) \psi (\lambda) \, d \lambda, & y<0, \\
 f(y), & y \geq 0,
\end{array}
\right.
$$
где $\psi (\lambda) = \dfrac{1}{\pi} \int\limits_0^{\infty} \sin
(x \sqrt{\lambda}) \psi_1 (x) \, dx,$ $\psi_1 (x) = \mbox{sh} (x)
\chi (x),$ $\chi$ "--- произвольная чётная функция из
$\mathring{C}^{\infty} \lr{E^1},$ равная единице в некоторой
окрестности нуля. Введённая функция $\psi$ ограничена, быстро
убывает на бесконечности, причём $\int\limits_0^{\infty} \lambda^n
\psi (\lambda) \, d \lambda = (-1)^n.$ Отсюда легко следует, что
оператор $\mbox{П}$ отображает пространство $\mathring{C}^{\infty}
(\ov{E_{+}^{1}})$ в $\mathring{C}^{\infty} \lr{E^1}.$

Пусть сначала $0 \leq \Re \nu < \dfrac{1}{2}.$ Тогда
из~\eqref{1.1.7}   для функции  $f \in \mathring{C}^{\infty}
(\ov{E_{+}^{1}})$ получаем
$$
P_{\nu}^{\frac{1}{2}-\nu} f(y) = \frac{\sqrt{\pi} }{\Gamma \lr{\nu+1}\, \Gamma \lr{\frac{1}{2}-\nu}} y^{- 2 \nu} \int\limits_1^{\infty} \lr{t^2-1}^{-\nu-\frac{1}{2}} f(yt) \, dt=
$$
$$
= \frac{y^{- 2 \nu}}{\sqrt{\pi}  \, \Gamma \lr{\nu+1}\, \Gamma \lr{\frac{1}{2}-\nu}} \int\limits_1^{\infty} \lr{t^2-1}^{-\nu-\frac{1}{2}} \int\limits_{-\infty}^{\infty} e^{i t \eta y} F \mbox{П} f (\eta) \, d \eta dt.
$$
Из наложенного на $\nu$ ограничения следует суммируемость подынтегральной функции. Применяя тогда теорему Фубини, получаем
$$
P_{\nu}^{\frac{1}{2}-\nu} f(y) =  \frac{y^{- 2 \nu}}{\sqrt{\pi}  \, \Gamma \lr{\nu+1}\, \Gamma \lr{\frac{1}{2}-\nu}} \int\limits_{-\infty}^{\infty} F \mbox{П} f (\eta) \int\limits_1^{\infty} \lr{t^2-1}^{-\nu-\frac{1}{2}}  e^{i t \eta y}  \, dt d \eta.
$$
Внутренний интеграл выражается через функцию Ханкеля первого рода
$H_{\nu}^{(1)}$ (функция Бесселя третьего рода) по
формуле~\cite[с.~95]{BE2}
$$
\int\limits_1^{\infty} \lr{t^2-1}^{-\nu-\frac{1}{2}}  e^{i t y \eta}  \, dt = i \sqrt{\pi} \, 2^{- \nu -1} (y \eta)^{\nu}\,
\Gamma \lr{\frac{1}{2}-\nu} H^{(1)}_{\nu} (y \eta).
$$
Подставляя эту формулу в предыдущую, находим
\begin{equation}
P_{\nu}^{\frac{1}{2}-\nu} f(y) =  \frac{i y^{-  \nu}}{2^{\nu+2}  \, \Gamma \lr{\nu+1}} \int\limits_{-\infty}^{\infty} H^{(1)}_{\nu} (y \eta) F \mbox{П} f (\eta) \, d \eta.
\label{1.3.3}
\end{equation}
Правая часть этой формулы суть аналитическая при $\Re \nu \geq 0$
функция. Этот факт легко усмотреть из представлений  функций
Ханкеля на той же странице цитированной книги~\cite{BE2}.
Следовательно, формула~\eqref{1.3.3} имеет место при $\Re \nu \geq
0.$

Из способа доказательства видно, что значение интеграла справа
в~\eqref{1.3.3} при $y>0$ не зависит от выбора оператора
продолжения. Этот факт может быть доказан и непосредственно,
исходя из самой формулы~\eqref{1.3.3} и используя одно обобщение
теоремы Пэли---Винера (см.~\cite{23}).

Сейчас мы дадим одно прямое следствие полученной формулы, которое будет использовано далее.

\begin{lemma}\label{lem:1.3.1}
Пусть функция $f \in \mathring{C}^{\infty} (\ov{E_{+}^{1}}).$
Тогда при $\Re \nu > 0$ справедливо соотношение
$$
\lim\limits_{y \to +0} y^{2 \nu} P_{\nu}^{\frac{1}{2}-\nu} f(y) = \frac{1}{2 \nu} f (0),
$$
а при $\nu=0$
$$
\lim\limits_{y \to +0} \frac{1}{\ln \frac{1}{y}} P_{0}^{\frac{1}{2}} f(y) = f (0).
$$
\end{lemma}

\begin{proof}
Доказательство легко следует из следующих (см.~\cite{BE2})
асимптотических формул для функций Ханкеля:
$$
z^{\nu}  H^{(1)}_{\nu}  (z) = - \frac{i 2^{\nu} \,
\Gamma(\nu)}{\pi} + o(1), \quad   \Re \nu > 0,\quad  z \to 0,
$$
$$
H^{(1)}_0  (z) = \frac{2 i}{\pi} \ln z + O(1), \quad   z \to 0,
$$
а также из ограниченности функций $H^{(1)}_{\nu}  (z)$ при $z \geq
\varepsilon > 0.$
\end{proof}

Рассмотрим случай вещественного $\nu \geq 0.$ Пусть функция $f \in
\mathring{C}^{\infty}_{-} (\ov{E_{+}^{1}}).$  Отсюда  следует, что
её продолжение по закону нечётности на всю прямую принадлежит
пространству $\mathring{C}^{\infty} \lr{E^1}.$ Выбирая в таком
случае в качестве оператора продолжения в~\eqref{1.3.3}
продолжение по закону нечётности, находим
$$
P_{\nu}^{\frac{1}{2}-\nu} f(y) =  \frac{y^{-  \nu}}{2^{\nu+1} \, \Gamma \lr{\nu+1}} \int\limits_{-\infty}^{\infty} H^{(1)}_{\nu} (y \eta) \eta^{\nu} F_{-} f (\eta) \, d \eta=
\frac{y^{-  \nu}}{2^{\nu+1} \, \Gamma \lr{\nu+1}}  \lr{\int\limits_{-\infty}^0+\int\limits_0^{\infty}} =
$$
$$
= \frac{y^{-  \nu}}{2^{\nu+1} \, \Gamma \lr{\nu+1}} \int\limits_0^{\infty} \lr{H^{(1)}_{\nu} (y \eta)- e^{i \pi \nu} H^{(1)}_{\nu} (-y \eta)} \eta^{\nu} F_{-} f (\eta) \, d \eta.
$$
Далее, поскольку (см.~\cite[с.~91]{BE2})
$$
- e^{i \pi \nu} H^{(1)}_{\nu} (-z) = H^{(2)}_{\nu} (z)
$$
 и (см.~\cite[с.~12]{BE2})
$$
J_{\nu} (z) = \frac{H^{(1)}_{\nu}(z)+H^{(2)}_{\nu}(z)}{2},
$$
то предыдущая формула принимает вид
$$
P_{\nu}^{\frac{1}{2}-\nu} f(y) =  \frac{1}{2^{\nu+1} \, \Gamma \lr{\nu+1}}\, y^{-  \nu} \int\limits_0^{\infty} J_{\nu} (y \eta)  \eta^{\nu} F_{-} f (\eta) \, d \eta.
$$
Заменяя здесь функцию $J_{\nu}$ на $j_{\nu}$ по
формуле~\eqref{1.3.2}, получаем
\begin{equation}
P_{\nu}^{\frac{1}{2}-\nu} f(y) =  \frac{1}{2^{2 \nu}  \, \Gamma^2 \lr{\nu+1}} \int\limits_{0}^{\infty} j_{\nu} (y \eta)  \eta^{2 \nu} F_{-} f (\eta) \, d \eta.
\label{1.3.4}
\end{equation}

\begin{lemma}\label{lem:1.3.2}
Пусть функция $f \in \mathring{C}^{\infty}_{-} (\ov{E_{+}^{1}})$ и
$\nu \geq 0.$ Тогда   имеет место представление
\begin{equation}
P_{\nu}^{\frac{1}{2}-\nu} f =  F_{\nu}^{-1} \lr{\frac{1}{\eta} F_{-} f}.
\label{1.3.5}
\end{equation}
Если $f \in \mathring{C}^{\infty}_{+} (\ov{E_{+}^{1}}),$ то
\begin{equation}
S_{\nu}^{\nu-\frac{1}{2}} f =  F_{\nu}^{-1} \lr{\eta F_{\nu} f}.
\label{1.3.6}
\end{equation}
Оператор $P_{\nu}^{\frac{1}{2}-\nu}$ взаимно однозначно отображает
пространство $\mathring{C}^{\infty}_{-} (\ov{E_{+}^{1}})$ на
$\mathring{C}^{\infty}_{+} (\ov{E_{+}^{1}}).$ Оператор
$S_{\nu}^{\nu-\frac{1}{2}}$ осуществляет обратное отображение.
\end{lemma}

\begin{proof}
Формула~\eqref{1.3.5} "--- это просто компактная запись
формулы~\eqref{1.3.4}. Формула~\eqref{1.3.6} доказывается
обращением предыдущей. Пусть функция $f \in
\mathring{C}^{\infty}_{-} (\ov{E_{+}^{1}}).$ Тогда $F_{-} f
(\eta)$ "--- нечётная быстроубывающая функция, а, следовательно,
$\dfrac{1}{\eta} F_{-} f (\eta)$ "--- чётная гладкая
быстроубывающая функция. Оператор $F_{\nu}$ является автоморфизмом
такого класса функций. Значит, функция  $F_{\nu}
\Big(\dfrac{1}{\eta} F_{-} f\Big) \in \mathring{C}^{\infty}_{+}
(\ov{E_{+}^{1}}),$ поскольку её финитность следует из
представления оператора по формуле~\eqref{1.1.7}. Точно так же
доказывается, что оператор $S_{\nu}^{\nu-\frac{1}{2}}$ отображает
пространство $\mathring{C}^{\infty}_{+} (\ov{E_{+}^{1}})$ в
$\mathring{C}^{\infty}_{-} (\ov{E_{+}^{1}}).$ Следовательно, эти
отображения сюръективны. Лемма доказана.
\end{proof}

Заметим, что лемму можно доказать иначе, опираясь на то, что
оператор $F_{\nu}$ отображает $\mathring{C}^{\infty}_{+}
(\ov{E_{+}^{1}})$ на множество чётных целых функций
экспоненциального типа~\cite{23}.

Отвлекаясь несколько в сторону от основного изложения, покажем,
как с помощью аналогов формул~\eqref{1.3.5} и~\eqref{1.3.6} можно
построить новый класс операторов преобразования~\cite{Kat1}. Пусть
$\nu \geq -\dfrac{1}{2}.$   Тогда   положим
\begin{equation}\label{preCM}
P_{\nu, \pm}^{(\varphi)}  =  F_{\nu}^{-1} \lr{ \varphi \lr{\eta} F_{\pm} }, \   S_{\nu, \pm}^{(\varphi)}  =  F_{\pm}^{-1} \lr{ \frac{1}{\varphi \lr{\eta}} F_{\nu} },
\end{equation}
где $\varphi (\eta)$ "---    некоторая   функция.    Нетрудно
заметить, что на подходящих областях определения операторы
$P_{\nu, \pm}^{(\varphi)}$ и $S_{\nu, \pm}^{(\varphi)}$
действительно будут операторами преобразования. Отметим наиболее
важные частные случаи. Если $\varphi (\eta) = \eta^{2 \nu+1},$ то
операторы $P_{\nu, +}^{(\varphi)}$ и $S_{\nu, +}^{(\varphi)}$
совпадают с классическими операторами Пуассона и Сонина. Если
$\varphi (\eta) \equiv 1,$ то $P_{\nu, +}^{(\varphi)}=
P_{\nu}^{-\frac{1}{2}-\nu}$ и  $S_{\nu, +}^{(\varphi)}=
S_{\nu}^{\nu+\frac{1}{2}}.$ Если же $\varphi (\eta) =
\eta^{\nu+\frac{1}{2}},$ то
\begin{multline}
P_{\nu, +}^{(\varphi)} f(y) = F^{-1}_{\nu}
\lr{\eta^{\nu+\frac{1}{2}} F_+ f} (y) =
\\
=\frac{-1}{2^{\nu+2}\,\Gamma (\nu+1)}
\left( \frac{\Gamma \lr{\nu+\frac{1}{2}}\,\cos (\pi \nu)}{2^{\nu}\, \Gamma (\nu+1)}  \int\limits_y^{\infty} {_2F_1} \lr{\frac{\nu}{2}+\frac{3}{4}, \frac{\nu}{2}+\frac{1}{4}; \nu+1; \frac{y^2}{t^2}} D f(t) \, dt + {} \right.\\
\left. {}+ \frac{\sqrt{2} \, \Gamma
\lr{\frac{\nu}{2}+\frac{3}{4}}}{\Gamma
\lr{\frac{\nu}{2}+\frac{1}{4}}} \frac{1}{y} \int\limits_0^y
{_2F_1} \lr{\frac{\nu}{2}+\frac{3}{4}, \frac{3}{4}-\frac{\nu}{2};
\frac{3}{2}; \frac{t^2}{y^2}} D f(t) \, dt \right),
 \label{1.3.8}
\end{multline}
\begin{multline*}
S_{\nu, +}^{(\varphi)} f(y) = F^{-1}_{+}
\lr{\eta^{-\nu-\frac{1}{2}} F_{\nu} f} (y) =
\\
=\frac{2^{\nu-1}\,\Gamma (\nu+1)}{\pi} \frac{\pr}{\pr y} \left(
\frac{\Gamma \lr{\nu+\frac{1}{2}}\,\cos (\pi \nu)}{2^{\nu}\,
\Gamma (\nu+1)}  y^{-\nu-\frac{1}{2}} \int\limits_0^y {_2F_1}
\lr{\frac{\nu}{2}+\frac{3}{4}, \frac{\nu}{2}+\frac{1}{4}; \nu+1;
\frac{t^2}{y^2}} t^{2 \nu +1} f(t) \, dt+{}\right.
\\
\left. {}+ \frac{\sqrt{2} \, \Gamma
\lr{\frac{\nu}{2}+\frac{3}{4}}}{\Gamma
\lr{\frac{\nu}{2}+\frac{1}{4}}} y \int\limits_y^{\infty} {_2F_1}
\lr{\frac{\nu}{2}+\frac{3}{4}, \frac{3}{4}-\frac{\nu}{2};
\frac{3}{2}; \frac{y^2}{t^2}} t^{\nu - \frac{1}{2}} f(t) \, dt
\right).
\end{multline*}
Операторы~\eqref{1.3.8}  названы в работе~\cite{Kat1}
изометрическими, поскольку они изометрически отображают
пространства $L_2$ в $L_{2, \nu}.$ В той же работе рассмотрены
некоторые их приложения к псевдодифференциальным операторам и к
спектральной теории. Отметим, что как показано в главе~\ref{ch3},
полученные представления можно упростить, выразив ядра операторов
преобразования через функции Лежандра. Подправленные на степенной
множитель, такие операторы можно свести к унитарным в одном
лебеговом пространстве на полуоси. Они по терминологии
главы~\ref{ch3} (см. ниже) являются комбинациями операторов
преобразования Бушмана---Эрдейи первого и второго родов. Свойство
унитарности делает их особенно полезными в приложениях, например,
они связывают решения возмущённых с операторами Бесселя и
невозмущённых со вторыми производными дифференциальных уравнений с
сохранением нормы. Этот важный класс операторов, введённых
В.\,В.~Катраховым, в~\cite{SitDis} предложено назвать операторами
Сонина---Катрахова и Пуассона---Катрахова (см. главу~\ref{ch3}).
Также отметим, что именно приведённая выше конструкция
В.\,В.~Катрахова~\eqref{preCM} послужила примером для разработки в
дальнейшем С.\,М.~Ситником общего композиционного метода
построения операторов преобразования различных классов, см.
главу~\ref{ch6} этой монографии.

Введем теперь операторы преобразования $P_{\nu, e}$ и $S_{\nu, e}$ по формулам
\begin{equation}
P_{\nu,e} = P_{\nu}^{\frac{1}{2}-\nu} I_e^{\nu - \frac{1}{2}}, \  S_{\nu,e} = I_e^{\frac{1}{2}-\nu} S_{\nu}^{\nu-\frac{1}{2}}.
\label{1.3.9}
\end{equation}
Их можно также выразить в виде
\begin{equation}
    P_{\nu,e} = P_{\nu} \mathcal{J}_{\nu - \frac{1}{2}, e}, \  S_{\nu,e} = \mathcal{J}_{\frac{1}{2}-\nu, e} S_{\nu},
    \label{1.3.10}
\end{equation}
где $\mathcal{J}_{\mu, e} = I_e^{\mu} I^{-\mu}.$ Использованные
здесь операторы $P_{\nu}^{\frac{1}{2}-\nu},$
$S_{\nu}^{\nu-\frac{1}{2}},$ $P_{\nu}$  определены в
пункте~\ref{sec4.1}, а операторы $I_e^{\mu}$ "--- в
пункте~\ref{sec4.2}.

Из результатов указанных разделов сразу же следует, что операторы $P_{\nu,e}$ и $S_{\nu,e}$ взаимно однозначно  отображают  пространство $\mathring{C}^{\infty}_{\{0\}} \lr{{E_{+}^{1}}}$ на себя, являются взаимно обратными и для них справедливы формулы
\begin{equation}
    B_{\nu} P_{\nu,e} = P_{\nu,e} D^2, \  D^2 S_{\nu,e} = S_{\nu, e} B_{\nu}.
    \label{1.3.11}
\end{equation}
Таким образом,  $P_{\nu,e}$ и $S_{\nu,e}$   также   являются операторами    преобразования.

Выясним связь   оператора $P_{\nu,e}$ с преобразованием Фурье.
Пусть функция $f \in \mathring{C}^{\infty} \lr{E^{1}}$ и $\mu$
"--- комплексное число. Тогда при $\Re \mu > 0$ по теореме
Фубини получаем
$$
F I_e^{\mu} f (\eta) = \frac{1}{\Gamma (\mu)} \int\limits_{-\infty}^{\infty} e^{- i y \eta} \int\limits_y^{\infty} (t-y)^{\mu -1 } e^{y-t} f(t) \, dt dy =
 \frac{1}{\Gamma (\mu)} \int\limits_{-\infty}^{\infty} f(t)  \int\limits_{-\infty}^t  e^{- i y \eta} (t-y)^{\mu -1 } e^{y-t} \, dy dt.
$$
Так как
$$
\int\limits_{-\infty}^t  e^{- i y \eta} (t-y)^{\mu -1 } e^{y-t} \, dy =  e^{- i t \eta} (1 - i \eta)^{- \mu} \Gamma(\mu),
$$
то
\begin{equation}
F I_e^{\mu} f (\eta) = (1 - i \eta)^{- \mu} F f (\eta).
\label{1.3.12}
\end{equation}
Пусть $\Re \mu \leq 0$    и $m = [-\Re \mu]+1.$ Тогда тем   же
методом получаем
$$
F I_e^{\mu} f (\eta) = \frac{(-1)^m}{\Gamma(m+\mu)} \int\limits_{-\infty}^{\infty} e^{- i t \eta} e^y  \int\limits_y^{\infty} (t-y)^{m+\mu -1 } D^m_t \lr{e^{-t} f(t)} \, dt dy =
$$
$$
= \frac{(-1)^m}{\Gamma(m+\mu)} \int\limits_{-\infty}^{\infty} D^m_t \lr{e^{-t} f(t)} \int\limits_{-\infty}^t  e^{- i y \eta} (t-y)^{m+\mu -1 } e^{y} \, dy dt =
$$
$$
= \frac{(-1)^m}{\Gamma(m+\mu)} \int\limits_{-\infty}^{\infty} D^m_t \lr{e^{-t} f(t)}  e^{- i t \eta + t}  (1 - i \eta)^{- (m+\mu)}  \Gamma(m+\mu) \, dt.
$$

Интегрирование по частям в последнем интеграле снова приводит к
формуле~\eqref{1.3.12}, которая, следовательно, доказана при всех
комплексных $\mu.$ Заметим, что здесь и ниже выбирается следующая
ветвь степенной функции:
$$
(1 - i \eta)^{- \mu} = \exp \lr{\frac{\mu}{2} \ln (1+ \eta^2) - i
\mu \, \mbox{arctg}\, \eta }.
$$

Учитывая формулы~\eqref{1.3.12},~\eqref{1.3.13}, при $\Re \mu
\geq 0$ получаем следующее представление для оператора    $P_{\nu,
e}$:
\begin{equation}
P_{\nu, e} f(y) = \frac{i y^{- \nu}}{ 2^{\nu+2}\, \Gamma (\nu+1)} \int\limits_{-\infty}^{\infty} H^1_{\nu} (y \eta) \eta^{\nu} (1-i \eta)^{\frac{1}{2}-\nu} F \mbox{П} f(\eta) \, d \eta,
 \label{1.3.13}
\end{equation}
где $f \in \mathring{C}^{\infty}_{\{0\}} \lr{{E_{+}^{1}}}.$

\subsection{Операторы преобразования и функциональные пространства (одномерная
теория)}\label{sec4.4}

В этом пункте рассматриваются оценки норм введённых выше
операторов преобразования. Эти результаты изложены по диссертации
В.\,В.~Катрахова~\cite{KatDis}. Применяемый метод, использующий по
существу преобразование Меллина и мультипликаторы преобразования
Меллина для оценок норм операторов преобразования, был предложен
С.\,М.~Ситником, ему же принадлежат основные результаты этого
пункта. Оригинальный подход С.\,М.~Ситника, основанный на теории
операторов преобразования Бушмана---Эрдейи и явном использовании
мультипликаторов Меллина, изложен в главе~\ref{ch3} данной
монографии, далее приводится другой вариант изложения этих
результатов из диссертации В.\,В.~Катрахова~\cite{KatDis}. На
самом деле разделение указанных результатов не имеет смысла, так
как они многократно и подробно обсуждались обоими авторами в
течение ряда лет в процессе поиска нужных подходов и техники
доказательств. Приложения указанных оценок норм операторов
преобразования  к весовым краевым задачам (глава~\ref{ch4}) и
краевым задачам с изолированными особенностями (глава~\ref{ch5})
полностью получены В.\,В.~Катраховым.

Обозначим через $L_{2, \nu} \lr{E_{+}^1},$ $\nu \geq -
\dfrac{1}{2},$ гильбертово    пространство функций $f(y),$ $y>0,$
для которых конечна норма
\begin{equation}
\| f \|_{L_{2, \nu} \lr{E^1_+}} = \lr{\int\limits_0^{\infty}
|f(y)|^2 y^{2 \nu +1} \, dy}^{\frac{1}{2}}. \label{1.4.1b}
\end{equation}
Хорошо известно, что преобразование Ханкеля $F_{\nu}$ отображает $L_{2, \nu} \lr{E_{+}^1}$ на себя, и справедливо равенство Парсеваля
\begin{equation}
\|F_{\nu} f \|_{L_{2, \nu} \lr{E^1_+}} = 2^{\nu} \, \Gamma (\nu
+1) \| f \|_{L_{2, \nu} \lr{E^1_+}}. \label{1.4.2b}
\end{equation}
Функциональное пространство $H^s_{\nu,+} \lr{E_{+}^1},$ $s \geq
0,$ $\nu \geq - \dfrac{1}{2},$   введённое в работе~\cite{Kip1},
определяется как замыкание по норме
\begin{equation}
    \| f \|_{H_{ \nu, +}^s \lr{E^1_+}} = \frac{1}{2^{\nu}\, \Gamma (\nu +1)} \| (1+\eta^2)^{\frac{s}{2}} F_{\nu} f \|_{L_{2, \nu} \lr{E^1_+}}
    \label{1.4.3}
\end{equation}
множества функций $\mathring{C}^{\infty}_{+} (\ov{E_{+}^{1}}).$
Предположение о чётности    здесь существенно, поскольку в
противном случае норма~\eqref{1.4.3} может быть равной
бесконечности.

Если замкнуть по той же норме~\eqref{1.4.3} множество
$\mathring{C}^{\infty}_{+} [0, R),$ $0<R<\infty,$ то получим
пространство $\mathring{H}^{s}_{\nu, +} (0, R),$  которое
непрерывно вложено в пространство $H^{s}_{\nu, +} \lr{E_{+}^1}.$ В
пространстве $\mathring{H}^{s}_{\nu, +} (0, R)$ выражение
\begin{equation}
\| f \|_{\mathring{H}_{ \nu, +}^s (0, R)} = \frac{1}{2^{\nu} \, \Gamma (\nu +1)} \| \eta^{s} F_{\nu} f \|_{L_{2, \nu} \lr{E^1_+}}.
\label{1.4.4}
\end{equation}
является нормой, эквивалентной норме~\eqref{1.4.3}. Из равенства
Парсеваля~\eqref{1.4.2b} легко следует, что при чётных $s \geq 0$
$$
\| f \|_{\mathring{H}_{ \nu, +}^s (0, R)} =  \|B_{\nu}^{\frac{s}{2}}  f \|_{L_{2, \nu} (0, R)}.
$$

Определим пространство С.\,Л.~Соболева $\mathring{H}^s (0, R),$ $s
\geq 0$ $0<R<\infty,$ как замыкание множества
$\mathring{C}^{\infty} [0, R)$ по норме
$$
\| f \|_{\mathring{H}^s (0, R)} =  \|D^s  f \|_{L_{2, \nu} (0, R)}.
$$

\begin{lemma}\label{lem:1.4.1}
При $\nu \geq 0$ оператор $S_{\nu}$ расширяется по непрерывности
до ограниченного оператора, отображающего пространство
$\mathring{H}_{\nu, +}^s (0, R)$ в  $\mathring{H}^s (0, R),$
причём справедлива оценка
\begin{multline}
\frac{2^{\nu+1}\, \Gamma (\nu +1)}{\sqrt{\pi}} \min\limits \lr{
\frac{1}{\sqrt{2}}, | \cos \frac{\pi (\nu+s)}{2}| }
\| f \|_{\mathring{H}_{ \nu, +}^s (0, R)} \leq   \\
\leq \|S_{\nu} f \|_{\mathring{H}^s (0, R)} \leq \frac{2^{\nu+1}\,
\Gamma (\nu +1)}{\sqrt{\pi}}  \min\limits \lr{ \frac{1}{\sqrt{2}},
| \cos \frac{\pi (\nu+s)}{2}| } \| f \|_{\mathring{H}_{ \nu, +}^s
(0, R)}. \label{1.4.5}
\end{multline}
\end{lemma}

\begin{proof}
Рассмотрим лишь случай $s=0,$ рассуждения при $s>0$ вполне
аналогичны. Докажем сначала справедливость формулы
\begin{equation}
S_{\nu} f(y) = \frac{2}{\pi} \int\limits_0^{\infty} \cos \lr{y
\eta - \frac{\pi \nu}{2} - \frac{\pi}{4}} \eta^{\nu+ \frac{1}{2}}
F_{\nu} f (\eta) \, d \eta. \label{1.4.6}
\end{equation}
для $\nu \geq 0$ и $f \in \mathring{C}^{\infty}_{+}
\lr{E_{+}^{1}}.$ На этом множестве ранее было получено следующее
представление оператора $S_{\nu}$:
\begin{equation}
S_{\nu} f = I^{\frac{1}{2}-\nu} F^{-1}_{-} \eta F_{\nu} f.
\label{1.4.7}
\end{equation}
Отсюда для полуцелых $\nu$ сразу же вытекает
справедливость~\eqref{1.4.6}, поскольку оператор
$I^{\frac{1}{2}-\nu}$ в этом случае есть просто оператор
дифференцирования $(-1)^{\nu-\frac{1}{2}} D^{\nu-\frac{1}{2}}.$
Пусть теперь $n \geq 0$ "--- чётное число и
$n-\dfrac{1}{2}<\nu<n+\dfrac{1}{2}.$ Так как функция $\eta^{1+n}
F_{\nu} f (\eta)$ суть гладкая, нечётная и быстроубывающая на
бесконечности, то по теореме Фубини имеем
\begin{multline}
S_{\nu} f(y) = I^{\frac{1}{2}-\nu} F^{-1}_{-} \eta F_{\nu} f (y) =
(-1)^n I^{\frac{1}{2}-\nu+n} D^n F^{-1}_{-} \eta F_{\nu} f (y) =
\\
=(-1)^{\frac{n}{2}} I^{\frac{1}{2}-\nu+n} F^{-1}_{-} \eta^{1+n}
F_{\nu} f (y) =
\\
= \frac{2 (-1)^{\frac{n}{2}}}{\pi \, \Gamma
\lr{\frac{1}{2}-\nu+n}} \int\limits_y^{\infty} (t-y)^{n - \nu
-\frac{1}{2}} \int\limits_0^{\infty} \sin (t \eta) \eta^{1+n}
F_{\nu} f (\eta) \, d \eta dt =
\\
= \frac{2 (-1)^{\frac{n}{2}}}{\pi \, \Gamma
\lr{\frac{1}{2}-\nu+n}} \lim\limits_{A \to \infty} \int\limits_y^A
(t-y)^{n - \nu -\frac{1}{2}} \int\limits_0^{\infty} \sin (t \eta)
\eta^{1+n} F_{\nu} f (\eta) \, d \eta dt =
\\
= \frac{2 (-1)^{\frac{n}{2}}}{\pi \, \Gamma
\lr{\frac{1}{2}-\nu+n}} \lim\limits_{A \to \infty}
\int\limits_0^{\infty} \eta^{1+n} f (\eta) \int\limits_y^A
(t-y)^{n - \nu -\frac{1}{2}} \sin (t \eta)  \, dt d \eta  =
\\
= \frac{2 (-1)^{\frac{n}{2}}}{\pi \, \Gamma
\lr{\frac{1}{2}-\nu+n}} \lim\limits_{A \to \infty}
\int\limits_0^{\infty} \eta^{1+n} f (\eta)
\lr{\int\limits_y^{\infty}-\int\limits_A^{\infty}} \, d \eta,
\label{1.4.8}
\end{multline}
где в последнем выражении два внутренних интеграла понимаются как
несобственные. Первый из них "--- табличный интеграл
\begin{equation}
\int\limits_y^{\infty} (t-y)^{n - \nu -\frac{1}{2}} \sin (t \eta)
\, dt = \eta^{ \nu-n -\frac{1}{2}} \Gamma \lr{n-\nu+\frac{1}{2}}
\cos \lr{y \eta - \frac{\pi (\nu-n)}{2} - \frac{\pi}{4}},
\label{1.4.9}
\end{equation}
а второй с помощью интегрирования по частям оценивается следующим
образом:
\begin{multline}
\left| \int\limits_A^{\infty} (t-y)^{n - \nu -\frac{1}{2}} \sin (t
\eta) \, dt \right| = \frac{1}{\eta} \left| \int\limits_A^{\infty}
(t-y)^{n - \nu -\frac{1}{2}} \frac{\pr}{\pr t}\cos (t \eta) \, dt
\right| \leq
\\
\leq  \frac{1}{\eta} \lr{(A-y)^{n - \nu -\frac{1}{2}} |\cos (A
\eta)| + \lr{\frac{1}{2}+\nu-n} \left| \int\limits_A^{\infty}
(t-y)^{n - \nu -\frac{3}{2}} \cos (t \eta) \, dt \right|} \leq
\\
\leq \frac{1}{\eta} \lr{(A-y)^{n - \nu -\frac{1}{2}} +
\lr{\frac{1}{2}+\nu-n}  \int\limits_A^{\infty} (t-y)^{n - \nu
-\frac{3}{2}} \, dt } = \frac{2}{\eta} (A-y)^{n - \nu
-\frac{1}{2}}. \label{1.4.10}
\end{multline}
Подставляя~\eqref{1.4.9},~\eqref{1.4.10} в~\eqref{1.4.8} мы
получаем формулу~\eqref{1.4.6}. Пусть теперь $n>0$ "--- нечётное и
$n-\dfrac{1}{2}<\nu<n+\dfrac{1}{2},$  тогда аналогично предыдущему
получаем
$$
S_{\nu} f(y) =  (-1)^n I^{\frac{1}{2}-\nu+n} D^n F^{-1}_{-} \eta
F_{\nu} f (y) =  (-1)^{\frac{n+1}{2}} I^{\frac{1}{2}-\nu+n}
F^{-1}_{+} \eta^{1+n} F_{\nu} f (y) =
$$
$$
= \frac{2 (-1)^{\frac{n+1}{2}}}{\pi \, \Gamma
\lr{\frac{1}{2}-\nu+n}} \int\limits_y^{\infty} (t-y)^{n - \nu
-\frac{1}{2}} \int\limits_0^{\infty} \cos (t \eta) \eta^{1+n}
F_{\nu} f (\eta) \, d \eta dt =
$$
$$
= \frac{2 (-1)^{\frac{n+1}{2}}}{\pi \, \Gamma
\lr{\frac{1}{2}-\nu+n}}  \int\limits_0^{\infty} \eta^{1+n} F_{\nu}
f (\eta) \int\limits_y^{\infty} (t-y)^{n - \nu -\frac{1}{2}} \cos
(t \eta)  \, dt d \eta.
$$
Так как
\begin{equation}
\int\limits_y^{\infty} (t-y)^{n - \nu -\frac{1}{2}} \cos (t \eta)
\, dt = - \eta^{ \nu-n -\frac{1}{2}} \Gamma \lr{\frac{1}{2}-\nu+n}
\sin \lr{y \eta - \frac{\pi (\nu-n)}{2} - \frac{\pi}{4}},
\label{1.4.11}
\end{equation}
то и в этом случае установлена справедливость~\eqref{1.4.6}.

Рассмотрим встретившийся в~\eqref{1.4.6} оператор $A_{\nu}$ вида
\begin{equation}
A_{\nu} f(y) = \int\limits_0^{\infty}  \cos \lr{y \eta - \frac{\pi
\nu}{2} - \frac{\pi}{4}} f (\eta)\, d \eta, \label{1.4.12}
\end{equation}
который является линейной комбинацией тригонометрических
преобразований Фурье и потому ограничен в $L_{2} \lr{E_{+}^1}.$
Введем несколько видоизмененный оператор Меллина $M$ по формуле
\begin{equation}
    M g (p) = \frac{1}{\sqrt{2 \pi}} \int\limits_0^{\infty} y^{ i p - \frac{1}{2}}  g(y) \, d y, \  p \in E^1.
    \label{1.4.13}
\end{equation}
Нетрудно заметить, что оператор $M$ изометрично отображает
пространство в $L_{2} \lr{E_{+}^1}$ на $L_{2} \lr{E_{+}^1}.$

Пусть для простоты функция  $f \in \mathring{C}^{\infty}
\lr{E_{+}^{1}},$   тогда с использованием формул~\eqref{1.4.9}
и~\eqref{1.4.11}, действуя как и при выводе формулы~\eqref{1.4.6},
получаем
$$
M A_{\nu} f (p) =  \frac{1}{\sqrt{2 \pi}} \int\limits_0^{\infty}
y^{ i p - \frac{1}{2}}  \int\limits_0^{\infty}  \cos \lr{y \eta -
\frac{\pi \nu}{2} - \frac{\pi}{4}} f (\eta)\, d \eta dy=
$$
$$
=  \frac{1}{\sqrt{2 \pi}} \int\limits_0^{\infty}  f (\eta)
\int\limits_0^{\infty} y^{ i p - \frac{1}{2}}   \cos \lr{y \eta -
\frac{\pi \nu}{2} - \frac{\pi}{4}} \, d \eta dy=
$$
$$
=  \frac{-1}{\sqrt{2 \pi}} \int\limits_0^{\infty}  f (\eta)
\eta^{ i p - \frac{1}{2}}   \Gamma \lr{i p +\frac{1}{2}} \sin
\frac{\pi (i p - \nu -1)}{2} \, d \eta=
$$
$$
=  -   \Gamma \lr{i p +\frac{1}{2}} \sin  \frac{\pi (i p - \nu
-1)}{2}  M f (-p) \equiv a_{\nu} (p) M f (-p).
$$
Отсюда и из свойства изометричности оператора $M$ следует
двусторонняя оценка для оператора~\eqref{1.4.12}
\begin{equation}
\inf\limits_{p \in E^1} |a_{\nu} (p)| \| f \|_{L_2 \lr{E^1_+}}
\leq \|A_{\nu} f \|_{L_2 \lr{E^1_+}} \leq \sup\limits_{p \in E^1}
|a_{\nu} (p)| \| f \|_{L_2 \lr{E^1_+}}, \label{1.4.14}
\end{equation}
постоянные в которой являются точными. Из известной формулы
$\left| \Gamma \lr{ip + \dfrac{1}{2}} \right|=\sqrt{\pi} /
\sqrt{\ch (\pi p)},$ а так же формулы  $\left|  \sin  \dfrac{\pi
(i p - \nu -1)}{2} \right|^2 = \lr{\sh \dfrac{\pi p}{2}}^2 +
\lr{\sin \dfrac{\pi (\nu+1)}{2}}^2,$ находим, что
$$
 \sup\limits_{p \in E^1} |a_{\nu} (p)| = \sqrt{\pi}   \sup\limits_{p \in E^1} \frac{\sqrt{\lr{\sh \frac{\pi p}{2}}^2 + \lr{\cos  \frac{\pi \nu}{2}}^2}}{\sqrt{\ch (\pi p)}} =  \sqrt{\pi} \lr{\sup\limits \frac{\lr{\sh \frac{\pi p}{2}}^2 + \lr{\cos  \frac{\pi \nu}{2}}^2}{1+2\lr{\sh \frac{\pi p}{2}}^2} }^{\frac{1}{2}} =
$$
$$
 =  \sqrt{\pi} \lr{\sup\limits_{t \geq 0} \frac{t + \lr{\cos  \frac{\pi \nu}{2}}^2}{1+2t} }^{\frac{1}{2}} = \sqrt{\pi} \lr{\max\limits \left\{ \left.\frac{t + \lr{\cos  \frac{\pi \nu}{2}}^2}{1+2t}\right|_{t=0},  \left.\frac{t + \lr{\cos  \frac{\pi \nu}{2}}^2}{1+2t}\right|_{t=+ \infty} \right\}}^{\frac{1}{2}} =
$$
$$
= \sqrt{\pi} \lr{\max\limits \left\{ \frac{1}{2},  \lr{\cos
\frac{\pi \nu}{2}}^2 \right\}}^{\frac{1}{2}} = \sqrt{\pi}
\max\limits \left\{ \frac{1}{\sqrt{2}},  |\cos  \frac{\pi \nu}{2}|
\right\}.
$$
Аналогичным образом получаем
$$
 \inf\limits_{p \in E^1} |a_{\nu} (p)| = \sqrt{\pi} \min\limits \left\{ \frac{1}{\sqrt{2}}, |\cos  \frac{\pi \nu}{2}|  \right\}.
$$

Таким образом оценка~\eqref{1.4.14} принимает вид
\begin{equation}
\sqrt{\pi} \min\limits \left\{ \frac{1}{\sqrt{2}}, |\cos \frac{\pi
\nu}{2}|  \right\} \| f \|_{L_2 \lr{E^1_+}}\leq \|A_{\nu} f
\|_{L_2 \lr{E^1_+}} \leq  \sqrt{\pi} \max\limits \left\{
\frac{1}{\sqrt{2}}, |\cos  \frac{\pi \nu}{2}|  \right\} \| f
\|_{L_2 \lr{E^1_+}}. \label{1.4.15}
\end{equation}

Заметим попутно, что оператор $\sqrt{\dfrac{2}{\pi}} A_{\nu}$
унитарен в $L_2 $ только при полуцелых $\nu$ ($A_{\nu}$ сводится в
этом случае к косинус- или синус-преобразованию), а имеет
ограниченный обратный в $L_2 \lr{E_{+}^1}$ только при $\nu \neq
\pm 1, \pm 3, \dots.$

Для завершения доказательства достаточно соединить
формулу~\eqref{1.4.6}, оценку~\eqref{1.4.15}, а также равенство
Парсеваля~\eqref{1.4.2b}. Лемма~\ref{lem:1.4.1} доказана.
\end{proof}

Отметим, что ниже в главе~\ref{ch3} дано более прозрачное
доказательство этой важной для дальнейшего леммы. Первоначальное
доказательство В.\,В.~Катрахова содержало неточное утверждение об
изометричности при всех значениях параметра, оно было исправлено
в~\cite{S66, S6}. Ещё раз отметим, что идея применять технику
преобразования Меллина для оценки норм операторов преобразования,
часто совместно с теоремой Слейтер, была впервые предложена
С.\,М.~Ситником также в работах~\cite{S66, S6}. Замечательная идея
В.\,В.~Катрахова <<подправить>> операторы Сонина и Пуассона, чтобы
они ограниченно действовали в одном пространстве, реализована в
общем виде в главе~\ref{ch3} введёнными там операторами
преобразования Бушмана---Эрдейи нулевого порядка гладкости.

\begin{corollary} \label{cor:1.4.1}
    В случае $s=0,$ $R=\infty$ постоянные в оценке~\eqref{1.4.5} являются точными.
\end{corollary}

Это утверждение вытекает из точности постоянных в
оценке~\eqref{1.4.15}.

\begin{corollary} \label{cor:1.4.2}
    Оператор $\sqrt{\pi} 2^{- \nu - \frac{1}{2}} \, \Gamma^{-1} (\nu+1) S_{\nu}$    при полуцелых
    $\nu > 0$   изометрично отображает пространство $L_{2, \nu} \lr{E_{+}^1}$   на $L_{2} \lr{E_{+}^1}.$
\end{corollary}

\begin{proof}
Изометричность следует из формулы~\eqref{1.4.5}, в которой при
$s=0$ можно положить $R=\infty.$ Из  леммы~\ref{lem:1.3.2}
вытекает, что множество   состоит из  всех    функций $f$ вида $f
= I^{\frac{1}{2} - \nu} g,$ $g \in \mathring{C}_{-}^{\infty}
(\ov{E_{+}^{1}}).$ Это множество всюду плотно в пространстве $L_2
\lr{E^1_+}.$ Следствие доказано.
\end{proof}

Введем новое функциональное пространство $\mathring{H}_{\nu}^s (0,
R).$ Обозначим через $\mathring{C}_{\nu}^{\infty} (0, R)$
множество всех функций $f,$ допускающих представление  $f =
P_{\nu} g,$ в котором $g \in \mathring{C}_{\nu}^{\infty} [0, R).$
Ясно, что $\mathring{C}_{\nu}^{\infty} (0, R) \subset
\mathring{C}_{\{0\}}^{\infty} (0, R) \subset
\mathring{C}_{\{0\}}^{\infty} \lr{E_{+}^1}.$ На
$\mathring{C}_{\nu}^{\infty} (0, R)$ введем  норму
\begin{equation}
\| f \|_{\mathring{H}_{ \nu}^s (0, R)} = \frac{\sqrt{\pi}} {2^{\nu+\frac{1}{2}}\, \Gamma (\nu +1)} \| S_{\nu} f \|_{\mathring{H}_{ \nu}^s (0, R)}.
\label{1.4.16}
\end{equation}
Это определение корректно, поскольку по условию $S_{\nu} f \in
\mathring{C}_{\nu}^{\infty} [0, R).$

Обозначим через $\mathring{H}^s_{loc} (0, R)$   множество всех
функций $f,$ равных нулю при $y \geq R,$ для которых при любом
$\varepsilon>0$ конечны полунормы
$$
P_{\varepsilon, s} (f) = \|D^s f \|_{L_2 (\varepsilon, \infty)}.
$$
$\mathring{H}^s_{loc} (0, R)$  является по топологии, порожденной этими полунормами, пространством Фреше.

Замыкание в $\mathring{H}^s_{loc} (0, R),$ линеала
$\mathring{C}_{\nu}^{\infty} (0, R)$ по норме~\eqref{1.4.16} мы
будем обозначать через $\mathring{H}^s_{\nu} (0, R).$ В
разделе~\ref{sec3} в более общей  ситуации доказано вложение
$\mathring{H}^s_{\nu} (0, R) \subset \mathring{H}^s_{loc} (0, R).$

Поскольку $\mathring{C}^{\infty}_{+} [0, R) \subset
\mathring{C}^{\infty}_{\nu} (0, R)$ при любых $\nu,$ то
пространство $\mathring{H}^s_{\nu, +} (0, R)$ по
лемме~\ref{lem:1.4.1} непрерывно вложено в $\mathring{H}^s_{\nu}
(0, R),$ причём индуцированная и собственная нормы
$\mathring{H}^s_{\nu, +} (0, R)$ эквивалентны друг другу при
$s+\nu \neq 1,3,5, \dots.$

Введем весовую функцию $\sigma_{\nu}$ по формуле
\begin{equation*}
\sigma_{\nu} (y) = \left\{
\begin{array}{lll}
y^{2 \nu}, & \mbox{если} \  \Re \nu>0, \\
\dfrac{1}{\ln y} & \mbox{если} \   \nu=0.
\end{array}
\right.
\end{equation*}
Весовым граничным значением, или, короче, весовым
$\sigma_{\nu}$-следом функции в точке $y=0$ называется предел
$$
\sigma_{\nu} f |_{y=0} = \lim\limits_{y \to +0} \sigma_{\nu} (y) f(y).
$$

\begin{theorem} \label{teo:1.4.1}
Пусть $\nu \geq 0,$ $s \geq 1$ и $s+\nu>1.$ Пусть $0<R<\infty.$
Тогда у любой функции~$f$ {\rm (}после исправления её на множестве
меры нуль, если это необходимо{\rm )} из пространства
$\mathring{H}^s_{\nu} (0, R)$ существует в точке $y=0$ весовой
$\sigma_{\nu}$-след. При этом справедливы неравенства{\rm :}
    \begin{equation}
    |\sigma_{\nu} f |_{y=0} | \leq \left\{
    \begin{array}{ll}
    \dfrac{2^{\nu-1} R^{s+\nu-1} \, \Gamma (\nu)}{ \sqrt{\pi} \, \Gamma \lr{s+\nu - \frac{1}{2}} \sqrt{s+\nu-1}}
    \|  f \|_{\mathring{H}_{ \nu}^s (0, R)}, & \mbox{если} \   \nu>0,  \\
    \dfrac{2 R^{s-1}}{ \sqrt{\pi} \, \Gamma \lr{s - \frac{1}{2}} \sqrt{s-1}}
    \|  f \|_{\mathring{H}_{0}^s (0, R)}, & \mbox{если} \   \nu=0,
    \end{array}
    \right.
    \label{1.4.17}
    \end{equation}
постоянные в котором являются точными.
\end{theorem}

\begin{proof}
Неравенство~\eqref{1.4.17} достаточно доказать для функций $f
\in\mathring{C}^{\infty}_{\nu} (0, R).$ Для таких функций по
лемме~\ref{lem:1.3.1} $\sigma_{\nu}$-след существует, причём имеет
место формула
\begin{equation}
\lim\limits_{y \to +0} \sigma_{\nu} (y) f(y) = \lim\limits_{y \to
+0} \sigma_{\nu} (y) P_{\nu}^{\frac{1}{2} - \nu}
S_{\nu}^{\nu-\frac{1}{2}} f(y)= \left\{
\begin{array}{ll}
\dfrac{1}{2 \nu} S_{\nu}^{\nu-\frac{1}{2}} f|_{y=0}, & \mbox{если} \ \nu > 0, \\
 S_{0}^{-\frac{1}{2}} f|_{y=0},  & \mbox{если} \ \nu = 0.
\end{array}
\right. \label{1.4.18}
\end{equation}
Пусть функция $f \in\mathring{C}^{\infty} [0, R).$ Тогда по
неравенству Коши---Буняковского при $\alpha > \dfrac{1}{2}$
получим
$$
|g(0)|= \left| I^{\alpha} I^{-\alpha} g|_{y=0} \right| =
\frac{1}{\Gamma (\alpha)} \left| \int\limits_0^R t^{\alpha -1}
I^{-\alpha} g(t) \, dt  \right| \leq
$$
$$
\leq \frac{1}{\Gamma (\alpha)} \left( \int\limits_0^R t^{2\alpha
-2}  \, dt  \right)^{\frac{1}{2}}
 \left( \int\limits_0^R \left| I^{-\alpha} g \right|^2 dt  \right)^{\frac{1}{2}} = \frac{R^{\alpha - \frac{1}{2}}}{ \Gamma (\alpha) \sqrt {2 \alpha-1}} \left\| I^{-\alpha} g \right\|_{L_2 (0, R)}.
$$
Подставляя в это неравенство вместо $g$ функцию
$S_{\nu}^{\nu-\frac{1}{2}} f = I^{\nu - \frac{1}{2}} S_{\nu} f,$
принадлежащую, очевидно, пространству  $\mathring{C}^{\infty} [0,
R),$ получим
$$
\left|S_{\nu}^{\nu-\frac{1}{2}} f(0) \right| \leq  \frac{R^{\alpha
- \frac{1}{2}}}{ \Gamma (\alpha) \sqrt {2 \alpha-1}} \left\|
I^{\nu-\alpha-\frac{1}{2}} S_{\nu} f \right\|_{L_2 (0, R)}.
$$
Заменяя здесь $\alpha$ на $s+\nu - \dfrac{1}{2}$ и   учитывая
формулу~\eqref{1.4.18} мы и приходим к неравенству~\eqref{1.4.17}.

Покажем, что постоянные в нем являются точными. Рассмотрим функцию
$$
f_0 (y) = \left\{
\begin{array}{ll}
\int\limits_y^R (t-y)^{s+\nu - \frac{3}{2}} t^{s+\nu-\frac{3}{2}} \, dt, & \mbox{если} \  0<y<R, \\
0,  & \mbox{если} \ y \geq R.
\end{array}
\right.
$$
Покажем, что функция $I^{\frac{1}{2}-\nu} f_0$ принадлежит
$\mathring{H}^s (0, R).$ Для этого достаточно показать, что
функция  $I^{\frac{1}{2}-s-\nu} f_0$ ограничена в окрестности
точки $y=R.$ При $y<R$ имеем
$$
I^{\frac{1}{2}-s-\nu} f_0 (y) =
\frac{(-1)^{[\nu-\frac{1}{2}]+s+1}}{\Gamma \lr{\frac{3}{2} -s-\nu
+[s+\nu-\frac{1}{2}]}} D_y^{[\nu-\frac{1}{2}]+s+1} \int\limits_y^R
(t-y)^{\frac{1}{2}-\nu+[\nu - \frac{1}{2}]} \int\limits_y^R (\tau-t)^{s+\nu - \frac{3}{2}}
\tau^{s+\nu-\frac{3}{2}} \, d \tau dt =
$$
$$
=\frac{(-1)^{[\nu-\frac{1}{2}]+s+1} \, \Gamma \lr{s+\nu -
\frac{1}{2}}}{\Gamma \lr{[\nu- \frac{1}{2}] +s+1}}  D_y^{[\nu- \frac{1}{2}] +s+1 } \int\limits_y^R
\tau^{s+\nu-\frac{3}{2}} (\tau -y)^{s+[\nu-\frac{1}{2}]} d \tau =
\Gamma \lr{s+\nu -\frac{1}{2} } y^{s+\nu-\frac{3}{2}}.
$$
Следовательно, $I^{\frac{1}{2}-\nu} f_0 \in \mathring{H}^s (0, R)$
и
\begin{equation}
\left\| I^{\frac{1}{2} - \nu} f_0 \right\|_{\mathring{H}^s \lr{0,
R}} = \frac{\Gamma \lr{s+\nu -\frac{1}{2} } }{\sqrt{2s+2 \nu - 2}}
R^{s+\nu-1}. \label{1.4.19}
\end{equation}
Тогда функция $f_1 (y) = P_{\nu} I^{\frac{1}{2} - \nu} f_0 \in
\mathring{H}^s_{\nu} \lr{0, R}.$ Кроме того, при $\nu>0$
$$
\sigma_{\nu} f_1 |_{y=0} = \frac{1}{2 \nu} S^{\nu -
\frac{1}{2}}_{\nu} f_1 (0) = \frac{1}{2 \nu}  f (0) =
\frac{R^{2s+2\nu-2}}{4 \nu (s+\nu-1)},
$$
из~\eqref{1.4.19} имеем
$$
\left\| f_1 \right\|_{\mathring{H}^s_{\nu} \lr{0, R}} =
\frac{\sqrt{\pi} } {2^{\nu+\frac{1}{2}} \, \Gamma (\nu+1)} \left\|
I^{\frac{1}{2} - \nu} f_0 \right\|_{\mathring{H}^s_{\nu} \lr{0,
R}}  =  \frac{\sqrt{\pi} \,  \Gamma \lr{s+\nu -\frac{1}{2}}
R^{s+\nu-1}} {2^{\nu+1} \, \Gamma (\nu+1) \sqrt{s+\nu-1}}.
$$
Две последние формулы показывают, что первое соотношение
в~\eqref{1.4.17} превращается для функции $f_1 \in
\mathring{H}^s_{\nu} (0, R)$    в   равенство.  Аналогичный
результат справедлив и для второго соотношения. Теорема доказана.
\end{proof}

Известна следующая формула~\cite{BE1}:
\begin{equation}
\lim\limits_{\nu \to +\infty} \frac{\Gamma (\nu + \alpha)}{\Gamma (\nu + \beta)} \, \nu^{\beta-\alpha} =1.
\label{1.4.20}
\end{equation}
Поэтому имеет место

\begin{corollary} \label{cor:1.4.3}
В условиях теоремы~\ref{teo:1.4.1} справедлива оценка
\begin{equation}
\left| \sigma_{\nu} f|_{y=0}  \right| \leq c (s,  R) 2^{\nu} R^{\nu} (\nu+1)^{-s}
\left\| f \right\|_{\mathring{H}^s_{\nu} \lr{0, R}},
\label{1.4.21}
\end{equation}
где $\nu \geq 0,$ причём постоянная $c (s, R)$ зависит лишь от $s$
и $R.$
\end{corollary}

\begin{corollary} \label{cor:1.4.4}
Пусть $\nu \geq 0,$ $s-2k- \dfrac{1}{2}>0,$ $s-2k+\nu>1,$
$k=0,1,\dots$ и пусть $0<R<\infty.$ Тогда для любой функции $f \in
\mathring{H}^s_{\nu} (0, R)$    {\rm (}после исправления на
множестве меры нуль{\rm )} существует $\sigma_{\nu}$-след функции
$B_{\nu}^k f.$ Справедлива оценка
\begin{equation}
\left| \sigma_{\nu} B^k_{\nu} f|_{y=0}  \right| \leq c (s, k, R) 2^{\nu} R^{\nu} (\nu+1)^{2k-s}
\left\| f \right\|_{\mathring{H}^s_{\nu} \lr{0, R}},
\label{1.4.22}
\end{equation}
где постоянная зависит лишь от $s,$ $k$ и $R.$
\end{corollary}

\section{Многомерные операторы преобразования}\label{sec5}

\subsection{Некоторые свойства пространства
С.\,Л.~Соболева}\label{sec5.1}

Приведем сначала некоторые известные результаты.

Пусть $E^n$ обозначает евклидово $n$-мерное пространство точек $x
= (x_1, \dots, x_n).$ Пусть $\Theta$ "--- единичная сфера в $E^n.$
Введем сферические координаты $r \geq 0,$ $\vartheta \in \Theta,$
где $r = |x|,$ $\vartheta = \dfrac{x}{|x|}.$ В угловых координатах
$\varphi_1, \dots, \varphi_{n-1}$ вектор $\vartheta =
(\vartheta_1, \dots, \vartheta_n)$ выражается по формулам
\begin{eqnarray*}
\vartheta_1 = \cos \varphi_1, \\
\vartheta_2 = \sin \varphi_1 \cos \varphi_2, \\
\dots \dots \dots \dots \dots \\
\vartheta_{n-1} = \sin \varphi_1 \sin \varphi_2 \dots \sin \varphi_{n-2}  \cos \varphi_{n-1}, \\
\vartheta_{n} = \sin \varphi_1 \sin \varphi_2 \dots \sin
\varphi_{n-2}  \sin \varphi_{n-1}.
\end{eqnarray*}
Оператор Лапласа
$$
\Delta = \sum\limits_{j=1}^n \frac{\pr^2}{\pr x^2_j}
$$
в сферических координатах имеет вид
$$
\Delta = B_{\frac{n}{2}-1} + \frac{1}{r^2} \Delta_{\Theta},
$$
где оператор Бесселя $B_{\nu} = \dfrac{\pr^2}{\pr r^2} + \frac{2
\nu +1}{r} \frac{\pr}{\pr r},$ $\nu = \dfrac{n}{2}-1,$ называется
радиальной частью оператора $\Delta,$ а оператор
$$
\Delta_{\Theta} = \frac{1}{\sin^{n-2} \varphi_1} \, \frac{\pr}{\pr \varphi_1} \lr{\sin^{n-2} \varphi_1  \, \frac{\pr}{\pr \varphi_1} } + \frac{1}{\sin^{2} \varphi_1 \sin^{n-3} \varphi_2}  \, \frac{\pr}{\pr \varphi_2} \lr{\sin^{n-3} \varphi_2 \, \frac{\pr}{\pr \varphi_2} } +
$$
$$
+ \dots + \frac{1}{\sin^{2} \varphi_1 \sin^{2} \varphi_2 \dots \sin^{2} \varphi_{n-2}}  \frac{\pr^2}{\pr \varphi^2_{n-1}}
$$
называется угловой  частью  оператора $\Delta.$

Сферической гармоникой порядка $k=0,1,2, \dots$ называется функция
$Y_k (\vartheta),$ удовлетворяющая уравнению
\begin{equation}
\Delta_{\Theta} Y_k + k (n+k-2) Y_k=0
\label{2.1.1}
\end{equation}
на сфере $\Theta.$ Это уравнение имеет $d_k =
\dfrac{(n+2k-2)(k+n-3)!}{k! (n-2)!}$ линейно независимых
ортонормированных в смысле пространства квадратично суммируемых на
сфере $\Theta$ функций $L_2 (\Theta)$  решений $Y_{k, l},$ $l=1,
\dots, d_k.$ Система функций $Y_{k, l},$ $k=0,1, \dots;$ $l=1,
\dots, d_k,$ образует ортонормированный базис в $L_2 (\Theta).$
Для функции $f$ коэффициенты её разложения в ряд по сферическим
гармоникам определяются по формуле
\begin{equation}
f_{k, l} (r) = \int\limits_{\Theta} f (r, \vartheta) Y_{k, l} (\vartheta) d \vartheta.
\label{2.1.2}
\end{equation}
Для функции $f \in \mathring{C}^{\infty} \lr{E^n}$  ряд по сферическим гармоникам
\begin{equation}
f (r) = \sum\limits_{k=0}^{\infty} \sum\limits_{l=1}^{d_k}  f_{k,
l} (r) Y_{k, l} (\vartheta) \label{2.1.3}
\end{equation}
сходится к ней абсолютно и равномерно.

Пусть $L_2 \lr{E^n},$ как обычно, обозначает пространство функций
$f(x)$ с конечной нормой
$$
\|f\|_{L_2 \lr{E^n}} = \lr{\int\limits_{E^n} |f(x)|^2 dx}^{\frac{1}{2}}.
$$
Через $L_{2, \nu} \lr{E^1_{+}},$ как и выше, обозначим  множество
функций с конечной нормой
$$
\|f\|_{L_{2, \nu} \lr{E^1_+}} = \lr{\int\limits_0^{\infty} |f(r)|^2 r^{2 \nu+1} dr}^{\frac{1}{2}}.
$$

\begin{theorem}\label{teo:2.1.1}
Пусть функция $f \in L_2 \lr{E^n}.$ Тогда ряд по сферическим
гармоникам~\eqref{2.1.3} сходится к $f$ по норме пространства $L_2
\lr{E^n}.$ При этом
\begin{equation}
\|f\|_{L_2 \lr{E^n}}^2 = \sum\limits_{k=0}^{\infty}
\sum\limits_{l=1}^{d_k}  \|r^{-k} f_{k, l} (r) \|^2_{L_{2,
\frac{n}{2}+k-1} (E^1_+)}. \label{2.1.4}
\end{equation}
Обратно, пусть функции $f_{k, l},$ $k=0,1, \dots,$ $l=1, \dots,
d_k,$ таковы, что $r^{-k} f_{k, l} \in L_{2, \frac{n}{2}+k-1}
\lr{E_{+}^1},$ и ряд справа в формуле~\eqref{2.1.4} является
сходящимся. Тогда ряд в формуле~\eqref{2.1.3} сходится по норме
пространства $L_2 \lr{E^n}$ к некоторой   функции $f \in L_2
\lr{E^n},$ норма которой вычисляется по формуле~\eqref{2.1.4}.
\end{theorem}

Приведем формулу действия $n$-мерного преобразования Фурье $F,$
имеющего вид
$$
F f(\xi) = \int\limits_{E^n} f(x) e^{-i \langle x, \xi \rangle } dx,
$$
где $\xi = (\xi_1, \dots, \xi_n) \in E^n,$ $\langle x, \xi
\rangle=x_1 \xi_1 + \dots + x_n \xi_n,$ на  функциях $f = g (r)
Y_k (\vartheta).$ Пусть $\rho = |\xi|,$ $\theta =
\dfrac{\xi}{|\xi|}$ "--- сферические    координаты  в двойственном
пространстве $E^n.$ Если $r^{-k} g(r) \in L_{2, \frac{n}{2}+k-1}
\lr{E_{+}^1},$ то имеет место формула
\begin{equation}
F \lr{g Y_k} (\rho, \theta) = (-i)^k (2 \pi )^{\frac{n}{2}} \rho^{\frac{2-n}{2}} Y_k (\theta) \int\limits_0^{\infty} J_{\frac{n}{2}+k-1} (r \rho) g(r) r^{\frac{n}{2}} \, dr.
 \label{2.1.5}
\end{equation}
Заменяя функцию Бесселя первого рода $J_{\nu}$ на нормированную
функцию по формуле $2^{\nu}\Gamma ({\nu+1}) J_{\nu} (z) = z^{\nu}
j_{\nu} (z)$ и используя одномерное преобразование Ханкеля
$F_{\nu},$ определяемое формулой~\eqref{1.3.1}, предыдущую формулу
можно переписать в виде
\begin{equation}
F \lr{g Y_k} (\rho, \theta) = \frac{ (-i)^k (2 \pi )^{\frac{n}{2}}}{2^{\frac{n}{2}+k-1} \, \Gamma \lr{\frac{n}{2}+k}} \rho^{k} Y_k (\theta) F_{k+\frac{n}{2}-1} \lr{r^{-k} g} (\rho).
\label{2.1.6}
\end{equation}

Все перечисленные выше факты общеизвестны. Их можно найти,
например, в книгах С.\,Л.~Соболева~\cite{77}, И.~Стейна,
Г.~Вейса~\cite{78}. Мы лишь сформулировали их в удобных для нас
обозначениях. Формулы~\eqref{2.1.5}--\eqref{2.1.6} получены
С.~Бохнером, это частный случай известной теоремы
Функа---Гекке~\cite{BE2,SKM}.

Пусть $U_R \subset E^n$ "--- открытый шар радиуса $R < \infty$ с
центром в начале координат. Пусть как обычно
$\mathring{C}^{\infty} \lr{U_R} = \{f: f \in \mathring{C}^{\infty}
\lr{E^n}, \  \supp f \subset U_R \}.$ Для целых $s \geq 0$
определим пространство $\mathring{H}^s \lr{U_R}$ как замыкание
множества функций $\mathring{C}^{\infty} \lr{U_R}$ по   норме
\begin{equation}
\| f \|_{\mathring{H}^s \lr{U_R}} = \lr{2 \pi}^{-\frac{n}{2}} \lr{\int\limits_{E^n} |F f (\xi)|^2 |\xi|^{2s} \, d \xi}^{\frac{1}{2}},
\label{2.1.7}
\end{equation}
где $|\xi| = \lr{\xi_1^2+\dots+\xi_n^2}^{1/2}.$

Обозначим через $\mathring{T}^{\infty}_{+} \lr{U_R}$ множество функций вида
\begin{equation}
f (r, \vartheta) = \sum\limits_{k=0}^{\mathcal{K}}
\sum\limits_{l=1}^{d_k} f_{k, l} (r) Y_{k, l}(\vartheta),
\label{2.1.8}
\end{equation}
где функции $r^{-k } f_{k, l} \in \mathring{C}^{\infty}_{+} [0,
R).$ Натуральное число $\mathcal{K}$ свое для каждой функции $f.$

\begin{lemma} \label{lem:2.1.1}
    Множество $\mathring{T}^{\infty}_{+} \lr{U_R}$ всюду плотно в пространствах  $\mathring{H}^s \lr{U_R}$ при $s \geq 0.$
\end{lemma}

\begin{proof}
Так как в $\mathring{H}^s \lr{U_R}$ всюду плотно множество
линейных комбинаций функций вида
\begin{equation}
f = \chi (r) Q(x), \label{2.1.9}
\end{equation}
где $\chi (r)$ "--- произвольная функция из пространства
$\mathring{C}^{\infty}_{+} [0, R),$ $Q (x)$ "--- однородные
полиномы $n$ переменных, то достаточно показать, что функции $f,$
определяемые формулой~\eqref{2.1.9}, принадлежат классу
$\mathring{T}^{\infty}_{+} \lr{U_R}.$ Обозначим через $q$ степень
полинома $Q.$ Тогда из представления Гаусса однородных полиномов
получаем разложение
$$
Q(x) = Q_0 (x) +|x|^2 Q_1 (x) + \dots + |x|^{2 l} Q_l (x),
$$
где $2 l \leq q$ и $Q_j$ "--- однородные гармонические полиномы
степени $q-2j,$ $j=0, \dots, l.$ Ввиду однородности получаем
$$
Q_j (x) = r^{q-2j} Q_j( \vartheta),
$$
причём функции $Q_j (\vartheta) = Y_{q-2j} (\vartheta)$ являются
сферическими гармониками порядка $q-2j.$ Комбинируя две последние
формулы, получим
$$
f(x) = \chi (r) Q(x) = \chi (r) \sum\limits_{j=0}^l r^{2j}
r^{q-2j} Y_{q-2j} (\vartheta).
$$
Это разложение эквивалентно разложению~\eqref{2.1.8}, поскольку
$r^{2 j} \chi(r) \in \mathring{C}^{\infty}_{+} [0, R)$ и функции
$Y_{q-2j} (\vartheta)$ могут быть представлены в виде линейной
комбинации ортонормированных гармоник  $Y_{q-2j, l},$ $l=1, \dots,
d_{q-2j}.$ Лемма доказана.
\end{proof}

\begin{lemma}\label{lem5.2}
Пусть функции $r^{-k} f_{k, l}  \in
\mathring{H}^s_{\frac{n}{2}+k-1} (0, R).$ Тогда функция
$f=\sum\limits_{k=0}^{\mathcal{K}} \sum\limits_{l=1}^{d_k} f_{k,
l} (r) Y_{k, l} (\vartheta),$ $\mathcal{K} < \infty,$ принадлежит
пространству $\mathring{H}^s \lr{U_R}$
 и имеет место формула
\begin{equation}
\|f\|^2_{\mathring{H}^s \lr{U_R}} =\sum\limits_{k=0}^{\mathcal{K}}
\sum\limits_{l=1}^{d_k} \|r^{-k} f_{k, l}
\|^2_{\mathring{H}^s_{\frac{n}{2}+k-1} (0, R)}. \label{2.1.10}
\end{equation}
\end{lemma}

\begin{proof}
В сферических координатах формула~\eqref{2.1.7} принимает вид
$$
\|f\|^2_{\mathring{H}^s \lr{U_R}} = \frac{1}{(2 \pi)^n}
\int\limits_0^{\infty} \int\limits_{\Theta} |F f(\rho,
\vartheta)|^2 d \vartheta \rho^{2s+n-1} \, d \rho.
$$
Отсюда, используя~\eqref{2.1.6}, получаем
$$
\|f\|^2_{\mathring{H}^s \lr{U_R}} = \frac{1}{(2 \pi)^n}
\int\limits_0^{\infty} \int\limits_{\Theta} \left|
\sum\limits_{k=0}^{\mathcal{K}} \sum\limits_{l=1}^{d_k}
\frac{(-i)^k (2 \pi)^{\frac{n}{2}}}{2^{\frac{n}{2}+k-1} \, \Gamma
\lr{\frac{n}{2}+k}}  \rho^k Y_{k, l} (\vartheta)
F_{\frac{n}{2}+k-1} (r^{-k} f_{k, l}) \right|^2 d \vartheta
\rho^{2s+n-1} d \rho =
$$
$$
= \sum\limits_{k=0}^{\mathcal{K}} \sum\limits_{l=1}^{d_k}
\frac{1}{2^{n+2k-2} \, \Gamma \lr{\frac{n}{2}+k}}
\int\limits_0^{\infty} \left| F_{\frac{n}{2}+k-1} (r^{-k} f_{k,
l}) \right|^2  \rho^{2k+2s+n-1} d \rho.
$$
Здесь была использована ортонормированность системы сферических
гармоник $Y_{k, l}$ в пространстве $L_2 (\Theta).$ Полученная
формула в сочетании с определением нормы пространства
$\mathring{H}^s_{\nu, +} (0, R)$ (см. пункт~\ref{sec4.4}) и
приводит к формуле~\eqref{2.1.10}. Лемма доказана.
\end{proof}

\begin{theorem} \label{teo:2.1.2}
Для принадлежности функции $f$ пространству $\mathring{H}^s
\lr{U_R}$ при некоторых $s \geq 0$ и $0 < R< \infty$ необходимо и
достаточно, чтобы функции $r^{-k} f_{k, l} (r)$ принадлежали
пространствам $\mathring{H}^s_{\frac{n}{2}+k-1, +} (0, R),$ и
чтобы числовой ряд
\begin{equation}
\sum\limits_{k=0}^{\infty} \sum\limits_{l=1}^{d_k} \|r^{-k} f_{k,
l} \|^2_{\mathring{H}^s_{\frac{n}{2}+k-1, +} (0, R)}
\label{2.1.11}
\end{equation}
был сходящимся. При этом функциональный ряд~\eqref{2.1.3} сходится
к функции $f$ по норме пространства $\mathring{H}^s \lr{U_R}$ и
\begin{equation}
\|f\|^2_{\mathring{H}^s \lr{U_R}} =\sum\limits_{k=0}^{\infty}
\sum\limits_{l=1}^{d_k} \|r^{-k} f_{k, l}
\|^2_{\mathring{H}^s_{\frac{n}{2}+k-1, +} (0, R)}. \label{2.1.12}
\end{equation}
\end{theorem}

\begin{proof}
{\it Необходимость}. Пусть функция $f \in \mathring{H}^s
\lr{U_R}.$ Тогда по теореме~\ref{teo:2.1.1} функции  $r^{-k} f_{k,
l} \in L_{2, \frac{n}{2}+k-1, +} (0, R)$ и ряд~\eqref{2.1.3}
сходится к $f$ по норме пространства $L_2 \lr{U_R}.$ Применяя
к~\eqref{2.1.3} преобразование Фурье и учитывая
формулу~\eqref{2.1.6}, получим
$$
F f (\rho, \vartheta) = \sum\limits_{k=0}^{\infty}
\sum\limits_{l=1}^{d_k} \frac{ (-i)^k (2 \pi
)^{\frac{n}{2}}}{2^{\frac{n}{2}+k-1} \, \Gamma \lr{\frac{n}{2}+k}}
\,  \rho^{k} \, Y_{k, l} (\vartheta) F_{k+\frac{n}{2}-1}
\lr{r^{-k} f_{k, l}}.
$$
Причём ряд справа сходится в $L_2 \lr{E^n}$ к   функции $F f \in
L_2 (E^n).$ Положим
$$
\widetilde{\sum\limits}_{\mathcal{K}}=
\sum\limits_{k=0}^{\mathcal{K}} \sum\limits_{l=1}^{d_k} \frac{
(-i)^k (2 \pi )^{\frac{n}{2}}}{2^{\frac{n}{2}+k-1} \, \Gamma
\lr{\frac{n}{2}+k}} \, \rho^{k} \, Y_{k, l} (\theta)
F_{k+\frac{n}{2}-1} \lr{r^{-k} f_{k, l}}.
$$
Тогда для почти всех $\rho > 0$ функции $F f (\rho, \vartheta)$ и
$\widetilde{\sum\limits}_{\mathcal{K}} (\rho, \vartheta)$
принадлежат пространствам $L_2 (\Theta).$ Учитывая
ортонормированность $Y_{k, l},$ найдем
$$
0 \leq \| F f(p, \cdot) - \widetilde{\sum\limits}_{\mathcal{K}}
(p, \cdot)   \|_{L_2 (\Theta)} = \| F f(p, \cdot)   \|_{L_2
(\Theta)} -
$$
$$
- \sum\limits_{k=0}^{\mathcal{K}} \sum\limits_{l=1}^{d_k} \frac{(2
\pi)^n}{ 2^{n+2k-2} \, \Gamma^2 \lr{\frac{n}{2}+k}} |\rho^k
F_{\frac{n}{2}+k-1} (r^{-k} f_{k, l})|^2.
$$
Умножим левую и правую части этого неравенства на функцию
$\rho^{2s+n-1}$ и проинтегрируем. Тогда
$$
\sum\limits_{k=0}^{\mathcal{K}} \sum\limits_{l=1}^{d_k}
\frac{1}{ 2^{n+2k-2} \, \Gamma^2 \lr{\frac{n}{2}+k}}
\int\limits_0^{\infty} \rho^{2k+2s+n-1} | F_{\frac{n}{2}+k-1}
(r^{-k} f_{k, l})|^2 \, d \rho \leq
$$
$$
\leq \frac{1}{(2 \pi)^n} \int\limits_0^{\infty}
\int\limits_{\Theta} |F f(\rho, \vartheta)|^2 \rho^{2s+n-1} \, d
\rho d \vartheta = \| f \|^2_{\mathring{H}^s \lr{U_R}}.
$$
Величина справа ограничена, значит каждое слагаемое слева также
ограничено и числовой ряд сходится. Это, в частности, означает,
что функции $r^{-k} f_{k, l} \in \mathring{H}^s_{\frac{n}{2}+k-1,
+} (0, R)$ и что последовательность
$\widetilde{\sum\limits}_{\mathcal{K}}$ фундаментальна в $F
\mathring{H}^s.$ Поскольку она сходится к функции $F f$ в смысле
$L_2 (E^n),$ то тогда она   сходится и в смысле пространства $F
\mathring{H}^s.$ А тогда ряд в~\eqref{2.1.3} сходится к $f$ в
пространстве $\mathring{H}^s (U_R),$ что приводит к
формуле~\eqref{2.1.12}.

{\it Достаточность}. Пусть функции $r^{-k} f_{k, l} \in
\mathring{H}^s_{\frac{n}{2}+k-1, +} (0, R)$ и числовой
ряд~\eqref{2.1.11} сходится. Тогда последовательность
$$
{\sum\limits}_{\mathcal{K}}= \sum\limits_{k=0}^{\mathcal{K}}
\sum\limits_{l=1}^{d_k} f_{k, l} (r) Y_{k,l} (\vartheta)
$$
фундаментальна в $\mathring{H}^s (0, R).$ Таким образом, найдется
функция $f \in \mathring{H}^s (0, R),$ для которой имеет место
формула~\eqref{2.1.3}, ряд в которой сходится по норме
пространства $\mathring{H}^s (U_R).$ Теорема доказана.
\end{proof}

Определим пространство $H^s (E^n)$ как замыкание по норме
$$
\| f \|^2_{H^s \lr{E^n}} = \frac{1}{(2 \pi)^n} \int\limits_{E^n} \| F f (\xi) \|^2 (1+ |\xi|^2)^s \, d \xi
$$
множества функций $\mathring{C}^{\infty} (E^n).$

Тогда аналогично теореме~\ref{teo:2.1.2} устанавливается
\begin{theorem} \label{teo:2.1.3}
Для того, чтобы функция $f \in H^s (E^n),$ необходимо и
достаточно, чтобы функции $r^{-k} f_{k, l} \in
\mathring{H}^s_{\frac{n}{2}+k-1} (E_{+}^1)$ и чтобы ряд
\begin{equation}
\sum\limits_{k=0}^{\infty} \sum\limits_{l=1}^{d_k} \|r^{-k} f_{k,
l} \|^2_{H^s_{\frac{n}{2}+k-1, +} (E_{+}^1)} \label{2.1.13}
\end{equation}
был сходящимся. При этом функциональный ряд~\eqref{2.1.3} сходится
к функции $f$ по норме   пространства $H^s (E^n)$    и   квадрат
нормы $\|f\|^2_{H^s (E^n)}$ равен   сумме   ряда~\eqref{2.1.12}.
\end{theorem}

\subsection{Определение многомерных операторов
преобразования}\label{sec5.2}

Обозначим через $\mathring{T}^{\infty}_{\{0\}} (U_{R, 0}),$ где
$U_{R, 0}=U_R \setminus 0,$ множество функций вида
\begin{equation}
f (r, \vartheta)= \sum\limits_{k=0}^{\mathcal{K}}
\sum\limits_{l=1}^{d_k} f_{k, l} (r) Y_{k,l} (\vartheta),
\label{2.2.1}
\end{equation}
где $ f_{k, l} (r) \in \mathring{C}^{\infty}_{\{0\}} (0, R),$
$\mathcal{K}=\mathcal{K}(f)$ "--- натуральное число.   Функции
$f_{k, l},$ а   вместе с ними и $f,$ могут  иметь в начале
координат произвольную особенность.

На пространстве $\mathring{T}^{\infty}_{\{0\}} (U_{R, 0})$  определим оператор $\mathfrak{G}_n$ по формуле
$$
\mathfrak{G}_n f (r, \vartheta)= \sum\limits_{k=0}^{\mathcal{K}}
\sum\limits_{l=1}^{d_k} \frac{\sqrt{\pi}
r^{\frac{1-n}{2}}}{2^{\frac{n}{2}+k-\frac{1}{2}} \, \Gamma
\lr{\frac{n}{2}+k}} S_{\frac{n}{2}+k-1} (r^{-k} f_{k, l}) Y_{k, l}
(\vartheta),
$$
где $S_{\nu}$ "--- операторы преобразования, введённые в
пункте~\ref{sec4.1}, и $f_{k, l}$ "--- коэффициент разложения
функции $f$ в ряд по сферическим гармоникам $Y_{k, l}.$

Выведем другое представление оператора  $\mathfrak{G}_n,$   не
использующее разложение по сферическим гармоникам.
Подставляя~\eqref{2.1.2} в~\eqref{2.2.1}, находим
$$
\mathfrak{G}_n f (r, \vartheta) = \int\limits_{\Theta}
\sum\limits_{k=0}^{\mathcal{K}} \frac{\sqrt{\pi}
r^{\frac{1-n}{2}}}{2^{\frac{n}{2}+k-\frac{1}{2}} \, \Gamma
\lr{\frac{n}{2}+k}} S_{\frac{n}{2}+k-1} (r^{-k} f (r, \vartheta'))
\sum\limits_{l=1}^{d_k} Y_{k,l} (\vartheta) Y_{k,l} (\vartheta')
\, d \vartheta'.
$$
По теореме сложения сферических гармоник~\cite[с.~235]{BE2}
получаем
$$
\sum\limits_{l=1}^{d_k} Y_{k,l} (\vartheta) Y_{k,l} (\vartheta') =
\frac{\Gamma \lr{\frac{n}{2}-1}}{2 \pi^{\frac{n}{2}}}
\lr{\frac{n}{2}+k-1} C_k^{\frac{n}{2}-1} (\gamma),
$$
где $\gamma = \langle \vartheta, \vartheta' \rangle$ обозначает
скалярное   произведение    векторов $\vartheta, \vartheta' \in
E^n.$ Так как $|\vartheta|=|\vartheta'|=1,$ то $\gamma$ "---
косинус угла между ними. Через $C_k^{\lambda}$ обозначены
ультрасферические многочлены Гегенбауэра. Отсюда и из
представления оператора $S_{\nu}$ по формуле~\eqref{1.1.23}  имеем
\begin{multline}\mathfrak{G}_n f (r, \vartheta) = \frac{\Gamma
\lr{\frac{n}{2}{-}1}}{4 \pi^{\frac{n}{2}}}  \int\limits_{\Theta}
\sum\limits_{k=0}^{\mathcal{K}} \frac{ 2 \sqrt{\pi}
}{2^{\frac{n}{2}{+}k{-}\frac{1}{2}} \, \Gamma
\lr{\frac{n}{2}{+}k}} r^{\frac{1{-}n}{2}} S_{\frac{n}{2}{+}k{-}1}
(r^{{-}k} f (r, \vartheta')) \lr{\frac{n}{2}{+}k-1}
C_k^{\frac{n}{2}-1} (\gamma) \, d \vartheta' =
\\
= \frac{ - \Gamma \lr{\frac{n}{2}-1}}{4 \pi^{\frac{n}{2}}}
r^{\frac{1-n}{2}} \frac{\pr}{\pr r} \left( r^{\frac{1+n}{2}}
\int\limits_{\Theta} \sum\limits_{k=0}^{\mathcal{K}} 2
\lr{\frac{n}{2}+k-1} C_k^{\frac{n}{2}-1} (\gamma) \int\limits_0^1
t^{- \frac{n+3}{2}} \cdot P^0_{\frac{n}{2}+k-\frac{3}{2}} (t) f
\lr{\frac{r}{t}, \vartheta'} \, dt d \vartheta' \right),
\label{2.2.3}
\end{multline}
где $P^0_{\nu}$ "--- функция Лежандра первого рода.

Дальнейшие выкладки зависят от чётности или нечётности размерности
$n.$ Пусть сначала $n \geq 3$ является нечётным числом. Приведем
некоторые известные соотношения для полиномов Лежандра и
Гегенбауэра, доказательство которых можно найти
в~\cite[гл.~10]{BE2}.

При целых $p>0$ справедлива формула
\begin{equation}
D_{\gamma}^p C_m^{\lambda} (\gamma) = \frac{2^p \, \Gamma (\lambda+p)}{\Gamma (\lambda)} C_{m-p}^{\lambda+p} (\gamma),
\label{2.2.4}
\end{equation}
полагая в которой $\lambda = \dfrac{1}{2},$ $m=k+\dfrac{n-3}{2},$
$p=\dfrac{n-3}{2}$ получим
$$
D_{\gamma}^{\frac{n-3}{2}} C_{\frac{n}{2}+k-\frac{3}{2}}^{\frac{1}{2}} (\gamma) = 2^{\frac{n-3}{2}}
\frac{\Gamma \lr{\frac{n}{2}-1}}{\Gamma \lr{\frac{1}{2}}} C_{k}^{\frac{n}{2}-1} (\gamma).
$$
Поскольку $C_{\nu}^{\frac{1}{2}} \equiv P_{\nu}^0,$ то
$$
D^{\frac{n-3}{2}} P_{\frac{n}{2}+k-\frac{3}{2}}^{0} (\gamma) = 2^{\frac{n-3}{2}}
\frac{\Gamma \lr{\frac{n}{2}-1}}{\sqrt{\pi}} C_{k}^{\frac{n}{2}-1} (\gamma).
$$

Подставляя теперь последнюю формулу в~\eqref{2.2.3}, найдем
\begin{multline}
\mathfrak{G}_n f (r, \vartheta) = - 2^{- \frac{1+n}{2}}
\pi^{\frac{1-n}{2}} r^{\frac{1-n}{2}} D_r \int\limits_{\Theta}
r^{\frac{1+n}{2}} D_{\gamma}^{\frac{n-3}{2}} \sum\limits_k 2
\lr{\frac{n}{2}+k-1}    P_{\frac{n}{2}+k-\frac{3}{2}}^{0} (\gamma)
\times{}
\\
{}\times\int\limits_0^1 t^{-\frac{n+3}{2}}  P_{\frac{n}{2}+k-\frac{3}{2}}^{0}  (t) f \lr{\frac{r}{t}, \vartheta'} \, dt d \vartheta' = \\
= - 2^{- \frac{1+n}{2}} \pi^{\frac{1-n}{2}} r^{\frac{1-n}{2}} D_r
\int\limits_{\Theta} r^{\frac{1+n}{2}}
D_{\gamma}^{r^{\frac{n-3}{2}}}
\sum\limits_{k=\frac{n-3}{2}}^{\infty} (2k+1) P_k^0 (\gamma)
\int\limits_0^1 t^{-\frac{n+3}{2}} P_{k}^{0}  (t) f
\lr{\frac{r}{t}, \vartheta'} \, dt d \vartheta', \label{2.2.5}
\end{multline}
где была произведена замена переменной суммирования. Так как
полином Лежандра $P_k^0$ имеет степень $k,$ то
$D_{\gamma}^{(n-3)/2} P_k^0 (\gamma) \equiv 0$  при $k<
\dfrac{n-3}{2}.$ Следовательно, не изменяя значения последнего
выражения, в последней сумме нижний предел суммирования можно
положить равным нулю. Далее, функция $t^{- \frac{n+3}{2}} f
\lr{\frac{r}{t}, \vartheta'}$ переменной $t$ бесконечно
дифференцируема при $0<t \leq 1$ и тождественно равна нулю
окрестности левого конца. Продолжая её нулём на отрезок $[-1, 0],$
мы получим функцию класса $C^{\infty} [-1, 1].$ Для любой функции
$g \in C^{\infty} [-1, 1]$ справедливо (см.~\cite[гл.~10]{BE2})
разложение в ряд Фурье по полиномам Лежандра
$$
g (\gamma) = \sum\limits_{k=0}^{\infty} (2k+1) P_{k}^{0}  (\gamma)
\int\limits_{-1}^1 g(t) P_{k}^{0}  (t) \, dt.
$$
Отсюда в нашем случае получим
$$
\sum\limits_{k=0}^{\infty} (2k+1) P_{k}^{0}  (\gamma)
\int\limits_{0}^1t^{-\frac{n+3}{2}}  P_{k}^{0}  (t) f
\lr{\frac{r}{t}, \vartheta'} \, dt =
 \left\{
\begin{array}{ll}
\gamma^{-\frac{n+3}{2}}   f \lr{\frac{r}{\gamma}, \vartheta'}, & \mbox{если} \  \gamma>0, \\
0, & \mbox{если} \  \gamma \leq 0.
\end{array}
\right.
$$
Подставляя это выражение в~\eqref{2.2.5}, находим окончательное
представление оператора $\mathfrak{G}_n$ в случае нечётного $n$:
\begin{equation}
\mathfrak{G}_n f (r, \vartheta) = - 2^{- \frac{1+n}{2}} \pi^{\frac{1-n}{2}} r^{\frac{1-n}{2}} D_r
\int\limits_{\langle \vartheta, \vartheta'\rangle > 0} r^{\frac{1+n}{2}}   \left. D_{\gamma}^{\frac{n-3}{2}} \lr{\gamma^{-\frac{n+3}{2}}   f \lr{\frac{r}{\gamma}, \vartheta'}} \right|_{\gamma=\langle \vartheta, \vartheta'\rangle} d \vartheta'.
\label{2.2.6}
\end{equation}

Проведем аналогичные построения в случае чётного $n.$ Пусть
сначала $n \geq 4.$ Снова воспользуемся формулой~\eqref{2.2.4}, в
которой на этот раз положим $\lambda=1,$ $p = \dfrac{n}{2}-2,$ $m
= \dfrac{n}{2}-2+k.$ Тогда
$$
D_{\gamma}^{\frac{n}{2}-2} C^1_{\frac{n}{2}+k-2} (\gamma) = 2^{\frac{n}{2}-2} \, \Gamma \lr{\frac{n}{2}-1} C_k^{\frac{n}{2}-1} (\gamma).
$$
Так как $C^1_{\frac{n}{2}+k-2}$ является многочленом Чебышева
второго рода, то (см.~\cite[с.~185]{BE2}):
$$
C^1_{\frac{n}{2}+k-2} (\gamma) =  \frac{\sin \lr{\lr{\frac{n}{2}+k-1}} \arccos \gamma}{\sin (\arccos \gamma)}.
$$
Отсюда
$$
C^1_{\frac{n}{2}+k-2} (\gamma) =  \frac{1}{\frac{n}{2}+k-1} D_{\gamma} \cos \lr{\lr{\frac{n}{2}+k-1} \arccos \gamma}.
$$
Следовательно,
$$
D_{\gamma}^{\frac{n}{2}-1} \cos \lr{\lr{\frac{n}{2}+k-1} \arccos \gamma} = 2^{\frac{n}{2}-2} \, \Gamma \lr{\frac{n}{2}-1} \lr{\frac{n}{2}+k-1} C_k^{\frac{n}{2}-1} (\gamma).
$$
Подставляя эту формулу в~\eqref{2.2.3}, находим
\begin{multline}
 \mathfrak{G}_n f (r, \vartheta) = - \frac{1}{\sqrt{2}} (2
\pi)^{\frac{3-n}{2}} r^{\frac{1-n}{2}} D_r
\int\limits_{\Theta} r^{\frac{1+n}{2}} D_{\gamma}^{\frac{n}{2}-1}  \int\limits_0^1 t^{- \frac{n+3}{2}} \times \\
 \times f \lr{\frac{r}{t}, \vartheta'} \sum\limits_{k=0}^{\infty} \cos
\lr{\lr{\frac{n}{2}+k-1} \arccos \gamma}
P_{\frac{n}{2}+k-\frac{3}{2}}^{0} (t)  \, dt d \vartheta' = \\
 = - 2^{\frac{2-n}{2}} \pi^{\frac{3-n}{2}} r^{\frac{1-n}{2}} D_r \int\limits_{\Theta} r^{\frac{1+n}{2}} D_{\gamma}^{\frac{n}{2}-1}  \int\limits_0^1 t^{- \frac{n+3}{2}}  f \lr{\frac{r}{t}, \vartheta'}  \sum\limits_{k=\frac{n}{2} -1}^{\infty} \cos \lr{k \arccos
\gamma} P_{k-\frac{1}{2}}^{0} (t) \, dt d \vartheta'.
\label{2.2.7}
\end{multline}
Многочлен Чебышева, первого рода $\cos \lr{k \arccos \gamma}$
имеет степень $k.$ Тогда $D_{\gamma}^{\frac{n}{2}-1} \lr{\cos
\lr{k \arccos \gamma}}$ при $0 \leq k < \dfrac{n}{2}-1.$ Имеет
место формула (см.~\cite[с.~166]{BE2})
\begin{equation}
P_{-\frac{1}{2}}^{0} (t) + 2 \sum\limits_{k=1}^{\infty} \cos \lr{k
\arccos \gamma} P_{k-\frac{1}{2}}^{0} (t) = \left\{
\begin{array}{ll}
\dfrac{1}{ \sqrt{2 (\gamma-t)}}, & \mbox{если} \  \gamma>t,  \\
0, &  \mbox{если} \  -1<\gamma<t.
\end{array}
\right.
\label{2.2.8}
\end{equation}
Учитывая вышесказанное, окончательно находим
\begin{equation}
\mathfrak{G}_n f (r, \vartheta) =  - 2^{\frac{1+n}{2}}
\pi^{-\frac{n}{2}} r^{\frac{1-n}{2}} D_r \int\limits_{\langle
\vartheta, \vartheta'\rangle > 0} r^{\frac{1+n}{2}}
D_{\gamma}^{\frac{n}{2}-1}  \int\limits_0^{\gamma} t^{-
\frac{n+3}{2}}    \left. \frac{f \lr{\frac{r}{t},
\vartheta'}}{\sqrt{\gamma-t}} dt \right|_{\gamma=\langle
\vartheta, \vartheta'\rangle} d \vartheta'. \label{2.2.9}
\end{equation}

Последняя формула доказана для чётных $n \geq 4.$ Предложенный
метод доказательства для $n=2$   непосредственно не  проходит. Но
сама эта формула верна и при $n=2.$ В самом деле, в полярных
координатах $x_1 = r \cos \varphi,$ $x_2 = r \sin \varphi,$ имеем
$$
f (r, \varphi) = f_0 (r) + \sum\limits_{k=1}^{\infty} \lr{f_{k, 1}
(r) \cos (k \varphi) + f_{k, 2} (r) \sin (k \varphi)},
$$
где
$$
f_{k, 1} (r) = \frac{1}{\pi} \int\limits_{- \pi}^{\pi} f (r, \varphi) \cos (k \varphi) \, d \varphi, \
f_{k, 2} (r) = \frac{1}{\pi} \int\limits_{- \pi}^{\pi} f (r, \varphi) \sin (k \varphi) \, d \varphi,
$$
$$
f_{0} (r) = \frac{1}{2 \pi} \int\limits_{- \pi}^{\pi} f (r, \varphi)  \, d \varphi.
$$
Поэтому
$$
\mathfrak{G}_2 f (r, \varphi) = \sqrt{\frac{\pi}{2 r}} \left( \frac{1}{2 \pi} \int\limits_{- \pi}^{\pi} S_0
f (r, \varphi') \, d \varphi' +  \right.
$$
$$
\left. + \sum\limits_{k=1}^{\infty} \frac{1}{\pi 2^k \, k!}
\int\limits_{- \pi}^{\pi} \lr{\cos (k \varphi)\cos (k \varphi')
+\sin (k \varphi) \sin (k \varphi')} S_k\lr {r^{-k} f (r,
\varphi')}  \, d \varphi' \right) =
$$
$$
=  \frac{-1}{2 \pi \sqrt{r}} D_r \int\limits_{- \pi}^{\pi}
r^{\frac{3}{2}}  \int\limits_0^1 f \lr{\frac{r}{t}, \varphi'}
r^{-\frac{5}{2}} \lr{P_{-\frac{1}{2}}^0 (t) + 2
\sum\limits_{k=1}^{\infty}  \cos (k (\varphi- \varphi'))
P_{k-\frac{1}{2}}^0 (t)} d t d \varphi'.
$$
Учитывая~\eqref{2.2.8}, получим
\begin{equation}
\mathfrak{G}_2 f (r, \varphi) =  \frac{-1}{2^{\frac{3}{2}} \pi \sqrt{r}} D_r \int\limits_{|\varphi-\varphi'|<\frac{\pi}{2}} r^{\frac{3}{2}}  \int\limits_0^{\cos  (\varphi- \varphi')}
\frac{t^{-\frac{5}{2}} f \lr{\frac{r}{t}, \varphi'}}{\sqrt{\cos  (\varphi- \varphi')-t}} \, d t d \varphi'.
\label{2.2.10}
\end{equation}
Это и заканчивает доказательство формулы~\eqref{2.2.9} в    случае
$n=2.$

Заметим, что доказательство
формул~\eqref{2.2.6},~\eqref{2.2.9},~\eqref{2.2.10} является не
совсем полным, поскольку не обоснована замена конечных сумм
рядами. На соответствующих выкладках, использующих стандартные
свойства разложений функций по полиномам Лежандра и Чебышева мы
здесь останавливаться не будем.

Определим оператор $\mathfrak{B}_n$ на множестве $\mathring{T}^{\infty}_{\{0\}} \lr{U_{R, 0}}$ по формуле
\begin{equation}
\mathfrak{B}_n f \lr{r, \vartheta} = \sum\limits_{k=0}^{\mathcal{K}} \sum\limits_{l=1}^{d_k} \frac{2^{\frac{n}{2}+k-\frac{1}{2}} \, \Gamma \lr{\frac{n}{2}+k} }{\sqrt{\pi}} r^k P_{\frac{n}{2}+k-\frac{1}{2}} \lr{r^{\frac{n-1}{2}} f_{k, l}} Y_{k, l} \lr{\vartheta},
\label{2.2.11}
\end{equation}
где $f_{k, l}$ определены в~\eqref{2.1.2}, а $P_{\nu}$ "---
оператор преобразования из пункта~\ref{sec4.1}. Для оператора
$\mathfrak{B}_n$ невозможно получить представление
типа~\eqref{2.2.6} или~\eqref{2.2.9}. Это обусловлено существом
дела, поскольку оператор $\mathfrak{G}_n$ устраняет особенности
функций в начале координат, а $\mathfrak{B}_n$ их порождает.
Оператор $\mathfrak{B}_n,$ например, функцию без особенностей
может отобразить в функцию с существенной особенностью в начале
координат.

\begin{theorem} \label{teo:2.2.1}
Операторы $\mathfrak{G}_n$ и $\mathfrak{B}_n$ отображают пространство $\mathring{T}^{\infty}_{\{0\}} \lr{U_{R, 0}}$ на себя и являются взаимно обратными. Для функции $f \in \mathring{T}^{\infty}_{\{0\}} \lr{U_{R, 0}}$ имеют место формулы
\begin{equation}
\mathfrak{G}_n  \Delta f = r^{\frac{1-n}{2}} D^2_r
\lr{r^{\frac{n-1}{2}} \mathfrak{G}_n f},\qquad \Delta
\mathfrak{B}_n f = \mathfrak{B}_n \lr{ r^{\frac{1-n}{2}} D^2_r
\lr{r^{\frac{n-1}{2}}  f}}, \label{2.2.12}
\end{equation}
\end{theorem}

\begin{proof}
Поскольку сферические гармоники $Y_{k, l}$  суть собственные
функции оператора $\Delta_{\Theta},$ то справедлива формула
$$
\Delta f =  \lr{D^2_r + \frac{n-1}{r} D_r + \frac{1}{r^2}
\Delta_{\Theta}} \sum\limits_{k=0}^{\mathcal{K}}
\sum\limits_{l=1}^{d_k}  f_{k, l} (r) Y_{k, l} \lr{\vartheta} =
$$
$$
= \sum\limits_{k=0}^{\mathcal{K}} \sum\limits_{l=1}^{d_k}
\lr{D^2_r + \frac{n-1}{r} D_r - \frac{k(n+k-2)}{r^2} } f_{k, l}
(r) Y_{k, l} \lr{\vartheta} =
 \sum\limits_{k=0}^{\mathcal{K}} \sum\limits_{l=1}^{d_k} r^k
B_{\frac{n}{2}+k-1}  \lr{r^{-k} f_{k, l} } Y_{k, l}
\lr{\vartheta},
$$
где через $B_{\nu},$ как обычно, обозначен оператор Бесселя,
действующий по радиальной переменной. Подставляя это выражение
в~\eqref{2.2.1}, получим
$$
\mathfrak{G} \Delta f = \sum\limits_{k=0}^{\mathcal{K}}
\sum\limits_{l=1}^{d_k} \frac{\sqrt{\pi}
r^{\frac{1-n}{2}}}{2^{\frac{n}{2}+k-1} \, \Gamma
\lr{\frac{n}{2}+k}} S_{\frac{n}{2}+k-1} B_{\frac{n}{2}+k-1}
\lr{r^{-k} f_{k, l} } Y_{k, l} \lr{\vartheta}.
$$
По теореме~\ref{theorem:1_1_2} имеем $S_{\nu} B_{\nu} = D_r^2
S_{\nu},$ последнее приводит к формуле~\eqref{2.2.12}. Второе
соотношение есть следствие первого. Теорема доказана.
\end{proof}

Теорема позволяет трактовать $\mathfrak{G}_n$ и $\mathfrak{B}_n$
как операторы преобразования, преобразующие многомерный оператор
Лапласа $\Delta$ в обыкновенный дифференциальный оператор вида
$r^{\frac{1-n}{2}} D^2_r r^{\frac{n-1}{2}}.$ Операторы
$r^{\frac{n-1}{2}} \mathfrak{G}_n$ и $\mathfrak{B}_n
r^{\frac{1-n}{2}}$ преобразуют    $\Delta$ в $D^2_r.$

Отметим, что интегральные представления для многомерных операторов
преобразования $\mathfrak{G}_n$ и $\mathfrak{B}_n$ были получены в
курсовой работе С.\,М.~Ситника в 1981~г. (задача была поставлена
В.\,В.~Катраховым). При этом был найден элементарный вывод с
использованием комбинаторных свойств вместо использования
спецфункций, это позволило получить строгое и относительно простое
доказательство интегральных
представлений~\eqref{2.2.6},~\eqref{2.2.9} и как следствие
теоремы~\ref{teo:2.2.1} при минимальных предположениях на
гладкость функции. Приведённое здесь доказательство
В.\,В.~Катрахова является более <<идейным>>, оно демонстрирует
связи с несколькими конструкциями теории операторов преобразования
и спецфункциями. Но при этом не удаётся, как указано выше,
получить полностью строгое доказательство и нужны избыточные
требования на гладкость функций.

\subsection{$L_2$-теория многомерных операторов
преобразования}\label{sec5.3}

В этом пункте будут приведены некоторые результаты об операторах
преобразования, связанные с пространствами $L_2.$ Основным
результатом здесь является

\begin{theorem} \label{teo:2.3.1}
Для функций $f \in \mathring{T}_{+}^{\infty} \lr{U_R}$ при
нечётном $n$ имеет место равенство
\begin{equation}
\|\mathfrak{G}_n f \|_{L_2 \lr{U_R}} =  \| f \|_{L_2 \lr{U_R}},
\label{2.3.1}
\end{equation}
а при чётном $n$ "--- оценка
\begin{equation}
\|\mathfrak{G}_n f \|_{L_2 \lr{U_R}} \leq \sqrt{2}  \| f \|_{L_2 \lr{U_R}},
\label{2.3.2}
\end{equation}
постоянная в которой точна при $R = \infty.$
\end{theorem}

\begin{proof}
Записывая $L_2$-норму в сферических координатах и
используя~\eqref{2.2.1}, получаем
$$
\|\mathfrak{G}_n f \|_{L_2 \lr{U_R}}^2 = \int\limits_0^R
\int\limits_{\Theta} | \mathfrak{G}_n f (r, \vartheta) |^2 d
\vartheta r^{n-1} dr =
$$
$$
= \int\limits_0^R \int\limits_{\Theta}  \left|
\sum\limits_{k=0}^{\mathcal{K}} \sum\limits_{l=1}^{d_k}
\frac{\sqrt{\pi} r^{\frac{n-1}{2}}}{2^{\frac{n}{2}+k-1} \, \Gamma
\lr{\frac{n}{2}+k}} S_{\frac{n}{2}+k-1} \lr{r^{-k} f_{k, l} }
Y_{k, l} \lr{\vartheta} \right|^2 \, d \vartheta dr.
$$
Ввиду ортонормированности системы сферических гармоник $Y_{k, l}$
в $L_2 \lr{\Theta}$ имеем
$$
\|\mathfrak{G}_n f \|_{L_2 \lr{U_R}}^2 =
\sum\limits_{k=0}^{\mathcal{K}} \sum\limits_{l=1}^{d_k}
\frac{\pi}{2^{n+2k-1} \, \Gamma \lr{\frac{n}{2}+k}} \left\|
S_{\frac{n}{2}+k-1} \lr{r^{-k} f_{k, l} }  \right\|^2_{L_2 (0,
R)}.
$$
По лемме~\ref{lem:1.4.1} при $n$ нечётном имеем отсюда
$$
\|\mathfrak{G}_n f \|_{L_2 \lr{U_R}}^2 =
\sum\limits_{k=0}^{\mathcal{K}} \sum\limits_{l=1}^{d_k}  \left|
r^{-k} f_{k, l}   \right|^2_{L_{2, \frac{n}{2}+k-1} (0, R)}.
$$
А тогда из теоремы~\ref{teo:1.4.1} и следует~\eqref{2.3.1}.
Оценка~\eqref{2.3.2} доказывается аналогично. Теорема доказана.
\end{proof}

Заметим, что конечность $R$ здесь не играет никакой роли. Поэтому
формулы~\eqref{2.3.1},~\eqref{2.3.2} остаются верными и при $R=
\infty$ для функций $f \in \mathring{T}_{+}^{\infty}
\lr{U_{\infty}} =  \mathring{T}_{+}^{\infty} \lr{E^n}.$ Исходя из
следствия~\ref{cor:1.4.2} теоремы~\ref{teo:2.3.1}, мы заключаем,
что область определения $\mathring{T}_{+}^{\infty} \lr{E^n}$
оператора $\mathfrak{G}_n$ плотна    в $L_2 (E^n).$ Поэтому
оператор $\mathfrak{G}_n$ допускает расширение по непрерывности до
ограниченного в $L_2 (E^n)$ оператора. Обозначим такое расширение
через $\mathfrak{G}_n^{L_2}.$ Пусть функции    $f, g \in
\mathring{T}_{+}^{\infty} \lr{E^n}.$  Тогда   из~\eqref{2.2.1}
получим
$$
\lr{\mathfrak{G}_n^{L_2} f, g }_{L_2 (E^n)} = \int\limits_0^{\infty} r^{n-1} \int\limits_{\Theta} \mathfrak{G}_n^{L_2} (r, \vartheta) \cdot g (r, \vartheta) \, d \vartheta dr =
$$
$$
 =\int\limits_0^{\infty} r^{n-1} \int\limits_{\Theta}  \sum\limits_{k=0}^{\mathcal{K}(f)} \sum\limits_{l=1}^{d_k} \lr{ \frac{\sqrt{\pi}}{2^{\frac{n}{2}+k-1} \, \Gamma \lr{\frac{n}{2}+k}}  r^{\frac{1-n}{2}}  S_{\frac{n}{2}+k-1}
 \lr{r^{-k} f_{k, l} } Y_{k, l} \lr{\vartheta} }   \sum\limits_{k=0}^{\mathcal{K}(g)} \sum\limits_{l=1}^{d_k} g_{k, l} (r)  Y_{k, l} \lr{\vartheta} \, d \vartheta d r =
$$
$$
= \sum\limits_{k=0}^{\min\limits \lr{\mathcal{K}(f),
\mathcal{K}(g)}} \sum\limits_{l=1}^{d_k}
\frac{\sqrt{\pi}}{2^{\frac{n}{2}+k-1} \, \Gamma
\lr{\frac{n}{2}+k}}   \int\limits_0^{\infty}  S_{\frac{n}{2}+k-1}
\lr{r^{-k} f_{k, l} } g_{k, l} (r) r^{\frac{n-1}{2}} \, d r.
$$
Воспользовавшись формулой~\eqref{1.1.23} и проинтегрировав по
частям, находим
$$
\lr{\mathfrak{G}_n^{L_2} f, g }_{L_2 (E^n)} =
\sum\limits_{k=0}^{\min\limits \lr{\mathcal{K}(f),
\mathcal{K}(g)}} \sum\limits_{l=1}^{d_k}  \int\limits_0^{\infty}
\int\limits_r^{\infty} \tau^{\frac{n-1}{2}} f_{k, l} (\tau)
P_{\frac{n}{2}+k-\frac{3}{2}}^0 \lr{\frac{r}{\tau}} \, d \tau \,
D_r \lr{ r^{\frac{n-1}{2}}  g_{k, l} (r)} dr =
$$
$$
=  \sum\limits_{k=0}^{\min\limits \lr{\mathcal{K}(f),
\mathcal{K}(g)}} \sum\limits_{l=1}^{d_k}  \int\limits_0^{\infty}
\tau^{\frac{n-1}{2}} f_{k, l} (\tau) \int\limits_0^{\tau} D_r \lr{
r^{\frac{n-1}{2}}  g_{k, l} (r)}   P_{\frac{n}{2}+k-\frac{3}{2}}^0
\lr{\frac{r}{\tau}} \, d r d \tau = \lr{f,
\lr{\mathfrak{G}_n^{L_2} }^* g }_{L_2 (E^n)},
$$
где
$$
\lr{\mathfrak{G}_n^{L_2} }^* g (r, \vartheta) = \sum\limits_{k=0}^{\mathcal{K}(g)} \sum\limits_{l=1}^{d_k}
r^{\frac{n-1}{2}} \int\limits_0^r  D_{\tau} \lr{ {\tau}^{\frac{n-1}{2}}  g_{k, l} (\tau)} P_{\frac{n}{2}+k-\frac{3}{2}}^0 \lr{\frac{\tau}{r}} \, d \tau Y_{k, l} \lr{\vartheta}.
$$
Последний ряд можно просуммировать, проводя такие же выкладки, как и при выводе соответствующего представления  оператора $\mathfrak{G}_n$: \\
при нечётном $n$
\begin{equation}
\lr{\mathfrak{G}_n^{L_2} }^* g (r, \vartheta) = 2^{- \frac{1+n}{2}} \pi^{\frac{1-n}{2}} \int\limits_{\langle \vartheta, \vartheta' \rangle > 0} D^{\frac{n-1}{2}}_{\gamma} \left. \lr{\gamma^{\frac{n-1}{2}} g(\gamma r, \vartheta')} \right|_{\gamma=\langle \vartheta, \vartheta' \rangle} d \vartheta'.
\label{2.3.3}
\end{equation}
а при чётном $n \geq 2$
\begin{equation}
\lr{\mathfrak{G}_n^{L_2} }^* g (r, \vartheta) = 2^{- \frac{1+n}{2}} \pi^{-\frac{n}{2}} \int\limits_{\langle \vartheta, \vartheta' \rangle > 0} D^{\frac{n-2}{2}}_{\gamma} \int\limits_0^{\gamma}  \left. D_{\tau} \lr{\tau^{\frac{n-1}{2}} g(\tau r, \vartheta')} \frac{d \tau}{\sqrt{\gamma-\tau}} \right|_{\gamma=\langle \vartheta, \vartheta' \rangle} d \vartheta'.
\label{2.3.4}
\end{equation}

Так как при нечётном $n$ оператор $\mathfrak{G}_n^{L_2}$ будет
унитарным, то формула~\eqref{2.3.3} дает представление также и
обратного к $\mathfrak{G}_n^{L_2}$ оператора.

Операторы $\mathfrak{G}_n^{L_2}$ можно также выразить через
преобразования Фурье. Пусть функция     $f \in
\mathring{T}_{+}^{\infty} \lr{E^n}.$ Тогда   по
лемме~\ref{lem:1.3.2} получим
$$
\mathfrak{G}_n^{L_2} f (r, \vartheta) = \sum\limits_{k=0}^{\mathcal{K}} \sum\limits_{l=1}^{d_k}
\frac{\pi r^{\frac{1-n}{2}} Y_{k, l} \lr{\vartheta}}{2^{\frac{n}{2}+k-\frac{1}{2}} \, \Gamma \lr{\frac{n}{2}+k}}
I^{\frac{3}{2}-\frac{n}{2}-k} F^{-1}_{-} \lr{\rho F_{\frac{n}{2}+k-1} (r^{-k} f_{k, l}) },
$$
где $I^{\mu}$ "--- лиувиллевский оператор, $F_{-}^{-1}$ "---
обратное синус-преобразование, $F_{\nu}$ "--- преобразование
Ханкеля. Функции $\rho F_{\frac{n}{2}+k-1} (r^{-k} f_{k, l})$ по
условию нечётные. Поэтому синус-преобразование можно заменить с
соответствующим множителем на обратное одномерное преобразование
Фурье по формуле $F_{-}^{-1} g = - 2 i F^{-1} g.$ Тогда
$$
\mathfrak{G}_n^{L_2} f (r, \vartheta) = (2 \pi)^{\frac{1-n}{2}} r^{\frac{1-n}{2}}   I^{\frac{1-n}{2}} F^{-1} \sum\limits_{k=0}^{\mathcal{K}} \sum\limits_{l=1}^{d_k}
\frac{(-i)^k (2 \pi)^{\frac{n}{2}} \rho^k}{2^{\frac{n}{2}+k-1} \, \Gamma \lr{\frac{n}{2}+k}}
 F_{\frac{n}{2}+k-1} (r^{-k} f_{k, l})  Y_{k, l} \lr{\vartheta}.
$$
Сравнивая это выражение с формулой~\eqref{2.1.6}, находим
$$
\mathfrak{G}_{n}^{L_2} f (r, \vartheta) = (2 \pi  r)^{\frac{1-n}{2}} I^{\frac{1-n}{2}} F^{-1} \lr{\rho^{\frac{n-1}{2}} F f},
$$
где $F$ "--- многомерное преобразование Фурье. Эта формула может
быть записана в следующем окончательном виде:
$$
\mathfrak{G}_{n}^{L_2} f (r, \vartheta) = (-2 \pi i r)^{\frac{1-n}{2}} F^{-1} \lr{\rho^{\frac{n-1}{2}} F f}.
$$
В последних формулах неявно присутствует оператор продолжения на отрицательное значения радиальной переменной, вид которого ясен из контекста.

\chapter{Теория операторов преобразования
Бушмана---Эрдейи}\label{ch3}

В этой главе систематически излагается теория операторов
преобразования Бушмана---Эрдейи. Отметим, что в предыдущей главе
рассматриваются некоторые специальные случаи интегральных
операторов этого класса.

Необходимо отметить, что следуя традиции теории операторов
преобразования и соответствующей литературы, мы зачастую
используем термин <<операторы>> там, где более точным был бы
термин <<дифференциальные выражения>>. В основных теоремах указаны
функциональные пространства, для которых они верны. Для
результатов, содержащих явные формулы, если конкретные классы
функций не указаны, то считается, что они сформулированы для
функций, финитных на положительной полуоси (бесконечно
дифференцируемых функциях, отличных от нуля на некотором отрезке
$(a,b), a>0, b<\infty$). Выкладки можно также проводить с
использованием очень полезных классов пробных функций МакБрайда,
рассматриваемых в~\cite{McB}.

\section{Интегральные операторы преобразования Бушмана---Эрдейи первого рода и нулевого порядка
гладкости}\label{sec6}

Теперь перейдём к описанию основных свойств важнейшего класса "---
операторов преобразования Бушмана---Эрдейи. Это класс ОП, который
при определённом выборе параметров является одновременным
обобщением ОП СПД и их сопряжённых, операторов дробного
интегродифференцирования Римана---Лиувилля и Эрдейи---Кобера, а
также интегральных преобразований Мелера---Фока. Интегральные
операторы указанного вида с функциями Лежандра в ядрах впервые
встретились в работах E.\,T.~Copson по уравнению
Эйлера---Пуассона---Дарбу в конце 1950-х годов~\cite{Cop1, Cop2,
Cop3}. Впервые подробное изучение разрешимости и обратимости
данных операторов было начато в 1960-х годах в работах
Р.~Бушмана~\cite{Bus1,Bus2} и А.~Эрдейи~\cite{Erd1, Erd2, Erd3,
Erd4, Erd5}. Операторы Бушмана---Эрдейи или их аналоги изучались
также в работах T.\,P.~Higgins~\cite{Hig1}, Ta Li~\cite{Ta1,Ta2},
E.\,R.~Love~\cite{Love1,Love2}, G.\,M.~Habibullah,
K.\,N.~Srivastava, Динь Хоанг Ань~\cite{Din}, В.\,И.~Смирнова,
Б.~Рубина, Н.\,А.~Вирченко, И.~Федотовой~\cite{Vir1},
А.\,А.~Килбаса, О.\,В.~Скором\-ник~\cite{KiSk2, KiSk1} и др. При
этом в основном изучались задачи о решении интегральных уравнений
с этими операторами, их факторизации и обращения. Эти результаты
частично упомянуты в монографии~\cite{SKM}, хотя случай выбранных
нами пределов интегрирования считается там особым и не
рассматривается, за исключением одного набора формул композиции.
Некоторые результаты для особого выбора пределов были добавлены в
английское расширенное издание монографии~\cite{SaKiMar}.

Термин   <<операторы Бушмана---Эрдейи>>  как наиболее исторически
оправданный был введён С.\,М.~Ситником в~\cite{S66, S6},
впоследствии он использовался и другими авторами. Ранее
в~\cite{SKM} встречался предложенный О.\,И.~Ма\-ричевым термин
<<операторы Бушмана>>. В теории преобразования Радона и
математической томографии также используется термин <<операторы
Чебышёва---Гегенбауэра>>, см.~\cite{Rub3, Rub1, Rub2, Rub4}.
Наиболее полное изучение операторов Бушмана---Эрдейи  было
проведено в работах С.\,М.~Ситника в 1980--1990-е
годы~\cite{S1,S70,S72,S2,S73,S4,S66,S65,S6,S5,S7}, и затем
продолжено в последующие годы в~\cite{S46,S14,S103,S400,SitDis,
SSfiz,S42,S94,S38,S401,S402}. При этом необходимо отметить, что
роль операторов Бушмана---Эрдейи как ОП до указанных работ  вообще
ранее нигде не отмечалась и не рассматривалась.

Из  относительно недавних работ, в которых изучались операторы
Бушмана---Эрдейи как интегральные операторы отметим работы
Н.\,А.~Вирченко~\cite{Vir1,Vir2}, А.\,А.~Килбаса, Б.~Рубина,
А.\,В.~Глушака и их учеников. Так в работах А.\,А.~Килбаса и
О.\,В.~Скоромник~\cite{KiSk1, KiSk2} рассматривается действие
операторов Бушмана---Эрдейи в весовых пространствах Лебега, а
также  многомерные обобщения в виде интегралов по пирамидальным
областям. В монографии Н.\,А.~Вирченко и И.~Федотовой~\cite{Vir1}
вводятся некоторые обобщения стандартных функций Лежандра, а затем
рассматриваются напоминающие операторы Бушмана---Эрдейи, но не
содержащие их как частные случаи, интегральные операторы с
введёнными авторами функциями в ядрах на всей положительной
полуоси (операторы Бушмана---Эрдейи определены на части
положительной полуоси). В работах Б.~Рубина среди других
результатов описаны множества определения и образы интегральных
операторов Бушмана---Эрдейи (Гегенбауэра---Чебышёва) в некоторых
функциональных пространствах~\cite{Rub3, Rub1, Rub2, Rub4} с
приложениями результатов к теории преобразования Радона и
томографии. Операторы преобразования Бушмана---Эрдейи используются
в недавних работах А.\,В.~Глушака~\cite{Glu2, Glu3}.

Дадим определение операторов Бушмана---Эрдейи первого рода.

\begin{definition}
{\it Операторами Бушмана---Эрдейи первого рода} называются
интегральные операторы
\begin{align}
\label{2BE1}
B_{0+}^{\nu,\mu}f&=\int\limits_0^x \left( x^2-t^2\right)^{-\frac{\mu}{2}}P_\nu^\mu \left(\frac{x}{t}\right)f(t)d\,t,\\
\label{2BE1a}
E_{0+}^{\nu,\mu}f&=\int\limits_0^x \left( x^2-t^2\right)^{-\frac{\mu}{2}}\mathbb{P}_\nu^\mu \left(\frac{t}{x}\right)f(t)d\,t,\\
\label{2BE2}
B_{-}^{\nu,\mu}f&=\int\limits_x^\infty \left( t^2-x^2\right)^{-\frac{\mu}{2}}P_\nu^\mu \left(\frac{t}{x}\right)f(t)d\,t,\\
\label{2BE2a} E_{-}^{\nu,\mu}f&=\int\limits_x^\infty \left(
t^2-x^2\right)^{-\frac{\mu}{2}}\mathbb{P}_\nu^\mu
\left(\frac{x}{t}\right)f(t)d\,t.
\end{align}
\end{definition}
Здесь $P_\nu^\mu(z)$ "--- функция Лежандра первого
рода~\cite{BE1}, $\mathbb{P}_\nu^\mu(z)$ "--- та же функция на
разрезе $-1\leq  t \leq 1,$ $f(x)$ "--- локально суммируемая
функция, удовлетворяющая некоторым ограничениям на рост при $x\to
0,x\to\infty.$ Параметры $\mu,\nu$ "--- комплексные числа, $\Re
\mu <1,$ можно ограничиться значениями $\Re \nu \geq -1/2.$

Интегральные операторы указанного вида с функциями Лежандра в
ядрах впервые встретились в работах E.\,T.~Copson по уравнению
Эйлера---Пуассона---Дарбу в конце 1950-х годов. А именно, в
работах~\cite{Cop1,Cop2} рассмотрено следующее утверждение,
которое мы назовём

\medskip

\textbf{Лемма Копсона.} {\it Рассмотрим дифференциальное уравнение
в частных производных с двумя переменными{\rm :}
\begin{equation} \label{C1}
\frac{\pd^2 u(x,y)}{\pd x^2}+\frac{2\alpha}{x}\frac{\pd u(x,y)}{\pd x}=
\frac{\pd^2 u(x,y)}{\pd y^2}+\frac{2\beta}{y}\frac{\pd u(x,y)}{\pd y}
\end{equation}
{\rm (}обобщённое уравнение Эйлера---Пуассона---Дарбу или
$B$-гиперболическое уравнение по терминологии
И.\,А.~Киприянова{\rm )} в открытой четверти плоскости $x>0,\ y>0$
при положительных параметрах $\beta>\alpha>0$ с краевыми условиями
на осях координат {\rm (}характеристиках{\rm )}
\begin{equation} \label{C2}
u(x,0)=f(x), u(0,y)=g(y), f(0)=g(0).
\end{equation}
Предполагается, что решение u(x,y) является непрерывно
дифференцируемым в замкнутом первом квадранте, имеет непрерывные
вторые производные в открытом квадранте, граничные функции $f(x),
g(y)$ являются непрерывно дифференцируемыми.

Тогда, если решение поставленной задачи существует, то для него
выполняются соотношения{\rm :}
\begin{equation}
\label{2Cop1}
\frac{\pd u}{\pd y}=0, y=0,  \frac{\pd u}{\pd x}=0, x=0,
\end{equation}
\begin{equation}
\label{2Cop2} 2^\beta \Gamma(\beta{+}\frac{1}{2})\int\limits_0^1
f(xt)t^{\alpha{+}\beta{+}1} \lr{1{-}t^2}^{\frac{\beta
{-}1}{2}}P_{{-}\alpha}^{1{-}\beta}(t) dt =2^\alpha
\Gamma(\alpha{+}\frac{1}{2})\int\limits_0^1
g(xt)t^{\alpha{+}\beta{+}1} \lr{1{-}t^2}^{\frac{\alpha
{-}1}{2}}P_{{-}\beta}^{1{-}\alpha}(t) dt,
\end{equation}
$$
\Downarrow
$$
\begin{equation}
\label{2Cop3}
g(y)=\frac{2\Gamma(\beta+\frac{1}{2})}{\Gamma(\alpha+\frac{1}{2})
\Gamma(\beta-\alpha)}y^{1-2\beta} \int\limits_0^y
x^{2\alpha-1}f(x) \lr{y^2-x^2}^{\beta-\alpha-1}x \,dx,
\end{equation}
где  $P_\nu^\mu(z)$ "--- функция Лежандра первого
рода~\cite{BE1}.}

\medskip

Соотношения~\eqref{2Cop1} были известны ранее до Копсона, они
очевидны. В его работе приводится нестрогий вывод~\eqref{2Cop2},
то есть получено, что граничные функции (или значения решения на
характеристиках) не могут быть произвольными, они связаны в
современной терминологии операторами Бушмана---Эрдейи. В этом
основное содержание леммы Копсона. Далее в этой главе
дополнительно доказано, что если две функции связаны операторами
Бушмана---Эрдейи указанных порядков, то на самом деле
выполняется~\eqref{2Cop3} "--- то есть они связаны более простыми
операторами Эрдейи---Кобера.

Но отсюда не следует, как иногда отмечается, что теперь можно
сразу  получить обращение соответствующего оператора
Бушмана---Эрдейи, хотя бы формально. Для этого произвольную
функцию в правой части соответствующего уравнения надо записать
также в виде оператора  Бушмана---Эрдейи соответствующего порядка,
чтобы подогнать под лемму Копсона. Но для этого уже надо уметь
оператор Бушмана---Эрдейи обращать "--- получается порочный круг.
Таким образом, неверно приписывать Копсону первый результат по
обращению операторов  Бушмана---Эрдейи, хотя насколько нам
известно в его работе эти операторы действительно встречаются в
явном виде впервые (этот абзац является отголоском сначала спора,
а затем продолжительного обсуждения с Анатолием Александровичем
Килбасом).

Доказательства в работах~\cite{Cop1,Cop2}  являются нестрогими
рассуждениями, скорее намечающими план доказательства, может быть
поэтому данный результат не включён Копсоном в его
монографию~\cite{Cop4}. Отметим также, что в знаменитой
монографии~\cite{SKM}  и других работах даётся не совсем точная
ссылка на работу~\cite{Cop1}, которую мы здесь исправили. Эта
первоначальная работа Копсона нашла продолжение в совместной
работе с Эрдейи~\cite{Cop3}, в которой уже даётся строгий вывод,
вводятся подходящие классы функций, явно озвучена связь с дробными
интегралами и операторами Кобера---Эрдейи.

Перейдём к изложению результатов  для ОП Бушмана---Эрдейи и их
приложений к дифференциальным уравнениям с особенностями в
коэффициентах.

Все рассмотрения ведутся ниже на полуоси. Поэтому будем обозначать
через $L_2$ пространство $L_2(0, \infty)$ и $L_{2, k}$ весовое
пространство $L_{2, k}(0, \infty).$

Вначале распространим определения~\eqref{2BE1}--\eqref{2BE2a} \ на
важный не исследованный ранее случай $\mu =1.$

\begin{definition}
 Введём при $\mu =1$ {\it операторы  Бушмана---Эрдейи нулевого порядка гладкости} по формулам
\begin{align}
\label{2BE01}
B_{0+}^{\nu,1}f&={_1 S^{\nu}_{0+}f}=\frac{d}{dx}\int\limits_0^x P_\nu \left(\frac{x}{t}\right)f(t)\,dt,\\
\label{2BE02}
E_{0+}^{\nu,1}f&={_1 P^{\nu}_{-}}f=\int\limits_0^x P_\nu \left(\frac{t}{x}\right)\frac{df(t)}{dt}\,dt,\\
\label{2BE03}
B_{-}^{\nu,1}f&={_1 S^{\nu}_{-}}f=\int\limits_x^\infty P_\nu \left(\frac{t}{x}\right)(-\frac{df(t)}{dt})\,dt,\\
\label{2BE04} E_{-}^{\nu,1}f&={_1
P^{\nu}_{0+}}f=(-\frac{d}{dx})\int\limits_x^\infty P_\nu
\left(\frac{x}{t}\right)f(t)\,dt,
\end{align}
где $P_\nu(z)=P_\nu^0(z)$ "--- функция Лежандра.
\end{definition}

Разумеется, при очевидных дополнительных условиях на функции
в~\eqref{2BE01}--\eqref{2BE04} можно продифференцировать под
знаком интеграла или проинтегрировать по частям.

\begin{theorem}\label{2fact1}
 Справедливы следующие формулы факторизации операторов
Бушмана---Эрдейи на подходящих функциях через дробные интегралы
Римана---Лиувилля и Бушмана---Эрдейи нулевого порядка
гладкости{\rm :}
\begin{equation}\label{1.9}{B_{0+}^{\nu,\,\mu} f=I_{0+}^{1-\mu}~ {_1 S^{\nu}_{0+}f},~B_{-}^{\nu, \,\mu} f={_1 P^{\nu}_{-}}~ I_{-}^{1-\mu}f,}\end{equation}
\begin{equation}\label{1.10}{E_{0+}^{\nu,\,\mu} f={_1
P^{\nu}_{0+}}~I_{0+}^{1-\mu}f,~E_{-}^{\nu, \, \mu} f=
I_{-}^{1-\mu}~{_1 S^{\nu}_{-}}f.}\end{equation}
\end{theorem}

\begin{proof}
Докажем первую формулу, остальные доказываются аналогично. С
учётом определений, финитности функции $f(x),$ согласно соглашению
в начале главы, и полугруппового свойства дробных интегралов
Римана---Лиувилля, получаем $$ B_{0+}^{\nu,\,\mu}
f=I_{0+}^{1-\mu}~ {_1 S^{\nu}_{0+}f}=
 I_{0+}^{-\mu}~\int\limits_0^t P_\nu \left(\frac{t}{y}\right)f(y)\,dy=\frac{1}{\Gamma(-\mu)}~\int\limits_0^x (x-t)^{-\mu-1} \left(\int\limits_0^t
P_\nu \left(\frac{t}{y}\right)f(y)\,dy \right)d\,t.
$$ Теперь
переставим пределы интегрирования, что возможно, ввиду финитности
функции, и для вычисления внутреннего интеграла
применим~\cite[т.~3, с.~163, формула~(7)]{PBM123}. Получим нужное
интегральное представление для оператора Бушмана---Эрдейи первого
рода, теорема доказана.
\end{proof}

Эти важные формулы позволяют <<разделить>> параметры $\nu$ и
$\mu.$ Мы докажем, что операторы~\eqref{2BE01}--\eqref{2BE04}
являются изоморфизмами пространств $L_2(0, \infty),$ если $\nu$ не
равно некоторым исключительным значениям. Поэтому
операторы~\eqref{2BE1}--\eqref{2BE2a} по действию в пространствах
типа $L_2$ в определённом смысле подобны операторам дробного
интегродиффенцирования  $I^{1-\mu},$ с которыми они совпадают при
$\nu=0.$ Далее операторы Бушмана---Эрдейи будут доопределены при
всех значениях $\mu.$ Исходя из этого, введём следующее

\begin{definition}
 Число $\rho=1-\Re \mu $ назовём {\it порядком гладкости} операторов Бушмана---Эрдейи~\eqref{2BE01}--\eqref{2BE04}.
\end{definition}

Таким образом, при $\rho > 0$ (то есть при $\Re \mu > 1$)
операторы Бушмана---Эрдейи являются сглаживающими, а при $\rho <
0$ (то есть при $\Re \mu < 1$) уменьшающими гладкость в
пространствах типа $L_2 (0, \infty).$
Операторы~\eqref{2BE01}--\eqref{2BE04}, для которых $\rho = 0,$
являются по данному определению операторами нулевого порядка
гладкости. Следует пояснить, что здесь под сглаживающими мы
понимаем операторы, которые представимы в виде $A=D^k B,$ где
$k>0,$ а оператор $B$ ограничен в $L_2 (0, \infty).$ Под
уменьшающими гладкость при этом понимаются операторы, которые
действуют из некоторого пространства $C^k(0, \infty), k>0$ в
пространство Лебега $L_2 (0, \infty).$

Перечислим основные свойства операторов Бушмана---Эрдейи первого
рода~\eqref{2BE1}--\eqref{2BE2a} с функцией Лежандра I рода в
ядре. Приведём эти свойства без доказательств, так как они все
следуют из основных свойств функций Лежандра. Будем обозначать
области определения операторов через $\mathfrak{D}(B_{0+}^{\nu, \,
\mu }),$ $\mathfrak{D}(E_{0+}^{\nu, \, \mu })$ и т.~д.

Простейшие свойства функций Лежандра приводят к тождествам,
выражающим симметрию по параметрам, соотношениям смежности и
свойству сдвига по параметрам операторов Бушмана---Эрдейи.
\begin{equation}\label{2.1}
\begin{aligned}
& B_{0+}^{\nu, \, \mu} f = B_{0+}^{-\nu-1,\mu} f, & E_{0+}^{\nu,\, \mu} f=E_{0+}^{-\nu-1, \, \mu}f,
\\
& B _{-}^{\nu, \, \mu} f=B_{-}^{-\nu-1, \, \mu}f, & E_{-}^{\nu,
\,\mu} f = E_{-}^{-\nu-1, \, \mu} f,
\end{aligned}
\end{equation}
\begin{equation}\label{2.2}
\begin{aligned}
& & (2 \nu +1)x \, B_{0+}^{\nu,\, \mu} {\frac{1}{x}f}=(\nu-\mu+1)B_{0+}^{\nu+1, \, \mu}
f +  (\nu+\mu)B_{0+}^{\nu-1, \, \mu}f, \\
& & (2 \nu +1) \frac{1}{x} \, B_{-}^{\nu,\mu} {x
f}=(\nu-\mu+1)B_{-}^{\nu+1,\mu} f +  (\nu+\mu)B_{-}^{\nu-1,\mu}f,
\end{aligned}
\end{equation}
\begin{equation}\label{2.3}
\begin{aligned}
& & B_{0+}^{\nu-1, \, \mu}f - B_{0+}^{\nu+1, \, \mu}f = -(2 \nu +1)B_{0+}^{0, \, \mu-1} { \frac{1}{x}f},\\
& & B_{-}^{\nu-1, \, \mu}f - B_{-}^{\nu+1, \, \mu}f = -(2 \nu
+1)\frac{1}{x}B_{-}^{\nu, \, \mu-1} f.
\end{aligned}
\end{equation}
Из формул разложения $L_{\nu}$ (см. далее~\eqref{275}) на
множители получаются тождества

\begin{eqnarray}
& & B_{0+}^{\nu, \, \mu -1}\left(\frac{d}{dx}- \frac{\nu}{x}\right) f = B_{0+}^{\nu-1, \, \mu } f,  \label{2.4}\\
& & B_{0+}^{\nu, \, \mu -1}\left(\frac{d}{dx}+
\frac{\nu}{x}\right) f = B_{0+}^{\nu+1,  \mu } f, \label{2.5}
\end{eqnarray}
справедливые при условиях $\Re \mu < 1,$ $\Re \nu > -
\dfrac{1}{2}$
$$
\lim\limits_{y \to 0} f(y) / y^{\nu}=0, ~ f \in
\mathfrak{D}(B_{0+}^{\nu \pm 1, \, \mu }),\quad
\left(\frac{d}{dx}\pm \frac{\nu}{x}\right)f \in
\mathfrak{D}(B_{0+}^{\nu,  \, \mu -1}).
$$

Формулы~\eqref{2.1} позволяют ограничиться случаем $\Re \nu \geq -
\dfrac{1}{2}.$ Функции, на которые действуют операторы, должны
принадлежать их областям определения. Для оператора
$E_{0+}^{\nu,\, \mu}$ справедливы те же формулы, что и для
$B_{0+}^{\nu,\, \mu}.$

\begin{theorem}
Операторы Бушмана---Эрдейи~\eqref{2BE1}--\eqref{2BE2a} определены,
если $\Re \mu < 1$ или $\mu \in \mathbb{N}$ и дополнительно
выполнены условия{\rm :}
\begin{enumerate}
\item[а)] для оператора  $B_{0+}^{\nu, \, \mu}$
$$
\int\limits_0^x \sq y \,|f(y) \ln y|\, dy < \infty,
$$
если $\nu=-\dfrac{1}{2},$ $\mu \neq \dfrac{1}{2},$ а во всех
остальных случаях
$$
\int\limits_0^x y^{-\Re \nu} \,|f(y)|\, dy < \infty;
$$

\item[б)] для оператора $E_{0+}^{\nu, \, \mu}$ дополнительные
условия не требуются{\rm ;}

\item[в)] для оператора $E_{-}^{\nu, \, \mu}$
$$
\int\limits_x^{\infty} y^{-\Re \nu} \,|f(y)|\, dy < \infty;
$$

\item[г)] для оператора $B_{-}^{\nu, \, \mu}$
$$
\int\limits_x^{\infty} y^{-\frac{1}{2}-\Re \nu} \,|\ln y \cdot
f(y)|\, dy < \infty,
$$
если $\nu=-\dfrac{1}{2},$ $\mu \neq \dfrac{1}{2},$ а во всех
остальных случаях
$$
\int\limits_x^{\infty} y^{\Re(\nu-\mu)} \,|f(y)|\, dy < \infty.
$$
\end{enumerate}
\end{theorem}

В этой теореме предполагается, что функция $f(x)$ является
локально суммируемой на $(0, \infty),$ $x$ "--- произвольное
положительное число. Доказательство следует из оценок интегралов
по модулю и использования известных асимптотик для функций
Лежандра в нуле и на бесконечности, см.\mbox{\cite{BE1, NIST}}.

Важно отметить, что при некоторых специальных значениях параметров
$\nu,~\mu$ операторы Бу\-шмана---Эрдейи сводятся к более простым.
Так при значениях $\mu=-\nu$ или $\mu=\nu+2$ они являются
операторами Эрдейи---Кобера; при $\nu = 0$ операторами дробного
интегродифференцирования $I_{0+}^{1-\mu}$ или $I_{-}^{1-\mu};$ при
$\nu=-\dfrac{1}{2},$ $\mu=0$ или $\mu=1$ ядра выражаются через
эллиптические интегралы; при  $\mu=0,$  $x=1,$ $v=it-\dfrac{1}{2}$
оператор $B_{-}^{\nu, \, 0}$ лишь на постоянную отличаются от
преобразования Мелера---Фока. Таким образом, операторы
Бушмана---Эрдейи первого рода являются обобщениями всех этих
указанных классов стандартных интегральных операторов.

\begin{theorem} Пусть или $\Re \mu < 0,$ или $\mu = m \in \mathbb{N},$
$-m \leq \nu \leq m-1,$ $\nu \in \mathbb{Z}.$ Тогда справедливы
тождества
\begin{eqnarray}
& \frac{d}{dx} B_{0+}^{\nu, \, \mu } f = B_{0+}^{\nu, \, \mu +1 } f, &
E_{0+}^{\nu, \, \mu }\frac{d \, f}{dx} = E_{0+}^{\nu, \, \mu +1 } f  \label{2.6},\\
& B_{-}^{\nu, \, \mu } \lr{-\frac{d \, f}{dx}} = B_{-}^{\nu, \,
\mu + 1}f, & \lr{-\frac{d}{dx}} E_{-}^{\nu, \, \mu } f =
E_{-}^{\nu, \, \mu +1 } f \label{2.7}.
\end{eqnarray}
если все указанные операторы определены.
\end{theorem}

Эта теорема позволяет доопределить операторы Бушмана---Эрдейи и на
значения $\Re \mu \geq 1,$ переопределив их для натуральных $\mu.$

\begin{definition}
 Пусть дано число $\sigma,$ $\Re \sigma \geq 1.$
Обозначим через $m$ наименьшее натуральное число, такое, что
$\sigma= \mu +m,$ $\Re \mu <1.$ Тогда операторы Бушмана---Эрдейи
доопределим по формулам
\begin{eqnarray}
& & B_{0+}^{\nu, \, \sigma}=B_{0+}^{\nu, \, \mu + m}=\lr{\frac{d}{dx}}^m \,
B_{0+}^{\nu, \, \mu}, \nonumber\\
& & E_{0+}^{\nu, \, \sigma}=E_{0+}^{\nu, \, \mu + m}=
E_{0+}^{\nu, \, \mu} \, \lr{\frac{d}{dx}}^m, \label{2.9}\\
& & B_{-}^{\nu, \, \sigma}=B_{-}^{\nu, \, \mu + m}=
B_{-}^{\nu, \, \mu} \lr{-\frac{d}{dx}}^m, \nonumber\\
& & E_{-}^{\nu, \, \sigma}=E_{-}^{\nu, \, \mu +
m}=\lr{-\frac{d}{dx}}^m E_{-}^{\nu, \, \mu}. \nonumber
\end{eqnarray}
\end{definition}

Отметим, что при натуральных $\mu$ операторы Бушмана---Эрдейи
определены и по формулам~\eqref{2BE1}--\eqref{2BE2a}. Мы
переопределяем их для этих значений $\mu$ по формуле~\eqref{2.9}.
Таким образом, символами $B_{0+}^{\nu, \, \mu},$ $E_{0+}^{\nu, \,
\mu},$ $B_{-}^{\nu, \, \mu},$ $E_{-}^{\nu, \, \mu}$ далее мы будем
обозначать операторы, определяемые по
формулам~\eqref{2BE1}--\eqref{2BE2a} при $\Re \mu < 1,$  и по
формулам~\eqref{2.9} при $\Re \mu \geq 1.$

Будем рассматривать наряду с оператором Бесселя также тесно связанный с ним дифференциальный оператор
\begin{equation}
\label{275}
L_{\nu}=D^2-\frac{\nu(\nu+1)}{x^2}=\left(\frac{d}{dx}-\frac{\nu}{x}\right)
\left(\frac{d}{dx}+\frac{\nu}{x}\right)=\left(\frac{d}{dx}+\frac{\nu+1}{x}\right)
\left(\frac{d}{dx}-\frac{\nu+1}{x}\right),
\end{equation}
который при $\nu \in \mathbb{N}$ является оператором углового момента из квантовой механики.
Их взаимосвязь устанавливают легко проверяемые формулы связи, приведём их.

 Пусть пара ОП $X_\nu, Y_\nu$ сплетают $L_{\nu}$ и вторую  производную:
\begin{equation}
\label{276}
X_\nu L_{\nu}=D^2 X_\nu, Y_\nu D^2 = L_{\nu} Y_\nu.
\end{equation}
Введём новую пару ОП по формулам
\begin{equation}
\label{277}
S_\nu=X_{\nu-1/2} x^{\nu+1/2}, P_\nu=x^{-(\nu+1/2)} Y_{\nu-1/2}.
\end{equation}
Тогда пара новых ОП $S_\nu, P_\nu$ сплетают оператор Бесселя и вторую производную:
\begin{equation}
\label{278}
S_\nu B_\nu = D^2 S_\nu, P_\nu D^2 = B_\nu P_\nu.
\end{equation}

Разумеется, по указанным формулам можно перейти и наоборот от ОП
для оператора Бесселя к ОП для оператора углового момента. А
именно, если дана пара ОП $S_\nu, P_\nu$ со сплетающим
свойством~\eqref{278}, то новая пара ОП, определяемых по формулам
\begin{equation}
\label{279}
X_\nu=S_{\nu+1/2} x^{-(\nu+1)}, Y_\nu=x^{\nu+1} P_{\nu+1/2},
\end{equation}
удовлетворяет соотношениям~\eqref{276}.

Преимуществом ОП, сплетающих вторую производную не с оператором
Бесселя, а с оператором углового момента, является тот факт, что
при определённых условиях они оказываются ограниченными в одном
пространстве, а не в паре разных пространств. Мы сохраним за ОП,
действующим по формулам~\eqref{276}, названия ОП типа Сонина и
Пуассона соответственно.

Перейдём к определению условий, при которых ОП Бушмана---Эрдейи
первого рода действительно являются операторами преобразования.

Определим класс $\Phi(B_{0+}^{\nu, \, \mu})$ как множество функций
таких, что
\begin{enumerate}
\item[1)] $f(x) \in \mathfrak{D}(B_{0+}^{\nu, \,
\mu})\bigcap\limits C^2(0, \infty),$

\item[2)] $\lim\limits_{y \to 0}\left| \dfrac{\ln y}{\sq y}
f(y)+\sq y \ln y \cdot f'(y) \right|=0,$ если $\nu=-\dfrac{1}{2},$
$ \mu \neq \dfrac{1}{2};$\\
$ \lim\limits_{y \to 0}{(\nu+1) y^{\nu} f(y)- y^{\nu+1} f'(y) }=0,
$ если $\mu= \nu +1,$ $ \Re \nu \neq - \dfrac{1}{2};$ и,
наконец,\\
$ \lim\limits_{y \to 0}\left( \nu
\dfrac{f(y)}{y^{\nu+1}}+\dfrac{f'(y)}{y^{\nu}} \right)=0 $ во всех
остальных случаях.
\end{enumerate}

\begin{theorem}
\label{2OPB0} Пусть $f(x) \in \Phi(B_{0+}^{\nu, \, \mu}),$ $\Re
\mu \leq 1.$ Тогда оператор $B_{0+}^{\nu, \, \mu}$ является
оператором преобразования типа Сонина и удовлетворяет
соотношению~\eqref{276} на функциях $f(x).$
\end{theorem}

Аналогичный результат справедлив и для других операторов
Бушмана---Эрдейи. При этом $E_{-}^{\nu, \, \mu}$ также является
оператором типа Сонина, а $E_{0+}^{\nu, \, \mu}$ и $B_{-}^{\nu, \,
\mu}$ "--- операторами типа Пуассона. Доказательство см.
в~\cite{S66}, оно основано на общих условиях для ядер сплетающих
операторов преобразования и асимптотиках для функций Лежандра.

Можно рассматривать случай, когда нижний предел в соответствующих
интегралах~\eqref{2BE1}--\eqref{2BE2a} равен произвольному числу
$a>0,$ или верхний предел в интегралах равен произвольному
конечному числу $b>0.$ При этом все результаты этого пункта
сохраняются, а их формулировки значительно упрощаются. В
частности, все операторы Бушмана---Эрдейи в этом случае определены
при единственном условии $\Re \mu < 1$ в
форме~\eqref{2BE1}--\eqref{2BE2a} и являются операторами
преобразования на функциях $f(x)$ таких, что $f(a)=f'(a)=0$ или
($f(b)=f'(b)=0$).

Теперь сделаем важное замечание. Из полученной теоремы следует,
что ОП Бушмана---Эрдейи связывают собственные функции операторов
Бесселя и второй производной. Таким образом, половина ОП
Бушмана---Эрдейи переводят тригонометрические или экспоненциальные
функции в приведённые функции Бесселя, а другая половина наоборот.
Эти формулы здесь не приводятся, их нетрудно выписать явно. Все
они являются обобщениями исходных формул Сонина и
Пуассона~\eqref{151}--\eqref{152} и представляют существенный
интерес. Ещё раз отметим, что подобные формулы являются
непосредственными следствиями доказанных сплетающих свойств ОП
Бушмана---Эрдейи, и могут быть непосредственно проверены при
помощи таблиц интегралов от специальных функций.

Перейдём к вопросу о различных факторизациях операторов
Бушмана---Эрдейи через операторы Эрдейи---Кобера и дробные
интегралы Римана---Лиувилля (см. определения в главе~\ref{ch1}).

\begin{theorem}
\label{2factBE} Справедливы следующие формулы факторизации
операторов Бушмана---Эрдейи через операторы дробного
интегродифференцирования и Эрдейи---Кобера{\rm :}
\begin{eqnarray}
& & B_{0+}^{\nu, \, \mu}=I_{0+}^{\nu+1-\mu} I_{0+; \, 2, \, \nu+ \frac{1}{2}}^{-(\nu+1)} {\lr{\frac{2}{x}}}^{\nu+1}\label{2.17}, \\
& & E_{0+}^{\nu, \, \mu}= {\lr{\frac{x}{2}}}^{\nu+1} I_{0+; \, 2,  - \frac{1}{2}}^{\nu+1} I_{0+}^{-(\nu+\mu)}  \label{2.18}, \\
& & B_{-}^{\nu, \, \mu}= {\lr{\frac{2}{x}}}^{\nu+1}I_{-; \, 2, \, \nu+ 1}^{-(\nu+1)} I_{-}^{\nu - \mu+2}  \label{2.19}, \\
& & E_{-}^{\nu, \, \mu}= I_{-}^{-(\nu+\mu)} I_{-; \, 2, \, 0}
^{\nu+1} {\lr{\frac{x}{2}}}^{\nu+1}  \label{2.20}.
\end{eqnarray}
\end{theorem}

Доказательство приведено в~\cite{S66}, оно аналогично~\cite{Kat2,
32, KatDis, Kat3}. Ещё проще эти формулы доказываются через
мультипликаторы Меллина, см. далее.

Многие основные свойства операторов Бушмана---Эрдейи как
интегральных операторов (но не как операторов преобразования!)
могут быть выведены из теоремы~\ref{2factBE}. Так можно получить,
что формально обратным к оператору Бушмана---Эрдейи с параметрами
($\nu,~\mu$) является тот же оператор с параметрами ($\nu, ~2 -
\mu$). При этом из двух операторов "--- прямого и обратного "---
всегда один будет иметь интегральное
представление~\eqref{2BE1}--\eqref{2BE2a}, а другой определяется
формулами~\eqref{2.9}; один будет обязательно иметь положительный
порядок гладкости, а другой отрицательный (кроме операторов
нулевого порядка гладкости). Кроме того, теорема~\ref{2factBE}
позволяет доопределить операторы Бушмана---Эрдейи на всю область
значений параметров. Такое доопределение согласуется
с~\eqref{2.9}. Отметим, что
факторизации~\eqref{2.17}--\eqref{2.20} являются новыми по
сравнению с факторизациями, приведёнными в~\cite{SKM}.

Рассмотрим связь между операторами Бушмана---Эрдейи и сплетающими
операторами Сонина---Пуассона---Дельсарта (СПД). Мы предпочтём
дать новые определения для них, чтобы сохранить однообразие
обозначений в этом пункте.

\begin{definition}
Переопределим {\it операторы преобразования
Сонина---Пуассона---Дель\-сар\-та} (см. главу~\ref{ch1}) по
формулам
\begin{eqnarray}
& & {_0S_{0+}^{\nu}}=B^{\nu, \, \nu+2}_{0+}= I_{0+; \, 2, \, \nu+ \frac{1}{2}}^{-(\nu+1)} {\lr{\frac{2}{x}}}^{\nu+1}, \label{2.21} \\
& & {_0P_{0+}^{\nu}}=E^{\nu,  - \nu}_{0+}=  {\lr{\frac{x}{2}}}^{\nu+1} I_{0+; \, 2,  - \frac{1}{2}}^{\nu+1}, \label{2.22} \\
& & {_0P_{-}^{\nu}}=B^{\nu, \, \nu+2}_{-}= {\lr{\frac{2}{x}}}^{\nu+1}I_{-; \, 2, \, \nu+ 1}^{-(\nu+1)}, \label{2.23} \\
& & {_0S_{-}^{\nu}}=E^{\nu,  - \nu}_{-}= I_{-; \, 2, \, 0}
^{\nu+1} {\lr{\frac{x}{2}}}^{\nu+1}. \label{2.24}
\end{eqnarray}
\end{definition}

Это определение в сочетании
с~\eqref{162},~\eqref{1.15}--\eqref{1.16} приводит к следующим
интегральным представлениям:
\begin{align*}
{_0S_{0+}^{\nu}} f &= \begin{cases}
\ds\frac{2^{\nu+2}}{\Gamma(-\nu-1)}x \int\limits_0^x
(x^2-t^2)^{-\nu-2}t^{\nu+1}f(t)\,dt, & \Re \nu < -1,
\\
\ds\frac{2^{\nu+1}}{\Gamma(-\nu)} \frac{d}{dx} \int\limits_0^x
(x^2-t^2)^{-\nu-1}t^{\nu+1}f(t)\,dt, & \Re \nu < 0;
\end{cases}
\\
{_0P_{0+}^{\nu}} f &= \begin{cases} \ds\frac{1}{2^{\nu}
\Gamma(\nu+1)}x^{-\nu} \int\limits_0^x (x^2-t^2)^{\nu}f(t)\,dt, &
\Re \nu > -1,
\\
\ds\frac{1}{2^{\nu} \Gamma(\nu+2)}\frac{1}{x^{\nu+1}} \frac{d}{dx}
\int\limits_0^x (x^2-t^2)^{\nu+1}f(t)\,dt, & \Re \nu > -2;
\end{cases}
\\
{_0P_{-}^{\nu}} f &= \begin{cases}
\ds\frac{2^{\nu+2}}{\Gamma(-\nu-1)}x^{\nu +1}
\int\limits_x^{\infty} (t^2-x^2)^{-\nu-2}tf(t)\,dt, & \Re \nu <
-1,
\\
\ds\frac{2^{\nu+1}}{\Gamma(-\nu)} x^{\nu} \lr{-\frac{d}{dx}}
\int\limits_x^{\infty} (t^2-x^2)^{-\nu-1}tf(t)\,dt, & \Re \nu <
0;
\end{cases}
\\
{_0S_{-}^{\nu}} f &= \begin{cases} \ds\frac{1}{2^{\nu}
\Gamma(\nu+1)}\int\limits_x^{\infty}
(t^2-x^2)^{\nu}t^{-\nu}f(t)\,dt, & \Re \nu > -1,
\\
\ds\frac{1}{2^{\nu+1} \Gamma(\nu+2)} \lr{-\frac{1}{x}
\frac{d}{dx}} \int\limits_x^{\infty}
(t^2-x^2)^{\nu+1}t^{-\nu}f(t)\,dt, & \Re \nu > -2.
\end{cases}
\end{align*}

Эти операторы являются сплетающими типа Сонина или Пуассона. Если построить новые ОП для оператора углового момента (см. выше), то получим операторы типа Сонина
$$
X_{\nu} f= {_0S_{0+}^{\nu- \frac{1}{2}}} x^{\nu} f=
\frac{2^{\nu+\frac{3}{2}}}{\Gamma(-\nu-\frac{1}{2})}x
\int\limits_0^x (x^2-t^2)^{-\nu-\frac{3}{2}}t^{2 \nu+1}f(t)\,dt,
$$
если $\Re \nu < -1/2,$ а если $\Re \nu < 1/2,$ то
\begin{equation}\label{2.25}{X_{\nu} f= S_{\nu} f=
\frac{2^{\nu+\frac{1}{2}}}{\Gamma(\frac{1}{2}-\nu)}\frac{d}{dx}
\int\limits_0^x (x^2-t^2)^{-\nu-\frac{1}{2}}t^{2
\nu+1}f(t)\,dt}.\end{equation}

Аналогично получим оператор типа Пуассона вида
\begin{equation}\label{2.26}{Y_{\nu} f= P_{\nu} f=
\frac{1}{2^{\nu}\Gamma(\nu+1)}\frac{1}{x^{2 \nu}} \int\limits_0^x
(x^2-t^2)^{\nu-\frac{1}{2}}f(t)\,dt}\end{equation} при условии $\Re \nu > -1/2.$

Отметим, что из теоремы~\ref{2factBE} выводятся и
формулы~\eqref{1.9}--\eqref{1.10}.

Перейдём теперь к изучению
операторов~\eqref{2BE01}--\eqref{2BE04}. Отметим, что если функция
$f(x)$ допускает дифференцирование под знаком интеграла или
интегрирование по частям, то
операторы~\eqref{2BE01}--\eqref{2BE04} принимают вид

\begin{eqnarray}
& & {_1S_{0+}^{\nu}}f=f(x)+\int\limits_0^x \frac{\partial}{\partial x}P_{\nu}\lr{\frac{x}{y}}f(y)dy, \label{2.31} \\
& & {_1P_{0+}^{\nu}}f=f(x)-\int\limits_0^x \frac{\partial}{\partial y}P_{\nu}\lr{\frac{y}{x}}f(y)dy, \label{2.32} \\
& & {_1P_{-}^{\nu}}f=f(x)+\int\limits_x^{\infty} \frac{\partial}{\partial y}P_{\nu}\lr{\frac{y}{x}}f(y)dy, \label{2.33} \\
& & {_1S_{-}^{\nu}}f=f(x)-\int\limits_x^{\infty}
\frac{\partial}{\partial x}P_{\nu}\lr{\frac{x}{y}}f(y)dy.
\label{2.34}
\end{eqnarray}
При этом для справедливости~\eqref{2.32} и~\eqref{2.33}
соответственно дополнительно необходимы условия
$$
\lim\limits_{x \to 0} P_{\nu}(0) f(x) = 0,~\lim\limits_{x \to \infty} P_{\nu}(x) f(x) = 0.
$$

 Несложными выкладками можно доказать, что при определённых условиях на функции операторы~\eqref{2BE01}--\eqref{2BE04} являются операторами преобразования. Они сплетают оператор углового момента и вторую производную.

Из теоремы~\ref{2factBE} вытекают следующие факторизации для
операторов  Бушмана---Эрдейи нулевого порядка гладкости:

\begin{corollary}
\begin{eqnarray}
& & {_1S_{0+}^{\nu}}= I_{0+}^{\nu+1} I_{0+; \, 2, \, \nu+ \frac{1}{2}}^{-(\nu+1)} {\lr{\frac{2}{x}}}^{\nu+1}, \label{2.35} \\
& & {_1P_{0+}^{\nu}}= {\lr{\frac{x}{2}}}^{\nu+1} I_{0+; \, 2,  - \frac{1}{2}}^{\nu+1} I_{0+}^{-(\nu+1)}, \label{2.36} \\
& & {_1P_{-}^{\nu}}= {\lr{\frac{2}{x}}}^{\nu+1}I_{-; \, 2, \, \nu+ 1}^{-(\nu+1)} I_{-}^{\nu+1}, \label{2.3.7} \\
& & {_1S_{-}^{\nu}}= I_{-}^{-(\nu+1)} I_{-; \, 2, \, 0} ^{\nu+1}
{\lr{\frac{x}{2}}}^{\nu+1}. \label{2.3.8}
\end{eqnarray}
\end{corollary}

Теперь рассмотрим более подробно свойства ОП Бушмана---Эрдейи
нулевого порядка гладкости, введённых по формулам~\eqref{2BE01}.
Подобный оператор был построен В.\,В.~Катраховым~\cite{Kat2, 32,
KatDis, Kat3} путём домножения стандартного ОП Сонина на обычный
дробный интеграл с целью взаимно компенсировать гладкость этих
двух операторов и получить новый, который бы действовал в одном
пространстве типа $L_2(0,\infty),$ см. главу~\ref{ch2}. Как
впоследствии оказалось, это можно сделать известными средствами,
так как ОП Сонина "--- это частный случай операторов
Эрдейи---Кобера. Существует замечательная теорема А.~Эрдейи,
позволяющая выделить стандартный дробный интеграл
Римана---Лиувилля из дробного интеграла по любой
функции~\cite{SKM}. В результате получается

\begin{theorem}\label{2tErd}
Рассмотрим оператор дробного интегродифференцирования
Эрдейи---Ко\-бе\-ра по функции $g(x)=x^2$
$$
I_{0+; \, x^2}^{\alpha} f = \frac{1}{\Gamma(\alpha)} \int\limits_0^x (x^2-t^2)^{\alpha-1} 2t \cdot f(t)\,dt
$$
при значениях $\Re \alpha > 0.$ Тогда при $0< \Re \alpha
<\dfrac{1}{2}$ на подходящих функциях справедливо представление
оператора Эрдейи---Кобера через дробный интеграл Римана---Лиувилля
и оператор Бушмана---Эрдейи нулевого порядка
гладкости~\eqref{2BE01} при $0<\Re \alpha <\dfrac{1}{2}$
\begin{equation}
\label{2799}
  I_{0+,x^2}^{\alpha}(f)(x)  =  I_{0+}^{\alpha} \left( \left( 2x\right)^\alpha f(x) +
\int\limits_0^x \lr{\frac{\pd}{\pd x} P_{-\alpha}\left(
\frac{x}{t} \right)} \left( 2t\right)^\alpha f(t)\,dt \right)=
B_{0+}^{\nu,1}\lr{\lr{2x}^\alpha f},
\end{equation}
где $I_{0+}^{\alpha}$ "--- обычный дробный интеграл
Римана---Лиувилля.
\end{theorem}

\begin{proof}
Из теоремы А.~Эрдейи~\cite{SKM}  мы получаем представление  вида
$$
(2x)^{\alpha} f(x) +\int\limits_0^x \frac{\partial}{\partial x}
\Phi(x, s)  f(s)\,ds.
$$
Для ядра $\Phi$ справедливо представление~\cite{SKM}

$$
\Phi (x, s)= \frac{\sin \pi \alpha}{\pi} 2s \cdot \int\limits_s^x
(x-u)^{-\alpha}(u-s)^{\alpha-1}
\frac{(u-1)^{1-\alpha}}{(u^2-s^2)^{1-\alpha}} du= \frac{\sin \pi
\alpha}{\pi} \cdot 2s \cdot \int\limits_s^x
(x-u)^{-\alpha}{(u^2-s^2)^{\alpha-1}} du.
$$
Интеграл вычисляется по формуле~I из~\cite[с.~301]{PBM123}.
Получаем
$$
\Phi (x, s)= \frac{\sin \pi \alpha}{\pi}\cdot 2s \cdot
(2s)^{\alpha-1}\frac{\pi}{\sin \pi \alpha}\,{_2F_1(\alpha,
1-\alpha; 1; \frac{1}{2}-\frac{1}{2}\frac{x}{s})}=
(2s)^{\alpha}{_2F_1(\alpha, 1-\alpha; 1;
\frac{1}{2}-\frac{1}{2}\frac{x}{s})}.
$$
Осталось воспользоваться формулой~(14) из~\cite[с.~129]{PBM123}.
\end{proof}

Операторы нулевого порядка гладкости выделяются тем, что только
для них можно доказать оценки в \textit{одном} пространстве типа
$L_p(0,\infty).$ При этом, учитывая структуру этих операторов,
удобно пользоваться техникой преобразования Меллина и теоремой
Слейтер (см. главу~\ref{ch1}).

\begin{theorem} \label{2tmult}~\par
\begin{enumerate}
\item Операторы Бушмана---Эрдейи нулевого порядка гладкости
действуют по правилу~\eqref{1712}, то есть в образах
преобразования Меллина их действие сводится к умножению на
некоторые мультипликаторы. Для их мультипликаторов справедливы
формулы{\rm :}
\begin{align}
m_{{_1S_{0+}^{\nu}}}(s)&=\frac{\Gamma(-\frac{s}{2}+\frac{\nu}{2}+1) \Gamma(-\frac{s}{2}-\frac{\nu}{2}+\frac{1}{2})}{\Gamma(\frac{1}{2}-\frac{s}{2})\Gamma(1-\frac{s}{2})}=  \nonumber  \\
&=\frac{2^{-s}}{\sq{ \pi}} \frac{\Gamma(-\frac{s}{2}-\frac{\nu}{2}+\frac{1}{2}) \Gamma(-\frac{s}{2}+\frac{\nu}{2}+1)}{\Gamma(1-s)},&& \Re s < \min\limits (2 + \Re \nu, 1- \Re \nu)\label{2.311}; \\
m_{{_1P_{0+}^{\nu}}}(s)&=\frac{\Gamma(\frac{1}{2}-\frac{s}{2})\Gamma(1-\frac{s}{2})}{\Gamma(-\frac{s}{2}+\frac{\nu}{2}+1) \Gamma(-\frac{s}{2}-\frac{\nu}{2}+\frac{1}{2})},&& \Re s < 1; \label{2.312} \\
m_{{_1P_{-}^{\nu}}}(s)&=\frac{\Gamma(\frac{s}{2}+\frac{\nu}{2}+1) \Gamma(\frac{s}{2}-\frac{\nu}{2})}{\Gamma(\frac{s}{2})\Gamma(\frac{s}{2}+\frac{1}{2})},&& \Re s > \max\limits(\Re \nu, -1-\Re \nu); \label{2.313} \\
m_{{_1S_{-}^{\nu}}}(s)&=\frac{\Gamma(\frac{s}{2})\Gamma(\frac{s}{2}+\frac{1}{2})}{\Gamma(\frac{s}{2}+\frac{\nu}{2}+\frac{1}{2})
\Gamma(\frac{s}{2}-\frac{\nu}{2})},&& \Re s >0. \label{2.314}
\end{align}

\item Кроме того, выполняются следующие соотношения для
мультипликаторов{\rm :}
\begin{eqnarray}
& & m_{{_1P_{0+}^{\nu}}}(s)=1/m_{{_1S_{0+}^{\nu}}}(s),~ m_{{_1P_{-}^{\nu}}}(s)=1/m_{{_1S_{-}^{\nu}}}(s), \label{2.315} \\
& & m_{{_1P_{-}^{\nu}}}(s)=m_{{_1S_{0+}^{\nu}}}(1-s),~
m_{{_1P_{0+}^{\nu}}}(s)=m_{{_1S_{-}^{\nu}}}(1-s). \label{2.316}
\end{eqnarray}

\item Справедливы следующие формулы для норм операторов
Бушмана---Эрдейи нулевого порядка гладкости в $L_2${\rm :}
\begin{eqnarray}
& & \| _1{S_{0+}^{\nu}} \| = \| _1{P_{-}^{\nu}}\|= 1/ \min\limits(1, \sq{1- \sin \pi \nu}), \label{2.322} \\
& & \| _1{P_{0+}^{\nu}}\| = \| _1{S_{-}^{\nu}}\|= \max\limits(1,
\sq{1- \sin \pi \nu}). \label{2.323}
\end{eqnarray}

\item Нормы операторов~\eqref{2BE01}--\eqref{2BE04} периодичны по
$\nu$ с периодом 2, то есть $\|x^{\nu}\|=\|x^{\nu+2}\|,$ где
$x^{\nu}$ "--- любой из операторов~\eqref{2BE01}--\eqref{2BE04}.

\item Нормы операторов ${_1 S_{0+}^{\nu}},$ ${_1 P_{-}^{\nu}}$ не
ограничены в совокупности по $\nu,$ каждая из этих норм не меньше
$1.$ Если $\sin \pi \nu \leq 0,$ то эти нормы равны $1.$ Указанные
операторы неограничены в $L_2$ тогда и только тогда, когда $\sin
\pi \nu = 1$ {\rm (}или $\nu=(2k) + 1/2,~k \in \mathbb{Z}${\rm )}.

\item Нормы операторов ${_1 P_{0+}^{\nu}},$ ${_1 S_{-}^{\nu}}$
ограничены в совокупности по $\nu,$ каждая из этих норм не больше
$\sq{2}.$ Все эти операторы ограничены в $L_2$ при всех $\nu.$
Если $\sin \pi \nu \geq 0,$ то их $L_2$-норма равна~$1.$
Максимальное значение нормы, равное $\sq 2,$ достигается тогда и
только тогда, когда $\sin \pi \nu = -1$ {\rm (}или $\nu=
-1/2+(2k),~k \in \mathbb{Z}${\rm )}.
\end{enumerate}
\end{theorem}

\begin{proof}
Будем доказывать требуемые утверждения для первого оператора, для
остальных доказательства аналогичны.

1. Вначале докажем формулу~\eqref{1712} с нужным
мультипликатором~\eqref{2.311}. Используя последовательно
формулы~\cite[(7), с.~130; (2), с.~129; (4), с.~130]{Marich1},
получим
$$
M\left[B_{0+}^{\nu,1}\right](s)=\frac{\Gamma(2-s)}{\Gamma(1-s)}\cdot
M\left[\int\limits_0^\infty \left\{H(\frac{x}{y}-1)P_\nu
(\frac{x}{y}) \right\} \left\{y f(y)\right\}\frac{dy}{y}
\right](s-1)=
$$
$$
=\frac{\Gamma(2-s)}{\Gamma(1-s)} \cdot M \left[(x^2-1)_+^0P_\nu^0
(x) \right] (s-1)\cdot M\left[f\right](s),
$$
где использованы обозначения из~\cite{Marich1} для функции
Хевисайда и усечённой степенной функции
$$
x_+^\alpha=\left\{
\begin{array}{rl}
x^\alpha, & \mbox{если } x\geqslant 0, \\
0, & \mbox{если } x<0, \\
\end{array}\right.
\qquad H(x)=x_+^0=\left\{
\begin{array}{rl}
1, & \mbox{если } x\geqslant 0, \\
0, & \mbox{если } x<0. \\
\end{array}\right.
$$
Далее, используя формулы~\cite[14(1), с.~234; (4),
с.~130]{Marich1}, получаем
$$
M\lrs{[(x-1)_+^0 P_\nu^0 (\sqrt x)}(s)=
\frac{\Gamma(\frac{1}{2}+\frac{\nu}{2}-s)\Gamma(-\frac{\nu}{2}-s)}
{\Gamma(1-s)\Gamma(\frac{1}{2}-s)},
$$
$$
M\left[(x^2-1)_+^0 P_\nu^0 (x) \right](s-1)=\frac{1}{2}\cdot \frac
{ \Gamma(\frac{1}{2}+\frac{\nu}{2}-\frac{s-1}{2})
\Gamma(-\frac{\nu}{2}-\frac{s-1}{2}) } {
\Gamma(1-\frac{s-1}{2})\Gamma(\frac{1}{2}-\frac{s-1}{2}) }=
\frac{1}{2}\cdot\frac { \Gamma(-\frac{s}{2}+\frac{\nu}{2}+1)
\Gamma(-\frac{s}{2}-\frac{\nu}{2}+\frac{1}{2}) }
{\Gamma(-\frac{s}{2}+\frac{3}{2})\Gamma(-\frac{s}{2}+1)}
$$
при условиях $\Re s<\min\limits(2+\Re\nu, \,1-\Re\nu).$ Отсюда
выводим формулу для мультипликатора
$$
M\lrs{B_{0+}^{\nu,1}}(s)=\frac{1}{2}\cdot\frac{\Gamma(2-s)}{\Gamma(1-s)}\cdot
\Gamma\lr{-\frac{s}{2}+\frac{3}{2}}\Gamma\lr{-\frac{s}{2}+1}.
$$
Применяя к $\Gamma\lr{2-s}$ формулу Лежандра удвоения аргумента
гамма-функции (см., например,~\cite{BE1}), получим
$$
M\lrs{B_{0+}^{\nu,1}}(s)=\frac{2^{-s}}{\sqrt\pi}\cdot \frac {
\Gamma(-\frac{s}{2}+\frac{\nu}{2}+1)
\Gamma(-\frac{s}{2}-\frac{\nu}{2}+\frac{1}{2}) } {\Gamma(1-s)}.
$$
Ещё одно применение формулы удвоения Лежандра к $\Gamma(1-s)$
приводит к нужной формуле для мультипликатора~\eqref{2.311}.

В работе~\cite{S66} показано, что за счёт рассмотрения подходящих
факторизаций, условия справедливости доказанной формулы,
приведённые в~\cite{Marich1} для более общего случая,  являются в
рассматриваемом нами случае  завышенными, их можно несколько
расширить. В частности, формула для мультипликатора справедлива
при условиях $0<\Re s<1$ при всех значениях параметра $\nu,$ что
проверяется непосредственно.

2. Теперь установим формулу для нормы~\eqref{2.322}. Из найденной
формулы для мультипликатора в силу теоремы~\ref{1tMel} получаем на
прямой $\Re s=1/2, s=i u+1/2$
$$
|M\lrs{B_{0+}^{\nu,1}}(i u+1/2)|=\frac{1}{\sqrt{2\pi}}\left|\frac
{ \Gamma(-i\frac{u}{2}-\frac{\nu}{2}+\frac{1}{4})
\Gamma(-i\frac{u}{2}+\frac{\nu}{2}+\frac{3}{4}) }
{\Gamma(\frac{1}{2}-iu)}\right|.
$$
Далее будем опускать у мультипликатора указание на порождающий его
оператор. Используем формулу для модуля комплексного числа
$|z|=\sqrt{z\bar{z}}$ и тождество для гамма-функции
$\ov{\Gamma(z)}=\Gamma(\bar z),$ вытекающее из её определения в
виде интеграла. Последнее равенство справедливо для класса так
называемых вещественно-аналитических функций, к которому относится
и гамма-функция. Тогда получим
$$
|M\lrs{B_{0+}^{\nu,1}}(i u+1/2)|=
\frac{1}{\sqrt{2\pi}}\left|\frac {
\Gamma(-i\frac{u}{2}-\frac{\nu}{2}+\frac{1}{4})
\Gamma(i\frac{u}{2}-\frac{\nu}{2}+\frac{1}{4})
\Gamma(-i\frac{u}{2}+\frac{\nu}{2}+\frac{3}{4})
\Gamma(i\frac{u}{2}+\frac{\nu}{2}+\frac{3}{4}) }
{\Gamma(\frac{1}{2}-iu)\Gamma(\frac{1}{2}+iu)}\right|.
$$
В числителе объединим крайние и средние сомножители, и три
образовавшиеся пары гамма-функций преобразуем по известной формуле
(см.~\cite{BE1})
$$
\Gamma(\frac{1}{2}+z)\  \Gamma(\frac{1}{2}-z)=\frac{\pi}{\cos \pi
z}.
$$
В результате получим
$$
|M\lrs{B_{0+}^{\nu,1}}(i u+1/2)|= \sqrt{ \frac{\cos(\pi i u)}
{2\cos\pi(\frac{\nu}{2}+\frac{1}{4}+i\frac{u}{2})
\cos\pi(\frac{\nu}{2}+\frac{1}{4}-i\frac{u}{2})} }=
\sqrt{ \frac{\ch(\pi  u)}{\ch\pi u-\sin\pi\nu} }.
$$
Далее обозначим $t=\ch\pi u, 1\le t <\infty.$ Отсюда, применяя
условие из теоремы~\ref{1tMel}, получаем
$$
\sup\limits_{u\in\R} |m(i u+\frac{1}{2})|=\sup\limits_{1\le t
<\infty} \sqrt{ \frac{t}{t-\sin\pi\nu} }.
$$
Поэтому, если $\sin\pi\nu\ge 0,$ то супремум достигается при
$t=1,$ и справедлива нужная формула  для нормы
$$
\|B_{0+}^{\nu,1}\|_{L_2}=\frac{1}{\sqrt{1-\sin\pi\nu}}.
$$
Если же $\sin\pi\nu\le 0,$ то супремум достигается при
$t\to\infty,$ и справедлива формула
$$
\|B_{0+}^{\nu,1}\|_{L_2}=1.
$$
Эта часть теоремы доказана.

Утверждения~3--6 теоремы теперь напрямую следуют из найденной
формулы для нормы и условий теоремы~\ref{1tMel}. Теорема полностью
доказана.
\end{proof}

Важнейшим свойством операторов Бушмана---Эрдейи нулевого порядка
гладкости является их унитарность при целых $\nu.$ Отметим, что
при интерпретации $L_{\nu}$ как оператора углового момента в
квантовой механике, параметр $\nu$ как раз и принимает целые
неотрицательные значения. Сформулируем один из основных
результатов данной главы.

\begin{theorem}\label{2tunit} Для унитарности в $L_2$ операторов~\eqref{2BE01}--\eqref{2BE04} необходимо и достаточно, чтобы число $\nu$ было целым. В этом случае пары операторов
$({_1 S_{0+}^{\nu}},$ ${_1 P_{-}^{\nu}})$ и  $({_1 S_{-}^{\nu}},$
${_1 P_{0+}^{\nu}})$ взаимно обратны.
\end{theorem}

\begin{proof}
При $\nu \in \mathbb{Z}$ получим $\sin\pi\nu=0$ и модуль
соответствующего мультипликатора в формуле~\eqref{2.311}
тождественно равен единице на нужной прямой $\Re s =\dfrac{1}{2}.$
Поэтому по свойству~г) теоремы~\ref{1tMel} данный оператор
является унитарным в  $L_2(0,\infty),$ как и его обратный. То, что
соответствующие пары операторов являются взаимно обратными, теперь
следует из того, что они являются сопряжёнными в $L_2(0,\infty).$
Теорема доказана.
\end{proof}

Эта теорема была первоначально сформулирована в~\cite{Kat1, 30,
Kat2}, доказательство содержало неточности (утверждалась
унитарность при всех $\nu$), затем скорректированные в~\cite{S1,
S70, S72, S2, S73}, см. также~\cite{S66, S6, S46, S14, S400, S42,
S94, S38, S401, S402, SitDis}.

Перед формулировкой частного случая как  следствия предположим,
что операторы~\eqref{2BE01}--\eqref{2BE04} заданы на таких
функциях $f(x),$ что справедливы
представления~\eqref{2.31}--\eqref{2.34} (для этого достаточно
предположить, что $x f(x) \to 0$ при $x \to 0$). Тогда при $\nu=1$
\begin{equation}\label{2.325}{_1{P_{0+}^{1}}f=(I-H_1)f,~_1{S_{-}^{1}}f=(I-H_2)f,}\end{equation} где
$H_1,$ $H_2$ "--- операторы Харди (см. главу~\ref{ch1})
\begin{equation}\label{2.326}{H_1 f = \frac{1}{x} \int\limits_0^x f(y) dy,~H_2 f =
\int\limits_x^{\infty}  \frac{f(y)}{y} dy,}\end{equation} $I$ "---
единичный оператор.

\begin{corollary} Операторы~\eqref{2.325} являются унитарными взаимно обратными в $L_2$ операторами. Они сплетают дифференциальные выражения $d^2 / d x^2$ и $d^2 / d x^2 - 2/ x^2.$
\end{corollary}

Кроме того, можно показать, что операторы~\eqref{2.325} являются
преобразованиями Кэли  от симметричных операторов $\pm 2 i (x
f(x))$ при соответствующем выборе областей определения.

В унитарном случае операторы Бушмана---Эрдейи нулевого порядка
гладкости образуют пару биортогональных преобразований Ватсона, а
их ядра образуют пары несимметричных ядер Фурье~\cite{Dzh1}.
Ограниченность операторов с подобными мультипликаторами изучалась
ещё Лесли Фоксом.

Отметим важность изучения унитарности для теории интегральных  уравнений. В этом случае обратный оператор необходимо искать в виде интеграла с другими, чем у исходного, пределами интегрирования.

Рассмотрим случай $\nu = i \alpha - \dfrac{1}{2},$ $\alpha \in
\mathbb{R},$ связанный с преобразованием Мелера---Фока.

\begin{theorem} Пусть $\nu = i \alpha - \dfrac{1}{2},$ $\alpha \in \mathbb{R}.$ Тогда операторы~\eqref{2BE01}--\eqref{2BE04} ограничены в $L_2$ при всех таких $\nu.$ Для их норм справедливы формулы
$$
\| _1{S_{0+}^{i \alpha - \frac{1}{2}}}\|=\| _1{P_{-}^{i \alpha -
\frac{1}{2}}}\|=1.
$$
\end{theorem}

Доказательство такое же, как и для вещественного $\nu.$

Далее перечислим некоторые общие свойства операторов, которые
действуют по правилу~\eqref{1712} как умножение на некоторый
мультипликатор в образах преобразования  Меллина и одновременно
являются сплетающими для второй производной и оператора углового
момента.

\begin{theorem} \label{2tOPmult} Пусть оператор $S_{\nu}$ действует по формулам~\eqref{1712}
и~\eqref{276}. Тогда
\begin{enumerate}
\item[а)] его мультипликатор удовлетворяет функциональному
уравнению
\begin{equation}\label{2.5.1}{m(s)=m(s-2)\frac{(s-1)(s-2)}{(s-1)(s-2)-\nu(\nu+1)};}\end{equation}

\item[б)] если функция $p(s)$ периодична с периодом~$2$ {\rm (}то
есть $p(s)=p(s-2)${\rm ),}  то функция $p(s)m(s)$ является
мультипликатором нового оператора  преобразования $S_2^{\nu},$
опять же сплетающего $L_{\nu}$ и вторую производную по
правилу~\eqref{276}.
\end{enumerate}
\end{theorem}

\begin{proof}
Второй пункт следует из первого. Уравнение~\eqref{2.5.1}
получается из~\eqref{276} применением преобразования Меллина и
использованием формул преобразования простейших
операций~\cite{Marich1}.
\end{proof}

Последняя теорема ещё раз показывает, насколько удобно изучение ОП в терминах мультипликаторов преобразования Меллина.

Определим преобразование Стилтьеса (см., например,~\cite{SKM}) по
формуле
$$
(S f)(x)= \int\limits_0^{\infty} \frac{f(t)}{x+t} dt.
$$
Этот оператор имеет вид~\eqref{1712} с мультипликатором  $p(s)=
\pi /sin (\pi s)$ и ограничен в $L_2.$ Очевидно, что
$p(s)=p(s-2).$ Поэтому из теоремы~\ref{2tOPmult} следует, что
композиция преобразования Стилтьеса с ограниченными сплетающими
операторами~\eqref{2BE01}--\eqref{2BE04} снова является оператором
преобразованием того же типа, ограниченным в $L_2.$

Отметим, что из предыдущего изложения  следует, что
$$
\|S\|_{L_r}=|\pi / \sin \frac{\pi}{r}|,\quad r>1.
$$
С другой стороны,
$$
\|S\|_{L_{2,\, k}}=|\pi / \sin \pi k |,\quad k \notin  \mathbb{Z}.
$$
Аналогично получаются оценки в весовых пространствах $L_{r, \, k}
~ r>0.$

Рассмотрим теперь оператор $H^{\nu}$ вида~\eqref{1712} с
мультипликатором
\begin{equation}
\label{2009} m(s)=\sq{\frac{\sin \pi s - \sin \pi \nu}{\sin \pi
s}}.
\end{equation}
Из теоремы~\ref{2tmult} получаем, что на прямой $\Re s =
\dfrac{1}{2}$ величина $m(s)$ по модулю равна единице. Тогда из
этой теоремы  следует, что
\begin{equation}\label{2010}{\|H^{\nu}\|_{L_2}=\|_1{P_{0+}^{\nu}}\|_{L_2}=\|_1{S_{-}^{\nu}}\|_{L_2}.}\end{equation}
Поэтому для оператора $H^{\nu}$ справедливо заключение
теоремы~\ref{2tmult}.  В частности $H^{\nu}$ ограничен в $L_2$ при
всех $\nu.$

Отметим, что формально этот оператор связан с преобразованием Стилтьеса формулой
\begin{equation}\label{2.53}{H^{\nu}=(1-\frac{\sin \pi \nu}{\pi}S)^{\frac{1}{2}}.}\end{equation}
Одновременно с $H^{\nu}$ введём оператор $\mathfrak{D}^{\nu}$ с
мультипликатором
$$
m_{\mathfrak{D}^{\nu}}(s)=\sq{\frac{\sin \pi s}{\sin \pi s - \sin \pi \nu}}.
$$
Отсюда получаем, что
\begin{equation}\label{2.54}{\|\mathfrak{D}^{\nu}\|_{L_2}=\|_1{S_{0+}^{\nu}}\|_{L_2}=\|_1{P_{-}^{\nu}}\|_{L_2}.}\end{equation}
и кроме того, оператор $\mathfrak{D}^{\nu}$ ограничен при $\sin
\pi \nu \neq 1.$

\begin{theorem} Рассмотрим композиции операторов
\begin{eqnarray}
& &  _3{S^{\nu}_{0+}}={_1S^{\nu}_{0+}}H^{\nu},~ _3{S^{\nu}_{-}}=\mathfrak{D}^{\nu} {_1S^{\nu}_{-}}, \label{2.55} \\
& &  _3{P^{\nu}_{0+}}=\mathfrak{D}^{\nu} {_1P^{\nu}_{0+}},~
_3{P^{\nu}_{-}}={_1P^{\nu}_{-}}H^{\nu}. \label{2.56}
\end{eqnarray}
Тогда операторы ${_3 S_{0+}^{\nu}},$ ${_3 S_{-}^{\nu}}$ являются
новыми операторами  преобразования типа Сонина, а ${_3
P_{0+}^{\nu}},$ ${_3 P_{-}^{\nu}}$ "--- типа Пуассона. Все эти
операторы унитарны в $L_2.$ Кроме того, если $\sin  \pi \nu \neq
1,$ то композиции~\eqref{2.55}--\eqref{2.56} можно вычислять в
любом порядке.
\end{theorem}

Доказательство этой теоремы очевидно и следует из перехода к
мультипликаторам. Аналогичная идея применена в~\cite{Lud} для
подправления преобразования Радона до изометрии.

Ниже будет получено явное интегральное представление операторов
преобразования, сплетающих $L_{\nu}$ и $d^2/dx^2,$ которые
\textit{являются унитарными при всех $\nu \in \mathbb{R}$}.

Изучим вопрос о взаимосвязи разносторонних операторов
Бушмана---Эрдейи. Полученные формулы аналогичны тем, которые
связывают лево- и правосторонние дробные интегралы
Римана---Лиувилля (см.~\cite[с.~163--171]{SKM}). С этой целью
введём оператор \begin{equation}\label{2.57}{C^{\nu} f=
f(x)-\frac{\sin \pi \nu}{\pi} S f,}\end{equation} где $S$ "---
преобразование Стилтьеса. Приведём без доказательства свойства
$C^{\nu}$:

\begin{enumerate}
\item[1)] $\|C^{\nu}\|_{L_2}=\min\limits (1, 1 - \sin \pi \nu)
\leq 1, ~ \nu \in \mathbb{R};$

\item[2)] $\|C^{\nu}\|_{L_2}=1+\ch \pi \alpha, ~  \nu=i \alpha -
\dfrac{1}{2},~ \alpha \in \mathbb{R}.$
\end{enumerate}

\begin{theorem} При $\nu \in \mathbb{R}$ справедливы тождества для композиций
\begin{eqnarray}
& &  C^{\nu}={_1S_{-}^{\nu}} \ {_1P_{0+}^{\nu}}={_1P_{0+}^{\nu}} \  {_1S_{-}^{\nu}}, \label{2.58} \\
& &  {_1S_{-}^{\nu}}={_1S_{0+}^{\nu}} \ C^{\nu}, ~   {_1P_{0+}^{\nu}}= {_1P_{-}^{\nu}} \ C^{\nu}, \label{2.59} \\
& &  {_1S_{-}^{\nu}}=C^{\nu} \  {_1S_{0+}^{\nu}}, ~
{_1P_{0+}^{\nu}}=C^{\nu} \  {_1P_{-}^{\nu}}, ~ \sin \pi \nu \neq
1. \label{2.510}
\end{eqnarray}
\end{theorem}

\section{Интегральные операторы преобразования Бушмана---Эрдейи второго рода и унитарные  операторы преобразования Сонина---Катрахова и
Пуассона---Катрахова}\label{sec7}
\sectionmarknum{Операторы преобразования Бушмана---Эрдейи второго рода,  Сонина---Катрахова и
Пуассона---Катрахова}

Теперь определим и изучим операторы Бушмана---Эрдейи второго рода.
В этом пункте для краткости некоторые доказательства будут
опущены, так как они в основном повторяют доказательства из
предыдущего пункта.

\begin{definition}
Введём новую пару {\it операторов Бушмана---Эрдейи} с функциями
Лежандра второго рода~\cite{BE1} в ядре:
\begin{equation}\label{2.61}{{_2S^{\nu}}f=\frac{2}{\pi} \left( - \int\limits_0^x (x^2-y^2)^{-\frac{1}{2}}Q_{\nu}^1 (\frac{x}{y}) f(y) dy  +
\int\limits_x^{\infty} (y^2-x^2)^{-\frac{1}{2}}\mathbb{Q}_{\nu}^1
(\frac{x}{y}) f(y) dy\right),}\end{equation}
\begin{equation}\label{2.62}{{_2P^{\nu}}f=\frac{2}{\pi} \left( - \int\limits_0^x
(x^2-y^2)^{-\frac{1}{2}}\mathbb{Q}_{\nu}^1 (\frac{y}{x}) f(y) dy -
\int\limits_x^{\infty} (y^2-x^2)^{-\frac{1}{2}}Q_{\nu}^1
(\frac{y}{x}) f(y) dy\right).}\end{equation}
\end{definition}
При $y \to x \pm 0$ интегралы понимаются в смысле главного
значения. Отметим без доказательства, что эти операторы определены
и являются сплетающими при некоторых условиях на функции $f(x)$
(при этом оператор~\eqref{2.61} будет типа Сонина,~\eqref{2.62}
"--- типа Пуассона).

\begin{theorem}\label{2rod}  Операторы~\eqref{2.61}-\eqref{2.62} представимы в виде~\eqref{1712}
с мультипликаторами
\begin{eqnarray}
& & m_{_2S^{\nu}}(s)=p(s) \ m_{_1S_{-}^{\nu}}(s), \label{2.63}\\
& & m_{_2P^{\nu}}(s)=\frac{1}{p(s)} \ m_{_1P_{-}^{\nu}}(s),
\label{2.64}
\end{eqnarray}
где мультипликаторы операторов ${_1S_-^{\nu}},$ ${_1P_-^{\nu}}$
определены формулами~\eqref{2.313}-\eqref{2.314}, а функция $p(s)$
{\rm (}с периодом~$2${\rm )} равна
\begin{equation}\label{2.65}{p(s)=\frac{\sin \pi \nu+ \cos \pi s}{\sin \pi \nu - \sin \pi s}.}\end{equation}
\end{theorem}

Вначале докажем лемму.
\begin{lemma} \label{2lem2rod} Рассмотрим более общий чем~\eqref{2.61} интегральный оператор при значениях $\Re \nu < 1${\rm :}
\begin{equation}\label{2.66}
{_3S^{\nu,\mu}}f=\frac{2}{\pi} \left( \int\limits_0^x (x^2+y^2)^{-\frac{\mu}{2}} e^{-\mu \pi i} Q_{\nu}^{\mu}( \frac{x}{y}) f(y)\, dy + \int\limits_x^{\infty} (y^2+x^2)^{-\frac{\mu}{2}}\mathbb{Q}_{\nu}^{\mu} (\frac{x}{y}) f(y)\, dy\right),
\end{equation}
где $Q_{\nu}^{\mu}(z)$ "--- функция Лежандра второго рода,
$\mathbb{Q}_{\nu}^{\mu}(z)$ "--- значение этой функции на разрезе.

Тогда на функциях из $C_0^{\infty}(0, \infty)$
оператор~\eqref{2.66} определён и действует по формуле
$$
M{{_3S^{\nu,\mu}}}(s)=m(s)\cdot M{x^{1-\mu} f}(s),
$$
\begin{equation}\label{2.67}
 m(s)=2^{\mu-1} \left( \frac{ \cos \pi(\mu-s) - \cos \pi \nu}{ \sin \pi(\mu-s) - \sin \pi \nu}  \right)  \left( \frac{\Gamma(\frac{s}{2})\Gamma(\frac{s}{2}+\frac{1}{2}))}{\Gamma(\frac{s}{2}+\frac{1-\nu-\mu}{2}) \Gamma(\frac{s}{2}+1+\frac{\nu-\mu}{2})} \right).
\end{equation}
\end{lemma}

\begin{proof}
Из асимптотики функций Лежандра~\cite{BE1} следует, что
оператор~\eqref{2.66} определён. Записывая его как свёртку
Меллина, и последовательно применяя формулы из~\cite{Marich1}: (2.50), с.~31;
(10), с.~283; 40(1), с.~251; (5), с.~130 "--- получим
выражение~\eqref{2.67} со значением мультипликатора
$$
\frac{2^{\mu{-}2}}{\pi^2}\cdot \frac{\sin \pi(\nu{-}\mu)}{\sin \pi
\mu}\cdot\Gamma(\frac{s}{2})\Gamma(\frac{s}{2}{+}\frac{1}{2})
\Gamma(\frac{1{+}\mu{+}\nu}{2}{-}\frac{s}{2})\Gamma(\frac{\mu{-}\nu}{2}{-}\frac{s}{2}){+}
\frac{2^{\mu{-}1}}{\sin \pi \mu}\cdot
\frac{\Gamma(\frac{1{+}\mu{+}\nu}{2}{-}\frac{s}{2})\Gamma(\frac{\mu{-}\nu}{2}{-}\frac{s}{2})}
{\Gamma(1{-}\frac{s}{2})\Gamma(\frac{1}{2}{-}\frac{s}{2})}+{}
$$
$$
{}+\frac{2^{\mu-1} \cos \pi \mu }{\sin \pi \mu}\cdot \frac{
\Gamma(\frac{s}{2})\Gamma(\frac{s}{2}+\frac{1}{2})}
{\Gamma(\frac{s}{2}+\frac{1-\nu-\mu}{2})
\Gamma(\frac{s}{2}+1+\frac{\nu-\mu}{2})}.
$$
В двух последних слагаемых <<поднимем>> гамма-функции в числитель
из знаменателя, переходя от $\Gamma (z)$ к  $\Gamma (1-z).$
Получим
$$
m(s)=\frac{2^{\mu-1}}{\pi^2 \sin \pi \mu} A(s) \cdot
\Gamma(\frac{s}{2})\Gamma(\frac{s}{2}+\frac{1}{2})\Gamma(\frac{\mu+\nu+1-s}{2})\Gamma(\frac{\mu-\nu-s}{2}).
$$
Выражение $A(s)$ последовательно преобразуем по элементарным
тригонометрическим формулам
$$
A(s)= \sin \pi (\nu - \mu) +
 2 \sin \frac{\pi s}{2} \cos \frac{\pi s}{2} - 2 \cos \pi \mu
\cdot \cos \frac{\pi}{2}(s-\nu-\mu) \sin \frac{\pi}{2}(s+\nu-\mu)=
$$
$$
=\sin \pi (\nu - \mu) + \sin \pi s - \cos \pi \mu (\sin \pi \nu  +
\sin \pi (s-\mu))=
 \sin \pi \mu \cdot (\cos \pi (s- \mu) - \cos \pi \nu).
$$
Подставим это выражение в $m(s).$ Перебрасывая  две последние
гамма-функции в знаменатель, получим~\eqref{2.67}. Применение
указанных формул из~\cite{Marich1} приводит к ограничениям на
значения переменной \begin{equation}\label{2.68}{0 < \Re
(\nu+\mu)< \min\limits (1+\Re(\nu+ \mu),~ \Re
(\mu-\nu)),}\end{equation} которые в нашем случае можно ослабить.
Лемма доказана.
\end{proof}

Перейдём к доказательству теоремы~\ref{2rod}. Из
формулы~\eqref{2.67} следует, что можно перейти к пределу при $\mu
\to 1-0.$ Отсюда получим~\eqref{2.63} и~\eqref{2.65}. Аналогичными
рассуждениями доказывается~\eqref{2.64}.

Отметим, что в процессе доказательства построено ещё одно
семейство операторов преобразования типа Сонина~\eqref{2.66}.

Доказательства следующих результатов полностью аналогичны
соответствующим утверждениям для операторов Бушмана---Эрдейи
первого рода.

\begin{theorem} Справедливы формулы для норм
\begin{eqnarray}
& &  \| {_2S^{\nu}} \|_{L_2}= \max\limits (1, \sq{1+\sin \pi \nu}), \label{2.69} \\
& &  \| {_2P^{\nu}} \|_{L_2}= 1 / {\min\limits (1, \sq{1+\sin \pi
\nu})}. \label{2.610}
\end{eqnarray}
\end{theorem}

Следствие.  Оператор ${_2S^{\nu}}$ ограничен при всех $\nu.$
Оператор ${_2P^{\nu}}$ не является непрерывным тогда и только
тогда, когда $\sin \pi \nu=-1.$

\begin{theorem} Для унитарности в $L_2$ операторов ${_2S^{\nu}}$ и ${_2P^{\nu}}$ необходимо и достаточно, чтобы параметр $\nu$ был целым
числом.
\end{theorem}

\begin{theorem} Пусть $\nu=i \beta+1/2,~\beta \in \mathbb{R}.$ Тогда
\begin{equation}\label{2.6.11}{\| {_2S^{\nu}} \|_{L_2}=\sq{1+\ch \pi \beta},~\| {_2P^{\nu}} \|_{L_2}=1.}\end{equation}
\end{theorem}

\begin{theorem} Справедливы представления
\begin{eqnarray}
& & {_2S^0} f = \frac{2}{\pi} \int\limits_0^{\infty} \frac{y}{x^2-y^2}f(y)\,dy, \label{2.612} \\
& & {_2S^{-1}} f = \frac{2}{\pi} \int\limits_0^{\infty}
\frac{x}{x^2-y^2}f(y)\,dy. \label{2.613}
\end{eqnarray}
\end{theorem}

Таким образом, в этом случае оператор ${_2S^{\nu}}$ сводится к
паре известных преобразований Гильберта на полуоси~\cite{SKM}.

Теперь мы можем решить следующую важную задачу. Перейдём к
построению операторов преобразования, унитарных при  всех $\nu.$
Такие операторы определяются по формулам:
\begin{eqnarray}
& & S_U^{\nu} f = - \sin \frac{\pi \nu}{2}\  {_2S^{\nu}}f+ \cos \frac{\pi \nu}{2}\  {_1S_-^{\nu}}f, \label{2.614} \\
& & P_U^{\nu} f = - \sin \frac{\pi \nu}{2}\  {_2P^{\nu}}f+ \cos
\frac{\pi \nu}{2}\  {_1P_-^{\nu}}f. \label{2.615}
\end{eqnarray}
Для любых значений $\nu \in \mathbb{R}$ они являются линейными
комбинациями операторов преобразования Бушмана---Эрдейи 1 и 2 рода
нулевого порядка гладкости. Их можно назвать операторами
Бушмана---Эрдейи третьего рода. В интегральной форме эти операторы
имеют вид:
\begin{multline} \label{2.616}
S_U^{\nu} f = \cos \frac{\pi \nu}{2} \left(- \frac{d}{dx} \right) \int\limits_x^{\infty} P_{\nu}\lr{\frac{x}{y}} f(y)\,dy + {}\\
+ \frac{2}{\pi} \sin \frac{\pi \nu}{2} \left(  \int\limits_0^x
(x^2-y^2)^{-\frac{1}{2}}Q_{\nu}^1 \lr{\frac{x}{y}} f(y)\,dy
\right. -  \int\limits_x^{\infty}
(y^2-x^2)^{-\frac{1}{2}}\mathbb{Q}_{\nu}^1 \lr{\frac{x}{y}}
f(y)\,dy \Biggl. \Biggr),
\end{multline}
\begin{multline}\label{2.617}
P_U^{\nu} f = \cos \frac{\pi \nu}{2}  \int\limits_0^{x} P_{\nu}\lr{\frac{y}{x}} \left( \frac{d}{dy} \right) f(y)\,dy -{}  \\
  {}-\frac{2}{\pi} \sin \frac{\pi \nu}{2} \left( - \int\limits_0^x (x^2-y^2)^{-\frac{1}{2}}\mathbb{Q}_{\nu}^1\lr{\frac{y}{x}} f(y)\,dy   \right.
 - \int\limits_x^{\infty} (y^2-x^2)^{-\frac{1}{2}} Q_{\nu}^1 \lr{\frac{y}{x}} f(y)\,dy \Biggl. \Biggr).
\end{multline}

\begin{theorem}\label{2unit} Операторы~\eqref{2.614}-\eqref{2.615} или~\eqref{2.616}-\eqref{2.617} при всех $\nu$ являются унитарными, взаимно сопряжёнными и обратными в $L_2.$ Они являются
сплетающими и действуют по формулам~\eqref{275}. При этом
$S_U^{\nu}$ является оператором типа Сонина {\rm
(}Сонина---Катрахова{\rm ),} а $P_U^{\nu}$ "--- типа Пуассона {\rm
(}Пуассона---Катрахова{\rm )}.
\end{theorem}

\begin{proof}
Проверим выполнение равенства~\eqref{1716} для одного из
мультипликаторов. Аналогично случаю теоремы~\ref{2tmult} с
использованием свёртки Меллина и формул для преобразования Меллина
специальных функций получается следующая формула для
мультипликатора:
$$
M\lrs{S_U^\nu}(s)= -\sin(\frac{\pi\nu}{2})\left(\frac {-\cos\pi
s-\cos\pi\nu}{\sin\pi s-\sin\pi\nu} \right)\cdot \left(
\frac{\Gamma(\frac{s}{2})\Gamma(\frac{s}{2}+\frac{1}{2})}
{\Gamma(\frac{s}{2}-\frac{\nu}{2})\Gamma(\frac{s}{2}+\frac{\nu}{2}+\frac{1}{2})}
\right)+
$$
$$
+\cos(\frac{\pi\nu}{2})\cdot\left(
\frac{\Gamma(\frac{s}{2})\Gamma(\frac{s}{2}+\frac{1}{2})}
{\Gamma(\frac{s}{2}-\frac{\nu}{2})
\Gamma(\frac{s}{2}+\frac{\nu}{2}+\frac{1}{2})} \right)=
$$
$$
=\left( -\sin(\frac{\pi\nu}{2}) \left(\frac {-\cos\pi
s-\cos\pi\nu}{\sin\pi s-\sin\pi\nu}\right)+ \cos(\frac{\pi\nu}{2})
\right)\cdot\left(
\frac{\Gamma(\frac{s}{2})\Gamma(\frac{s}{2}+\frac{1}{2})}
{\Gamma(\frac{s}{2}-\frac{\nu}{2})
\Gamma(\frac{s}{2}+\frac{\nu}{2}+\frac{1}{2})} \right).
$$
Полностью вычисления проведены в~\cite{S66}. Далее рассматриваем в
соответствии с теоремой~\ref{1tMel} с учётом выкладок из
теоремы~\ref{2tmult} величину
$$
\left|M\lrs{S_U^\nu}(i u+\frac{1}{2})\right|=
\left|-\sin(\frac{\pi\nu}{2}) \left(\frac {-\cos\pi (i
u+\frac{1}{2})-\cos\pi\nu}{\sin\pi (i
u+\frac{1}{2})-\sin\pi\nu}\right)+ \cos(\frac{\pi\nu}{2})
\right|\cdot \sqrt{\frac{\cos\pi u -\sin\pi \nu}{\cos\pi u}}=
$$
$$
=\left|
\frac{\sin(\frac{\pi\nu}{2})-\cos\left(\frac{\pi\nu}{2}+\pi i
u\right)} {\sin\pi\nu-\cos\pi i u} \right| \cdot
\sqrt{\frac{\cos\pi u -\sin\pi \nu}{\cos\pi u}}=
\left| \frac{\sin\frac{\pi\nu}{2}-\cos\frac{\pi\nu}{2}\ch\pi u +
i\sin\frac{\pi\nu}{2}\sh\pi u} {\sqrt{\ch\pi u(\ch\pi-\sin\pi
\nu)}} \right|.
$$
Вычисляя модуль и заменяя затем тригонометрический и
гиперболический синусы на косинусы, получим:
$$
\left|M\lrs{S_U^\nu}(i u+\frac{1}{2})\right|=
\sqrt{\frac{\left(\sin\frac{\pi\nu}{2}-\cos\frac{\pi\nu}{2}\ch\pi
u\right)^2+ \left(\sin\frac{\pi\nu}{2}\sh\pi u\right)^2} {\ch\pi
u(\ch\pi-\sin\pi \nu)}}=
\sqrt{\frac{\ch^2\pi u-\sin\pi\nu\ch\pi u}{\ch\pi
u(\ch\pi-\sin\pi \nu)}} =1.
$$
Унитарность этим доказана. Взаимная сопряжённость следует из
определений~\eqref{2.614}-\eqref{2.615}, если рассматривать
операторы  как расширения с множества финитных функций.
Следовательно, эти унитарные операторы и взаимно обратны. То, что
они являются ОП типа Сонина и Пуассона, вытекает из проверки
условий для мультипликаторов теоремы~\ref{1tMel}. Теорема
доказана.
\end{proof}

Этим результатом завершается история построения унитарных ОП типа
Сонина и Пуассона. Унитарные операторы преобразования тесно
связаны с унитарностью оператора рассеяния в задачах квантовой
механики~\cite{AM, Fad1,Fad2, ShSa}.

Интересно рассмотреть частный случай, который вытекает из
теоремы~\ref{2unit} при $\nu=1.$ Получаем пару очень простых
операторов
\begin{equation}
\label{282} B_{0+}^{1,1}f=f(x)-\frac{1}{x}\int\limits_0^x
f(y)\,dy,\  B_{-}^{1,1}f=f(x)-\int\limits_x^\infty
\frac{f(y)}{y}\,dy,
\end{equation}
связанных со знаменитыми операторами Харди
\begin{equation}
\label{283} H_1f=\frac{1}{x}\int\limits_0^x f(y)\,dy,
H_2f=\int\limits_x^\infty \frac{f(y)}{y}\,dy.
\end{equation}
По поводу теории неравенств Харди см.~\cite{OpKu}. Из наших
результатов следует

\begin{theorem} \label{2tHar}
Операторы~\eqref{282} образуют пару взаимнообратных унитарных в
$L_2(0,\infty)$ операторов. Они сплетают как ОП
$\dfrac{d^2}{dx^2}$ и $\dfrac{d^2}{dx^2}-\dfrac{2}{x^2}.$
\end{theorem}

Как следует из~\eqref{283}, операторы Бушмана---Эрдейи могут
рассматриваться как \textit{обобщения операторов Харди}, а
неравенства для их норм являются \textit{определёнными обобщениями
неравенств Харди}, что позволяет взглянуть на этот класс
операторов под новым интересным углом зрения. Кроме того, можно
показать, что операторы~\eqref{282} являются преобразованиями Кэли
от симметричных операторов $\pm 2i (xf(x))$ при соответствующем
выборе областей определения. Их спектром является единичная
окружность. В~\cite{S6, S66} эти вопросы рассмотрены и для
пространств со степенным весом.

Результат об унитарности из теоремы~\ref{2tHar}   рассматривался
Куфнером, Перссоном и Малиграндой~\cite{KMP}, давшими его
элементарное доказательство и приложения. Теорема~\ref{2unit}
позволяет выписать ещё несколько пар унитарных в $L_2(0,\infty)$
операторов очень простого вида, которые являются частными случаями
операторов Бушмана---Эрдейи при целых $\nu$ (см.~\cite{S6, S66}):
\begin{align*}
\label{284} U_3f&= f+\int\limits_0^x f(y)\,\frac{dy}{y}, &  U_4f&=
f+\frac{1}{x}\int\limits_x^\infty f(y)\,dy,\\\nonumber U_5f&=
f+3x\int\limits_0^x f(y)\,\frac{dy}{y^2}, &  U_6f&=
f-\frac{3}{x^2}\int\limits_0^x y f(y)\,dy,\\\nonumber
U_7f&=f+\frac{3}{x^2}\int\limits_x^\infty y f(y)\,dy, & U_8f&=f-3x
\int\limits_x^\infty f(y)\frac{dy}{y^2},\\\nonumber
U_9f&=f+\frac{1}{2}\int\limits_0^x
\left(\frac{15x^2}{y^3}-\frac{3}{y}\right)f(y)\,dy, &
U_{10}f&=f+\frac{1}{2}\int\limits_x^\infty
\left(\frac{15y^2}{x^3}-\frac{3}{x}\right)f(y)\,dy.\\\nonumber
\end{align*}

Этот перечень можно продолжить дальше. Интересно отметить, что
приведённые выше примеры противоречат утверждению, высказанному
В.\,А.~Марченко в~\cite{Mar9}: <<{\it Но среди вольтерровских
операторов только единичный оператор унитарен}>> "--- имеются в
виду интегральные операторы Вольтерра второго рода, форму которых
имеют большинство классических операторов преобразования.

ОП в форме подобной~\eqref{2.616}-\eqref{2.617}, но только с
ядрами, выражающимися через общую гипергеометрическую функцию
Гаусса, были впервые построены в 1980~г. В.\,В.~Катраховым, в его
работах ядра выражались через гипергеометрическую функцию Гаусса.
Поэтому можно предложить для них названия: операторы
преобразования Сонина---Катрахова и Пуассона---Катрахова. Мы
получили для них выражения  через функции Лежандра первого и
второго родов, кроме того их удаётся включить в общую схему
построения операторов преобразования композиционным методом,
изложенным ниже в главе~\ref{ch6}. При этом основными становятся
наиболее простые  формулы факторизации
вида~\eqref{2.614}-\eqref{2.615}. На этом пути построение подобных
операторов перестаёт быть специальным искусным приёмом, а
встраивается в общую методику построения целых классов подобных
операторов преобразования композиционным методом.

\section[Приложения операторов преобразования Бушмана---Эрдейи, Сонина---Катрахова и Пуассона---Катрахова к дифференциальным уравнениям с особенностями в
коэффициентах]{Приложения операторов преобразования Бушмана---Эрдейи, Сонина---Катрахова и Пуассона---Катрахова к дифференциальным уравнениям с особенностями в
коэффициентах}\label{sec8}
\sectionmarknum{Приложения операторов преобразования к дифференциальным уравнениям с особенностями в
коэффициентах}

\subsection{Приложения операторов преобразования Бушмана---Эрдейи к  задачам для уравнения Эйлера---Пуассона---Дарбу и лемме
Копсона}\label{sec8.1}

   Рассмотрим теперь приложения операторов преобразования Бушмана---Эрдейи первого рода и нулевого порядка гладкости к обобщениям и уточнению леммы Копсона.

\begin{theorem} Рассмотрим задачу Дирихле в четверти плоскости для уравнения Эйлера---Пуассона---Дарбу в условиях леммы Копсона~\eqref{C1}-\eqref{C2}. Тогда между данными задачи Дирихле справедливы соотношения, выражающиеся через операторы  Бушмана---Эрдейи первого рода и нулевого порядка гладкости{\rm :}
\begin{equation}
\frac{c_\beta}{x^\beta} E_{0+}^{-\alpha,1-\beta}(y^{\alpha+\beta+1}f(y))=
\frac{c_\alpha}{x^\alpha} E_{0+}^{-\beta,1-\alpha}(y^{\alpha+\beta+1}g(y)), c_\beta=2\Gamma(\beta+1/2),
\end{equation}
\begin{equation}
\frac{c_\beta}{x^\beta} {_1 P_{0+}^{-\alpha}}
I_{0+}^{\beta}(y^{\alpha+\beta+1}f(y))= \frac{c_\alpha}{x^\alpha}
{_1 P_{0+}^{-\beta}} I_{0+}^{\alpha}(y^{\alpha+\beta+1}g(y)).
\end{equation}
\end{theorem}

Данные формулы являются прямыми следствиями теоремы~\ref{2fact1}.
C другой стороны, применяя теперь полученную выше
теорему~\ref{2factBE}, в которой операторы Бушмана---Эрдейи
нулевого порядка гладкости были профакторизованы через дробные
интегралы Римана---Лиувилля и Эрдейи---Кобера, мы получаем такой
результат.

\begin{theorem} Между данными задачи Дирихле в условиях леммы Копсона справедливы соотношения, выражающиеся через операторы Эрдейи---Кобера{\rm :}
\begin{equation}
x^{\alpha+\beta+1}f(x)=\frac{c_\alpha}{c_\beta} I_{0+;2;-1/2}^{\alpha-\beta}(y^{\alpha+\beta+1}g(y)).
\end{equation}
\end{theorem}

Последнее соотношение уточняет соответствующий результат из
оригинальной работы Копсона, который, по-видимому, содержит
неточности в коэффициентах.

Связь между данными на осях, выраженная через операторы
Бушмана---Эрдейи  нулевого порядка гладкости, позволяет установить
дополнительные результаты.

\begin{theorem} Рассмотрим пространство с весом $L_{2,x^k}(0,\infty),$ и пусть $\alpha,\beta$ "--- целые числа. Тогда для весовых норм данных Дирихле справедливы соотношения
\begin{equation}
c_\beta \|I_{0+}^{\beta}(y^{\alpha+\beta+1}f(y))\|=c_\alpha \|I_{0+}^{\alpha}(y^{\alpha+\beta+1}g(y))\|.
\end{equation}
\end{theorem}

Этот результат следует из условия унитарности операторов
Бушмана---Эрдейи  нулевого порядка гладкости в
теореме~\ref{2tunit}. Кроме того, в этом случае можно из той же
теоремы выразить данные не через оператор Эрдейи---Кобера, а через
обратный оператор Бушмана---Эрдейи  нулевого порядка гладкости,
получив новое соотношение. Для произвольных, не являющихся целыми
значений параметров, из теоремы~\ref{2tmult} также получается
соотношение между весовыми нормами данных Дирихле.

   Рассмотрим приложения операторов преобразования рассматриваемых нами классов к постановкам задачи Коши для уравнения Эйлера---Пуассона---Дарбу (ЭПД).

Рассмотрим уравнение ЭПД в полупространстве
$$
B_{\alpha,\, t} u(t,x)= \frac{ \pd^2 \, u }{\pd t^2} + \frac{2 \alpha+1}{t} \frac{\pd u}{\pd t}=\Delta_x u+F(t, x),
$$
где $t>0,~x \in \mathbb{R}^n.$ Дадим  описание процедуры,
позволяющей получать различные постановки начальных условий при
$t=0$ единым методом. Образуем по формулам~\eqref{275} операторы
преобразования $X_{\alpha, \, t}$ и $Y_{\alpha, \, t}.$
Предположим, что существуют выражения $X_{\alpha, \, t} u=v(t,x),$
$X_{\alpha, \, t} F=G(t,x).$ Пусть обычная (несингулярная) задача
Коши \begin{equation}\label{2.7.28}{\frac{ \pd^2 \, v }{\pd t^2}
=\Delta_x v+G,~ v|_{t=0}=\varphi (x),~ v'_t|_{t=0}=\psi
(x)}\end{equation} корректно разрешима в полупространстве. Тогда
получаем следующие начальные условия для уравнения ЭПД:
\begin{equation}\label{2.7.29}{X_{\alpha} u|_{t=0}=a(x),~(X_{\alpha}
u)'|_{t=0}=b(x).}\end{equation} При этом различному выбору
операторов преобразования  $X_{\alpha, t}$ (операторы Сонина,
Бушмана---Эрдейи, Бушмана---Эрдейи нулевого порядка гладкости 1
или 2 рода, унитарные операторы третьего рода~\eqref{2.616},
обобщённые операторы Бушмана---Эрдейи) будут соответствовать
различные начальные условия. Следуя изложенной методике в каждом
конкретном случае их можно привести к более простым аналитическим
формулам. При этом с использованием интегральных ОП различных
типов для каждого конкретного ОП будут получатся некоторые
нелокальные начальные условия, в том числе и с возможностью
рассмотрения решений с особенностями.

Данная схема обобщается на дифференциальные уравнения с большим
числом переменных, по которым могут действовать операторы Бесселя
с различными параметрами, а также уравнения других типов.
Применение операторов преобразований позволяет сводить сингулярные
(или иначе вырождающиеся) уравнения с операторами Бесселя по одной
или нескольким переменным (уравнения ЭПД, сингулярное уравнение
теплопроводности, $B$-эллиптические уравнения по определению
И.\,А.~Кипри\-янова, уравнения обобщённой осесимметрической теории
потенциала "--- теории GASPT ({\it Generalized  Axially  Symmetric
Potential  Theory}) "--- А.~Вайнстейна и др.) к несингулярным. При
этом априорные оценки для сингулярного случая получаются как
следствия соответствующих априорных оценок для регулярных
уравнений, если только удалось оценить сами операторы
преобразования  в нужных функциональных пространствах.
Значительное число подобных оценок было приведено выше.

\subsection{Приложения операторов преобразования Бушмана---Эрдейи, Сонина---Катрахова и Пуассона---Катрахова к установлению формул связи между решениями дифференциальных
уравнений}\label{sec8.2}

Начнём с простого, но важного замечания. Любую из рассмотренных в этой главе пар операторов преобразования, сплетающих операторы Бесселя и вторую производную, можно использовать для установления связей между решениями дифференциального уравнения с особенностями в коэффициентах с операторами Бесселя вида
$$
\sum\limits a_k B_{\nu_k,x_k} u(x) =f(x)
$$
и невозмущённого уравнения с постоянными коэффициентами
$$
\sum\limits a_k \frac{\pr^2 v(x)}{\pr x_k^2} =g(x).
$$
При этом если пары взаимообратных ОП действуют по каждой переменной по формулам
\begin{equation}
S_\nu B_\nu=D^2 S_\nu,\quad  P_\nu D^2=B_\nu P_\nu,
\end{equation}
то решения возмущённого и невозмущённого уравнений связаны
соотношениями $$ u(x)=\prod\limits_k S_{\nu_k} v(x),\ \
v(x)=\prod\limits_k P_{\nu_k} u(x)
$$ (операторы типа Сонина и
Пуассона). При этом результаты об ограниченности, оценках норм,
унитарности операторов преобразований приводят автоматически к
соответствующим утверждениям для пар решений дифференциальных
уравнений. Мы ограничимся этой схемой, не выписывая
соответствующих формул связи и оценок для решений дифференциальных
уравнений с особенностями в коэффициентах.

Рассмотрим подробно один пример, позволяющий использовать
построенные ОП различных классов и их оценки для построения
решений одного из нелинейных уравнений. В работах А.\,В.~Бицадзе и
В.\,И.~Пашковского~\cite{Bitz1,Bitz12}  указано, что в
математической физике ряд задач сводятся к решению нелинейных
уравнений Максвелла---Эйнштейна вида
\begin{equation}\label{ur1}
\Delta u +\frac{1}{x}u_x-\frac{1}{u}(1-\frac{u^2}{A^2-u^2})(u_x^2+u_y^2)=0,
\end{equation}
\begin{equation}\label{ur2}
\Delta u +\frac{1}{x}u_x-\frac{1}{u}(u_x^2+u_y^2)=0.
\end{equation}

Используя результаты этой работы А.\,В.~Бицадзе и развитую нами
технику ОП, получаем следующее приложение рассмотренных классов ОП
к нелинейным уравнениям математической физики
Максвелла---Эйнштейна.

\begin{theorem}
Пусть $P$ "--- произвольный оператор преобразования типа Пуассона,
удовлетворяющий сплетающему свойству на гладких функциях
\begin{equation}
P\ D^2=(\frac{d^2}{d x^2}+\frac{1}{x}\frac{d}{dx})\ P,
\end{equation}
$g(x,y)$ "--- произвольная гармоническая функция. Тогда функция
$u_1(x,y)=\dfrac{A}{\ch(a P_x g(x,y))}$ является решением
нелинейного уравнения~\eqref{ur1}, а функция  $u_2(x,y)=\exp(b P_x
g(x,y))$ является решением нелинейного уравнения~\eqref{ur2},
$a,b$ "--- произвольные постоянные.
\end{theorem}

Нами построены различные классы ОП типа Пуассона: Бушмана---Эрдейи
первого, второго родов, нулевого порядка гладкости,
Сонина---Катрахова и Пуассона---Катрахова и т.~д. Теперь мы можем
использовать их в приведённой теореме для получения представлений
решений нелинейных уравнений Максвелла---Эйнштейна через
гармонические функции.

\subsection[Приложения операторов преобразования  Сонина---Катрахова и Пуассона---Катрахова к решению некоторых интегродифференциальных
уравнений]{\hspace{-5pt}Приложения операторов преобразования
Сонина---Катрахова и Пуассона---Катрахова к решению некоторых
интегродифференциальных уравнений}\label{sec8.3}

Рассмотрим приложения унитарных операторов Сонина---Катрахова и
Пуассона---Катрахова к решению соответствующих
интегродифференциальных уравнений.

\begin{theorem}
Пусть функции $f(x),\ g(x) \in L_2(0,\infty)$ и непрерывно
дифференцируемы на полуоси. Тогда следующие
интегро-дифференциальные уравнения взаимно обратны и решаются по
приведённым формулам{\rm :}
\begin{multline}
g(x) = \cos \frac{\pi \nu}{2} \left(- \frac{d}{dx} \right) \int\limits_x^{\infty} P_{\nu}\lr{\frac{x}{y}} f(y)\,dy +  {}\\
{}+ \frac{2}{\pi} \sin \frac{\pi \nu}{2} \left(  \int\limits_0^x
(x^2-y^2)^{-\frac{1}{2}}Q_{\nu}^1 \lr{\frac{x}{y}} f(y)\,dy
\right.
-  \int\limits_x^{\infty} (y^2-x^2)^{-\frac{1}{2}}\mathbb{Q}_{\nu}^1 \lr{\frac{x}{y}} f(y)\,dy \Biggl. \Biggr),  \\
\end{multline}
\begin{multline}
f(x) = \cos \frac{\pi \nu}{2}  \int\limits_0^{x} P_{\nu}\lr{\frac{y}{x}} \left( \frac{d}{dy} \right) g(y)\,dy - {} \\
{}-\frac{2}{\pi} \sin \frac{\pi \nu}{2} \left( - \int\limits_0^x
(x^2-y^2)^{-\frac{1}{2}}\mathbb{Q}_{\nu}^1\lr{\frac{y}{x}}
g(y)\,dy   \right. - \int\limits_x^{\infty}
(y^2-x^2)^{-\frac{1}{2}} Q_{\nu}^1 \lr{\frac{y}{x}} g(y)\,dy
\Biggl. \Biggr).
\end{multline}
При этом в указанном пространстве нормы решений и правых частей равны.
\end{theorem}

Это приложение теоремы об унитарности ОП Сонина---Катрахова и
Пуассона---Катрахова. Отметим, что при специальных значениях
параметра $\nu,$ при которых функции Лежандра выражаются через
более простые функции, получается список конкретных
интегро-дифференциальных уравнений, для которых получены явные
решения с  оценками их норм. Мы не будем приводить здесь этот
список.

\section{Приложения операторов преобразования Бушмана---Эрдейи к установлению эквивалентности норм пространств И.\,А.~Киприянова и весовых пространств
С.\,Л.~Соболева}\label{sec9}
\sectionmarknum{\sП\sр\sи\sл\sо\sж\sе\sн\sи\sя\s{ }\sо\sп\sе\sр\sа\sт\sо\sр\sо\sв\s{ }\sБ\sу\sш\sм\sа\sн\sа\s---\sЭ\sр\sд\sе\sй\sи\s{ }\sк\s{
}\sэ\sк\sв\sи\sв\sа\sл\sе\sн\sт\sн\sо\sс\sт\sи\s{ }\sн\sо\sр\sм\s{
}\sп\sр\sо\sс\sт\sр\sа\sн\sс\sт\sв\s{
}\sК\sи\sп\sр\sи\sя\sн\sо\sв\sа\s{ }\sи\s{
}\sв\sе\sс\sо\sв\sы\sх\s{ }\sп\sр\sо\sс\sт\sр\sа\sн\sс\sт\sв{
}\sС\sо\sб\sо\sл\sе\sв\sа}

Пространства И.\,А.~Киприянова, как было показано в его работах и
в работах представителей его научной школы, идеально подходят для
изучения $B$-эллиптических уравнений с частными
производными~\cite{Kip1,Kip2}. Поэтому свойства этих пространств
играют существенную роль при изучении дифференциальных уравнений с
особенностями в коэффициентах. В данном пункте доказано, что в
одномерном случае нормы в пространствах И.\,А.~Киприянова
эквивалентны нормам в весовых пространствах С.\,Л.~Соболева. Этот
результат справедлив и для некоторых модельных областей в
многомерных пространствах. При этом соответствующие результаты по
существу являются переформулировками результатов предыдущего
параграфа об условиях ограниченности и унитарности в пространствах
Лебега на полуоси операторов преобразования Бушмана---Эрдейи
нулевого порядка гладкости. Таким образом, указанные результаты об
условиях ограниченности и унитарности в пространствах Лебега на
полуоси операторов преобразования Бушмана---Эрдейи нулевого
порядка гладкости имеют важное значение для теории уравнений с
частными производными с операторами Бесселя, находя в этой области
свои как стандартные, так и неожиданные приложения. Получаемые
результаты относятся также к энергетическим пространствам для
соответствующих дифференциальных уравнений.

В данном пункте мы приведём список основных результатов без
доказательств. Доказательства по существу вытекают из приведённых
ранее свойств операторов Бушмана---Эрдейи, однако аккуратные
определения и выкладки требуют немалого объёма текста.

 Определим множество функций $\mathfrak{D}(0, \infty).$ Если $f(x) \in \mathfrak{D}(0, \infty),$ то $f(x) \in C^{\infty}(0, \infty),~ f(x)$ "--- финитна на бесконечности. На этом множестве функций введём полунормы
\begin{eqnarray}
& & \|f\|_{h_2^{\alpha}}=\|I_-^{\alpha}f\|_{L_2(0, \infty)}, \label{2.71} \\
& &  \|f\|_{\widehat{h}_2^{\alpha}}=\|x^{\alpha}
(-\frac{1}{x}\frac{d}{dx})^{\alpha}f\|_{L_2(0, \infty)},
\label{2.72}
\end{eqnarray}
где $I_-^{\alpha}$ "--- дробная производная Римана---Лиувилля,
оператор в~\eqref{2.72} определяется по формуле
\begin{equation}\label{2.73}{(-\frac{1}{x}\frac{d}{dx})^{\beta}=2^{\beta}I_{-; \, 2,
\,0}^{-\beta}x^{-2 \beta},}\end{equation} $I_{-; 2, \,
0}^{-\beta}$ "--- оператор Эрдейи---Кобера, $\alpha$ "---
произвольное действительное число. При $\beta = n \in
\mathbb{N}_0$ выражение~\eqref{2.73} понимается в обычном смысле,
что согласуется с определениями~\eqref{162}--\eqref{1.16}
главы~\ref{ch1}.

\begin{lemma}
 Пусть $f(x) \in \mathfrak{D}(0, \infty).$ Тогда справедливы тождества{\rm :}
\begin{eqnarray}
& & I_-^{\alpha}f={_1S_-^{\alpha-1}} {x^{\alpha} (-\frac{1}{x}\frac{d}{dx})^{\alpha}} f, \label{2.74} \\
& & x^{\alpha}
(-\frac{1}{x}\frac{d}{dx})^{\alpha}f={_1P_-^{\alpha-1}}
I_-^{\alpha}f. \label{2.75}
\end{eqnarray}
\end{lemma}

Таким образом, операторы Бушмана---Эрдейи нулевого порядка
гладкости первого рода осуществляют связь между дифференциальными
операторами (при $\alpha \in \mathbb{N}$) из определений
полунорм~\eqref{2.71} и~\eqref{2.72}.

\begin{lemma}
 Пусть $f(x) \in \mathfrak{D}(0, \infty).$ Тогда справедливы неравенства между полунормами
\begin{eqnarray}
& &  \|f\|_{h_2^{\alpha}} \leq \max\limits (1, \sq{1+\sin \pi \alpha}) \|f\|_{\widehat{h}_2^{\alpha}}, \label{2.77}\\
& & \|f\|_{\widehat{h}_2^{\alpha}} \leq \frac{1}{\min\limits (1,
\sq{1+\sin \pi \alpha})} \|f\|_{h_2^{\alpha}}, \label{2.78}
\end{eqnarray}
где $\alpha$ "--- любое действительное число, $\alpha \neq
-\dfrac{1}{2}+2k,~k \in \mathbb{Z}.$
\end{lemma}

Постоянные в неравенствах~\eqref{2.77}-\eqref{2.78} не меньше
единицы, что будет далее использовано. В случае $\sin \pi \alpha =
-1 $ или $\alpha = -\dfrac{1}{2}+2k,~k \in \mathbb{Z},$
оценка~\eqref{2.78} не имеет места.

Приведённые леммы "--- это переформулировки результатов о формулах
факторизации и оценках норм для операторов Бушмана---Эрдейи
нулевого порядка гладкости.

Введём на $\mathfrak{D} (0, \infty )$ соболевскую норму
\begin{equation}\label{2.79}{\|f\|_{W_2^{\alpha}}=\|f\|_{L_2 (0, \infty)}+\|f\|_{h_2^{\alpha}}.}\end{equation}
Введём также другую норму
\begin{equation}\label{2.710}{\|f\|_{\widehat{W}_2^{\alpha}}=\|f\|_{L_2 (0,
\infty)}+\|f\|_{\widehat{h}_2^{\alpha}}}.\end{equation}
Пространства $W_2^{\alpha},~ \widehat{W}_2^{\alpha}$ определим как
замыкания $\mathfrak{D}(0, \infty)$ по нормам~\eqref{2.79}
и~\eqref{2.710} соответственно.

\begin{theorem}~\par
\begin{enumerate}
\item[а)] При всех $\alpha \in \mathbb{R}$ пространство
$\widehat{W}_2^{\alpha}$ непрерывно вложено в $W_2^{\alpha},$
причём
\begin{equation}\label{2.711}{\|f\|_{W_2^{\alpha}}\leq A_1
\|f\|_{\widehat{W}_2^{\alpha}},}\end{equation} где
$A_1=\max\limits (1, \sq{1+\sin \pi \alpha}).$

\item[б)] Пусть $\sin \pi \alpha \neq -1$ или $\alpha \neq
-\dfrac{1}{2} + 2k, ~ k \in \mathbb{Z}.  $ Тогда справедливо
обратное вложение $W_2^{\alpha}$  в $\widehat{W}_2^{\alpha},$
причём
\begin{equation}\label{2.712}{\|f\|_{\widehat{W}_2^{\alpha}}\leq A_2
\|f\|_{W_2^{\alpha}},}\end{equation} где $A_2 =1/  \min\limits (1,
\sq{1+\sin \pi \alpha}).$

\item[в)] Пусть $\sin \pi \alpha \neq -1,$ тогда пространства
$W_2^{\alpha}$  и $\widehat{W}_2^{\alpha}$ изоморфны, а их нормы
эквивалентны.

\item[г)] Константы в неравенствах
вложений~\eqref{2.711}-\eqref{2.712} точные.
\end{enumerate}
\end{theorem}

Эта теорема фактически является следствием результатов по
ограниченности операторов Бушмана---Эрдейи нулевого порядка
гладкости в $L_2,$ а именно теоремы~\ref{2tmult}. В свою очередь,
из теоремы об унитарности этих операторов~\ref{2tunit} вытекает
результат об эквивалентной нормировке рассматриваемых нами
вариантов пространств С.\,Л.~Соболева.

\begin{theorem} Нормы
\begin{eqnarray}
& & \|f\|_{W_2^{\alpha}} = \sum\limits_{j=0}^s \| \mathfrak{D}_-^j f\|_{L_2}, \label{2.713} \\
& & \|f\|_{\widehat{W}_2^{\alpha}}=\sum\limits_{j=0}^s \|
x^j(-\frac{1}{x}\frac{d}{dx})^j f \|_{L_2} \label{2.714}
\end{eqnarray}
задают эквивалентные нормировки в пространстве Соболева при целых
$s \in \mathbb{Z}.$ Кроме того, каждое слагаемое в~\eqref{2.713}
тождественно равно соответствующему слагаемому в~\eqref{2.714} с
тем же индексом $j.$
\end{theorem}

Как уже было отмечено, И.\,А.~Киприянов ввёл в~\cite{Kip1,Kip2}
шкалу пространств, которые оказали существенное влияние на теорию
уравнений в частных производных с оператором Бесселя по одной или
нескольким переменным. Эти пространства можно определить следующим
образом. Рассмотрим подмножество чётных функций на
$\mathfrak{D}(0, \infty),$ у которых все производные нечётного
порядка равны нулю при $x=0.$ Обозначим это множество
$\mathfrak{D}_c (0, \infty)$ и введём на нём норму
\begin{equation}\label{2.715}{\|f\|_{\widetilde{W}_{2, k}^s} =
\|f\|_{L_{2, k}}+\|B_k^{\frac{s}{2}}\|_{L_{2, k}}},\end{equation}
где $s$ "--- чётное натуральное число, $B^{s/2}_k$ "--- итерация
оператора Бесселя. Пространство И.\,А.~Киприянова при чётных $s$
определяется как замыкание $\mathfrak{D}_c (0, \infty)$ по
норме~\eqref{2.715}. Известно, что эквивалентная~\eqref{2.715}
норма может быть задана по формуле~\cite{Kip2}
\begin{equation}\label{2.716}{\|f\|_{\widetilde{W}_{2, k}^s} =
\|f\|_{L_{2, k}}+\|x^s(-\frac{1}{x}\frac{d}{dx})^s f\|_{L_{2,
k}}}.\end{equation} Это позволяет доопределить норму в
$\widetilde{W}_{2, \, k}^s$ для всех $s.$ Отметим, что по существу
этот подход совпадает с одним из принятых в~\cite{Kip2}, другой
подход основан на использовании преобразования Фурье---Бесселя.
Далее будем считать, что $\widetilde{W}_{2, k}^s$ нормируется по
формуле~\eqref{2.716}.

Введём весовую соболевскую норму
\begin{equation}\label{2.717}{\|f\|_{W_{2, k}^s} = \|f\|_{L_{2, k}}+\|\mathfrak{D}_-^s f\|_{L_{2, k}}}\end{equation}
и пространство $W_{2, \, k}^s$ как замыкание $\mathfrak{D}_c (0,
\infty)$ по этой норме.

\begin{theorem} \label{2tvloz1}~\par
\begin{enumerate}
\item[а)] Пусть $k \neq -n, ~ n \in \mathbb{N}.$ Тогда
пространство $\widetilde{W}_{2, \, k}^s$ непрерывно вложено в
$W_{2, \, k}^s,$ причём существует постоянная $A_3>0$ такая, что
\begin{equation}\label{2.718}{\|f\|_{W_{2, k}^s}\leq A_3 \|f\|_{\widetilde{W}_{2,
k}^s}.}\end{equation}

\item[б)] Пусть $k+s \neq -2m_1-1, ~ k-s \neq -2m_2-2, ~ m_1 \in
\mathbb{N}_0, ~ m_2 \in \mathbb{N}_0.$ Тогда справедливо обратное
вложение $W_{2, \, k}^s$ в $\widetilde{W}_{2, \, k}^s,$ причём
существует постоянная $A_4>0,$ такая, что
\begin{equation}\label{2.719}{\|f\|_{\widetilde{W}_{2, k}^s}\leq A_4 \|f\|_{W_{2,
k}^s}.}\end{equation}

\item[в)] Если указанные условия не выполняются, то
соответствующие вложения не имеют места.
\end{enumerate}
\end{theorem}

\begin{corollary}
Пусть выполнены условия{\rm :} $k \neq -n,$ $n \in \mathbb{N};$
$k+s \neq -2m_1-1,$ $m_1 \in \mathbb{N}_0;$ $k-s \neq -2m_2-2,$
$m_2 \in \mathbb{N}_0.$ Тогда пространства И.\,А.~Киприянова можно
определить как замыкание $\mathfrak{D}_c (0, \infty)$ по весовой
соболевской норме~\eqref{2.717}.
\end{corollary}

\begin{corollary}  Точные значения постоянных в неравенствах вложения~\eqref{2.718}-\eqref{2.719} есть
$$
A_3 = \max\limits (1, \|{_1S_-^{s-1}} \| _ {L_{2, k}}), ~
A_4=\max\limits(1, \|{_1P_-^{s-1}}\|_{L_{2, k}}).
$$
\end{corollary}

Очевидно, что приведённая теорема и следствия из неё вытекают из
приведённых выше результатов для операторов Бушмана---Эрдейи.
Отметим, что  нормы операторов Бушмана---Эрдейи нулевого порядка
гладкости в $L_{2, \, k}$  дают значения точных постоянных в
неравенствах вложения~\eqref{2.718}-\eqref{2.719}. Оценки норм
операторов Бушмана---Эрдейи в банаховых пространствах $L_{p,
\alpha}$ позволяют рассматривать вложения для соответствующих
функциональных пространств.

Неравенство для полунорм, приводящее к вложению~\eqref{2.718} ($s$
"--- целое число), получено в~\cite{Lis}. Вложения типа полученных
в теореме~\ref{2tvloz1} ранее изучались в~\cite{Lei1,Lei2}. В
последней работе рассматривался случай $k>-1/2,$ $s \in
\mathbb{N},$ пространства $W_{p,\,k}^s.$ Мы рассматриваем
гильбертовы пространства $W_{2,\,k}^s,$ $k$ и $s$ "--- любые
действительные числа. Кроме того, в теореме~\ref{2tvloz1} уточнены
условия отсутствия вложений из~\cite{Lei1,Lei2}, которые содержали
ошибки (в качестве контрпримеров для некоторых вложений
приводилась функция синуса, которая считалась чётной).   Отметим,
что в теореме~\ref{2tvloz1} фактически установлены более точные
неравенства между соответствующими полунормами, чем в
предшествующих работах~\cite{Lei1,Lei2}. Это стало возможным
благодаря применению подробно изученных выше ОП Бушмана---Эрдейи.

Перейдём к рассмотрению правосторонних сплетающих
операторов~\eqref{2BE01}--\eqref{2BE04}. Мы покажем, что в общем
случае они осуществляют изоморфизм пространства С.\,Л.~Соболева и
И.\,А.~Киприянова. Определим пространства $H^{2s},$
$H_{\alpha}^{2s}$ и $K_{\alpha}^{2s}$ как замыкания множества
функций $\mathfrak{D} (0, \infty)$ по нормам
\begin{eqnarray}
& & \|f\|_{H^{2s}} = \|f\|_{L_2}+\|I_-^{2s} f\|_{L_2}, \label{2.720} \\
& & \|f\|_{H^{2s}_{\alpha}} = \|f\|_{L_{2, {\alpha}}}+\|I_-^{2s} f\|_{L_{2, {\alpha}}}, \label{2.721} \\
& & \|f\|_{K^{2s}_{\alpha}} = \|f\|_{L_{2,
{\alpha}}}+\|B_{\alpha}^s f\|_{L_{2, {\alpha}}}, \label{2.722}
\end{eqnarray}
$s$ "--- натуральное число, $\alpha \in \mathbb{R}.$ Определим
также пару операторов типа~\eqref{276}
\begin{equation}\label{2.723}{{_1X_-^{\alpha}}={_1S_-^{\alpha-\frac{1}{2}}}
x^{\alpha+\frac{1}{2}}, ~ {_1Y_-^{\alpha}}=
x^{-(\alpha+\frac{1}{2})} {_1P_-^{\alpha-\frac{1}{2}}}.
}\end{equation}

\begin{theorem} \label{2tvloz2}
Пусть $\alpha \in \mathbb{R},$ $s \in \mathbb{N}.$ Тогда оператор
${_1X^2_-}$ действует непрерывно из $H^{2s}_{\alpha}$ в
$K^{2s}_{\alpha},$ причём
\begin{equation}\label{2.724}{\|{_1X_-^{\alpha}}f\|_{H^{2s}_{\alpha}} \leq A_5
\|f\|_{K^{2s}_{\alpha}},}\end{equation} где
$A_5=\|{_1X_-^{\alpha}}\|_{H^{2s}_{\alpha} \to
K^{2s}_{\alpha}}=\|{_1S_-^{\alpha-\frac{1}{2}}}\|_{L_2}=
\max\limits (1, \, \sq{1+\cos \pi \alpha}).$

Пусть $s \in \mathbb{N}, \alpha \neq 2k+1, ~ k \in \mathbb{Z}$
{\rm (}или $\cos \pi \alpha \neq -1${\rm )}. Тогда оператор
${_1Y_-^{\alpha}}$ действует непрерывно из в $K^{2s}_{\alpha}$ в
$H^{2s}_{\alpha},$ причём справедливо неравенство
\begin{equation}\label{2.725}{\|{_1Y_-^{\alpha}}f\|_{K^{2s}_{\alpha}} \leq A_6
\|f\|_{H^{2s}_{\alpha}},}\end{equation} с постоянной
$$
A_6=\|{_1Y_-^{\alpha}}\|_{K^{2s}_{\alpha} \to
H^{2s}_{\alpha}}=\|{_1P_-^{\alpha-\frac{1}{2}}}\|_{L_2}=1/
\max\limits (1, \, \sq{1+\cos \pi \alpha}).
$$
\end{theorem}

Все утверждения теоремы вновь следуют из результатов для свойств
норм ОП Бушмана---Эрдейи нулевого порядка гладкости.

Оператор Бесселя является радиальной частью лапласиана в
$\mathbb{R}^n.$ При такой интерпретации этого оператора в случае
нечётномерных пространств  будет выполнено условие
теоремы~\ref{2tvloz2}.

\begin{theorem} \label{2tvloz3}
Пусть выполнены условия{\rm :} $\alpha \neq 2k+1,$ $k \in
\mathbb{Z};$ $\alpha \neq -n,$ $n \in \mathbb{N};$ $\alpha+2s \neq
-2m_1-1,$ $m_1 \in \mathbb{N}_0;$ $\alpha-2s \neq -2m_2-2,$ $m_2
\in \mathbb{N}_0.$ Тогда операторы~\eqref{2.723} осуществляют
топологический изоморфизм пространств Соболева $H^{2s}$ и весового
пространства Соболева $H^{2s}_{\alpha}.$
\end{theorem}

Очевидно, что все условия теоремы~\ref{2tvloz3} выполнены при
полуцелых $\alpha \in \mathbb{R}.$ Поэтому справедлива

\begin{theorem} Пусть $s \in \mathbb{N},~\alpha - \dfrac{1}{2} \in \mathbb{Z}.$ Тогда операторы~\eqref{2.723} осуществляют топологический изоморфизм пространств Соболева $H^{2s}$ и весовых пространств Соболева $H^{2s}_{\alpha}.$
\end{theorem}

Аналогично можно ввести по формулам типа~\eqref{2.723} операторы
${_1X^{\alpha}}_{0+}$ и ${_1Y^{\alpha}}_{0+}.$ В качестве
приложения приведённых выше результатов можно также рассмотреть
действие операторов~\eqref{2.723} в пространствах с
нормами~\eqref{2.720}--\eqref{2.722} при произвольных весах, не
согласованных с постоянной $\alpha$ в операторе Бесселя
$B_{\alpha}.$

Полученные в этом пункте результаты для одномерного случая очевидным образом переносятся на многомерный случай для области, которая состоит из декартового произведения положительных полуосей или отрезков по каждой переменной. Например, в двумерном случае полученные результаты сразу могут быть применены для оценок норм и доказательству вложений в первом положительном квадранте или лежащем в нём прямоугольнике.

Таким образом в этом пункте с помощью ОП Бушмана---Эрдейи нулевого
порядка гладкости дан положительный ответ на вопрос, который давно
обсуждался в устном <<фольклоре>> "--- \textit{пространства
Киприянова изоморфны весовым пространствам Соболева}.  Приведённые
результаты ни в коем случае не умаляют ни существенного значения,
ни необходимости использования  пространств И.\,А.~Киприянова для
подходящего круга задач в теории дифференциальных уравнений с
частными производными с операторами Бесселя. Принципиальная
важность пространств И.\,А.~Киприянова для теории уравнений в
частных производных различных типов с операторами Бесселя отражает
общий методологический подход, который автор услышал в виде
красивого афоризма на пленарной лекции чл.-корр. РАН
Л.\,Д.~Кудрявцева: <<{\it Каждое дифференциальное уравнение должно
изучаться в своём собственном пространстве}!>>

Полученные результаты также подтверждают полезность и
эффективность для теории дифференциальных уравнений специального
класса ОП "--- Бушмана---Эрдейи нулевого порядка гладкости, по
существу введённого В.\,В.~Катраховым в 1980-х годах.

\chapter{Общие весовые краевые задачи для сингулярных эллиптических
уравнений}\label{ch4}

Эта глава распадается на две естественные части. В первой части
вводятся и изучаются новые функциональные пространства сначала для
случая полупространства, а затем для ограниченной области. При
этом в первом случае мы используем метод введённых в первом
параграфе операторов преобразования, а во втором случае
пространства вводятся с помощью локальных карт. Приводятся все
необходимые для дальнейшего сведения, в том числе внутренние
теоремы вложения, теоремы о весовых следах и соответствующие
теоремы о мультипликаторах.

Во второй части сначала строится регуляризатор для несингулярного
уравнения с постоянными коэффициентами с нелокальными граничными
условиями, в которых присутствует оператор лиувиллевского типа,
этот результат имеет вспомогательное значение. Основные результаты
главы сосредоточены в заключительном пункте, в котором ставится и
изучается общая весовая краевая задача для сингулярных
эллиптических уравнений высшего порядка. При этом в случае
<<постоянных>> коэффициентов применяется метод операторов
преобразования, который сводит весовую задачу к указанной выше
несингулярной задаче. В случае переменных коэффициентов
применяется один из вариантов классической эллиптической техники,
восходящей к Шаудеру. Мы строим двусторонний регуляризатор во
введённых ранее пространствах, откуда легко следуют все основные
положения о нётеровости изучаемой  весовой краевой задачи.

\section{Функциональные пространства $H_{\nu}^s(E_{+}^{n+1})$}\label{sec10}

\subsection{Определения и внутренние теоремы
вложения}\label{sec10.1}

Пусть $E^{n+1}$ "--- евклидово $(n+1)$-мерное пространство точек
$x= (x', y)= (x_1, \dots, x_n, y),$ где $x' \in E^n,$ $y \in E^1.$
Через $E^{n+1}_{+}$ обозначим полупространство $y>0,$ а через
$\ov{E^{n+1}_{+}}$ "--- его замыкание. На протяжении всей главы мы
будем использовать операторы преобразования $P_{\nu, e}$ и
$S_{\nu, e},$ введённые в пункте~\ref{sec4.3}. Здесь они будут
действовать по последней переменной, то есть по $y.$

Введём некоторые обозначения. Пусть $\mathring{C}^{\infty}
(\ov{E^{n+1}_{+}})$ является множеством бесконечно
дифференцируемых функций в $\ov{E^{n+1}_{+}}$ с компактным в
$\ov{E^{n+1}_{+}}$ носителем. Через $\mathring{C}^{\infty}_{\nu}
(E^{n+1}_{+})$ обозначим множество функций, допускающих
представление $f = P_{\nu, e}\, g,$  где $g \in
\mathring{C}^{\infty} (\ov{E^{n+1}_{+}}).$ В символической записи
$\mathring{C}^{\infty}_{\nu} (E^{n+1}_{+}) =  P_{\nu, e}\,
\mathring{C}^{\infty} (\ov{E^{n+1}_{+}}).$

Введем для целых $s \geq 0$ на пространстве
$\mathring{C}^{\infty}_{\nu} (E^{n+1}_{+})$ следующую норму:
\begin{equation}
\| f \|_{H^s_{\nu} \lr{E^{n+1}_{+}}} = \| S_{\nu, e} f \|_{H^s \lr{E^{n+1}_{+}}}.
    \label{3.1.1}
\end{equation}
Здесь и в дальнейшем через $H^s (\Omega) \equiv W^s_2 (\Omega)$ будем   обозначать пространства С.\,Л.~Соболева. Определение нормы~\eqref{3.1.1} корректно, поскольку $S_{\nu, e} = \lr{P_{\nu, e}}^{-1},$    а   значит $S_{\nu, e}\, f \in \mathring{C}^{\infty} (\ov{E^{n+1}_{+}}).$ Обозначим через $H^s_{loc} (E^{n+1}_{+})$ множество функций, определенных  в $E^{n+1}_{+},$ которые в каждом слое\\
${E^{n+1}_{a, b} = \{x=(x', y): x' \in E^n, a<y<b \}}$ при
$0<a<b<\infty$ принадлежит пространству $H^s (E^{n+1}_{a, b}).$
Система полунорм $p_{s, a, b} (f) = \| f\|_{H^s (E^{n+1}_{a, b})}$
превращает его в пространство Фреше.

\begin{lemma} \label{lem: 3.1.1}
Пусть чётное $s \geq 0$ и $\Re \nu \geq 0.$   Тогда   для любой
функции $ f \in \mathring{C}^{\infty}_{\nu} \lr{E^{n+1}_{+}}$
имеет   место оценка
\begin{equation}
p_{s, a, b} (f) \leq c \,  \|  f \|_{H^s_{\nu} \lr{E^{n+1}_{+}}},
\label{3.1.2}
\end{equation}
где $0<a<b<\infty$ и постоянная $c>0$ не    зависит от $f.$
\end{lemma}

\begin{proof}
Функция $g=S_{\nu} f$ по условию принадлежит пространству
$\mathring{C}^{\infty} (\ov{E^{n+1}_{+}})$ и ${P_{\nu, e} g =f}.$
Применим    к функции формулу (см.~\eqref{1.3.13})
\begin{equation}
P_{\nu, e} \, g (x', y) = c_{\nu}\, y^{-\nu}
\int\limits_{-\infty}^{\infty} H^{(1)}_{\nu} (y \eta) \eta^{\nu}
(1 - i \eta)^{\frac{1}{2}- \nu} F \mbox{П} g (x', \eta)   \, d
\eta, \label{3.1.3}
\end{equation}
где $F$ "--- преобразование Фурье по последней переменной, a
$\mbox{П}$ "--- оператор продолжения по Уитни из $E^{n+1}_{+}$ в
$E^{n+1}.$ Пусть сначала $s=0.$ По  неравенству
Коши---Буняковского получаем
\begin{multline}
\| P_{\nu, e}\, g \|^2_{L_2 \lr{E^{n+1}_{a, b}}} \leq c
\int\limits_{E^{n+1}_{a, b}} \int\limits_{|\eta| < 1}
|H^{(1)}_{\nu} (y \eta) \eta^{\nu} (1 - i \eta)^{\frac{1}{2}- \nu}
|^2   \, d \eta  \int\limits_{|\eta| < 1} | F \mbox{П} g |^2 \, d
\eta   \, d x + {}
\\
{} + c  \int\limits_{E^{n+1}_{a, b}}  | \int\limits_{|\eta| > 1}
H^{(1)}_{\nu} (y \eta) \, \eta^{\nu} (1 - i \eta)^{\frac{1}{2}-
\nu} F \mbox{П} g (x', \eta) \, d \eta|^2 \, dx. \label{3.1.4}
\end{multline}
Первый внутренний интеграл в первом слагаемом справа
в~\eqref{3.1.4} ограничен, поскольку по асимптотическим формулам
для функций Ханкеля~\cite{BE2} подынтегральная функция либо
ограничена, либо имеет логарифмическую особенность в точке
$\eta=0.$ Поэтому, расширяя пределы интегрирования и учитывая
равенство Парсеваля, мы оценим первое слагаемое величиной  $c\,
\|\mbox{П} g \|_{L_2 (E^{n+1}_{+})}.$    Для оценки  второго
слагаемого нужно воспользоваться асимптотикой функции Ханкеля
$H_{\nu}^{(1)} (z)$ при $z \to \infty.$ Она   имеет   вид
(см.~\cite[с.~98]{BE2})
\begin{equation}
H_{\nu}^{(1)} (z) = \frac{e^{i z}}{\sqrt{z}} \lr{c+ O(|z|^{-1})},
\label{3.1.5}
\end{equation}
где $c$ "--- некоторая постоянная. Используя равенство Парсеваля,
имеем
$$
\int\limits_a^b \left| \int\limits_{|\eta| > 1} e^{i z}  (y
\eta)^{-\frac{1}{2}} \, \eta^{\nu} (1 - i \eta)^{\frac{1}{2}- \nu}
F \mbox{П} g (x', \eta)   \, d \eta\right|^2 \, dy \leq
 c \int\limits_{-\infty}^{\infty} | F \mbox{П} g |^2 \, d \eta
\leq c \int\limits_{- \infty}^{\infty} |  \mbox{П} g (x', y) |^2
\, d y.
$$
Таким образом, получена нужная оценка выражения, соответствующего
первому слагаемому из~\eqref{3.1.5}. Оценивая выражение,
соответствующее остаточному члену, по неравенству
Коши---Буняковского, получим
$$
\int\limits_a^b | \int\limits_{|\eta| > 1}   (y
\eta)^{-\frac{1}{2}} \eta^{\nu}  (1 - i \eta)^{\frac{1}{2}- \nu}
O(|y \eta|^{-1}) F \mbox{П} g (x', \eta)   \, d \eta|^2 \, dy \leq
$$
$$
\leq c \int\limits_a^b  \int\limits_{|\eta| > 1} | y \eta |^{-2}
\, d \eta \int\limits_{|\eta| > 1}  | F \mbox{П} g |^2 \, d \eta
\,  dy \leq c \int\limits_{- \infty}^{\infty} |  \mbox{П} g (x',
y) |^2 \, d y.
$$
Следовательно, так как $\|\mbox{П} g\|_{L_2 (E^n)} \leq c \,
\|g\|_{L_2 (E^{n+1}_{+})},$    то  неравенство~\eqref{3.1.2} при
$s=0$ установлено.

Рассмотрим случай чётного $s>0.$    Заметим сначала, что
$$
p^2_{s, a, b} (f) \leq \sum\limits_{|\alpha'|+2 \alpha_{n+1} \leq
s} \|D^{\alpha'}_{x'} \, D^{2 \alpha_{n+1}}_y \, f \|^2_{L_2
(E^{n+1}_{a, b})} \leq
 c \sum\limits_{|\alpha'|+2 \alpha_{n+1} \leq s} \|D^{\alpha'}_{x'}
\, B^{\alpha_{n+1}}_y \, f \|^2_{L_2 (E^{n+1}_{a, b})},
$$
где использованы обычные обозначения:
$$
D_{x'}^{\alpha'} = \frac{\pr^{|\alpha'|}}{\pr x_1^{\alpha_1} \pr
x_2^{\alpha_2} \dots \pr x_n^{\alpha_n}},~ |\alpha'| = \alpha_1 +
\dots + \alpha_n,~B_y = \frac{\pr^2}{\pr y^2} + \frac{2 \nu +1}{y}
\frac{\pr }{\pr y}.
$$

Учитывая  доказанное  и формулу для операторов преобразования
$S_{\nu, e} B = D^2 S_{\nu, e},$ получаем
$$
p^2_{s, a, b} (f) \leq c \sum\limits_{|\alpha'|+2 \alpha_{n+1}
\leq s} \|S_{\nu}\, D^{\alpha'}_{x'} \, B^{ \alpha_{n+1}}_y \, f
\|^2_{L_2 (E^{n+1}_{+})} \leq
 c\,  \| S_{\nu, e} \, f \|_{H^s (E^{n+1}_{+})} = c \,
\|f\|^2_{H^s_{\nu} (E^{n+1}_{+})}.
$$
Лемма доказана.
\end{proof}

Пространство $\mathring{C}^{\infty}_{\nu} ({E}^{n+1}_{+}),$
наделенное нормой~\eqref{3.1.1}, в силу леммы непрерывно вложено в
полное пространство $H^s_{loc} ({E}^{n+1}_{+}).$ Замыкая
$\mathring{C}^{\infty} ({E}^{n+1}_{+})$ по норме~\eqref{3.1.1} при
$s \geq 0,$ $\Re \nu \geq 0,$ мы получаем пространство,
обозначаемое через  $H^s_{\nu} ({E}^{n+1}_{+}).$

Из леммы~\ref{lem: 3.1.1} непосредственно вытекает

\begin{corollary} {\label{cor: 3.1.2}}
     При $s \geq 0,$ $\Re \nu \geq 0$  гильбертово пространство $H^s_{\nu} ({E}^{n+1}_{+})$ непрерывно вложено в пространство Фреше    $H^s_{loc} ({E}^{n+1}_{+}).$
\end{corollary}

\begin{lemma}{\label{lem:3.1.2}}
Пусть чётное $s \geq 0,$ $\Re \nu \geq 0$ и $0<a<b< \infty.$ Пусть
функция $f \in H^s_{\nu} ({E}^{n+1}_{+})$ и её носитель заключён в
слое $E^{n+1}_{a, b}.$ Тогда $f \in H^s ({E}^{n+1}_{+})$ и
справедливо неравенство
\begin{equation}
c' \, \| f \|_{H^s (E^{n+1}_{+})} \leq \| f \|_{H^s_{\nu} (E^{n+1}_{+})} \leq c'' \, \| f \|_{H^s (E^{n+1}_{+})},
\label{3.1.6}
\end{equation}
постоянные $c', c'' > 0$ в котором не зависят от $f.$
\end{lemma}

\begin{proof}
То, что $f \in H^s ({E}^{n+1}_{+}),$ а также левое неравенство
в~\eqref{3.1.6}, установлено по сути дела при доказательстве
леммы~\ref{lem: 3.1.1}. Докажем правое неравенство. Напомним, что
\begin{equation}
S_{\nu, e} = \mathcal{J}_{\nu - \frac{1}{2}, e}\, S_{\nu},
\label{3.1.7}
\end{equation}
где
\begin{equation}
 \mathcal{J}_{\nu - \frac{1}{2}, e} \,f (x', y) =  f (x', y) + \lr{\frac{1}{2}-\nu} \int\limits_y^{\infty}   f (x', t) \Phi \lr{\nu + \frac{1}{2}, 2; y-t} dt,
\label{3.1.8}
\end{equation}
\begin{equation}
S_{\nu} f (x', y) =  c_{\nu} \lr{y^{\nu + \frac{1}{2}} f (x', y) -
\int\limits_y^{\infty} t^{\nu + \frac{1}{2}}  f (x', t) \frac{\pr
}{\pr y} P_{\nu - \frac{1}{2}}^0 \lr{\frac{y}{t}} dt}
\label{3.1.9}
\end{equation}
и $\Phi$ "--- вырожденная гипергеометрическая функция, а
$P^0_{\mu}$ "--- функция Лежандра первого рода. Эти функции,
являющиеся ядрами интегральных  операторов, суть гладкие функции.
Поэтому, поскольку $f(x', y)=0$ при $y<a$ и при $y>b,$ то  норма
функции $S_{\nu, e} f$ в пространстве $H^s (E^{n+1}_{+}),$
оценивается через норму функции $f$ в том же пространстве. Этим
заканчивается доказательство леммы.
\end{proof}

Введем  пространство $H^s_{\nu} (E^{n+1}_{0, b})$ для слоя
$E^{n+1}_{0, b},$ $0<b<\infty$ как замыкание по
норме~\eqref{3.1.1} множества всех функций из
$\mathring{C}^{\infty}_{\nu} ({E}^{n+1}_{+}),$  носители которых
заключены в $\ov{E}^{n+1}_{0, b}.$

\begin{lemma}{\label{lem: 3.1.2}} Пусть чётное $s \geq 0,$ $\Re \nu \geq 0$ и $0<b< \infty.$ Тогда норма~\eqref{3.1.1} и норма
\begin{equation}
 \| f \|_{H^s (E^{n+1}_{0, b})} = \| S_{\nu} f \|_{H^s (E^{n+1}_{0, b})}
\label{3.1.10}
\end{equation}
на пространстве $H^s_{\nu} (E^{n+1}_{0, b})$ эквивалентны.
\end{lemma}

\begin{proof}
Пусть функция  $f \in \mathring{C}^{\infty}_{\nu} ({E}^{n+1}_{+})$
и $\supp f \subset \ov{E^{n+1}_{+}}.$ По определению это означает,
что $S_{\nu, e}\, f \in \mathring{C}^{\infty}_{\nu}
(\ov{E^{n+1}_{+}}).$ А тогда и функция $S_{\nu} f =
\mathcal{J}_{\frac{1}{2}-\nu, e}\, S_{\nu, e}\, f \in
\mathring{C}^{\infty}_{\nu} (\ov{E^{n+1}_{+}}).$ Последнее легко
следует из~\eqref{3.1.7} и~\eqref{3.1.8}. Заметим, что
формула~\eqref{3.1.8} верна при всех комплексных $\nu.$ Операторы
$\mathcal{J}_{\mu, e}$   и $\mathcal{J}_{-\mu, e}$ обратны друг
другу на пространстве $\mathring{C}^{\infty}_{\nu}
(\ov{E^{n+1}_{+}}).$ Кроме того, из гладкости функции $\Phi$ легко
получаем оценку
$$
c' \, \|  \mathcal{J}_{\mu, e} g \|_{H^s (E^{n+1}_{+})} \leq \| g
\|_{H^s (\ov{E^{n+1}_{+}})} \leq c'' \, \| \mathcal{J}_{\mu, e} g
\|_{H^s (E^{n+1}_{+})},
$$
которая верна для функций $g \in \mathring{C}^{\infty}
(\ov{E^{n+1}_{+}})$ таких, что $\supp g \subset \ov{E}^{n+1}_{0,
b}.$ Постоянные $c'$ и $c''$ зависят от $b,$ но не зависят от
функции $g.$ Таким образом, эквивалентность норм~\eqref{3.1.1}
и~\eqref{3.1.10} доказана на плотном множестве $f \in
\mathring{C}^{\infty}_{\nu} ({E}^{n+1}_{+}),$ $\supp f \subset
\ov{E}^{n+1}_{0, b}.$ Лемма доказана.
\end{proof}

Сделаем одно замечание. Нами доказана эквивалентность
норм~\eqref{3.1.1} и~\eqref{3.1.10} для пространств $H^s_{\nu}
(E^{n+1}_{0, b})$ в случае конечного $b.$ При $b=\infty$ это уже
не имеет места. Можно показать, что   норма~\eqref{3.1.10} при
$b=\infty$ подчинена норме~\eqref{3.1.1}. Поэтому, если в
определении пространства $H^s (E^{n+1}_{+})$ заменить
норму~\eqref{3.1.1} на норму~\eqref{3.1.10}, то получим более
широкое пространство, которое оказывается не вложенным в
$H^s_{loc} (E^{n+1}_{+}).$ Это и послужило причиной использования
нормы~\eqref{3.1.1}, которая порождена оператором преобразования
$S_{\nu, e},$ а не оператором $S_{\nu}.$

Установим соотношение введённых пространств и функциональных
пространств, изученных И.\,А.~Киприяновым~\cite{Kip1}. Пусть $F$
$(n+1)$ "--- мерное преобразование Фурье:
$$
F f(\xi)  = \int\limits_{E^{n+1}} f (x', y) e^{- i \langle \xi', x' \rangle - i y \eta} dx,
$$
где $\xi = (\xi', \eta) \in E^{n+1},$ $\xi' \in E^n,$ $\langle
\xi', x' \rangle = \xi_1 x_1 + \dots +  \xi_n x_n.$ Через
$F_{\nu}$ обозначим преобразование Фурье---Бесселя (или
Фурье---Ханкеля)
$$
F_{\nu} f(\xi)  = \int\limits_{E^{n+1}_{+}} f (x', y) e^{- i \langle \xi', x' \rangle} j_{\nu} (y \eta) y^{2 \nu +1 } \, dx.
$$
Пространства    $H^s_{\nu, +} (E^{n+1}_{0, b}),$ $s \geq 0,$ $\nu
\geq -\dfrac{1}{2}$ определяются как замыкание множества чётных
функций $f \in \mathring{C}^{\infty} (\ov{E^{n+1}_{+}}),$ для
которых $\supp f \subset \ov{E^{n+1}_{0, b}}$ (множество  таких
функций обозначим   через $\mathring{C}^{\infty}_{+}
(\ov{E^{n+1}_{0, b}})$) по норме
\begin{equation}
\| f \|_{H^s_{\nu, +} (E^{n+1}_{+})} =  \lr{ ~ \int\limits_{E^{n+1}_{+}} |F_{\nu} f (\xi) |^2 (1+|\xi|^2)^s {\eta}^{2 \nu +1 } \, d \xi}^{\frac{1}{2}}.
\label{3.1.11}
\end{equation}

Многократно применяя лемму~\ref{lem: 3.1.2}, мы приходим к
следующему утверждению.

\begin{lemma}{\label{lem: 3.1.4}}
Пусть чётное $s \geq 0$ и вещественное $\nu \geq 0.$ Тогда
пространство $H^s_{\nu, +} (E^{n+1}_{0, b})$ непрерывно  вложено в
$H^s_{\nu} (E^{n+1}_{0, b}).$ При $\nu \neq 1,3,5,\dots$
$H^s_{\nu, +} (E^{n+1}_{0, b})$ образует собственное
подпространство  пространства $H^s_{\nu} (E^{n+1}_{0, b}),$ причём
индуцированная и собственная нормы на $H^s_{\nu, +} (E^{n+1}_{0,
b})$ эквивалентны.
\end{lemma}

Рассмотрим теперь внутренние теоремы вложения введённых
пространств $H^s_{\nu}.$

\begin{theorem} \label{teo: 3.1.1}
    Пусть $\Re \nu \geq 0$ и $s'>s \geq 0.$ Тогда пространство    $H^{s'}_{\nu} (E^{n+1}_{+})$ непрерывно вложено в $H^s_{\nu} (E^{n+1}_{+})$  и справедливо  соответствующее неравенство между нормами.
\end{theorem}

Это утверждение есть прямое следствие определения пространств
$H^s_{\nu}$ и теорем вложения пространств $H^s.$

\begin{theorem}{\label{teo: 3.1.2}}
    Пусть $\Re \nu \geq 0$ и $s'>s \geq 0.$ Пусть $f_k,$ $k=1,2,\dots$ ограниченная последовательность функций из $H^{s'}_{\nu} (E^{n+1}_{+}),$ причём $\supp f_k \subset \mathcal{Y},$  где $\mathcal{Y}$ "--- ограниченное множество в $E^{n+1}_{+}.$ Тогда существует подпоследовательность, сходящаяся по норме пространства $H^s_{\nu} (E^{n+1}_{+}).$
\end{theorem}

\begin{proof}
Можно построить последовательность функций  $\varphi_k (x) \in
\mathring{C}^{\infty}_{\nu} ({E}^{n+1}_{+})$  таких, что $\supp
\varphi_k  \subset \mathcal{Y}'$ и  $\|f_k - \varphi_k
\|_{H^{s'}_{\nu} (E^{n+1}_{+})} \to 0$ при $k \to \infty.$ Здесь
$\mathcal{Y}'$ "--- некоторый компакт в $\ov{E^{n+1}_{+}}.$ Тогда
последовательность норм $\|\varphi_k  \|_{H^{s'}_{\nu}
(E^{n+1}_{+})} = \|S_{\nu, e} \varphi_k  \|_{H^{s'}_{\nu}
(E^{n+1}_{+})}$ ограничена. По теоремам о полной непрерывности
вложения пространств $H^s$    (см. С.\,М.~Никольский~\cite{66})
найдется функция $g \in H^s (E^{n+1}_{+}),$ к которой сходится
некоторая подпоследовательность $S_{\nu, e} \varphi_{k_m}.$
Значит, она фундаментальна в    $H^s (E^{n+1}_{+}),$ а  тогда
фундаментальна и последовательность $f_{k_m}$ в $H^s_{\nu}
(E^{n+1}_{+}).$ В силу полноты пространства $H^s_{\nu}$ эта
последовательность будет сходящейся. Теорема доказана.
\end{proof}

\subsection{Некоторые результаты о
мультипликаторах}\label{sec10.2}

В этом пункте будут выяснены достаточные условия на функцию
$a(x),$ при которых отображение $f \to a f$ будет непрерывным в
пространствах $H^s_{\nu} (E^{n+1}_{+}).$ Кроме того, будут даны
оценки нормы этого отображения.

\begin{lemma}{\label{lem: 3.2.1}}
    Пусть $a(x) \in C^{\infty}  (\ov{E^{n+1}_{+}})$  и выполнено условие
\begin{equation}
D^k_y a(x) = 0  \  \mbox{при} \  y = 0, \  k =1,2, \dots.
\label{3.2.1}
\end{equation}
Тогда при $\Re \nu \geq 0$ функция $a f$ принадлежит множеству
$\mathring{C}^{\infty}_{\nu} ({E}^{n+1}_{+}),$ если $f \in
\mathring{C}^{\infty}_{\nu} ({E}^{n+1}_{+}).$
\end{lemma}

\begin{proof}
Покажем, что если функция $g \in \mathring{C}^{\infty}
(\ov{E^{n+1}_{+}}),$  то функция $S_{\nu, e}\, (a P_{\nu, e}\, g)$
принадлежит тому же классу. Так как $S_{\nu, e} =
I^{\frac{1}{2}-\nu}_{e} S^{\nu- \frac{1}{2}}_{\nu}$ и $P_{\nu, e}
= P^{\frac{1}{2}-\nu}_{\nu} I^{\nu- \frac{1}{2}}_{e}$ и операторы
$I^{\mu}_e$ отображают $\mathring{C}^{\infty} (\ov{E^{n+1}_{+}})$
на себя при любых   комплексных $\mu,$  то достаточно показать,
что $S^{\nu- \frac{1}{2}}_{\nu} (a P^{\frac{1}{2}-\nu}_{\nu} g)
\in \mathring{C}^{\infty} (\ov{E^{n+1}_{+}}).$   Из определения
операторов преобразования $S^{\nu- \frac{1}{2}}_{\nu}$ и
$P^{\frac{1}{2}-\nu}_{\nu},$ где $N$ "--- натуральное число,
получаем
$$
S_{\nu}^{\nu - \frac{1}{2}} \lr{a P_{\nu}^{ \frac{1}{2} - \nu} f}
(x', y) = \frac{(-1)^{N+1} \, 2^{1-N}}{\Gamma \lr{\nu +
\frac{1}{2}} \,  \Gamma \lr{N-\nu + \frac{1}{2}}}\, \frac{\pr}{\pr
y} \int\limits_y^{\infty} (t^2-y^2)^{\nu - \frac{1}{2}} a (x', t)\, t \times{}
$$
$$
{}\times \int\limits_t^{\infty} (\tau^2-t^2)^{N-\nu - \frac{1}{2}} \, \lr{\frac{\pr}{\pr \tau} \frac{1}{\tau}}^N g(x', \tau) \, d \tau d t =
$$
$$
= \frac{(-1)^{N+1} \, 2^{1-N}}{\Gamma \lr{\nu + \frac{1}{2}} \,
\Gamma \lr{N-\nu + \frac{1}{2}}}\, \frac{\pr}{\pr y}
\int\limits_y^{\infty} \lr{\lr{\frac{\pr}{\pr \tau}
\frac{1}{\tau}}^N g(x', \tau)} \,
 \int\limits_y^{\tau} (\tau^2-y^2)^{\nu - \frac{1}{2}} \, t \, (\tau^2-t^2)^{N-\nu - \frac{1}{2}}   a (x', t)\, d t d \tau.
$$
Во внутреннем интеграле справа сделаем замену переменной по
формуле $t = \sqrt{y^2 + z (\tau^2-y^2)}.$ Тогда
$$
\int\limits_y^{\tau} (\tau^2-y^2)^{\nu - \frac{1}{2}} \, t \,
(\tau^2-t^2)^{N-\nu - \frac{1}{2}}   a (x', t)\, d t =
\frac{(\tau^2-y^2)^N }{2}
 \int\limits_0^1 z^{\nu - \frac{1}{2}}  z^{N-\nu -
\frac{1}{2}} a (x', \sqrt{y^2 + z (\tau^2-y^2)}) \, d z.
$$
Отсюда интегрированием по частям находим
$$
S_{\nu}^{- \frac{1}{2} + \nu} \lr{a P_{\nu}^{ \frac{1}{2} - \nu}
g} (x', y) = \frac{-1}{\Gamma \lr{\nu + \frac{1}{2}} \,  \Gamma
\lr{N-\nu + \frac{1}{2}}}\, \frac{\pr}{\pr y}
\int\limits_y^{\infty} g(x', \tau)  \times
$$
$$
\times \frac{\pr^N}{(\pr \tau^2)^N} \left(  (\tau^2-y^2)^N
\int\limits_0^1 z^{\nu - \frac{1}{2}}  z^{N-\nu - \frac{1}{2}} a
(x', \sqrt{y^2 + z (\tau^2-y^2)}) \, d z \right)  d \tau   =
$$
$$
= \frac{-1}{\Gamma \lr{\nu + \frac{1}{2}} \,  \Gamma \lr{N-\nu +
\frac{1}{2}}}\, \frac{\pr}{\pr y} \int\limits_y^{\infty} g(x',
\tau) \int\limits_0^1 \sum\limits_{k=0}^{N} 2^{k-N} \binom{N}k
\frac{N!}{(N-k)!}  \times
$$
$$
\times (\tau^2-y^2)^{N-k} z^{\nu + N - k - \frac{1}{2}}
(1-z)^{N-\nu - \frac{1}{2}}\left. \lr{\frac{\pr}{\lambda \pr
\lambda}}^{N-k} a(x', \lambda)\right|_{\lambda= \sqrt{y^2 + z
(\tau^2-y^2)}} dz d \tau.
$$
Здесь через $\binom{N}k$ обозначены биномиальные коэффициенты.
Дифференцируя по параметру $y,$ получим
\begin{equation}
S_{\nu}^{\nu- \frac{1}{2}} \lr{a P_{\nu}^{ \frac{1}{2} - \nu} g}
(x', y) =  a (x', y) g (x', y) + \int\limits_y^{\infty} g(x',
\tau) a_{\nu} (x', y, \tau) \, d \tau, \label{3.2.2}
\end{equation}
где положено
\begin{multline}
    a_{\nu} (x', y, \tau) = \frac{-1}{\Gamma \lr{\nu +\frac{1}{2}}
    \, \Gamma \lr{N-\nu +\frac{1}{2}}} \frac{\pr}{\pr y}
    \sum\limits_{k=0}^{N}   \binom{N}k \frac{2^{k-N} \, N!}{(N-k)!}
    (\tau^2-y^2)^{N-k} \times \\
\left. \times \int\limits_0^1
    z^{\nu + N - k - \frac{1}{2}} (1-z)^{N-\nu - \frac{1}{2}}
    \lr{\frac{\pr}{\lambda \pr \lambda}}^{N-k} a(x',
    \lambda)\right|_{\lambda= \sqrt{y^2 + z (\tau^2-y^2)}} dz.
    \label{3.2.3}
\end{multline}
Из~\eqref{3.2.3} и из условия~\eqref{3.2.1} следует, что функция
$a_{\nu} (x', y, \tau)$ бесконечно дифференцируема при $x' \in
E^n$ и $y, \tau \geq 0.$ Поэтому левая часть формулы~\eqref{3.2.2}
также бесконечно дифференцируема и, очевидно, финитна. Лемма
доказана.
\end{proof}

Наложим на функцию $a(x)$ еще некоторые ограничения. Пусть при некоторых $R < \infty$
\begin{equation}
D_y \, a (x', y) = 0, \  y \geq  R,
\label{3.2.4}
\end{equation}
и при всех мультииндексах $\alpha$
\begin{equation}
\sup\limits_{x \in E_{+}^{n+1}} \left| D^{\alpha'}_{x'}  \lr{\frac{1}{y} D_y}^{\alpha_{n+1}} a(x)\right| = M^{\alpha', \alpha_{n+1} } = M^{\alpha}  < \infty.
\label{3.2.5}
\end{equation}

Уточним в этих условиях некоторые свойства функции $a_{\nu} (x',
y, \tau).$ Имеет место оценка
\begin{equation}
 \left| D^{\alpha'}_{x'}  \lr{\frac{1}{y} D_y}^{l} \lr{\frac{1}{\tau} D_{\tau}}^{m} a_{\nu} (x', y, \tau) \right| \leq c \, (1+\tau)^{- 2 \Re \nu -1-m} \max\limits_{k \leq N+l+m+1}  M^{\alpha', k},
\label{3.2.6}
\end{equation}
где постоянная $c>0$ не зависит от $x',$ $y,$ $\tau$ и функции
$a.$ В самом деле, в правую часть формулы~\eqref{3.2.3}  входит не
сама функция $a,$ а только производные от неё по последней
переменной. Тогда из~\eqref{3.2.4} имеем
\begin{equation}
a_{\nu} (x, y, \tau) = 0, \  \tau \geq y \geq  R.
\label{3.2.7}
\end{equation}
Следовательно, оценка~\eqref{3.2.6} тем более справедлива при
$\tau \geq y\geq R.$ Кроме того, она, очевидно, справедлива и при
$0 \leq y \leq \tau \leq R.$    Пусть теперь $\tau>R>y,$    тогда
ввиду~\eqref{3.2.4}   имеем
$$
\left| \left. \int\limits_0^1 z^{\nu + N - k - \frac{1}{2}} (1-z)^{N-\nu - \frac{1}{2}} \lr{\frac{\pr}{\lambda \pr \lambda}}^{N-k} a(x', \lambda)\right|_{\lambda= \sqrt{y^2 + z (\tau^2-y^2)}} dz \right|=
$$
$$
= \left| \left. \int\limits_0^{\frac{R^2-y^2}{\tau^2-y^2}} z^{\nu + N - k - \frac{1}{2}} (1-z)^{N-\nu - \frac{1}{2}} \lr{\frac{\pr}{\lambda \pr \lambda}}^{N-k} a(x', \lambda)\right|_{\lambda= \sqrt{y^2 + z (\tau^2-y^2)}} dz \right| \leq
$$
$$
\leq \frac{1}{ |\nu + N - k + \frac{1}{2}|} \lr{\frac{R^2-y^2}{\tau^2-y^2}}^{\Re \nu + N - k + \frac{1}{2}} \sup\limits_{y \geq 0} \left|  \lr{\frac{\pr}{y \pr y}}^{N-k} a(x', y) \right|.
$$
Это соотношение и приводит к неравенству~\eqref{3.2.6}.

Формула~\eqref{3.2.3} позволяет получить аналогичное представление
оператора $S_{\nu, e}\, a\, P_{\nu, e}.$ Так как $S_{\nu, e}=
I^{\frac{1}{2} - \nu} S^{\nu - \frac{1}{2}}_{\nu},$ $P_{\nu, e}=
P_{\nu}^{\frac{1}{2} - \nu} I^{\nu - \frac{1}{2}}_{e},$ то
из~\eqref{3.2.3} получаем
\begin{equation}
S_{\nu, e} \lr{a \, P_{\nu, e} \, g} (x', y) = I_{e}^{\frac{1}{2}
- \nu} \lr{a \, I_{e}^{\nu - \frac{1}{2} } g} (x', y) +
I_{e}^{\frac{1}{2} - \nu}  \int\limits_y^{\infty} a_{\nu} (x', y,
\tau) I_{e}^{\nu - \frac{1}{2} } g (x', \tau)\, d \tau,
\label{3.2.8}
\end{equation}
где функция $g \in \mathring{C}^{\infty} (\ov{E^{n+1}_{+}}).$
Первое слагаемое справа рассматривалось в лемме~\ref{lem:1.2.3}.
Поэтому нам достаточно оценить второе слагаемое. При $\dfrac{1}{2}
< \Re \nu < N+\dfrac{1}{2}$ имеем
\begin{multline}
I_{e}^{\frac{1}{2} - \nu}  \int\limits_y^{\infty} a_{\nu} (x', y,
\tau) I_{e}^{\nu - \frac{1}{2} } g (x', \tau)\, d \tau
= (-1)^N e^y I^{\frac{1}{2} - \nu+N} D^N_y \int\limits_y^{\infty} e^{-y} a_{\nu} (x', y, \tau) I_{e}^{\nu - \frac{1}{2} } g (x', \tau)\, d \tau = \\
\left. = \sum\limits_{k=0}^{N-1} (-1)^k  I^{\frac{3}{2} - \nu+k} \lr{e^y D^k_y \lr{e^{-\tau} a_{\nu} (x', y, \tau)}\right|_{\tau=y}  I^{\nu - \frac{1}{2} } g (x', y)}+ \\
+ (-1)^N  I_e^{\frac{1}{2} - \nu+N} \int\limits_y^{\infty} e^{y}
D^N_y \lr{e^{-y} a_{\nu} (x', y, \tau)} I_{e}^{\nu - \frac{1}{2} }
g (x', \tau)\, d \tau. \label{3.2.9}
\end{multline}
Из формулы~\eqref{3.2.3} нетрудно усмотреть, что функция
$$
e^y D^k_y \left. \lr{e^{-y} a_{\nu} (x', y, \tau)}\right|_{\tau=y}, \  k = 0, \dots, N-1,
$$
удовлетворяет всем условиям следствия~\ref{cor:1.2.2} из
главы~\ref{ch2}. Поэтому
\begin{multline}
\| \left. \sum\limits_{k=0}^{N-1} (-1)^k  I^{\frac{3}{2} - \nu+k} \lr{e^y D^k_y \lr{e^{-y} a_{\nu} (x', y, \tau)} \right|_{\tau=y}  I_e^{\nu - \frac{1}{2} } g (x', y)} \|_{H^s (E_{+}^{n+1})} \leq \\
\leq c \sum\limits_{\substack{|\alpha'| \leq  s \\
k \leq 3 N +s+1}} M^{\alpha', k} \|g\|_{H^s (E^{n+1}_{+})},
\label{3.2.10}
\end{multline}
где постоянная $c>0$ не зависит от функций $a$ и $g.$

Для оценки последнего интеграла справа в~\eqref{3.2.9} заметим,
что операторы $I_e^{\frac{1}{2}-\nu+N}$ и $I_e^{\nu-\frac{1}{2}}$
принадлежат классу $L (H^s (E^{n+1}_{+}), H^s  (E^{n+1}_{+})),$
так как $\Re \Big(\dfrac{1}{2}-\nu+N\Big)>0,$ $\Re \nu -
\dfrac{1}{2}
> 0$  (см. лемму~\ref{lem:1.2.2}). Кроме того, оценка~\eqref{3.2.6} показывает, что интегральный оператор
$$
w \to \int\limits_y^{\infty} w(x', \tau) \, e^y \, D_y^N  \lr{e^{-y} a_{\nu} (x', y, \tau)} d \tau
$$
также принадлежит классу ${L} \lr{ H^s (E^{n+1}_{+}), H^s (E^{n+1}_{+})}$ при любых $s$ и его норма не превосходит величины
$$
c \sum\limits_{  k \leq 3 N +s+1,~ |\alpha'| \leq  s} M^{\alpha', k},
$$
где постоянная $c>0$ не зависит от функции $a.$

Такие же оценки справедливы и в случае  $0 \leq \Re \nu \leq
\dfrac{1}{2}.$ Для их доказательства достаточно оценить норму
второго слагаемого в формуле~\eqref{3.2.8}. Интегрируя по частям,
находим
$$
I_e^{\frac{1}{2} {-} \nu}\!\! \int\limits_y^{\infty}  I_e^{\nu {-} \frac{1}{2}} g (x'\!, \tau)  a_{\nu} (x'\!, y\!, \tau)  d \tau =  I_e^{\frac{1}{2} {-} \nu}  a_{\nu} (x'\!, y\!, y) I_e^{\nu {+} \frac{1}{2}} g {+}
 I_e^{\frac{1}{2} {-} \nu}\!\! \int\limits_y^{\infty}  e^{{-} \tau} D_{\tau} \lr{ e^{ \tau} a_{\nu} (x'\!, y\!, \tau)} I_e^{\nu {+} \frac{1}{2}} g (x'\!, \tau)  d \tau.
$$
Из неравенства~\eqref{3.2.6} теперь легко следует, что $ H^s
(E^{n+1}_{+})$ "--- норма этого выражения не превосходит величины
$c \sum\limits_{k \leq 4+s,~|\alpha'|\leq s} M^{\alpha', k},$ где
постоянная $c>0$ не зависит от функции $a.$

Итак, доказана

\begin{theorem}\label{teo: 3.2.1}
Пусть функция $a \in C^{\infty} (\ov{E^{n+1}_{+}})$ и
удовлетворяет условиям~\eqref{3.2.1}, \eqref{3.2.4},
\eqref{3.2.5}. Тогда для любых комплексных $\nu$ из полуплоскости
$\Re \nu \geq 0,$ любых $s \geq 0$ и любых ${R>0},$ не
превосходящих некоторого числа $R_0>0,$ оператор $S_{\nu, e} \, a
\, P_{\nu, e}$ непрерывно отображает пространство $ H^s
(E^{n+1}_{+})$ в себя. Справедлива оценка
\begin{equation}
\| S_{ \nu, e} \, a \, P_{ \nu, e} f \|_{H^s (E^{n+1}_{+})} \leq c \, \|f\|_{H^s (E^{n+1}_{+})} \sum\limits_{\substack{|\alpha'| \leq  s \\
\alpha_{n+1} \leq 3 N +s+1}} M^{\alpha', \alpha_{n+1} },
\label{3.2.11}
\end{equation}
где $N$ наименьшее натуральное число, для которого $\Re \nu <
N+1.$ Постоянная $c>0$ зависит от $\nu, n, s, R_0,$ но не зависит
от функций $a, f,$    а также от $R.$
\end{theorem}

\begin{corollary}\label{cor: 3.2.1}
    В условиях теоремы справедлива оценка
    \begin{equation}
    \|  a  f \|_{H^s_{\nu} (E^{n+1}_{+})} \leq c  \sum\limits_{\substack{|\alpha'| \leq  s \\
\alpha_{n+1} \leq 3 N +s+1}} M^{\alpha', \alpha_{n+1} }
\|f\|_{H^s_{\nu} (E^{n+1}_{+})}.
    \label{3.2.12}
    \end{equation}
\end{corollary}

\begin{corollary}\label{cor: 3.2.2}
Пусть выполнены условия теоремы и пусть
$$
\sup\limits_{x \in E^{n+1}_{+}} \left| D_{x'}^{\alpha'} \lr{\frac{1}{y} D_y}^{\alpha_{n+1}} \frac{a(x', y)}{y} \right| = \widetilde{M}^{\alpha', \alpha_{n+1} } < \infty.
$$
Тогда справедлива оценка
\begin{equation}
\|  a \, D_y \,  f \|_{H^s_{\nu} (E^{n+1}_{+})} \leq c  \sum\limits_{\substack{|\alpha'| \leq  s \\
\alpha_{n+1} \leq 3 N +s+1}} \widetilde{M}^{\alpha', \alpha_{n+1}
} \|f\|_{H^{s+1}_{\nu} (E^{n+1}_{+})}. \label{3.2.13}
\end{equation}
\end{corollary}

\begin{proof}
Формула
$$
y \, D_y \, P_{\nu}^{\frac{1}{2} - \nu} = P_{\nu}^{\frac{1}{2} -
\nu} \lr{y D_y-2 \nu}
$$
легко следует из определения оператора преобразования
(см.~\eqref{1.1.7},~\eqref{1.1.8}. В связи с этим для функции
$S_{\nu, e} f$ по аналогии с~\eqref{3.2.8} имеем
$$
S_{\nu, e} \lr{a \, D_y \, P_{\nu, e}\, g } = I_e^{\frac{1}{2} -
\nu } \lr{a (D_y - \frac{2 \nu}{y}) I^{\nu - \frac{1}{2}} } +
 I_e^{\frac{1}{2} - \nu} \int\limits_y^{\infty} \lr{\tau D_{\tau}
- 2 \nu}   \widetilde{a}_{\nu}  (x', y, \tau) I_e^{\nu -
\frac{1}{2}}  g (x', \tau)  \, d \tau.
$$
Здесь функция $\widetilde{a}_{\nu}$ определяется по
формуле~\eqref{3.2.3} с заменой в последней функции $a$ на функцию
$\dfrac{1}{y} \, a (x', y).$ После интегрирования по частям
последняя формула приводится к виду
$$
S_{\nu, e} \lr{a \, D_y \, P_{\nu, e} \, g } = I_e^{\frac{1}{2} -
\nu }  \lr{a \, I_e^{\nu - \frac{1}{2}} \, D_y \, g  } - 2 \nu
I_e^{\frac{1}{2} - \nu } \lr{ \frac{a}{y} I^{\nu - \frac{1}{2}} g}
-
 2 \nu \, I_e^{\frac{1}{2} - \nu} \int\limits_y^{\infty} I^{\nu -
\frac{1}{2}}  g (x', \tau)  \widetilde{a}_{\nu}  (x', y, \tau) \,
d \tau -{}
$$
$$
{}- I_e^{\frac{1}{2} - \nu} \lr{ y \, \widetilde{a}_{\nu}  (x', y,
\tau) I_e^{\nu - \frac{1}{2}} g} -
  I_e^{\frac{1}{2} - \nu} \int\limits_y^{\infty} \tau
\widetilde{a}_{\nu}  (x', y, \tau) I^{\nu - \frac{1}{2}}  D_{\tau}
\, g (x', \tau) \, d \tau.
$$
Каждое слагаемое, как нетрудно заметить, удовлетворяет нужным
условиям для возможности применения схемы рассуждений из
доказательства предыдущей теоремы. Отсюда и следует справедливость
оценки~\eqref{3.2.13}. Следствие доказано.
\end{proof}

\subsection{Весовые следы}\label{sec10.3}

В этом пункте будет введено понятие весового следа и доказаны прямые теоремы о весовых следах.

Определим как и раньше весовую функцию $\sigma_{\nu} (y)$ по формулам
\begin{equation*}
\sigma_{\nu} (y) = \left\{
\begin{array}{lll}
y^{2 \nu}, & \mbox{если} \  \Re \nu>0, \\
 \dfrac{1}{ - \ln y}, & \mbox{если} \   \nu=0, \\
1, & \mbox{если} \  \Re \nu<0.
\end{array}
\right.
\end{equation*}

Случай мнимых значений параметра $\nu$ считается особым и будет рассмотрен ниже отдельно.

В классическом смысле весовым $\sigma_{\nu}$ "--- следом функции
называется следующий предел:
$$
\left. \sigma_{\nu} f \right|_{y=0} = \lim\limits_{y \to + 0 } \sigma_{\nu} (y) f (x', y) = \psi(x').
$$

Покажем, что у функций  $f \in \mathring{C}_{\nu}^{\infty}
(E_{+}^{n+1})$ $\sigma_{\nu}$ "--- след существует и принадлежит
пространству $\mathring{C}^{\infty} (E^{n}).$ Если   $f \in
\mathring{C}_{\nu}^{\infty} (E_{+}^{n+1}),$ то по определению
этого пространства, существует функция $g \in
\mathring{C}_{\nu}^{\infty} (\ov{E_{+}^{n+1}})$ такая, что
$f=P_{\nu, e} g.$ Тогда $g = S_{\nu, e} f.$ В пункте~\ref{sec4.3}
в одномерном случае была доказана формула
\begin{equation}
\lim\limits_{y \to + 0 } \sigma_{\nu} (y) P_{\nu}^{\frac{1}{2} - \nu } f_1 (y) = \left\{
\begin{array}{ll}
\dfrac{1}{2 \nu} f_1 (0), & \mbox{если} \  \Re \nu>0, \\[1.5ex]
f_1 (0), & \mbox{если} \   \nu=0.
\end{array}
\right.
\label{3.3.1}
\end{equation}
Такой же результат, очевидно, верен и в многомерном случае.
Используя формулу $P_{\nu, e}= P_{\nu}^{\frac{1}{2}-\nu} I_e^{\nu
- \frac{1}{2}},$ отсюда получаем
\begin{equation}
\left. \sigma_{\nu} f \right|_{y=0} = \left\{
\begin{array}{ll}
\dfrac{1}{2 \nu}  I_e^{\nu - \frac{1}{2}} g (x', y) \Big|_{y=0}, & \mbox{если} \  \Re \nu>0, \\
 I_e^{ - \frac{1}{2}} g (x', y) \Big|_{y=0}, & \mbox{если} \
\nu=0.
\end{array}
\right.
\label{3.3.2}
\end{equation}

Таким образом, $\left. \sigma_{\nu} f \right|_{y=0} \in
\mathring{C}^{\infty} (E^n),$ если $f \in
\mathring{C}_{\nu}^{\infty} (E_{+}^{n+1}).$

В связи с формулой~\eqref{3.3.2} нам   надлежит изучить следы
дробных  интегралов  $I^{\mu}_e$ в пространствах $H^s.$  В таких
результатах мы будем нуждаться и в пятом параграфе.

Итак, пусть функция $g \in \mathring{C}^{\infty}
(\ov{E^{n+1}_{+}})$ и  пусть $\mbox{П} g$ "---  её  продолжение по
Уитни на всё пространство $E^{n+1}.$ Тогда для любого комплексного
$\mu$ и натурального $k$ справедлива  формула
$$
\left. D_y^k \, I^{\mu}_e \, g (x', y) \right|_{y=0} =  \left. D^k \, I^{\mu}_e \, \mbox{П}  g (x', y) \right|_{y=0} = \psi_k (x'),
$$
где $\psi_k \in \mathring{C}^{\infty} (E^n).$ В образах Фурье это
равенство имеет вид
$$
\frac{1}{2 \pi} \int\limits_{- \infty}^{\infty} (1 - i \eta)^{- \mu} i^k \eta^k F \mbox{П}  g (\xi', \eta) \, d \eta = F' \psi_k (x'),
$$
где через $F'$ обозначено преобразование Фурье по первым $n$
переменным. По неравенству Коши---Буняковского получаем
\begin{equation}
\| (1 + |\xi'|^2)^s F'  \psi_k (x')  \|_{L_2(E^n)}^2 \leq c \int\limits_{E^n} (1 + |\xi'|^2)^s   \int\limits_{- \infty}^{\infty} \frac{(1+\eta^2)^{- \Re \mu} \eta^{2 k}}{ (1 + |\xi|^2)^{s'}} \, d \eta \int\limits_{- \infty}^{\infty}
(1 + |\xi'|^2)^{s'} \left| F \mbox{П}  g (\xi) \right|^2 \, d \eta d \xi'.
\label{3.3.3}
\end{equation}

\begin{lemma} \label{lem: 3.3.1.}
Пусть $s'>k- \Re \mu + \dfrac{1}{2} > 0.$ Тогда
\begin{equation}
\int\limits_{- \infty}^{\infty} \frac{(1+\eta^2)^{- \Re \mu} \eta^{2 k}}{ (1 + |\xi|^2)^{s'}} \, d \eta \leq c \, (1 + |\xi'|^2)^{k-s'+\frac{1}{2} - \Re \mu},
\label{3.3.4}
\end{equation}
где постоянная $c>0$ не зависит от $\xi' \in E^n.$
\end{lemma}

\begin{proof}
Разобьем интеграл в~\eqref{3.3.4} на два. Для первого из них имеем
$$
 \int\limits_{|\eta|>1} \frac{(1+\eta^2)^{- \Re \mu} \eta^{2 k}}{ (1 + |\xi|^2)^{s'}} \, d \eta \leq c_1  \int\limits_{- \infty}^{\infty} \frac{ \eta^{2 k - \Re \mu}}{ (1 + |\xi|^2)^{s'}} \, d \eta = c_2 \, (1 + |\xi'|^2)^{k-s' - \Re \mu +\frac{1}{2}}.
$$
Второй интеграл допускает элементарную оценку
$$
\int\limits_{|\eta|<1} \frac{(1+\eta^2)^{- \Re \mu} \eta^{2 k}}{
(1 + |\xi|^2)^{s'}} \, d \eta \leq c_3 \, (1 + |\xi'|^2)^{-s'}.
$$
Лемма доказана.
\end{proof}

Соединяя неравенства~\eqref{3.3.3} и~\eqref{3.3.4}, получаем
$$
\| \psi_k \|^2_{H^s(E^n)} \leq c \int\limits_{E^{n+1}}  (1 + |\xi'|^2)^{s-s' + k - \Re \mu +\frac{1}{2}} (1 + |\xi|^2)^{-s'} \left| F \mbox{П}  g (\xi) \right|^2 \, d \xi \leq
$$
$$
\leq c \int\limits_{E^{n+1}}  (1 + |\xi|^2)^{s + k - \Re \mu +\frac{1}{2}}  \left| F \mbox{П}  g (\xi) \right|^2 \, d \xi \leq c \, \| g \|_{H^{s + k - \Re \mu +\frac{1}{2}} (E^{n+1}_{+})}.
$$
Следовательно, доказана

\begin{theorem} \label{teo: 3.3.1}
Пусть $s>k- \Re \mu + \dfrac{1}{2} > 0.$ Тогда для любой функции
$g \in H^s (E_{+}^{n+1})$ на гиперплоскости $y=0$ существует след
функции $D^k_y I_e^{\mu} g,$ который принадлежит пространству
$H^{s - k + \Re \mu - \dfrac{1}{2}} (E^{n}).$ Справедлива оценка
\begin{equation}
\|\left. D^k_y I^{\mu}_e g \right|_{y=0} \|_{H^{s - k + \Re \mu - \frac{1}{2}} (E^{n})} \leq c \, \| g \|_{H^{s} (E^{n+1}_{+})},
\label{3.3.5}
\end{equation}
где постоянная $c>0$ не зависит от функции $g.$
\end{theorem}

Вернемся к изучению весовых следов. Пусть функция $f \in
\mathring{C}_{\nu}^{\infty} (E_{+}^{n+1}).$ Выпишем
формулу~\eqref{3.3.1} для функции $B^k f,$ также принадлежащей
$\mathring{C}_{\nu}^{\infty} (E_{+}^{n+1})$:
\begin{equation}\label{3.3.6}
\begin{aligned}
\left. \sigma_{\nu} B^k f \right|_{y=0} &= \left\{
\begin{array}{ll}
\dfrac{1}{2 \nu} I_e^{\nu - \frac{1}{2}} S_{\nu, e} \, B^k f \Big|_{y=0}, & \mbox{если} \  \Re \nu>0, \\
 I_e^{ - \frac{1}{2}} S_{0, e} B^k f \Big|_{y=0}, &
\mbox{если} \   \nu=0
\end{array}
\right.= \\
&= \left\{
\begin{array}{ll}
\dfrac{1}{2 \nu}  D_y^{2 k} I_e^{\nu - \frac{1}{2}} S_{\nu, e} f \Big|_{y=0}, & \mbox{если} \  \Re \nu>0, \\
 D_y^{2 k} I_e^{ - \frac{1}{2}} S_{0, e}  f \Big|_{y=0},
& \mbox{если} \   \nu=0.
\end{array}
\right.
\end{aligned}
\end{equation}
Здесь была использована формула $S_{\nu, e} B^k= D_y^{2 k} S_{\nu,
e}.$ Соединяя~\eqref{3.3.1} и~\eqref{3.3.6} с теоремой~\ref{teo:
3.3.1} и определением нормы пространства $H^s_{\nu}
(E^{n+1}_{+}),$ приходим к следующей теореме о весовых следах.

\begin{theorem} \label{teo: 3.3.2}
Пусть $\Re \nu > 0$ или $\nu = 0$ и $s> 2 k - \Re \nu +1 >0.$
Тогда отображение $\left. f \to \sigma_{\nu} B^k f \right|_{y=0},$
определяемое по формуле~\eqref{3.3.6} для функций $f \in
\mathring{C}_{\nu}^{\infty} (E_{+}^{n+1}),$ расширяется по
непрерывности до линейного ограниченного отображения пространства
$H^s_{\nu} (E_{+}^{n+1})$ в пространство $H^{s-2 k + \Re \nu
-1}_{\nu} (E^{n}).$ Справедливо неравенство
$$
\|\left. \sigma_{\nu} B^k f \right|_{y=0} \|_{H^{s - 2k + \Re \nu - 1} (E^{n})} \leq c \, \| g \|_{H^{s}_{\nu} (E^{n+1}_{+})},
$$
где постоянная $c>0$ не зависит от функции $f.$
\end{theorem}

Рассмотрим случай мнимых значений параметра $\nu.$  Формула
типа~\eqref{3.3.1} здесь уже не имеет места. В самом деле, для
любой функции $g \in \mathring{C}^{\infty} (\ov{E_{+}^{n+1}})$ при
$\Re \geq 0$ имеем (см.~\eqref{1.3.3})
\begin{equation}
P_{\nu}^{\frac{1}{2}-\nu} g (x', y) = \frac{i y^{-\nu}}{2^{\nu+2}\, \Gamma (\nu+1)} \int\limits_{-\infty}^{\infty} H^{(1)}_{\nu} (y \eta) \eta^{\nu} F \mbox{П} g (x', \eta) \, d \eta.
\label{3.3.7}
\end{equation}
Функция Ханкеля $ H^{(1)}_{\nu}(z)$ в случае $\Re \nu =0, \nu \neq
0,$ в окрестности точки $z=0$ имеет вид~\cite{BE2}
$$
H^{(1)}_{\nu}  (z) = {c'}_{\nu} z^{- \nu} \lr{1+o(1)}+{c''}_{\nu} z^{- \nu} \lr{1+o(1)},
$$
где $c'_{\nu}$ и $c''_{\nu}$ "--- вполне определенные постоянные.
Функции $z^{- \nu}$ и $z^{\nu}$ осциллируют вблизи начала
координат: $z^{- \nu} = \cos (\mu \ln z) - i \sin (\mu \ln z),$
$z^{\nu} = \cos (\mu \ln z) + i \sin (\mu \ln z).$   Поэтому не
существует степенной функции $\varphi (z),$ для которой
$\lim\limits_{z \to 0} H^{(1)}_{\nu}  (z) =1.$ Вследствие этого и
у функции $P_{\nu}^{\frac{1}{2}- \nu} g(x', y),$ вообще говоря, не
существует весовых следов. Поэтому мы будем рассматривать следы не
самой функции, а её производной.

Воспользуемся следующей рекуррентной формулой для функции
$H^{(1)}_{\nu}$ (см.~\cite[с.~20]{BE2}):
$$
\frac{\pr}{z \pr z} \lr{z^{- \nu} H^{(1)}_{\nu}  (z)} = - z^{-\nu-1} H^{(1)}_{\nu+1}  (z).
$$
Тогда из~\eqref{3.3.7} получим
$$
\frac{\pr}{y \pr y} P_{\nu}^{\frac{1}{2}-\nu} g (x', y) = \frac{-i y^{-\nu-1}}{2^{\nu+2} \Gamma (\nu+1)} \int\limits_{-\infty}^{\infty} H^{(1)}_{\nu+1} (y \eta) \eta^{\nu+1} F \mbox{П} g (x', \eta) \, d \eta.
$$
Следовательно,
$$
\frac{\pr}{y \pr y} P_{\nu}^{\frac{1}{2}-\nu} g (x', y) = -2 (\nu+1) P_{\nu}^{-\frac{1}{2}-\nu} g (x', y).
$$
Отсюда и из~\eqref{3.3.1}  получаем
\begin{equation}
\lim\limits_{y \to +0}  y^{2 \nu +1} \frac{\pr}{ \pr y} P_{\nu}^{\frac{1}{2}-\nu} g (x', y) = - g(x', 0).
\label{3.3.8}
\end{equation}

Формула~\eqref{3.3.8}, справедливая при всех $\nu$ из
полуплоскости $\Re \nu \geq 0,$ в   том числе   и для мнимых
$\nu,$ служит как и~\eqref{3.3.1} основой получения следующей
теоремы о весовых следах:

\begin{theorem} \label{teo: 3.3.3}
    Пусть $\Re \nu \geq 0$  и $s> 2 k - \Re \nu +1 >0.$

    Тогда отображение   $\left. f \to y^{2 \nu +1} D_y B^k f \right|_{y=0},$ определяемое по формуле
    $$
    y^{2 \nu +1} \left. D_y B^k f \right|_{y=0} = - D^{2 k}_y I_e^{\nu - \frac{1}{2}} \left. S_{\nu, e} f \right|_{y=0}
    $$
    для функций $f \in \mathring{C}_{\nu}^{\infty} (E_{+}^{n+1}),$ расширяется по непрерывности до  ограниченного отображения пространства $H^s_{\nu} (E_{+}^{n+1})$ в пространство  $H^{s-2 k + \Re \nu -1} (E^{n}).$ Справедливо неравенство
    $$
    \|y^{2 \nu +1} \left. D_y B^k f \right|_{y=0} \|_{H^{s - 2k + \Re \nu - 1} (E^{n})} \leq \|f\|_{H_{\nu}^s (E^{n+1}_{+})},
    $$
    где постоянная $c>0$ не зависит от функции $f.$
\end{theorem}

\section{Функциональные пространства $H^{s}_{\nu}
(\Omega)$}\label{sec11}

\subsection{Разбиение единицы и определения функциональных
пространств}\label{sec11.1}

Пусть $\Omega$ "--- ограниченная область полупространства
$E^{n+1}_{+}$ с границей класса $C^{\infty}.$ Пусть область
$\Omega_0,$ расположенная строго внутри $\Omega$ (это означает,
что $\ov{\Omega}_0 \subset \Omega$ или, что то же самое,
$\ov{\Omega}_0 \cap \pr \Omega = \varnothing$), и области
$\Omega_l,$ $l =1, 2, \dots, \ov{l},$ имеющие непустое пересечение
с границей $\pr \Omega,$ образуют покрытие области $\ov{\Omega}.$
Положим $\Omega_l^{+} = \Omega_l \cap \Omega,$ $l =1, 2, \dots,
\ov{l}.$ Пусть существует диффеоморфизмы $\varkappa_l$ класса
$C^{\infty},$ отображающие области $\Omega_l,$ $l =1, 2, \dots,
\ov{l},$ в области $\omega_l,$  расположенные в пространстве
$E^{n+1}.$ При этом через $\varkappa_0$ обозначим тождественное
преобразование. Пусть область $\Omega_l^{+}$ отображается  в
$\omega_l^{+} = \omega_l \cap E^{n+1}_{+}$   и часть границ
$\Omega_l \cap \pr \Omega$ отображается в часть гиперплоскости
$\omega_l \cap \{ y=0 \}.$ Кроме того, мы предполагаем, что если
$\Omega_l \cap  \Omega_{l'} \neq \varnothing,$ то естественное
отображение $\varkappa_{l'} \varkappa_{l}^{-1}$ части
$\varkappa_{l } \lr{\Omega_l^+ \cap  \Omega_{l'}^+ }$ области
$\omega_l^+$ в часть $\varkappa_{l '} \lr{\Omega_l^+ \cap
\Omega_{l'}^+ }$ области $\omega_{l'}^+$ является   невырожденным
преобразованием с положительным якобианом, осуществляемым по
формулам $x'=x' (\widetilde{x}'),$ $y = \widetilde{y}$ в некоторой
окрестности границы.

В дальнейшем указанное покрытие области считается фиксированным, а все другие необходимые нам покрытия будут получаться из него измельчением.

\begin{lemma}\label{lem: 4.1.1.}
 Существуют функции $h_{l} \in \mathring{C}^{\infty} (E^{n+1})$ обладающие следующими свойствами{\rm :}
    \begin{enumerate}
    \item[1)] $h_{l} (x)=0,$  если  $x \in E^{n+1} \setminus \Omega_{l},$ ${l}=0, \dots, \ov{{l}};$
    \item[2)] $0 \leq h_{l} (x) \leq 1$  при всех $x \in E^n;$
    \item[3)] $h_0 (x)+ \ldots + h_{\ov{l}} (x)=1$  при $x \in \Omega;$
    \item[4)] $D_y h_{l} (x) =0,$ ${l}=0, \dots, \ov{l},$  в локальных координатах в некоторой окрестности гиперплоскости $y=0.$
    \end{enumerate}

\end{lemma}

\begin{proof}
Обозначим через $U_R (x)$ открытый  шар с центром в точке $x \in
E^{n+1}$ радиуса $R> 0.$   Пусть
$$
\omega_{l, \varepsilon} = \{x \in \omega_{l}: U_{\varepsilon}
\subset \omega_{l} \}.
$$
Легко видеть, что при достаточно малом $\varepsilon > 0$ области
$\varkappa_{l}^{-1} \omega_{ {l}, \varepsilon}$ образуют покрытие
области $\ov{\Omega}.$

Пусть $\varepsilon > 0$  как раз таково, что области
$\varkappa_{l}^{-1} \omega_{ {l}, \varepsilon},$ $l = 0, \dots,
\ov{l},$ образуют покрытие области~$\ov{\Omega}.$ Рассмотрим
какую-либо область $\omega_{ {l}, \varepsilon}.$ Обозначим через
$\mathcal{K}_{\delta} (x)$ открытый $n+1$-мерный куб с центром в
точке $x$ с длиной ребра $2 \delta$ и с гранями, параллельными
координатным гиперплоскостям. Поскольку каждая точка $x \in
\omega_{ {l}, \varepsilon} \cap \{y=0\}$ входит в $\omega_l$
вместе с шаром радиуса   $\varepsilon,$ то существует конечное
покрытие её кубами вида $\mathcal{K}_{\varepsilon_p/2} (x^p),$
$p=1, \dots, \ov{p} < \infty,$ $\varepsilon_p > 0,$ с центрами,
соответственно, в точках $x^p \in \omega_{ {l}, \varepsilon} \cap
\{y=0\},$ причём такое, что  $\mathcal{K}_{\varepsilon_p} (x^p)
\subset \omega_l.$ Покроем оставшуюся часть области $\omega_{ l,
\varepsilon}$ кубами    вида  $\mathcal{K}_{\varepsilon_p/2}
(x^p),$ $p= \ov{p}+1, \dots, \ov{\ov{p}} < \infty,$ и такими,
чтобы $\mathcal{K}_{\varepsilon_p} (x^p) \subset \omega_l$ и чтобы
$\ov{\mathcal{K}}_{\varepsilon_p} (x^p)$ не имело общих точек с
гиперплоскостью $\{y=0\}.$

Обозначим через $\varphi_{\delta} (t)$ функцию одной переменной
такую, что $\varphi_{\delta} \in \mathring{C}^{\infty} (E^1),$
$\varphi_{\delta} (t) \geq 0$ при $t \in E^1,$ $\varphi_{\delta}
(t)=1$ при  $|t| \leq \delta,$  $\varphi_{\delta} (t)=0$ при $|t|
\geq 2 \delta.$

Введем теперь функцию $\psi_{l}(x),$ $x= (x', y)$ по    формуле
$$
\psi_{l}(x) = \sum\limits_{p=1}^{\ov{p}} \varphi_{\varepsilon_p}
(y) \prod\limits_{q=1}^{n} \varphi_{\varepsilon_p} (x_q-x_q^p) +
\sum\limits_{p=\ov{p}+1}^{\ov{\ov{p}}} \psi_{{l}, p} (x),
$$
где $x_q^p$ обозначает  $q$-ю координату точки $x^p,$ а через
$\psi_{l, p}(x)$ обозначены функции, обладающие свойствами:
$\psi_{{l}, p} (x) \in \mathring{C}^{\infty} (E^{n+1}),$ $\psi_{l,
p}(x) \geq 0$ при $x \in E^{n+1},$ $\psi_{l, p}(x)=1$ при $x \in
\mathcal{K}_{\varepsilon_p/2} (x^p),$ $\psi_{l, p}(x)=0$ при  $x
\in E^{n+1} \setminus \mathcal{K}_{\varepsilon_p} (x^p),$ $p=
\ov{p}+1, \dots, \ov{\bar{p}}.$

Доказательство существования функций $\varphi_{\delta}$ и
$\psi_{l, p},$ удовлетворяющих всем перечисленным условиям, можно
найти, например, в книге~\cite{65}. Тогда построенная функция
$\psi_l (x)$ обладает свойствами: $\psi_{l} (x) \in
\mathring{C}^{\infty} (E^{n+1}),$ $\psi_{l} (x) \geq 0$ при $x \in
E^{n+1},$ $\psi_{l} (x) \geq 1$ при $x \in \omega_{l,
\varepsilon}$ и $\supp \psi_l \subset \omega_l.$ Нетрудно  также
заметить, что $D_y \psi_l (x) \equiv 0$ при $|y| < \delta,$ где
$\delta$ "--- некоторое положительное число.

Искомые функции $h_l$ определим по формуле
$$
h_{l} (x) = \frac{\varkappa_{l}^{-1}
\psi_{l}}{\sum\limits_{m=0}^{\ov{l}} \varkappa_{l}^{-1} \psi_m}.
$$
Проверка выполнения условий 1)--4) не вызывает затруднения. Лемма
доказана.
\end{proof}

 В дальнейшем мы будем говорить, что набор функций $\{h_l\}$ образует разбиение единицы, {\it подчиненное покрытию} $\{\Omega_l\},$ если выполнены условия 1)--4).

Введем теперь пространство $H^s_{\nu} (\Omega),$ $s \geq 0,$ $\Re \nu \geq 0$ как множество функций $f,$ определенных в области
$\Omega,$ таких, что функции $\varkappa_l (h_l f)$ принадлежат
пространству $H^s_{\nu} (E^{n+1}_{+}).$ Это пространство
становится гильбертовым по норме
\begin{equation}
\| f \|_{H^s_{\nu} (\Omega)} = \lr{ \sum\limits_{l=0}^{\ov{l}} \| \varkappa_l (h_l f) \|^2_{H^s_{\nu} (E^{n+1}_{+})}}^{\frac{1}{2}}.
\label{4.1.1}
\end{equation}
Здесь $\{h_l\}$ "--- разбиение единицы, подчиненное покрытию
$\{\Omega_l\}.$

Покажем, что нормы~\eqref{4.1.1} при различных выборах разбиений
единицы эквивалентны. Получим даже более общий результат. Пусть
$\{{\Omega'}_{l'}\}$ другое покрытие области $\ov{\Omega}.$ Пусть
для каждого $l'=0, \dots, \ov{l'}$ найдется такое    $l=0, \dots,
\ov{l},$ что ${\Omega'}_{l'} \subset {\Omega}_{l}.$ При выполнении
этого условия покрытие $\{{\Omega'}_{l'}\}$ будем называть
измельчением покрытия $\{{\Omega}_{l}\}.$ Обозначим через
$\{h'_{l'}\}$ разбиение единицы, подчиненное покрытию
$\{{\Omega'}_{l'}\}.$ Тогда можно ввести следующую норму:
\begin{equation}
\| f \|_{H^s_{\nu} (\Omega)} = \lr{ \sum\limits_{l'=0}^{\ov{l'}} \| \varkappa'_{l'} (h'_{l'} f) \|^2_{H^s_{\nu} (E^{n+1}_{+})}}^{\frac{1}{2}},
\label{4.1.2}
\end{equation}
где отображение $\varkappa'_{l'}=\varkappa_l,$ $l=l(l'),$ причём
$l$ выбирается таким образом, чтобы $\supp h'_{l'} \subset
\Omega_l.$ Если   таких $l$ несколько, то выбирается любое из них
(из свойств отображения $\varkappa_l$ следует, что соответствующие
слагаемые будут эквивалентными нормами).

Оценим одно из слагаемых справа в~\eqref{4.1.2}. Имеем
$$
\| \varkappa'_{l'} (h'_{l'} f) \|_{H^s_{\nu} (E^{n+1}_{+})} = \| \sum\limits_{l=0}^{\ov{l}} \varkappa'_{l'} (h'_{l'} h_{l} f) \|_{H^s_{\nu} (E^{n+1}_{+})} \leq
  \sum\limits_{l=0}^{\ov{l}} \| \varkappa'_{l'} (h'_{l'} h_{l} f) \|_{H^s_{\nu} (E^{n+1}_{+})} \leq
$$
$$
\leq c_1 \sum\limits_{l=0}^{\ov{l}} \| \varkappa_{l} (h'_{l'} h_{l} f) \|_{H^s_{\nu} (E^{n+1}_{+})}
\leq c_2 \sum\limits_{l=0}^{\ov{l}} \| \varkappa_{l} (h_{l} f) \|_{H^s_{\nu} (E^{n+1}_{+})}.
$$
При этом в предпоследнем неравенстве мы воспользовались свойствами
отображения $\varkappa'_{l'} \varkappa^{-1}_{l}$ (см. начало
пункта), а в последнем неравенстве "--- результатами о
мультипликаторах из пункта~\ref{sec3.2}. Здесь и потребовалось
свойство~4 функций ${h'}_{l'}.$ Обратная оценка доказывается
аналогично. Таким образом, эквивалентность норм~\eqref{4.1.1}
и~\eqref{4.1.2} установлена.

Рассмотрим еще один вопрос, касающийся структуры пространства
$H^s_{\nu} (\Omega).$ Пусть функция $f \in H^s_{\nu} (\Omega),$ $a
\in \mathring{C}^{\infty} (\Omega).$ Тогда из результатов
пункта~\ref{sec3.1} следует, что функция $a f$ принадлежит
пространству $H^s (\Omega).$ Более того, справедлива оценка
\begin{equation}
c_1  \| a f \|_{H^s (\Omega)} \leq   \| a f \|_{H^s_{\nu} (\Omega)} \leq c_2  \| a f \|_{H^s (\Omega)},
\label{4.1.3}
\end{equation}
где постоянные $c_1, c_2 > 0$ не зависят от функции $f.$

\begin{lemma} \label{lem: 4.1.2}
    Пусть $\widetilde{\Omega}$ "--- строго внутренняя подобласть области $\Omega.$ Тогда сужение функции $f \in H^s_{\nu} (\Omega)$ на
    $\widetilde{\Omega}$ принадлежит  пространству $H^s (\widetilde{\Omega})$ и справедлива оценка
    \begin{equation}
        \| f \|_{H^s (\widetilde{\Omega})} \leq c \,  \|  f \|_{H^s_{\nu} (\Omega)},
        \label{4.1.4}
    \end{equation}
    в которой постоянная $c > 0$ не зависит от функции $f.$
\end{lemma}

\begin{proof}
Введём  функцию  $a \in \mathring{C}^{\infty} (\Omega)$ такую, что
$a (x) = 1$ при  $x \in \widetilde{\Omega}.$ Тогда из
неравенства~\eqref{4.1.3} имеем
$$
\| f \|_{H^s (\widetilde{\Omega})} \leq   \| a f \|_{H^s (\Omega)}
\leq  c \, \| a f \|_{H^s_{\nu} (\Omega)}.
$$
По теореме~\ref{teo: 3.2.1} о мультипликаторах получаем оценку
$$
\| a f \|_{H^s_{\nu} (\Omega)} \leq  c \, \|  f \|_{H^s_{\nu}
(\Omega)},
$$
где $c > 0$ не зависит от функции $f.$  Лемма доказана.
\end{proof}

Пространства $H^s_{\nu, +} (\Omega),$ $s \geq 0,$ $\nu \geq 0,$
введённые И.\,А.~Киприяновым~\cite{Kip1}, определяются как
замыкание множества всех функций из $C^{\infty} (\ov{\Omega}),$
для которой в каждой л.с.к. $D_y^{2 k +1} f = 0$ при $y=0,$ $k=0,
1, \dots$ по норме
\begin{equation}
\| f \|_{H^s_{\nu, +} (\Omega)} =\lr{ \sum\limits_{l=0}^{\ov{l}} \| \varkappa_{l} (h_{l} f) \|^2_{H^s_{\nu, +} (E^{n+1}_{+})}}^{\frac{1}{2}},
\label{4.1.5}
\end{equation}
где нормы в $H^s_{\nu, +} (E_{+}^{n+1})$ определены в
пункте~\ref{sec3.1}.

\begin{lemma} \label{lem: 4.1.3}
    Пусть $s \geq 0$ и вещественное $\nu \geq 0,$ $\nu \neq 1, 3, 5, \dots.$ Тогда пространство $H^s_{\nu, +} (\Omega)$ является собственным пространством пространства $H^s_{\nu} (\Omega),$ причём нормы~\eqref{4.1.1} и~\eqref{4.1.5} на $H^s_{\nu, +} (\Omega)$ эквивалентны. При $\nu = 1, 3,  \dots$ пространство $H^s_{\nu, +} (\Omega)$ вложено в $H^s_{\nu} (\Omega).$
\end{lemma}

\begin{proof}
Ввиду финитности функций $h_l f,$ участвующих в определении нормы,
мы можем применить лемму~\ref{lem: 3.1.4}. Лемма доказана.
\end{proof}

\subsection{Теоремы вложения}\label{sec11.2}

Рассмотрим сначала внутренние теоремы вложения. Следующее
утверждение просто вытекает из теоремы~\ref{teo: 3.1.1}.

\begin{theorem} \label{teo: 4.2.1}
    Пусть $s'>s \geq 0$ и $\Re \nu \geq 0.$ Тогда пространство $H^{s'}_{\nu} (\Omega)$  вложено в $H^{s}_{\nu} (\Omega)$ и справедливо неравенство
    $$
        \|  f \|_{H^s_{\nu} (\Omega)} \leq  c \, \|  f \|_{H^{s'}_{\nu} (\Omega)},
    $$
    где  $f \in H^{s'}_{\nu} (\Omega)$  и постоянная $c>0$ не зависит от $f.$
\end{theorem}

\begin{theorem} \label{teo: 4.2.2}
    При $s'>s \geq 0$ и $\Re \nu \geq 0$ оператор вложения пространства $H^{s'}_{\nu} (\Omega)$ в $H^{s}_{\nu} (\Omega)$ вполне непрерывен.
\end{theorem}

\begin{proof}
Пусть последовательность функций $f_j,$ $j=1,2, \dots$ ограничена
по норме пространства $H^{s}_{\nu} (\Omega).$ Обозначим через
$\mathsf{y}_j,$ $j=1,2, \dots$ расширяющуюся систему строго
внутренних подобластей области $\Omega,$ объединение которых
совпадает с $\Omega.$ По лемме~\ref{lem: 4.1.2} последовательность
$f_j$ ограничена в смысле нормы каждого пространства $H^{s'}
(\mathsf{y}_j),$ поэтому из неё можно выделить сходящуюся по норме
пространства $H^{s} (\mathsf{y}_1)$ последовательность $f_{1j}.$
Обозначим предел этой подпоследовательности через $f^1.$ Из
$f_{1j}$ выделим подпоследовательность $f_{2j},$ сходящуюся к
$f^2$ по норме $H^{s} (\mathsf{y}_2).$ Ясно, что $f^2|_\mathsf{y}
= f^1.$ Продолжим этот процесс далее. Тогда мы найдем функцию $f$
такую, что при каждом $j$ сужение $f|_{\mathsf{y}_j} \in H^s
(\mathsf{y}_j)$ и   $f|_{\mathsf{y}_j}.$ Образуем диагональную
последовательность $f_{jj}.$ Для неё
$$
\| \left. (f-f_{jj}) \right|_{\mathsf{y}_p} \|_{H^s
(\mathsf{y}_p)} \to 0, ~ j \to \infty, ~  p  = 1,2, \dots.
$$
Не ограничивая общности будем считать, что сама $fj$ обладает этим
свойством.

Далее, последовательность $\{\varkappa_l (h_l f_j) \} \subset
H^{s'}_{\nu} (\omega_l^{+})$ и ограничена по норме этого
пространства. Тогда по теореме~\ref{teo: 3.1.2} найдется функция
$g_l \in H^{s}_{\nu} (\omega_l^{+}),$ к которой сходится некоторая
подпоследовательность последовательности $\varkappa_l (h_l f_j).$
Опять не ограничивая общности считаем, что этим свойством при всех
$l=0,1,\dots, \ov{l}$ обладает сама последовательность $f_j.$
Пусть $\varepsilon>0.$ Обозначим через $\omega^{+}_{l,
\varepsilon}$ часть области $\omega_l^{+},$ расположенную в
полупространстве $y > \varepsilon.$ Тогда по  лемме~\ref{lem:
4.1.2} последовательность $\varkappa_l (h_l f_j)|_{\{y >
\varepsilon\}}$ сходится к функции $\varkappa_l (h_l f)|_{\{y >
\varepsilon\}}$    по  норме пространства $H^s_{\nu}
(\omega^{+}_{l, \varepsilon}).$ Ввиду произвольности $\varepsilon$
получаем, что $g_l= \varkappa_l (h_l f)$   при всех $l.$ А тогда
$f \in H^s_{\nu} (\Omega),$ что и завершает доказательство
теоремы.
\end{proof}

Рассмотрим теперь утверждение о весовых следах. Нам потребуются
граничные пространства $H^s (\pr \Omega),$ описание которых можно
найти в книге С.\,М.~Никольского~\cite{66}. Одна из эквивалентных
норм пространства $H^s (\pr \Omega)$ может быть определена по
формуле
$$
\| \varphi \|^2_{H^s(\pr \Omega)} = \sum\limits_{l=0}^{\ov{l}} \| \varkappa_{l} (h_{l} f) \|_{H^s (E^{n})}.
$$
Обозначим через $\sigma_{\nu}$ функцию, которая в каждой л.с.к.
совпадает с ранее введённой весовой функцией $\sigma_{\nu} (y).$

Будем говорить, что функция $H^s_{\nu} (\Omega)$ обладает весовым
следом $\left. \sigma_{\nu} f \right|_{\pr \Omega} = \varphi$ на
границе  $\pr \Omega,$ если в каждой л.с.к. функции $h_l f,$
$l=0,1,\dots, \ov{l}$ обладают весовым $\sigma_{\nu}$-следом.

Обозначим через $B_{\nu}$ дифференциальный оператор, который в каждой л.с.к. совпадает с оператором Бесселя.

\begin{theorem}\label{teo: 4.2.3}
    Пусть $\Re \nu > 0$ или $\nu=0,$  $s>2k+1- \Re \nu > 0.$ Тогда для любой функции $f \in H^s_{\nu} (\Omega)$ существует  весовой след функции $B_{\nu}^k f$ и имеет место неравенство{\rm :}
$$
\| \left. \sigma_{\nu} B^k_{\nu} f \right|_{\pr \Omega} \|_{H^{s-2k-1+ \Re \nu} (\pr \Omega)} \leq c \, \| f \|_{H^s_{\nu} (\Omega)},
$$
где постоянная $c>0$ не зависит от функции $f.$

Если $\Re \nu \geq 0$ и $s>2k+1- \Re \nu > 0,$ то для любой
функции  $f \in H^s_{\nu} (\Omega)$ существуют весовые следы
$\left. \sigma_{\nu + \frac{1}{2}} D B^k_{\nu} f \right|_{\pr
\Omega}$ функций $D B^k f$ и выполнены неравенства{\rm :}
$$
\| \left. \sigma_{\nu + \frac{1}{2}} D B^k_{\nu} f \right|_{\pr \Omega} \|_{H^{s-2k-1+ \Re \nu} (\pr \Omega)} \leq c \, \| f \|_{H^s_{\nu} (\Omega)},
$$
где постоянная $c > 0$ не зависит от функции $f$ и через $D$
обозначен оператор, совпадающий с $\dfrac{\pr }{\pr y}$ в каждой
л.с.к. вблизи границы.
\end{theorem}

Это утверждение есть следствие теорем~\ref{teo: 3.3.2} и~\ref{teo:
3.3.3}.

\section{Эллиптическая краевая задача в полупространстве с нелокальными краевыми условиями дробного
порядка}\label{sec12}

\subsection{Постановка краевой задачи}\label{sec12.1}

Рассмотрим следующий однородный оператор в частных производных с
постоянными коэффициентами:
$$
A \lr{\frac{1}{i} D} = \sum\limits_{|\alpha|=2 m} a_{\alpha}
\lr{\frac{1}{i} D}^{\alpha},
$$
где $\alpha = (\alpha_1, \dots, \alpha_{n+1}) = (\alpha',
\alpha_{n+1}),$ $\alpha_j \geq 0$ "--- целые числа, $|\alpha| =
\alpha_1+ \dots+ \alpha_{n+1},$
$$
\lr{\frac{1}{i} D}^{\alpha} = \frac{\pr^{|\alpha|}}{(i \pr
x_1)^{\alpha_1} \dots (i \pr x_{n+1})^{\alpha_{n+1}}}.
$$

Оператор $A$ называется эллиптическим, если
$$
A(\xi) = \sum\limits_{|\alpha|=2 m} a_{\alpha}  \xi^{\alpha} \neq
0
$$
при всех $\xi \in E^{n+1},$ $\xi \neq 0.$  Здесь $\xi^{\alpha} =
\xi^{\alpha_1}_1 \dots \xi^{\alpha_{n+1}}_{n+1}.$

Пусть выполнено условие собственной эллиптичности:
характеристический полином $A (\xi) = A (\xi', \eta)$ комплексного
переменного $\eta$ для всякого $\xi' \in E^n,$ $\xi \neq 0,$ имеет
$m$ корней с учётом их кратностей $\eta_j^{+},$ $j=1, \dots, m,$ с
положительной мнимой частью.

Пусть $G_j \Big(\dfrac{1}{i} D\Big),$ $j=0, \dots, m-1$ граничные
операторы с постоянными коэффициентами, причём $G_j (\xi)$ "---
однородный полином степени $m_j \geq 0.$

Пусть выполнено условие Шапиро---Лопатинского: для любого $0 \neq
\xi' \in E^n$ полиномы $G_j (\xi', \eta)$ переменного $\eta$
линейно независимы по модулю полинома $A^{+} (\xi', \eta) = (\eta
- \eta_1^{+} (\xi ')) \dots (\eta - \eta_m^{+} (\xi ')).$

При выполнении перечисленных условий будем говорить, что операторы
$\{A, G_j\},$ а также соответствующие им полиномы (символы),
образуют эллиптический набор.

Рассмотрим краевую задачу следующего вида:
\begin{equation}
\begin{cases}
A \lr{\dfrac{1}{i} D} u(x) = f (x), &  x \in E^{n+1}_{+}, \\
G_j \lr{\dfrac{1}{i} D} \left.  D^k_y I^{\mu}_e u \right|_{y=0} =
g_j (x'), &  x' \in E^n, \  j=0,\dots, m-1.
\end{cases}
\label{5.1.1}
\end{equation}
Здесь $k$ "--- целое неотрицательное число, a $I^{\mu}_e$
обозначает как обычно оператор лиувиллевского типа, введённый в
пункте~\ref{sec4.2} и действующий по последней переменной $y.$ Его
<<порядок>> "--- $\mu$ может быть комплексным числом.

После применения преобразования Фурье по первым $n$ переменным, определяемого по формуле
$$
\widetilde{f} (\xi', y) = F' f (\xi', y) = \int\limits_{E^n} e^{- i \langle x', \xi' \rangle} f (x', y) \, d x',
$$
краевая задача~\eqref{5.1.1} принимает вид
\begin{equation}
\begin{cases}
A \lr{\xi', \dfrac{1}{i} D_y}  \widetilde{u} (\xi', y) = \widetilde{f} (\xi', y), &  \xi' \in E^{n}, \  y>0, \\
G_j \lr{\xi', \dfrac{1}{i} D_y}  \left.  D^k I^{\mu}_e
\widetilde{u} \right|_{y=0} = g_j (\xi'), &  \xi' \in E^n, \
j=0,\dots, m-1.
\end{cases}
\label{5.1.2}
\end{equation}

\begin{lemma} \label{lem:5.1.1}
 Пусть операторы $A$ и $G_j,$ $j=0,\dots, m-1,$ образуют эллиптический набор. Пусть выполнены следующие соотношения{\rm :}
\begin{equation}
k - \Re \mu + \frac{1}{2} > 0, \ 2 m +s -k + \Re \mu - \frac{1}{2} - \max\limits_j m_j >0.
\label{5.1.3}
\end{equation}
Пусть функция $\widetilde{f} (\xi', y)$ по последней переменной
принадлежит пространству $H^s (E^1_{+}).$ Тогда краевая
задача~\eqref{5.1.2} при любом, $\xi' \neq 0,$ $\xi' \in E^n$
имеет в пространстве $H^{2m+s} (E_{+}^1)$ единственное решение и
для него справедлива оценка
\begin{multline}
c \, \sum\limits_{l=0}^{2 m + s} |\xi'|^{2(2m+s-l)} \| D^{l}_y \widetilde{u} (\xi', y) \|_{L_2 (E^1_+)} \leq \\
\leq  \sum\limits_{l=0}^{s} |\xi'|^{2(s-l)} \| D^{l}_y
\widetilde{f} (\xi', y) \|_{L_2 (E^1_+)} + \sum\limits_{j=0}^{ m -
1} |\xi'|^{2(2m+s-m_j-k + \Re \mu - \frac{1}{2})}
|\widetilde{g}_j (\xi')|^2, \label{5.1.4}
\end{multline}
где постоянная $c>0$ не зависит от $f$ и $g_j.$ Если же $|\xi'|
\geq 1,$ то эта постоянная не зависит от $\xi'.$
\end{lemma}

Эту лемму можно считать известной, хотя авторы затрудняются дать
конкретную ссылку. Доказательство её стандартно и легко может быть
восстановлено по следующей схеме. После нахождения частного
решения уравнения с помощью продолжения и применения
преобразования Фурье по последней переменной, нам останется
рассмотреть лишь краевую задачу с однородным уравнением. Применяя
к уравнению оператор $D_y^k I_e^{\mu}$ для функции $D_y^k
I_e^{\mu} \widetilde{u}$ получим обычную краевую задачу, решение
которой дается с помощью контурного интеграла. Затем применяется
обратный оператор ($D_y^k I_e^{\mu}$ на пространстве
быстроубывающих функций обратим) и дается прямая оценка
полученного выражения.

\subsection{Регуляризатор и априорные оценки}\label{sec12.2}

Краевая задача~\eqref{5.1.1}, порождает следующий оператор:
\begin{equation*}
\mathfrak{U}: u \to \mathfrak{U} u = \left\{ A  u, G_0  D_y^k I_e^{\mu} u |_{y=0}, \dots, G_{m-1}  D_y^k I_e^{\mu} u |_{y=0}  \right\}.
\end{equation*}
Пусть
\begin{equation*}
\mathcal{H}^s  \lr{E_{+}^{n+1}, E^n, m} = H^s (E_{+}^{n+1})
\times \prod\limits_{j=0}^{m-1} H^{s+2m-m_j-k+\Re \mu -
\frac{1}{2}} \lr{E^n}.
\end{equation*}

В $\mathcal{H}^s$ зададим   топологию прямого произведения.

Из теоремы~\ref{3.3.1} следует, что оператор $\mathfrak{U}$
непрерывно отображает пространство $H^s (E_{+}^{n+1})$   в
$\mathcal{H}^s \lr{E_{+}^{n+1}, E^n, m}.$

Пусть $\Phi = \{f(x), g_0 (x'), \dots,  g_{m-1} (x')\}$ элемент
пространства $\mathcal{H}^s.$ Левым регуляризатором для
$\mathfrak{U}$ называется оператор $\mathfrak{R}_{\mbox{л}}:
\mathcal{H}^s \to  H^{s+m}$ для которого имеет место формула
$$
\mathfrak{R}_{\mbox{л}} \mathfrak{U} u = u + T_{\mbox{л}} u, \  u \in H^{s+2 m } (E^{n+1}_+),
$$
в которой ${T}_{\mbox{л}}$ "---  сглаживающий оператор,
${T}_{\mbox{л}}: H^{s+2m} \to  H^{s+2m+1}.$

Оператор $\mathfrak{R}_{\mbox{п}}: \mathcal{H}^s \to  H^{s+2m}$ называется правым регуляризатором, если
$$
  \mathfrak{U}  \mathfrak{R}_{\mbox{п}} \Phi = \Phi + T_{\mbox{п}} \Phi,
$$
где ${T}_{\mbox{п}}: \mathcal{H}^s \to  \mathcal{H}^{s+1}$ "---
сглаживающий оператор.

Оператор $\mathfrak{R}$ называется двухсторонним регуляризатором или просто регуляризатором, если он одновременно и правый и левый регуляризатор.

Основной целью настоящего пункта является построение
регуляризатора для введённого выше оператора $\mathfrak{U}.$ При
этом мы несколько изменили схему построения регуляризатора,
примененную Л.~Хермандером~\cite{83}. Искомый регуляризатор мы
построим следующим образом. Пусть функция $U (\xi', y)$ при
$|\xi'| \geq 1$ является единственным решением краевой задачи
\begin{equation}
    \begin{cases}
        A \lr{\xi', \dfrac{1}{i} D_y}  {U} (\xi', y) = \widetilde{f} (\xi', y),  &  y>0, \\
        G_j \lr{\xi', \dfrac{1}{i} D_y}  \left.  D^k I^{\mu}_e {U} \right|_{y=0} = \widetilde{g}_j (\xi'), &  j=0,\dots, m-1, \\
        \widetilde{f} (\xi', y) = F' f, \quad  \widetilde{g}_j (\xi') = F' g_j,
    \end{cases}
    \label{5.2.1}
\end{equation}
где через $F'$ обозначен оператор Фурье по первым $n$ переменным. При $|\xi'|<1 $ функция $U (\xi', y)$ является решением задачи
\begin{equation}
\begin{cases}
A \lr{\xi'_0, \dfrac{1}{i} D_y}  {U} (\xi', y) = \widetilde{f} (\xi', y),  &  y>0, \\
G_j \lr{\xi'_0, \dfrac{1}{i} D_y}  \left.  D^k I^{\mu}_e {U}
\right|_{y=0} = \widetilde{g}_j (\xi'), &  j=0,\dots, m-1,
\end{cases}
\label{5.2.2}
\end{equation}
причём $\xi'_0 \in E^n$ "--- произвольная фиксированная  точка
единичной сферы, $|\xi'_0|=1.$ Заметим, что отличие от схемы
Л.~Хермандера построения регуляризатора заключается как раз в
способе задания функции $U (\xi', y)$ при $|\xi'| < 1.$

Теперь определим оператор $\mathfrak{R}$ по формуле
\begin{equation}
\mathfrak{R} \Phi = \lr{F'}^{-1} {U},
\label{5.2.3}
\end{equation}
где $\Phi = \{f, g_0, \dots, g_{m-1}\}.$

Покажем, что введённый оператор $\mathfrak{R}$ непрерывно
отображает пространство $\mathcal{H}^s (E_{+}^{n+1}, E^n, m)$ в
$H^{s+2m} (E_{+}^{n+1}).$ Имеем
$$
\| \mathfrak{R} \Phi  \|_{H^{s+2m} (E^{n+1}_{+})} = \| (F')^{-1}  U  \|^2_{H^{s+2m} (E^{n+1}_{+})} \leq
$$
$$
\leq c \lr{ \| (F')^{-1}  U  \|_{L_2 (E^{n+1}_{+})} +
\sum\limits_{l=0}^{2m+s} \sum\limits_{|\alpha'|=2m+s-l}
\|D_{x'}^{\alpha'} D_y^{l} (F')^{-1}  U  \|_{L_2 (E^{n+1}_{+})}^2}
\leq
$$
$$
\leq c \lr{ \|  U (\xi', y)  \|^2_{L_2 (E^{n+1}_{+})} +
\sum\limits_{l=0}^{2m+s} \sum\limits_{|\alpha'|=2m+s-l}
\|{\xi'}^{\alpha'} D_y^{l}   U (\xi', y)  \|_{L_2
(E^{n+1}_{+})}^2}.
$$
Поскольку функция $U (\xi', y)$ определяется различно при $|\xi'|
\geq 1$ и $|\xi'| < 1,$ то рассмотрим сначала первый случай. Имеем
$$
I_1 = \int\limits_{|\xi'|>1} \int\limits_0^{\infty} \lr{ |U (\xi',
y)|^2+ \sum\limits_{l=0}^{2m+s} \sum\limits_{|\alpha'|=2m+s-l}
|\xi_1^{\alpha_1} \dots \xi_n^{\alpha_n}  D_y^{l}   U (\xi',
y)|^2} \, d \xi' dy \leq
$$
$$
\leq c \sum\limits_{l=0}^{2m+s} \int\limits_{|\xi'|>1}
\int\limits_0^{\infty} (1+|\xi'|^2)^{2m+s-l}  | D_y^{l}   U (\xi',
y)|^2 \, dy d \xi'.
$$
Применяя лемму~\ref{lem:5.1.1}, получаем
$$
I_1 \leq c \int\limits_{|\xi'|>1} \left( \sum\limits_{l=0}^{s}
(1+|\xi'|^2)^{s-l}  \int\limits_0^{\infty} | D_y^{l}   F' f (\xi',
y)|^2 \, dy  \right. +
$$
$$
\left. + \sum\limits_{j=0}^{m-1} (1+|\xi'|^2)^{2m+s-m_j-k + \Re \mu - \frac{1}{2}} |F' g_j (\xi')|^2  \right) \, d \xi' \leq c \,
\|\Phi \|^2_{\mathcal{H}^s},
$$
где постоянная $c>0$ не зависит от $\Phi.$

Рассмотрим случай $|\xi'| < 1.$ Из определения функции
(см.~\eqref{5.2.2}) получаем
$$
I_2  = \int\limits_{|\xi'|<1} \int\limits_0^{\infty} \lr{ |U
(\xi', y)|^2+ \sum\limits_{l=0}^{2m+s}
\sum\limits_{|\alpha'|=2m+s-l}\hspace{-3ex}  |{\xi'}^{\alpha'}
D_y^{l} U (\xi', y)|^2} \, dy d \xi'  \leq
 c \sum\limits_{l=0}^{2m+s} \int\limits_{|\xi'|<1} \int\limits_0^{\infty} | D_y^{l}   U (\xi', y)|^2 \, dy d \xi'.
$$
Применяя еще раз лемму~\ref{lem:5.1.1}, получим
$$
I_2 \leq c \sum\limits_{l=0}^{s} \int\limits_{|\xi'|<1} \lr{
\int\limits_0^{\infty}  | D_y^{l}   F' (\xi', y) |^2 \, dy +
\sum\limits_{l=0}^{m-1} |F' g_j (\xi') |^2 } d \xi'.
$$
Отсюда следует оценка $I_2 \leq c \, \| \Phi
\|^2_{\mathcal{H}^s}.$ Соединяя оценки для $I_1$ и $I_2,$ получаем
$$
\| \mathfrak{R} \Phi \|_{H^{s+2m} (E_{+}^{n+1})} \leq c \, \| \Phi \|_{\mathcal{H}^s (E_{+}^{n+1}, E^n, m)}.
$$
Таким образом, непрерывность оператора $\mathfrak{R}$ доказана.

Проверим, что оператор $\mathfrak{R}$ является  левым регуляризатором. Имеем
\begin{equation}
\mathfrak{R} \mathfrak{U} u =u + \frac{1}{(2 \pi)^n} \int\limits_{|\xi'|<1} e^{i \langle x', \xi' \rangle} V_{(\xi', y)} \, d \xi' = u + T_{\mbox{л}} u,
\label{5.2.4}
\end{equation}
где функция $V (\xi', y)$ является решением краевой задачи
\begin{equation}
    \begin{cases}
        A \lr{\xi'_0, \dfrac{1}{i} D_y}  V (\xi', y) {=} \left[ A \lr{\xi'_0, \dfrac{1}{i} D_y} -A \lr{\xi', \dfrac{1}{i} D_y}   \right] F'u(\xi', y), &  y>0, \\
        G_j \lr{\xi'_0, \dfrac{1}{i} D_y}  \left.  D^k_y I^{\mu}_e  V \right|_{y=0} {=}\left[ G_j \lr{\xi'_0, \dfrac{1}{i} D_y} {-} G_j \lr{\xi', \dfrac{1}{i} D_y}   \right]  D^k_y I^{\mu}_e   \left. F'u(\xi', y) \right|_{y=0}, & j=\overline{0,\ m-1}.
    \end{cases}
    \label{5.2.5}
\end{equation}
Здесь $T_{\mbox{л}} u = \lr{F'}^{-1} V_1,$ где  $V_1 (\xi', y) =
V(\xi', y),$ если $|\xi'|<1$ и $V_1 (\xi', y)=0,$ если $|\xi'|
\geq 1.$ Отсюда
$$
\| T_{\mbox{л}} u \|_{H^{s+2m} (E^{n+1}_{+})} \leq c
\sum\limits_{l=0}^{s+2m} \int\limits_{|\xi'|<1}
\int\limits_0^{\infty} | D^{l}_y V(\xi', y)|^2 \, dy d\xi'.
$$
Применяя оценку леммы~\ref{lem:5.1.1} к функции $V$ и замечая, что
операторы, стоящие в квадратных скобках в~\eqref{5.2.5}, имеют
одинаковую главную часть, получаем
$$
\| T_{\mbox{л}} u \|_{H^{s+2m} (E^{n+1}_{+})} \leq c \left(
\sum\limits_{l=0}^{s} \int\limits_{|\xi'|<1}  \| D^{l+2m-1}_y  F'u
(\xi', y) \|^2_{L_2 (E^1_{+})} \,  d\xi'+ \right.
$$
$$
+ \left. \sum\limits_{l=0}^{m_j-1} \sum\limits_{j=0}^{s-1} \int\limits_{|\xi'|<1} |\left.  D^{l+k}_y I^{\mu}_e  F' u (\xi', y) \right|_{y=0}|^2  \, d \xi' \right).
$$
Первая сумма, очевидно, оценивается нормой $\| u \|_{H^{s+2m+1}
(E^{n+1}_{+})}.$ Вторая сумма также может быть оценена той же
нормой, как это следует из теоремы~\ref{teo: 3.3.1}. Таким
образом, доказана оценка
$$
\| T_{\mbox{л}} u \|_{H^{s+2m} (E^{n+1}_{+})} \leq c \, \| u \|_{H^{s+2m+1} (E^{n+1}_{+})},
$$
в которой постоянная не зависит от функции $u.$ Следовательно,
оператор $T_{\mbox{л}}$ "--- сглаживающий.

То, что оператор $\mathfrak{R}$ является правым регуляризатором доказывается вполне аналогично.

Итак, нами доказана

\begin{theorem} \label{teo: 5.2.1}
 Пусть операторы $A$ и $G_j,$ $j=0, \dots, m-1$ образуют эллиптический набор. Пусть выполнены соотношения
 $$
 k - \Re \mu + \frac{1}{2}>0, \  s+2m-k+ \Re \mu -  \frac{1}{2} - \max\limits_{j} m_j > 0.
 $$
Тогда для оператора $\mathfrak{U}$ существует регуляризатор
$\mathfrak{R},$ принадлежащий пространствам  $L \lr{\mathcal{H}^s
(E^{n+1}_{+}, E^n, m), H^{s+2m} (E^{n+1}_{+})}.$ Если функция $u
\in H^{2m+s} (E^{n+1}_{+}),$ а $\mathfrak{U} u \in H^{s+p}$ при
некотором $p>0,$ то $u \in H^{2m+s+p} (E^{n+1}_{+})$ и справедлива
априорная оценка
$$
c \, \| u \|_{H^{s+2m} (E^{n+1}_{+})} \leq \| A u \|_{H^{s} (E^{n+1}_{+})} + \sum\limits_{j=0}^{m-1} \| G_j D_y^k I^{\mu}_e  u|_{y=0}  \|_{H^{s_j} (E^n)} +
 \| u \|_{H^{s+2m-1} (E^{n+1}_{+})},
$$
в которой $c>0$ не зависит от функции $u,$  $s_j = s + 2m-m_j - k
+ \Re \mu - \dfrac{1}{2}.$
\end{theorem}

\section{Общие весовые краевые задачи для сингулярных эллиптических
уравнений}\label{sec13}

\subsection{Весовая краевая задача в полупространстве. Постоянные
коэффициенты}\label{sec13.1}

Рассмотрим в полупространстве $E_{+}^{n+1}$ уравнение
\begin{equation}
    A \lr {\frac{1}{i} D_{x'}, \frac{1}{i^2} B_y}  u \equiv \sum\limits_{|\alpha'| + 2 \alpha_{n+1} = 2 m}  a_{\alpha} \lr{\frac{1}{i} D_{x'}}^{\alpha'} \lr{\frac{1}{i^2} B_y}^{\alpha_{n+1}} u(x) = f(x),
    \label{6.1.1}
\end{equation}
где $B$ "--- оператор Бесселя с комплексным параметром $\nu.$ На
гиперплоскости $y = 0$ рассмотрим два типа краевых условий:
\begin{equation}
\sigma_{\nu} G_j \lr {\frac{1}{i} D_{x'}, \frac{1}{i^2} B_y} B_y^k \left. u \right|_{y=0} = g_j (x'), \  j=0, \dots, m-1,
\label{6.1.2}
\end{equation}
если $\Re \nu > 0$ или $\nu = 0,$
\begin{equation}
y^{2 \nu +1} D_y G_j \lr {\frac{1}{i} D_{x'}, \frac{1}{i^2} B_y} B_y^k \left. u \right|_{y=0} = g_j (x'), \  j=0, \dots, m-1,
\label{6.1.3}
\end{equation}
если $\Re \nu \geq 0.$ Операторы $G_j,$ $j=0, \dots, m-1$ имеют
следующий вид:
$$
G_j \lr {\frac{1}{i} D_{x'}, \frac{1}{i^2} B_y}  = \sum\limits_{|\alpha'| + 2 \alpha_{n+1} = m_j}  g_{\alpha j} \lr{\frac{1}{i} D_{x'}}^{\alpha'} \lr{\frac{1}{i^2} B_y}^{\alpha_{n+1}}.
$$
Указанные краевые задачи порождают операторы $\mathfrak{U}_{\nu}$ и $\mathfrak{U'}_{\nu}$ определяемые по формулам
$$
\mathfrak{U}_{\nu} u = \left\{ Au, \sigma_{\nu} \left. G_0 B^k_y u \right|_{y=0}, \dots,  \sigma_{\nu} \left. G_{m-1} B^k_y u \right|_{y=0}   \right\},
$$
$$
\mathfrak{U'}_{\nu} u = \left\{ Au, y^{2 \nu +1} D_y \sigma_{\nu} \left. G_0 B^k_y u \right|_{y=0}, \dots,  y^{2 \nu +1} D_y  \left. G_{m-1} B^k_y u \right|_{y=0}   \right\}.
$$
Нашей задачей является построение регуляризаторов для
$\mathfrak{U}_{\nu}$ и $\mathfrak{U'}_{\nu}.$ Это построение будет
осуществлено методом операторов преобразования.

Введем пространство
$$
\mathcal{H}^s_{\nu} \lr{E^{n+1}_{+}, E^n, m} = H^s_{\nu}
\lr{E^{n+1}_{+}} \times \prod\limits_{j=0}^{m-1} H^{2m+s-m_j - 2 k
- 1 + \Re \nu} (E^n)
$$
и снабдим его топологией прямого произведения. Предположим, что
выполнены соотношения
\begin{equation}
2 k + 1 - \Re \nu > 0, \  s+2 m -2 k- 1 + \Re \nu -
\max\limits_j m_j > 0. \label{6.1.4}
\end{equation}

Тогда из результатов раздела~\ref{sec3} следует, что операторы
$\mathfrak{U}_{\nu}$ и $\mathfrak{U'}_{\nu}$ непрерывно отображают
пространства $H^{s+2 m}_{\nu} \lr{E^{n+1}_{+}}$ в  $H^{s}_{\nu}
\lr{E^{n+1}_{+}, E^n, m}.$

Ниже мы изучим подробно лишь оператор $\mathfrak{U}_{\nu},$
поскольку оператор $\mathfrak{U'}_{\nu}$ исследуется аналогично.

Наряду с оператором $\mathfrak{U}_{\nu}$ рассмотрим оператор
$\mathfrak{U},$ определяемый по формуле
$$
\mathfrak{U} u = \left\{ A \lr{ \frac{1}{i} D_{x'},  \frac{1}{i^2} D^2_{y}} u, G_0  \lr{ \frac{1}{i} D_{x'},  \frac{1}{i^2} D^2_{y}} D_y^{2 k } \left. I^{\nu -\frac{1}{2}}_e u \right|_{y=0}, \dots,  \right.
$$
$$
\left. G_{m-1}  \lr{ \frac{1}{i} D_{x'},  \frac{1}{i^2} D^2_{y}} D_y^{2 k } \left. I^{\nu -\frac{1}{2}}_e u \right|_{y=0} \right\}.
$$
Этот оператор был изучен нами в предыдущем параграфе, где было показано, что он непрерывно отображает пространство $H^{s+2 m} \lr{E^{n+1}_{+}}$ в пространство
$$
\mathcal{H}^s \lr{E^{n+1}_{+}, E^n, m} =  H^{s} \lr{E^{n+1}_{+}} \times \prod\limits_{j=0}^{m-1} H^{s+2m-m_j - 2 k + \Re \nu -1} \lr{E^{n}},
$$
наделенное топологией прямого произведения.

Введем операторы преобразования $\mathfrak{B}_{\nu, e}$ и
$\mathfrak{G}_{\nu, e},$   связанные с изучаемой краевой задачей.
При этом мы лишь расширим теперь области определения операторов
преобразования и $P_{\nu, e}$ и $S_{\nu, e},$ введённых и
изученных в первом и третьем параграфах. Обозначим через $\Phi$
следующий набор функций:
$$
\Phi=\{f, g_0, \dots, g_{m-1}\}.
$$
Положим тогда
$$
\mathfrak{B}_{\nu, e} \Phi=\{P_{\nu, e} f, c_{\nu} g_0, \dots, c_{\nu} g_{m-1}\},
$$
$$
\mathfrak{G}_{\nu, e} \Phi=\{S_{\nu, e} f, \frac{1}{c_{\nu}} g_0, \dots, \frac{1}{c_{\nu}} g_{m-1}\},
$$
где постоянная $c_{\nu} = 2 \nu$ при $\Re \nu > 0$ и $c_0=1$ при
$\nu=0,$ а операторы  $P_{\nu, e},$ $S_{\nu, e}$ определены в
разделе~\ref{sec1}. Оператор $\mathfrak{B}_{\nu, e}$ изоморфно
отображает пространство $\mathcal{H}^s$ в $\mathcal{H}^s_{\nu},$ а
$\mathfrak{G}_{\nu, e} $ осуществляет    обратное отображение.

Заметим теперь, что операторы $\mathfrak{U}$ и $\mathfrak{U}_{\nu}$ связаны следующими простыми соотношениями:
$$
\mathfrak{U}_{\nu} =  \mathfrak{B}_{\nu, e} \mathfrak{U} S_{\nu, e}, \  \mathfrak{U} =  \mathfrak{G}_{\nu, e} \mathfrak{U}_{\nu} P_{\nu, e}.
$$
Наличие таких соотношений позволяет построить регуляризатор
сингулярного оператора $\mathfrak{U}_{\nu}$ по формуле
$\mathfrak{R}_{\nu} =  P_{\nu, e} \mathfrak{R} \mathfrak{G}_{\nu,
e},$ где $\mathfrak{R}$ регуляризатор оператора $\mathfrak{U},$
построенный в разделе~\ref{sec5} при условии, что операторы  $A
\Big( \dfrac{1}{i} D_{x'},  \dfrac{1}{i^2} D^2_{y}\Big),$ $G_j
\Big( \dfrac{1}{i} D_{x'}, \dfrac{1}{i^2} D^2_{y}\Big)$ образуют
эллиптический набор. Покажем, что оператор $\mathfrak{R}_{\nu}$
действительно является регуляризатором.

Имеем
$$
\mathfrak{U}_{\nu} \mathfrak{R}_{\nu} =  \mathfrak{B}_{\nu, e}
\mathfrak{U} S_{\nu, e}  P_{\nu, e}  \mathfrak{R}
\mathfrak{G}_{\nu, e} =\mathfrak{B}_{\nu, e} \mathfrak{U}
\mathfrak{R} \mathfrak{G}_{\nu, e} = \mathfrak{B}_{\nu, e}  (I +
T_{\mbox{п}}) \mathfrak{G}_{\nu, e} = I + T_{\mbox{п}, \nu},
$$
где $T_{\mbox{п}, \nu } = \mathfrak{B}_{\nu, e} T_{\mbox{п}}
\mathfrak{G}_{\nu, e}.$ Так как $ T_{\mbox{п}} \in
L(\mathcal{H}^s, \mathcal{H}^{s+1}),$ то
 $T_{\mbox{п}, \nu } \in L(\mathcal{H}^s_{\nu}, \mathcal{H}^{s+1}_{\nu}).$

Аналогично получаем
$$
 \mathfrak{R}_{\nu} \mathfrak{U}_{\nu} = I + T_{\mbox{л}, \nu},  \    T_{\mbox{л}, \nu} = P_{\nu, e}  T_{\mbox{л}} S_{\nu, e},
$$
и так как $T_{\mbox{л}} \in L(H^{s+2 m} (E^{n+1}_{+}), H^{s+2m+1}
(E^{n+1}_{+})),$ то $T_{\mbox{л}, \nu} \in  L (H^{s+2 m}_{\nu}
(E^{n+1}_{+}), H^{s+2m+1}_{\nu} (E^{n+1}_{+})).$

Таким образом, доказана

\begin{lemma} \label{teo: 6.1.1}
Пусть операторы $A \Big( \dfrac{1}{i} D_{x'},  \dfrac{1}{i^2}
D^2_{y}\Big),$ $G_j  \Big( \dfrac{1}{i} D_{x'},  \dfrac{1}{i^2}
D^2_{y}\Big)$ образуют эллиптический набор и пусть выполнены
соотношения~\eqref{6.1.4}. Тогда операторы $\mathfrak{U}_{\nu}$ и
$\mathfrak{U}_{\nu}'$  имеют регуляризаторы, непрерывно
отображающие пространства $\mathcal{H}^s_{\nu}$ в $H^{s+2 m}_{\nu}
(E^{n+1}_{+})$  при всех    допустимых~$\nu.$
\end{lemma}

\subsection{ Весовая краевая задача в полупространстве. Маломеняющиеся
коэффициенты}\label{sec13.2}

Рассмотрим в полупространстве $E^{n+1}_{+}$ следующее сингулярное эллиптическое уравнение в частных производных:
\begin{equation}
    A \lr {x, \frac{1}{i} D_{x'}, \frac{1}{i^2} B_y}  u \equiv \sum\limits_{|\alpha'| + 2 \alpha_{n+1} \leq 2 m}  a_{\alpha} (x) \lr{\frac{1}{i} D_{x'}}^{\alpha'} \lr{\frac{1}{i^2} B_y}^{\alpha_{n+1}} u = f(x).
    \label{6.2.1}
\end{equation}

В случае $\Re \nu >0$ или $\nu=0$ присоединим к
уравнению~\eqref{6.2.1} следующие краевые условия:
\begin{multline}
\sigma_{\nu} (y) G_j \lr {\frac{1}{i} D_{x'}, \frac{1}{i^2} B_y}
B_y^k \left. u \right|_{y=0}
\equiv \sigma_{\nu}(y) \hspace{-3ex}\sum\limits_{|\alpha'| + 2 \alpha_{n+1} \leq m_j}\hspace{-3ex}  g_{\alpha j} (x') \lr{\frac{1}{i} D_{x'}}^{\alpha'} \lr{\frac{1}{i^2} B_y}^{\alpha_{n+1}}  B_y^k \left. u \right|_{y=0}  = g_j (x'), \\
 j=0, \dots, m-1.
\label{6.2.2}
\end{multline}
При $\Re \nu \geq 0$ будем рассматривать также краевые условия
вида:
\begin{equation}
y^{2 \nu +1} D_y G_j \lr {x', \frac{1}{i} D_{x'}, \frac{1}{i^2} B_y} B_y^k \left. u \right|_{y=0} = g_j (x'), \  j=0, \dots, m-1,
\label{6.2.3}
\end{equation}

Относительно коэффициентов  $a_{\alpha}$ и $g_{\alpha j}$ будем
предполагать,  что они бесконечно дифференцируемы и кроме того
пусть $\left. D^p_y a_{\alpha}\right|_{y=0} = 0$ при $p= 1, 2,
\dots.$ Разумеется, $\dfrac{1}{i^2}=-1,$ но запись выше
используется для единообразия.

Пусть выполнены следующие условия:
\begin{enumerate}
\item[1)] многочлены $A (0, \xi', \eta^2)$ и $G_j (0, \xi',
\eta^2)$ образуют  эллиптический набор, (т.~е. удовлетворяют
условию Шапиро---Лопатинского);

\item[2)]  при $|\alpha'| + 2 \alpha_{n+1} = 2 m$ и $x \in
\ov{E_{+}^{n+1}}$ выполняются неравенства $|a_{\alpha} (x) -
a_{\alpha} (0)| < \varepsilon,$ $\Big|\Big(\dfrac{1}{y}
D_{y}\Big)^p a_{\alpha} (x) \Big| < \varepsilon,$ $1 \leq p \leq 3
\Re \nu +s+1;$

\item[3)] существует такое $R>0,$ что при $|x| \geq R,$ и при
$|\alpha'| + 2 \alpha_{n+1} = 2 m$  $a_{\alpha} (x) = a_{\alpha}
(0);$

\item[4)] при $|\alpha'| + 2 \alpha_{n+1} < 2 m$ коэффициенты
$a_{\alpha} (x) \in \mathring{C}^{\infty} (\ov{E^{n+1}_{+}});$

\item[5)] при $j = 0, \dots, m-1,$ $|\alpha'| + 2 \alpha_{n+1} =
m_j$ и при всех $x' \in E^n$ $|g_{\alpha_j } (x) - g_{\alpha_j}
(0)| < \varepsilon;$

\item[6)] при $j = 0, \dots, m-1,$ $|\alpha'| + 2 \alpha_{n+1} =
m_j$ и при $|x'| \geq R$ $g_{\alpha_j } (x')=g_{\alpha_j} (0);$

\item[7)] при $j = 0, \dots, m-1,$ $|\alpha'| + 2 \alpha_{n+1} <
m_j$ коэффициенты $g_{\alpha_j } (x') \in \mathring{C}^{\infty}
(E^n).$
\end{enumerate}

Определим операторы $\mathfrak{U}_{\nu}$ и $\mathfrak{U}_{\nu}'$ по формулам
$$
\mathfrak{U}_{\nu} u = \left\{ Au, \sigma_{\nu}  G_0 \left. B^k_y u  \right|_{y=0}, \dots, \sigma_{\nu}  G_{m-1} \left. B^k u  \right|_{y=0}  \right\},
$$
$$
\mathfrak{U}_{\nu}' u = \left\{ Au, y^{2 \nu + 1}   \left. D_y G_0 B^k u  \right|_{y=0}, \dots, y^{2 \nu + 1}   \left. D_y G_{m-1} B^k u  \right|_{y=0}  \right\}.
$$

При выполнении условий 1)--7) операторы $\mathfrak{U}_{\nu}$ и
$\mathfrak{U}_{\nu}'$ будем называть сингулярными эллиптическими
операторами с $(\varepsilon, s)$ "--- маломеняющимися
коэффициентами. Ниже при достаточно малых $\varepsilon$ будут
построены для них регуляризаторы. При этом мы изучим подробно
только оператор $\mathfrak{U}_{\nu},$ поскольку оператор
$\mathfrak{U}_{\nu}'$ изучается аналогично.

Имея в виду применение метода возмущения, разобьем оператор
$\mathfrak{U}_{\nu}$ на сумму трех слагаемых: оператора с
постоянными коэффициентами  $\mathfrak{U}_{\nu, 0},$ оператора с
малой нормой $\mathfrak{U}_{\nu, 1}$ и оператора
$\mathfrak{U}_{\nu, 2},$ содержащего лишь младшие члены. Таким
образом, мы полагаем
$$
\mathfrak{U}_{\nu, 0} u = \left\{ A_0 \lr{0, \frac{1}{i} D_{x'}, \frac{1}{i^2} B_y} u, \sigma_{\nu} (y) G_{0, 0} \lr{0, \frac{1}{i} D_{x'}, \frac{1}{i^2} B_y} \left. B^k_y u  \right|_{y=0}, \dots, \right.
$$
$$
\left. \sigma_{\nu} (y) G_{m-1, 0} \lr{0, \frac{1}{i} D_{x'}, \frac{1}{i^2} B_y} \left. B^k_y u  \right|_{y=0} \right\},
$$
где
$$
 A_0 \lr{0, \frac{1}{i} D_{x'}, \frac{1}{i^2} B_y}  = \sum\limits_{|\alpha'| + 2 \alpha_{n+1} = 2 m} a_{\alpha} (0)   \lr{\frac{1}{i} D_{x'}}^{\alpha'} \lr{\frac{1}{i^2} B_y}^{\alpha_{n+1}},
$$
$$
G_{j, 0} \lr{0, \frac{1}{i} D_{x'}, \frac{1}{i^2} B_y}  = \sum\limits_{|\alpha'| + 2 \alpha_{n+1} =  m_j} g_{\alpha j} (0)   \lr{\frac{1}{i} D_{x'}}^{\alpha'} \lr{\frac{1}{i^2} B_y}^{\alpha_{n+1}},
$$
$$
j = 0, \dots, m-1.
$$
Оператор $\mathfrak{U}_{\nu, 1}$ имеет вид
$$
\mathfrak{U}_{\nu, 1} u = \left\{ A_1 \lr{x, \frac{1}{i} D_{x'}, \frac{1}{i^2} B_y} u, \sigma_{\nu} (y) G_{0, 1} \lr{x', \frac{1}{i} D_{x'}, \frac{1}{i^2} B_y} \left. B^k_y u  \right|_{y=0},  \right.
$$
$$
\left. \dots, \sigma_{\nu} (y) G_{m-1, 1} \lr{x', \frac{1}{i} D_{x'}, \frac{1}{i^2} B_y} \left. B^k_y u  \right|_{y=0}\right\},
$$
где
$$
A_1 \lr{x, \frac{1}{i} D_{x'}, \frac{1}{i^2} B_y}  = \sum\limits_{|\alpha'| + 2 \alpha_{n+1} = 2 m} (a_{\alpha} (x)- a_{\alpha} (0))   \lr{\frac{1}{i} D_{x'}}^{\alpha'} \lr{\frac{1}{i^2} B_y}^{\alpha_{n+1}},
$$
$$
G_{j, 1} \lr{x', \frac{1}{i} D_{x'}, \frac{1}{i^2} B_y}  = \sum\limits_{|\alpha'| + 2 \alpha_{n+1} =  m_j} (g_{\alpha j} (x')-g_{\alpha j} (0))   \lr{\frac{1}{i} D_{x'}}^{\alpha'} \lr{\frac{1}{i^2} B_y}^{\alpha_{n+1}},
$$
$$
j = 0, \dots, m-1.
$$
Оператор $\mathfrak{U}_{\nu, 2}$ определяется по формуле
$$
\mathfrak{U}_{\nu, 2} u = \left\{ A_2 \lr{x, \frac{1}{i} D_{x'}, \frac{1}{i^2} B_y} u, \sigma_{\nu} (y) G_{0, 2} \lr{x', \frac{1}{i} D_{x'}, \frac{1}{i^2} B_y} \left. B^k_y u  \right|_{y=0},  \right.
$$
$$
\dots, \left. \sigma_{\nu} (y) G_{m-1, 2} \lr{x', \frac{1}{i} D_{x'}, \frac{1}{i^2} B_y} \left. B^k_y u  \right|_{y=0} \right\},
$$
где
$$
A_2 \lr{x, \frac{1}{i} D_{x'}, \frac{1}{i^2} B_y}  = \sum\limits_{|\alpha'| + 2 \alpha_{n+1} \leq 2 m-1} a_{\alpha} (x)  \lr{\frac{1}{i} D_{x'}}^{\alpha'} \lr{\frac{1}{i^2} B_y}^{\alpha_{n+1}},
$$
$$
G_{j, 2} \lr{x', \frac{1}{i} D_{x'}, \frac{1}{i^2} B_y}  = \sum\limits_{|\alpha'| + 2 \alpha_{n+1} \leq m_j-1} g_{\alpha j} (x')   \lr{\frac{1}{i} D_{x'}}^{\alpha'} \lr{\frac{1}{i^2} B_y}^{\alpha_{n+1}},
$$
$$
j = 0, \dots, m-1.
$$
причём, в последнем соотношении, предполагается, что $m_j \geq 1.$
Если же $m_j=0$ при некотором номере $j,$ то мы полагаем $G_{j,
2}=0.$

Следовательно,
\begin{equation}
\mathfrak{U}_{\nu} =  \mathfrak{U}_{\nu, 0} + \mathfrak{U}_{\nu, 1} +\mathfrak{U}_{\nu, 2}.
\label{6.2.4}
\end{equation}
Во всём дальнейшем предполагается выполнение
соотношений~\eqref{6.1.4}.

Заметим теперь, что оператор $\mathfrak{U}_{\nu, 0}$ порожден
краевой задачей с однородным оператором с постоянными
коэффициентами как в уравнении, так и в граничных условиях. Для
таких операторов в предыдущем пункте был построен регуляризатор,
который здесь мы обозначим через $\mathfrak{R}_{\nu, 0}.$ По
доказанному оператор   $\mathfrak{R}_{\nu, 0} \in L
\lr{\mathcal{H}^s_{\nu} (E^{n+1}_{+}, E^n, m), H^{s+2m}
(E^{n+1}_{+})}.$ Для него имеют место формулы
\begin{equation}
\mathfrak{U}_{\nu, 0}  \mathfrak{R}_{\nu, 0}= I +   T_{\mbox{п}, 0}, \  \mathfrak{R}_{\nu, 0} \mathfrak{U}_{\nu, 0}  = I +   T_{\mbox{л}, 0},
\label{6.2.5}
\end{equation}
где через $I$ обозначен тождественный оператор, а операторы
$T_{\mbox{п}, 0}$ и $T_{\mbox{л}, 0}$ сглаживающие, то есть
$T_{\mbox{п}, 0} \in L \lr{H^{s+2m}_{\nu}, H^{s+2m+1}_{\nu} },$
$T_{\mbox{л}, 0} \in L \lr{ \mathcal{H}^{s}_{\nu},
\mathcal{H}^{s+1}_{\nu} }.$

Отсюда и из разложения~\eqref{6.2.4} получаем формулы вида:
\begin{equation}
\mathfrak{U}_{\nu}  \mathfrak{R}_{\nu, 0}= I +   T_{\mbox{п}, 0} +  \mathfrak{U}_{\nu, 1} \mathfrak{R}_{\nu, 0}    +\mathfrak{U}_{\nu, 2} \mathfrak{R}_{\nu, 0},
\label{6.2.6}
\end{equation}
\begin{equation}
 \mathfrak{R}_{\nu, 0} \mathfrak{U}_{\nu} = I +   T_{\mbox{л}, 0} + \mathfrak{R}_{\nu, 0} \mathfrak{U}_{\nu, 1} + \mathfrak{R}_{\nu, 0} \mathfrak{U}_{\nu, 2}.
\label{6.2.7}
\end{equation}

Ниже мы покажем, что операторы $I + \mathfrak{U}_{\nu, 1}
\mathfrak{R}_{\nu, 0}$ и $I +  \mathfrak{R}_{\nu, 0}
\mathfrak{U}_{\nu, 1}$ имеют ограниченные обратные в
соответствующих пространствах при достаточно малых $\varepsilon.$
Тогда, вводя обозначение
\begin{equation}
\mathfrak{R}_{\nu, 0} \lr{I + \mathfrak{U}_{\nu, 1} \mathfrak{R}_{\nu, 0}}^{-1} = \mathfrak{R}_{\nu}
\label{6.2.8}
\end{equation}
и замечая, что
\begin{equation}
\mathfrak{R}_{\nu, 0} \lr{I + \mathfrak{U}_{\nu, 1} \mathfrak{R}_{\nu, 0}}^{-1} = \lr{I + \mathfrak{R}_{\nu, 0}  \mathfrak{U}_{\nu, 1}}^{-1}  \mathfrak{R}_{\nu, 0},
\label{6.2.9}
\end{equation}
формулы~\eqref{6.2.6}-\eqref{6.2.7} преобразуем к виду
\begin{equation}
\mathfrak{U}_{\nu}  \mathfrak{R}_{\nu}= I +   \lr{T_{\mbox{п}, 0} +  \mathfrak{U}_{\nu, 2} \mathfrak{R}_{\nu, 0}}   \lr{I + \mathfrak{U}_{\nu, 1} \mathfrak{R}_{\nu, 0}}^{-1} \equiv I +  {T_{\mbox{п}}},
\label{6.2.10}
\end{equation}
\begin{equation}
\mathfrak{U}_{\nu}  \mathfrak{R}_{\nu}= I +  \lr{I + \mathfrak{R}_{\nu, 0} \mathfrak{U}_{\nu, 1} }^{-1}  \lr{T_{\mbox{л}, 0} +   \mathfrak{R}_{\nu, 0} \mathfrak{U}_{\nu, 2}}    \equiv I +  {T_{\mbox{л}}},
    \label{6.2.11}
\end{equation}
причём операторы $T_{\mbox{п}}$ и $T_{\mbox{л}}$  оказываются
сглаживающими.

Для построения обратного к оператору $I + \mathfrak{U}_{\nu, 1}
\mathfrak{R}_{\nu, 0}$ достаточно доказать сходимость в
соответствующей операторной топологии ряда Неймана, стоящего в
правой части следующей формулы:
\begin{equation}
\lr{I + \mathfrak{U}_{\nu, 1} \mathfrak{R}_{\nu, 0} }^{-1}  =
\sum\limits_{l=0}^{\infty} (-1)^l \lr{\mathfrak{U}_{\nu, 1}
\mathfrak{R}_{\nu, 0}}^l.
    \label{6.2.12}
\end{equation}

Пусть $\Phi = \{f, g_0, \dots, g_{m-1} \} \in
\mathcal{H}^s_{\nu}.$ Положим временно $u= \mathfrak{R}_{\nu, 0}
\Phi.$ Тогда  $u \in H^{s+2m}_{\nu} (E_{+}^{n+1})$ и
\begin{equation}
    \left\| \mathfrak{U}_{\nu, 1} \mathfrak{R}_{\nu, 0} \Phi \right\|^2_{\mathcal{H}^s_{\nu}}  =  \left\| \mathfrak{U}_{\nu, 1}  u \right\|^2_{\mathcal{H}^s_{\nu}}
    =  \left\| A_{ 1}  u \right\|^2_{H^s_{\nu} (E_{+}^{n+1})} + \sum\limits_{j=0}^{m-1}  \left\| \sigma_{\nu} G_{j, 1}  B^k_y u|_{y=0} \right\|^2_{H^{s+2m - m_j - 2 k + \Re \nu -1} (E^n)}.
    \label{6.2.13}
\end{equation}

Дадим оценку каждого слагаемого последней формулы, причём
условимся через $c_j,$ $j=1, 2, \dots$ обозначать постоянные,
которые не зависят от $\mathfrak{U}_{\nu, 1}$ и $\Phi.$ Из
следствия~\ref{cor: 3.2.1} и условий 2), 3) непосредственно
получаем
\begin{equation}
    \left\| A_{ 1}  u \right\|^2_{H^s_{\nu} (E_{+}^{n+1})} \leq c_1 \varepsilon \| u \|_{ H^{2m}_{\nu} (E_{+}^{n+1})}.
    \label{6.2.14}
\end{equation}
Для $s> 0$ по формуле Лейбница находим
$$
    \left\| A_{ 1}  u \right\|^2_{H^s_{\nu} (E_{+}^{n+1})} = \sum\limits_{|\alpha'| + \alpha_{n+1} \leq s} \left| D_{x'}^{\alpha'} D^{\alpha_{n+1}}_y S_{\nu, e} A_1 u \right|^2_{L_2 (E_{+}^{n+1})}  \leq
$$
$$
\leq c_2 \lr{\sum\limits_{|\alpha'| + \alpha_{n+1} \leq s} \left|  D^{\alpha_{n+1}}_y S_{\nu, e} \lr{ A_1 D_{x'}^{\alpha'} u} \right|^2_{L_2 (E_{+}^{n+1})} + \sum\limits_{\substack{|\alpha'|>0 \\
|\alpha'| + \alpha_{n+1} \leq s}} \left|  D^{\alpha_{n+1}}_y
S_{\nu, e} \lr{ A_1^{\alpha'} u} \right|^2_{L_2 (E_{+}^{n+1})}},
$$
где
$$
A^{(\alpha')}_1 u = D_{x'}^{\alpha'} A_1 u - A_1  D_{x'}^{\alpha'} =
 \sum\limits_{\substack{|\gamma'| + |\delta'|= |\alpha'|\\
|\delta'|>0}}  \sum\limits_{|\beta'| + 2 \beta_{n+1} = 2 m} c
(\delta', \gamma') \lr{D^{\delta'}_{x'} a_{\beta} (x)}
D^{\beta'+\gamma'}_{x'} B_y^{\beta_{n+1}} u,
$$
$c (\delta', \gamma')$ "--- вполне   определенные постоянные.
Порядок оператора $A^{(\alpha')}_1$ не превосходит  $2m +
|\alpha'| -1,$  поэтому оценивая слагаемые в соответствии с
теоремой~\ref{teo: 3.2.1}, получаем
$$
\left\|  D^{\alpha_{n+1}}_y S_{\nu, e}  A_1 D_{x'}^{\alpha'} u \right\|_{L_2 (E_{+}^{n+1})} \leq
 c_3 \, \| u  \|_{H^{s+2 m}_{\nu} (E_{+}^{n+1})} \max\limits_{l \leq |\alpha_{n+1}| +3N+1} \sup\limits_x \left| \lr{\frac{1}{y} D_y}^l (a_{\alpha} (x)-a_{\alpha} (0 )) \right| \leq
$$
$$
\leq c_3 \,  \varepsilon \, \| u  \|_{H^{s+2 m}_{\nu} (E_{+}^{n+1})}
$$
и
$$
\left\|  D^{\alpha_{n+1}}_y S_{\nu, e}  A_1^{\alpha'} u \right\|^2_{L_2 (E_{+}^{n+1})} \leq
 c_4 \, \| u  \|_{H^{s+2 m-1}_{\nu} (E_{+}^{n+1})} \max\limits_{l \leq 3N+ |\alpha_{n+1}| +1} \sup\limits_x \left| D_{x'}^{\delta'} \lr{\frac{1}{y} D_y}^l a_{\alpha} (x) \right|.
$$
Тогда
\begin{equation}
\left\| A_{ 1}  u \right\|^2_{H^s_{\nu} (E_{+}^{n+1})} \leq c_5 \,
\varepsilon \, \| u \|_{ H^{2m}_{\nu} (E_{+}^{n+1})} + c_5 \, M_s
(a) \| u \|_{ H^{s+2m+1}_{\nu} (E_{+}^{n+1})} + c_6 \, M_s (a) \|
u \|_{ H^{s+2m+1}_{\nu} (E_{+}^{n+1})}, \label{6.2.15}
\end{equation}
где положено
\begin{equation}
M_s (a) = \max\limits_{|\beta| \leq 3N + s +1}
\max\limits_{|\alpha|=2m} \sup\limits_{x \in E^{n+1}_{+}} \left|
D_{x'}^{\beta'} \lr{\frac{1}{y} D_y}^{\beta_{n+1}} a_{\alpha} (x)
\right| \label{6.2.16}
\end{equation}
и натуральное $N$ таково, что  $\Re \nu < N + \dfrac{1}{2}.$

По теореме~\ref{teo: 3.3.2}  о весовых следах существует не
зависящая от операторов  $G_{j, 1}$  и функции $u(x)$ постоянная
$c_7 > 0$ такая, что
\begin{equation}
 \left\| \sigma_{\nu} \left. G_{j, 1} B^k u \right|_{y=0} \right\|_{H^{s+2 m -m_j-2 k + \Re \nu -1} (E^n)} \leq c_7  \left\| G_{j, 1} u \right\|_{H^{s+2 m -m_j}_{\nu} (E^{n+1}_{+})}.
\label{6.2.17}
\end{equation}

Коэффициенты операторов $G_{j, 1}$  не зависят от последней переменной. Тогда из условий 5), 6) немедленно получаем
$$
\| G_{j, 1} u \|_{H^0_{\nu} (E_{+}^{n+1})} \leq c_8 \, \varepsilon \, \| u \|_{H^{m_j}_{\nu} (E_{+}^{n+1})}.
$$
Отсюда при $s > m_j -2 m$ по формуле Лейбница имеем
\begin{equation}
\left\|  G_{j, 1}  u \right\|^2_{H^{s+2m-m_j}_{\nu} (E_{+}^{n+1})}
\leq c_9 \, \varepsilon \, \| u \|_{ H^{s+2m}_{\nu} (E_{+}^{n+1})}
+ c_{10} \, \widetilde{M}_s (\widetilde{g}_j) \| u \|_{
H^{s+2m-1}_{\nu} (E_{+}^{n+1})}, \label{6.2.18}
\end{equation}
где
\begin{equation}
\widetilde{M}_s (\widetilde{g}_j)  = \max\limits_{|\beta'| \leq s+2 m - m_j|} \max\limits_{|\alpha'| +\alpha_{n+1} = m_j} \sup\limits_{x' \in E^n} \left| D^{\beta'}_{x'} g_{\alpha j} (x') \right|.
\label{6.2.19}
\end{equation}
Вспоминая, что $u = \mathfrak{R}_{\nu, 0} \Phi,$ из
оценок~\eqref{6.2.15},~\eqref{6.2.17},~\eqref{6.2.18} получаем
$$
\left\| \mathfrak{U}_{\nu, 1} \mathfrak{R}_{\nu, 0} \Phi \right\|^2_{\mathcal{H}^s_{\nu} (E_{+}^{n+1}, E^n, m)} \leq c_{11} \, \varepsilon \, \left\|  \mathfrak{R}_{\nu, 0} \Phi \right\|_{H^{s+2m}_{\nu} (E_{+}^{n+1})}+
 c_{12}  \, M_s \left\|  \mathfrak{R}_{\nu, 0} \Phi \right\|_{H^{s+2m-1}_{\nu} (E_{+}^{n+1})},
$$
где $M_s = \max \left\{M_s (a), \widetilde{M}_s (\widetilde{g}_1),
\dots, \widetilde{M}_s (\widetilde{g}_{m-1}) \right\}.$ Так как
$\mathfrak{R}_{\nu, 0}$ ограниченный оператор, то последнее
неравенство принимает вид
\begin{equation}
\left\| \mathfrak{U}_{\nu, 1} \mathfrak{R}_{\nu, 0} \Phi \right\|_{\mathcal{H}^s_{\nu}} \leq c_{13} \, \varepsilon \, \left\| \Phi \right\|_{\mathcal{H}^{s}_{\nu}} +  c_{14} \, M_s  \left\| \Phi \right\|_{\mathcal{H}^{s-1}_{\nu}}.
\label{6.2.20}
\end{equation}
Заменяя в этом неравенстве $\Phi$ на $ \mathfrak{U}_{\nu, 1}
\mathfrak{R}_{\nu, 0} \Phi,$ получим
\begin{equation}
\left\| (\mathfrak{U}_{\nu, 1} \mathfrak{R}_{\nu, 0})^2 \Phi  \right\|_{\mathcal{H}^s_{\nu}} \leq c_{13} \, \varepsilon \,  \left\| \mathfrak{U}_{\nu, 1} \mathfrak{R}_{\nu, 0} \Phi \right\|_{\mathcal{H}^{s}_{\nu}} +  c_{14} \, M_s  \left\| \mathfrak{U}_{\nu, 1} \mathfrak{R}_{\nu, 0} \Phi \right\|_{\mathcal{H}^{s-1}_{\nu}}.
\label{6.2.21}
\end{equation}
Оценим последнее слагаемое в правой части. Не ограничивая
общности, будем считать, что $\varepsilon \leq 1$ $M_s \geq 1.$
Используя неравенство Эрлинга---Ниренберга и
неравенства~\eqref{6.2.14},~\eqref{6.2.17}, находим
$$
\left\| \mathfrak{U}_{\nu, 1} \mathfrak{R}_{\nu, 0} \Phi \right\|_{H^{s-1}_{\nu}} \leq \left\| A_{1} \mathfrak{R}_{\nu, 0} \Phi \right\|_{H^{s-1}_{\nu} (E^{n+1}_{+})} + c_{15} \sum\limits_{j=0}^{m-1} \left\| G_{j,1} \mathfrak{R}_{\nu, 0} \Phi \right\|_{H^{s+2 m -m_j-1}_{\nu} (E^{n+1}_{+})} \leq
$$
$$
\leq \varepsilon_1 \left\| A_{1} \mathfrak{R}_{\nu, 0} \Phi \right\|_{H^{s}_{\nu} (E^{n+1}_{+})} + c (\varepsilon_1) \left\| A_{1} \mathfrak{R}_{\nu, 0} \Phi \right\|_{H^{0}_{\nu} (E^{n+1}_{+})} +
$$
$$
 + \varepsilon_1 \sum\limits_{j=0}^{m-1} \left\| G_{j,1} \mathfrak{R}_{\nu, 0} \Phi \right\|_{H^{s+2 m -m_j}_{\nu} (E^{n+1}_{+})} + c (\varepsilon_1)   \sum\limits_{j=0}^{m-1} \left\| G_{j,1} \mathfrak{R}_{\nu, 0} \Phi \right\|_{H^{0}_{\nu} (E^{n+1}_{+})} \leq
 c_{16} \lr{\varepsilon_1 + \varepsilon  c (\varepsilon_1) } M_s \left\|  \Phi \right\|_{H^s_{\nu}},
$$
где $\varepsilon_1$ "---  произвольное положительное число, а $c
(\varepsilon_1)$ зависит от $\varepsilon_1.$ В сочетании с
неравенствами~\eqref{6.2.20} и~\eqref{6.2.21} это приводит нас к
оценке
\begin{equation}
\left\| (\mathfrak{U}_{\nu, 1} \mathfrak{R}_{\nu, 0})^2 \Phi  \right\|_{\mathcal{H}^s_{\nu}} \leq c_{17} M_s^2 \lr{\varepsilon_1 + \varepsilon  c (\varepsilon_1) }  \left\|  \Phi \right\|_{\mathcal{H}^{s}_{\nu}}.
\label{6.2.22}
\end{equation}
Выберем $\varepsilon_1$ так, чтобы $c_{17} M^2_s \varepsilon_1 =
\dfrac{1}{2}$ и предположим, что
$$
c (\varepsilon_1) c_{17} M_s^2  \varepsilon < \frac{1}{2}.
$$
Если $\varepsilon$ удовлетворяет этому условию, то
$$
\left\| (\mathfrak{U}_{\nu, 1} \mathfrak{R}_{\nu, 0})^2 \Phi  \right\|_{\mathcal{H}^s_{\nu}} \leq q  \left\|  \Phi \right\|_{\mathcal{H}^{s}_{\nu}},
$$
где $0<q<1.$ Тогда оператор  $ (\mathfrak{U}_{\nu, 1}
\mathfrak{R}_{\nu, 0})^2$ является сжимающим. А этого достаточно
для сходимости ряда Неймана~\eqref{6.2.12} в банаховом
пространстве линейных ограниченных в $\mathcal{H}^s_{\nu}$
операторов. Следовательно, оператор  $ (I+ \mathfrak{U}_{\nu, 1}
\mathfrak{R}_{\nu, 0})^{-1}$  определён и непрерывен в
пространствах  $\mathcal{H}^s_{\nu}.$ Осталось    только заметить,
что оператор  $\mathfrak{R}_{\nu},$ определенный по формуле
$$
\mathfrak{R}_{\nu} = \mathfrak{R}_{\nu, 0} (I+ \mathfrak{U}_{\nu, 1} \mathfrak{R}_{\nu, 0})^{-1},
$$
в связи с формулами~\eqref{6.2.9}--\eqref{6.2.11} является
регуляризатором для оператора  $\mathfrak{R}_{\nu},$ так как
операторы   $T_{\mbox{п}}$ и $T_{\mbox{л}},$ очевидно,
сглаживающие.

Стало быть, нами доказана
\begin{theorem} \label{teo: 6.2.1}
    Пусть операторы $\mathfrak{U}_{\nu}$ и  $\mathfrak{U'}_{\nu}$ являются сингулярными эллиптическими операторами с  $(\varepsilon, s)$ "--- маломеняющимися коэффициентами и при этом выполнены соотношения
    $$
    2k+1-\Re \nu > 0,
    $$
    $$
    s+2m-2k-2+ \Re \nu - \max\limits_j m_j > 0.
    $$
Тогда при достаточно малых $\varepsilon$ оператор
$\mathfrak{U}_{\nu},$ если  $\Re \nu > 0$ или  $\nu=0,$ и оператор
$\mathfrak{U'}_{\nu},$ если  $\Re \nu \geq 0,$ имеют
регуляризаторы, принадлежащие пространствам  $L
\big(\mathcal{H}^s_{\nu} (E^{n+1}_{+}, E^n, m),$ $H^{s+2m}_{\nu}
(E^{n+1}_{+})\big).$
\end{theorem}

\subsection{Весовая краевая задача в ограниченной
области}\label{sec13.3}

Пусть $\Omega$ "---  ограниченная область с гладкой границей
полупространства $E^{n+1}_{+}.$  Предполагается, что выполнены все
необходимые для неё условия из раздела~\ref{sec4}.

Обозначим через $A$ эллиптический внутри области $\Omega$ оператор
чётного порядка $2 m$  с бесконечно дифференцируемыми в $\Omega$
комплекснозначными коэффициентами. Пусть вблизи границы в каждой
л.с.к. оператор $A$ допускает следующее представление:
\begin{equation}
A= A \lr{x, \frac{1}{i} D_{x'},  \frac{1}{i^2} B_y} = \sum\limits_{|\alpha'|+2 \alpha_{n+1} \leq 2m} a_{\alpha} (x) \lr{\frac{1}{i} D_{x'}}^{\alpha'} \lr{\frac{1}{i^2} B_y}^{\alpha_{n+1}}.
\label{6.3.1}
\end{equation}
Коэффициенты $a_{\alpha} (x)$  предполагаются бесконечно дифференцируемыми вплоть до гиперплоскости  $y=0$  и удовлетворяющими условию
\begin{equation}
D_y^p a_{\alpha} (x) =0, \  y=0, \  p=1,2, \dots.
\label{6.3.2}
\end{equation}

Зададим $m$ граничных операторов $G_j,$ $j=0, 1, \dots, m-1,$ в
каждой л.с.к. по формулам вида
$$
G_j \lr{x', \frac{1}{i} D_{x'},  \frac{1}{i^2} B_y} = \sum\limits_{|\alpha'|+2 \alpha_{n+1} \leq m_j} g_{\alpha j} (x') \lr{\frac{1}{i} D_{x'}}^{\alpha'} \lr{\frac{1}{i^2} B_y}^{\alpha_{n+1}}.
$$
Здесь $g_{\alpha j}$ "--- бесконечно дифференцируемые коэффициенты
переменного $x' \in E^n.$

Предположим, что многочлены
$$
A_{0} (x, \xi', \eta^2) = \sum\limits_{|\alpha'|+2 \alpha_{n+1} \leq 2m} a_{\alpha} (x) {\xi'}^{\alpha'} \eta^{2 \alpha_{n+1}}
$$
и
$$
G_{j, 0} (x, \xi', \eta^2) = \sum\limits_{|\alpha'|+2 \alpha_{n+1} \leq m_j} g_{\alpha j} (x') {\xi'}^{\alpha'} \eta^{2 \alpha_{n+1}},
$$
где $\xi' \in E_n,$ $j=0, \dots, m-1,$ образуют эллиптический
набор для каждой фиксированной точки границы.

Рассмотрим краевую задачу вида
\begin{equation}
\begin{cases}
A u = f, & \\
\left. \sigma_{\nu} \widetilde{G}_j u \right|_{\pr \Omega} = g_j,
& j=0, \dots, m-1,
\end{cases}
\label{6.3.3}
\end{equation}
где операторы $\widetilde{G}_j$ в каждой л.с.к. имеет вид
$$
\widetilde{G}_j =  G_j \lr{x', \frac{1}{i} D_{x'},  \frac{1}{i^2} B_y} B_y^k.
$$
В дальнейшем, как и ранее, предполагается выполнение соотношений
\begin{equation}
2 k +1 - \Re \nu >0, \  s+2 m - 2 k-1 - \Re \nu - \max\limits_j m_j > 0.
\label{6.3.4}
\end{equation}
Основной целью этого пункта и всей  главы является доказательство нётеровости поставленной краевой задачи (в конце мы рассмотрим еще одну краевую задачу). Для этого нам необходимо построить регуляризатор для оператора
$$
\mathfrak{U}_{\nu}: u \to  \mathfrak{U} u = \{Au, \left. \sigma_{\nu} \, \widetilde{G}_0 \right|_{\pr \Omega}, \dots, \left. \sigma_{\nu}\, \widetilde{G}_{m-1} \right|_{\pr \Omega}  \}.
$$

Из результатов раздела~\ref{sec4} следует, что оператор
$\mathfrak{U}_{\nu}$ непрерывно ото6ражает пространство
$H^{s+2m}_{\nu} (\Omega)$ в пространство
$$
\mathcal{H}^{s}_{\nu} (\Omega, \pr \Omega, m) = H^{s}_{\nu}
(\Omega)  \times \prod\limits_{j=0}^{m-1} H^{s+2m-2k -1 + \Re \nu - m_j} (\pr \Omega),
$$
наделенное топологией прямого произведения.

Построение регуляризатора будем осуществлять локально. Для этого нам потребуется некоторое специальное продолжение оператора $\mathfrak{U}_{\nu}$ вида
\begin{equation}
\mathfrak{U}_{\nu} u = \{Au, \left. \sigma_{\nu} \, \widetilde{G}_0 \right|_{y=0}, \dots, \left. \sigma_{\nu} \, \widetilde{G}_{m-1} \right|_{y =0}  \},
\label{6.3.5}
\end{equation}
заданного по указанным выше формулам в л.с.к., скажем, в некоторой
окрестности точки $x^0$ гиперплоскости $y=0.$ Пусть для
определенности оператор~\eqref{6.3.5} (точнее его коэффициенты)
определены в некотором полукубе $K^{+}_{\delta_0} (x^0),$
$\delta_0 >0,$ где
$$
K^{+}_{\delta_0} = \{x = (x', y):x \in E^{n+1}_{+},~|x_p - x^0_p| < \delta,~p=1,\dots,n+1 \}.
$$

Введем по главной части оператора~\eqref{6.3.5}  для произвольного
$\varepsilon \in (0, 1)$ оператор $\mathfrak{U}_{\nu}^{x^0,
\varepsilon}$ вида
\begin{equation}
\mathfrak{U}_{\nu, 0}^{x^0, \varepsilon} u = \{A^{x^0, \varepsilon} u, \left. \sigma_{\nu} \, \widetilde{G}_0^{x^0, \varepsilon} \right|_{y=0}, \dots, \left. \sigma_{\nu} \, \widetilde{G}_{m-1}^{x^0, \varepsilon} \right|_{y =0}  \},
\label{6.3.6}
\end{equation}
где
\begin{multline*}
A^{x^0, \varepsilon} u(x) = \hspace{-3ex}\sum\limits_{|\alpha'|+ 2
\alpha_{n+1} = 2 m }\hspace{-3ex} \varphi_{n+1} \lr{
\frac{x-x^0}{\delta_0}}
 \lr{  \lr{a_{\alpha} (\varepsilon (x- x^0) + x^0) -a_{\alpha} (x^0)} \lr{\frac{1}{i} D_{x'}}^{\alpha'}}  \lr{\frac{1}{i^2} B_{y}}^{\alpha_{n+1}}\hspace{-2ex} u +{}
\\
{}+  \sum\limits_{|\alpha'|+ 2 \alpha_{n+1} = 2 m } a_{\alpha}
(x^0) \lr{\frac{1}{i} D_{x'}}^{\alpha'}  \lr{\frac{1}{i^2}
B_{y}}^{\alpha_{n+1}} u,
\end{multline*}
\begin{multline*}
\widetilde{G}_j^{x^0, \varepsilon} u(x) = \sum\limits_{|\alpha'|+
2 \alpha_{n+1} = m_j } (-1)^k  \varphi_{n} \lr{
\frac{x-x^0}{\delta_0}}  \times
\\
\times \lr{  \lr{g_{\alpha j} (\varepsilon (x- x^0) + x^0)
-g_{\alpha j} (x^0)}} t \lr{\frac{1}{i} D_{x'}}^{\alpha'}
\lr{\frac{1}{i^2} B_{y}}^{\alpha_{n+1} +k} u +{}
\\
{}+
 \sum\limits_{|\alpha'|+ 2 \alpha_{n+1} = m_j } (-1)^k g_{\alpha j} (x^0) \lr{\frac{1}{i} D_{x'}}^{\alpha'}  \lr{\frac{1}{i^2} B_{y}}^{\alpha_{n+1}+k} u,
\quad j=0, \dots, m-1.
\end{multline*}
причём в граничных операторах под $x^0$
понимается точка  $(x_1^0, \dots, x_n^0) \in E^n$ для упрощения обозначений. Через
$\varphi_p,$ $p=n, n+1,$ обозначена функция вида $\varphi_p (x) =
\varphi (x_1) \cdot \dots  \, \cdot \varphi (x_p),$ где    функция
одной переменной $ \varphi (t)\in \mathring{C}^{\infty} (E^1),$
причём $ \varphi (t)=1$ при $t \leq \dfrac{1}{2},$ $ \varphi
(t)=0$ при $t \geq 1.$

Таким образом, оператор $\mathfrak{U}_{\nu, 0}^{x^0, \varepsilon}$
определён уже во всём полупространстве  $E^{n+1}_{+}.$ Отметим
некоторые свойства его бесконечно дифференцируемых коэффициентов.
Из~\eqref{6.3.2} получаем
$$
\left|  D_{x'}^{\beta'} \lr{ \frac{1}{y} D_y}^{\beta_{n+1}} \lr{ \varphi_{n+1} \lr{{ \frac{x-x^0}{\delta_0}}} \lr{a_{\alpha} (\varepsilon (x- x^0) + x^0) -a_{\alpha} (x^0)} } \right| \leq c_{\beta} \, \varepsilon,
$$
$$
\left|  D_{x'}^{\beta'} \lr{ \varphi_{n+1} \lr{{ \frac{x'-x^0}{\delta_0}}} \lr{g_{\alpha j} (\varepsilon (x- x^0) + x^0) -g_{\alpha j} (x^0)} } \right| \leq c_{\beta'} \, \varepsilon,
$$
где постоянные $c_{\beta}, c_{\beta'} > 0$ не зависят от
$\varepsilon \in (0, 1)$ и от $x \in \ov{E_{+}^{n+1}}.$
Следовательно, в соответствии с теоремой~\ref{teo: 6.2.1} оператор
$\mathfrak{U}_{\nu, 0}^{x^0, \varepsilon}$ является сингулярным
оператором с $(\varepsilon, s)$ "--- маломеняющимися
коэффициентами и поэтому при достаточно малом $\varepsilon =
\varepsilon_0,$ для него существует регуляризатор
$\mathfrak{R}_{\nu, 0}^{\varepsilon} \in L \lr{
\mathcal{H}^s_{\nu} (E^{n+1}_{+}, E^n, m), H^{s+2m}_{\nu}
(E^{n+1}_{+})},$ то есть справедливы формулы
\begin{equation}
\begin{array}{l}
\mathfrak{U}_{\nu, 0}^{x^0, \varepsilon} \mathfrak{R}_{\nu}^{x^0, \varepsilon} \Phi = \Phi + T_{\mbox{п}, \nu, 0} \Phi, \  \Phi \in
\mathcal{H}^s_{\nu} (E^{n+1}_{+}, E^n, m), \\
\mathfrak{R}_{\nu}^{x^0, \varepsilon} \mathfrak{U}_{\nu, 0}^{x^0,
\varepsilon} u = u + T_{\mbox{л}, \nu, 0} u, \  u \in
H^{s+2m}_{\nu} (E^{n+1}_{+}),
\end{array}
\label{6.3.7}
\end{equation}
в которых $T_{\mbox{п}, \nu, 0},$  $T_{\mbox{л}, \nu, 0}$
сглаживающие операторы.

Введем теперь оператор растяжения $Q_{\varepsilon}$ по формуле
$$
Q_{\varepsilon}: u \to Q_{\varepsilon} u(x) = u \lr{\frac{x-x^0}{\varepsilon}+x^0}.
$$
Обратным к $Q_{\varepsilon}$ будет оператор
$Q_{\frac{1}{\varepsilon}}.$ Нам потребуется также оператор
$$
Q_{\varepsilon, m}: \Phi = \{f(x), g_0 (x'), \dots, g_{m-1}(x') \} \to Q_{\varepsilon, m} \Phi =
$$
$$
= \left\{ \varepsilon^{-2m } f_0 \lr{\frac{x-x^0}{\varepsilon}+x^0},\ \varepsilon^{-2 \nu - m_0 - 2k } g_0 \lr{\frac{x-x^0}{\varepsilon}+x^0},\
\varepsilon^{-2 \nu - m_{m-1} - 2k } g_{m-1} \lr{\frac{x-x^0}{\varepsilon}+x^0} \right\}.
$$
Обратным к $Q_{\varepsilon, m}$ будет оператор
$Q_{\frac{1}{\varepsilon}, m}.$ Операторы $Q_{\varepsilon},$
$Q_{\varepsilon, m}$ при $\varepsilon>0$ отображают изоморфно
пространства $H_{\nu}^s  (E_{+}^{n+1}),$ $\mathcal{H}_{\nu}^s
(E_{+}^{n+1}, E^n, m),$ соответственно, на себя.

Применим к первой из формул~\eqref{6.3.7} оператор
$Q_{\varepsilon_0, m}$ и заменим в ней $\Phi$ на
$Q_{\frac{1}{\varepsilon}, m} \Phi.$ Тогда после простых
преобразований получим формулу вида
\begin{equation}
\mathfrak{U}_{\nu, 0}^{x^0} \mathfrak{R}_{\nu}^{x^0} \Phi = \Phi + T_{\mbox{п}, \nu, 0}^{x^0} \Phi, \  \Phi \in
\mathcal{H}^s_{\nu} (E^{n+1}_{+}, E^n, m),
\label{6.3.8}
\end{equation}
в которой оператор $T_{\mbox{п}, \nu, 0}^{x^0}  =
Q_{\varepsilon_0, m} T_{\mbox{п}, \nu, 0}^{x^0, \varepsilon_0}
Q_{\frac{1}{\varepsilon_0}, m} $ будет сглаживающим, как и $T_{\mbox{п}, \nu,
0}^{x^0, \varepsilon_0}.$ Оператор
$\mathfrak{R}_{\nu}^{x^0}$ имеет вид $\mathfrak{R}_{\nu}^{x^0}  =
Q_{\varepsilon_0} \mathfrak{R}_{\nu}^{x^0, \varepsilon_0}
Q_{\frac{1}{\varepsilon_0}, m} $ и принадлежит тому же классу $L
\big(\mathcal{H}^s_{\nu} (E^{n+1}_{+}, E^n, m),$ $H^{s+2 m}_{\nu}
(E^{n+1}_{+})\big),$ что и оператор $\mathfrak{R}_{\nu}^{x^0,
\varepsilon_0}.$ Оператор $\mathfrak{U}_{\nu, 0}^{x^0}$ имеет вид
$$
\mathfrak{U}_{\nu, 0}^{x^0} u  = \left\{ \sum\limits_{|\alpha'|+ 2 \alpha_{n+1} = 2 m } \lr{ \varphi_{n+1} \lr{ \frac{x-x^0}{\varepsilon_0 \delta_0}} \lr{a_{\alpha} (x) - a_{\alpha} (x^0)} -a_{\alpha} (x^0)} \right. \times
$$
$$
\times
\lr{\frac{1}{i} D_{x'}}^{\alpha'}  \lr{\frac{1}{i^2} B_{y}}^{\alpha_{n+1}} u,
\sigma_{\nu} (y) \sum\limits_{|\alpha'|+ 2 \alpha_{n+1} = m_0 } (-1)^k  \left( \varphi_{n} \lr{ \frac{x-x^0}{\varepsilon_0 \delta_0}} \times \right.
$$
$$
\times \left.
\lr{g_{\alpha 0} (x) - g_{\alpha 0} (x^0)} -g_{\alpha 0} (x^0)\right)\lr{\frac{1}{i} D_{x'}}^{\alpha'} \left. \lr{\frac{1}{i^2} B_{y}}^{\alpha_{n+1}+k} u \right|_{y=0},
$$
$$
\sigma_{\nu} (y) \sum\limits_{|\alpha'|+ 2 \alpha_{n+1} = m_{m-1} } (-1)^k  \lr{ \varphi_{n} \lr{ \frac{x-x^0}{\varepsilon_0 \delta_0}} \lr{g_{\alpha m-1} (x) - g_{\alpha m-1} (x^0)} -g_{\alpha m-1} (x^0)}
\times
$$
$$
\left.
\times \lr{\frac{1}{i} D_{x'}}^{\alpha'}
\left. \lr{\frac{1}{i^2} B_{y}}^{\alpha_{n+1}+k} u \right|_{y=0} \right\}.
$$

Аналогичным образом трансформируется и вторая из
формул~\eqref{6.3.7}. Отметим, что в некоторой окрестности точки
$x^0,$ а именно в полукубе $K^{+}_{\frac{\varepsilon_0
\delta_0}{2}} (x^0),$ коэффициенты оператора $\mathfrak{U}_{\nu,
0}^{x^0}$ совпадают с коэффициентами главной части оператора
$\mathfrak{U}_{\nu}$ (см.~\eqref{6.3.5}).

Определим теперь оператор $\mathfrak{U}_{\nu}^{x^0}$ по формуле
$\mathfrak{U}_{\nu}^{x^0}= \mathfrak{U}_{\nu, 0}^{x^0} +
\mathfrak{U}_{\nu, 1}^{x^0},$ в которой оператор
$\mathfrak{U}_{\nu, 1}^{x^0}$ имеет вид
$$
\mathfrak{U}_{\nu, 1}^{x^0} u  = \left\{ \sum\limits_{|\alpha'|+ 2 \alpha_{n+1} < 2 m } \varphi_{n+1} \lr{ \frac{x-x^0}{\delta_0}} a_{\alpha} (x)  \lr{\frac{1}{i} D_{x'}}^{\alpha'}  \lr{\frac{1}{i^2} B_{y}}^{\alpha_{n+1}} u, \right.
$$
$$
\sigma_{\nu} (y) \sum\limits_{|\alpha'|+ 2 \alpha_{n+1} < m_0 } (-1)^k  \varphi_{n} \lr{ \frac{x'-x^0}{ \delta_0}} g_{\alpha 0} (x)  \lr{\frac{1}{i} D_{x'}}^{\alpha'} \left. \lr{\frac{1}{i^2} B_{y}}^{\alpha_{n+1}+k} u \right|_{y=0},
$$
$$
 \sigma_{\nu} (y) \sum\limits_{|\alpha'|+ 2 \alpha_{n+1} < m_{m-1} } (-1)^k  \varphi_{n} \lr{ \frac{x'-x^0}{ \delta_0}} g_{\alpha m-1} (x)  \lr{\frac{1}{i} D_{x'}}^{\alpha'}
\left.
\left. \lr{\frac{1}{i^2} B_{y}}^{\alpha_{n+1}+k} u \right|_{y=0}\right\}.
$$
Тогда регуляризатор $\mathfrak{R}_{\nu}^{x^0}$ оператора
$\mathfrak{U}_{\nu}^{x^0}$ будет, очевидно, регуляризатором
оператора $\mathfrak{U}_{\nu}^{x^0}.$

Итак, для оператора $\mathfrak{U}_{\nu}^{x^0}$ существует
регуляризатор  $\mathfrak{R}_{\nu}^{x^0} \in L
\lr{\mathcal{H}^s_{\nu} (E^{n+1}_{+}, E^n, m), H^{s+2 m}_{\nu}
(E^{n+1}_{+})}$ и сглаживающие в соответствующих пространствах
операторы  $T_{\mbox{п}, \nu}^{x^0},$ $T_{\mbox{л}, \nu}^{x^0}$
такие, что
\begin{equation}
\mathfrak{U}_{\nu}^{x^0}  \mathfrak{R}_{\nu}^{x^0} = I + T_{\mbox{п}, \nu}^{x^0},
\label{6.3.9}
\end{equation}
\begin{equation}
\mathfrak{R}_{\nu}^{x^0} \mathfrak{U}_{\nu}^{x^0}   = I +
T_{\mbox{л}, \nu}^{x^0}. \label{6.3.10}
\end{equation}
При этом в некоторое полукубе $K^{+}_{\frac{\varepsilon_0
\delta_0}{2}} (x^0)$ коэффициенты оператора
$\mathfrak{U}_{\nu}^{x^0}$ совпадают с коэффициентами
оператора~\eqref{6.3.5}. Искомое продолжение оператора
$\mathfrak{U}_{\nu}$ из~\eqref{6.3.5} построено.

Такие рассуждения можно провести для любой граничной точки $x$ и
таким образом построить покрытие границы прообразами полукубов
$K^{+}_{\delta} (x),$ $x \in \{y=0\},$ $\delta= \delta(x)>0,$
причём в полукубе $K^{+}_{\delta / 2} (x)$ коэффициенты операторов
$\mathfrak{U}_{\nu}^{x}$ и $\mathfrak{U}_{\nu}$ совпадают (в
л.с.к.). Из указанного покрытия выделим конечное подпокрытие
границы, а вместе с ней и некоторой её окрестности. Обозначим
прообразы центров полукубов, образующих это конечное покрытие,
через $\widetilde{x}^q,$ $q=1, \dots, \ov{q},$ а их образы через
$x^q = \varkappa_{l_q} \widetilde{x}^q.$ Ребро  $q$-го куба
обозначим через $\delta_q.$

Рассмотрим область $\mathsf{y} = \Omega \setminus
\bigcup\limits_{q=1}^{\ov{q}}  \varkappa_{l_q}^1 K^{+}_{\delta_q}
(x^q),$ которая лежит строго внутри области  $\Omega.$ Для области
$\mathsf{y}$ нужно проделать те же построения, что и выше. Нам нет
нужды этим заниматься, поскольку внутри области $\Omega$ оператор
$A$ эллиптичен, причём имеет гладкие коэффициенты. Поэтому мы
можем сослаться, например, на книгу Л.~Хермандера~\cite{83}, в
которой показано существование конечного покрытия области
$\mathsf{y},$ например, кубами $K_{\delta_q} (x^q),$ $q= \ov{q}+1,
\dots, \ov{\ov{q}},$ расположенными  вместе с замыканиями в
$\Omega.$ Причём это покрытие обладает следующим свойством: для
каждого куба  $K_{\delta_q} (x^q),$ $q> \ov{q},$ существует
сохраняющий гладкость оператор продолжения $\mbox{П}^{x_q}$ (его,
впрочем,  можно построить по указанной выше схеме) такой, что для
оператора
$$
A^{x_q} = \sum\limits_{|\alpha| \leq 2 m} \lr{\mbox{П}^{x_q} b_{\alpha} (x) } \lr{\frac{1}{i} D_x}^{\alpha},
$$
где $b_{\alpha}$ "--- коэффициенты оператора $A,$ определенные в
области $\Omega,$ существует регуляризатор $\mathfrak{R}^{x^q} \in
L \lr{H^s (E^{n+1}),  H^{s+2m} (E^{n+1}) }.$ Это означает, что
справедливы формулы
$$
A^{x^q}  \mathfrak{R}^{x^q} = I + T^{x^q}_{\mbox{п}}, \    \mathfrak{R}^{x^q}  A^{x^q} = I + T^{x^q}_{\mbox{л}},
$$
где $ T^{x^q}_{\mbox{п}} \in L \lr{H^s (E^{n+1}), H^{s+1}
(E^{n+1})},$   $ T^{x^q}_{\mbox{л}} \in L \lr{H^{s+2m} (E^{n+1}),
H^{s+2m+1} (E^{n+1})}.$ Весьма важно то, что коэффициенты
оператора $A^{x^q}$ совпадают с коэффициентами оператора $A$ в
кубе $K_{\delta/2} (x^q),$ $q > \ov{q}.$

Итак, мы имеем покрытие области $\ov{\Omega}$ конечным множеством
областей $\varkappa_{l_q}^{-1} K_{\delta} (x^q)$ при этом
диффеоморфизмы $x_{l_q}$  при $q > \ov{q}$ считаются
тождественными отображениями. Не ограничивая общности,
предположим, что области $\varkappa_{l_q}^{-1} K_{\delta_q/4}
(x^q)$  также образуют покрытие $\ov{\Omega}.$ Обозначим через
$\{h_q\}_{q=1}^{\ov{\ov{q}}}$ разбиение единицы, подчиненное этому
последнему покрытию. Введем еще множество финитных бесконечно
дифференцируемых функций $\{\psi_q\}_{q=1}^{\ov{\ov{q}}}$ таких,
что
\begin{equation}
\supp \psi_q \subset \varkappa_{l_q}^{-1} K_{\delta_q / 2} (x^q),
\label{6.3.11}
\end{equation}
\begin{equation}
\psi_q (x) h_q (x) \equiv h_q (x),
\label{6.3.12}
\end{equation}
где $q = 1, \dots, \ov{\ov{q}}.$ Пусть также в каждой л.с.к.
\begin{equation}
D_y \psi_q = 0
\label{6.3.13}
\end{equation}
в некоторой окрестности гиперплоскости $y=0.$ Существование такого
набора функций по существу доказано при выводе леммы~\ref{lem:
4.1.1.}.

Пусть $\Phi = \{f, g_0, \dots, g_{m-1} \} \in \mathcal{H}^s
(\Omega, \pr \Omega, m).$ Тогда искомый регуляризатор
$\mathfrak{R}_{\nu}$ оператора $\mathfrak{U}_{\nu}$ определим по
формуле
$$
\mathfrak{R}_{\nu} \Phi = \sum\limits_{q=1}^{\ov{q}} h_q
\varkappa_{l_q}^{-1}  \mathfrak{R}_{\nu}^{x^q} \varkappa_{l_q}
\lr{\psi_q \Phi} +  \sum\limits_{q=\ov{q}+1}^{\ov{\ov{q}}} h_q
\mathfrak{R}^{x^q} (\psi_q f).
$$
Здесь $\psi_q \Phi = \{\psi_q f, \psi_q|_{\pr \Omega} \, g_0,
\dots, \psi_q|_{\pr \Omega} \, g_{m-1}  \}.$ Используя
теорему~\ref{teo: 3.2.1} и утверждения об эквивалентности норм,
нетрудно усмотреть, что  $ \mathfrak{R}_{\nu} \in L
\lr{\mathcal{H}_{\nu}^{s} (\Omega, \pr \Omega, m), H^{s+2m}
(\Omega)}.$

Покажем, что оператор $\mathfrak{R}_{\nu}$ действительно является
регуляризатором. Пусть функция $u \in H^{s+2m}_{\nu} (\Omega).$
Тогда
\begin{multline}
\mathfrak{R}_{\nu} \mathfrak{U}_{\nu} u = \sum\limits_{q=1}^{\ov{q}} h_q \varkappa_{l_q}^{-1}  \mathfrak{R}_{\nu}^{x^q} \varkappa_{l_q} \lr{\psi_q \mathfrak{U}_{\nu}  u} + \sum\limits_{q= \ov{q}+ 1}^{\ov{\ov{q}}} h_q \mathfrak{R}^{x^q}  \lr{\psi_q A  u} = \\
= \sum\limits_{q=1}^{\ov{q}} h_q \varkappa_{l_q}^{-1}
\mathfrak{R}_{\nu}^{x^q} \varkappa_{l_q} \lr{\psi_q
\mathfrak{U}_{\nu}^{x^q}  u} + \sum\limits_{q= \ov{q}+
1}^{\ov{\ov{q}}} h_q \mathfrak{R}^{x^q}  \lr{\psi_q A^{x^q}  u}.
\label{6.3.14}
\end{multline}
Здесь мы воспользовались тем, что коэффициенты операторов
$A^{x^q},$ $G_j^{x^q}$ совпадают с коэффициентами операторов $A,$
$G_j,$ соответственно, в прообразах кубов  $K_{\delta_q / 2}
(x^q),$ следовательно, и на носителях функций  $\psi_q.$ Далее, по
формуле Лейбница имеем
\begin{equation}
\varkappa_{l_q}  \lr{ \psi_q \mathfrak{U}_{\nu}^{x^q} u } = \varkappa_{l_q}  \lr{ \mathfrak{U}_{\nu}^{x^q} (\psi_q u) + \mathfrak{\hat{U}}_{\nu}^{x^q} u  },~ q = 1, \dots, \ov{q},
\label{6.3.15}
\end{equation}
\begin{equation}
\psi_{q}  A^{x^q} u  = A^{x^q} (\psi_q u) + \hat{A}^{x^q} u,~ q = \ov{q} + 1, \dots, \ov{\ov{q}}.
\label{6.3.16}
\end{equation}

Заметим, что формула Лейбница для операторов Бесселя имеет вид
$$
B (\psi v) = v B \psi + \psi B v + 2 \lr{D_y \psi} D_y v.
$$

Отсюда индукцией по $r$ выводится следующая формула:
$$
B^r (\psi v) =  \psi B^r v + \sum\limits c_{\nu} (r_1, \dots, r_5)
y^{-r_1} \lr{ D_y^{r_2} B^{r_3} \psi } D_y^{r_4} B^{r_5} v,
$$
где суммирование справа осуществляется по всем неотрицательным
целым индексам $r_1, \dots, r_5$ таким, что $r_1+r_2+2 r_3 + r_4
+2 r_5 = 2 r-1,$ $r_2+2 r_3 > 0,$ $r_{3, 4} =0,1.$ Через $c_{\nu}
(r_1, \dots, r_5)$ обозначены вполне определенные постоянные (они
совпадают с биномиальными коэффициентами при $\nu = -
\dfrac{1}{2}$). Сочетание этой формулы и обычной формулы Лейбница
по остальным переменным и приводит нас к
формулам~\eqref{6.3.15},~\eqref{6.3.16} в которых операторы
$\mathfrak{\hat{U}}_{\nu}$ и $\hat{A}$ имеют порядок на единицу
меньший, чем порядки операторов $\mathfrak{U}_{\nu}$ и $A,$
соответственно.

В силу свойств функций  $\psi$ (см.
формулы~\eqref{6.3.11}--\eqref{6.3.13}) по следствиям~\ref{cor:
3.2.1},~\ref{cor: 3.2.2} получаем, что
$$
\varkappa_{l_q}   \mathfrak{\hat{U}}_{\nu}^{x^q} \in L \lr{H^{s+2m}_{\nu} (E^{n+1}_{+}), \mathcal{H}^{s+2m}_{\nu} (E^{n+1}_{+}, E^n, m)},
$$
$$
\hat{A}^{x^q} \in L \lr{H^{s+2m} (E^{n+1}), H^{s+1} (E^{n+1})}.
$$
Так как $\mathfrak{R}_{\nu}^{x^q}$ являются регуляризаторами, то имеет место формула
\begin{equation}
\mathfrak{R}_{\nu} \mathfrak{U}_{\nu} u =
\sum\limits_{q=1}^{\ov{q}} h_q  \psi_q u + T_{\mbox{л}, \nu} u = u
+ T_{\mbox{л}, \nu} u, \label{6.3.17}
\end{equation}
в которой $T_{\mbox{л}, \nu}$ оказывается по доказанному сглаживающим оператором.

Таким образом, доказано, что оператор $\mathfrak{R}_{\nu}$ является левым регуляризатором. Совершенно аналогично показывается, что он будет и правым регуляризатором.

\begin{theorem}\label{teo: 6.3.1}
 Пусть $\Re \nu > 0$ или $\nu=0$ и выполнены соотношения~\eqref{6.3.4}. Тогда краевая задача~\eqref{6.3.3} нётерова. Если $u \in H^{s+2m}_{\nu} (\Omega),$ $ \mathfrak{U}_{\nu} u  \in \mathcal{H}_{\nu}^{s+p} (\Omega, \pr \Omega, m)$ при некотором $p>0,$ то $u \in H^{s+2m+p}_{\nu} (\Omega),$  и имеет место оценка
 \begin{equation}
 c \, \| u \|_{ H^{s+2m+p}_{\nu} (\Omega)} \leq  \|A u \|_{ H^{s+2m+p}_{\nu} (\Omega)}
 + \sum\limits_{j=0}^{m-1} \|\left. \sigma_{\nu} \widetilde{G}_j u\right|_{\pr \Omega} \|_{H^{s+p+2m-m_j-2k-1 + \Re \nu} (\pr \Omega)} + \| u \|_{ H^{s+2m}_{\nu} (\Omega)},
 \label{6.3.18}
 \end{equation}
 где постоянная  $c>0$ не зависит $u.$
\end{theorem}

\begin{proof}
Нётеровость краевой задачи эквивалентна существованию
регуляризатора (см., например, по этому поводу абстрактные
результаты из книги З.~Пресдорфа~\cite{71}). Априорная
оценка~\eqref{6.3.18} следует из того, что для оператора $R_{\nu}$
справедлива формула~\eqref{6.3.17}, в которой $T_{\mbox{л}, \nu} $
"--- сглаживающий оператор. Утверждение о повышении гладкости
также есть следствие этой формулы. В самом деле, если    $u \in
H^{s+2m}_{\nu} (\Omega),$ то  $T_{\mbox{л}, \nu} u \in
H^{s+2m+1}_{\nu} (\Omega).$ Так как $\mathfrak{U}_{\nu} u  \in
\mathcal{H}^{s+p},$ то $\mathfrak{R}_{\nu} \mathfrak{U}_{\nu} u
\in  H^{s+2m+p}_{\nu} (\Omega).$ Тогда    из~\eqref{6.3.16}
следует, что  $u \in H^{s+2m+1}_{\nu} (\Omega).$ Повторяя те же
рассуждения, мы в конце концов покажем, что $u \in
H^{s+2m+p}_{\nu} (\Omega).$ Теорема доказана.
\end{proof}

Рассмотрим теперь сингулярное эллиптическое
уравнение~\eqref{6.3.1} с весовыми краевыми условиями вида
\begin{equation}
\sigma_{\nu + \frac{1}{2}} \widetilde{G}'_j u|_{\pr \Omega} = g_j, \  j = 0, \dots, m-1,
\label{6.3.19}
\end{equation}
где $\Re \nu > 0$ и операторы $\widetilde{G}'_j$ в каждой л.с.к. имеют вид
$$
\widetilde{G}'_j = D_y G_j \lr{x', \frac{1}{i} D_{x'},  \frac{1}{i^2} B_y } B_y^k.
$$
При выполнении всех предшествующих условий справедлива
\begin{theorem}\label{teo: 6.3.2}
Пусть $\Re \nu \geq 0$ и выполнены соотношения~\eqref{6.3.4}.
Тогда краевая задача~\eqref{6.3.1},~\eqref{6.3.18} нётерова. Если
$u \in H^{s+2m}_{\nu} (\Omega),$ а  $ \mathfrak{U}_{\nu}' u  \in
\mathcal{H}_{\nu}^{s+p} (\Omega, \pr \Omega, m),$ $p>0,$ то $u \in
H^{s+2m+p}_{\nu} (\Omega)$ и имеет место оценка
$$
c \, \| u \|_{ H^{s+2m+p}_{\nu} (\Omega)} \leq  \|A u \|_{ H^{s+p}_{\nu} (\Omega)} +
 \sum\limits_{j=0}^{m-1} \|\left. \sigma_{\nu+ \frac{1}{2}} \widetilde{G}'_j u\right|_{\pr \Omega} \|_{H^{s+p+2m-m_j-2k-1 + \Re \nu} (\pr \Omega)} + \| u \|_{ H^{s+2m}_{\nu} (\Omega)},
$$
где постоянная $c>0$ не зависит от $u.$
\end{theorem}

Доказательство этой теоремы полностью аналогично доказательству
теоремы~\ref{teo: 6.3.1} и поэтому не приводится.

\chapter{Новые краевые задачи для уравнения Пуассона с особенностями в изолированных
точках}\label{ch5}

В этой главе рассматриваются новые краевые задачи для уравнения
Пуассона, решения которых могут иметь особенности в изолированных
внутренних точках, которые считаются граничными, причём эти
особенности могут быть не только степенными типа полюса, но и
существенными бесконечного порядка.  Вводятся и изучаются новые
функциональные пространства типа Фреше в ограниченной области с
гладкой границей. Теория этих пространств характеризуется тремя
моментами. Они, во-первых, шире, чем пространства С.\,Л.~Соболева.
Во-вторых, они содержат все гармонические функции, имеющие
произвольные особенности в конечном числе фиксированных внутренних
точек (не ограничивая общность, рассмотрен случай одной такой
точки). И, в-третьих, вне особой точки локально они совпадают с
пространствами С.\,Л.~Соболева. Сочетания этих свойств в одном
пространстве удалось получить благодаря использованию многомерных
операторов преобразования из второй главы. Кроме того, вводится
понятие в определенном смысле нелокального следа в особой точке,
который становится нетривиальным лишь для сингулярных в этой точке
функций. Здесь же доказываются прямая и обратная теоремы о следах.
В терминах указанного следа  дается классификация изолированных
особых точек гармонических функций и доказывается основной
результат главы об однозначной разрешимости соответствующей
краевой задачи для уравнения Пуассона. Ещё раз отметим, что у
решения и правой части уравнения допускается в особой точке рост
произвольного порядка.

В особых точках необходимо использование нового нелокального
краевого условия, которое названо далее сигма "--- следом. Можно
также предложить название $K$-след в честь В.\,В.~Катрахова,
которым это краевое условие было введено и подробно изучено.

\section{Функциональные пространства}\label{sec14}

\subsection{Определение и теоремы вложения}\label{sec14.1}

Пусть $\Omega$ "--- ограниченная область в пространстве $E^n$ с
гладкой границей $\pr \Omega.$ Пусть начало координат $\mathbf{0}$
принадлежит $\Omega.$ Обозначим через $\Omega_0$ область $\Omega
\setminus \mathbf{0}.$ Пусть $R_0>0$ такое число, что шар
$\ov{U}_{2 R_0} \subset \Omega.$

Обозначим через $X$ множество функций $\chi (r) \in
\mathring{C}^{\infty} (\ov{E_{+}^1})$ таких, что $\chi (r) = 1$
при $0 \leq r \leq 1,$  $\chi (r) = 0$ при  $r \geq 2$ и $\chi
(r)$ монотонно убывают при $1<r<2.$ По функции $\chi$ функция
$\chi_R,$ $R> 0,$ определяется по формуле $\chi_R (r) = \chi
\Big(\dfrac{r}{R}\Big).$ Через  $\mathring{T}^{\infty} (\Omega_0)$
обозначим множество функций $f \in C^{\infty} \lr{\ov{\Omega}
\setminus \mathbf{0}}$ таких, что
\begin{equation}
\chi_R (r)  f (r, \vartheta) = \sum\limits_{k=0}^{\mathcal{K}}
\sum\limits_{l=1}^{d_k} \chi_R (r) f_{k, l} (r) Y_{k, l}
(\vartheta), \label{7.1.1}
\end{equation}
где функции $f_{k,l}$ таковы, что $r^{-k} f_{k, l} \chi_R \in
\mathring{C}^{\infty}_{\frac{n}{2} +k -1} (E_{+}^1).$ Натуральное
число $\mathcal{K}$ свое для каждой функции $f.$ Напомним (см.
главу~\ref{ch2}), что множество $\mathring{C}^{\infty}_{\nu}
(E_{+}^1)$ состоит из функций $h(r)$ вида $h = P_{\nu} g,$ где $ g
\in \mathring{C}^{\infty}_{\nu}  (\ov{E_{+}^1})$

Для каждого $s \geq 0,$ $0 <R < R_0$ и каждой функции $\chi \in X$
зададим на $\mathring{T}^{\infty} (\Omega_0)$ нормы
\begin{equation}
 \| f \|_{ s, R} = \left(  \sum\limits_{k=0}^{\mathcal{K}} \sum\limits_{l=1}^{d_k}  \|r^{-k} \chi_R  f_{k, l} \|^2_{\mathring{H}^{s}_{\frac{n}{2} +k -1} (0, 2R)} + \|(1 - \chi_R) f \|^2_{H^s (\Omega)} \right)^{1/2}.
\label{7.1.2}
\end{equation}
Здесь и всюду ниже будут использованы пространства
$\mathring{H}^s_{\nu} (0, R),$ введённые в пункте~\ref{sec4.4}. Из
результатов главы~\ref{ch2} следует, что при чётных  $s$
справедлива формула
\begin{equation}
 \| f \|_{ s, R}^2 = \| \mathfrak{G}_n \Delta^{\frac{s}{2} }  (\chi_R f) \|^2_{L_2 (U_{2R})} + \|  (1-\chi_R) f \|^2_{H^s (\Omega)}.
\label{713}
\end{equation}

Система норм~\eqref{7.1.2} при $R \in (0, R_0)$ и при
фиксированном $s \geq 0$ задает на линеале $\mathring{T}^{\infty}
(\Omega_0)$ топологию.

Введем пространство  $H^s_{loc} (\Omega_0),$ состоящее из функций
$f$ таких, что при любом  $R \in (0, R_0)$ функция $ (1-\chi_R) f
\in H^s (\Omega).$ Наделим это пространство топологией,
определяемой семейством полунорм
$$
p_{s, R} (f) = \|  (1-\chi_R) f \|^2_{H^s (\Omega)}, \  0 <R < R_0.
$$
Легко видеть, что введённая топология превращает $H^s_{loc}
(\Omega_0)$ в полное топологическое векторное пространство.
Очевидно, что $\mathring{T}^{\infty} (\Omega_0) \subset H^s_{loc}
(\Omega_0)$ и для любой функции $f \in \mathring{T}^{\infty}$
имеет место оценка
$$
 \| f \|_{ s, R} \geq p_{s, R} (f).
$$
Таким образом, это вложение будет и топологическим.

Определим пространство $M^s (\Omega_0)$ как замыкание пространства
$\mathring{T}^{\infty} (\Omega_0)$ по
топологии, порожденной системой норм~\eqref{7.1.2}. Следовательно, элементами
пространства  $M^s (\Omega_0)$ являются обычные функции, которые
принадлежат классу $H^s_{loc} (\Omega_0).$

Установим соотношение между пространствами $M^s (\Omega_0)$ и $H^s
(\Omega).$ Обозначим через $\mathring{T}^{\infty}_{+} (\Omega)$
подмножество функций $f$ из $\mathring{T}^{\infty} (\Omega_0),$
для которых функции  $f_{k, l}$ из разложения~\eqref{7.1.1}
удовлетворяют условию $r^{-k} f_{k, l} \chi_R \in
\mathring{C}^{\infty}_{+} (\ov{E_{+}^1}).$

Пусть функция $f \in \mathring{T}^{\infty}_{+} (\Omega).$ Тогда
при любом $R \in (0, R_0)$ $\chi_R f \in \mathring{T}^{\infty}_{+}
(U_{2 R}).$ По лемме~\ref{lem5.2} и из результатов
пункта~\ref{sec4.4} следует, что
$$
 \| \chi_R f \|_{H^s(U_{2, R})} = \sum\limits_{k=0}^{\mathcal{K}} \sum\limits_{l=1}^{d_k}  \|r^{-k} f_{k, l} \chi_R \|^2_{\mathring{H}^{s}_{\frac{n}{2} +k -1, +} (0, 2R)}  \geq
  \sum\limits_{k=0}^{\mathcal{K}} \sum\limits_{l=1}^{d_k}  \|r^{-k} f_{k, l} \chi_R \|^2_{\mathring{H}^{s}_{\frac{n}{2} +k -1} (0, 2R)}.
$$
Следовательно,
$$
\| f \|^2_{s, R} \leq c \lr{ \| \chi_R f \|_{\mathring{H}^s(U_{2R})}^2 + \| (1-\chi_R) f \|_{H^s(\Omega)}^2}.
$$

Выражение справа является при любом $R \in (0, R_0)$ квадратом
одной из эквивалентных норм пространства $H^s(\Omega).$ Отсюда и
из того, что множество $\mathring{T}^{\infty}_{+} (\Omega)$ всюду
плотно в $H^s(\Omega)$ по лемме~\ref{lem:2.1.1} вытекает

\begin{theorem} \label{teo: 7.1.1}
    Для любого $s \geq 0$ имеет место вложение
$$
H^s(\Omega) \subset  M^s (\Omega_0) \subset H^s_{loc} (\Omega_0),
$$
понимаемое в топологическом смысле, причём в случае нечётного $n$
индуцированная левым вложением и собственная топология
пространства $H^s(\Omega)$ равносильны. Для любого  $R \in (0,
R_0)$ топологии подпространства $\widetilde{H}^s (\Omega \subset
\ov{U}_R)$ пространства $H^s(\Omega),$ состоящего из функций $f
\in H^s(\Omega),$ равных нулю в шаре $U_R,$ индуцированная
вложениями   $\widetilde{H}^s (\Omega \setminus \ov{U}_R) \subset
H^s(\Omega)$ и $\widetilde{H}^s (\Omega \setminus \ov{U}_R)
\subset M^s(\Omega_0),$ равносильны.
\end{theorem}
\begin{corollary} \label{cor: 7.1.1}
    Пространство $H^s(\Omega)$ не всюду плотно в $M^s(\Omega_0).$
\end{corollary}

Пусть $\chi$ $\chi'$ "--- две функции из множества $X.$ Покажем,
что порождаемые ими топологии равносильны.

Пусть $f \in \mathring{T}^{\infty} (\Omega_0)$  и $0<R'<R<R_0.$
Тогда
\begin{multline}
\sum\limits_{k=0}^{\mathcal{K}} \sum\limits_{l=1}^{d_k}  \|\chi_R r^{-k}   f_{k, l} \|^2_{\mathring{H}^{s}_{\frac{n}{2} +k -1} (0, 2R)}  \leq \\
\leq 2 \sum\limits_{k=0}^{\infty} \sum\limits_{l=1}^{d_k}
\|\chi'_{R'} r^{-k}   f_{k, l} \|^2_{\mathring{H}^{s}_{\frac{n}{2}
+k -1} (0, 2R)} + 2 \sum\limits_{k=0}^{\infty}
\sum\limits_{l=1}^{d_k} \|(\chi_{R} - \chi'_{R'}) r^{-k}   f_{k,
l} \|^2_{\mathring{H}^{s}_{\frac{n}{2} +k -1} (0, 2R)}.
\label{7.1.3}
\end{multline}
Так как $\chi$ и $\chi'$ принадлежат множеству $X,$ то $\chi_R (r)
= \chi'_{R'} (r) = 1$ при $r< R'.$  Тогда  ${\chi_R (r) -
\chi'_{R'} (r) = 0}$ при  $0 \leq r \leq R',$ а также и при $r \geq
2 R.$ Следовательно, $(\chi_R (r) - \chi'_{R'} (r)) r^{-k} f_{k,
l} \in \mathring{C}_{+}^{\infty} [0, 2R)$ и по
теореме~\ref{teo:2.1.2} вторая сумма справа в~\eqref{7.1.3}
оценивается нормой
$$
\| (\chi_R - \chi'_{R'}) f \|_{\mathring{H}^s (U_{2 R_0})} \leq c \, \| (\chi_R (r) - \chi'_{R'}) f \|_{{H}^s (\Omega)}.
$$
Отсюда получаем
\begin{equation}
\| f \|_{ s, R}^2 \leq c  \left(  \sum\limits_{k=0}^{\infty}
\sum\limits_{l=1}^{d_k}  \|\chi'_{R'} r^{-k}   f_{k, l}
\|^2_{\mathring{H}^{s}_{\frac{n}{2} +k -1} (0, 2R')}  + \|(1 -
\chi_R) f \|^2_{H^s (\Omega)} +  \|(1 - \chi'_{R'}) f \|^2_{H^s
(\Omega)} \right). \label{7.1.4}
\end{equation}

Из монотонности функций $\chi$ и $\chi'$  и из условия  $R' < R$
следует, что $1-  \chi'_{R'} \geq \delta > 0$ на носителе функции
$1 - \chi_R.$ Значит, функция  $(1 - \chi_R) / (1-  \chi'_{R'})
\in C^{\infty} (\ov{\Omega}).$ Тогда имеет место оценка
$$
\|(1 - \chi_R) f \|_{H^s (\Omega)}  \leq c \,  \|(1 - \chi'_{R'}) f \|_{H^s (\Omega)}.
$$
Подставляя это соотношение в~\eqref{7.1.4}, окончательно получаем
$$
\| f \|_{ s, R}^2 \leq c  \left(  \sum\limits_{k=0}^{\infty}
\sum\limits_{l=1}^{d_k}  \|\chi'_{R'} r^{-k}   f_{k, l}
\|^2_{\mathring{H}^{s}_{\frac{n}{2} +k -1} (0, 2R')}  +  \|(1 -
\chi'_{R'}) f \|^2_{H^s (\Omega)} \right) =  c \, \| f \|_{ s,
R'}^2.
$$
Здесь слева стоит норма, порождаемая функцией $\chi,$ справа "---
функцией $\chi'.$ Меняя местами $\chi$ и $\chi',$ получаем
противоположную оценку. Следовательно, нами доказана

\begin{lemma}\label{lem: 7.1.1}
Топологии, порождаемые в пространстве $M^s (\Omega_0)$ различными
функциями из множества $X,$ равносильны.
\end{lemma}
\begin{lemma}\label{lem: 7.1.2}
    Пространство $M^s (\Omega_0)$ суть полное счётно-нормируемое топологическое пространство, то есть пространство Фреше.
\end{lemma}

\begin{proof}
Рассмотрим счётное множество норм $\| \|_{s, R_0/m},$ $m=1, 2
\dots.$ Пусть произвольное число $R \in (0, R_0).$ Тогда найдутся
такие   натуральные числа $m_1$ и $m_2,$ что $\dfrac{R_0}{m_1}<R<
\dfrac{R_0}{m_2}.$ В таком случае по лемме~\ref{lem: 7.1.1}
найдутся постоянные $c_1, c_2 > 0$ такие, что для любой функции $f
\in M^s (\Omega_0)$ справедлива оценка
$$
c_2 \|f \|_{s, R_0/ m_2} \leq \| f\|_{s, R} \leq  c_1 \|f \|_{s,
R_0/ m},
$$
которая и доказывает наше утверждение.
\end{proof}

\begin{theorem}\label{teo: 7.1.2}
    Пусть $0 \leq s_1 < s_2.$ Тогда пространство $M^{s_1} (\Omega_0)$ непрерывно вложено в $M^{s_2} (\Omega_0).$
\end{theorem}

\begin{proof}
Поскольку для любой функции $g \in \mathring{C}^{\infty} [0, R)$
$$
\| g \|_{\mathring{H}^{s_1} (0, R)} = \| D^{s_1} g\|_{L_2 (0, R)}
\leq c \, (s_1 - s_2, R) \| D^{s_2} g\|_{L_2 (0, R)}=
  c\, (s_1 - s_2, R) \|  g\|_{\mathring{H}^{s_2} (0, R)},
$$
то
$$
\|\chi_{R} \, r^{-k}   f_{k, l}
\|^2_{\mathring{H}^{s_1}_{\frac{n}{2} +k -1} (0, 2R)} =
\frac{\sqrt{\pi}}{2^{\frac{n}{2} +k -1} \, \Gamma \lr{\frac{n}{2}
+k}}  \| S_{\frac{n}{2} +k -1} (\chi_{R} \, r^{-k}   f_{k, l})
\|_{\mathring{H}^{s_1} (0, R)} \leq
$$
$$
\leq c \, (s_1 - s_2, R)
\|\chi_{R} \, r^{-k}   f_{k, l}
\|^2_{\mathring{H}^{s_2}_{\frac{n}{2} +k -1} (0, 2R)}.
$$
Суммируя эти неравенства по $k$ и $l$ и замечая, что
$$
\| (1-\chi_R) f \|_{H^{s_1} (\Omega)} \leq  \| (1-\chi_R) f
\|_{H^{s_2} (\Omega)},
$$
получаем неравенство $ \|f\|_{s_1, R} \leq c \, \|f\|_{s_2, R},$ в
котором постоянная не зависит от функции  $f \in M^{s_2}
(\Omega_0).$ Теорема    доказана.
\end{proof}

Перейдем к изучению следов функций из $M^s (\Omega_0)$ на границе
области $\Omega_0,$ которая состоит из точки $\mathbf{0}$ и
поверхности $\pr \Omega.$

Отметим, что так как пространство $M^s$ вне любой окрестности
точки $\mathbf{0}$ устроено так же как и пространство $H^s,$ то
справедлива
\begin{theorem}\label{teo: 7.1.3}
Пусть $s > \dfrac{1}{2}.$ Тогда отображение  $f \to f|_{\pr
\Omega}$  непрерывно действует из пространства  $M^s (\Omega_0)$ в
$H^{s- \frac{1}{2}} (\pr \Omega).$ Существует оператор, непрерывно
отображающий $H^{s- \frac{1}{2}} (\pr \Omega)$ в $M^s (\Omega_0),$
такой, что $f|_{\pr \Omega}=g,$ где $f$ "--- образ функции  $g$
при этом отображении.
\end{theorem}

\subsection{Прямая и обратная теоремы о $\sigma$-следах
($K$-следах)}\label{sec14.2}

Дадим определение следа в граничной точке $\mathbf{0} \in \pr
\Omega_0.$ Легко увидеть, что каким бы большим не было $s,$
обычного следа в точке~$\mathbf{0}$ у функции $f \in M^s
(\Omega_0),$ вообще говоря, не существует. Поэтому мы вводим
некоторое новое понятие следа, который носит в определенном смысле
нелокальный характер.

Обозначим через $A(\Theta)$  множество определенных на сфере
$\Theta$ функций $\psi\lr{\vartheta},$ которые
вещественно-аналитичны на $\Theta$ и для которых при каждом $h \in
(0, 1)$ конечны нормы
\begin{equation}
\| \psi \|_h = \lr{  \sum\limits_{k=0}^{\infty}
\sum\limits_{l=1}^{d_k} \|\psi_{k, l} \|^2 h^{-2 k}
}^{\frac{1}{2}}, \label{7.2.1}
\end{equation}
где через  $\psi_{k, l}$ обозначены коэффициенты разложения
функции $\psi$ по сферическим гармоникам $Y_{k, l}.$

\begin{lemma}\label{lem: 7.2.1}
    Пространство $A(\Theta)$ является полным счётно-нормируемым топологическим пространством.
\end{lemma}

\begin{proof}
Докажем полноту. Пусть $\{ \psi^m \} \subset A(\Theta)$ "---
фундаментальная последовательность. Тогда для любого $h \in (0,
1)$ $\| \psi^{m_1} - \psi^{m_2} \|_h \to 0$ при $m_1, m_2 \to
\infty.$ Поэтому для каждого $h$ найдется такая функция $g_h
(\vartheta),$ что $\| g_h \|_h < \infty$ и $\| \psi^m - g_h \|_h
\to 0$ при $m \to \infty.$ Из~\cite[с.~496]{77}, следует, что
функция $g_h$ аналитична на $\Theta.$ Так как  $\| \psi \|_h \leq
\| \psi \|_{h'},$ если $h' < h,$ то из условия $\|g_{h'}- \psi^m
\|_{h'}  \to 0$ следует, что $\|g_{h'}- \psi^m \|_{h}  \to 0$ при
$m \to \infty.$ Следовательно, все функции $g_h$ совпадают друг с
другом. Эта единственная функция принадлежит, очевидно,
пространству $A(\Theta)$ и является пределом последовательности
$\{ \psi^m \}$ по любой норме~\eqref{7.2.1}. Полнота доказана.

Рассмотрим счётное множество норм  $\| \|_{1/m},$ $m=2,3, \dots.$
Легко видеть, что определяемая ими топология равносильна исходной.
Таким образом, пространство  $A(\Theta)$ счётно-нормируемо. Лемма
доказана.
\end{proof}

Отметим, что пространство $A(\Theta)$ состоит из тех
вещественно-аналити\-ческих функций на сфере  $\Theta,$ которые
допускают гармоническое продолжение на всё пространство $E^n.$ По
заданной функции  $\psi \in A(\Theta)$ такое продолжение
определяется по формуле
$$
\Psi (r, \vartheta) = \sum\limits_{k=0}^{\infty}
\sum\limits_{l=1}^{d_k} r^k \psi_{k, l} Y_{k, l} (\vartheta).
$$
Ясно, что функция $\Psi$ гармоническая во всём пространстве $E^n.$

Пусть функция $f \in \mathring{T}^{\infty} (\Omega_0).$ Это
означает, что при $r < 2 R_0$ имеет место формула
\begin{equation}
f (r, \vartheta) = \sum\limits_{k=0}^{\mathcal{K}}
\sum\limits_{l=1}^{d_k} f_{k, l} (r) Y_{k, l} (\vartheta),
\label{7.2.2}
\end{equation}
причём функция $r^{-k} \chi_R f_{k, l} \in
\mathring{C}^{\infty}_{\frac{n}{2} +k -1} (0, 2R),$ $0<R<R_0.$
Заметим, что сами функции $f_{k, l} (r)$ могут иметь особенность в
точке  $r=0.$ При $r < R_0$ и $n \geq 3$ определим операцию
усреднения $\sigma$ на классе функций $\mathring{T}^{\infty}
(\Omega_0)$ по формуле
\begin{equation}
\sigma f (r, \vartheta) = \sum\limits_{k=0}^{\mathcal{K}}
\sum\limits_{l=1}^{d_k} r^{n+k-2} f_{k, l} (r) Y_{k, l}
(\vartheta). \label{7.2.3}
\end{equation}

Назовем {\it $\sigma$-следом} функции  $f \in
\mathring{T}^{\infty} (\Omega_0)$ предел
$$
\lim\limits_{r \to + 0} \sigma f (r, \vartheta) = \sigma f|_0.
$$
Из условий на функцию  $f$ по теореме~\ref{teo:1.4.1} следует, что
её $\sigma$-след существует и равен функции
$$
\psi (\vartheta) = \sum\limits_{k=0}^{\mathcal{K}}
\sum\limits_{l=1}^{d_k} \psi_{k, l} (r) Y_{k, l} (\vartheta),
$$
где
$$
\psi_{k, l} = \lim\limits_{r \to + 0}  r^{n+k-2} f_{k, l} (r).
$$
Ясно также, что $\psi \in A(\Theta).$

Дадим теперь другое, но эквивалентное прежнему, определение
$\sigma$-следа, не использующее разложения по сферическим
гармоникам. Пусть $n \geq 3,$ функции $f_{k, l}$ из~\eqref{7.2.2}
определяются по формуле
$$
f_{k, l} (r) = \int\limits_{\Theta} f (r, \vartheta) Y_{k, l} (\vartheta) \,  d \vartheta.
$$
Подставляя это выражение в~\eqref{7.2.3}, находим
$$
\sigma f (r, \vartheta) = \sum\limits_{k=0}^{\infty}
\sum\limits_{l=1}^{d_k} r^{n+k-2} Y_{k, l} (\vartheta)
\int\limits_{\Theta} f (r, \vartheta') Y_{k, l} (\vartheta') \,  d
\vartheta' =
 \int\limits_{\Theta} f (r, \vartheta')  \sum\limits_{k=0}^{\infty} \sum\limits_{l=1}^{d_k} r^{n+k-2} Y_{k, l} (\vartheta)  Y_{k, l} (\vartheta') \,  d \vartheta'.
$$
Известно (см.например,~\cite{BE2}), что при $r < 1$ имеет место
формула
$$
\sum\limits_{k=0}^{\infty} \sum\limits_{l=1}^{d_k} r^{k} Y_{k, l}
(\vartheta)  Y_{k, l} (\vartheta') = \frac{\Gamma
\lr{\frac{n}{2}}}{2 \pi^{\frac{n}{2}}} \frac{1-r^2}{|r \vartheta -
\vartheta' |^n}=K_n (r \vartheta, \vartheta'),
$$
где функция $K_n (x, y)$ называется ядром Пуассона для единичной
сферы $\Theta.$ Следовательно,
$$
\sigma f|_0 =  \lim\limits_{r \to + 0}  r^{n-2}  \int\limits_{\Theta} f (r, \vartheta') K_n (r \vartheta, \vartheta') \, d \vartheta'.
$$
Это и есть искомое явное определение $\sigma$-следа.

Для размерности $n=2$  $\sigma$-след определяется несколько иначе.
Пусть $r>0,$ $|\varphi|<\pi$ "--- полярные координаты на
плоскости. Вместо~\eqref{7.2.3} в этом случае мы полагаем
$$
\sigma f (r, \varphi) = \frac{f_0 (r)}{\ln r} +
\sum\limits_{k=1}^{\mathcal{K}} r^k \lr{f_{k, 1} (r) \cos (k
\varphi)+ f_{k, 2} (r) \sin (k \varphi)},
$$
где
\begin{align*}
f_0 (r) &= \frac{1}{2 \pi} \int\limits_{- \pi}^{\pi}  f (r, \varphi)  \, d \varphi,
\\
f_{k,1} (r) &= \frac{1}{\pi} \int\limits_{- \pi}^{\pi}  f (r,
\varphi) \cos(k \varphi) \, d  \varphi,
\\
f_{k,2} (r) &= \frac{1}{2 \pi} \int\limits_{- \pi}^{\pi}  f (r,
\varphi) \sin (k \varphi)  \, d \varphi.
\end{align*}
Отсюда имеем
$$
\sigma f (r, \varphi) = \frac{f_0 (r)}{\ln r} +  \frac{1}{\pi}
\int\limits_{- \pi}^{\pi} f (r, \varphi')
\sum\limits_{k=1}^{\infty} r^k  \cos (k( \varphi-\varphi')) \, d
\varphi'.
$$
Так как при $r<1$ имеет место формула
$$
\sum\limits_{k=1}^{\infty} r^k  \cos (k( \varphi-\varphi')) =
\frac{r \cos( \varphi-\varphi') -r^2 }{1-2 r \cos(
\varphi-\varphi') +r^2},
$$
то, стало быть, $\sigma$-след при  $n=2$ можно определить и по
формуле
$$
\sigma f|_0 =  \lim\limits_{r \to + 0}  \frac{1}{2 \pi} \int\limits_{- \pi}^{\pi} f (r, \varphi') \lr{ \frac{2 r \cos( \varphi-\varphi') - 2 r^2 }{1-2 r \cos( \varphi-\varphi') +r^2} + \frac{1}{\ln r}}  \, d \varphi'.
$$

Распространим понятие $\sigma$-следа на функции $f \in M^s
(\Omega_0).$ Для любой функции $f \in M^s (\Omega_0)$ существует
последовательность функций $f^j \in \mathring{T}^{\infty}
(\Omega_0),$ $j = 1, 2, \dots,$ сходящихся к $f$ при $j \to
\infty$ по топологии пространства $M^s (\Omega_0).$ Для каждой
функции $f^j$ $\sigma$-след, по доказанному, существует и
принадлежит пространству $A(\Theta).$ Если последовательность
$\psi^j = \sigma f^j|_0$ сходится при $j \to \infty$ к функции
$\psi \in A(\Theta)$ по топологии пространства $A(\Theta),$ причём
функция $\psi$ не зависит от выбора последовательности $f^j,$ то
она и называется $\sigma$-следом функции $f.$

Имеет место следующая прямая теорема о $\sigma$-следах.

\begin{theorem}\label{teo: 7.2.1}
    Пусть $s \geq 1,$ если  $n \geq 3,$ и $s \geq 2,$ если  $n=2.$ Тогда для каждой функции $f \in M^s (\Omega_0)$ существует $\sigma$-след $\sigma f|_0 \in A(\Theta).$  При этом оператор $f \to \sigma f|_0$ непрерывно отображает пространство  $ M^s (\Omega_0)$ в $A(\Theta).$
\end{theorem}

\begin{proof}
Достаточно показать, что указанный оператор непрерывно отображает
пространство $\mathring{T}^{\infty} (\Omega_0)$ с индуцированной
пространством $M^s (\Omega_0)$ топологией в пространство
$A(\Theta).$ Для этого  необходимо доказать, что для любого $h \in
(0, 1)$ существует такое число $R \in (0, R_0)$ и такая постоянная
$c>0,$ что для любой функции $f \in \mathring{T}^{\infty}
(\Omega_0)$ справедлива оценка
\begin{equation}
\| \sigma f|_0 \|_h \leq c \, \| f \|_{s, R}. \label{7.2.4}
\end{equation}
Пусть $f_{k, l}$ "--- коэффициенты разложения функций $f \in
\mathring{T}^{\infty} (\Omega_0)$ по сферическим гармоникам $Y_{k,
l}.$ Тогда функции $r^{-k} \chi_R f_{k, l}$ принадлежат
пространству $\mathring{C}^{\infty}_{\frac{n}{2}+k-1} (E^1_{+}),$
а значит, и пространству  $\mathring{H}^{s}_{\frac{n}{2}+k-1} (0,
2R).$

По следствию~\ref{cor:1.4.3} для любой функции $g \in
\mathring{H}^{s}_{\nu} (0, 2R)$ справедлива оценка
$$
|\left. \sigma_{\nu} (r) g(r)\right|_{r=0} | \leq c \, (s, R)
(4R)^{\nu} (\nu+1)^{-s} \| g\|_{\mathring{H}^{s}_{\nu} (0, 2R)},
$$
если $s \geq 1,$ $\nu \geq 0,$ $s+\nu > 1.$ Здесь $\sigma_{\nu}
(r) = r^{2 \nu}$ при $\nu > 0$ и $\sigma_0 (r) = \Big(\ln
\dfrac{1}{r}\Big)^{-1}.$ Полагая в этом неравенстве $\nu =
\dfrac{n}{2}+k-1,$ $ g = \chi_R r^{-k} f_{k, l},$ получим
$$
|\left. \sigma_{\frac{n}{2}+k-1} (r)  r^{-k} f_{k, l}
(r)\right|_{r=0} | \leq c \, (s, n,  R) (4R)^{k} (k+1)^{-s}  \| \chi_R \, r^{-k} f_{k,
l}\|_{\mathring{H}^{s}_{\frac{n}{2}+k-1} (0, 2R)}.
$$
Далее, так как
$$
\sigma f |_0 = \sum\limits_k \sum\limits_l
\sigma_{\frac{n}{2}+k-1} r^{-k} f_{k, l} (r)|_{r=0} Y_{k, l}
(\vartheta),
$$
то для любого $h \in (0, 1)$ из предыдущей оценки получаем
$$
\| \sigma f |_0 \|_h^2 = \sum\limits_k \sum\limits_l |
\sigma_{\frac{n}{2}+k-1} r^{-k} f_{k, l} (r)|_{r=0} |^2 h^{-2 k}
\leq
$$
$$
\leq c \, (s, R, n) \sum\limits_k \sum\limits_l (4 R)^{2 k}
(k+1)^{-2 s} h^{-2 k}  \|  \chi_R \, r^{-k} f_{k, l}
\|^2_{\mathring{H}^{s}_{\frac{n}{2}+k-1} (0, 2R)} \leq
$$
$$
\leq c \, (s, R, n) \sum\limits_k \sum\limits_l  \|  \chi_R \,
r^{-k} f_{k, l} \|^2_{\mathring{H}^{s}_{\frac{n}{2}+k-1} (0, 2R)},
$$
где положено $R = \dfrac{1}{4} h.$ Теорема доказана.
\end{proof}

Рассмотрим обратную теорему о $\sigma$-следах.

\begin{theorem}\label{teo: 7.2.2}
    Для любой функции $\psi \in A (\Theta)$ существует функция $f,$ принадлежащая пространствам $M^s (\Omega_0)$ при всех $s \geq 0,$ для которой  $\psi$  является её $\sigma$-следом. При этом для любого $R \in (0, R_0)$ существует число $h \in (0, 1)$ такое,    что имеет место неравенство
\begin{equation}
\| f \|_{s, R} \leq c \, \| \psi \|_h
\label{7.2.5}
\end{equation}
с постоянной $c>0,$ не зависимой от выбора функции $\psi.$
\end{theorem}

\begin{proof}
Пусть функция $\psi \in A (\Theta)$ и $\psi_{k, l}$ "--- её
коэффициенты разложения в ряд по сферическим гармоникам. Функцию
$f,$ о которой говорится в теореме, определим по формуле
\begin{equation}
f (r, \vartheta) = \sum\limits_{k=0}^{\infty}
\sum\limits_{l=1}^{d_k} r^{2-n - k}  \psi_{k, l}  Y_{k, l}
(\vartheta) \label{7.2.6}
\end{equation}
в случае, если  $n \geq 3.$ Если же $n=2,$ то положим
\begin{equation}
f (r, \vartheta) = \psi_0 \ln r + \sum\limits_{k=1}^{\infty}
\sum\limits_{l=1}^{2} r^{ - k}  \psi_{k, l}  Y_{k, l} (\vartheta).
\label{7.2.7}
\end{equation}

Отметим, что так как $\psi \in A (\Theta),$ то ряды~\eqref{7.2.6}
и~\eqref{7.2.7} сходятся абсолютно и равномерно вне любого шара с
центром в  $\mathbf{0}$ к аналитической функции. По той же причине
функция $f$ будет даже гармонической во всём пространстве $E^n$ за
исключением точки $\mathbf{0}.$

Покажем справедливость оценки~\eqref{7.2.5}. Пусть сначала $s \geq
2$ чётное число. Тогда, используя свойства операторов
преобразования (см. главу~\ref{ch2}), получаем
\begin{multline}
\| \chi_R \, r^{2-n-2k} \|^2_{\mathring{H}^{s}_{\frac{n}{2}+k-1}
(0, 2R)} \leq \frac{2 \pi}{ 2^{n+2k-1}\, \Gamma^2
\lr{\frac{n}{2}+k} }
\| D^s S_{\frac{n}{2}+k-1} (\chi_R \, r^{2-n-2k}) \|^2_{L_2 (0, 2R)} = \\
= \frac{2 \pi}{ 2^{n+2k-1}\, \Gamma^2 \lr{\frac{n}{2}+k} } \|
S_{\frac{n}{2}+k-1} B^{\frac{s}{2}}_{\frac{n}{2}+k-1} (\chi_R
r^{2-n-2k}) \|^2_{L_2 (0, 2R)}, \label{7.2.8}
\end{multline}
где $S_{\nu}$ "--- операторы преобразования из
пункта~\ref{sec4.1}.

Применяя соотношение Дарбу---Вайнштейна  в~\eqref{7.2.8} получаем
\begin{equation}
\| \chi_R \, r^{2-n-2k} \|^2_{\mathring{H}^{s}_{\frac{n}{2}+k-1}
(0, 2R)} \leq \frac{2 \pi}{ 2^{n+2k-1}\, \Gamma^2
\lr{\frac{n}{2}+k} }
 \| S_{\frac{n}{2}+k-1} (r^{2-n-2k}
B^{\frac{s}{2}}_{\frac{n}{2}+k-1} \chi_R) \|^2_{L_2 (0, 2R)}.
\label{7.2.9}
\end{equation}
Так как функция $\chi_R (r) =1$ при $0 \leq r \leq R,$ то
$$
B_{- \nu}  \chi_R (r) = r^{2 \nu -1} \frac{\pr }{\pr r} \lr{r^{1 -
2 \nu}  \frac{\pr \chi_R}{\pr r}} =0
$$
при тех же значениях $r.$ Значит и $B_{- \nu}^{\frac{s}{2}} \chi_R
(r) = 0$ при $r \leq R,$ а также и при $r \geq 2R,$ поскольку
$\chi_R (r) =0$ при  $r \geq 2R.$ Отсюда следует, что функция
$B^{\frac{s}{2}}_{\frac{n}{2}+k-1} \chi_R$ бесконечно
дифференцируема, имеет компактный носитель и обращается в нуль
вблизи начала. Следовательно, $B^{\frac{s}{2}}_{\frac{n}{2}+k-1}
\chi_R \in \mathring{C}_{+}^{\infty} (\ov{E_{+}^1}).$ Тогда   на
основании леммы~\ref{lem:1.4.1} из формулы~\eqref{7.2.9} получаем
\begin{equation}
\| r^{2-n-2k}  \chi_R \|^2_{\mathring{H}^{s}_{\frac{n}{2}+k-1} (0,
2R)} \leq 2 \int\limits_{R}^{2 R}
|B^{\frac{s}{2}}_{1-\frac{n}{2}-k} \chi_R (r) |^2  r^{3-n-2k} \,
dr. \label{7.2.10}
\end{equation}
Так как $\chi (t) \in C^{\infty} [1, 2],$ то имеет место оценка
$$
| B^{\frac{s}{2}}_{-\nu} \chi (t) | = \left| \lr{D^2_t + \frac{1-2
\nu}{t} D_t}^{\frac{s}{2}} \chi (t) \right| \leq c \, (s) (1+
\nu)^{\frac{s}{2}}
$$
с постоянной $c \, (s)>0,$ не зависящей от $\nu.$ Отсюда и
из~\eqref{7.2.10} следует неравенство
\begin{equation}
\| r^{2-n-2k}  \chi_R \|^2_{\mathring{H}^{s}_{\frac{n}{2}+k-1} (0,
2R)} \leq 2 R^{4 - s- n-2 k}
  \int\limits_{1}^{2} |B^{\frac{s}{2}}_{1-\frac{n}{2}-k}
\chi (t) |^2 \,  t^{3-n-2k} \, d t
\leq c (s, n, R) \, R^{- 2 k} (1+k)^{s-1}. \label{7.2.11}
\end{equation}
Таким образом,
\begin{equation}
\sum\limits_k \sum\limits_l \| r^{2-n-2k} \psi_{k, l} \chi_R
\|^2_{\mathring{H}^{s}_{\frac{n}{2}+k-1} (0, 2R)}
 \leq c  \sum\limits_k \sum\limits_l | \psi_{k, l}|^2  R^{- 2 k} (1+k)^{s-1}.
\label{7.2.12}
\end{equation}

Оценим теперь выражение $\| (1 - \chi_R) f \|_{H^s (\Omega)}.$

Обозначим через $\ov{R}$ диаметр области $\Omega.$ Тогда функция
$\chi_{\ov{R}} (r) = 1$ на $\Omega$ и поэтому
$$
\|  (1 - \chi_R) f \|_{H^s (\Omega)}^2 =  \| \chi_{\ov{R}} (1 -
\chi_R) f \|_{H^s (\Omega)} \leq
  \| \chi_{\ov{R}} (1 - \chi_R) f \|_{\mathring{H}^s (U_{2
\ov{R}})},
$$
где $U_{2 \ov{R}}$ "--- шар с центром в $\mathsf{0}$ радиуса ${2
\ov{R}}.$ Из теоремы~\ref{teo:2.1.2} вытекает формула
$$
 \| \chi_{\ov{R}} (1 - \chi_R) f \|_{\mathring{H}^s (U_{2 \ov{R}})}^2 = \sum\limits_{k=0}^{\infty} \sum\limits_{l=1}^{d_k}  \| \chi_{\ov{R}} (1 - \chi_R) f_{k, l} \|_{\mathring{H}^{s}_{\frac{n}{2}+k-1, +} (0, 2 \ov{R})}^2 =
$$
$$
=  \sum\limits_{k=0}^{\infty} \sum\limits_{l=1}^{d_k} |\psi_{k,
l}|^2  \| \chi_{\ov{R}} (1 - \chi_R) r^{2-n-2k}
\|_{\mathring{H}^{s}_{\frac{n}{2}+k-1, +} (0, 2 \ov{R})}^2 =
$$
$$
=  \sum\limits_{k=0}^{\infty} \sum\limits_{l=1}^{d_k} |\psi_{k,
l}|^2 \int\limits_R^{2 \ov{R}} |B^{\frac{s}{2}}_{\frac{n}{2}+k-1}
(\chi_{\ov{R}} (1 - \chi_R) r^{2-n-2k}) |^2   r^{2k+n-1} \, dr =
$$
$$
=  \sum\limits_{k=0}^{\infty} \sum\limits_{l=1}^{d_k} |\psi_{k,
l}|^2 \int\limits_R^{2 \ov{R}} |B^{\frac{s}{2}}_{1-\frac{n}{2}-k}
(\chi_{\ov{R}} (1 - \chi_R)) |^2   r^{3-n-2k} \, dr,
$$
где мы еще раз использовали соотношение Дарбу---Вайнштейна.

Нетрудно убедиться в справедливости следующей оценки:
$$
\max\limits_{r} |B^{\frac{s}{2}}_{1-\frac{n}{2}-k}  (\chi_{\ov{R}} (1 - \chi_R)) | \leq c \, (1+k)^s,
$$
где постоянная $c>0$ зависит от $s,$ $n,$ $R,$ $\ov{R},$ $\chi,$
но не зависит от $k.$ Объединяя две последние формулы, получаем
\begin{equation}
 \| \chi_{\ov{R}} (1 - \chi_R) f \|_{\mathring{H}^s (U_{2 \ov{R}})}^2 \leq c \sum\limits_{k=0}^{\infty} \sum\limits_{l=1}^{d_k} |\psi_{k, l}|^2 (1+k)^s \int\limits_R^{2 \ov{R}} r^{3-n-2k} \, dr
 \leq c' \, \sum\limits_{k=0}^{\infty} \sum\limits_{l=1}^{d_k}
|\psi_{k, l}|^2  (1+k)^{s-1}  R^{- 2 k}, \label{7.2.13}
\end{equation}
где постоянная $c' > 0$ не зависит от $\psi.$ При $h < \min\limits
(1, R)$ из~\eqref{7.2.12} и~\eqref{7.2.13} теперь получаем
$$
\| f \|^2_{s, R} =  \sum\limits_{k=0}^{\infty}
\sum\limits_{l=1}^{d_k} \| \chi_R\, r^{2-n-2k} \psi_{k, l}
\|^2_{\mathring{H}^{s}_{\frac{n}{2}+k-1} (0, 2R)} +
 \| (1 - \chi_R) f\|_{H^s (\Omega)}^2 \leq c \sum\limits_{k=0}^{\infty}
\sum\limits_{l=1}^{d_k} |\psi_{k, l}|^2 h^{-2 k} = c \, \| \psi
\|_h^2,
$$
где постоянная $c>0$ не зависит от $\psi.$ Таким образом,
оценка~\eqref{7.2.5} доказана при всех чётных $s \geq 2.$ При
доказательстве теоремы~\ref{teo: 7.1.2} было установлено
неравенство $\| f \|_{s', R}  \leq c  \| f \|_{s, R},$
справедливое при  $s'<s.$ Отсюда следует справедливость
оценки~\eqref{7.2.5} при всех $s \geq 0.$ Теорема доказана.
\end{proof}

В заключение этого раздела покажем, что вложение пространства $M^s
(\Omega_0)$ в $M^{s'} (\Omega_0)$ при $s>s',$ непрерывность
которого устанавливает теорема~\ref{teo: 7.1.2}, не является
вполне непрерывным.

В самом деле, пусть  $V$ "--- произвольное ограниченное множество
в пространстве  $A (\Theta).$ По обратной теореме о
$\sigma$-следах существует ограниченное множество $W,$
принадлежащее как пространству  $M^s (\Omega_0),$ так и
пространству $M^{s'} (\Omega_0),$ множество $\sigma$-следов
функций из которого совпадает с $V.$ Если указанный оператор
вложения пространства  $M^s (\Omega_0)$ в $M^{s'} (\Omega_0)$ был
бы вполне непрерывным, то множество $W$ было бы относительно
компактно в пространстве $M^{s'} (\Omega_0).$ Тогда по прямой
теореме о $\sigma$-следах его образ (множество $\sigma$-следов
функций из $W$) был бы относительно компактным в пространстве $A
(\Theta).$ Таким образом, любое ограниченное множество было бы
относительно компактным в пространстве $A (\Theta).$ Но это
возможно лишь в конечномерных пространствах, каковым $A (\Theta),$
очевидно, не является. Полученное противоречие и доказывает наше
утверждение.

\section{Новая краевая задача для уравнения Пуассона}\label{sec15}

\subsection[Постановка\! краевой\! задачи\! и\! изолированные\! особые\! точки\! гармонических\!
функций]{Постановка краевой задачи и изолированные особые точки гармонических
функций}\label{sec15.1}

Основной целью далее является изучение следующей краевой задачи. В
области $\Omega_0$ рассмотрим уравнение Пуассона
\begin{equation}
\Delta u= f(x), \  x \in \Omega_0. \label{8.1.1}
\end{equation}
Присоединим к этому уравнению краевое условие на $\pr \Omega_0$:
\begin{equation}
 u|_{\pr \Omega}  = g(x), \  x \in \pr \Omega,
\label{8.1.2}
\end{equation}
\begin{equation}
\sigma u|_{0} = \psi (\vartheta), \  \vartheta \in \Theta,
\label{8.1.3}
\end{equation}
где через $\sigma u |_0$  обозначен, как и ранее, $\sigma$-след
функции и в точке $\mathsf{0} \in \pr \Omega_0.$ Краевая
задача~\eqref{8.1.1}--\eqref{8.1.3} порождает оператор следующего
вида:
$$
\Lambda: u \to  \Lambda u = \{ \Delta u, u|_{\pr \Omega}, \sigma u|_0 \}.
$$

Снабдим пространство $M^s =  M^s (\Omega_0) \times H^{s+ \frac{3}{2}} (\pr \Omega) \times A (\Theta)$ топологией прямого произведения.

Из результатов предыдущего параграфа легко следует, что оператор
$\Lambda$ непрерывно отображает пространство $M^{s+2} (\Omega_0)$
в $M^s$ при любых $s \geq 0.$

Сформулируем основной результат этой главы.

\begin{theorem}\label{teo: 8.1.1}
    Оператор $\Lambda$ имеет обратный оператор $\Lambda^{-1},$ непрерывно отображающий пространство  $M^s$ на пространство  $M^{s+2} (\Omega_0)$ при любых чётных $s \geq 0.$
\end{theorem}

Доказательство теоремы состоит из нескольких этапов.

Сначала рассмотрим одно свойство особых точек гармонических
функций, представляющее, на наш взгляд, самостоятельный интерес.
Оказывается, что $\sigma$-след полностью характеризует поведение
сингулярной части гармонической функции в окрестности
изолированной особой точки. Более точно это означает следующее.

\begin{lemma}\label{lem: 8.1.1}
Для гармонической в $\Omega_0$ функции $u(x)$ существует
$\sigma$-след в точке $\mathsf{0},$ который принадлежит
пространству $A(\Theta).$ Обратно, для любой функции  $\psi \in
A(\Theta)$ существует единственная {\rm (}с точностью до
гармонической в окрестности $\mathsf{0}$ функции{\rm )}
гармоническая в $\Omega_0$ функция $u(x),$ для которой $\sigma u
|_0 = \psi.$
\end{lemma}

\begin{proof}
Пусть $u(x)$ "--- гармоническая в области $\Omega_0$ функция.
Тогда (см., например,~\cite{77}) при $n \geq 3$ имеет место
разложение
\begin{equation}
u (r, \vartheta) = \sum\limits_{k=0}^{\infty}
\sum\limits_{l=1}^{d_k} a_{k, l} r^k Y_{k, l} (\vartheta) +
\sum\limits_{k=0}^{\infty} \sum\limits_{l=1}^{d_k} b_{k, l}
r^{2-n-k} Y_{k, l} (\vartheta) =u_1+ u_2, \label{8.1.4}
\end{equation}
где первый ряд сходится абсолютно и равномерно в любом шаре с
центром в $\mathsf{0},$ лежащем вместе с замыканием в области
$\Omega,$ а второй "--- вне любого шара с центром в точке
$\mathsf{0}.$ Тогда $\sigma$-след функции $u_1,$ очевидно, равен
нулю. Для функции $u_2$ при $r> 0 $ имеем
$$
\sigma u_2 (r, \vartheta) = r^{n-2} \int\limits_{\Theta}    u_2
(r, \vartheta') \mathcal{K}_n  u_2 (r \vartheta, \vartheta') \, d
\vartheta'=
 \sum\limits_{k=0}^{\infty} \sum\limits_{l=1}^{d_k} b_{k, l} r^{-k}
\int\limits_{\Theta}   Y_{k, l} (\vartheta') \mathcal{K}_n   (r
\vartheta, \vartheta') \, d \vartheta',
$$
где $K_n (x, y)$ "--- ядро Пуассона для сферы $\Theta.$ Последний
интеграл есть интеграл Пуассона. Он определяет гармоническую в
единичном шаре функцию, которая на границе шара, то есть на сфере~$\Theta,$ равна функции $Y_{k, l} (\vartheta).$ Такой единственной
функцией будет функция   $r^k Y_{k, l} (\vartheta).$
Следовательно,
$$
\sigma u_2 (r, \vartheta) =  \sum\limits_{k=0}^{\infty}
\sum\limits_{l=1}^{d_k} b_{k, l} Y_{k, l} (\vartheta).
$$
Правая часть этого соотношения не зависит от переменной $r.$
Поэтому у неё существует предел при $r \to + 0,$ который по
определению и является $\sigma$-следом. Следовательно,
\begin{equation}
\sigma u|_0 = \sigma u_2|_0 =  \sum\limits_{k=0}^{\infty}
\sum\limits_{l=1}^{d_k} b_{k, l} Y_{k, l} (\vartheta).
\label{8.1.5}
\end{equation}

Далее, поскольку функция $u_2$ гармонична во всём пространстве
$E^n$ за исключением точки $\mathsf{0}$ и так как $b_{k, l}
r^{2-n-k}$ является коэффициентами разложения функции $u_2$ в ряд
по сферическим гармоникам, то для любого $r>0$ имеем формулу
$$
b_{k, l}  r^{2-n-k} = \int\limits_{\Theta}    u_2 (r, \vartheta')
Y_{k, l} (\vartheta) d \vartheta.
$$
Отсюда по неравенству Коши---Буняковского получаем оценку
$$
|b_{k, l}|  r^{2-n-k} \leq \lr{ \int\limits_{\Theta}    |u_2 (r,
\vartheta)|^2 \, d \vartheta}^{\frac{1}{2}}.
$$
Таким образом, $\sigma u_2|_0  \in A (\Theta).$ Первая часть леммы
доказана. Для доказательства второй части предположим, что функция
$\psi  \in A (\Theta).$ Тогда имеет место формула
$$
\psi (\vartheta) = \sum\limits_{k=0}^{\infty}
\sum\limits_{l=1}^{d_k} \psi_{k, l} Y_{k, l} (\vartheta),
$$
в которой коэффициенты $\psi_{k, l}$ допускают оценку $|\psi_{k,
l}| \leq c\, h^k,$ где $h$ "--- любое положительное число и
постоянная $c$ зависит от $h,$ но не зависит от параметров  $k$ и
$l.$ Среди функций вида~\eqref{8.1.4} требуется найти такую
функцию $u,$ чтобы $\sigma u_0 |_0= \psi,$  $\sigma$-след функции
$u$ совпадал с $\sigma$-следом функции $u_2$ и вычислялся по
формуле~\eqref{8.1.5}. В силу однозначности разложений по
сферическим гармоникам, в формуле~\eqref{8.1.4} необходимо
положить $b_{k, l}= \psi_{k, l}.$ Таким образом, мы показали, что
искомая функция необходимо должна иметь следующий вид:
$$
u (r, \vartheta) = u_1 (r, \vartheta) + \sum\limits_{k=0}^{\infty}
\sum\limits_{l=1}^{d_k} \psi_{k, l}  r^{2-n-k}  Y_{k, l}
(\vartheta),
$$
где $u_1$ "--- произвольная гармоническая в окрестности точки
$\mathsf{0}$ функция. Для сферических гармоник хорошо известна
оценка:
$$
| Y_{k, l} (\vartheta) | \leq c_n \, k^{\frac{n}{2}-1}, \  k = 0,
1, \dots; \  l = 1, \dots, d_k.
$$
Отсюда вытекает, что функциональный ряд в предыдущей формуле
сходится абсолютно и равномерно вне любого шара с центром в точке
$\mathsf{0}.$ Каждое слагаемое этого ряда представляет собой
гармоническую функцию, поэтому и сумма ряда есть гармоническая вне
точки $\mathsf{0}$ функция.

Для случая $n=2$ рассуждения вполне аналогичны с той лишь
разницей, что в этом случае в разложении в ряд по сферическим
гармоникам возникает одно логарифмическое слагаемое. Лемма
доказана.
\end{proof}

Отметим, что в термина $\sigma$-следа можно дать классификацию
изолированных особых точек гармонических функций. Так, если
$\left. \sigma u \right|_0 = 0,$  то $\mathsf{0}$ является
устранимой особой точкой. Если в разложении функции $\psi = \left.
\sigma u \right|_0 $ в ряд по сферическим гармоникам содержится
конечное число слагаемых, то функция $u$ имеет особенность типа
полюса. В противном случае $u$ имеет особенность типа существенной
особой точки.

\begin{lemma} \label{lem: 8.1.2}
    Пусть $f \in M^s (\Omega_0),$ $g \in H^{s+ \frac{3}{2}} (\pr \Omega),$ $\psi \in A (\Theta)$ и $s \geq 0.$
    Тогда у краевой задачи~\eqref{8.1.1}--\eqref{8.1.3} существует не более одного решения в пространстве $M^{s+2} (\Omega_0).$
\end{lemma}

\begin{proof}
Решение краевой задачи $\Delta v = 0,$ $\left. v \right|_{\pr
\Omega}=0,$ $\left. \sigma u \right|_0 = 0$ имеет в $\mathsf{0}$
устранимую особенность. Тогда $v \in H^{s+2} (\Omega)$ и
однородная краевая задача имеет лишь тривиальное решение.
\end{proof}

\subsection{Существование и априорная оценка
решения}\label{sec15.2}

В этом пункте будет завершено доказательство теоремы~\ref{teo:
8.1.1}.

\begin{lemma} \label{lem: 8.2.1}
Пусть $U_{\ov{R}, 0} = U_{\ov{R}} \setminus \mathsf{0},$ где
$U_{\ov{R}}$ "--- шар в $E^n$ радиуса $\ov{R}$ с центром в точке
$\mathsf{0}.$ Пусть $s \geq 0$ чётное число, функция $f \in M^s
(U_{\ov{R}, 0})$ обращается в нуль вблизи границы шара
$U_{\ov{R}}.$ Тогда существует функция $v \in M^{s+2} (U_{\ov{R},
0})$ такая,   что
\begin{equation}
\Delta v= f(x), \  x \in U_{\ov{R}, 0},
\label{8.2.1}
\end{equation}
\begin{equation}
\left. \sigma v \right|_0 = 0.
\label{8.2.2}
\end{equation}
При этом оператор $f \to v$ непрерывно действует из пространства
$M^s (U_{\ov{R}, 0})$ в пространство $M^{s+2} (U_{\ov{R}, 0}).$
\end{lemma}

\begin{proof}
Рассмотрим лишь случай  $n\geq 3,$  поскольку для $n=2$
рассуждения вполне аналогичны. Пусть сначала функция $f \in
\mathring{T}^{\infty} (U_{\ov{R}, 0}).$ Напомним, что это
означает, что имеет место разложение
\begin{equation}
    f (r, \vartheta) = \sum\limits_{k=0}^{\mathcal{K}} \sum\limits_{l=1}^{d_k} f_{k, l} (r) Y_{k, l} (\vartheta),
    \label{8.2.3}
\end{equation}
в котором $\mathcal{K}=\mathcal{K} (f)$ натуральное число, а
функция $r^{-k} f_{k, l} \in
\mathring{C}^{\infty}_{\frac{n}{2}+k-1} (0, \ov{R})$ (см.
главу~\ref{ch2}). Искомое решение $v$ в этом случае имеет вид
\begin{equation}
v (r, \vartheta) = - \sum\limits_{k=0}^{\mathcal{K}}
\sum\limits_{l=1}^{d_k} Y_{k, l} (\vartheta) r^k
\int\limits_r^{\ov{R}} t^{1-2k-n} \int\limits_0^t \tau^{n+k-1}
f_{k, l} (\tau) \, d \tau dt. \label{8.2.4}
\end{equation}
Функции из $ \mathring{C}^{\infty}_{\nu} (0, \ov{R})$ имеют
особенность в нуле не выше степенной порядка $- 2 \nu.$ Отсюда
следует, что повторный интеграл в~\eqref{8.2.4} есть величина
$O(r^{4-2k-n})$ и, следовательно,  выполнено~\eqref{8.2.2}.
Проверка~\eqref{8.2.1} осуществляется прямым дифференцированием
в~\eqref{8.2.4}.

Покажем, что оператор $f \to v,$ определенный
формулой~\eqref{8.2.4}, непрерывен в соответствующих топологиях.

Обозначим через $v_{k, l}$ функцию  вида
$$
v_{k, l} = - r^k \int\limits_r^{\ov{R}} t^{1-2k-n} \int\limits_0^t
\tau^{n+k-1} f_{k, l} (\tau) \, d \tau dt.
$$
Пусть $2 R < \ov{R}.$ Тогда
$$
v_{k, l} (r) = - r^k \int\limits_r^{\ov{R}} t^{1-2k-n}
\int\limits_0^t \tau^{n+k-1} \chi_{R/4} f_{k, l} (\tau) \, d \tau
dt -{}
$$
$$
{}- r^k \int\limits_r^{\ov{R}} t^{1-2k-n} \int\limits_0^t
\tau^{n+k-1} (1-\chi_{R/4}) f_{k, l} (\tau) \, d \tau dt = v_{k,
l}^1+v_{k, l}^2.
$$
Введем функции $v^j,$ $j=1,2$ по формуле
$$
v^j (r, \vartheta) = \sum\limits_{k=0}^{\mathcal{K}}
\sum\limits_{l=1}^{d_k} v^j_{k, l} (r) Y_{k, l} (\vartheta)
$$
и оценим в отдельности каждую из них. Рассмотрим сначала функцию
$v^1 (r, \vartheta).$  По  формуле Лейбница    имеем
\begin{equation}
B_{\frac{n}{2}+k-1} \lr{\chi_R \, r^{-k} v^1_{k,l}} = \chi_R
B_{\frac{n}{2}+k-1} \lr{r^{-k} v^1_{k, l}} + 2 \frac{\pr
\chi_R}{\pr r} \frac{\pr r^{-k} v^1_{k, l}}{\pr r} + r^{-k}
v^1_{k, l} B_{\frac{n}{2}+k-1} \chi_R. \label{8.2.5}
\end{equation}
Для первого слагаемого справа имеем следующую формулу:
$$
\chi_R  B_{\frac{n}{2}+k-1} \lr{ r^{-k} v^1_{k,l}} = \chi_R \,
r^{-k} \chi_{R/4} f_{k, l} (r) =
 \chi_{R/4} (r) r^{-k} f_{k, l} (r),
$$
так как $\chi_R  \chi_{R/4} =  \chi_{R/4}.$ Учитывая, что $D
\chi_R (r) =0$ при $0 \leq r \leq R,$ а  $\chi_{R/4} (r) = 0$ при
$r \geq \dfrac{R}{2},$ для второго слагаемого в~\eqref{8.2.5}
получаем  выражение
\begin{equation}
2 \frac{\pr \chi_R}{\pr r} \frac{\pr (r^{-k} v^1_{k, l})}{\pr r} =
2  \frac{\pr \chi_R}{\pr r} r^{1-2k-n} \int\limits_0^{R/2}
\tau^{n+k-1} \chi_{R/4} f_{k, l} \, d \tau. \label{8.2.6}
\end{equation}
Из тех же соображений получаем представление третьего слагаемого
\begin{multline}
r^{-k} v^1_{k, l} B_{\frac{n}{2}+k-1} \chi_R (r)
= - \lr{B_{\frac{n}{2}+k-1} \chi_R} \int\limits_r^{\ov{R}} t^{1-2k-n} \int\limits_0^t \tau^{n+k-1} \chi_{R/4} f_{k, l} (\tau) \, d \tau dt =\\
= \lr{B_{\frac{n}{2}+k-1} \chi_R}
\frac{r^{2-2k-n}-\ov{R}^{2-2k-n}}{2-2k-n} \int\limits_0^{R/2}
\tau^{n+k-1} \chi_{R/4} f_{k, l} (\tau) \, d \tau. \label{8.2.7}
\end{multline}
Учитывая, что
$$
\frac{2}{r} D \chi_R +  \frac{1}{2-2k-n} B_{\frac{n}{2}+k-1}
\chi_R = \frac{-1}{n+2k-2} B_{1-\frac{n}{2}-k} \, \chi_R,
$$
из~\eqref{8.2.5}--\eqref{8.2.7} находим
\begin{multline}
 B_{\frac{n}{2}+k-1} \lr{ \chi_R (r) r^{-k} v^1_{k, l}} = \chi_{R/4} r^{-k} f_{k, l} - \\
- \frac{1}{n+2k-2}  \lr{B_{1-\frac{n}{2}-k} \chi_R} r^{2-2k-n} \int\limits_0^{R/2} \tau^{n+k-1} \chi_{R/4} f_{k, l} (\tau) \, d \tau +\\
+ \frac{1}{n+2k-2} \ov{R}^{2-2k-n}  \lr{B_{\frac{n}{2}+k-1} \chi_R
}  \int\limits_0^{R/2} \tau^{n+k-1} \chi_{R/4} f_{k, l} (\tau) \,
d \tau. \label{8.2.8}
\end{multline}
Теперь в силу формулы~\eqref{7.1.2} получаем следующую оценку
функции $v^1$:
\begin{multline}\label{8.2.9}
\| v^1 \|^2_{s+2, R} = \sum\limits_k \sum\limits_l \| \chi_R \,
r^{-k} v^1_{k,l} \|^2_{\mathring{H}^{s+2}_{\frac{n}{2}+k-1 (0,
2R)}} + \| (1-\chi_R)  v^1 \|^2_{{H}^{s+2} (U_{\ov{R}})} \leq
\\
\leq 3  \sum\limits_k \sum\limits_l \| \chi_{R/4} \, r^{-k}
f_{k,l} \|^2_{\mathring{H}^{s}_{\frac{n}{2}+k-1} (0, 2R)} +
3  \sum\limits_k \sum\limits_l \frac{1}{(n+2k-2)^2} \left(  \| r^{2-2k-n} B_{1-\frac{n}{2}-k} \chi_R \|^2_{\mathring{H}^{s}_{\frac{n}{2}+k-1 (0, 2R)}} + \right.  \\
+ \left. \ov{R}^{4-4k-2n} \|B_{\frac{n}{2}+k-1} \chi_R
\|^2_{\mathring{H}^{s}_{\frac{n}{2}+k-1} (0, 2R)}   \right) \left|
\int\limits_0^{R/2} \tau^{n+k-1} \chi_{R/4} f_{k, l} \, d \tau
\right|^2 + \|(1-\chi_R) v^1 \|^2_{H^{s+2} (U_{\ov{R}})}=
\\
=3 \Sigma_1 +3 \Sigma_2 + \|(1-\chi_R) v^1 \|^2_{H^{s+2}
(U_{\ov{R}})}.
\end{multline}

Оценим каждое слагаемое справа в последней формуле. Слагаемое
$\Sigma_1$ допускает очевидную оценку
\begin{equation}
\Sigma_1 \leq \| f \|^2_{s, R/4}. \label{8.2.10}
\end{equation}

Рассмотрим интеграл, присутствующий во втором слагаемом. Введем
обозначение $\nu = \dfrac{n}{2}+k-1,$ $\omega (r) =
\chi_{\frac{R}{4}} r^{-k} f_{k, l}.$ Тогда указанный интеграл
примет вид
$$
Q_{\nu} (\omega, R) = \int\limits_0^{R/2} \tau^{2 \nu +1}
\omega(\tau) \, d \tau =  \int\limits_0^{R/2} \tau^{2 \nu +1}
P_{\nu}^{\frac{1}{2}-\nu} S_{\nu}^{\nu-\frac{1}{2}} \omega \, d
\tau,
$$
где $P_{\nu}^{\frac{1}{2}-\nu}$ и $S_{\nu}^{\nu-\frac{1}{2}}$ "---
операторы преобразования из пункта~\ref{sec4.1}. Так как $r^{-k}
f_{k, l} \in \mathring{C}^{\infty}_{\frac{n}{2}+k-1} (0, \ov{R}),$
то по определению класса $\mathring{C}^{\infty}_{\nu} (0, \ov{R})$
функции $S_{\nu}^{\nu-\frac{1}{2}} \omega  =
S_{\nu}^{\nu-\frac{1}{2}} (\chi_R \, r^{-k} f_{k, l})$ принадлежит
пространству $\mathring{C}^{\infty} [0, R).$ Пусть
$\widetilde{\omega} = S_{\nu}^{\nu-\frac{1}{2}} \omega$ и пусть
$\nu< N + \dfrac{1}{2},$ где $N$ "--- натуральное число. Тогда по
определению операторов преобразования получаем
$$
P_{\nu}^{\frac{1}{2}-\nu} \widetilde{\omega} (\tau) = \frac{(-1)^N
2^{-N} \sqrt{\pi} \, \tau^{2(N-\nu)}}{\Gamma(\nu+1) \,
\Gamma(N-\nu+\frac{1}{2})}  \lr{\frac{\pr}{\pr \tau}
\frac{1}{\tau}}^N \int\limits_{\tau}^{\infty} \tau^{2 \nu}
(t^2-\tau^2)^{N-\nu-\frac{1}{2}} t^{-2 N} \widetilde{\omega} (t)
\, d t.
$$
Отсюда интегрированием по частям получаем
$$
Q_{\nu} = \int\limits_0^{R/2} \tau^{2 \nu +1} \omega(\tau) \, d
\tau =
$$
$$
=  \frac{(-1)^N 2^{-N} \sqrt{\pi}}{\Gamma(\nu+1) \,
\Gamma(N-\nu+\frac{1}{2})} \int\limits_0^{R/2} \tau^{2 \nu +1}
\lr{\frac{\pr}{\pr \tau} \frac{1}{\tau}}^N
\int\limits_{\tau}^{\infty} \tau^{2 \nu}
(t^2-\tau^2)^{N-\nu-\frac{1}{2}} t^{-2 N} \widetilde{\omega} (t)
\, d t d \tau =
$$
$$
= \frac{ \sqrt{\pi} \, \Gamma(N+\frac{3}{2})}{\Gamma(\nu+1) \,
\Gamma(N-\nu+\frac{1}{2})} \int\limits_0^{R/2} \tau^{2 \nu +1}
\int\limits_{\tau}^{\infty}  (t^2-\tau^2)^{N-\nu-\frac{1}{2}}
t^{-2 N} \widetilde{\omega} (t) \, d t d \tau =
$$
$$
= \frac{ \sqrt{\pi} \, \Gamma(N+\frac{3}{2})}{\Gamma(\nu+1) \,
\Gamma(N-\nu+\frac{1}{2}) \, \Gamma(\frac{3}{2})}
\int\limits_0^{R/2}    \widetilde{\omega} (t) t^{-2 N}
\int\limits_{0}^{t} \tau^{2 \nu +1}
(t^2-\tau^2)^{N-\nu-\frac{1}{2}}  \, d \tau d t.
$$
Внутренний интеграл выражается через функции Эйлера:
$$
\int\limits_{0}^{t} \tau^{2 \nu +1}
(t^2-\tau^2)^{N-\nu-\frac{1}{2}}  \, d \tau = t^{2 N+1}
\frac{\Gamma(\nu+1) \, \Gamma(N-\nu+\frac{1}{2})}{ 2 \,
\Gamma(N+\frac{3}{2})}.
$$
Следовательно,
$$
Q_{\nu} = \int\limits_0^{R/2} t \widetilde{\omega}(t) \, d t =
\int\limits_0^{R/2} t S_{\nu}^{\nu-\frac{1}{2}}   \omega(t) \, d
t.
$$
Далее, поскольку $S_{\nu} = I^{\frac{1}{2} - \nu}
S_{\nu}^{\nu-\frac{1}{2}},$ где $I^{\mu}$ "--- лиувиллевский
оператор, то из предыдущей формулы получаем
$$
Q_{\nu} = \int\limits_0^{R/2} t  I^{s+\nu-\frac{1}{2} } I^{-s}
S_{\nu} \, \omega (t)  \, d t =
 \frac{1}{\Gamma\lr{s+\nu-\frac{1}{2}}} \int\limits_0^{R/2} t
\int\limits_t^{R/2} (\tau-t)^{s+\nu-\frac{3}{2}} I^{-s} S_{\nu} \,
\omega (t)  \, d \tau d t =
$$
$$
= \frac{1}{\Gamma\lr{s+\nu-\frac{1}{2}}} \int\limits_0^{R/2}
\lr{I^{-s} S_{\nu} \, \omega(\tau)} \int\limits_0^{\tau} t
(\tau-t)^{s+\nu-\frac{3}{2}}   \,  d t d \tau.
$$
Так как
$$
 \int\limits_0^{\tau} t (\tau-t)^{s+\nu-\frac{3}{2}}   \,  d t  =\tau^{s+\nu+\frac{1}{2}} \frac{\Gamma(s+\nu-\frac{1}{2})}{\Gamma(s+\nu+\frac{3}{2})},
$$
то
$$
Q_{\nu} =  \frac{1}{\Gamma(s+\nu+\frac{3}{2})} \int\limits_0^{R/2}
\tau^{s+\nu+\frac{1}{2}} I^{-s} S_{\nu} \, \omega(\tau) \, d \tau.
$$

Теперь по неравенству Коши---Буняковского получаем оценку
$$
|Q_{\nu}| \leq  \frac{1}{\Gamma(s+\nu+\frac{3}{2})} \lr{
\int\limits_0^{R/2} \tau^{2 s+2\nu+1} \, dt}^{1/2} \|D^s S_{\nu}\,
\omega \|_{L_2 (0, R/2)} =
$$
$$
= \frac{R^{s+\nu+1}}{2^{s+\nu+\frac{3}{2}}\sqrt{s+\nu+1} \, \Gamma
(s+\nu+\frac{3}{2})}  \|S_{\nu} \, \omega \|_{\mathring{H}^s (0,
R/2)} =
 \frac{R^{s+\nu+1} \, \Gamma(\nu+1)}{2^{s+1} \sqrt{s+\nu+1} \,
\Gamma (s+\nu+\frac{3}{2})}  \| \omega \|_{\mathring{H}^s_{\nu}
(0, R/2)}.
$$
Из соотношения~\eqref{1.4.20} для гамма-функций отсюда следует
неравенство
$$
|Q_{\nu}| \leq c (s, R) R^{\nu} \nu^{-1-s} \| \omega
\|_{\mathring{H}^s_{\nu} (0, R/2)}.
$$
Возвращаясь к старым обозначениям, выпишем полученную
окончательную оценку интеграла $Q_{\nu}$ в следующем виде:
\begin{equation}
|Q_{\frac{n}{2}+k-1} (\chi_{R/4} \, r^{-k} f_{k, l}, R)| \leq c
(n, s, R)  \frac{R^k}{(k+1)^{1+s}} \| \chi_{R/4} \, r^{-k} f_{k,
l} \|_{\mathring{H}^s_{\frac{n}{2}+k-1} (0, R/2)}. \label{8.2.11}
\end{equation}
Из результатов пункта~\ref{sec4.1} следует, что
$$
\| B_{\frac{n}{2}+k-1} \chi_R
\|^2_{\mathring{H}^s_{\frac{n}{2}+k-1} (0, 2 R)} \leq 2 \|
B_{\frac{n}{2}+k-1} \chi_R \|^2_{\mathring{H}^s_{\frac{n}{2}+k-1,
+} (0, 2 R)} =
 2 \int\limits_0^{2R} |B_{\frac{n}{2}+k-1}^{\frac{s}{2}+1} \chi_R
|^2 \, r^{2 (\frac{n}{2}+k-1)+1} dr =
$$
$$
= 2 R^{n+2k-s-2} \int\limits_0^{2}
|B_{\frac{n}{2}+k-1}^{\frac{s}{2}+1} \chi(t) |^2 \, t^{n+2k-1} dt
=
 c (s, n, k) \, R^{n+2k-s-2} \leq c (s, n) \, 2^{2 k} \,
R^{n+2k-s-2} (k+1)^{s+1},
$$
а при доказательстве теоремы~\ref{teo: 7.2.2} было получено
неравенство
$$
\| r^{2-2k-n} B_{1-\frac{n}{2}-k} \chi_R
\|_{\mathring{H}^s_{\frac{n}{2}+k-1} (0, 2 R)} \leq c (s, n, R) \,
R^{-2 k} (k+1)^{s+1}.
$$
Два последних неравенства с учётом~\eqref{8.2.11} приводят к
следующей оценке слагаемого из формулы~\eqref{8.2.9}:
$$
\Sigma_2 \leq c \sum\limits_k \sum\limits_l R^{2k} (k+1)^{-4-2s}
\left( R^{-2 k} (k+1)^{s+1} + {}\right.
$$
$$
\left.{} + \ov{R}^{4-4k-2n} R^{n+2k-s-2} \, 2^{2k} (k+1)^{s+1}
\right)   \| \chi_{R/4} \, r^{-k} f_{k, l}
\|^2_{\mathring{H}^s_{\frac{n}{2}+k-1} (0, R/2)}.
$$
Так как $2 R < \ov{R},$ то мы приходим к окончательной оценке
слагаемого $\Sigma_2$:
\begin{equation}
\Sigma_2 \leq c \sum\limits_k \sum\limits_l \| \chi_{R/4} \,
r^{-k} f_{k, l} \|^2_{\mathring{H}^s_{\frac{n}{2}+k-1} (0, R/2)}
\leq c \, \|f \|^2_{s, R/4}, \label{8.2.12}
\end{equation}
где постоянная $c>0$ не зависит от $f.$

Для завершения оценки функции $v^1$ осталось рассмотреть последнее
слагаемое в~\eqref{8.2.9}. Поскольку $\chi_{\ov{R}} (r) = 1$ в
шаре $U_{\ov{R}},$ то
\begin{equation}
\| (1 - \chi_R) v^1 \|_{H^{s+2} (U_{\ov{R}}) } \leq c  \|
\chi_{\ov{R}} (1 - \chi_R) v^1 \|_{\mathring{H}^{s+2} (U_{2
\ov{R}}) }. \label{8.2.13}
\end{equation}
По аналогии с формулами~\eqref{8.2.5}--\eqref{8.2.8} имеем
$$
B_{\frac{n}{2}+k-1} (\chi_{\ov{R}} (1 - \chi_R) r^{-k} v^1_{k, l})
=
$$
$$
= \frac{1}{2k+n-2}  \lr{B_{\frac{n}{2}+k-1} (\chi_{\ov{R}} (1 -
\chi_R)) \ov{R}^{2-2k-n} - B_{1-\frac{n}{2}-k} (\chi_{\ov{R}} (1 -
\chi_R) r^{2-2k-n}}   \int\limits_{0}^{R/2} \tau^{n+k-1} \chi_{R/4} f_{k, l} \, d
\tau.
$$
Следовательно,
$$
\| \chi_{\ov{R}} (1 - \chi_R) v^1 \|^2_{\mathring{H}^{s+2} (U_{2
\ov{R}}) } \leq \sum\limits_k \sum\limits_l
\frac{|Q_{\frac{n}{2}+k-1}|^2}{(2k+n-2)^2}  \left( \| r^{2-2k-n}
B_{1-\frac{n}{2}-k}^{\frac{s+2}{2}} (\chi_{\ov{R}} (1 - \chi_R))
\|^2_{L_2,  \frac{n}{2}+k-1}  +{} \right.
$$
$$
\left.{} +  \ov{R}^{4-4k-2n}  \|
B_{\frac{n}{2}+k-1}^{\frac{s+2}{2}} (\chi_{\ov{R}} (1 - \chi_R))
\|_{L_2,  \frac{n}{2}+k-1}  \right).
$$
Нормы из последней суммы уже оценивались нами, в частности, при
доказательстве теоремы~\ref{teo: 7.2.2}. Откуда имеем
$$
\| \chi_{\ov{R}} (1 - \chi_R) v^1 \|^2_{\mathring{H}^{s+2} (U_{2
\ov{R}}) } \leq c\, \sum\limits_k \sum\limits_l (k+1)^{-s-1}  \|
\chi_{R/4} \, r^{-k} f_{k, l}
\|^2_{\mathring{H}^s_{\frac{n}{2}+k-1} (0, R/2)} \leq c \, \| f
\|^2_{s, R/4}.
$$
Следовательно,
\begin{equation}
\| (1 - \chi_R) v^1 \|_{{H}^{s+2} (U_{ \ov{R}}) } \leq c \, \| f
\|^2_{s, R/4}. \label{8.2.14}
\end{equation}
Стало быть, (см.~\eqref{8.2.10},~\eqref{8.2.12},~\eqref{8.2.14})
нами доказана следующая оценка функции $v^1$:
\begin{equation}
\| v^1 \|_{s+2, R} \leq c \, \| f \|_{s, R/4}, \label{8.2.15}
\end{equation}
постоянная $c>0$ в которой не зависит от функции $f.$

Рассмотрим функцию $v^2.$ Нетрудно  заметить, что она принадлежит
пространству ${H}^{s+2} (U_{ \ov{R}})$ и является решением
следующей краевой задачи:
$$
\Delta v^2 = (1 - \chi_{R/4}) f(x), \  x \in U_{ \ov{R}},
$$
$$
\left. v^2 \right|_{\pr U_{\ov{R}}} = 0.
$$
Причём так как $f \in M^s (U_{ \ov{R}, 0}),$ то, очевидно, $(1 -
\chi_{R/4}) f \in {H}^{s} (U_{ \ov{R}}).$ Хорошо известно, что
решение этой краевой задачи Дирихле для уравнения Пуассона
единственно и для него справедлива оценка
$$
\| v^2 \|_{{H}^{s+2} (U_{ \ov{R}}) } \leq c \, \| (1 - \chi_{R/4})
f \|_{{H}^{s} (U_{ \ov{R}}) }.
$$
Теперь из определения норм $\| \|_{s, R}$ и из теоремы~\ref{teo:
7.1.1} получаем
\begin{equation}
\| v^2 \|_{s+2, R} \leq c \, \| f \|_{s, R/4}, \label{8.2.16}
\end{equation}
где постоянная не зависит от функции $f.$

Таким образом, нами получены оценки каждого слагаемого разложения
$v=v^1+v^2.$ Следовательно, для любого $R \in \Big(0,
\dfrac{\ov{R}}{2}\Big)$ и $s \geq 0$ найдётся такая постоянная
$c>0,$ что для любой функции $f \in \mathring{T}^{\infty}
(U_{\ov{R}, 0})$ имеет место неравенство
$$
\| v \|_{s+2, R} \leq c \, \| f \|_{s, R/4}.
$$

Завершим доказательство леммы предельным переходом. Пусть функция
$f \in M^s (U_{\ov{R}, 0})$ и удовлетворяет условию леммы. Тогда
найдётся последовательность функций $f^m \in \mathring{T}^{\infty}
(U_{\ov{R}, 0})$ сходящаяся к $f$ по топологии этого пространства.
Для каждой функции $f^m$ определим функции $v^m$ по
формуле~\eqref{8.2.4}. Тогда
\begin{equation}
\Delta v^m = f^m \to f, \  m \to \infty. \label{8.2.17}
\end{equation}
По доказанному, отображение $f^m \to v^m$ непрерывно действует из
пространства $M^s (U_{\ov{R}, 0})$ в $M^{s+2} (U_{\ov{R}, 0}).$
Тогда последовательность $v^m$ фундаментальна в  $M^{s+2}
(U_{\ov{R}, 0}).$ Отсюда, в силу полноты пространства  $M^{s+2}
(U_{\ov{R}, 0})$ существует функция  $v \in M^{s+2} (U_{\ov{R},
0}),$ являющаяся пределом последовательности функций $v^m$ в
смысле топологии этого пространства.

Оператор $\Delta$ непрерывно отображает пространство $M^{s+2}
(U_{\ov{R}, 0})$ в  $M^{s} (U_{\ov{R}, 0}).$  Поэтому $\Delta v^m
\to \Delta v,$ $m \to \infty,$ в смысле пространства $M^s
(U_{\ov{R}}).$    Тогда   из~\eqref{8.2.17} следует, что $\Delta v
= f.$ По
 прямой теореме о $\sigma$-следах имеем
$$
\left. \lim\limits_{m \to \infty} \sigma v^m \right|_0 = \left.
\sigma v \right|_0,
$$
где предельный переход понимается в смысле пространства $A
(\Theta).$ Поскольку $\left. \sigma v^m \right|_0=0,$ то тогда и
$\left. \sigma v \right|_0=0.$  Лемма~\ref{lem: 8.2.1} полностью
доказана.
\end{proof}

\begin{proof}
[Доказательство теоремы~\ref{teo: 8.1.1}] Утверждение о
единственности решения краевой задачи~\eqref{8.1.1}--\eqref{8.1.3}
установлено в лемме~\ref{lem: 8.1.2}. Существование решения
докажем следующим образом. Пусть $\ov{R}$ "---  диаметр области
$\Omega$ и пусть функция $u^1 \in M^{s+2} (U_{\ov{R}, 0})$
является решением краевой задачи
\begin{equation}
\begin{cases}
\Delta u^1 =  \chi_{R_0}\, f(x), &  x \in U_{ \ov{R}, 0},\\
\left. \sigma u^1 \right|_{0} = 0,&
\end{cases}
\label{8.2.18}
\end{equation}
построенном в лемме~\ref{lem: 8.1.1}. Пусть функция $u^2$ является
решением краевой задачи
\begin{equation}
\begin{cases}
\Delta u^2 = (1- \chi_{R_0}) f(x), &  x \in \Omega,\\
\left. u^2 \right|_{\pr \Omega} = g - \left. u^1 \right|_{\pr
\Omega} - \left. u^3 \right|_{\pr \Omega},&
\end{cases}
\label{8.2.19}
\end{equation}
где через $u^3$ обозначена гармоническая в $E^n \setminus
\mathsf{0}$ функция, построенная в теореме~\ref{teo: 7.2.2},
такая, что
$$
\left. \sigma u^3 \right|_{0} = \psi.
$$

Тогда функция $u = u^1+u^2+u^3$ является решением краевой
задачи~\eqref{8.1.1}--\eqref{8.1.3}.

Функция $u$ принадлежит пространству $M^{s+2} (\Omega_0).$ Это
следует из того, что  $u^1 \in M^{s+2} (\Omega_0)$ по
лемме~\ref{lem: 8.2.1}, $u^3 \in M^{s+2} (\Omega_0)$ по
теореме~\ref{teo: 7.2.2},  $u^2 \in H^{s+2} (\Omega) \subset
M^{s+2} (\Omega_0)$ по общей теории эллиптических    краевых задач
(см., например,~\cite{51, 62}). Из этих результатов вытекает и
непрерывность разрешающего оператора $\Lambda^{-1}.$
Теорема~\ref{teo: 8.1.1} доказана.
\end{proof}

\newpage

\chapter{Композиционный метод построения операторов
преобразования}\label{ch6}

\section{Общая схема композиционного метода построения операторов
преобразования}\label{sec16}

Исследованные ранее классы ОП строились каждый своими методами.
Поэтому возникает необходимость в разработке общей схемы
построения ОП. Такая схема "--- метод факторизации или
композиционный метод (КМ) "--- предлагается в настоящей главе.
Метод основан на представлении ОП в виде композиции интегральных
преобразований. Композиционный метод даёт алгоритмы не только для
построения новых ОП, но содержит как частные случаи ОП СПД,
Векуа---Эрдейи---Лаундеса, Бушмана---Эрдейи различных типов,
унитарные ОП Сонина---Катрахова и Пуассона---Катрахова, обобщённые
операторы Эрдейи---Кобера, а также введённые
Р.~Кэрролом~\cite{Car1, Car2, Car3} классы эллиптических,
гиперболических и параболических ОП. В этой главе вводятся их
обобщения: классы $B$-эллиптических, $B$-гиперболических и
$B$-параболических ОП.

Композиционный метод построения операторов преобразования был
разработан и последовательно развивался в работах
С.\,М.~Ситника~\cite{S6,S66,S7,S5,S46,S14,S400,S42,S38,S401,SitDis,FJSS}.
Отметим, что первоначальная идея, близкая к этому методу, была
применена В.\,В.~Катраховым для построения одного специального
класса операторов преобразования, которые в этой книге названы
операторами преобразования Сонина---Катрахова и
Пуассона---Катрахова, см. главу~\ref{ch3}.

Для перевода на английский язык названия композиционного метода
предложен такой вариант: {\it Integral Transforms Composition Method}
(ITCM), см.~\cite{FJSS}.

В этом пункте приведены основные определения и небольшая часть результатов, которые могут быть получены с использованием КМ. Так как сами формулировки результатов достаточно громоздки, они приводятся без доказательств. Полное изложение композиционного метода с многочисленными примерами планируется в дальнейшем.

Общая схема композиционного метода следующая. На вход подаётся
пара операторов произвольного вида $A,B,$ а также связанные с ними
обобщённые преобразования Фурье $F_A, F_B,$ которые обратимы и
действуют по формулам
\begin{equation}
\label{4301} F_A A =g(t) F_A,\qquad  F_B B= g(t) F_B,
\end{equation}
где $t$ "--- двойственная переменная, $g(t)$ "--- произвольная
подходящая функция, самый очевидный выбор "--- это $g(t)=-t^2,$
как для классических интегральных преобразований.

\medskip

\textit{Целью композиционного метода является формальное
построение на выходе пары взаимно обратных операторов
преобразования  $P$ и $S$ по следующим формулам{\rm :}
\begin{equation}\label{4302}
S=F^{-1}_B \frac{1}{w(t)} F_A,\qquad P=F^{-1}_A w(t) F_B
\end{equation}
с произвольной весовой функцией $w(t).$ Тогда $P$ и $S$ являются
взаимно обратными операторами преобразования, сплетающими исходные
операторы $A$ и $B${\rm :}}
\begin{equation}\label{Inter}
SA=BS,\qquad PB=AP.
\end{equation}

\medskip

Формальная проверка получается прямой подстановкой. Основную трудность представляет вычисление введённых композиций в явном интегральном виде, а также нахождение соответствующих областей определений операторов.

Перечислим основные достоинства композиционного метода.

\begin{itemize}
\item{Простота метода "--- многочисленные операторы преобразования
строятся как из строительных блоков по простому правилу из
элементарных кирпичиков: классических интегральных
преобразований.} \item{Метод позволяет построить единым способом
все известные ранее в явном виде операторы преобразования.}
\item{Метод позволяет построить единым способом многие новые
классы операторов преобразования.} \item{Метод позволяет просто
получать обратные преобразования в том же композиционном виде.}
\item{Метод позволяет выводить оценки норм прямых и обратных ОП,
используя известные оценки норм для классических интегральных
преобразований в различных функциональных пространствах.} \item{КМ
позволяет получать явные формулы, связывающие решения сплетаемых
дифференциальных уравнений.}
\end{itemize}

Можно указать и один недостаток КМ: действие кирпичиков "---
интегральных преобразований "--- известно как правило в
классических пространствах, а в применении к дифференциальным
уравнениям или дифференциальным операторам  для конкретных задач
оценки и действие на функциях требуются на других классах,
например с занулениями в начале координат или на бесконечности,
или и там и там. В этих ситуациях КМ можно применять для получения
явной формы ОП, а затем распространять его на нужные пространства.

Подчеркнём, что формулы вида~\eqref{4301}-\eqref{4302},
разумеется, сами по себе не являются новыми ни для теории
интегральных преобразований, ни для приложений к дифференциальным
уравнениям. \textit{Но композиционный метод является новым для
приложений в теории операторов преобразований}\,!

В других задачах для интегральных преобразований и связанных с
ними дифференциальных уравнений формулы композиционного
метода~\eqref{4302} при выборе классического преобразования Фурье
приводят к классу псевдодифференциальных операторов с символами
$w(t), \dfrac{1}{w(t)},$ см., например,~\cite{OlRa, Agr}. При
выборе  опять же в качестве интегрального преобразования
классического преобразования Фурье $A=B=D^2,\ \ F_A=F_B=F$ и
весовой функции $w(t)=(\pm it)^{-s}$ мы получаем дробные интегралы
Римана---Лиувилля на всей действительной прямой, при выборе
весовой функции  $w(t)=|t|^{-s}$ получаются потенциалы М.~Рисса,
при выборе  $w(t)=(1+t^2)^{-s}$ по формулам~\eqref{4302}
получаются потенциалы Бесселя, а при $w(t)=(1\pm it)^{-s}$ "---
модифицированные потенциалы Бесселя, см.~\cite{SKM}.

При выборе в качестве интегрального преобразования классического
преобразования Ханкеля и весовых функций
\begin{equation}\label{trans} A=B=B_\nu,\quad  F_A=F_B=H_\nu,\quad
g(t)=-t^2,\quad  w(t)=j_\nu(st)
\end{equation}
получаются операторы обобщённого сдвига (ООС) Дельсарта,
см.~\cite{Lev2, Lev3, Lev4, Mar9}; напомним обозначения для
оператора Бесселя $B_\nu,$ преобразования Ханкеля $H_\nu,$
нормированной или <<малой>> функции Бесселя $j_\nu(\cdot),$ см.
вводную главу~\ref{ch1}. Отметим, что в общем случае при выборе
данных для композиционного метода $A=B, F_A=F_B$ с произвольными
весовыми функциями $g(t), w(t)$ мы будем получать операторы
преобразования, коммутирующие с данным оператором $A,$ так же, как
и ООС коммутирует с оператором Бесселя.

Интересным вариантом приложения композиционного  метода является
выбор в качестве оператора $A$ выражения для квантового
осциллятора, а в качестве соответствующего интегрального
преобразования $F_A$ выбор квадратичного преобразования Фурье
(КПФ, дробное преобразование Фурье, преобразование
Фурье---Френеля, преобразование Вейерштрасса),
см.~\cite{OZK,AbOs,Os}. Это важное интегральное преобразование
недостаточно широко известно (пока), оно возникло из предложения
Френеля заменить стандартные плоские волны с линейными аргументами
в экспонентах на более общие волны с квадратичными аргументами в
экспонентах, что позволило полностью объяснить парадоксы со
спектральными линиями при дифракции Фраунгофера. Математически
операторы КПФ являются дробными степенями $F^\alpha$ обычного
преобразования Фурье, достраивая его до полугруппы по параметру
$\alpha,$ они были определены Н.~Винером и А.~Вейлем. В теории
всплесков, в которой принято каждую формулу считать новой и
называть по-новому давно известные вещи, КПФ называется
преобразованием Габора. Интересные приложения этого интегрального
преобразования к группе Гейзенберга, квантовым осцилляторам и
всплескам недавно получены в~\cite{AlKi1}. Изложенный выше
композиционный метод позволяет с помощью этого преобразования
строить ОП для одномерного оператора Шрёдингера "--- квантового
осциллятора~\cite{S45,S42}. При этом может быть использовано и
более общее квадратичное преобразование Ханкеля~\cite{S60}.

Применение композиционного метода вместо классических подходов
позволяет также строить в явном интегральном виде различные формы
дробных степеней оператора Бесселя,
см.~\cite{S135,S133,S127,S123,S18,S42,S700,SS,FJSS}. Такое
построение невозможно методами теории полугрупп или спектральными
методами.

Распространение композиционного метода на многомерный случай
очевидно, при этом $t$ является вектором, а $g(t),w(t)$ в
формулах~\eqref{4301}-\eqref{4302} становятся вектор-функциями. К
сожалению, в этом случае известны или могут быть введены явно лишь
немногие классы ОП. Но к их числу относятся хорошо известные
классы потенциалов. Например, в том случае, когда в композиционном
методе используются классическое преобразование Фурье и в качестве
весовой функции $w(t)$ в~\eqref{4302} выбирается положительно
определённая квадратичная форма, получаются эллиптические
потенциалы М.~Рисса~\cite{Riesz,SKM}; если $w(t)$ "---
знаконеопределённая квадратичная форма, то получаются
гиперболические потенциалы М.~Рисса,~\cite{Riesz,SKM,Nogin}; при
выборе $w(x,t)=(|x|^2-it)^{-\alpha/2}$ "--- параболические
потенциалы~\cite{SKM}. При использовании в композиционном методе в
формулах~\eqref{4301}-\eqref{4302} преобразования Ханкеля и в
качестве $w(t)$ квадратичной формы мы приходим к
эллиптическим~\cite{Lyah3,Gul1} или гиперболическим~\cite{ShiE2}
$B$-потенциалам М.~Рисса. В этих случаях применяется теория
обобщённых функций и их свёрток, а для обращения таких потенциалов
нужны специальные процедуры осреднения и аппроксимаций,
см.~\cite{Nogin,ShiE2}, а в качестве пространств используются
пространства Шварца или Лизоркина пробных функций и двойственные
пространства для распределений.

Одним из важных применений композиционного метода являются формулы связи между решениями возмущённых и невозмущённых дифференциальных уравнений, которые удалось связать с помощью ОП. Например, для уравнений с операторами Бесселя вида
\begin{equation}\label{GenEPD}
\sum\limits_{k=1}^n A_k\left(   \frac{\partial^2 u}{\partial x_k^2}+\frac{\nu_k}{x_k}\frac{\partial u}{\partial x_k}\right) \pm \lambda^2 u=0,
\end{equation}
$$
x_k>0,\qquad A_k=\const,\qquad \nu_k=\const,\qquad \lambda=\const,
$$
композиционным методом устанавливаются формулы связи с решениями невозмущённого уравнения
\begin{equation}\label{GenWave}
\sum\limits_{k=1}^n A_k  \frac{\partial^2 v}{\partial x_k^2}\pm \lambda^2 v=0.
\end{equation}
В частности, к этому классу соотношений относятся формулы связи
решений для уравнений с операторами Бесселя с различными
параметрами, они называются операторами сдвига по
параметру~\cite{FJSS}. В теории уравнений дробного порядка
используется аналогичная идея, которая называется <<принципом
субординации>>, по существу состоящая в построении ОП,
осуществляющих сдвиг по порядку уравнения, что, в частности,
позволяет выразить решения уравнений дробного порядка через
решения классических уравнений целого порядка, см.,
например,~\cite{Jan,EiIvKoch,Bajlekova0,Bajlekova1,FJSS}.

Таким образом, можно сделать вывод, что композиционный метод для построения операторов преобразования является эффективным, он связан с рядом других известных методов и задач, с его помощью получаются все известные в литературе явные представления для операторов преобразования, а также могут быть построены многочисленные новые классы операторов преобразования. Приложения алгоритма композиционного метода проводятся в три этапа или шага.

\begin{itemize}
\item Шаг 1. Для данной пары операторов $A,B$ и связанных с ними
обобщённых преобразований Фурье $F_A, F_B$ определяем и вычисляем
пару взаимно обратных операторов преобразования $P,S$ по основным
формулам метода~\eqref{4301}-\eqref{4302}.

\item Шаг 2. Находим точные условия и определяем классы функций,
для которых операторы преобразования, построенные на шаге 1,
удовлетворяют сплетающим свойствам~\eqref{Inter}.

\itemШаг 3. Применяем корректно определённые на шагах 1 и 2
операторы преобразования на соответствующих пространствах функций
к решению задач для дифференциальных уравнений, например, к
установлению формул соответствия для решений возмущённых и
невозмущённых дифференциальных уравнений.
\end{itemize}

Далее на основе композиционного метода строятся гиперболические,
эллиптические и параболические по Р.~Кэрролу операторы
преобразования, обобщённые операторы Эрдейи---Кобера и другие.
Отметим, что в форме представления решений одних абстрактных
уравнений через другие ОП указанных типов рассматривались
в~\cite{Lav1, BD1, BD2, BD3}. В частности, применяя этот вариант
метода ОП получены несколько необычные формулы, которые выражают
решения волнового уравнения через решения уравнения
теплопроводности и наоборот. Для случая абстрактного оператора
Бесселя подобные формулы получены в~\cite{Glu13}.

Мы используем классификацию операторов преобразования,
предложенную Р.~Кэрролом и связанную с типом уравнения в частных
производных, которому удовлетворяет ядро этого оператора. Мы также
вводим по аналогии и используем естественные обобщения этих
понятий: $B$-эллиптические, $B$-гиперболические и
$B$-параболические операторы преобразования. При этом используются
определения для соответствующих классов сингулярных уравнений в
частных производных, которые были введены И.\,А.~Киприяновым.
Например, классические сплетающие операторы
Сонина---Пуассона---Дельсарта по этой терминологии являются
$B$-гиперболическими.

\subsection{$B$-гиперболические операторы
преобразования}\label{sec16.1}

Так будем называть операторы сдвига по параметру оператора Бесселя, удовлетворяющие соотношению
\begin{equation}\label{449}{T B_{\nu} = B_{\mu} T.}\end{equation}
Будем искать такие операторы в факторизованном виде
\begin{equation}\label{4410}{T_{\nu, \, \mu}^{(\varphi)} = F_{\mu}^{-1}
\lr{\varphi(t) F_{\nu}}.}\end{equation}

Если $\nu=-\dfrac{1}{2}$ или $\mu=-\dfrac{1}{2},$ то такие
операторы сводятся к уже изученным. Будем полагать $\varphi (t) =
t^{\alpha},$ $T^{(\varphi)}=T^{(\alpha)}.$

\begin{theorem} Пусть выполнены условия
    $$
    -2-2 \Re \mu < \Re \alpha < \Re(\nu-\mu).
    $$
    Тогда справедливо интегральное представление
    \begin{multline}\label{4411}
    \lr{T^{(\alpha)}_{\nu,\, \mu} f}(x) = \frac{C_1}{x^{2 \mu + \alpha +2}}  \int\limits_0^{x} y^{2 \nu +1}\, {_2F_1}\lr{\mu+\frac{\alpha}{2}+1, \frac{\alpha}{2}+1; \nu+1; \frac{y^2}{x^2}} f(y)  \, dy + {} \\
    {}+ C_2 \int\limits_x^{\infty}  y^{-2\mu+2 \nu-\alpha -1}    {_2F_1}\lr{\mu+\frac{\alpha}{2}+1, \mu-\nu+\frac{\alpha}{2}+1; \mu+1; \frac{x^2}{y^2}} f(y) \,dy,
    \end{multline}
    $$
    C_1=   \frac{2^{ \mu-\nu+\alpha+1} \Gamma \lr{\mu+\frac{\alpha}{2}+1}}{\Gamma \lr{-\frac{\alpha}{2}} \Gamma \lr{\nu+1}},\qquad C_2 = \frac{2^{ \mu-\nu+\alpha+1} \Gamma \lr{\mu+\frac{\alpha}{2}+1}}{\Gamma \lr{\nu-\mu-\frac{\alpha}{2}} \Gamma \lr{\mu+1}},  \nonumber
   $$
    где ${_2F_1}$ "--- гипергеометрическая функция Гаусса.
\end{theorem}

Рассмотрим несколько частных случаев оператора~\eqref{4411}.

\medskip

\noindent а) Пусть $\alpha=-1-2\mu+\nu,$ $\Re \nu > -1,$ $\Re \mu
> -1.$ Тогда
\begin{multline*}
    \lr{T^{(-1-2\mu+\nu)}_{\nu,\, \mu} f}(x) = 2^{- \mu} \frac{1}{x^{ \mu}}  \int\limits_x^{\infty} y^{\nu} (y^2-x^2)^{\frac{\mu}{2}} P_{\frac{\nu}{2}-\frac{1}{2}}^{- \mu} \lr{1-2 \frac{x^2}{y^2}}  f(y)  \, dy + {} \\
    {}+2^{1- \mu} e^{i \mu \pi}  \frac{\Gamma \lr{\mu+\frac{\nu}{2}+\frac{1}{2}}}{\Gamma \lr{\mu-\frac{\nu}{2}+\frac{1}{2}}} \frac{1}{x^{ \mu}}  \int\limits_0^{x} y^{\nu} (x^2-y^2)^{\frac{\mu}{2}} Q_{\frac{\nu}{2}-\frac{1}{2}}^{ \mu} \lr{2 \frac{x^2}{y^2}-1} f(y) \,dy.
\end{multline*}
б) Пусть $\alpha=0;$ $-1< \Re \mu < \Re \nu.$ В этом случае
получается замечательный оператор <<спуска>> по параметру, который
не зависит от начального и конечного значений параметров $\nu,$
$\mu,$ а зависит лишь от величины <<спуска>>  $\gamma = \nu -
\mu$:
$$
\lr{T^{(0)}_{\nu,\, \mu} f}(x) = \frac{2^{1-(\nu-\mu)}}{\Gamma(\nu-\mu)} \int\limits_x^{\infty} y (y^2-x^2)^{\nu-\mu-1} f(y)  \, dy.
$$
Этот оператор лишь числовым множителем отличается от дробного
интеграла $I_{-, x^2}$ по функции $g(x)=x^2.$ В такой форме он был
открыт А. Эрдейи, это частный случай операторов Эрдейи---Кобера
или Дж. Лаундеса.

\medskip

\noindent в) Пусть $\alpha=2 \nu,$ $-1< \Re (\nu + \mu) < 0.$
Тогда
$$
\lr{T^{(2 \nu)}_{\nu,\, \mu} f}(x) = \frac{\sin (\pi \mu)}{\pi} 2^{\mu+\nu+1} \Gamma(\nu+\mu+1) \int\limits_x^{\infty} y^{2 \nu +1} (y^2-x^2)^{-\mu-\nu-1} f(y)  \, dy-
$$
$$
- \frac{\sin (\pi \nu)}{\pi} 2^{\mu+\nu+1} \Gamma(\nu+\mu+1) \int\limits_0^x y^{2 \nu +1} (x^2-y^2)^{-\mu-\nu-1} f(y)  \, dy.
$$
г) Пусть  $\mu = \nu,$ $-2 \Re \nu -2 < \Re \alpha < 0.$ Тогда
получаем семейство операторов перестановочных с $B_{\nu}$:
$$
\lr{T^{(\alpha)}_{\nu,\, \nu} f}(x) =
\frac{2^{\alpha+2}}{\sqrt{\pi}} \frac{e^{-i \frac{\pi \alpha}{2}}}{\Gamma\lr{-\frac{1}{2}-\frac{\alpha}{2}}}\frac{1}{x^{\nu+\frac{1}{2}}} \int\limits_0^{\infty} y^{\nu + \frac{3}{2}}|x^2-y^2|^{-\frac{\alpha}{2}-1} Q_{\nu - \frac{1}{2}}^{\frac{\alpha}{2}+1} \lr{\frac{x^2+y^2}{2xy}} f(y) \, dy.
$$
д) Пусть $\mu = - \nu,$ $2 \Re \nu -2 < \Re \alpha < \Re 2
\nu.$ Тогда
$$
\lr{T^{(\alpha)}_{\nu, -\nu} f}(x) = 2^{-2\nu+\alpha+1}
\frac{\Gamma\lr{- \nu+\frac{\alpha}{2} +1}}{\Gamma\lr{2
\nu-\frac{\alpha}{2}}} \, x^{\nu} \int\limits_0^{\infty}
|y^2-x^2|^{\nu-\frac{\alpha}{2}-1} y^{\nu + 1} P_{\nu -
\frac{\alpha}{2}-1}^{\nu} \lr{\frac{x^2+y^2}{|x^2-y^2|}} f(y) \,
dy.
$$

Другой подход  к построению операторов сдвига по параметру
заключается в следующем. Пусть по приведённым выше формулам
построен оператор типа Сонина, для которого выбрано
$\varphi=\varphi_1,$ и  оператор типа Пуассона при значениях
$\nu=\mu,$ $\varphi=\varphi_2.$ Тогда их композиция $T=P_{\mu}
S_{\nu}$ и будет искомым оператором. При этом в случае, если
одновременно выбраны или синус или косинус "--- преобразования
Фурье, получится в точности конструкция~\eqref{4410}, где
$\varphi=\dfrac{\varphi_2}{\varphi_1}.$ Поэтому вопрос об
ограниченности такого оператора в лебеговых пространствах со
степенным весом $L_{p, \gamma}(0, \infty)$ сводится к задаче о
возможности деления в пространствах мультипликаторов. В случае
если выбраны разные преобразования, мы приходим к следующей
факторизации: \begin{equation}\label{4412}{T^{(\varphi_1,
\varphi_2)}_{\nu,\, \mu} = F^{-1}_{\mu} \cdot \varphi_2(t)\cdot
F_{\scriptsize\left\{
\begin{matrix} s
\\
c \end{matrix} \right\}} \cdot F^{-1}_{\scriptsize\left\{
\begin{matrix} c
\\
s \end{matrix} \right\}} \frac{1}{\varphi_1(t)}
F_{\nu}.}\end{equation} Несложно показать, что композиция
преобразований Фурье сводится к так называемым преобразованиям
Гильберта на полуоси
\begin{eqnarray}
& & (F_s F_c f)(x) = \int\limits_0^{\infty} \frac{x}{x^2-y^2} f(y) \, dy,  \label{4413}\\
& & (F_c F_s f)(x) = \int\limits_0^{\infty} \frac{y}{y^2-x^2} f(y)
\, dy, \label{4414}
\end{eqnarray}
где интеграл понимается в смысле главного значения, причём
указанные преобразования, домноженные на нужные постоянные,
унитарны в $L_2 (0, \infty).$ Таким образом, в рассматриваемом
случае оператор~\eqref{4412} факторизуется через два
преобразования Фурье---Бес\-селя и одно из двухвесовых
преобразований
$$
(A_1 f)(x) = x\, \varphi_2 (x) \int\limits_0^{\infty} \frac{f(y)}{(x^2-y^2)} \cdot \frac{dy}{\varphi_1 (x)},  \qquad
(A_2 f)(x) = \varphi_2 (x) \int\limits_0^{\infty} \frac{y}{\varphi_1 (x)} \cdot \frac{f(y)}{(y^2-x^2)} \,  dy.
$$
Вопрос об ограниченности~\eqref{4412} в этом случае сводится к
задаче об оценках с двумя весами для полуосевых преобразований
Гильберта в соответствующих пространствах.

Другие классы $B$-гиперболических операторов преобразования можно
построить, если использовать вместо преобразования Фурье---Бесселя
$F_{\nu}$ преобразование с функцией Неймана $Y_{\nu}(z)$ в ядре.

Отметим, что операторы сдвига по параметрам для операторов Бесселя
находят важные приложения в теории сингулярных дифференциальных
уравнений, см.~\cite{FJSS}, а также работы
Ш.\,Т.~Каримова~\cite{KarST, KarST2, KarST3, KarST4}.

\subsection{$B$-эллиптические операторы
преобразования}\label{sec16.2}

Эти операторы удовлетворяют соотношению
\begin{equation}\label{4415}{T B_{\nu} = - D^2 T.}\end{equation} Этот необычный
класс ОП осуществляет связь между решениями $B$-эллиптических и
$B$-гиперболических дифференциальных уравнений, выражая их решения
друг через друга.

Построение таких операторов основано на замене в предыдущих
факторизациях синус- и косинус-преобразований Фурье на
преобразование Лапласа или на замене преобразования
Фурье---Бесселя на одно из преобразований с функциями Макдональда
и Неймана
\begin{eqnarray}
& &  (K_{\nu} f) (t) = \frac{1}{t^{\nu}} \int\limits_0^{\infty} y^{\nu+1} K_{\nu}(t y) f(y) \, dy, \label{4416} \\
& &  (Y_{\nu} f) (t) = \frac{1}{t^{\nu}} \int\limits_0^{\infty}
y^{\nu+1} Y_{\nu}(t y) f(y) \, dy. \label{4417}
\end{eqnarray}
Возможно также использование преобразования с функцией Струве в ядре.

\begin{theorem}  Пусть $|\Re \nu| + \Re (\alpha+\nu) < 1.$ Тогда оператор
    $$
    \lr{A_{\nu}^{\alpha} f} (x) = F_c^{-1} t^{- \alpha} K_{\nu} f
    $$
    является $B$-эллиптическим, удовлетворяющим соотношению~\eqref{4415}. Для него справедливо интегральное представление
    $$
    (A^{\alpha}_{\nu} f) (x) = \frac{\pi \Gamma(1-\alpha)}{4 \sin \frac{\pi}{2}  (1-\alpha-2 \nu)} \int\limits_0^{\infty} y^{\nu+1} (x^2+y^2)^{\frac{\alpha+\nu-1}{2}} \times
    $$
    $$
    \times \left[ P_{-\alpha-\nu}^{-\nu} \lr{\frac{x}{\sqrt{x^2+y^2}}} + P_{-\alpha-\nu}^{-\nu} \lr{ - \frac{x}{\sqrt{x^2+y^2}}}  \right] f(y) \, dy.
    $$
\end{theorem}

Определим оператор, удовлетворяющий соотношению~\eqref{4415}, по
формуле
$$
(C^{\alpha}_{\nu} f) (x) = L (t^{-\alpha} F_{\nu} f),
$$
где $L$ "--- преобразование Лапласа.

\begin{theorem} Пусть $\Re \alpha < 1.$ Тогда справедливо интегральное представление
    $$
    (A^{\alpha}_{\nu} f) (x) = \Gamma(1-\alpha) \int\limits_0^{\infty} y^{\nu+1} (x^2+y^2)^{\frac{\alpha+\nu-1}{2}} P_{-\alpha-\nu}^{-\nu} \lr{\frac{x}{\sqrt{x^2+y^2}}} f(y) \, dy.
    $$
\end{theorem}

Аналогичное представление при $|\Re \nu| + \Re (\alpha+\nu) < 1$ справедливо и для оператора
$$
L (t^{-\alpha} Y_{\nu} f)(x) =
 - \frac{2}{\pi} \Gamma(1-\alpha) \int\limits_0^{\infty} y^{\nu+1} (x^2+y^2)^{\frac{\alpha+\nu-1}{2}} Q_{-\alpha-\nu}^{-\nu} \lr{\frac{x}{\sqrt{x^2+y^2}}} f(y) \, dy.
$$

Подобные формулы выведены и для более широкого класса
$B$-эллипти\-ческих операторов преобразования, сплетающих
$B_{\nu}$ и $\lr{-B_{\mu}}.$

Рассмотрим простейшие из введённых выше операторов
$A_{\nu}^{\alpha}$ и $C_{\nu}^{\alpha}$ при значениях $\nu = \pm
\dfrac{1}{2}.$ Эти операторы будут определяться по формулам
$$
\lr{A^{\beta} f}(x) = \lr{F^{-1}_{\left\{ \begin{matrix} s \\
c \end{matrix} \right\}} t^{\beta} L} f, \qquad \lr{C^{\beta}
f}(x) = \lr{L  t^{-\beta} F_{\left\{ \begin{matrix} s \\
c
\end{matrix} \right\}} } f,
$$
где $L$ "--- преобразование Лапласа. Они сплетают $D^2$ и $-D^2.$
В этом случае проще вычислить ядра интегральных операторов
непосредственно. Это приводит к формулам
$$
\lr{C^{\beta} f}(x) = \sqrt{\frac{2}{\pi}} \Gamma (1-\beta) \int\limits_0^{\infty} \frac{ f(y)}{\lr{x^2+y^2}^{\frac{1-\beta}{2}}} \left\{ \begin{matrix} \sin \\
\cos \end{matrix} \right\}  \Big[(1-\beta) \arctg \frac{y}{x}\Big]
dy,
$$

$\Re \beta < 1 + \delta,$ $\delta=1$ при выборе синуса,
$\delta=0$ при выборе косинуса.
$$
\lr{A^{\beta} f}(x) = \sqrt{\frac{2}{\pi}} \Gamma (1+\beta) \int\limits_0^{\infty} \frac{ f(y)}{\lr{x^2+y^2}^{\frac{1+\beta}{2}}} \left\{ \begin{matrix} \sin \\
\cos \end{matrix} \right\}  \Big[(1+\beta) \arctg \frac{x}{y}\Big]
dy,
$$
где $\Re \beta > - \delta - 1,$ $\delta$ определенно выше. В
частности, при $\beta=0$ получаем пару операторов преобразования,
связанную с интегралами Пуассона для полупространства:
$$
\lr{C^{0} f}(x) = \sqrt{\frac{2}{\pi}} \int\limits_0^{\infty} \frac{y f(y)}{x^2+y^2} \, dy,
\qquad
\lr{A^{0} f}(x) = \sqrt{\frac{2}{\pi}} \int\limits_0^{\infty} \frac{x f(y)}{x^2+y^2}\, dy.
$$
Эти операторы и операторы $C^{\beta},$ $A^{\beta}$ для частных
значений $\beta \in \mathbb{N}$ построены в~\cite{Car1}.

\subsection{$B$-параболические операторы
преобразования}\label{sec16.3}

Этот необычный класс ОП позволяет выражать решения параболических уравнений через гиперболические и наоборот.

Введем интегральные преобразования по формулам
$$
\lr{F'_c f} (x) =  \lr{F_c f} (\sqrt{x}), ~  \lr{F'_s f} (x) =  \lr{F_s f} (\sqrt{x}),
$$
$$
\lr{P f} (x) = \lr{L \varphi (t) F'_{\scriptsize\left\{ \begin{matrix} s \\
c \end{matrix} \right\}}}(x).
$$

Тогда на финитных функциях оператор $P$ сплетает вторую и первую производные по формулам
$$
P D^2 f = D P f.
$$
Таким образом, этот оператор является параболическим по
терминологии Р.~Кэррола.

\subsection{Операторы сдвига по спектральному параметру типа
Лаундеса}\label{sec16.4}

Операторы этого типа появились при установлении формул, выражающих
решения уравнений Гельмгольца через гармонические функции. Их
изучение, начатое в работах И.\,Н.~Векуа  и А.~Эрдейи  было
продолжено в работах Дж.~Лаундеса. Поэтому С.\,М.~Ситником было
предложено название для этого класса "--- <<операторы
преобразования Векуа---Эрдейи---Лаундеса>>,
см.~\cite{S66,S6,S46,S59,S125}. Подобные операторы также
применялись к решению сингулярных дифференциальных уравнений в
работах Ш.\,Т.~Каримова, см., например,~\cite{KarST}.

Рассмотрим оператор
$$
T_1 = F^{-1}_{\nu} \lr{ \varphi (t) F'_{\mu}},
$$
где введено преобразование с параметром $\lambda$
$$
\lr{F'_{\mu} f} (t) = \frac{1}{t^{\nu}} \int\limits_0^{\infty} y^{\nu+1} J_{\nu}(y \sqrt{t^2+\lambda^2}) f(y) \, dy.
$$
Оператор $T_1$ удовлетворяет соотношению
\begin{equation}\label{4420}{T_1 B_{\mu} = (B_{\nu} - \lambda^2) T_1.}\end{equation}

\begin{theorem} При условиях $-1< \Re \nu<  \Re \mu$ и выборе  $\varphi (x) = x^{\mu} (x^2+\lambda^2)^{-\frac{\mu}{2}}$ справедливо интегральное представление
    $$
    \lr{T_1 f} (x) = \lambda^{1+\nu-\mu} \int\limits_0^{\infty} y (y^2-x^2)^{\frac{\mu-\nu-1}{2}} J_{\nu-\mu-1}(\lambda \sqrt{y^2-x^2}) f(y) \, dy.
    $$
\end{theorem}

Приведем ряд других операторов, получаемых в факторизованном виде.

\medskip

\noindent а) Пусть $\nu=1,$ $\varphi (x) = x^{\mu-2}
(x^2+\lambda^2)^{\frac{\mu}{2}}.$

Тогда оператор $T_2= F_1^{-1} \varphi (t) F'_{\mu},$
удовлетворяющий соотношению
$$
T_2 B_{\mu} = (B_1 -\lambda^2)\, T_2,
$$
при $\Re \lambda > 0,$  $\Re \mu > -1$ может быть представлен
в виде
$$
(T_2 f)(x)= \frac{1}{x^2 \lambda^2} \int\limits_0^{\infty} y^{\nu+1} J_{\mu}(\lambda y) f(y) \, dy -
 \frac{1}{x^2 \lambda^2} \int\limits_x^{\infty} y (y^2-x^2)^{\frac{\mu}{2}} J_{\mu}(\lambda \sqrt{y^2-x^2}) f(y) \, dy.
$$
б) Пусть  $\varphi (x) = x^{\mu} (x^{\frac{\mu}{2}}+\lambda^2)^2,$
$\Re \lambda > 0,$ и пусть $T_3 = F^{-1}_{\nu} \varphi
F'_{\mu}.$ Тогда если $-1< \Re \nu< - \Re \mu,$ то
$$
(T_3 f) (x) =   \frac{2 \sin (\pi \mu)}{\pi} \lambda^{\mu+\nu+1}  \int\limits_0^x y^{2 \mu + 1} (x^2-y^2)^{-\frac{\mu+\nu+1}{2}} K_{\mu+\nu+1}(\lambda \sqrt{x^2-y^2}) f(y) \, dy +
$$
$$
+ \lambda^{\nu+\mu+1} \int\limits_x^{\infty} y^{2 \mu + 1} (y^2-x^2)^{-\frac{\mu+\nu+1}{2}}  \left[ \sin (\pi \nu) Y_{\mu+\nu+1}(\lambda \sqrt{y^2-x^2}) -
 \cos (\pi \nu) J_{\nu+\mu+1}(\lambda \sqrt{y^2-x^2}) \right]  f(y) \, dy.
$$
в) Пусть  $\varphi (x) = x^{\mu-1}/(x^2+\lambda^2) (x^2+\frac{\lambda}{2}^2)^{\frac{\mu}{2}}.$\\
Положим $T_4 = Y^{-1}_{\nu} \varphi F'_{\mu}.$ Тогда если $\Re \lambda > 0,$ $-\dfrac{1}{2}< \Re \nu< 3+ \Re \mu,$ то
справедливо соотношение
$$
(T_4 f) (x) =  - \frac{\lambda^{\nu-\mu-1}}{2^{- \frac{\mu}{2}}} \cdot \frac{K_{\nu} (\lambda x)}{x^{\nu}}  \int\limits_0^{\infty} y^{\mu+1} J_{\mu}\lr{\frac{\lambda y}{\sqrt{2}} } f(y) \, dy.
$$
Приведенные примеры демонстрируют важность свободы выбора функции
$\varphi.$

Аналогично осуществляется построение операторов $T,$
удовлетворяющих соотношению \begin{equation}\label{4421}{T B_{\mu}
= (B_{\mu} + \lambda^2) \, T.}\end{equation} Например, оператор
вида
$$
T_5 = (F'_{\mu})^{-1} \varphi F_{\nu}
$$
для одной из функций $\varphi$  является также оператором Дж. Лаундеса
$$
(T_4 f)(x) = \frac{\lambda^{\nu-\mu+1}}{x^{2 \mu}} \int\limits_0^{\infty} y^{2 \nu +1} (x^2-y^2)^{\frac{\mu-\nu-1}{2}} J_{\mu-\nu-1}(\lambda \sqrt{y^2-x^2}) f(y) \, dy.
$$

Отметим, что соотношения~\eqref{4420}-\eqref{4421} ввиду
линейности входящих в них операторов могут быть переписаны в форме
$$
T_1 (B_{\mu} + \lambda^2) = B_{\nu} T,\quad T (B_{\nu} -
\lambda^2) = B_{\mu} T.
$$
Кроме того, во всех операторах можно обосновать замену $\lambda
\to i \lambda.$ Соотношение для наиболее общих операторов
преобразования вида \begin{equation}\label{22}{T (B_{\nu}+\alpha)
= (B_{\mu} + \beta) T}\end{equation} эквивалентно уже
рассмотренным
$$
T(B_{\nu}+\alpha-\beta) = B_{\mu}T,\quad  (B_{\mu}+\beta-\alpha)T
=T B_{\nu}.
$$
Впрочем, операторы, для которых выполнены предыдущие соотношения,
могут быть получены и непосредственно. Таким же образом строятся и
$B$-эллиптические операторы типа Лаундеса, удовлетворяющие
соотношению
$$
T (B_{\nu}+\lambda) = (-B_{\mu} + \beta)\, T.
$$

При выборе значений параметров $\nu = \mu = - \dfrac{1}{2}$
получаем операторы, сплетающие $D^2$ и $D^2 \pm \lambda^2.$

Укажем на возможность применения КМ при выборе в качестве
интегрального преобразования квадратичного (дробного)
преобразования Фурье или Ханкеля. Для этих целей наиболее
интересны соотношения~\eqref{Ht^2D}-\eqref{x^2HD}. По общей схеме
КМ мы можем сконструировать ОП, сплетающие дифференциальные
операторы $D^2$ и
$$
\left(\sin^{2}{\frac{\alpha}{2}}L_{\nu} - \frac{1}{2}i\sin{\alpha}\left(XD+DX\right) - X^{2}\cos^{2}{\frac{\alpha}{2}}
\right),
$$
где
$$
L_{\nu}= -\frac{1}{4}D^2-\frac{\nu^2-1/4}{x^2} + \frac{1}{4}x^2 - \frac{\nu+1}{2}
$$
с произвольными параметрами $\alpha, \nu.$ Для этого в наших
обозначениях  надо положить
$$
A=\left(\sin^{2}{\frac{\alpha}{2}}L_{\nu} - \frac{1}{2}i\sin{\alpha}\left(XD+DX\right) - X^{2}\cos^{2}{\frac{\alpha}{2}}\right),\quad B=D^2,
$$
$$
F(A)=H_\nu^\alpha, F(B)=F_c, g(t)=-t^2,
$$
где $H_\nu^\alpha$ "--- квадратичное преобразование
Фурье---Френеля, $F_c$ "--- косинус преобразование Фурье.

Отметим, что КМ может быть также с успехом применён к построению
дробных степеней оператора Бесселя. С его помощью может быть
получено решение многих интегродифференциальных уравнений,
см.~\cite{S46, SitDis}.

\chapter{Приложения метода операторов преобразования к оценкам решений для дифференциальных уравнений с переменными коэффициентами и задаче
Е.\,М.~Ландиса}\label{ch7}

\chaptermarknum{Приложения метода операторов преобразования}

\section{Приложения метода операторов преобразования для возмущённого уравнения Бесселя с переменным
потенциалом}\label{sec17}
\sectionmarknum{\sП\sр\sи\sл\sо\sж\sе\sн\sи\sя\s{ }\sм\sе\sт\sо\sд\sа\s{ }\sо\sп\sе\sр\sа\sт\sо\sр\sо\sв\s{
}\sп\sр\sе\sо\sб\sр\sа\sз\sо\sв\sа\sн\sи\sя\s{ }\sд\sл\sя\s{
}\sв\sо\sз\sм\sу\sщ\sё\sн\sн\sо\sг\sо\s{
}\sу\sр\sа\sв\sн\sе\sн\sи\sя\s{ }\sБ\sе\sс\sс\sе\sл\sя\s{ }\sс\s{
}\sп\sе\sр\sе\sм\sе\sн\sн\sы\sм{
}\sп\sо\sт\sе\sн\sц\sи\sа\sл\sо\sм}

Рассматривается задача о построении интегральной формулы для решений дифференциального уравнения с определённой асимптотикой
\begin{equation}\label{4.2.1}
B_{\alpha} g(x) - q(x) g(x)=\lambda^2 g(x),
\end{equation}
где  $B_\alpha$ "--- оператор Бесселя, который в данном пункте нам
удобно определить в таком виде
\begin{equation}\label{4.2.2}
B_{\alpha}g =g''(x)+\frac{2\alpha}{x}g'(x),~ \alpha>0.
\end{equation}

Заметим, что в работе используются разные формы записи постоянных в операторе Бесселя в зависимости от решаемых задач, это объясняется более простыми формулами для записи получаемых результатов; это не должно привести к путанице, так как рассматриваемые задачи в разных местах работы не имеют пересечений.

Данная задача решается методом операторов преобразования. Для этого достаточно построить пару взаимно обратных операторов преобразования,  первый из которых $S_\alpha$ вида
\begin{equation}\label{4.2.3}
S_\alpha h(x)=h(x)+\int\limits_x^{\infty}S(x,t)h(t)\,dt,
\end{equation}
сплетает операторы
$B_\alpha - q(x)$ и $B_\alpha$ по формуле
\begin{equation}\label{4.2.4}
S_{\alpha}(B_{\alpha}-q(x))h=B_{\alpha}S_{\alpha}h,
\end{equation} а второй $P_\alpha,$ обратный к первому, должен
быть построен в виде интегрального с ядром $P(x,t)$
\begin{equation}\label{4.2.3*}
P_\alpha h(x)=h(x)+\int\limits_x^{\infty}P(x,t)h(t)\,dt \end{equation}
и действовать по формуле
$$
P_\alpha B_{\alpha}h=(B_{\alpha}-q(x))P_\alpha h
$$
на функциях $h \in C^2(0, \infty).$

В результате на решениях дифференциального уравнения~\eqref{4.2.1}
функция $S_{\alpha}u=v$ будет выражаться через решения
невозмущённого уравнения, получаемого отбрасыванием слагаемого с
потенциалом в~\eqref{4.2.1}, то есть фактически через функции
Бесселя, а функция $u=P_\alpha v$ будет являться решением
исходного возмущённого уравнения~\eqref{4.2.1}. При этом для
решения будет получено интегральное представление~\eqref{4.2.3*} с
явным описанием ядра $P(x,t).$ Эта методика отражает одно из
основных применений ОП "--- выражение решений более сложных
дифференциальных уравнений через подобные более простые, что уже
несколько раз отмечалось ранее. Также заметим, что в силу
использования линейных ОП одна и та же пара взаимно обратных ОП
позволяет получить как представления решений дифференциального
уравнения~\eqref{4.2.1} со спектральным параметром, так и
представление решений в более простом случае однородного уравнения
$$
B_{\alpha} h(x) - q(x) h(x)=0.
$$
При этом, если ставится задача о нахождении представления для
решений возмущённого уравнения~\eqref{4.2.1}, то построение
прямого оператора преобразования можно пропустить, а сразу перейти
к построению обратного и нахождению интегрального представления
для искомого решения вида~\eqref{4.2.3*}.

Оригинальная методика для решения подобных задач была разработана
В.\,В.~Сташевской~\cite{Sta1,Sta2}, что позволило ей включить в
рассмотрение сингулярные потенциалы с предельно точной оценкой в
нуле $|q(x)| \leq c x^{- 3/2+\varepsilon},$ $\varepsilon
> 0$ при целых $\alpha.$ Эта методика, основанная на применении
обобщённых теорем Пэли---Винера, получила широкое развитие и
признание. Случай непрерывной~$q$ при $\alpha>0$  рассмотрен  в
работах А.\,С.~Сохина~\cite{Soh1, Soh2, Soh3, Soh4}, а
также~\cite{Volk}. При этом в работах В.\,В.~Сташевской и
В.\,Я.~Волка строились операторы преобразования типа Повзнера с
интегрированием по конечному промежутку, а в работах А.\,С.~Сохина
"--- типа Левина с интегрированием по бесконечному промежутку.
Далее автором предлагается новый модифицированный метод,
позволяющий скомбинировать оба эти подхода.

Среди многих работ по получению  представлений решений для
возмущённого уравнения Бесселя~\eqref{4.2.1}-\eqref{4.2.2} отметим
те, в которых решение ищется в виде рядов специального вида "---
это работы А.~Фитоухи с соавторами~\cite{CFH, FH} и
В.\,В.~Кравченко с соавторами~\cite{CKT1,CKT2, Krav1, Krav2,
Krav3, Krav4, Krav5, Krav6, Krav7, Krav8, Krav9, Krav10}.
Критический разбор ряда результатов по этой задаче недавно
проведён в~\cite{Hol}.

Вместе с тем во многих математических и физических задачах
необходимо рассматривать сильно сингулярные потенциалы, например,
допускающие произвольную степенную особенность в нуле. В настоящей
работе сформулированы результаты по интегральному представлению
решений уравнений с подобными сингулярными потенциалами. От
потенциала требуется лишь мажорируемость определенной функцией,
суммируемой на бесконечности. В частности, к классу допустимых в
данной работе относятся  сингулярный потенциал $q=x^{-2},$ сильно
сингулярный потенциал со степенной особенностью
$q=x^{-2-\varepsilon},$ $\varepsilon > 0,$ потенциалы Юкавы типа
$q=e^{-\alpha x}/x,$ потенциалы Баргмана и
Батмана---Шадана~\cite{ShSa} и ряд других. При этом на функцию
$q(x)$ не накладывается никаких дополнительных условий типа
быстрой осцилляции в начале координат или знакопостоянства, что
позволяет изучать притягивающие и отталкивающие потенциалы единым
методом.

Следует отметить, что в данной книге построены операторы
преобразования специального вида, отличающиеся от ранее известных
некоторыми деталями. До этого рассматривались лишь случаи
одинаковых пределов (оба вида $[0;a]$ или $[a;\infty]$) в основном
интегральном уравнении для ядра оператора преобразования. В данной
работе  показано, что можно рассматривать случай различных
пределов в основном интегральном уравнении. Именно такая
расстановка пределов и позволила охватить более широкий класс
потенциалов с особенностями в нуле. Кроме того, по сравнению с
рассуждениями по образцу классической  работы
Б.\,М.~Левитана~\cite{Lev7} мы вносим  усовершенствование в эту
схему. Используемую в доказательствах функцию Грина как оказалось
можно выразить не только через общую гипергеометрическую функцию
Гаусса, но и более конкретно через функцию Лежандра, зависящую от
меньшего числа параметров, что позволяет избавиться от
неопределённых постоянных в оценках из предыдущих работ.

Из-за ограниченности объёма книги в этом пункте приводятся только
постановка задачи, сводка основных результатов и следствий без
доказательств, подробное изложение см.
в~\cite{S8,S63,S46,S19,S43,S4}.

\subsection{Решение основного интегрального уравнения для ядра оператора
преобразования}\label{sec17.1}

Введем новые переменные и функции по формулам:
$$
\xi=\frac{t+x}{2}, ~ \eta=\frac{t-x}{2}, ~ \xi \geq \eta > 0;
$$
\begin{equation}\label{4.2.5}
K(x, t)= \left(\frac{x}{t}\right)^\alpha P(x, t), ~ u(\xi, \eta)= K(\xi-\eta, \xi+\eta).
\end{equation}
Обозначим $\nu=\alpha-1.$ Таким образом, для обоснования
представления~\eqref{4.2.3*} для решения  уравнения~\eqref{4.2.1}
достаточно определить функцию $u(\xi, \eta ).$
Известно~\cite{Soh1, Soh2, Soh3, Soh4}, что если существует дважды
непрерывно дифференцируемое решение $u(\xi, \eta)$ интегрального
уравнения
$$
u(\xi, \eta)=-\frac{1}{2}\int\limits_\xi^{\infty}R_\nu(s, 0; \xi, \eta) q(s) \, ds -
\int\limits_\xi^{\infty} ds \int\limits_0^{\eta} q(s+\tau) R_{\nu}(s, \tau; \xi, \eta) u(s, \tau) \, d \tau,
$$
при условиях $ 0< \tau < \eta < \xi < s,$ то искомая функция $P(x,
t)$ определяется по формулам~\eqref{4.2.5} через это решение
$u(\xi, \eta).$ Функция $R_{\nu}=R_{\alpha-1}$ является функцией
Римана, возникающей при решении некоторой задачи Гурса для
сингулярного гиперболического уравнения
$$
\frac{\partial^2 u(\xi, \eta)}{\partial \xi \partial \eta}+ \frac{4 \alpha(\alpha-1) \xi \eta}{(\xi^2-\eta^2)^2} u(\xi, \eta)=q(\xi+\eta)u(\xi, \eta).
$$
Эта функция известна в явном виде, см.~\cite{Soh1, Soh2, Soh3,
Soh4}, она выражается через гипергеометрическую функцию Гаусса
$_2{F_1}$ по формуле
\begin{equation}\label{4.2.6}
R_\nu=\left(\frac{s^2-\eta^2}{s^2-\tau^2}\cdot \frac{\xi^2-\tau^2}{\xi^2-\eta^2}\right)^{\nu}{_2F_1} \left(-\nu, -\nu; 1; \frac{s^2-\xi^2}{s^2-\eta^2}\cdot \frac{\eta^2-\tau^2}{\xi^2-\tau^2}\right).
\end{equation}
Это выражение упрощено в~\cite{S8}, где показано, что функция
Римана в рассматриваемом случае выражается через функцию Лежандра
по формуле
\begin{equation}\label{4.2.7}
R_\nu (s, \tau, \xi, \eta)=P_\nu \left(\frac{1+A}{1-A}\right),~A=\frac{\eta^2-\tau^2}{\xi^2-\tau^2}\cdot \frac{s^2-\xi^2}{s^2-\eta^2}.
\end{equation}

Основное содержание этого пункта составляет

\begin{theorem}  Пусть функция $q(r)\in C^1 (0,\infty)$  удовлетворяет условию
    \begin{equation}\label{4.2.8}
    |q(s+\tau)|\leq |p(s)|, ~ \forall s, \forall \tau, ~ 0< \tau <s,\ \int\limits_\xi^\infty |p(t)| \, dt < \infty, \forall \xi>0.
    \end{equation}
    Тогда существует интегральное представление вида~\eqref{4.2.3*}, ядро которого удовлетворяет оценке
    $$
    \gathered
    |P(r, t)| \leq \left(\frac{t}{r}\right)^ \alpha \frac {1}{2} \int\limits_{\frac{t+r}{2}}^\infty P_{\alpha-1}\left(\frac{y^2(t^2+r^2)-(t^2-r^2)}{2try^2}\right)|p(y)|\, dy \times
    \\
    \times\exp \left[ \left(\frac{t-r}{2}\right) \frac{1}{2} \int\limits_{\frac{t+r}{2}}^\infty P_{\alpha-1}\left(\frac{y^2(t^2+r^2)-(t^2-r^2)}{2try^2}\right)|p(y)|\, dy  \right].
    \endgathered
    $$
    При этом ядро оператора преобразования $P(x,t),$ а также  решение уравнения~\eqref{4.2.1} являются дважды непрерывно дифференцируемыми на  $(0,\infty)$ функциями по своим аргументам.
\end{theorem}

Перечислим классы потенциалов, для которых выполнены
условия~\eqref{4.2.8}. Если $|q(x)|$ монотонно убывает,  то можно
принять $p(x)=|q(x)|.$ Для потенциалов с произвольной особенностью
в начале координат и возрастающих при $0<x<M$ (например,
кулоновских $q=-\dfrac{1}{x}$), которые обрезаны нулём на
бесконечности, $q(x)=0,$ $x>M,$ можно принять $p(x)=|q(M)|,$
$x<M,$ $p(x)=0,$ $x \geq M.$ Условию~\eqref{4.2.8}  будут также
удовлетворять потенциалы с оценкой $q(x+\tau) \leq
c|q(x)|=|p(x)|.$ На возможность подобного усиления
теоремы~\ref{4.2.1} указал В.\,В.~Катрахов.

В частности, приведенным условиям удовлетворяют следующие
потенциалы, встречающиеся в приложениях: сильно сингулярный
потенциал со степенной особенностью вида
$q(x)=x^{-2-\varepsilon},$ различные потенциалы Баргмана
$$
q_1 (x) = - \frac{e^{-ax}}{(1+\beta e^{-ax})^2}, ~ q_2 (x)= \frac{c_2}{(1+c_3 x)^2}, ~ q_3 (x)= \frac{c_4}{ch^2(c_5 x)},
$$
и Юкавы
$$
q_4 (x) = - \frac{e^{-ax}}{x},~ q_5 (x) = \int\limits_x^\infty e^{-at} \, d c(t).
$$
(см., например,~\cite{ShSa}).

\begin{remark}
Фактически при доказательстве приведенной теоремы  не требуется
явный вид функции Римана~\eqref{4.2.7}.  Используется только
существование функции Римана, её положительность и некоторое
специальное свойство монотонности. Эти факты являются довольно
общими, поэтому полученные результаты можно распространить на
достаточно широкий класс дифференциальных уравнений.
\end{remark}

Оценку из теоремы~\ref{4.2.1} для потенциалов общего вида можно
преобразовать в более грубую, но зато и более обозримую.

\begin{theorem}
    Пусть выполнены условия теоремы~\ref{4.2.1}. Тогда для ядра оператора преобразования $P(x, t)$ справедлива оценка
    $$
    |P(x, t)| \leq \frac{1}{2} \left(\frac{t}{x}\right)^{\alpha} P_{\alpha-1} \left(\frac{t^2+x^2}{2tx}\right) \int\limits_x^\infty |p(y)|\, dy\,  \exp \left[ \frac{1}{2} \left(\frac{t-x}{2}\right) P_{\alpha-1} \left(\frac{t^2+x^2}{2tx}\right) \int\limits_x^\infty |p(y)| \, dy \right].
    $$
\end{theorem}

Отметим, что при $x \to 0$ ядро интегрального представления может иметь экспоненциальную особенность.

Для класса потенциалов со степенной сингулярностью вида
\begin{equation}\label{4.2.16}
q(x)=x^{-(2\beta +1 )},~ \beta > 0
\end{equation}
полученные оценки можно упростить, не снижая их точности. Ограничение на $\beta$ вызвано условием суммируемости на бесконечности.

\begin{theorem} Рассмотрим потенциал вида~\eqref{4.2.16}. Тогда теорема~\ref{4.2.1} выполняется с оценкой
    $$
    |P(x, t)|  \leq \left(\frac{t}{x}\right)^{\alpha}  \frac{\Gamma(\beta)4^{\beta-1}}{(t^2-x^2)^\beta} \cdot P_{\alpha-1}^{- \beta} \left(\frac{t^2+x^2}{2tx}\right) \exp\left[\left(\frac{t-x}{x}\right)  \frac{\Gamma(\beta)4^{\beta-1}}{(t^2-x^2)^\beta} P_{\alpha-1}^{- \beta} \left(\frac{t^2+x^2}{2tx}\right)   \right],
    $$
    где $P_\nu^\mu(\cdot)$ "--- функция Лежандра, величина $\beta$ определяется из~\eqref{4.2.16}, а величина $\alpha$ "--- из~\eqref{4.2.2}.
\end{theorem}

Отметим, что данная оценка получается после довольно длинных
вычислений с использованием знаменитой теоремы
Слейтер---Маричева~\cite{Marich1}, которая помогает вычислить в
терминах гипергеометрических функций необходимые интегралы после
их сведения к свёртке Меллина.

Простейшая подобная оценка была получена в работе~\cite{S8} для
потенциала $q(x)=cx^{-2},$ для которого $\beta=\dfrac{1}{2}.$ Как
следует из~\cite{BE1}, в этом случае функция Лежандра
$P_{\nu}^{-\frac{1}{2}}(z)$ может быть выражена через элементарные
функции. Поэтому и соответствующая оценка   может быть выражена
через элементарные функции.

Другим потенциалом,  для которого полученная  оценка   может быть
ещё упрощена и выражена через элементарные функции, является
потенциал вида $q(x)=x^{-(2 \beta + 1)},$ когда параметры связаны
соотношением $\beta=\alpha-1.$

\begin{corollary} Пусть выполнено соотношение между параметрами $\beta=\alpha-1.$ Тогда оценка из теоремы~\ref{4.2.3} принимает вид
    $$
    |P(x, t)|  \leq \left(\frac{t}{x}\right)^{\beta+1} \frac{2^{\beta-2}}{\beta} \left[\frac{t^2+x^2}{2tx}\right]^\beta  \exp \left[ \left(\frac{t-x}{2}\right) \frac{2^{\beta-2}}{\beta} \left[\frac{t^2+x^2}{2tr}\right]^\beta \right]=
    $$
    \begin{equation}\label{4.2.20}
    = \frac{1}{4 \beta} \frac{1}{x^{2 \beta +1}} (t^2+x^2)^{\beta} \exp \left[  \frac{2^{\beta-2}}{\beta} \left(\frac{t-x}{2}\right) \left(\frac{t^2+x^2}{2tx}\right)^\beta \right].
    \end{equation}
\end{corollary}

Отметим, что при $\alpha=0$ в
формулах~\eqref{4.2.1}--\eqref{4.2.3*}, теорема~\ref{4.2.1}
сводится к известным оценкам для ядра интегрального представления
решений Йоста для уравнения Штурма---Лиувилля.

Изложенная техника полностью переносится и на задачу о построении
неклассических операторов обобщённого сдвига. Эта задача по
существу эквивалентна выражению решений уравнения
\begin{equation}\label{4.2.22}
B_{\alpha,x} u(x,y) - q(x) u(x,y)=B_{\beta,y} u(x,y)
\end{equation}
через решения невозмущённого  уравнения Эйлера---Пуассона---Дарбу
(в несингулярном случае "--- волнового) при наличии дополнительных
условий, обеспечивающих корректность. Такие представления
получаются уже из факта существования операторов преобразования и
изучались для несингулярного случая ($\alpha=\beta=0$)
в~\cite{Lev1, Lev2, Lev3} как следствия теории обобщённого сдвига,
см. также~\cite{Mar9}. Интересная оригинальная методика для
получения подобных представлений также в несингулярном случае
разработана в работах А.\,В.~Боровских~\cite{Bor1,Bor2}. Из
результатов настоящей работы следуют интегральные представления
некоторого подкласса решений уравнения~\eqref{4.2.22} в общем
сингулярном случае для достаточно произвольных потенциалов с
особенностями в начале координат. При этом оценки для решений не
содержат никаких неопределенных постоянных, а для ядер
интегральных представлений в явном виде выписываются интегральные
уравнения, которым они удовлетворяют.

\section{Приложение метода операторов преобразования к задаче
Е.\,М.~Ландиса}\label{sec18}

В заметке Е.\,М.~Ландиса~\cite{Lan} поставлена следующая задача:
доказать, что решение стационарного уравнения Шрёдингера с
ограниченным потенциалом
\begin{eqnarray}
& &  \Delta u(x) - q(x) u(x) =0,~x \in \R^n,~|x| \geq R_0 > 0,   \label{l1} \\
& & |q(x)| \leq \lambda^2,~ \lambda > 0,~u(x) \in C^2 \lr{|x| \geq
R_0}, \nonumber
\end{eqnarray}
удовлетворяющее оценке вида
$$
|u(x)| \leq \const\cdot \, e^{- (\lambda + \varepsilon)|x|},~ \varepsilon > 0,
$$
тождественно равно нулю.

В.\,З.~Мешковым в известных работах~\cite{Mesh1,Mesh2} был дан
отрицательный ответ на данный вопрос. При этом было доказано
существование контрпримеров с решениями, которые являются
комплексными функциями. Более того, было показано, что если
усилить оценку в гипотезе Е.\,М.~Ландиса до следующей:
$$
|u(x)| \leq \const\cdot \, e^{- (\lambda + \varepsilon)|x|^{4/3}},~ \varepsilon > 0,
$$
то ответ будет положительным "--- таких ненулевых решений не
существует. В последнее время интерес к этим результатам не
пропал, тематика, связанная с гипотезой Е.\,М.~Ландиса и
результатами В.\,З.~Мешкова, активно развивается, причём в том
числе и ведущими математиками в области дифференциальных уравнений
"--- Ж.~Бургейном, К.~Кёнигом и рядом других,
см.~\cite{BuKe1,Ken1,Ken2,DKW1,Rossi}. Основным вопросом остаётся
исследование гипотезы Е.\,М.~Ландиса для действительных решений,
причём ответ на этот вопрос до сих пор не удалось получить. В
связи с вышеизложенным представляется обоснованным название
\textit{задача Ландиса---Мешкова} в следующей формулировке.

\medskip

\textbf{Задача Ландиса---Мешкова.}
{\it Верно ли, что для заданных области $D$ и положительных функций
$r(x), s(x)$ только нулевое классическое решение стационарного
уравнения Шрёдингера
\begin{equation}\label{LM1}
\Delta u(x) - q(x) u(x) =0,\ \ x \in \mathbb{R}^n, \ \ |q(x)|\leq r(x),
\end{equation}
удовлетворяет оценке}
\begin{equation}\label{LM2}
|u(x)| \leq s(x)?
\end{equation}

\medskip

Из результатов В.\,З.~Мешкова следует отрицательный ответ в этой
задаче в случае комплексных решений,
 $D$ "--- внешность некоторого круга, $q(x)=\lambda^2,$ $s(x)=e^{- (\lambda + \varepsilon)|x|},$ $\varepsilon > 0,$ и положительный ответ в этой задаче в случае комплексных решений,
 $D$ "--- внешность некоторого круга, $q(x)=\lambda^2,$
 $s(x)=e^{- (\lambda + \varepsilon)|x|^{4/3}},\ \varepsilon > 0.$ Для действительных решений даже в этих частных случаях ответы неизвестны.

Далее мы показываем, что несмотря на общее отрицательной решение
В.\,З.~Мешкова для первоначальной формулировки задачи
Е.\,М.~Ландиса, тем не менее для некоторых классов потенциалов
проблема решается положительно, причём для действительных решений.
При этом используется метод операторов преобразования специального
вида~\cite{S3, S71, S75}.

Далее эта задача решена для случая потенциала, зависящего только
от одной переменной: $q(x)= q(x_i),$ где $1 \leq i \leq n;$ далее
для определённости считается, что $i=1.$ Для этого случая в работе
доказано обобщение утверждения~\eqref{l1} для уравнения
\begin{equation}\label{l2}{\Delta u - q(x_1) u = 0,}\end{equation} в котором потенциал $q(x_1)$
ограничен произвольной неубывающей функцией. Решение основано на
использовании операторов преобразования, сводящих
уравнение~\eqref{l2} к уравнению Лапласа.

\medskip

{\bf 1.} Условия задачи~\eqref{l1}  выполнены в полупространстве
$x_1 \geq R_0$ и инвариантны относительно замены переменных $z=x_1
- R_0.$ Поэтому мы будем рассматривать задачу~\eqref{l1} в
полупространстве $z \geq 0$ или, сохраняя для переменной $x_1$
прежнее обозначение, $x_1 \geq 0.$ Будет доказано, что решение
задачи~\eqref{l1} равно нулю в полупространстве $x_1 \geq 0,$ а
тогда в силу теоремы Кальдерона о единственности продолжения
(см.~\cite[гл.~6, с.~14]{Miran}) такое решение тождественно равно
нулю во всём пространстве  $\R^n.$

Обозначим через $T\lr{\delta}$ множество функций, удовлетворяющих
в полупространстве\\
$\R^n_+=\{ x \in \R^n,$ $x_1 \geq 0 \}$ следующим
условиям~\eqref{l3}--\eqref{l5}:
\begin{eqnarray}
& & u(x) \in C^2 \lr{\R^n_+}, \label{l3} \\
& & \left| u(x)\right| \leq c_1 \, e^{- \delta |x|},~\delta>0, \label{l4} \\
& &  \left| \frac{\pr u }{\pr x_1}\right| \leq c_2 \, e^{- \delta
|x|}.\label{l5}
\end{eqnarray}
Построим для функций из $T \lr{\lambda+\varepsilon}$ оператор
преобразования вида (см.~\cite{S3, S71, S75})
\begin{equation}\label{l6}{ S u(x) = u(x) +
\int\limits_{x_1}^{\infty} K (x_1, t) u (t,
x^1)\,dt,}\end{equation} чтобы выполнялось равенство
\begin{equation}\label{l7}{S \lr{\frac{\pr^2 u}{\pr x_1^2} -
q(x_1) u} = \frac{\pr^2 }{\pr x_1^2} S u,~|q(x_1)|\leq
\lambda^2,}\end{equation} где, как обычно, через $(x_1, x^1)$
обозначено $(x_1, x_2, \dots, x_n).$ Подстановка
выражения~\eqref{l6} в формулу~\eqref{l7} приводит к равенствам
\begin{equation}\label{l8}{ \frac{\pr^2 K}{\pr t^2} - \frac{\pr^2
K}{\pr x_1^2} = q(t) K,}\end{equation}
\begin{equation}\label{l9}{3 \frac{\pr K(x_1, x_1)}{\pr x_1} =
q(x_1),}\end{equation} \begin{equation}\label{l10}{\lim\limits_{t
\to \infty} K(x_1, t) \frac{\pr u (t, x^1)}{\pr t} -
\lim\limits_{t \to \infty} \frac{\pr  K(x_1, t)}{\pr t} u (t, x^1)
= 0.}\end{equation} Выполняя стандартную замену переменных
$w=\dfrac{t+x_1}{2},~v=\frac{t-x_1}{2},$ мы сводим
систему~\eqref{l8}-\eqref{l9} к более простой (выполнение
условия~\eqref{l10} на решениях задачи~\eqref{l1} будет показано
позже)
\begin{eqnarray}
& & \frac{\pr^2 K}{\pr w \pr v} = q (w+v) K, \label{l11} \\
& & K(w, 0) =  \frac{1}{3} \int\limits_0^w q(s) \, ds, \label{l12}
\end{eqnarray}
которая, в свою очередь, является следствием одного интегрального уравнения
\begin{equation}
K(w,v) = \frac{1}{3} \int\limits_0^w q(s) \, ds +
 \int\limits_0^w d\alpha  \int\limits_0^v q(\alpha+\beta) K (\alpha, \beta) \, d \beta,\quad |q| \leq \lambda^2,\quad w \geq v \geq 0. \label{l13}
\end{equation}

Уравнение~\eqref{l13} отличается от обычно используемого при
рассмотрении операторов преобразования на бесконечном интервале
интегрального уравнения изменением области интегрирования с
полуоси $\lr{w, \infty}$ на отрезок $\lr{0, w},$ что влечёт
экспоненциальный рост ядра $K \lr{x_1, t}.$ Далее доказывается,
что такое ядро существует и оператор преобразования с таким
ядром~\eqref{l6} определён на множестве $T
\lr{\lambda+\varepsilon}.$ Возможность сведения
задачи~\eqref{l8}--\eqref{l10} к неэквивалентным интегральным
уравнениям вытекает из недоопределённости задачи
Коши~\eqref{l11}-\eqref{l12}.

Следует отметить, что интегральное уравнение~\eqref{l13} должно
быть решено в более широкой области без ограничений $w \geq v,$
иначе не будет определено ядро под знаками интегралов.
Доказательство существования решения в этой более широкой области
проводится так же, как приведённое ниже доказательство. На этот
нюанс при доказательстве существования решения интегрального
уравнения~\eqref{l13} обычно не обращают внимания (замечание
А.\,В.~Боровских).

\begin{lemma}\label{L18.1} Существует единственное непрерывное решение уравнения~\eqref{l13}, удовлетворяющее неравенству
\begin{equation}\label{l14}{|K (w, v)| \leq \frac{\lambda}{3}  \sqrt{\frac{w}{v}}\,
I_1 \lr{2\lambda \, \sqrt{wv}},}\end{equation} где $I_1 (x)$ "---
модифицированная функция Бесселя. При этом на допустимом
потенциале $q(x_1) \equiv \lambda^2$ в~\eqref{l14} достигается
знак равенства.
\end{lemma}

\begin{remark} В дальнейшем символом $c$ обозначаются абсолютные положительные постоянные, величина которых не играет роли.
\end{remark}

\begin{proof}
Введём обозначения
$$
K_0 (w, v) = \frac{1}{3} \int\limits_0^w q(s) \, ds,
$$
$$
P K (w, v) = \int\limits_0^w d\alpha  \int\limits_0^v
q(\alpha+\beta) K (\alpha+\beta) \, d \beta.
$$
Тогда уравнение~\eqref{l13} запишется в виде $K = K_0 +P K.$ Будем
искать его решение в виде ряда Неймана
\begin{equation}\label{l15}{K = K_0 +P K_0 + P^2 K_0 + \dots }\end{equation} Для
слагаемых ряда~\eqref{l15} с учётом условия $|q(x_1)| \leq
\lambda^2$ получаем \begin{equation}\label{l16}{\left| P^n K_0
(w_0 v) \right| \leq \frac{1}{3} \,  \lr{\lambda^2}^{n+1}
\frac{w^{n+1}}{(n+1)!} \frac{v^n}{n!},~n=0, 1, 2,
\dots}\end{equation} Отсюда вытекает неравенство~\eqref{l14}, если
использовать представление функции $I_1 (x)$ в виде ряда
$$
I_1 (x) = \sum\limits_{k=0}^{\infty} \frac{(x/2)^{2 k+1}}{k!
(k+1)!}.
$$
Оценка~\eqref{l14} является точной, так как при $q(x_1) \equiv
\lambda^2,$  неравенства~\eqref{l16} превращаются в равенства для
всех целых $n \geq 0.$ Лемма доказана.
\end{proof}

\begin{lemma}\label{L18.2} В терминах переменных  $x_1,$ $t$ справедлива оценка
$$
\left| K(x_1, t) \right| \leq c \, t \,e^{\lambda t }.
$$
\end{lemma}

\begin{proof}
Рассмотрим неравенство
$$
\left| \frac{1}{x} I_1 (x) \right| \leq c \, e^x,~x \geq 0,
$$
для проверки истинности которого надо разобрать случаи а) $x \geq
1,$ б) $0 \leq x  \leq 1$ и использовать известную асимптотику
функции $I_1 (x)$ при $x \to \infty$ и $x \to +0$
(см.~\cite{BE2}). Отсюда с помощью очевидных неравенств
$$
 \frac{x_1+t}{2} \leq t, ~ 2 \sqrt{w v} = \sqrt{t^2 -x_1^2} \leq t
$$
из оценки~\eqref{l14} следует утверждение леммы.
\end{proof}

Из леммы  следует, что выражение~\eqref{l6} определено на функциях
из $T \lr{\lambda+\varepsilon}.$ Покажем, что выражение~\eqref{l6}
в действительности задаёт оператор преобразования на  $T
\lr{\lambda+\varepsilon}.$ Для этого осталось проверить
соотношение~\eqref{l10}. Из того, что $u(x) \in T
\lr{\lambda+\varepsilon}$  и из леммы  вытекает равенство
$$
\lim\limits_{t \to \infty} K (x_1, t) \frac{\pr u(t, x^1)}{\pr t} = 0.
$$
Поэтому осталось доказать, что если $u(x) \in T
\lr{\lambda+\varepsilon},$ то
$$
\lim\limits_{t \to \infty}  \frac{\pr K (x_1, t)}{\pr t}  u(x_1, t) = 0.
$$
Последнее соотношение следует из оценки
\begin{equation}\label{l17}{\left|  \frac{\pr K (x_1, t)}{\pr t}
\right| \leq c\, t\, e^{\lambda t }. }\end{equation} Для
доказательства неравенства~\eqref{l17} нужно перейти к переменным
$w,$ $v$ и с использованием уже установленных оценок для ядра
$K(x_1, t)$ оценить производные $\dfrac{\pr K}{\pr w},$
$\dfrac{\pr K}{\pr v},$дифференцируя уравнение~\eqref{l13}. Так
как
$$
\frac{\pr K}{\pr t} = \frac{1}{2} \lr{\frac{\pr K}{\pr
w}+\frac{\pr K}{\pr v}},
$$
то мы придём к~\eqref{l17}.

\medskip

{\bf 2.} Покажем, что любое решение задачи~\eqref{l1} принадлежит
$T \lr{\lambda+\varepsilon}$ и, следовательно, на таких решениях
определён оператор~\eqref{l6}. Для этого необходимо проверить
выполнение условия~\eqref{l5}.

\begin{lemma}\label{L18.3} Пусть функция $u(x) \in C^2 \lr{|x| \geq R_0}$ есть решение задачи~\eqref{l1}. Тогда найдётся такая постоянная $c>0,$ что
$$
\left| \frac{\pr u}{\pr x_1} \right| \leq c \, e^{-(\lambda+\varepsilon)|x|}.
$$
\end{lemma}

\begin{proof}
В силу априорных оценок Шаудера, в замкнутом шаре $B(x, 1)$
единичного радиуса с центром в точке $x,$ $|x| \geq R_0+1,$ имеем
(см.~\cite[теорема~33, II]{Miz})
$$
u_1 \leq c \lr{{u_{1, \lambda_1}}^{\frac{1}{1+\lambda_1}} \cdot
{u_{0}}^{\frac{\lambda_1}{\lambda_1+1}} + u_0 },
$$
где обозначено
$$
u_0= \| u(x) \|_{C^0 \lr{B(x, 1)}},~u_1 = \| u(x) \|_{C^1 \lr{B(x,
1)}};
$$
$u_{1, \lambda_1}$ есть сумма коэффициентов Гёльдера функции
$u(x)$ и её производных первого порядка $\dfrac{\pr u}{\pr x_i},$
$1 \leq  i \leq n.$ Отсюда следует, что
\begin{equation}\label{l18}{\left|\frac{\pr u(x)}{\pr x_1}\right|
\leq c \lr{{u_{1, \lambda_1}}^{\frac{1}{1+\lambda_1}} \cdot
{u_{0}}^{\frac{\lambda_1}{\lambda_1+1}} + u_0 }.}\end{equation}
Отметим, что так как выполнены все условия~\cite[утверждение~33,
V]{Miz}, то константа $c$ в формуле~\eqref{l18} не зависит от $x.$

Из результатов Морри (см.~\cite[теорема~39, IV]{Miran}) вытекает
оценка для величины  $u_{1, \lambda_1}$
\begin{equation}\label{l19}{ u_{1, \lambda_1} \leq c \left[ \| u
\|_{L_2 \lr{B(x, 1)}}+ \|q\, u \|_{L_2 \lr{B(x,
1)}}\right],}\end{equation} причём постоянная в~\eqref{l19}
по-прежнему не зависит от $x.$ Из условий задачи $|q(x_1)| \leq
\lambda^2,$ следовательно, с помощью теоремы о  среднем получаем
из~\eqref{l19}
$$
u_{1, \lambda_1}  \leq c \lr{\int\limits_{B(x, 1)} |u(y)|^2 \,
dy}^{1/2} \leq {c^1}\, e^{- \lr{\lambda+\varepsilon}|x|}.
$$
Подставляя последнее неравенство  в~\eqref{l18}, получаем
$$
\left| \frac{\pr u}{\pr x_1} \right| \leq c \left[ \lr{ e^{-
\lr{\lambda+\varepsilon}|x|}}^{\frac{1}{1+\lambda}+\frac{\lambda}{1+\lambda}}
+
 e^{- \lr{\lambda+\varepsilon}|x|}\right] \leq c \, e^{- \lr{\lambda+\varepsilon}|x|}.
$$

Таким образом, требуемое неравенство установлено для $|x_1| \geq
R_0 +1.$ Так как множество $R_0 \leq |x| \leq R_{0}+1$ является
компактом в $\R^n,$ то это неравенство справедливо и при $|x| \geq
R_0.$ Лемма~\ref{L18.3} доказана.
\end{proof}

Выполняя опять замену координат $z=x_1-R_0,$ получаем, что
лемма~\ref{L18.3} справедлива в полупространстве $x_1\geq 0$ (мы
переобозначим $z$ через $x_1$).

\medskip

{\bf 3.} Применим к уравнению~\eqref{l2} оператор $S.$ Из
тождества~\eqref{l7} и перестановочности $S$ с производными
$\dfrac{\pr^2 u}{\pr x_i^2},$ $2 \leq i \leq n,$ получаем, что в
полупространстве $\R_{+}^n$
$$
S \lr{\Delta u - q(x_1) u} = \Delta S u = 0.
$$

Обозначим функцию $Su$ через $v.$ Из~\eqref{l6},~\eqref{l13}
следует, что если $u(x) \in C^2 \lr{\R_{+}^n},$ $q(x) \in C
\lr{\R_{+}^n},$ то $v(x) \in C^2 \lr{\R_{+}^n}.$ Покажем, что
$v(x)$ экспоненциально убывает в $\R_{+}^n$ при $|x| \to \infty$
и, следовательно, равна нулю.

\begin{lemma} Пусть $u(x) \in T \lr{\lambda+ \varepsilon}.$ Тогда для $x \in \R_{+}^n$
$$
|v|= \left| S u \right|  \leq c\, |x| \, e^{-\varepsilon |x|},~\varepsilon>0.
$$
\end{lemma}

\begin{proof}
Из представления~\eqref{l6} и леммы~\ref{L18.3} получаем
$$
\left| S u \right|  \leq |u(x)| + \int\limits_{x_1}^{\infty} t \,
e^{\lambda t} c\, e^{-(\lambda+\varepsilon) \sqrt{t^2+|x^1|^2}} \,
dt \leq
 c \lr{ e^{-(\lambda + \varepsilon) |x|} +
\int\limits_{x_1}^{\infty} t \, e^{-(\lambda+\varepsilon)
\sqrt{t^2+|x^1|^2}} \, dt}.
$$
Вычисляя интеграл с помощью замены переменных по формуле $y =
\sqrt{t^2 + |x^1|^2}$ с последующим интегрированием по частям,
получаем требуемую оценку. Лемма доказана.
\end{proof}

Итак, $v(x)=0$ в $\R_{+}^n.$ Определим на $T \lr{\lambda+
\varepsilon}$ обратный к $S$ оператор $P$ по формуле
$$
P u(x) = u(x)+ \int\limits_{x_1}^{\infty} N(x_1, t) u(t, x^1) \, dt.
$$
Тогда для ядра $N \lr{x_1, t}$ справедливы утверждения
лемм~\ref{L18.1}--\ref{L18.3}. Кроме того, если $S u \in T
\lr{\lambda+ \varepsilon},$ то \begin{equation}\label{l20}{P S
u(x) = u(x).}\end{equation}

Так как, очевидно, $0 \in T \lr{\lambda+ \varepsilon},$ то
применяя~\eqref{l20} к обеим частям установленного в $\R_{+}^n$
равенства $S u = 0,$ получим $u=0$ в $\R_{+}^n.$ Выше было
показано, что это влечёт $u \equiv 0$ во всём $\R^n.$

\begin{remark} Переход к полупространству $\R_{+}^n$ использовался при доказательстве потому, что выражение~\eqref{l6} не определено в области, получаемой пересечением шара $|x|\leq R_0$ и бесконечного полуцилиндра $\{|x^1| \leq R_0,~|x_1| \leq R_0 \}.$
\end{remark}

Итак, доказана

\begin{theorem} \label{Theo1} Любое решение  $u(x) \in C^2 \lr{|x|>R_0}$ стационарного уравнения Шрёдингера с ограниченным потенциалом
$$
\Delta u(x) - q(x_1) u = 0, ~x \in \R^n,~ |x| \geq R_0 > 0,
$$
$$
q(x_1) \in C \lr{|x| \geq R_0},~|q(x_1) | \leq \lambda^2,~\lambda>0,
$$
удовлетворяющее  оценке
$$
|u(x)| \leq \const e^{- \lr{\lambda+\varepsilon}|x|},~\varepsilon>0,
$$
есть тождественный нуль.
\end{theorem}

{\bf 4.} Использованная техника операторов преобразования
позволяет усилить полученный результат. Будем обозначать через
$L_{2,\, loc} \lr{x_1 \geq R_0}$ множество функций, для которых
при любом $x_1 \geq R_0$ конечен интеграл $\int\limits_{R_0}^{x_1}
\psi^2(s)\,ds.$ Пусть далее, задана неотрицательная функция
$g(x),$ для которой интеграл $\int\limits_{x_1}^{\infty} t\, g(t,
x^1) \, dt = p(x)$ конечен при любом $x_1 \geq R_0$ и для
некоторой постоянной $\alpha > 0$
$$
|p(x)| \leq c \cdot \exp \lr{- \alpha |x|^{\delta}},~\delta>0.
$$
Тогда по схеме доказательства предыдущей теоремы может быть установлена

 \begin{theorem} \label{Theo2} Пусть $\psi(x_1) \in L_{2,\, loc} \lr{x_1 \geq R_0},$  $\psi(x_1)$ "--- неубывающая функция, функция $g(x)$ удовлетворяет перечисленным выше требованиям. Тогда любое решение уравнения
$$
\Delta u(x) - q(x_1) u = 0, ~x \in \R^n,~ |x| \geq R_0 > 0,
$$
$$
|q(x_1)| \leq \psi^2 (x_1),
$$
для которого выполнено неравенство
$$
\psi (x_1) |u(x)| \leq \const e^{-\psi (x_1) |x|} g(x),~g(x) \geq 0,
$$
есть тождественный нуль.
\end{theorem}

В условиях теоремы~\ref{Theo1} нужно положить $g(x) =
e^{-\varepsilon |x|}.$ Примером другой допустимой $g(x)$ является
функция $g(x)= \exp \lr{-\varepsilon |x|^{\delta} },$ $0< \delta
<1.$ Этот случай является также примером обобщённой задачи
Ландиса---Мешкова~\eqref{LM1}-\eqref{LM2}.

По аналогичной схеме может быть также рассмотрен  случай
потенциала, зависящего только от радиальной переменной. Ответ в
первоначальной формулировке задачи Е.\,М.~Ландиса здесь тоже
положительный, после перехода к сферическим координатам нужно
использовать операторы преобразования  для возмущённого оператора
Бесселя, подобные тем, что были рассмотрены в предыдущем пункте
этой книги, см.~\cite{S3, S71, S75}.

Возможно рассмотрение обобщений задачи Е.\,М.~Ландиса на случай
более общих дифференциальных уравнений и соответствующих оценок
роста решений. Например, представляет интерес исследование
поставленных вопросов для нелинейного уравнения
$p$-Лапласиана~\cite{Lind, DKN} (эта задача возникла во время
обсуждения доклада автора на семинаре кафедры дифференциальных
уравнений МГУ с проф. А.\,А.~Коньковым).


\chapter*{Благодарности}
\addcontentsline{toc}{chapter}{Благодарности}
\sectionmark{\sc Благодарности}

В заключение С.\,М.~Ситник благодарит своих коллег, которые
просмотрели текст книги и сделали ряд полезных замечаний,
исправлений и дополнений: А.\,В.~Глушака, Д.\,Б.~Карпа,
В.\,В.~Кравченко, А.\,Б.~Муравника, Э.\,Л.~Шишкину. Автор
благодарит Е.\,М.~Варфоломеева за большой труд по редактированию
этой книги. Также выражаю благодарность С.\,Н.~Ушакову  за
квалифицированный компьютерный набор рукописи диссертации
В.\,В.~Катрахова.

\newpage

\chapter*{Валерий Вячеславович Катрахов\\(1949--2010):\\краткая биографическая справка}
\addcontentsline{toc}{chapter}{Валерий Вячеславович
Катрахов: краткая биографическая справка}
\sectionmark{\sc Валерий Вячеславович Катрахов (1949--2010): краткая биографическая справка}

\epigraph{<<Написать его биографию было бы делом его друзей; но
замечательные люди исчезают у нас, не оставляя по себе следов. Мы
ленивы и нелюбопытны...>>}{А.\,С.~Пушкин. <<Путешествие в Арзрум>>
(1836).}

\begin{figure}[H]
\centering
\begin{subfigure}[t]{0.45\textwidth}
\centering
\includegraphics[width=0.7\textwidth]{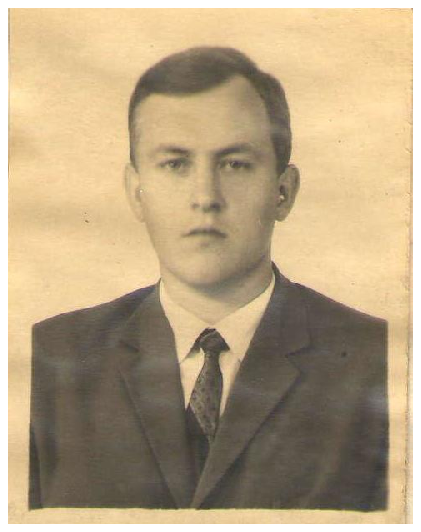}
\subcaption{\hfil Валерий Вячеславович Катрахов\hfil\\
\hphantom{\hspace{1.4cm}}(фото из личного дела 1972~г.)}
\label{pic3}
\end{subfigure}
\qquad
\begin{subfigure}[t]{0.45\textwidth}
\centering
\includegraphics[width=0.7\textwidth]{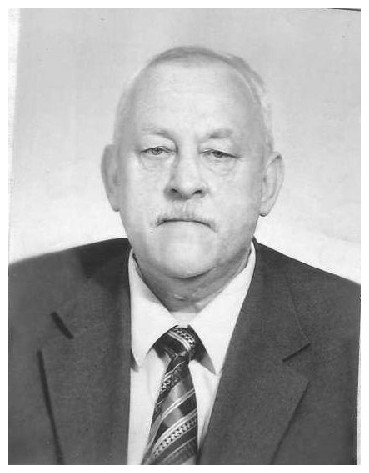}
\subcaption{\hfil Валерий Вячеславович Катрахов\hfil\\
\hphantom{\hspace{1.4cm}}(фото из личного дела 2006~г.)}
\label{pic4}
\end{subfigure}
\end{figure}

Валерий Вячеславович Катрахов родился 4~августа 1949~г. на
Сахалине, отец "--- Катрахов Вячеслав Тимофеевич,  офицер
Советской армии, впоследствии работал на Воронежском авиационном
заводе, мать "--- Катрахова Клавдия Васильевна, врач, работала в
поликлинике.

В 1966~г. закончил математический класс школы \No~15 г.~Воронежа.
С 1966 по 1971~гг. учился на математическом факультете
Воронежского государственного университета (ордена Ленина, им.
Ленинского Комсомола), который окончил с отличием. В 1971-1972~гг.
обучался в аспирантуре Воронежского государственного университета.

Дальнейшая трудовая деятельность В.\,В.~Катрахова во многом была
связана с факультетом прикладной математики и механики
Воронежского государственного университета. Факультет был основан
в 1969~г., когда произошло разделение математико-механического
факультета на математический факультет и факультет прикладной
математики и механики (ПММ). Факультет был образован, главным
образом, благодаря усилиям известного ученого в области механики,
профессора Геннадия Ивановича Быковцева, который стал его первым
деканом.

С 1972~г. В.\,В.~Катрахов работал в Воронежском университете на
кафедре дифференциальных уравнений факультета прикладной
математики и механики  ассистентом (1972--1974~гг.), старшим
преподавателем (1974--1980~гг.), доцентом (1980--1983~гг.),
заведующим кафедрой вычислительной математики факультета
прикладной математики и механики (1983--1987~гг.). В 1977~г. его
работа <<Сингулярные краевые задачи и их приложение к оптимальному
управлению>> была удостоена премии Воронежского Комсомола в
области науки и техники, в течение ряда лет он руководил научным
студенческим обществом факультета прикладной математики и
механики, работал на общественных началах заместителем декана по
научной работе, являлся членом совета факультета ПММ и
специализированного совета ВГУ по защите кандидатских диссертаций.

С 1987~г. по 2006~г. В.\,В.~Катрахов работал во Владивостоке.
Сначала в 1987--1995~гг. в институте прикладной математики ДВО АН
СССР (затем ДВО РАН) в должностях старшего научного сотрудника,
заведующего лабораторией, заместителя директора,  а  в
1995--2006~гг. профессором Дальневосточного коммерческого
института (ДВКИ, сейчас "--- Тихоокеанский государственный
экономический университет (ТГЭУ)), а также в ряде других
дальневосточных вузов.

В 2006~г. вернулся в Воронеж, где работал профессором на кафедрах
дифференциальных уравнений, а также программного обеспечения и
администрирования систем факультета прикладной математики и
механики Воронежского государственного университета до 2008~г.

В 1974~г. защитил кандидатскую диссертацию в Воронежском
государственном университете на тему <<О спектральной функции
сингулярных дифференциальных операторов>> под руководством
профессора Ивана Александровича Киприянова, а в  1989~г. "---
докторскую диссертацию в Новосибирском государственном
университете на тему <<Сингулярные эллиптические краевые задачи.
Метод операторов преобразования>>.

В.\,В.~Катрахов получил ряд значительных результатов в различных
областях математики, среди них "--- теория сингулярных и
вырождающихся дифференциальных уравнений, операторы преобразования
и интегральные преобразования, специальные функции и операторы
дробного интегродифференцирования, теория функций и функциональный
анализ, спектральная теория и псевдодифференциальные операторы,
численные методы и оптимальное управление, математическая физика и
моделирование систем. Выделим несколько важных результатов, для
получения которых В.\,В.~Катраховым были использованы новые идеи и
методы: решение систем сингулярных дифференциальных и
псевдодифференциальных уравнений с весовыми краевыми условиями,
постановка и решение новых краевых задач для решений с
существенными особенностями, введение нового свёрточного
нелокального краевого условия <<$K$-следа>> для решений с
особенностями, введение нового класса счётно-нормированных
пространств Фреше на основе операторов преобразований для
постановки задач и доказательства корректности для эллиптических
уравнений с существенными сингулярностями в изолированных особых
точках, введение и применение новых классов операторов
преобразования и дробного интегродифференцирования, новые
постановки краевых задач для дифференциальных уравнений с
особенностями в угловых областях и в пространстве Лобачевского,
матричный метод исследования модели Изинга. В своих основных
работах по теории дифференциальных уравнений использовал и
развивал метод операторов преобразования.

Среди учеников В.\,В.~Катрахова четыре доктора наук
(Н.\,И.~Головко, А.\,Б.~Муравник, И.\,П.~Половинкин,
С.\,М.~Ситник) и около десяти кандидатов наук.

В.\,В.~Катрахов является автором около 150 научных работ в ведущих
российских журналах, таких как <<Доклады АН СССР>>, <<Доклады
РАН>>, <<Дифференциальные уравнения>>, <<Математический сборник>>,
<<Сибирский математический журнал>> и др., а также автором семи
монографий:
\begin{itemize}
\item Катрахов~В.\,В., Мазелис~Л.\,С. <<Непрерывность, пополнение,
замыкание в метрических пространствах>>,
    2000, Владивосток, изд. ДВО, 112~с.

\item  Головко~Н.\,И., Катрахов~В.\,В. <<Анализ систем массового
обслуживания, функционирующих в случайной среде>>,      2000,
Владивосток, Изд-во Дальневосточной государственной академии
экономики и управления (ДВГАЭУ), 144~с.

\item Дмитриев~А.\,А., Катрахов~В.\,В., Харченко~Ю.\,Н. <<Корневые
трансфер-матрицы в моделях Изинга>>,
    2004, Москва, Наука, 192~с.

\item Катрахов~В.\,В., Рыжков~Д.\,Е. <<Введение в
функционально-аналитический метод в динамической теории массового
обслуживания>>, 2004, Владивосток, изд-во ДВГУ, 102~с.

\item Катрахов~В.\,В., Головко~Н.\,И., Рыжков~Д.\,Е. <<Введение в
теорию марковских дважды стохастических систем массового
обслуживания>>,
    2005, Владивосток, Изд-во ДВГУ,  212~с.

\item Головко~Н.\,И., Катрахов~В.\,В. <<Применение моделей СМО в
информационных системах>>,
    2008, Владивосток, Изд-во ТГЭУ, 272~с.

\item Катрахов~В.\,В., Ситник~С.\,М. <<Метод операторов
преобразования и краевые задачи для сингулярных эллиптических
уравнений>>, журнал <<Современная математика. Фундаментальные
направления>>, 2018, т.~64, \No~2, 216~с.

\end{itemize}

Был женат, жена "--- Катрахова Алла Анатольевна, преподаватель
математики, дочь "--- Катрахова Алла Валерьевна, врач.

В.\,В.~Катрахов умер в 2010~г., похоронен на Коминтерновском
кладбище г.~Воронежа.

\def\numberline{}

\engthebibliography{999}{

\bibitem{eng_AbOs}  Z.~L.~Abzhandadze and V.~F.~Osipov, {\it
Preobrazovanie Fur'e---Frenelya i nekotorye ego prilozheniya}
[Fourier--Fresnel Transform and some Applications], S.-Peterb.
Univ., SPb., 2000 (in Russian).

\bibitem{eng_AbSi}  M.~Ablowitz  and H.~Segur, {\it Solitony i metod
obratnoy zadachi} [Solitons and the Inverse Scattering Transform],
Mir, Moscow, 1979 (Russian translation).

\bibitem{eng_AS}  M.~Abramowitz and I.~Stegun, {\it Spravochnik po
spetsial'nym funktsiyam} [Handbook Of Mathematical Functions],
Nauka, Moscow, 1979 (Russian translation).

\bibitem{eng_1}  S.~Agmon, A.~Douglis, and L.~Nirenberg, {\it Otsenka resheniy
ellipticheskikh uravneniy vblizi granitsy} [Estimates Near the
Boundary for Solutions of Elliptic Partial Differential
Equations], Inostr. lit., Moscow, 1962 (Russian translation).

\bibitem{eng_AM}
Z.~S.~Agranovich and V.~A.~Marchenko, {\it Obratnaya zadacha
teorii rasseyaniya} [Inverse Scattering Problem], KHGU, Khar'kov,
1960 (in Russian).

\bibitem{eng_Agr}  M.~S.~Agranovich, {\it Ellipticheskie psevdodifferentsial'nye
operatory. Ch.~1,~2} [Elliptic Pseudodifferential Operators.
V.~1,~2], 2003, Moscow, 2004 (in Russian).

\bibitem{eng_AMR}  N.~V.~Azbelev,
V.~P.~Maksimov, and L.~F.~Rakhmatullina, {\it Elementy sovremennoy
teorii funktsional'no-differentsial'nykh uravneniy. Metody i
prilozheniya} [Elements of the Modern Functional-Differential
Equations Theory. Methods and Applications], Inst. Komp. Issl.,
Moscow, 2002 (in Russian).

\bibitem{eng_2}   Sh.~A.~Alimov, ``Drobnye stepeni ellipticheskikh operatorov i
izomorfizm klassov differentsiruemykh funktsiy'' [Fractional
powers of elliptic operators and isomorphism of classes of
differentiable functions], {\it Diff. uravn.} [Differ. Equ.],
1972, {\bf 8}, No.~9, 1609--1626 (in Russian).

\bibitem{eng_Arsh}   E.~A.~Arshava, ``Obrashchenie
integral'nykh operatorov metodom operatornykh tozhdestv''
[Inversion of integral operators by the method of integral
identities], {\it Nauch. vedom. Belgorod. gos. un-ta. Ser. Mat.
Fiz.} [Sci. Bull. Belgorod State Univ. Ser. Math. Phys.], 2009,
{\bf 17/2}, No.~13 (68), 18--29 (in Russian).

\bibitem{eng_Ahi1}  N.~I.~Akhiezer, ``K teorii sparennykh integral'nykh
uravneniy'' [To the theory of paired integral equations], {\it
Uch. zap. Khar'kov. gos. un-ta.} [Sci. Notes Khar'kov State
Univ.], 1957, {\bf 80}, 5--21 (in Russian).

\bibitem{eng_Ahi2}  N.~I.~Akhiezer, {\it Lektsii ob integral'nykh
preobrazovaniyakh} [Lectures on Integral Transforms], Vishcha
shkola. Khar'kov State Univ., Khar'kov, 1984 (in Russian).

\bibitem{eng_3}  V.~V.~Babikov, {\it Metod fazovykh funktsiy v
kvantovoy mekhanike} [A Method of Phase Functions in Quantum
Mechanics], Nauka, Moscow, 1976 (in Russian).

\bibitem{eng_BMYa}
I.~I.~Bavrin, V.~L.~Matrosov, and  O.~E.~Yaremko, {\it Operatory
preobrazovaniya dlya kraevykh zadach, integral'nykh predstavleniy
i vosstanovleniya zavisimostey} [Transmutation Operators for
Boundary--Value Problems, Integral Representations and Dependence
Recovering], Prometey, Moscow, 2015 (in Russian).

\bibitem{eng_BaSa}  V.~G.~Bagrov and  B.~A.~Samsonov,
``Preobrazovanie Darbu uravneniya Shredingera'' [Darboux transform
for the Shr\"odinger equation], {\it Fiz. elem. chastits i atom.
yadra} [Phys. Elem. Part. Atom Kern.], 1997,  {\bf 28}, No.~4,
951--1012 (in Russian).

\bibitem{eng_Baid}  A.~N.~Baydakov, ``Apriornye otsenki
gel'derovykh norm resheniy kvazilineynykh $B$-ellipticheskikh
uravneniy'' [A priori estimates of H\"older norms of solutions to
quasilinear $B$-elliptic equations], {\it Diff. uravn.} [Differ.
Equ.], 1987, {\bf 23},  No.~11, 1923--1930 (in Russian).

\bibitem{eng_Bas}  A.~G.~Baskakov, {\it Garmonicheskiy analiz
lineynykh operatorov} [Harmonic Analysis of Linear Operators],
VGU, Voronezh, 1987 (Russian translation).

\bibitem{eng_BE1}  H.~Bateman and  A.~Erd\'elyi, {\it Vysshie transtsendentnye
funktsii. T.~1} [Higher Transcendental Functions. V.~1], Nauka,
Moscow, 1966 (Russian translation).

\bibitem{eng_BE2}  H.~Bateman and  A.~Erd\'elyi, {\it Vysshie transtsendentnye funktsii. T.~2}
[Higher Transcendental Functions. V.~2], Nauka, Moscow, 1966 (in
Russian).

\bibitem{eng_BE3}  H.~Bateman and  A.~Erd\'elyi, {\it Vysshie transtsendentnye
funktsii. T.~3} [Higher Transcendental Functions. V.~3], Nauka,
Moscow, 1967 (Russian translation).

\bibitem{eng_BB}  E.~Beckenbach and  R.~Bellman, {\it Neravenstva}
[Inequalities], Mir, Moscow, 1965 (Russian translation).

\bibitem{eng_Berg}  S.~Bergman, {\it Integral'nye operatory v teorii uravneniy s chastnymi
proizvodnymi} [Integral Operators in the Theory of Linear Partial
Differential Equations], Mir, Moscow, 1964 (Russian translation).

\bibitem{eng_6}  Yu.~M.~Berezanskiy, {\it Razlozheniya po
sobstvennym funktsiyam samosopryazhennykh operatorov} [Expansions
in Eigenfunctions of Self-Adjoint Operators], Naukova Dumka, Kiev,
1965 (in Russian).

\bibitem{eng_Berk}  L.~M.~Berkovich, {\it Faktorizatsiya i preobrazovaniya differentsial'nykh
uravneniy. Metody i prilozheniya} [Factorization and
Transformations of Differential Equations. Methods and
Applications], RKHD, Moscow, 2002 (in Russian).

\bibitem{eng_7}   L.~Bers, {\it Matematicheskie voprosy dozvukovoy
i okolozvukovoy gazovoy dinamiki} [Mathematical Aspects of
Subsonic and Transonic Gas Dynamics], Inostr. lit., Moscow, 1961
(Russian translation).

\bibitem{eng_8}  O.~V.~Besov, V.~P.~Il'in, and S.~M.~Nikol'skiy, {\it Integral'nye
predstavleniya funktsiy i teoremy vlozheniya} [Integral
Representations of Functions and Embedding Theorems], Nauka,
Moscow, 1975 (in Russian).

\bibitem{eng_Bitz2}
A.~V.~Bitsadze, {\it Uravneniya smeshannogo tipa} [Equations of
Mixed Type], Izd. AN SSSR, Moscow, 1959 (in Russian).

\bibitem{eng_9}  A.~V.~Bitsadze, {\it Nekotorye klassy uravneniy v
chastnykh proizvodnykh} [Some Classes of Partial Differential
Equations], Nauka, Moscow, 1981 (in Russian).

\bibitem{eng_Bitz1}  A.~V.~Bitsadze and V.~I.~Pashkovskiy, ``K teorii uravneniy
Maksvella---Eynshteyna'' [To the theory of Maxwell--Einstein
equations], {\it Dokl. AN SSSR} [Rep. Acad. Sci. USSR], 1974, {\bf
216}, No.~2, 9-10 (in Russian).

\bibitem{eng_Bitz12}
A.~V.~Bitsadze and V.~I.~Pashkovskiy, ``O nekotorykh klassakh
resheniy uravneniya Maksvella---Eynshteyna'' [On some classes of
Maxwell--Einstein equation], {\it Tr. MIAN} [Proc. Math. Inst.
Russ. Acad. Sci.], 1975, {\bf 134}, 26--30 (in Russian).

\bibitem{eng_Bloh}  A.~Sh.~Blokh, ``Ob opredelenii differentsial'nogo operatora po ego
spektral'noy matritse-funktsii'' [On defining of a differential
operator by its spectral matrix-function], {\it Dokl. AN SSSR}
[Rep. Acad. Sci. USSR], 1953,  {\bf 92}, No.~2, 209--212 (in
Russian).

\bibitem{eng_Bor1}  A.~V.~Borovskikh, ``Formula rasprostranyayushchikhsya voln dlya
odnomernoy neodnorodnoy sredy'' [A formula of wave extending for
one-dimensional nonhomogeneous medium], {\it Diff. uravn.}
[Differ. Equ.], 2002,  {\bf 38}, No.~6, 758--767 (in Russian).

\bibitem{eng_Bor2}  A.~V.~Borovskikh, ``Metod rasprostranyayushchikhsya
voln'' [A method of extending waves], {\it Tr. sem. im.
I.~G.~Petrovskogo} [Proc. Petrovskii Semin.], 2004, {\bf 24},
3--43 (in Russian).

\bibitem{eng_Boyar}   B.~Boyarskiy, ``Obobshchennye resheniya sistemy differentsial'nykh
uravneniy pervogo poryadka ellip\-ti\-ches\-kogo tipa s razryvnymi
koeffitsientami'' [Generalized solutions to a first order system
of dif\-fe\-ren\-tial equations of elliptic type with discontinuous
coefficients], {\it Mat. sb.} [Math. Digest], 1957,  {\bf 43},
No.~4, 451--503 (in Russian).

\bibitem{eng_11}  Yu.~A.~Brychkov and
A.~P.~Prudnikov, {\it Integral'nye preobrazovaniya obobshchennykh
funktsiy} [Integral Transforms of Generalized Functions], Nauka,
Moscow, 1977 (in Russian).

\bibitem{eng_Burb}  N.~Bourbaki, {\it Funktsii deystvitel'nogo
peremennogo} [Functions of a Real Variable], Nauka, Moscow, 1965
(Russian translation).

\bibitem{eng_Bur1}  V.~I.~Burenkov, {\it Funktsional'nye
prostranstva} [Functional Spaces], RUDN, Moscow, 1989 (in
Russian).

\bibitem{eng_Bur2}  V.~I.~Burenkov and  M.~L.~Gol'dman, {\it Metodicheskie
rekomendatsii k izucheniyu kursa ``Funktsional'nye prostranstva''}
[Methodical Materials to the Course of Functional Spaces], RUDN,
Moscow, 1989 (in Russian).

\bibitem{eng_But}  S.~A.~Buterin, ``O vosstanovlenii svertochnogo vozmushcheniya operatora
Shturma---Liuvillya po spektru'' [On a recovery of a convolution
perturbation of Sturm--Lioville operator by spectrum], {\it Diff.
uravn.} [Differ. Equ.], 2010,  {\bf 46}, No.~1, 146--149 (in
Russian).

\bibitem{eng_Val}   Yu.~N.~Valitskiy, ``Ob operatore preobrazovaniya dlya
integro-differentsial'nykh operatorov tipa Vol'terra'' [On a
transmutation operator for integral-differential operators of
Volterra type], In: {\it Matematicheskaya fizika} [Mathematical
Physics], Naukova Dumka, Kiev, 1965, pp.~23--36 (in Russian).

\bibitem{eng_VaRo}
E.~M.~Varfolomeev and L.~E.~Rossovskiy, {\it
Funktsional'no-differentsial'nye uravneniya i ikh prilozheniya k
issledovaniyu neyronnykh setey i peredache informatsii nelineynymi
lazernymi sistemami s obratnoy svyaz'yu. Ucheb. posobie}
[Functional-Differential Equations and Their Applications to the
Study of Neuron Nets and Information Transmission by Nonlinear
Laser Systems with Feedback. Textbook], RUDN, Moscow, 2008 (in
Russian).

\bibitem{eng_Wat}
G.~N.~Watson, {\it Teoriya besselevykh funktsiy. T.~1} [A Treatise
on the Theory of Bessel Functions. V.~1], Inostr. lit., Moscow,
1949 (Russian translation).

\bibitem{eng_12}    A.~A.~Vasharin and P.~I.~Lizorkin, ``Nekotorye kraevye zadachi dlya
ellipticheskikh uravneniy s sil'nym vyrozhdeniem na granitse''
[Some boundary-value problems for elliptic equations with strong
degeneration on the boundary], {\it Dokl. AN SSSR} [Rep. Acad.
Sci. USSR], 1961, {\bf 137}, No.~5, 1015--1018 (in Russian).

\bibitem{eng_Vek1}   I.~N.~Vekua, ``O resheniyakh
uravneniya $\Delta u + \lambda^2 u$'' [On solutions to the
equation $\Delta u + \lambda^2 u$], {\it Soobshch. AN Gruz. SSR}
[Bull. Acad. Sci. Georgian SSR], 1942, {\bf 3}, No.~4, 307--314
(in Russian).

\bibitem{eng_Vek2}   I.~N.~Vekua, ``Obrashchenie odnogo integral'nogo
preobrazovaniya i ego nekotorye primeneniya'' [Inversion of one
integral transform and some of its applications], {\it Soobshch.
AN Gruz. SSR} [Bull. Acad. Sci. Georgian SSR], 1945,  {\bf 6},
No.~3, 177--183 (in Russian).

\bibitem{eng_13}   I.~N.~Vekua, ``Ob odnom obobshchenii integrala
Puassona dlya ploskosti'' [On a generalization of the Poisson
integral on the plane], {\it Dokl. AN SSSR} [Rep. Acad. Sci.
USSR], 1947,  {\bf 56}, No.~2, 229--231 (in Russian).

\bibitem{eng_Vek3}   I.~N.~Vekua, {\it Novye metody resheniya
ellipticheskikh uravneniy} [New Methods for Solving Elliptic
Equations], GITTL, Moscow--Leningrad, 1948 (in Russian).

\bibitem{eng_Vek4}   I.~N.~Vekua, {\it Obobshchennye analiticheskie
funktsii} [Generalized Analytic Functions], Nauka, Moscow, 1988
(in Russian).

\bibitem{eng_ViGa}   N.~A.~Virchenko and V.~Gaydey, {\it Klassicheskie
i obobshchennye mnogoparametricheskie funktsii} [Classical and
Generalized Multiparameter Functions], Kiev, 2008 (in Ukrainian).

\bibitem{eng_ViRy}   N.~A.~Virchenko and V.~Ya.~Rybak, {\it Osnovy drobnogo
integrodifferentsirovaniya} [Foundations of Fractional
Integrodifferentiation], Kiev, 2007 (in Ukrainian).

\bibitem{eng_15}
M.~I.~Vishik and V.~V.~Grushin, ``Kraevye zadachi dlya
ellipticheskikh uravneniy, vyrozhdayushchikhsya na granitse
oblasti'' [Boundary-value problems for elliptic equations
degenerated at the boundary of domain], {\it Mat. sb.} [Math.
Digest], 1969,  {\bf 80}, 455--491 (in Russian).

\bibitem{eng_16}   B.~C.~Vladimirov, {\it Obobshchennye funktsii v
matematicheskoy fizike} [Generalized Functions in Mathematical
Physics], Nauka, Moscow, 1979 (in Russian).

\bibitem{eng_Volk}   V.~Ya.~Volk, ``O formulakh obrashcheniya dlya differentsial'nogo uravneniya
s osobennost'yu pri $x=0$'' [On inversion formulas for a
differential equation with a singularity at $x=0$], {\it Usp. mat.
nauk} [Progr. Math. Sci.], 1953,  {\bf 111}, No.~4, 141--151 (in
Russian).

\bibitem{eng_VoKa}   I.~K.~Volkov and A.~N.~Kanatnikov, {\it Integral'nye
preobrazovaniya i operatsionnoe ischislenie} [Integral Transforms
and Operational Calculus], MGTU im. N.~E.~Baumana, Moscow, 2002
(in Russian).

\bibitem{eng_VZ}  V.~F.~Volkodavov and V.~N.~Zakharov, {\it Tablitsy funktsiy Rimana i
Rimana---Adamara dlya nekotorykh differentsial'nykh uravneniy v
$n$-mernykh evklidovykh prostranstvakh} [Tables of Riemann and
Riemann--Hadamard Functions for Some Differential Equations in
$n$-Dimensional Euclidean Spaces], Samara, 1994 (in Russian).

\bibitem{eng_VNN}   V.~F.~Volkodavov,
M.~E.~Lerner, N.~Ya.~Nikolaev, and V.~A.~Nosov, {\it Tablitsy
nekotorykh funktsiy Rimana, integralov i ryadov} [Tables of Some
Riemann Functions, Integrals and Series], Kuybyshev. Gos. Ped.
Inst., Kuybyshev, 1982 (in Russian).

\bibitem{eng_VoNi}
V.~F.~Volkodavov and N.~Ya.~Nikolaev, {\it Kraevye zadachi dlya
uravneniya Eylera---Paussona---Darbu} [Boundary--Value Problems
for the Euler--Poisson--Darboux equation], Izd. Kuybyshev. Gos.
Ped. Inst., Kuybyshev, 1984 (in Russian).

\bibitem{eng_VN}   V.~F.~Volkodavov and N.~Ya.~Nikolaev, {\it Integral'nye uravneniya
Vol'terra pervogo roda s nekotorymi spetsial'nymi funktsiyami v
yadrakh i ikh prilozheniya} [Integral Equations of Volterra of the
First Kind with Some Kernel Special Functions and Applications],
Samara Univ., Samara, 1992 (in Russian).

\bibitem{eng_Vrag}  V.~N.~Vragov, {\it Kraevye
zadachi dlya neklassicheskikh uravneniy matematicheskoy fiziki}
[Boundary-Value Problems for Nonclassical Equations of
Mathematical Physics], NGU, Novosibirsk, 1983 (in Russian).

\bibitem{eng_Gan}   F.~R.~Gantmacher, {\it Teoriya matrits} [Matrix Theory],
Nauka, Moscow, 1988 (in Russian).

\bibitem{eng_17}   I.~M.~Gel'fand, M.~I.~Graev, and N.~Ya.~Vilenkin, {\it Integral'naya
geometriya i svyazannye s ney voprosy teorii predstavleniy}
[Integral Geometry and Connected Problems of Representation
Theory], Gos. Fiz.-Mat. Lit., Moscow, 1962 (in Russian).

\bibitem{eng_Glu2}   A.~V.~Glushak, ``O vozmushchenii abstraktnogo uravneniya
Eylera---Puassona---Darbu'' [On perturbation of abstract
Euler--Poisson--Darboux equation], {\it Mat. zametki} [Math.
Notes], 1996,  {\bf 60}, No.~3, 363--369 (in Russian).

\bibitem{eng_Glu4}   A.~V.~Glushak, ``O stabilizatsii resheniya zadachi Dirikhle
dlya odnogo ellipticheskogo uravneniya v banakhovom prostranstve''
[On stabilization of solution of Dirichlet problem for an elliptic
equation in the Banach space], {\it Diff. uravn.} [Differ. Equ.],
1997, {\bf 33}, No.~4, 433--437 (in Russian).

\bibitem{eng_Glu3}   A.~V.~Glushak, ``Operatornaya funktsiya
Besselya'' [Bessel Operator Function], {\it Dokl. RAN} [Rep. Russ.
Acad. Sci.], 1997, {\bf 352}, No.~5, 587--589 (in Russian).

\bibitem{eng_Glu5}   A.~V.~Glushak, ``Operatornaya funktsiya Besselya i svyazannye s neyu
polugruppy i modifitsirovannoe preobrazovanie Gil'berta''
[Operator Bessel function and connected semigroups and modified
Hilbert transform], {\it Diff. uravn.} [Differ. Equ.], 1999,  {\bf
35}, No.~1, 128--130 (in Russian).

\bibitem{eng_Glu55}   A.~V.~Glushak, ``Regulyarnoe i
singulyarnoe vozmushcheniya abstraktnogo uravneniya
Eylera---Puassona---Darbu'' [Regular and singular perturbations
for an abstract Euler--Poisson--Darboux equations], {\it Mat.
zametki} [Math. Notes], 1999,  {\bf 66}, No.~3, 364--371 (in
Russian).

\bibitem{eng_Glu6}   A.~V.~Glushak, ``Operatornaya funktsiya
Lezhandra'' [Legendre operator function], {\it Izv. RAN. Ser.
Mat.} [Bull. Russ. Acad. Sci. Ser. Math.], 2001,  {\bf 65}, No.~6,
3--14 (in Russian).

\bibitem{eng_Glu7}   A.~V.~Glushak, ``Zadacha tipa Koshi dlya abstraktnogo
differentsial'nogo uravneniya s drobnymi proizvodnymi'' [A problem
of Cauchy-type for abstract equation with fractional derivatives],
{\it Mat. zametki} [Math. Notes], 2005,  {\bf 77}, No.~1, 28--41
(in Russian).

\bibitem{eng_Glu8}   A.~V.~Glushak, ``O svyazi
prointegrirovannoy kosinus-operator-funktsii s operatornoy
funktsiey Besselya'' [On connection of integrated cosine operator
function with Bessel operator function], {\it Diff. uravn.}
[Differ. Equ.], 2006, {\bf 42}, No.~5, 583--589 (in Russian).

\bibitem{eng_Glu10}
A.~V.~Glushak, ``Nachal'naya zadacha dlya slabo nagruzhennogo
uravneniya Eylera---Puassona---Darbu'' [Initial problem for weakly
loaded Euler--Poisson--Darboux equation], Math. Int. Sci. Conf.
{\it ``Aktual'nye problemy teorii uravneniy v chastnykh
proizvodnykh''} [Actual problems in the theory of partial
differential equations] MGU, Moscow, 2016, p.~101 (in Russian).

\bibitem{eng_Glu11}
A.~V.~Glushak, ``Nelokal'naya zadacha dlya abstraktnogo uravneniya
Eylera---Puassona---Darbu'' [Nonlocal problem for abstract
Euler--Poisson--Darboux equation], {\it Izv. vuzov. Ser. Mat.}
[Bull. Higher Edu. Inst. Ser. Math.], 2016,  No.~6, 1--9 (in
Russian).

\bibitem{eng_Glu12}
A.~V.~Glushak, ``Abstraktnaya zadacha Koshi dlya uravneniya
Besselya---Struve'' [Abstract problem for the Bessel--Struve
equation], {\it Diff. uravn.} [Differ. Equ.], 2017,  {\bf 53},
No.~7, 891--905 (in Russian).

\bibitem{eng_Glu1}   A.~V.~Glushak,
V.~I.~Kononenko, and S.~D.~Shmulevich, ``Ob odnoy singulyarnoy
abstraktnoy zadache Koshi'' [On a singular abstract Cauchy
problem], {\it Izv. vuzov. Ser. Mat.} [Bull. Higher Edu. Inst.
Ser. Math.], 1986, No.~6, 55--56 (in Russian).

\bibitem{eng_Glu9}
A.~V.~Glushak and O.~A.~Pokruchin, ``Kriteriy razreshimosti
zadachi Koshi dlya abstraktnogo uravneniya
Eylera---Puassona---Darbu'' [Criterion of solvability for Cauchy
problem for abstract Euler--Poisson--Darboux equations], {\it
Diff. uravn.} [Differ. Equ.], 2016, {\bf 52}, No.~1, 41--59 (in
Russian).

\bibitem{eng_Glu13}   A.~V.~Glushak and T.~G.~Romanchenko, ``Formuly svyazi mezhdu
resheniyami abstraktnykh singulyarnykh differentsial'nykh
uravneniy'' [Connection formulas for solutions of abstract
singular differential equa\-tions], {\it Nauch. vedom. Belgorod.
gos. un-ta. Ser. Mat. Fiz.} [Sci. Bull. Belgorod State Univ. Ser.
Math. Phys.], 2016, {\bf 42}, No.~6, 36--39 (in Russian).

\bibitem{eng_Glushko}   V.~P.~Glushko, {\it Lineynye vyrozhdayushchiesya differentsial'nye
uravneniya} [Linear Degenerate Dif\-fe\-ren\-tial Equations],
 VGU, Voronezh, 1972 (in Russian).

\bibitem{eng_GKE}  L.~S.~Gnoenskiy,
G.~A.~Kamenskiy, and L.~E.~El'sgol'ts, {\it Matematicheskie osnovy
teorii upravlyaemykh sistem} [Mathematical Foundations of the
Theory of Controllable Systems], Nauka, Moscow, 1969 (in Russian).

\bibitem{eng_Grin}   G.~A.~Grinberg, {\it Izbrannye voprosy
matematicheskoy teorii elektricheskikh i magnitnykh yavleniy}
[Selected Problems of Mathematical Theory of Electric and Magnetic
Phenomena], AN SSSR, Moscow, 1948 (in Russian).

\bibitem{eng_Gul2}  V.~S.~Guliev, {\it Integral'nye operatory, funktsional'nye prostranstva
i voprosy approksimatsii na gruppe Geyzenberga} [Integral
Operators, Functional Spaces and Approximation Problems for the
Heisenberg Group], ELM, Baku, 1996 (in Russian).

\bibitem{eng_Gul1}   V.~S.~Guliev, {\it Funktsional'nye prostranstva, integral'nye
operatory i dvukhvesovye otsenki na odnorodnykh gruppakh.
Nekotorye prilozheniya} [Functional Spaces, Integral Operators and
Two-Weighted Estimates on Homogeneous Groups. Some Applications],
Chashyogly, Baku, 1999 (in Russian).

\bibitem{eng_Gur}   M.~I.~Gurevich, {\it Teoriya struy ideal'noy
zhidkosti} [Jet Theory for an Ideal Liquid], Nauka, Moscow, 1979
(in Russian).

\bibitem{eng_Gus1}   I.~M.~Guseynov, ``Ob odnom operatore
preobrazovaniya'' [On one transmutation operator], {\it Mat.
zametki} [Math. Notes], 1997,  {\bf 62}, No.~2, 206--215 (in
Russian).

\bibitem{eng_Gus2}   I.~M.~Guseynov, A.~A.~Nabiev, and R.~T.~Pashaev, ``Operatory
preobrazovaniya i asimptoticheskie formuly dlya sobstvennykh
znacheniy polinominal'nogo puchka operatorov Shturma---Liuvillya''
[Transmutations and asymptotic formulas for eigenvalues of a
polynomial pencil of Sturm--Liouville operators], {\it Sib.
mat.~zh.} [Siberian Math.~J.], 2000, {\bf 41}, No.~3, 554--566 (in
Russian).

\bibitem{eng_Dza2}   G.~V.~Dzhayani, {\it Reshenie nekotorykh zadach dlya odnogo
vyrozhdayushchegosya ellipticheskogo uravneniya i ikh prilozheniya
k prizmaticheskim obolochkam} [Solution of Some Problems for a
Degenerated Elliptic Equations and Applications to Prismatic
Shells], Tbilisi Univ., Tbilisi, 1982 (in Russian).

\bibitem{eng_Dza1}   G.~V.~Dzhayani, {\it Uravnenie
Eylera---Paussona---Darbu} [Euler--Poisson--Darboux Equation],
 Tbilisi Univ., Tbilisi, 1984 (in Russian).

\bibitem{eng_Gini}   C.~Gini, {\it Srednie velichiny} [Means],
Statistika, Moscow, 1970 (Russian translation).

\bibitem{eng_Dzh1}
M.~M.~Dzhrbashyan, {\it Integral'nye preobrazovaniya i
predstavleniya funktsiy v kompleksnoy oblasti} [Integral
Transforms and Representations of Functions in the Complex
Domain], Nauka, Moscow, 1966 (in Russian).

\bibitem{eng_Din}   Kh.~A.~Din', ``Integral'nye uravneniya
s funktsiey Lezhandra v yadrakh v osobykh sluchayakh'' [Integral
equations with Legendre function as kernels in special cases],
{\it Dokl. AN Belorus. SSR} [Rep. Acad. Sci. Belorus. SSR], 1989,
{\bf 33}, No.~7, 591--594 (in Russian).

\bibitem{eng_EPP}   I.~E.~Egorov, S.~G.~Pyatkov, and
S.~V.~Popov, {\it Neklassicheskie differentsial'no-operatornye
uravneniya} [Nonclassical Differential-Operator Equations], Nauka,
Novosibirsk, 2000 (in Russian).

\bibitem{eng_EgFed}   I.~E.~Egorov and V.~Evs.~Fedorov, {\it Neklassicheskie
uravneniya matematicheskoy fiziki vysokogo poryadka} [Nonclassical
Higher Order Equations of Mathematical Physics],  VTS SO RAN,
Novosibirsk, 1995 (in Russian).

\bibitem{eng_ZhM1}   V.~I.~Zhegalov and A.~N.~Mironov, {\it Differentsial'nye uravneniya
so starshimi chastnymi proizvodnymi} [Differential Equations with
Highest Partial Derivatives], Kazan. Mat. Ob-vo, Kazan', 2001 (in
Russian).

\bibitem{eng_ZhM2}   V.~I.~Zhegalov, A.~N.~Mironov, and E.~A.~Utkina, {\it Uravneniya s
dominiruyushchey chastnoy proizvodnoy} [Equations with Dominating
Partial Derivative], Kazan. Univ, Kazan', 2014 (in Russian).

\bibitem{eng_Zhit}
Ya.~I.~Zhitomirskiy, ``Zadacha Koshi dlya sistem lineynykh
uravneniy v chastnykh proizvodnykh s differentsial'nymi
operatorami tipa Besselya'' [A Cauchy problem for for systems of
partial differential equations with Bessel-type differential
operators], {\it Mat. sb.} [Math. Digest], 1955,  {\bf 36}, No.~2,
299--310 (in Russian).

\bibitem{eng_ZhuSi1}   N.~V.~Zhukovskaya and S.~M.~Sitnik,
``Differentsial'nye uravneniya tipa Eylera drobnogo poryadka''
[Differential equations of fractional order of Euler type], {\it
Mat. zametki SVFU} [Math. Notes North-East Fed. Univ.], 2018, {\bf
25}, No.~2, 27--39 (in Russian).

\bibitem{eng_Zhu}   V.~M.~Zhuravlev, {\it Nelineynye
volny. Tochno reshaemye zadachi} [Nonlinear Waves. Exactly
Solvable Problems], Ul'yanovsk, 2001 (in Russian).

\bibitem{eng_Zai1}  V.~A.~Zaytsev, ``O printsipe Gyuygensa dlya
nekotorykh uravneniy s osobennostyami'' [On Huygens principle for
some equations with singularities], {\it Dokl. AN SSSR} [Rep.
Acad. Sci. USSR], 1978,  {\bf 242}, No.~1, 28--31 (in Russian).

\bibitem{eng_Zai2}  V.~A.~Zaytsev, ``Slabye
lakuny dlya odnomernykh strogo giperbolicheskikh uravneniy s
postoyannymi koef\-fi\-tsi\-en\-tami'' [Weak lacunas for one-dimensional
strictly hyperbolic equations with constant coefficients], {\it
Sib. mat.~zh.} [Siberian Math.~J.], 1984,  {\bf 25}, No.~4, 54--62
(in Russian).

\bibitem{eng_ZMNP}
V.~E.~Zakharov,  S.~V.~Manakov, S.~P.~Novikov, and
L.~P.~Pitaevskiy, {\it Teoriya solitonov: metod obratnoy zadachi}
[Soliton Theory: the Method of Inverse Problem], Nauka, Moscow,
1980 (in Russian).

\bibitem{eng_ZKNE}  A.~M.~Zverkin,
G.~A.~Kamenskiy, S.~B.~Norkin, and L.~E.~El'sgol'ts,
``Differentsial'nye uravneniya s otklonyayushchimsya argumentom''
[Differential Equations with Deviating Argument], {\it Usp. mat.
nauk} [Progr. Math. Sci.], 1962, {\bf 17}, No.~2, 77--164 (in
Russian).

\bibitem{eng_18}   A.~I.~Ibragimov, ``O povedenii v okrestnosti granichnykh
tochek i teoremy ob ustranimykh mnozhestvakh dlya ellipticheskikh
uravneniy vtorogo poryadka s nepreryvnymi koeffitsientami'' [On
behavior near boundary points and theorems on removable
singularities for the second order elliptic equations with
continuous derivatives], {\it Dokl. AN SSSR} [Rep. Acad. Sci.
USSR], 1980,  {\bf 250}, No.~1, 25--28 (in Russian).

\bibitem{eng_19}  V.~M.~Ivakin, ``Vidoizmenennaya zadacha Dirikhle dlya
vyrozhdayushchikhsya na granitse uravneniy i sistem'' [Modified
Dirichlet problem for equations and systems degenerated at the
boundary], {\it Diff. uravn.} [Differ. Equ.], 1982,  {\bf 18},
No.~2, 319--324 (in Russian).

\bibitem{eng_Iva1}  L.~A.~Ivanov, ``O zadache Koshi dlya
operatorov, raspadayushchikhsya na mnozhiteli
Eylera---Puassona---Darbu'' [On Cauchy problem for operators
disintegrated into products of Euler--Poisson-Darboux
multipliers], {\it Diff. uravn.} [Differ. Equ.], 1978, {\bf 14},
No.~4, 736--739 (in Russian).

\bibitem{eng_Iva2}  L.~A.~Ivanov, ``Zadacha Koshi
dlya nekotorykh operatorov s osobennostyami'' [A Cauchy problem
for some operators with singularities], {\it Diff. uravn.}
[Differ. Equ.], 1982, {\bf 18}, No.~6, 1020--1028 (in Russian).

\bibitem{eng_20}  V.~A.~Il'in, ``Yadra drobnogo poryadka'' [Kernels of fractional
order], {\it Mat. sb.} [Math. Digest], 1957,  {\bf 41}, No.~4,
459--480 (in Russian).

\bibitem{eng_21}  K.~S.~Kazaryan, ``O
zadache Dirikhle v vesovoy metrike'' [On a Dirichlet problem in
weighted metrics], In: {\it Primenenie metodov teorii funktsiy i
funktsional'nogo analiza k zadacham matematicheskoy fiziki}
[Application of Method of Function Theory and Functional Analysis
to Problems of Mathematical Physics],  EGU, Erevan, 1982,
pp.~134--136 (in Russian).

\bibitem{eng_KaSk}  G.~A.~Kamenskiy and A.~L.~Skubachevskiy, {\it Lineynye
kraevye zadachi dlya differentsial'no-raznostnykh uravneniy}
[Linear Boundary-Value Problems for Differential--Difference
equations], MAI, Moscow, 1992 (in Russian).

\bibitem{eng_Kap1}  O.~V.~Kaptsov, {\it Metody integrirovaniya uravneniy s chastnymi
proizvodnymi} [Methods of Integration of Partial Differential
Equations], Fizmatlit, Moscow, 2009 (in Russian).

\bibitem{eng_KaSa}  N.~K.~Karapetyants and
S.~G.~Samko, {\it Uravneniya s involyutivnymi operatorami i ikh
prilozheniya} [Equations with Involution Operators and
Applications],  Rostov. un-ta, Rostov-na-Donu, 1988 (in Russian).

\bibitem{eng_KarST3}  Sh.~T.~Karimov, ``Mnogomernyy operator
Erdeyi---Kobera i ego prilozhenie k resheniyu zadachi Koshi dlya
trekhmernogo giperbolicheskogo uravneniya s singulyarnymi
koeffitsientami'' [Multidimensional Erd\'elyi--Kober operator and
its applications to solving of Cauchy problem for
three--dimensional hyperbolic equation with singular
coefficients], {\it Uzb. mat.~zh.} [Uzb. Math.~J.], 2013,  No.~1,
70--80 (in Russian).

\bibitem{eng_KarST4}  Sh.~T.~Karimov, ``Ob odnom metode resheniya zadachi Koshi
dlya obobshchennogo uravneniya Eylera---Puassona---Darbu'' [On a
method of solution of the Cauchy problem for a generalized
Euler--Poisson--Darboux equation], {\it Uzb. mat.~zh.} [Uzb.
Math.~J.], 2013, No.~3, 57--69 (in Russian).

\bibitem{eng_KarST2}  Sh.~T.~Karimov, ``Reshenie zadachi Koshi dlya
trekhmernogo giperbolicheskogo uravneniya s sin\-gu\-lyar\-ny\-mi
koeffitsientami i so spektral'nym parametrom'' [Solution of the
Cauchy problem for three-dimensional hyperbolic equation with
singular coefficients and spectral parameter], {\it Uzb. mat.~zh.}
[Uzb. Math.~J.], 2014,  No.~2, 55--65 (in Russian).

\bibitem{eng_KarST1}
Sh.~T.~Karimov, ``Ob odnom metode resheniya zadachi Koshi dlya
odnomernogo polivolnovogo uravneniya s singulyarnym operatorom
Besselya'' [On a method of solution of the Cauchy problem for
one-dimensional polywave equation with singular Bessel operator],
{\it Izv. vuzov. Ser. Mat.} [Bull. Acad. Sci. USSR. Ser. Math.],
2017,  No.~8, 27--41 (in Russian).

\bibitem{eng_Kar1}   D.~B.~Karp, ``Prostranstva s gipergeometricheskimi
vosproizvodyashchimi yadrami i drobnye pre\-o\-bra\-zo\-va\-niya tipa
Fur'e'' [Spaces with hypergeometric reproducing kernels and
fractional Fourier transforms], {\it PhD Thesis}, Vladivostok,
2000 (in Russian)

\bibitem{eng_S142}   D.~B.~Karp and S.~M.~Sitnik, ``Drobnoe
preobrazovanie Khankelya i ego prilozheniya'' [Fractional Hankel
transforms and its applications], {\it Abstr. of Voronezh. Spring
Math. School $(17--23$  Apr. $1996)$}, VGU, Voronezh, 1996, p.~92
(in Russian).

\bibitem{eng_22}   V.~V.~Katrakhov, ``O zadache na sobstvennye
znacheniya dlya singulyarnykh ellipticheskikh operatorov'' [On an
eigenvalue problem for singular elliptic operators], {\it Dokl. AN
SSSR} [Rep. Acad. Sci. USSR], 1972,  {\bf 207}, No.~2, 284--287
(in Russian).

\bibitem{eng_23}
V.~V.~Katrakhov, ``K teorii uravneniy s chastnymi proizvodnymi s
singulyarnymi koeffitsientami'' [On the theory of partial
differential equations with singular coefficients], {\it Dokl. AN
SSSR} [Rep. Acad. Sci. USSR], 1974,  {\bf 218}, No.~1, 17--20 (in
Russian).

\bibitem{eng_24}   V.~V.~Katrakhov, ``Spektral'naya funktsiya nekotorykh singulyarnykh
differentsial'nykh operatorov'' [Spectral function for some
singular differential operators], {\it Diff. uravn.} [Differ.
Equ.], 1976, {\bf 12}, No.~7, 1256--1266 (in Russian).

\bibitem{eng_25}   V.~V.~Katrakhov, ``Operatory preobrazovaniya v teorii
odnomernykh psevdodifferentsial'nykh operatorov'' [Transmutations
in the theory of one-dimensional pseudodifferential operators],
In: {\it Primenenie metodov teorii funktsiy i funktsional'nogo
analiza k zadacham matematicheskoy fiziki} [Application of Method
of Function Theory and Functional Analysis to Problems of
Mathematical Physics], IM SO AN USSR, Novosibirsk, 1979,
pp.~72--75 (in Russian).

\bibitem{eng_26}   V.~V.~Katrakhov, ``Operatory
preobrazovaniya i psevdodifferentsial'nye operatory''
[Transmutations and pseudodifferential operators], {\it Sib.
mat.~zh.} [Siberian Math.~J.], 1980,  {\bf 21}, No.~1, 86--97 (in
Russian).

\bibitem{eng_Kat1}   V.~V.~Katrakhov, ``Izometricheskie operatory
preobrazovaniya i spektral'naya funktsiya dlya odnogo klassa
odnomernykh singulyarnykh psevdodifferentsial'nykh operatorov''
[Isometric transmutations and a spectral function for a class of
one-dimensional singular pseudodifferential operators], {\it Dokl.
AN SSSR} [Rep. Acad. Sci. USSR], 1980,  {\bf  251}, No.~5,
1048--1051 (in Russian).

\bibitem{eng_28}   V.~V.~Katrakhov, ``Obshchie kraevye zadachi dlya odnogo
klassa singulyarnykh i vyrozhdayushchikhsya uravneniy'' [General
boundary--value problems for a class of singular and degenerated
equations], {\it Dokl. AN SSSR} [Rep. Acad. Sci. USSR], 1980, {\bf
251}, No.~6, 1296--1300 (in Russian).

\bibitem{eng_29}   V.~V.~Katrakhov, ``Obshchie kraevye zadachi
dlya odnogo klassa singulyarnykh i vyrozhdayushchikhsya
ellipticheskikh uravneniy'' [General boundary--value problems for
a class of singular and degenerated elliptic equations], {\it Mat.
sb.} [Math. Digest], 1980,  {\bf 112}, No.~3, 354--379 (in
Russian).

\bibitem{eng_Kat2}   V.~V.~Katrakhov, ``Ob odnoy kraevoy zadache
dlya uravneniya Puassona'' [On a boundary-value problem for the
Poisson equation], {\it Dokl. AN SSSR} [Rep. Acad. Sci. USSR],
1981, {\bf 259}, No.~5, 1041--1045 (in Russian).

\bibitem{eng_30}   V.~V.~Katrakhov, ``Singulyarnye kraevye
zadachi i operatory preobrazovaniya'' [Singular boundary-value
problems and transmutations], In: {\it Korrektnye kraevye zadachi
dlya neklassicheskikh urav\-ne\-niy matematicheskoy fiziki}
[Well-Posed Boundary-Value Problems for Nonclassical Equations of
Mathematical Physics], IM SO AN USSR, Novosibirsk, 1981,
pp.~87--91 (in Russian).

\bibitem{eng_32}   V.~V.~Katrakhov, ``Metod operatorov preobrazovaniya v teorii obshchikh
vesovykh kraevykh zadach dlya singulyarnykh i vyrozhdayushchikhsya
ellipticheskikh uravneniy s parametrom'' [Transmutation method in
the theory of general weighted boundary-value problems for
singular and degenerated elliptic equations with a parameter],
{\it Dokl. AN SSSR} [Rep. Acad. Sci. USSR], 1982,  {\bf  266},
No.~5, 1037--1040 (in Russian).

\bibitem{eng_KatDis}   V.~V.~Katrakhov, ``Singulyarnye ellipticheskie kraevye zadachi.
Metod operatorov preobrazovaniya'' [Singular elliptic
boundary--value problems.  Method of transmutations], {\it
Doctoral Thesis}, Novosibirsk, 1989 (in Russian)

\bibitem{eng_Kat3}   V.~V.~Katrakhov, ``Ob odnoy singulyarnoy kraevoy zadache dlya
uravneniya Puassona'' [On a  boundary-value problem for the
Poisson equation], {\it Mat. sb.} [Math. Digest], 1991,  {\bf
182}, No.~6, 849--876 (in Russian).

\bibitem{eng_Kat4}   V.~V.~Katrakhov, ``Singulyarnye kraevye zadachi dlya nekotorykh
ellipticheskikh uravneniy v oblastyakh s uglovymi tochkami''
[Singular boundary-value problems for some elliptic equations in
corner domains], {\it Dokl. AN SSSR} [Rep. Acad. Sci. USSR], 1991,
{\bf 316}, No.~5, 1047--1050 (in Russian).

\bibitem{eng_KaKa}   V.~V.~Katrakhov and
A.~A.~Katrakhova, ``Formula Teylora s operatorom Besselya dlya
funktsiy odnoy i dvukh peremennykh'' [The Taylor formula with
Bessel operator for functions in one or two variables], {\it Dep.
VINITI}, Voronezh, 1982, (in Russian)

\bibitem{eng_33}   V.~V.~Katrakhov and
I.~A.~Kipriyanov, ``Stepeni singulyarnogo ellipticheskogo
operatora'' [Powers of a singular elliptic operator], In: {\it
Teoriya kubaturnykh formul i prilozheniya funktsional'nogo analiza
k zadacham matematicheskoy fiziki} [Theory of Cubature Formulas
and Applications of Functional Analysis to Problems of
Mathematical Physics], Novosibirsk, 1980, pp.~60--80 (in Russian).

\bibitem{eng_S1}   V.~V.~Katrakhov and
S.~M.~Sitnik, ``Kraevaya zadacha dlya  statsionarnogo  uravneniya
Shredingera  s singulyarnym potentsialom'' [Boundary-value problem
for the steady-state Shr\"odinger equation with a singular
potential], {\it Dokl. AN SSSR} [Rep. Acad. Sci. USSR], 1984, {\bf
278}, No.~4, 797--799 (in Russian).

\bibitem{eng_S5}   V.~V.~Katrakhov and S.~M.~Sitnik, ``Metod faktorizatsii v teorii
operatorov preobrazovaniya'' [Factorization method in the
transmutation theory], In: {\it Memorial'nyy sbornik pamyati
Borisa Alekseevicha Bubnova: neklassicheskie uravneniya i
uravneniya smeshannogo tipa} [Memorial Digest to the Memory of
Boris A.~Bubnov: Nonclassical Equations and Mixed-Type Equations],
Novosibirsk, 1990, pp.~104--122 (in Russian).

\bibitem{eng_S7}   V.~V.~Katrakhov and S.~M.~Sitnik, ``Kompozitsionnyy
metod postroeniya $B$-ellipticheskikh, $B$-parabolicheskikh i
$B$-giperbolicheskikh operatorov preobrazovaniya'' [Composition
method for constructing $B$-elliptic, $B$-hyperbolic and
$B$-parabolic transmutations], {\it Dokl. RAN} [Rep. Russ. Acad.
Sci.], 1994, {\bf 337}, No.~3, 307--311 (in Russian).

\bibitem{eng_S8}   V.~V.~Katrakhov and S.~M.~Sitnik, ``Otsenki resheniy Yosta  dlya
odnomernogo  uravneniya  Shredingera s singulyarnym potentsialom''
[Estimates for Jost solutions to one--dimensional Shr\"odinger
equation with a singular potential], {\it Dokl. RAN} [Rep. Russ.
Acad. Sci.], 1995, {\bf 340}, No.~1, 18--20 (in Russian).

\bibitem{eng_Kach1}   A.~P.~Kachalov and
Ya.~V.~Kurylev, ``Metod operatorov preobrazovaniya v obratnoy
zadache rasseyaniya, odnomernyy Shtark-effekt'' [Transmutation
method for inverse scattering problem, one-dimensional Stark
effect], {\it Zap. nauch. sem. LOMI} [Notes Sci. Semin. Leningrad
Dept. Math. Inst.  Acad. Sci.], 1989,  {\bf 179}, 73--87 (in
Russian).

\bibitem{eng_35}   B.~V.~Kvyadaras, ``Reshenie zadachi Dirikhle dlya vyrozhdennogo
ellipticheskogo uravneniya'' [A solution of the Dirichlet problem
for a degenerated elliptic equation], In: {\it Differentsial'nye
uravneniya s chastnymi proizvodnymi: trudy konferentsii po
differentsial'nym uravneniyam i vychislitel'noy matematike}
[Partial Differential Equations: Proceedings of the Conference on
Differential Equations and Computational Mathematics], Nauka,
Novosibirsk, 1980, pp.~35--36 (in Russian).

\bibitem{eng_37}   M.~V.~Keldysh, ``O
razreshimosti i ustoychivosti zadachi Dirikhle'' [On solvability
and stability of the Dirichlet problem], {\it Usp. mat. nauk}
[Progr. Math. Sci.], 1941,  {\bf  8}, 171--292 (in Russian).

\bibitem{eng_Kel}   M.~V.~Keldysh, ``O nekotorykh sluchayakh vyrozhdeniya
uravneniy ellipticheskogo tipa na granitse oblasti'' [On some
cases of degeneration of elliptic equations on the boundary of a
domain], {\it Dokl. AN SSSR} [Rep. Acad. Sci. USSR], 1951,  {\bf
77}, No.~1, 181--183 (in Russian).

\bibitem{eng_KSZ}    A.~A.~Kilbas, M.~Saygo, and
V.~A.~Zhuk, ``O kompozitsii operatorov obobshchennogo drobnogo
integrirovaniya s differentsial'nym operatorom osesimmetricheskoy
teorii potentsiala'' [On a composition of generalized fraction
integration operators with a differential operator of axisymmetric
potential theory], {\it Diff. uravn.} [Differ. Equ.], 1991,  {\bf
27}, No.~9, 1640--1642 (in Russian).

\bibitem{eng_KiSk2}   A.~A.~Kilbas and O.~V.~Skoromnik, ``Reshenie
mnogomernogo integral'nogo uravneniya pervogo roda s funktsiey
Lezhandra po piramidal'noy oblasti'' [Solution of the
multidimensional integral equation of the first kind with Legendre
function at a pyramidal domain], {\it Dokl. RAN} [Rep. Russ. Acad.
Sci.], 2009,  {\bf 429}, No.~4, 442--446 (in Russian).

\bibitem{eng_Kip2}   I.~A.~Kipriyanov, ``Preobrazovaniya Fur'e---Besselya i teoremy
vlozheniya dlya vesovykh klassov'' [Fourier--Bessel transforms and
embedding theorems for weighted classes], {\it Tr. MIAN} [Proc.
Math. Inst. Russ. Acad. Sci.], 1967,  {\bf  89}, 130--213 (in
Russian).

\bibitem{eng_39}   I.~A.~Kipriyanov, ``Kraevye zadachi dlya singulyarnykh ellipticheskikh
operatorov v chastnykh pro\-iz\-vod\-nykh'' [Boundary--value problems
for singular elliptic partial differential operators], {\it Dokl.
AN SSSR} [Rep. Acad. Sci. USSR], 1970,  {\bf  195}, No.~1, 32--35
(in Russian).

\bibitem{eng_40}   I.~A.~Kipriyanov, ``Ob odnom klasse singulyarnykh
ellipticheskikh operatorov'' [On a class of singular elliptic
operators], {\it Diff. uravn.} [Differ. Equ.], 1971,  {\bf  7},
No.~11, 2065--2077 (in Russian).

\bibitem{eng_41}   I.~A.~Kipriyanov, ``Ob odnom klasse singulyarnykh
ellipticheskikh uravneniy'' [On a class of singular elliptic
equations], {\it Sib. mat.~zh.} [Siberian Math.~J.], 1973, {\bf
14}, No.~3, 560--568 (in Russian).

\bibitem{eng_Kip1}   I.~A.~Kipriyanov, {\it Singulyarnye ellipticheskie kraevye
zadachi} [Singular Elliptic Boundary-Value Problems],
Nauka---Fizmatlit, Moscow, 1997 (in Russian).

\bibitem{eng_KipIv1}   I.~A.~Kipriyanov and
L.~A.~Ivanov, ``O lakunakh dlya nekotorykh klassov uravneniy s
osobennostyami'' [On lacunas for some classes of singular
equations], {\it Mat. sb.} [Math. Digest], 1979,  {\bf  110},
No.~2, 235--250 (in Russian).

\bibitem{eng_42}   I.~A.~Kipriyanov and
L.~A.~Ivanov, ``Uravnenie Eylera---Puassona---Darbu v rimanovom
prostranstve'' [Euler--Poisson--Darboux equation in the Riemann
space], {\it Dokl. AN SSSR} [Rep. Acad. Sci. USSR], 1981,  {\bf
260}, No.~4, 790--794 (in Russian).

\bibitem{eng_KipIv3}
I.~A.~Kipriyanov and L.~A.~Ivanov, ``Zadacha Koshi dlya uravneniya
Eylera---Puassona---Darbu v odnorodnom simmetricheskom rimanovom
prostranstve.~I'' [Cauchy problem for the  Euler--Poisson--Darboux
equation in homogeneous symmetric Riemann space.~I], {\it Tr.
MIAN} [Proc. Math. Inst. Russ. Acad. Sci.], 1984,  {\bf 170},
139--147 (in Russian).

\bibitem{eng_KipIv2}   I.~A.~Kipriyanov and L.~A.~Ivanov, ``Zadacha
Koshi dlya uravneniya Eylera---Puassona---Darbu v simmetricheskom
prostranstve'' [Cauchy problem for the Euler--Poisson--Darboux
equation in symmetric space], {\it Mat. sb.} [Math. Digest], 1984,
{\bf  124}, No.~1, 45--55 (in Russian).

\bibitem{eng_KipIv4}   I.~A.~Kipriyanov and
L.~A.~Ivanov, ``Predstavlenie Dalambera i ravnoraspredelenie
energii'' [The D'Alambert representation and energy
equipartition], {\it Diff. uravn.} [Differ. Equ.], 1990,  {\bf
26}, No.~3, 458--464 (in Russian).

\bibitem{eng_KiKa1}   I.~A.~Kipriyanov and
V.~V.~Katrakhov, ``Ob odnom klasse mnogomernykh singulyarnykh
psev\-do\-dif\-fe\-ren\-tsi\-al'\-nykh operatorov'' [On a class of
multidimensional singular pseudodifferential operators], {\it Mat.
sb.} [Math. Digest], 1977,  {\bf  104}, No.~1, 49--68 (in
Russian).

\bibitem{eng_44}
  I.~A.~Kipriyanov and V.~V.~Katrakhov, ``Kraevaya zadacha dlya ellipticheskikh uravneniy
vtorogo poryadka pri nalichii osobennostey v izolirovannykh
granichnykh tochkakh'' [Boundary--value problem for elliptic
equations of the second order with singularities at isolated
boundary points], {\it Dokl. AN SSSR} [Rep. Acad. Sci. USSR],
1984,  {\bf 276}, No.~2, 274--276 (in Russian).

\bibitem{eng_KiKa2}
I.~A.~Kipriyanov and V.~V.~Katrakhov, ``Ob odnoy singulyarnoy
ellipticheskoy kraevoy zadache v oblastyakh na sfere'' [On a
singular elliptic boundary--value problem at the sphere domain],
{\it Preprint IPM DVO RAN}, 1989 (in Russian).

\bibitem{eng_KiKa3}
I.~A.~Kipriyanov and V.~V.~Katrakhov, ``Singulyarnye kraevye
zadachi dlya nekotorykh ellipticheskikh uravneniy vysshikh
poryadkov'' [Singular boundary-value problems for some elliptic
higher order equations], {\it Preprint IPM DVO RAN}, 1989 (in
Russian).

\bibitem{eng_KiKa4}
I.~A.~Kipriyanov and V.~V.~Katrakhov, ``Ob odnoy kraevoy zadache
dlya ellipticheskikh uravneniy vtorogo poryadka v oblastyakh na
sfere'' [On a  boundary-value problem for elliptic equations of
the second order at the sphere domain], {\it Dokl. AN SSSR} [Rep.
Acad. Sci. USSR], 1990, {\bf  313}, No.~3, 545--548 (in Russian).

\bibitem{eng_45}   I.~A.~Kipriyanov and M.~I.~Klyuchantsev, ``O yadrakh
Puassona dlya kraevykh zadach s differentsial'nym operatorom
Besselya'' [On Poisson kernels for boundary-value problems with
Bessel differential operator], In: {\it Differentsial'nye
uravneniya s chastnymi proizvodnymi} [Partial Differential
Equations], Moscow, 1970, pp.~119--134 (in Russian).

\bibitem{eng_KipKo1}   I.~A.~Kipriyanov and V.~I.~Kononenko, ``O fundamental'nykh resheniyakh
uravneniy v chastnykh proizvodnykh s differentsial'nym operatorom
Besselya'' [On fundamental solutions for partial differential
equations with Bessel differential operator], {\it Dokl. AN SSSR}
[Rep. Acad. Sci. USSR], 1966,  {\bf 170}, No.~2, 261--264 (in
Russian).

\bibitem{eng_KipKo2}
I.~A.~Kipriyanov and V.~I.~Kononenko, ``Fundamental'nye resheniya
$B$-ellipticheskikh uravneniy'' [Fundamental solutions for
$B$-elliptic equations], {\it Diff. uravn.} [Differ. Equ.], 1967,
{\bf 3}, No.~1, 114--129 (in Russian).

\bibitem{eng_KipKo3}
I.~A.~Kipriyanov and V.~I.~Kononenko, ``O fundamental'nykh
resheniyakh nekotorykh singulyarnykh uravneniy v chastnykh
proizvodnykh'' [On fundamental solutions for some singular partial
differential equations], {\it Diff. uravn.} [Differ. Equ.], 1969,
{\bf 5}, No.~8, 1470--1483 (in Russian).

\bibitem{eng_KipKu1}   I.~A.~Kipriyanov and A.~A.~Kulikov,
``Fundamental'nye resheniya $B$-gipoellipticheskikh uravneniy''
[Fundamental solutions for $B$-hypoelliptic equations], {\it Diff.
uravn.} [Differ. Equ.], 1991, {\bf 27}, No.~8, 1387--1395 (in
Russian).

\bibitem{eng_Kli2}   S.~B.~Klimentov, ``Klassy
Khardi obobshchennykh analiticheskikh funktsiy'' [Hardy classes of
generalized analytic functions], {\it Izv. vuzov. Sev.-Kavkaz.
reg. Ser. Estestv. nauki} [Bull. Higher Edu. Inst. North Caucas.
Reg. Nat. Sci.], 2003, No.~3, 6--10 (in Russian).

\bibitem{eng_Kli3}   S.~B.~Klimentov, ``Klassy Smirnova obobshchennykh
analiticheskikh funktsiy'' [Smirnov classes of generalized
analytic functions], {\it Izv. vuzov. Sev.-Kavkaz. reg. Ser.
Estestv. nauki} [Bull. Higher Edu. Inst. North Caucas. Reg. Nat.
Sci.], 2005, No.~1, 13--17 (in Russian).

\bibitem{eng_Kli5}
S.~B.~Klimentov, ``Klassy VMO obobshchennykh analiticheskikh
funktsiy'' [BMO classes of generalized analytic functions], {\it
Vladikavkaz. mat.~zh.} [Vladikavkaz. Math.~J.], 2006,  {\bf  8},
No.~1, 27--39 (in Russian).

\bibitem{eng_Kli1}   S.~B.~Klimentov, {\it Granichnye svoystva
obobshchennykh analiticheskikh funktsiy} [Boundary Properties of
Generalized Analytic Functions], Izd. Yuzhnogo mat. in-ta VNTS RAN
i RSO-A, Vladikavkaz, 2014 (in Russian).

\bibitem{eng_Kly1}   M.~I.~Klyuchantsev, ``O postroenii
$r$-chetnykh resheniy singulyarnykh differentsial'nykh uravneniy''
[On construction of $r$-even solutions to singular differential
equations], {\it Dokl. AN SSSR} [Rep. Acad. Sci. USSR], 1975, {\bf
224}, No.~5, 1004--1007 (in Russian).

\bibitem{eng_Kly2}   M.~I.~Klyuchantsev, ``Integraly drobnogo
poryadka i singulyarnye kraevye zadachi'' [Integrals of fractional
order and singular boundary-value problems], {\it Diff. uravn.}
[Differ. Equ.], 1976, {\bf  12}, No.~6, 983--990 (in Russian).

\bibitem{eng_KF}   A.~Kolmogorov~N and S.~V.~Fomin, {\it Elementy
teorii funktsiy i funktsional'nogo analiza} [Elements of Function
Theory and Functional Analysis], Nauka, Moscow, 1981 (in Russian).

\bibitem{eng_Kol1}   D.~Colton and R.~Kress, {\it Metody integral'nykh uravneniy v teorii
rasseyaniya} [Integral Equation Methods in Scattering Theory],
Mir, Moscow, 1987 (Russian translation).

\bibitem{eng_Kor1}   Yu.~F.~Korobeynik, {\it Operatory
sdviga na chislovykh semeystvakh} [Shift Operators on Numerical
Families], Rostov. Univ., Rostov-na-Donu, 1983 (in Russian).

\bibitem{eng_Kor2}   Yu.~F.~Korobeynik, {\it O razreshimosti v kompleksnoy oblasti
nekotorykh obshchikh klassov lineynykh integral'nykh uravneniy}
[On Solvability in Complex Domain of Some General Classes of
Linear Integral Equations], Rostov. Univ., Rostov-na-Donu, 2005
(in Russian).

\bibitem{eng_Koch1}   A.~N.~Kochubey, ``Zadacha Koshi dlya evolyutsionnykh
uravneniy drobnogo poryadka'' [Cauchy problem of evolutionary
equations of fractional order], {\it Diff. uravn.} [Differ. Equ.],
1989, {\bf  25}, No.~8, 1359--1369 (in Russian).

\bibitem{eng_Koch2}   A.~N.~Kochubey, ``Diffuziya drobnogo
poryadka'' [Diffusion of fractional order], {\it Diff. uravn.}
[Differ. Equ.], 1990,  {\bf  26}, No.~4, 660--770 (in Russian).

\bibitem{eng_KGS}
N.~S.~Koshlyakov, E.~B.~Gliner, and M.~M.~Smirnov, {\it Uravneniya
v chastnykh proizvodnykh matematicheskoy fiziki} [Partial
Differential Equations of Mathematical Physics], Vysshaya shkola,
Moscow, 1962 (in Russian).

\bibitem{eng_Kra2}
V.~F.~Kravchenko $($ed.$)$, {\it Tsifrovaya obrabotka signalov i
izobrazheniy v radiofizicheskikh prilozheniyakh} [Digital Signal
and Image Processing in Radio Physics Applications], Fizmatlit,
Moscow, 2007 (in Russian).

\bibitem{eng_Kra1}
V.~F.~Kravchenko and V.~L.~Rvachev, {\it Algebra logiki, atomarnye
funktsii i veyvlety v fizicheskikh prilozheniyakh} [Algebra of
Logic, Atomistic Functions and Wavelets in Physical Applications],
Fizmatlit, Moscow, 2006 (in Russian).

\bibitem{eng_KPS}   S.~G.~Kreyn,
Yu.~I.~Petunin, and E.~M.~Semenov, {\it Interpolyatsiya lineynykh
operatorov} [Interpolation of Linear Operators], Nauka, Moscow,
1978 (in Russian).

\bibitem{eng_48}   L.~D.~Kudryavtsev, {\it Pryamye i obratnye teoremy vlozheniya. Prilozheniya k
resheniyu variatsionnym metodom ellipticheskikh uravneniy} [Direct
and Inverse Embedding Theorems. Applications to Solving Elliptic
Equations by Variational Method], Nauka, Moscow, 1959 (in
Russian).

\bibitem{eng_KN}   L.~D.~Kudryavtsev and S.~M.~Nikol'skiy, ``Prostranstva differentsiruemykh
funktsiy mnogikh peremennykh i teoremy vlozheniya'' [Spaces of
multivariate differentiable functions and embedding theorems],
{\it Itogi nauki i tekhn. Sovrem. probl. mat. Fundam. napravl.}
[Totals Sci. Tech. Contemp. Probl. Math. Fundam. Directions],
VINITI, Moscow, 1988, {\bf 26}, 5--157 (in Russian).

\bibitem{eng_Kuz3} N.~V.~Kuznetsov, ``O sobstvennykh funktsiyakh odnogo
integral'nogo uravneniya'' [On eigenfunctions of one integral
equation], {\it Zap. nauch. sem. LOMI} [Notes Sci. Semin.
Leningrad Dept. Math. Inst.  Acad. Sci.], 1970,  {\bf  17}, No.~3,
66--149 (in Russian).

\bibitem{eng_Kuz2}
N.~V.~Kuznetsov, ``Gipoteza Petersona dlya parabolicheskikh form
vesa nul' i gipoteza Linnika. Summy summ Kloostermana''
[Peterson's hypothesis for parabolic forms of zero weight and
Linnik's hypothesis. Kloostermann sums], {\it Mat. sb.} [Math.
Digest], 1980, {\bf 111}, No.~3, 334--383 (in Russian).

\bibitem{eng_Kuz1}   N.~V.~Kuznetsov, {\it Formuly sleda i
nekotorye ikh prilozheniya v teorii chisel} [Trace Formulas and
Some of Its Applications in Number Theory], Dal'nauka,
Vladivostok, 2003 (in Russian).

\bibitem{eng_49}   R.~Courant, {\it Uravneniya s
chastnymi proizvodnymi} [Partial Differentail Equations], Mir,
Moscow, 1979 (Russian translation).

\bibitem{eng_Kus}   A.~G.~Kusraev, {\it Mazhoriruemye operatory}
[Majorizable Operators], Nauka, Moscow, 2003 (in Russian).

\bibitem{eng_Lav1}  M.~M.~Lavrent'ev, {\it Odnomernye obratnye
zadachi matematicheskoy fiziki} [One-Dimensional Inverse Problems
of Mathematical Physics], Nauka, Novosibirsk, 1982 (in Russian).

\bibitem{eng_51}   O.~A.~Ladyzhenskaya and N.~N.~Ural'tseva, {\it Lineynye i kvazilineynye
uravneniya ellipticheskogo tipa} [Linear and Quasilinear Equations
of Elliptic Type], Nauka, Moscow, 1973 (in Russian).

\bibitem{eng_Laks1}   P.~Lax and
R.~Phillips, {\it Teoriya rasseyaniya} [Scattering Theory], Mir,
Moscow, 1971 (Russian translation).

\bibitem{eng_Laks2}   P.~Lax and
R.~Phillips, {\it Teoriya rasseyaniya dlya avtomorfnykh funktsiy}
[A Scattering Theory for Automorphic Functions], Mir, Moscow, 1979
(Russian translation).

\bibitem{eng_53}
E.~M.~Landis, {\it Uravnenie vtorogo poryadka ellipticheskogo i
parabolicheskogo tipov} [The Second Order Equations of Elliptic
and Parabolic Types], Nauka, Moscow, 1971 (in Russian).

\bibitem{eng_Lan}   E.~M.~Landis, ``Zadachi E.~M.~Landisa''
[Problems of E.\,M.~Landis], {\it Usp. mat. nauk} [Progr. Math.
Sci.], 1982,  {\bf  37}, No.~6, 278--281 (in Russian).

\bibitem{eng_54}
N.~S.~Landkof, {\it Osnovy sovremennoy teorii potentsiala}
[Foundations of Modern Potential Theory], Nauka, Moscow, 1966 (in
Russian).

\bibitem{eng_Lar6}   A.~A.~Larin, ``O spektral'nykh razlozheniyakh, otvechayushchikh
samosopryazhennym rasshireniyam nekotorykh singulyarnykh
ellipticheskikh operatorov'' [On spectral expansions corresponding
to self-adjoint extensions of some singular elliptic operators],
{\it Dokl. AN SSSR} [Rep. Acad. Sci. USSR], 1987,  {\bf 293},
No.~2, 309--312 (in Russian).

\bibitem{eng_Lar5}
A.~A.~Larin, ``O svoystvakh sobstvennykh funktsiy nekotorykh
singulyarnykh ellipticheskikh operatorov'' [On eigenfunction
properties of some singular elliptic operators], {\it Diff.
uravn.} [Differ. Equ.], 1991,  {\bf  27}, No.~5, 849--856 (in
Russian).

\bibitem{eng_Lar4}   A.~A.~Larin, ``Ob ogranichennosti stepeney samosopryazhennykh
rasshireniy singulyarnykh ellipticheskikh operatorov,
deystvuyushchikh v vesovykh klassakh'' [On boundedness of powers
of self-adjoint extensions of singular elliptic operators acting
on weighted classes], {\it Diff. uravn.} [Differ. Equ.], 1992,
{\bf 28}, No.~3, 528--529 (in Russian).

\bibitem{eng_Lar3}
A.~A.~Larin, ``O predstavlenii resheniy odnogo singulyarnogo
ellipticheskogo uravneniya vtorogo poryadka v okrestnosti uglovoy
tochki'' [On representation of solutions to a singular elliptic
equation of the second order nearby the corner point], {\it Diff.
uravn.} [Differ. Equ.], 2000,  {\bf  36}, No.~4, 566--568 (in
Russian).

\bibitem{eng_Lar2}   A.~A.~Larin, ``Ob odnoy kraevoy zadache v ploskom ugle
dlya singulyarnogo ellipticheskogo uravneniya vtorogo poryadka''
[On one boundary-value problem in a plane corner for a
second-order singular elliptic equation], {\it Diff. uravn.}
[Differ. Equ.], 2000, {\bf 36}, No.~12, 1687--1694 (in Russian).

\bibitem{eng_Lar1}   A.~A.~Larin, ``O teoreme
suzheniya na sfericheskuyu poverkhnost' dlya preobrazovaniy
Fur'e---Besselya'' [On a restriction to spherical surface theorem
for Fourier--Bessel transform], {\it Dokl. Adygskoy
$($Cherkesskoy$)$ mezhd. akad. nauk.} [Rep. Adyg. (Cherkess.) Int.
Acad. Sci.], 2014,  {\bf  16}, No.~3, 22--29 (in Russian).

\bibitem{eng_Levin2}   B.~Ya.~Levin, ``Preobrazovaniya tipa Fur'e i Laplasa pri pomoshchi
resheniy differentsial'nogo uravneniya vtorogo poryadka'' [Fourier
and Laplace type transforms by means of solutions to the second
order differential equations], {\it Dokl. AN SSSR} [Rep. Acad.
Sci. USSR], 1956,  {\bf 106}, No.~2, 187--190 (in Russian).

\bibitem{eng_Levin1}   B.~Ya.~Levin, {\it Raspredelenie
korney tselykh funktsiy} [Distribution of Zeroes of Entire
Functions], GITTL, Moscow, 1956 (in Russian).

\bibitem{eng_Lev1}   B.~M.~Levitan, {\it Razlozhenie po sobstvennym
funktsiyam differentsial'nykh uravneniy vtorogo poryadka}
[Expansions in Eigenfunctions Of Second-Order Differential
Equations], Gostekhizdat, Moscow, 1950 (in Russian).

\bibitem{eng_Lev7}   B.~M.~Levitan, ``Razlozheniya po funktsiyam
Besselya v ryady i integraly Fur'e'' [Expansions in Bessel
function series and Fourier integrals], {\it Usp. mat. nauk}
[Progr. Math. Sci.], 1951,  {\bf 6}, No.~2, 102--143 (in Russian).

\bibitem{eng_Lev4}   B.~M.~Levitan, {\it Pochti-periodicheskie
funktsii} [Almost Periodic Functions], GITTL, Moscow, 1953 (in
Russian).

\bibitem{eng_Lev2}   B.~M.~Levitan, {\it Operatory obobshchennogo sdviga i
nekotorye ikh primeneniya} [Generalized Translations and Some
Applications], GIFML, Moscow, 1962 (in Russian).

\bibitem{eng_Lev3}   B.~M.~Levitan, {\it Teoriya operatorov
obobshchennogo sdviga} [Theory of Generalized Translation
Operators], Nauka, Moscow, 1973 (in Russian).

\bibitem{eng_Lev5}
B.~M.~Levitan, {\it Obratnye zadachi Shturma---Liuvillya} [Inverse
Sturm--Lioville Problems], Nauka, Moscow, 1984 (in Russian).

\bibitem{eng_Lev8}   B.~M.~Levitan and A.~Ya.~Povzner, ``Differentsial'nye uravneniya
Shturma---Liuvillya na poluosi i teorema Plansherelya''
[Sturm--Lioville differential equations on semiaxis and the
Plancherel theorem], {\it Dokl. AN SSSR} [Rep. Acad. Sci. USSR],
1946, {\bf 52}, No.~6, 483--486 (in Russian).

\bibitem{eng_Lev6}   B.~M.~Levitan and I.~S.~Sargsyan, {\it
Operatory Shturma---Liuvillya i Diraka} [Sturm--Lioville and Dirac
Operators], Nauka, Moscow, 1988 (in Russian).

\bibitem{eng_Lei1}
M.~A.~Leyzin, ``K teoremam vlozheniya dlya odnogo klassa
singulyarnykh differentsial'nykh operatorov v poluprostranstve''
[On Embedding theorems for a class of singular differential
equations in a half-space], {\it Diff. uravn.} [Differ. Equ.],
1976, {\bf 12}, No.~6, 1073--1083 (in Russian).

\bibitem{eng_Lei2}   M.~A.~Leyzin, ``O vlozhenii nekotorykh vesovykh
klassov'' [On embedding of some weighted classes], In: {\it Metody
resheniy operatornykh uravneniy} [Methods of Solution for Operator
Equations], VGU, Voronezh, 1978, pp.~96--103 (in Russian).

\bibitem{eng_Leo}   A.~F.~Leont'ev, ``Otsenka rosta resheniya
odnogo differentsial'nogo uravneniya pri bol'shikh po modulyu
znacheniyakh parametra i ee primeneniya k nekotorym voprosam
teorii funktsiy'' [Growth estimate of one differential equation
for large in moduli values and its application to some function
theory problems], {\it Sib. mat.~zh.} [Siberian Math.~J.], 1960,
{\bf  1}, No.~3, 456--487 (in Russian).

\bibitem{eng_Ler}  M.~E.~Lerner, {\it Printsipy maksimuma dlya
uravneniy giperbolicheskogo tipa i novye svoystva funktsii Rimana}
[Maximum Principles for Hyperbolic Type Equations and New
Properties of the Riemann Function], Samar. gos. tekh. un-t,
Samara, 2001 (in Russian).

\bibitem{eng_56}   P.~I.~Lizorkin, ``Obobshchennoe liuvillevskoe
differentsirovanie i metod mul'tiplikatorov v teorii vlozheniy
klassov differentsiruemykh funktsiy'' [The generalized Liouville
differentiation and multiplier method in the theory of class
embeddings of differentiable functions], {\it Tr. MIAN} [Proc.
Math. Inst. Russ. Acad. Sci.], 1969,  {\bf  105}, 89--167 (in
Russian).

\bibitem{eng_Lis}
P.~I.~Lizorkin, ``Klassy funktsiy, postroennye na osnove
usredneniy po sferam. Sluchay prostranstv Soboleva'' [Function
classes constructed by means over spheres. Case of Sobolev
spaces], {\it Tr. MIAN} [Proc. Math. Inst. Russ. Acad. Sci.],
1990, No.~192, 122--139 (in Russian).

\bibitem{eng_57}   P.~I.~Lizorkin and S.~M.~Nikol'skiy, ``Ellipticheskoe uravnenie s
vyrozhdeniem. Variatsionnyy metod'' [An elliptic degenerated
equation. Variational method], {\it Dokl. AN SSSR} [Rep. Acad.
Sci. USSR], 1981,  {\bf  257}, No.~1, 42--45 (in Russian).

\bibitem{eng_58}    P.~I.~Lizorkin and S.~M.~Nikol'skiy, ``Ellipticheskie
uravneniya s vyrozhdeniem. Differentsial'nye svoystva resheniy''
[An elliptic degenerated equation. Differential properties of
solutions], {\it Dokl. AN SSSR} [Rep. Acad. Sci. USSR], 1981, {\bf
257}, No.~2, 278--282 (in Russian).

\bibitem{eng_59}    P.~I.~Lizorkin and
S.~M.~Nikol'skiy, ``Koertsitivnye svoystva ellipticheskogo
uravneniya s sil'nym vyrozhdeniem $($sluchay obobshchennykh
resheniy$)$'' [Coercive properties of elliptic equations with
strong degeneracy (case of generalized solutions)], {\it Dokl. AN
SSSR} [Rep. Acad. Sci. USSR], 1981, {\bf  259}, No.~1, 28--30 (in
Russian).

\bibitem{eng_62}  J.-L.~Lions and  E.~Magenes, {\it Neodnorodnye granichnye zadachi i
ikh prilozheniya} [Non-Homogeneous Boundary Value Problems and
Applications], Mir, Moscow, 1971 (Russian translation).

\bibitem{eng_Lit1}  G.~S.~Litvinchuk, {\it Kraevye zadachi
i singulyarnye integral'nye uravneniya so sdvigom} [Kraevye
zadachi i singulyarnye integral'nye uravneniya so sdvigom], Nauka,
Moscow, 1977 (in Russian).

\bibitem{eng_Luke3}  Yu.~Luke, {\it Spetsial'nye matematicheskie
funktsii i ikh approksimatsii} [Mathematical Functions and Their
Approximations], Mir, Moscow, 1980 (Russian translation).

\bibitem{eng_Lyah3}  L.~N.~Lyakhov, ``Obrashchenie $B$-potentsialov''
[Inversion of $B$-potentials], {\it Dokl. AN SSSR} [Rep. Acad.
Sci. USSR], 1991,  {\bf 321}, No.~3, 466--469 (in Russian).

\bibitem{eng_Lyah1}
L.~N.~Lyakhov, {\it Vesovye sfericheskie funktsii i potentsialy
Rissa, porozhdennye obobshchennym sdvigom} [Weighted Spherical
Functions and Riesz Potentials Generated by Generalized
Translations], VGTA, Voronezh, 1997 (in Russian).

\bibitem{eng_Lyah2}  L.~N.~Lyakhov, {\it $B$-gipersingulyarnye integraly i ikh prilozheniya k opisaniyu funktsional'nykh klassov
Kipriyanova i k integral'nym uravneniyam s $B$-potentsial'nymi
yadrami} [$B$-Hypersingular Integrals and Their Applications to
Description of Kipriyanov's Functional Classes and to Integral
Equations with $B$-potential Kernels],  LGPU, Lipetsk, 2007 (in
Russian).

\bibitem{eng_LPSh1}   L.~N.~Lyakhov, I.~P.~Polovinkin, and E.~L.~Shishkina, ``Ob
odnoy zadache I.~A.~Kipriyanova dlya singulyarnogo
ul'tragiperbolicheskogo uravneniya'' [On a problem of
I.\,A.~Kipriyanov for a singular ultrahyperbolic equation], {\it
Diff. uravn.} [Differ. Equ.], 2014,  {\bf  50}, No.~4, 516--528
(in Russian).

\bibitem{eng_LPSh2}
L.~N.~Lyakhov, I.~P.~Polovinkin, and E.~L.~Shishkina, ``Formuly
resheniya zadachi Koshi dlya singulyarnogo volnovogo uravneniya s
operatorom Besselya po vremeni'' [Formulas for solution to Cauchy
problem for singular wave equation with Bessel operators in time],
{\it Dokl. RAN} [Rep. Russ. Acad. Sci.], 2014,  {\bf  459}, No.~5,
533--538 (in Russian).

\bibitem{eng_LShFrac}   L.~N.~Lyakhov and
E.~L.~Shishkina, {\it Drobnye proizvodnye i integraly i ikh
prilozheniya} [Fractional Derivatives and Integrals and Their
Applications],
 VGU, Voronezh, 2011 (in Russian).

\bibitem{eng_Mal1}   M.~M.~Malamud, ``Ob operatorakh preobrazovaniya dlya
obyknovennykh differentsial'nykh uravneniy vysshikh poryadkov''
[On transmutations for ordinary differential equations of higher
orders], In: {\it Matematicheskiy analiz i teoriya veroyatnostey}
[Mathematical Analysis and Probability Theory], Naukova dumka,
Kiev, 1978, pp.~108--111 (in Russian).

\bibitem{eng_Mal2}   M.~M.~Malamud, ``Neobkhodimye usloviya sushchestvovaniya operatora
preobrazovaniya dlya uravneniy vysshikh poryadkov'' [Necessary
conditions for existence of transmutations for equations of higher
orders], {\it Funkts. analiz i ego prilozh.} [Funct. Anal. Appl.],
1982, {\bf  16}, No.~3, 74--75 (in Russian).

\bibitem{eng_Mal3}   M.~M.~Malamud, ``K voprosu ob
operatorakh preobrazovaniya'' [On a question  of transmutations],
{\it Preprint IM AN USSR}, Kiev, 1984 (in Russian).

\bibitem{eng_Mal4}
M.~M.~Malamud, ``Operatory preobrazovaniya dlya  uravneniy
vysshikh poryadkov'' [Transmutations for  equations of higher
orders], {\it Mat. fiz. i nelin. mekh.} [Math. Phys. Nonlinear
Mech.], 1986, No.~6, 108--111 (in Russian).

\bibitem{eng_Mal5}
M.~M.~Malamud, ``K voprosu ob operatorakh preobrazovaniya dlya
obyknovennykh differentsial'nykh uravneniy'' [On a question  of
transmutations for ordinary differential equations], {\it Tr.
Mosk. mat. ob-va} [Proc. Moscow Math. Soc.], 1990,  {\bf 53},
68--97 (in Russian).

\bibitem{eng_Marich1}  O.~I.~Marichev, {\it Metod vychisleniya integralov ot spetsial'nykh
funktsiy} [Method of Calculation of Integrals of Special
Functions], Nauka i Tekhnika, Minsk, 1978 (in Russian).

\bibitem{eng_MaKiRe}  O.~I.~Marichev,
A.~A.~Kilbas, and O.~A.~Repin, {\it Kraevye zadachi dlya uravneniy
v chastnykh proizvodnykh s razryvnymi koeffitsientami}
[Boundary-Value Problems for Partial Differential Equations with
Discontinuous Coefficients], Samar. Gos. Ekonom. Univ., Samara,
2008 (in Russian).

\bibitem{eng_Mar3}  V.~A.~Marchenko, ``Nekotorye voprosy teorii differentsial'nogo
operatora vtorogo poryadka'' [Some questions of the theory of
second-order differential operators], {\it Dokl. AN SSSR} [Rep.
Acad. Sci. USSR], 1950,  {\bf 72}, No.~3, 457--460 (in Russian).

\bibitem{eng_Mar4}  V.~A.~Marchenko, ``Operatory preobrazovaniya''
[Transmutation operators], {\it Dokl. AN SSSR} [Rep. Acad. Sci.
USSR], 1950,  {\bf  74}, No.~2, 185--188 (in Russian).

\bibitem{eng_Mar5}
V.~A.~Marchenko, ``O formulakh obrashcheniya, porozhdaemykh
lineynym differentsial'nym operatorom vtorogo poryadka'' [On
inversion formulas generated by the second order differential
operator], {\it Dokl. AN SSSR} [Rep. Acad. Sci. USSR], 1950, {\bf
74}, No.~4, 657--660 (in Russian).

\bibitem{eng_Mar6}  V.~A.~Marchenko, ``Nekotorye voprosy teorii odnomernykh
differentsial'nykh operatorov vtorogo poryadka.~I'' [Some
questions of the theory of one-dimensional  second-order
differential operators.~I], {\it Tr. Mosk. mat. ob-va} [Proc.
Moscow Math. Soc.], 1952, {\bf 1}, 327--420 (in Russian).

\bibitem{eng_Mar7}  V.~A.~Marchenko, ``Nekotorye voprosy
teorii odnomernykh differentsial'nykh operatorov vtorogo
poryadka.~II'' [Some questions of the theory of one-dimensional
second-order differential operators.~II], {\it Tr. Mosk. mat.
ob-va} [Proc. Moscow Math. Soc.], 1953,  {\bf 2}, 3--82 (in
Russian).

\bibitem{eng_Mar1}  V.~A.~Marchenko, {\it Spektral'naya teoriya
operatorov Shturma---Liuvillya} [Spectral Theory of
Sturm--Liouville Operators], Naukova Dumka, Kiev, 1972 (in
Russian).

\bibitem{eng_Mar2}  V.~A.~Marchenko, {\it Operatory
Shturma---Liuvillya i ikh prilozheniya} [Sturm--Liouville
Operators and Their Applications], Naukova Dumka, Kiev, 1977 (in
Russian).

\bibitem{eng_Mar8}  V.~A.~Marchenko, {\it Nelineynye uravneniya i
operatornye algebry} [Nonlinear Equations and Operator Algebras],
Naukova Dumka, Kiev, 1986 (in Russian).

\bibitem{eng_Mar9}
V.~A.~Marchenko, ``Obobshchennyy sdvig, operatory preobrazovaniya
i obratnye zadachi'' [Generalized tran\-s\-la\-tion, transmutations and
inverse problems], In: {\it Matematicheskie sobytiya KHKH veka}
[Mathe\-ma\-ti\-cal Events in XX Century], Fazis, Moscow, 2003 (in
Russian).

\bibitem{eng_Mat1}
M.~I.~Matiychuk, {\it Parabolichni singulyarni krayovi zadachi}
[Parabolic Singular Boundary-Value Problems], In-t matematiki NAN
Ukra\"\i ni, Ki\"\i v, 1999 (in Russian).

\bibitem{eng_Mat2}   M.~I.~Matiychuk, {\it Parabolіchnі ta elіptichnі
krayovі zadachі z osoblivostyami} [Parabolic and Elliptic
Boundary-Value Problems with Singularities], Prut, Chernіvtsі,
2003 (in Russian).

\bibitem{eng_Mats}
V.~I.~Matsaev, ``O sushchestvovanii operatora preobrazovaniya dlya
differentsial'nykh uravneniy vysshikh poryadkov'' [On an existence
of transmutations for higher order differential equations], {\it
Dokl. AN SSSR} [Rep. Acad. Sci. USSR], 1960,  {\bf 130}, No.~3,
499--502 (in Russian).

\bibitem{eng_S24}  Kh.~Mekhrez and S.~M.~Sitnik, ``Monotonnost' otnosheniy nekotorykh
gipergeometricheskikh funktsiy'' [Monotonicity of ratios for some
hypergeometric functions], {\it Sib. elektron. mat. izv.}
[Siberian Electron Math. Bull.], 2016,  {\bf  13}, 260--268 (in
Russian).

\bibitem{eng_Mesh1}   V.~Z.~Meshkov, ``Vesovye differentsial'nye neravenstva
i ikh primenenie dlya otsenok skorosti ubyvaniya na beskonechnosti
resheniy ellipticheskikh uravneniy vtorogo poryadka'' [Weighted
differential inequalities and its applications for estimates of
decay rate for second-order elliptic equations], {\it Tr. MIAN}
[Proc. Math. Inst. Russ. Acad. Sci.], 1989,  {\bf 190}, 139--158
(in Russian).

\bibitem{eng_Mesh2}   V.~Z.~Meshkov, ``O vozmozhnoy skorosti ubyvaniya na beskonechnosti
resheniy uravneniy v chastnykh proizvodnykh vtorogo poryadka'' [On
a possible rate of decay at infinity for solutions to partial
differential second order equations], {\it Mat. sb.} [Math.
Digest], 1991,  {\bf 182}, No.~3, 364--383 (in Russian).

\bibitem{eng_Me1}  V.~V.~Meshcheryakov, ``Differentsial'no-raznostnye operatory, assotsiirovannye s sistemami
korney kokseterovskogo tipa'' [Differential-difference operators
associated with kernel systems of Coxeter type], {\it PhD Thesis},
Kolomna, 2008 (in Russian).

\bibitem{eng_Miz}  S.~Mizohata, {\it Teoriya
uravneniy s chastnymi proizvodnymi} [The Theory of Partial
Differential Equations], Mir, Moscow, 1977 (Russian translation).

\bibitem{eng_Miran}
C.~Miranda, {\it Uravneniya s chastnymi proizvodnymi
ellipticheskogo tipa} [Partial Differential Equations of Elliptic
Type], Mir, Moscow, 1957 (Russian translation).

\bibitem{eng_Moi}  E.~I.~Moiseev, {\it Uravneniya smeshannogo tipa so spektral'nym
parametrom} [Mixed Type Equations with Spectral Parameter],  MGU,
Moscow, 1988 (in Russian).

\bibitem{eng_Mur9}   A.~B.~Muravnik, ``O
stabilizatsii resheniy nekotorykh singulyarnykh kvazilineynykh
parabolicheskikh zadach'' [On stabilization of solutions of some
singular quasilinear parabolic problems], {\it Mat. zametki}
[Math. Notes], 2003,  {\bf  74}, No.~6, 858--865 (in Russian).

\bibitem{eng_Mur8}
A.~B.~Muravnik, ``O stabilizatsii resheniy singulyarnykh
ellipticheskikh uravneniy'' [On stabilization of solutions of
singular elliptic equations], {\it Fundam. i prikl. mat.} [Fundam.
Appl. Math.], 2006,  {\bf  12}, No.~4, 169--186 (in Russian).

\bibitem{eng_Mur}
A.~B.~Muravnik, ``Funktsional'no-differentsial'nye parabolicheskie
uravneniya: integral'nye predstavleniya i kachestvennye svoystva
resheniy zadachi Koshi'' [Functional-differential parabolic
equations: integral representations and quantitative properties of
solutions to Cauchy problems], {\it Sovrem. mat. Fundam. napravl.}
[Contemp. Math. Fundam. Directions], 2014,  {\bf 52}, 3--141 (in
Russian).

\bibitem{eng_Mys}   A.~D.~Myshkis, {\it Lineynye
differentsial'nye uravneniya s zapazdyvayushchim argumentom}
[Linear Differential Equations with Delayed Argument],
Gostekhizdat, Moscow--Leningrad, 1951 (in Russian).

\bibitem{eng_Nai}  M.~A.~Naymark, {\it Lineynye differentsial'nye
operatory} [Linear Differential Operators], Nauka, Moscow, 1969
(in Russian).

\bibitem{eng_65}  R.~Narasimhan, {\it Analiz na deystvitel'nykh i kompleksnykh
mnogoobraziyakh} [Analysis on Real and Complex Manifolds], Mir,
Moscow, 1971 (Russian translation).

\bibitem{eng_Nat1}  F.~Natterer, {\it Matematicheskie
aspekty komp'yuternoy tomografii} [The Mathematics of Computerized
Tomography], Mir, Moscow, 1990 (Russian translation).

\bibitem{eng_Nah1}
A.~M.~Nakhushev, {\it Uravneniya matematicheskoy biologii}
[Equations of Mathematical Biology], Vysshaya Shkola, Moscow, 1995
(in Russian).

\bibitem{eng_Nah2}  A.~M.~Nakhushev, {\it Elementy drobnogo
ischisleniya i ikh primenenie} [Elements of Fractional Calculus
and Their Application], Nal'chik, 2000 (in Russian).

\bibitem{eng_Nah3}  A.~M.~Nakhushev, {\it Drobnoe ischislenie i ego
primenenie} [Fractional Calculus and Its Applications], Fizmatlit,
Moscow, 2003 (in Russian).

\bibitem{eng_Nah4}  A.~M.~Nakhushev, {\it Nagruzhennye
uravneniya i ikh primeneniya} [Loaded Equations and Their
Applications], Nauka, Moscow, 2012 (in Russian).

\bibitem{eng_Nizh1}  L.~P.~Nizhnik, {\it Obratnaya nestatsionarnaya
zadacha teorii rasseyaniya} [Inverse Nonstationary Scattering
Theory Problem], Naukova Dumka, Kiev, 1973 (in Russian).

\bibitem{eng_Nizh2}  L.~P.~Nizhnik, {\it Obratnye zadachi rasseyaniya dlya giperbolicheskikh
uravneniy} [Inverse Scattering  Problems for Hyperbolic
Equations], Naukova dumka, Kiev, 1990 (in Russian).

\bibitem{eng_66}   S.~M.~Nikol'skiy, {\it Priblizhenie funktsiy
mnogikh peremennykh i teoremy vlozheniya} [Approximation of
Multivariate Functions and Embedding Theorems], Nauka, Moscow,
1977 (in Russian).

\bibitem{eng_67}   S.~M.~Nikol'skiy, ``Variatsionnaya problema dlya uravneniya
ellipticheskogo tipa s vyrozhdeniem na granitse'' [Variational
problem for  an elliptic equation with degeneration at the
boundary], {\it Tr. MIAN} [Proc. Math. Inst. Russ. Acad. Sci.],
1979,  {\bf  150}, 212--238 (in Russian).

\bibitem{eng_68}   S.~M.~Nikol'skiy and P.~I.~Lizorkin, ``O nekotorykh
neravenstvakh dlya funktsiy iz vesovykh klassov i kraevykh
zadachakh s sil'nym vyrozhdeniem na granitse'' [On some
inequalities for functions from weighted classes and
boundary-value problems with strong degeneracy at the boundary],
{\it Dokl. AN SSSR} [Rep. Acad. Sci. USSR], 1964,  {\bf 159},
No.~3, 512--515 (in Russian).

\bibitem{eng_Nov}   O.~A.~Novozhenova, {\it Biografiya i nauchnye trudy Alekseya Nikiforovicha
Gerasimova. O lineynykh operatorakh, uprugo-vyazkosti, elevteroze
i drobnykh proizvodnykh} [Biography and Scientific Works of Alexey
Nikiforovich Gerasimov. On Linear Operatos, Viscoelasticity,
Eleutherosis and Fractional Derivatives], Pero, Moscow, 2018 (in
Russian).

\bibitem{eng_69}   A.~A.~Novruzov, ``O zadachakh Dirikhle dlya
ellipticheskikh uravneniy vtorogo poryadka'' [On Dirichlet
problems for second-order elliptic equations], {\it Dokl. AN SSSR}
[Rep. Acad. Sci. USSR], 1979,  {\bf  246}, No.~1, 11--14 (in
Russian).

\bibitem{eng_Nogin}  V.~A.~Nogin and E.~V.~Sukhinin, ``Obrashchenie i
opisanie giperbolicheskikh potentsialov s $L_p$-plotnostyami''
[Inversion and description of hyperbolic potentials with
$L_p$-densities], {\it Dokl. RAN} [Rep. Russ. Acad. Sci.], 1993,
{\bf 329}, No.~5, 550--552 (in Russian).

\bibitem{eng_New}   A.~Newell, {\it Solitony v
matematike i fizike} [Solitons in Mathematics and Physics], Mir,
Moscow, 1989 (Russian translation).

\bibitem{eng_OlRa}   O.~A.~Oleynik and
E.~V.~Radkevich, {\it Uravneniya s neotritsatel'noy
kharakteristicheskoy formoy} [Equations with Nonnegative
Characteristic Form], MGU, Moscow, 2010 (in Russian).

\bibitem{eng_Ome}   A.~V.~Omel'chenko, {\it Metody integral'nykh preobrazovaniy v
zadachakh matematicheskoy fiziki} [Methods of Integral Transforms
in Problems of Mathematical Physics], MTSNMO, Moscow, 2010 (in
Russian).

\bibitem{eng_Os}
V.~F.~Osipov, {\it Pochti periodicheskie funktsii Bora---Frenelya}
[Almost Periodic Bohr--Fresnel Functions], S.-Peterb. Univ.,
Saint-Petersburg, 1992 (in Russian).

\bibitem{eng_Pas1}   A.~E.~Pasenchuk, {\it Abstraktnye singulyarnye operatory}
[Abstract Singular Operators], Novocherkassk, 1993 (in Russian).

\bibitem{eng_Plat2}   S.~S.~Platonov, ``Obobshchennye sdvigi
Besselya i nekotorye zadachi teorii priblizheniya funktsiy v
metrike $L_2$.~1'' [Generalized Bessel shifts and some problems of
function approximation theory in metric $L_2$.~1], {\it Tr.
PetrGU. Ser. Mat.} [Proc. Petr. State Univ. Ser. Math.], 2000,
{\bf 7}, 70--82 (in Russian).

\bibitem{eng_Plat3}
S.~S.~Platonov, ``Obobshchennye sdvigi Besselya i nekotorye
zadachi teorii priblizheniya funktsiy v metrike $L_2$.~2''
[Generalized Bessel shifts and some problems of function
approximation theory in metric $L_2$.~2], {\it Tr. PetrGU. Ser.
Mat.} [Tr. PetrGU. Ser. Mat.], 2001,  {\bf 8}, 20--36 (in
Russian).

\bibitem{eng_Plat1}   S.~S.~Platonov, ``Garmonicheskiy analiz Besselya i priblizhenie
funktsiy na polupryamoy'' [Bessel harmonic analysis and
approximation of functions on a semiaxis], {\it Izv. RAN. Ser.
Mat.} [Bull. Russ. Acad. Sci. Ser. Math.], 2007,  {\bf  71},
No.~5, 149--196 (in Russian).

\bibitem{eng_Povz}   A.~Ya.~Povzner, ``O differentsial'nykh uravneniyakh tipa
Shturma---Liuvillya na poluosi'' [On Sturm--Liouville differential
equations on a semiaxis], {\it Mat. sb.} [Math. Digest], 1948,
{\bf 23}, No.~1, 3--52 (in Russian).

\bibitem{eng_Pol1}   G.~N.~Polozhiy, {\it Uravneniya matematicheskoy fiziki}
[Equations of Mathematical Physics], Vysshaya shkola, Moscow, 1964
(in Russian).

\bibitem{eng_Pol2}   G.~N.~Polozhiy, {\it Obobshchenie teorii analiticheskikh
funktsiy kompleksnogo peremennogo. $P$-ana\-li\-ti\-ches\-kie i
$(P,Q)$-analiticheskie funktsii i nekotorye ikh primeneniya} [A
Generalization of the Theory of Analytic Functions of Complex
Variable. $P$-Analytic and $(P,Q)$-Analytic Functions and Some
Their Applications], KGU, Kiev, 1965 (in Russian).

\bibitem{eng_Pol3}
G.~N.~Polozhiy, {\it Teoriya i primenenie $p$-analiticheskikh
funktsiy} [Theory and Applications of $p$-Analytic Functions],
Naukova dumka, Kiev, 1973 (in Russian).

\bibitem{eng_ZaiPol}   A.~D.~Polyanin and V.~F.~Zaytsev, {\it Cpravochnik po nelineynym
uravneniyam matematicheskoy fiziki: tochnye resheniya} [Hand-Book
on Nonlinear Equations of Mathematical Physics: Exact Solutions],
Fizmatlit, Moscow, 2002 (in Russian).

\bibitem{eng_Pol}   V.~T.~Polyatskiy, ``O svoystvakh resheniy nekotorogo uravneniya''
[On properties of solution of some equation], {\it Usp. mat. nauk}
[Progr. Math. Sci.], 1965, {\bf 17}, No.~4, 119--124 (in Russian).

\bibitem{eng_71}
Z.~Presdorf, {\it Nekotorye klassy singulyarnykh uravneniy} [Some
Classes of Singular Equations], Mir, Moscow, 1979 (in Russian).

\bibitem{eng_PBM123}   A.~P.~Prudnikov, Yu.~A.~Brychkov, and
O.~I.~Marichev, {\it Integraly i ryady. T.~1,~2,~3} [Integrals and
Series. V.~1,~2,~3], Nauka, Moscow, 1981, 1983, 1986 (in Russian).

\bibitem{eng_PBM}
A.~P.~Prudnikov, Yu.~A.~Brychkov, and O.~I.~Marichev,
``Vychislenie integralov i preobrazovanie Mellina'' [Calculation
of integrals and the Mellin transform], {\it Itogi nauki i tekhn.
Mat. analiz.} Totals Sci. Tech. Math. Anal.], 1989, {\bf 27},
3--146 (in Russian).

\bibitem{eng_Pshu2}   A.~V.~Pskhu, ``Integral'nye preobrazovaniya s funktsiey Rayta v
yadre'' [Integral transforms with the Wright function in the
kernel], {\it Dokl. Adygskoy (Cherkesskoy) mezhd. akad. nauk}
[Rep. Adyg. (Cherkess.) Int. Acad. Sci.], 2002, {\bf 6}, No.~1,
35--47 (in Russian).

\bibitem{eng_Pshu1}   A.~V.~Pskhu, {\it Kraevye zadachi dlya differentsial'nykh
uravneniy s chastnymi proizvodnymi drobnogo i kontinual'nogo
poryadka} [Boundary-Value Problems for Differential Equations with
Partial Derivatives of Fractional and Continual Order], Nal'chik,
2005 (in Russian).

\bibitem{eng_Pul1}   S.~P.~Pul'kin, ``Nekotorye kraevye zadachi dlya uravneniya
$u_{xx}\pm u_{yy}+\frac{p}{x}u_x$'' [Some boundary-value problems
for the equation $u_{xx}\pm u_{yy}+\frac{p}{x}u_x$], {\it Uch.
zap. Kuybyshev. ped. in-ta} [Sci. Notes Kuibyshev Ped. Univ.],
1958, {\bf 21}, 3--54 (in Russian).

\bibitem{eng_Pul2}   S.~P.~Pul'kin, {\it Izbrannye trudy} [Selected Works],
Univers Grupp, Samara, 2007 (in Russian).

\bibitem{eng_Pulkina}    L.~S.~Pul'kina, ``Ob odnoy neklassicheskoy
zadache dlya vyrozhdayushchegosya giperbolicheskogo uravneniya''
[On one nonclassical problem for degenerating hyperbolic
equation], {\it Izv. vuzov. Ser. Mat.} [Bull. Higher Edu. Inst.
Ser. Math.], 1991, No.~11, 48--51 (in Russian).

\bibitem{eng_Ram1}   A.~G.~Ramm, {\it Mnogomernye obratnye
zadachi teorii rasseyaniya} [Multidimensional Inverse Scattering
Problems], Mir, Moscow, 1994 (Russian translation).

\bibitem{eng_Rv3}  V.~A.~Rvachev, ``Finitnye resheniya
funktsional'no-differentsial'nykh uravneniy i ikh primeneniya''
[Finite solutions of functional differential equations and their
applications], {\it Usp. mat. nauk} [Progr. Math. Sci.], 1990,
{\bf 45}, No.~1, 77--103 (in Russian).

\bibitem{eng_Rv1}
V.~L.~Rvachev and V.~A.~Rvachev, {\it Teoriya priblizheniy i
atomarnye funktsii} [Approximation Theory and Atomic Functions],
Znanie, Moscow, 1978 (in Russian).

\bibitem{eng_Rv2}   V.~L.~Rvachev and V.~A.~Rvachev, {\it Neklassicheskie metody teorii
priblizheniy v kraevykh zadachakh} [Nonclassical Methods of
Approximation Theory in Boundary-Value Problems], Naukova dumka,
Kiev, 1979 (in Russian).

\bibitem{eng_Rep}
O.~A.~Repin, {\it Kraevye zadachi so smeshcheniem dlya uravneniy
giperbolicheskogo i smeshannogo tipov} [Boundary-Value Problems
with Shift for hyperbolic and mixed equations], Samara, 1992 (in
Russian).

\bibitem{eng_RiNa}   F.~Riesz and B.~Sz.-Nagy, {\it Lektsii po
funktsional'nomu analizu} [Functional Analysis], Mir, Moscow, 1979
(Russian translation).

\bibitem{eng_72}   Ya.~A.~Roytberg and
Z.~G.~Sheftel', ``Ob obshchikh ellipticheskikh zadachakh s sil'nym
vyrozhdeniem'' [On general elliptic problems with strong
degeneration], {\it Dokl. AN SSSR} [Rep. Acad. Sci. USSR], 1980,
{\bf 254}, No.~6, 1336--1341 (in Russian).

\bibitem{eng_Ros3}
 L.~E.~Rossovskiy, ``Ellipticheskie funktsional'no-differentsial'nye uravneniya so szhatiem
i rastyazheniem argumentov neizvestnoy funktsii'' [Elliptic
functional differential equations with contractions and extensions
of independent variables of the unknown function], {\it Sovrem.
mat. Fundam. napravl.} [Contemp. Math. Fundam. Directions], 2014,
{\bf  54}, 3--138 (in Russian).

\bibitem{eng_Ross1}  L.~E.~Rossovskiy and
L.~E.~Skubachevskiy, ``Razreshimost' i regulyarnost' resheniy
nekotorykh klassov ellipticheskikh
funktsional'no-differentsial'nykh uravneniy'' [Solvability and
regularity of solutions for some classes of elliptic functional
differential equations], {\it Itogi nauki i tekhn. Ser. Sovrem.
mat. i ee pril.} [Totals Sci. Tech. Ser. Contemp. Math. Appl.],
1999, {\bf 66}, 114--192 (in Russian).

\bibitem{eng_Rut2}  S.~Rutkauskas, ``Zadachi Dirikhle s asimptoticheskimi usloviyami dlya vyrozhdayushcheysya v tochke ellipticheskoy
sistemy.~I'' [Dirichlet problems with asymptotic conditions for an
elliptic system degenerated at a point.~I], {\it Diff. uravn.}
[Differ. Equ.], 2002, {\bf 38}, No.~3, 385--392 (in Russian).

\bibitem{eng_Rut3}  S.~Rutkauskas, ``Zadachi Dirikhle s asimptoticheskimi usloviyami dlya
vyrozhdayushcheysya v tochke ellipticheskoy sistemy.~II''
[Dirichlet problems with asymptotic conditions for an elliptic
system degenerated at a point.~II], {\it Diff. uravn.} [Differ.
Equ.], 2002, {\bf 38}, No.~5, 681--686 (in Russian).

\bibitem{eng_Rut1}  S.~Rutkauskas, ``O zadache tipa Dirikhle
dlya ellipticheskikh sistem s vyrozhdeniem na pryamoy'' [On the
Dirichlet-type problem for elliptic systems with degeneration on
an axis], {\it Mat. zametki} [Math. Notes], 2016,  {\bf 100},
No.~2, 270--278 (in Russian).

\bibitem{eng_Ryko1}   V.~S.~Ryko, ``Kompozitsionnaya
struktura ryadov Fur'e, diskretnykh preobrazovaniy Fur'e i Mellina
i vychislenie ikh summ'' [Compositional structure of the Fourier
series, discrete Fourier and Mellin transforms, and calculation of
their sums], {\it Dep. v VINITI AN SSSR} [Dep. VININI], Minsk,
1987, No.~2826--V~87 (in Russian).

\bibitem{eng_Ryko2}
V.~S.~Ryko, ``Metod summirovaniya i uluchsheniya skhodimosti
funktsional'nykh ryadov'' [A method of summation and improvement
of convergence of functional series], {\it Dep. v VINITI AN SSSR}
[Dep. VININT], Vologda, 1988, No.~3542--V~88 (in Russian).

\bibitem{eng_SaIl}
K.~B.~Sabitov and R.~R.~Il'yasov, ``Reshenie zadachi Trikomi dlya
uravneniya smeshannogo tipa s singulyarnym koeffitsientom
spektral'nym metodom'' [Solution of the Tricomi problem for a
mixed-type equation with a singular coefficient by a spectral
method], {\it Izv. vuzov. Ser. Mat.} [Bull. Higher Edu. Inst. Ser.
Math.], 2004, No.~2, 64--71 (in Russian).

\bibitem{eng_SKM}   S.~G.~Samko, A.~A.~Kilbas, and
O.~I.~Marichev, {\it Integraly i proizvodnye drobnogo poryadka i
nekotorye ikh prilozheniya} [Integrals and Derivatives of
Fractional Order and Some Their Applications], Nauka i Tekhnika,
Minsk, 1987 (in Russian).

\bibitem{eng_Sah1}   L.~A.~Sakhnovich, ``Spektral'nyy
analiz vol'terrovskikh operatorov i obratnye zadachi'' [Spectral
analysis of the Volterra operators and inverse problems], {\it
Dokl. AN SSSR} [Rep. Acad. Sci. USSR], 1957, {\bf 115}, No.~4,
666--669 (in Russian).

\bibitem{eng_Sah2}   L.~A.~Sakhnovich, ``Obratnaya zadacha
dlya differentsial'nykh operatorov poryadka $n>2$ s
analiticheskimi koeffitsientami'' [Inverse problem for
differential operators of order $n>2$ with analytic coefficients],
{\it Mat. sb.} [Math. Digest], 1958, {\bf 46}, No.~1, 61--76 (in
Russian).

\bibitem{eng_Sah3}
L.~A.~Sakhnovich, ``Neobkhodimye usloviya nalichiya operatorov
preobrazovaniya dlya uravneniya chetvertogo poryadka'' [Necessary
conditions of existence of transmutations for a fourth-order
equation], {\it Usp. mat. nauk} [Progr. Math. Sci.], 1961,  {\bf
16}, No.~5, 199--205 (in Russian).

\bibitem{eng_S75}   S.~M.~Sitnik, ``O skorosti ubyvaniya resheniy statsionarnogo uravneniya
Shredingera s potentsialom, zavisyashchim ot odnoy peremennoy''
[On the decay rate of solutions of steady-state Schr\"odinger
equation with a potential depending on one variable], In: {\it
Kraevye zadachi dlya neklassicheskikh uravneniy matematicheskoy
fiziki} [Boundary-Value Problems for Nonclassical Equations of
Mathematical Physics], Novosibirsk, 1985, pp.~139--147 (in
Russian).

\bibitem{eng_S71}   S.~M.~Sitnik, ``O skorosti ubyvaniya resheniy nekotorykh
ellipticheskikh i  ul'traellipticheskikh uravneniy'' [On the decay
rate of solutions of some elliptic and ultraelliptic equations],
{\it Dep. v VINITI} [Dep. VINITI], VGU, Voronezh, 13.11.1986,
No.~7771--V86 (in Russian).

\bibitem{eng_S72}   S.~M.~Sitnik, ``Ob unitarnykh
operatorakh preobrazovaniya'' [On unitary transmutations], {\it
Dep. v VINITI} [Dep. VINITI], VGU, Voronezh, 13.11.1986,
No.~7770--V86 (in Russian).

\bibitem{eng_S70}   S.~M.~Sitnik, ``Operatory preobrazovaniya dlya
differentsial'nogo vyrazheniya Besselya'' [Transmutations for the
Bessel differential expression], {\it Dep. v VINITI} [Dep.
VINITI], VGU, Voronezh, 23.01.1987, No.~535--V87 (in Russian).

\bibitem{eng_S2}   S.~M.~Sitnik, ``Ob odnoy pare operatorov
preobrazovaniya'' [On one pair of transmutations], In: {\it
Kraevye zadachi dlya neklassicheskikh uravneniy matematicheskoy
fiziki} [Boundary-Value Problems for Nonclassical Equations of
Mathematical Physics], Novosibirsk, 1987, pp.~168--173 (in
Russian).

\bibitem{eng_S73}   S.~M.~Sitnik, ``Metod operatorov preobrazovaniya dlya statsionarnogo uravneniya
Shredingera'' [Tran\-s\-mu\-ta\-tion method for the steady-state
Schr\"odinger equation], {\it PhD Thesis}, Voronezh, 1987 (in
Russian).

\bibitem{eng_S3}       S.~M.~Sitnik, ``O skorosti ubyvaniya resheniy
nekotorykh ellipticheskikh  i  ul'traellipticheskikh uravneniy''
[On the decay rate of solutions of some elliptic and ultraelliptic
equations], {\it Diff. uravn.} [Differ. Equ.], 1988, {\bf  24},
No.~3, 538--539 (in Russian).

\bibitem{eng_S4}   S.~M.~Sitnik, ``Operatory
preobrazovaniya dlya singulyarnykh differentsial'nykh uravneniy s
operatorom Besselya'' [Transmutations for singular differential
equations with the Bessel operator], In: {\it Kraevye zadachi dlya
neklassicheskikh uravneniy matematicheskoy fiziki} [Boundary-Value
Problems for Nonclassical Equations of Mathematical Physics],
Novosibirsk, 1989, pp.~179--185 (in Russian).

\bibitem{eng_S66}   S.~M.~Sitnik, ``Unitarnost' i
ogranichennost' operatorov Bushmana---Erdeyi nulevogo poryadka
gladkosti'' [Unitarity and boundedness of the Buschman---Erd\'elyi
operators of zero order of smoothness], {\it Preprint In-ta
avtomatiki i protsessov upravl. DVO RAN} [Preprint Inst. Automat.
Control Proc. RAS], Vladivostok, 1990 (in Russian).

\bibitem{eng_S6}
S.~M.~Sitnik, ``Faktorizatsiya i otsenki norm v vesovykh
lebegovykh prostranstvakh operatorov Bushmana---Erdeyi''
[Factorization and estimates of norms of the Buschman---Erd\'elyi
operators in weighted Lebesgue spaces], {\it Dokl. AN SSSR} [Rep.
Acad. Sci. USSR], 1991, {\bf 320}, No.~6, 1326--1330 (in Russian).

\bibitem{eng_S63}   S.~M.~Sitnik, ``Operator
preobrazovaniya  i  predstavlenie Yosta dlya uravneniya s
singulyarnym potentsialom'' [Transmutation and Jost representation
for an equation with singular potential], {\it Preprint In-ta
avtomatiki i protsessov upravl. DVO RAN} [Preprint Inst. Automat.
Control Proc. RAS], Vladivostok, 199 (in Russian).

\bibitem{eng_S61}   S.~M.~Sitnik, ``Neravenstva dlya polnykh ellipticheskikh
integralov Lezhandra'' [Inequalities for full elliptic Legendre
integrals], {\it Preprint In-ta avtomatiki i protsessov upravl.
DVO RAN} [Preprint Inst. Automat. Control Proc. RAS], Vladivostok,
1994 (in Russian).

\bibitem{eng_S9}
S.~M.~Sitnik, ``Neravenstva dlya funktsiy Besselya'' [Inequalities
for Bessel functions], {\it Dokl. RAN} [Rep. Russ. Acad. Sci.],
1995,  {\bf 340}, No.~1, 29--32 (in Russian).

\bibitem{eng_S140}   S.~M.~Sitnik, ``Metod polucheniya
posledovatel'nykh utochneniy neravenstva Koshi---Bunyakovskogo i
ego primeneniya k otsenkam spetsial'nykh funktsiy'' [A method of
consecutive refinements of the Cauchy--Bunyakovskii inequality and
its application to estimates of special functions], {\it Abstr. of
Voronezh. Vesen. Mat. Shkola ``Sovremennye metody v teorii
kraevykh zadach. Pontryaginskie chteniya-VII,''} VGU, Voronezh,
1996, p.~164 (in Russian).

\bibitem{eng_S140p}   S.~M.~Sitnik, ``Formula Teylora dlya
operatorov tipa Besselya'' [Taylor formula for Bessel-type
operators], {\it Abst. of Voronezh. Vesen. Mat. Shkola
``Sovremennye metody v teorii kraevykh zadach. Pontryaginskie
chteniya-VII,''} VGU, Voronezh, 1996, p.~102 (in Russian).

\bibitem{eng_S135}   S.~M.~Sitnik, ``O nekotorykh obobshcheniya drobnogo
integro-differen\-tsirovaniya'' [On some gene\-ra\-li\-za\-tions of
fractional integrodifferentiation], {\it Proc. Int. Rus.-Uzb.
Symp. ``Uravneniya smeshannogo tipa i rodstvennye problemy analiza
i informatiki,''} Nal'chik, 2003, pp.~86--88 (in Russian).

\bibitem{eng_S133}   S.~M.~Sitnik, ``Drobnoe integrodifferentsirovanie dlya
differentsial'nogo operatora Besselya'' [Fractional
integrodifferentiation for the Bessel differential operator], {\it
Proc. Int. Rus.-Kazakh. Symp. ``Uravneniya smeshannogo tipa i
rodstvennye problemy analiza i informatiki,''} Nal'chik, 2004,
pp.~163--167 (in Russian).

\bibitem{eng_S53}    S.~M.~Sitnik, ``Obobshcheniya neravenstv Koshi---Bunyakovskogo metodom
srednikh znacheniy i ikh prilozheniya'' [Generalization of the
Cauchy--Bunyakovskii inequalities by the means method and their
applications], {\it Chernozemnyy al'manakh nauch. issl. Ser.
Fundam. mat.} [Chernozem. Digest Sci. Study. Ser. Fundam. Math.],
2005, No.~1 (1), 3--42 (in Russian).

\bibitem{eng_S127}   S.~M.~Sitnik, ``Ob
obobshchenii formuly Khille---Tamarkina dlya rezol'venty na
sluchay operatorov drobnogo integrirovaniya Besselya'' [On a
generalization of the Hille--Tamarkin formula for the resolvent in
the case of Bessel fractional differentiation operators], {\it
Abstr. of III Int. Conf. ``Nelokal'nye kraevye zadachi i
rodstvennye problemy matematicheskoy biologii, informatiki i
fiziki,''} Nal'chik, 2006, pp.~269--270 (in Russian).

\bibitem{eng_S123}   S.~M.~Sitnik, ``Operatory drobnogo
integro-differentsirovaniya dlya differentsial'nogo operatora
Besselya'' [Operators of fractional integrodifferentiation for the
Bessel differential operator], {\it Proc. Fourth All-Rus. Sci.
Conf. ``Matematicheskoe modelirovanie i kraevye zadachi,'' V.~3},
Samara, 2007, pp.~158--160 (in Russian).

\bibitem{eng_S125}   S.~M.~Sitnik, ``Postroenie
operatorov preobrazovaniya Vekua---Erdeyi---Laundesa''
[Construction of Vekua--Erd\'elyi--Lowndes transmutations], {\it
Abstr. of Int. Conf. decicated to 100th anniversary since birthday
of Acad. I.~N.~Vekua, ``Differentsial'nye uravneniya, teoriya
funktsiy i prilozheniya,''} Novosibirsk, 2007, pp.~469--470 (in
Russian).

\bibitem{eng_S46}   S.~M.~Sitnik, ``Operatory preobrazovaniya i ikh prilozheniya'' [Transmutations and their applications], In: {\it Issledovaniya po
sovremennomu analizu i matematicheskomu modelirovaniyu} [Studies
in Contemporary Analysis and Mathematical Modelling], Vladikavkaz.
Nauch. Tsentr RAN i RSO-A, Vladikavkaz, 2008, pp.~226--293 (in
Russian).

\bibitem{eng_S14}   S.~M.~Sitnik, ``Metod faktorizatsii operatorov
preobrazovaniya v teorii differentsial'nykh uravneniy'' [The
method of factorization of transmutations in the theory of
differential equations], {\it Vestn. Samar. gos. un-ta.
Estestvennonauch. ser.} [Bull. Samar. State Univ. Ser. Nat. Sci.],
2008, No.~8/1 (67), 237--248 (in Russian).

\bibitem{eng_S45}   S.~M.~Sitnik, ``Utochneniya i obobshcheniya klassicheskikh neravenstv''
[Refinements and generalizations of classic inequalities], In:
{\it Itogi nauki. Yuzhnyy federal'nyy okrug. Ser. Mat. forum.
T.~3. Issl. po mat. anal.} [Totals Sci. Southern Fed. Distr. Ser.
Math. Forum. V.~3. Math. Anal.], Yuzhnyy Mat. Inst. VNTS RAN i RSO
Alaniya, Vladikavkaz, 2009, pp.~221--266 (in Russian).

\bibitem{eng_S18}    S.~M.~Sitnik, ``O yavnykh
realizatsiyakh drobnykh stepeney differentsial'nogo operatora
Besselya i ikh prilozheniyakh k differentsial'nym uravneniyam''
[On explicit realizations of fractional powers of the Bessel
differential operator and their applications to differential
equations], {\it Dokl. Adygskoy (Cherkesskoy) mezhd. akad. nauk}
[Rep. Adyg. (Cherkess.) Int. Acad. Sci.], 2010, {\bf 12}, No.~2,
69--75 (in Russian).

\bibitem{eng_S19}   S.~M.~Sitnik, ``O
predstavlenii v integral'nom vide resheniy odnogo
differentsial'nogo uravneniya s osobennostyami v koeffitsientakh''
[On integral representation of solutions of one differential
equation with singularities in coefficietns], {\it Vladikavkaz.
mat.~zh.} [Vladikavkaz. Math.~J.], 2010, {\bf 12}, No.~4, 73--78
(in Russian).

\bibitem{eng_S43}   S.~M.~Sitnik, ``Operator preobrazovaniya
spetsial'nogo vida dlya differentsial'nogo operatora s
sin\-gu\-lyar\-nym v nule potentsialom'' [Transmutation of special form
for a differential operator with singular at zero potential], In:
{\it Neklassicheskie uravneniya matematicheskoy fiziki}
[Nonclassical Equations of Mathematical Physics], Inst. Mat. im.
S.~L.~Soboleva SO RAN, Novosibirsk, 2010, pp.~264--278 (in
Russian).

\bibitem{eng_S103}  S.~M.~Sitnik, ``Razlichnye klassy operatorov preobrazovaniya Bushmana---Erdeyi''
[Different classes of Buschman---Erd\'elyi transmutations], {\it
Abstr. of Int. Conf. ``Differentsial'nye uravneniya i smezhnye
voprosy,''} MGU, Moskva, 2011, pp.~344--345 (in Russian).

\bibitem{eng_S92}   S.~M.~Sitnik, ``Novye kraevye zadachi s
$K$-sledom dlya resheniy s sushchestvennymi osobennostyami v
rabotakh V.~V.~Katrakhova'' [New boundary-value problems with
$K$-trace for solutions with strong singularities in works by
V.~V.~Katrakhov], {\it Abstr. of Int. Conf. ``Obratnye i
nekorrektnye zadachi matematicheskoy fiziki,''} Sib. nauch. ,
Novosibirsk, 2012, p.~439 (in Russian).

\bibitem{eng_S95}   S.~M.~Sitnik, ``Raboty V.~V.~Katrakhova po
teorii operatorov preobrazovaniya'' [Works by V.~V.~Katrakhov on
the transmutation theory], {\it Proc. of the Second Rus.-Uzb.
Symp. ``Uravneniya smeshannogo tipa i rodstvennye problemy analiza
i informatiki,''} KBR, Nal'chik, 2012, pp.~241--243 (in Russian).

\bibitem{eng_S400}   S.~M.~Sitnik, ``Obzor osnovnykh
svoystv operatorov preobrazovaniya Bushmana---Erdeyi'' [Survey of
basic properties of Buschman---Erd\'elyi transmutations], {\it
Chelyabinsk. fiz.-mat.~zh.} [Chelyabinsk Phys.-Math.~J.], 2016,
{\bf 1}, No.~4, 63--93 (in Russian).

\bibitem{eng_SitDis}  S.~M.~Sitnik, ``Primenenie
operatorov preobrazovaniya Bushmana---Erdeyi i ikh obobshcheniy v
teorii differentsial'nykh uravneniy s osobennostyami v
koeffitsientakh'' [Application of Buschman---Erd\'elyi
transmutations and their generalizations in the theory of
differential equations with singularities in coefficients], {\it
Doctoral Thesis}, Voronezh, 2016 (in Russian).

\bibitem{eng_S62}   S.~M.~Sitnik and D.~B.~Karp, ``Formuly kompozitsiy dlya integral'nykh
preobrazovaniy s  funktsiyami Besselya v yadrakh'' [Composition
formulas for integral transforms with Bessel functions in
kernels], {\it Preprint In-ta avtomatiki i protsessov upravl. DVO
RAN} [Preprint Inst. Automat. Control Proc. RAS], Vladivostok,
1993 (in Russian).

\bibitem{eng_S60}   S.~M.~Sitnik and D.~B.~Karp, ``Drobnoe preobrazovanie Khankelya i ego
prilozheniya v matematicheskoy fizike'' [Fractional Hankel
transform and its applications in mathematical physics], {\it
Preprint In-ta avtomatiki i protsessov upravl. DVO RAN} [Preprint
Inst. Automat. Control Proc. RAS], Vladivostok, 1994 (in Russian).

\bibitem{eng_S65}   S.~M.~Sitnik and G.~V.~Lyakhovetskiy, ``Formuly
kompozitsiy dlya operatorov Bushmana---Erdeyi'' [Composition
formulas for Buschman---Erd\'elyi transmutations], {\it Preprint
In-ta avtomatiki i protsessov upravl. DVO RAN} [Preprint Inst.
Automat. Control Proc. RAS], Vladivostok, 1991 (in Russian).

\bibitem{eng_S59}   S.~M.~Sitnik and
G.~V.~Lyakhovetskiy, ``Operatory preobrazovaniya
Vekua---Erdeyi---Laundesa'' [Vekua--Erd\'elyi--Lowndes
transmutations], {\it Preprint In-ta avtomatiki i protsessov
upravl. DVO RAN} [Preprint Inst. Automat. Control Proc. RAS],
Vladivostok, 1994 (in Russian).

\bibitem{eng_SitShishSemi} S.~M.~Sitnik and E.~L.~Shishkina, ``Ob odnom tozhdestve dlya
iterirovannogo vesovogo sfericheskogo srednego i ego
prilozheniyakh'' [On one identity for an iterated weighted
spherical mean and its applications], {\it Sib. elektron. mat.
izv.} [Siberian Electron Math. Bull.], 2016, {\bf 13}, 849--860
(in Russian).

\bibitem{eng_SSfiz}  S.~M.~Sitnik and
E.~L.~Shishkina, {\it Metod operatorov preobrazovaniya dlya
differentsial'nykh uravneniy s operatorami Besselya} [The
Transmutation Method for Differential Equations with Bessel
Operators], Fizmatlit, Moscow, 2018 (in Russian).

\bibitem{eng_S700}  S.~M.~Sitnik and
E.~L.~Shishkina, ``O drobnykh stepenyakh operatora Besselya na
poluosi'' [On fractional powers of the Bessel operator on a
semiaxis], {\it Sib. elektron. mat. izv.} [Siberian Electron Math.
Bull.], 2018, {\bf 15}, 1--10 (in Russian).

\bibitem{eng_Sku1}   A.~L.~Skubachevskiy, ``Neklassicheskie kraevye zadachi.~I'' [Nonclassical boundary-value problems.~I], {\it Sovrem.
mat. Fundam. napravl.} [Contemp. Math. Fundam. Directions], 2007,
{\bf 26}, 3--132 (in Russian).

\bibitem{eng_Sku4}   A.~L.~Skubachevskiy, {\it Nelokal'nye kraevye zadachi i ikh prilozheniya k issledovaniyu mnogomernykh diffuzionnykh
protsessov i protsessov termoregulyatsii zhivykh kletok. Ucheb.
posobie} [Nonlocal Boundary-Value Problems and Their Applications
to Study of Multidimensional Diffusion Processes and Processes of
Thermoregulation in Living Cells. Textbook], RUDN, Moscow, 2008
(in Russian).

\bibitem{eng_Sku2}    A.~L.~Skubachevskiy, ``Neklassicheskie kraevye zadachi.~II''
[Nonclassical boundary-value problems.~II], {\it Sovrem. mat.
Fundam. napravl.} [Contemp. Math. Fundam. Directions], 2009, {\bf
33}, 3--179 (in Russian).

\bibitem{eng_Sku3}
A.~L.~Skubachevskiy and P.~L.~Gurevich, {\it Primenenie metodov
nelineynogo funktsional'nogo analiza k nelokal'nym problemam
protsessov raspredeleniya tepla. Ucheb. posobie} [Application of
Methods of Nonlinear Functional Analysis to Nonlocal Problems of
Heat Conduction. Textbook], RUDN, Moscow, 2008 (in Russian).

\bibitem{eng_76}   M.~M.~Smirnov, {\it Vyrozhdayushchiesya ellipticheskie i
giperbolicheskie uravneniya} [Degenerating Elliptic and Hyperbolic
Equations], Nauka, Moscow, 1966 (in Russian).

\bibitem{eng_Smi}   M.~M.~Smirnov, {\it Uravneniya smeshannogo tipa}
[Equations of Mixed Type], Nauka, Moscow, 1970 (in Russian).

\bibitem{eng_77}   S.~L.~Sobolev, {\it Vvedenie v teoriyu
kubaturnykh formul} [Introduction to Theory of Cubature Formulas],
Nauka, Moscow, 1974 (in Russian).

\bibitem{eng_Sob}
S.~L.~Sobolev, {\it Uravneniya matematicheskoy fiziki} [Equations
of Mathematical Physics], Nauka, Moscow, 1992 (in Russian).

\bibitem{eng_Sol}   A.~P.~Soldatov, {\it Odnomernye singulyarnye
operatory i kraevye zadachi teorii funktsiy} [One-Dimensional
Singular Operators and Boundary-Value Problems of Function
Theory], Vysshaya shkola, Moscow, 1991 (in Russian).

\bibitem{eng_Son1}   N.~Ya.~Sonin, {\it Issledovaniya o
tsilindricheskikh  funktsiyakh i spetsial'nykh polinomakh}
[Investigations on Cylinder Functions and Special Polynomials],
Gostekhteoretizdat, Moscow, 1954 (in Russian).

\bibitem{eng_Soh1}   A.~S.~Sokhin, ``Ob odnom klasse operatorov preobrazovaniya''
[On one class of transmutations], {\it Tr. fiz.-tekh. in-ta nizk.
temp. AN USSR} [Proc. Phys.-Tech. Inst. Low Temp. Acad. Sci. Ukr.
SSR], 1969, {\bf 1}, 117--125 (in Russian).

\bibitem{eng_Soh2}   A.~S.~Sokhin, ``Obratnye zadachi rasseyaniya dlya uravneniy s
osobennost'yu'' [Inverse scattering problems for equations with a
singularity], {\it Tr. fiz.-tekh. in-ta nizk. temp. AN USSR}
[Proc. Phys.-Tech. Inst. Low Temp. Acad. Sci. Ukr. SSR], 1971,
{\bf 2}, 182--233 (in Russian).

\bibitem{eng_Soh3}   A.~S.~Sokhin, ``Obratnye zadachi rasseyaniya dlya
uravneniy s osobennostyami spetsial'nogo vida'' [Inverse
scattering problems for equations with singularities of special
form], {\it Teor. funktsiy, funkts. analiz i ikh prilozh.} [Funct.
Theory. Funct. Anal. Appl.], 1973, {\bf 17}, 36--64 (in Russian).

\bibitem{eng_Soh4}   A.~S.~Sokhin, ``O
preobrazovanii operatorov  dlya uravneniy s osobennost'yu
spetsial'nogo vida'' [On transforms of operators for equations
with singularity of special form], {\it Vestn. Khar'kov. un-ta}
[Bull. Kharkov Univ.], 1974, {\bf 113}, 36--42 (in Russian).

\bibitem{eng_Sta1}
V.~V.~Stashevskaya, ``Metod operatorov preobrazovaniya'' [The
transmutation method], {\it Dokl. AN SSSR} [Rep. Acad. Sci. USSR],
1953, {\bf 113}, No.~3, 409--412 (in Russian).

\bibitem{eng_Sta2}  V.~V.~Stashevskaya, ``Ob obratnoy
zadache spektral'nogo analiza dlya differentsial'nogo operatora s
osobennost'yu v nule'' [On an inverse problem of spectral analysis
for differential operator with singularity at zero], {\it Uch.
zap. Khar'kov. mat. ob-va} [Sci. Notes Kharkov Math. Soc.], 1957,
No.~5, 49--86 (in Russian).

\bibitem{eng_78}
E.~Stein and G.~Weiss, {\it Vvedenie v garmonicheskiy analiz na
evklidovykh prostranstvakh} [Introduction to Fourier Analysis on
Euclidean Spaces], Mir, Moscow, 1974 (Russian translation).

\bibitem{eng_Ter1}   S.~A.~Tersenov, {\it Vvedenie v teoriyu uravneniy,
vyrozhdayushchikhsya na granitse} [Introduction to Theory of
Equations Degenerating at the Boundary], NGU, Novosibirsk, 1973
(in Russian).

\bibitem{eng_Tit1}
E.~Titchmarsh, {\it Vvedenie v teoriyu integralov Fur'e}
[Introduction to the Theory of Fourier Integrals], GITTL,
Moscow--Leningrad, 1948 (Russian translation).

\bibitem{eng_79}   A.~N.~Tikhonov and A.~A.~Samarskiy, {\it Uravneniya matematicheskoy fiziki}
[Equations of Mathematical Physics], Nauka, Moscow, 1972 (in
Russian).

\bibitem{eng_Trib1}   H.~Triebel, {\it Teoriya interpolyatsii.
Funktsional'nye prostranstva. Differentsial'nye operatory}
[Inter\-po\-la\-tion Theory, Function Spaces, Differential Operators],
Mir, Moscow, 1980 (Russian translation).

\bibitem{eng_Tricomi1}   F.~Tricomi, {\it Lektsii po uravneniyam v
chastnykh proizvodnykh} [Lectures on Partial Differential
Equations], Inostr. lit., Moscow, 1957 (Russian translation).

\bibitem{eng_Tyr}   E.~E.~Tyrtyshnikov, {\it Matrichnyy analiz i
lineynaya algebra} [Matrix Analysis and Linear Algebra],
Fizmatlit, Moscow, 2007 (in Russian).

\bibitem{eng_Wit1}   E.~Whittaker and
G.~Watson, {\it Kurs sovremennogo analiza. Ch.~2.
Transtsendentnye funktsii} [A Course of Modern Analysis. V.~2.
Transcendental Functions], GIFML, Moscow, 1963 (Russian
translation).

\bibitem{eng_Usp}   S.~V.~Uspenskiy, ``O teoremakh vlozheniya dlya
vesovykh klassov'' [On embedding theorems for weight classes],
{\it Tr. MIAN} [Proc. Math. Inst. Russ. Acad. Sci.], 1961, {\bf
60}, 282--303 (in Russian).

\bibitem{eng_Fage3}
M.~K.~Fage, ``Postroenie operatorov preobrazovaniya i reshenie
odnoy problemy momentov dlya obyk\-no\-ven\-nykh lineynykh
differentsial'nykh uravneniy proizvol'nogo poryadka''
[Construction of trans\-mu\-ta\-tions and solution of one problem on
moments for ordinary linear differential equations of arbit\-rary
order], {\it Usp. mat. nauk} [Progr. Math. Sci.], 1957, {\bf 12},
No.~1, 240--245 (in Russian).

\bibitem{eng_Fage1}   M.~K.~Fage, ``Operatorno-analiticheskie
funktsii odnoy nezavisimoy pe\-re\-men\-noy'' [Operator-analytic
func\-tions of one independent variable], {\it Dokl. AN SSSR} [Rep.
Acad. Sci. USSR], 1957, {\bf 112}, No.~6, 1008--1011 (in Russian).

\bibitem{eng_Fage2}   M.~K.~Fage, ``Integral'nye predstavleniya
operatorno-analiticheskikh funktsiy odnoy nezavisimoy pe\-re\-men\-noy''
[Integral representations of operator-analytic functions of one
independent variable], {\it Dokl. AN SSSR} [Rep. Acad. Sci. USSR],
1957, {\bf 115}, No.~5, 874--877 (in Russian).

\bibitem{eng_Fage4}
M.~K.~Fage, ``Operatorno-analiticheskie funktsii odnoy nezavisimoy
peremennoy'' [Operator-analytic func\-tions of one independent
variable], {\it Tr. Mosk. mat. ob-va} [Proc. Moscow Math. Soc.],
1958,  {\bf  7}, 227--268 (in Russian).

\bibitem{eng_Fage5}
M.~K.~Fage, ``Integral'nye predstavleniya
operatorno-analiticheskikh funktsiy odnoy nezavisimoy pe\-re\-men\-noy''
[Integral representations of operator-analytic functions of one
independent variable], {\it Tr. Mosk. mat. ob-va} [Proc. Moscow
Math. Soc.], 1958,  {\bf 8}, 3--48 (in Russian).

\bibitem{eng_Fage6}   M.~K.~Fage, {\it Operatorno-analіtichnі funktsії odnієї nezalezhnoї
zmіnnoї} [Operator-Analytic Functions of One Independent
Variable], L'vov Univ., L'vov, 1959 (in Russian).

\bibitem{eng_FN}   D.~K.~Fage and
N.~I.~Nagnibida, {\it Problema ekvivalentnosti obyknovennykh
differentsial'nykh operatorov} [Problem of Equivalency of Ordinary
Differential Operators], Nauka, Novosibirsk, 1977 (in Russian).

\bibitem{eng_Fad1}   L.~D.~Faddeev, ``Obratnaya zadacha
kvantovoy teorii rasseyaniya.~1'' [Inverse problem of quantum
scattering theory.~1], {\it Usp. mat. nauk} [Progr. Math. Sci.],
1959, {\bf 14}, No.~4, 57--119 (in Russian).

\bibitem{eng_81}   L.~D.~Faddeev, ``Razlozhenie po sobstvennym funktsiyam
operatora Laplasa na fundamental'noy oblasti diskretnoy gruppy na
ploskosti Lobachevskogo'' [Expansion in eigenfunctions of the
Laplace operator on a fundamental domain of a discrete group on
the Lobachevskiy plane], {\it Tr. Mosk. mat. ob-va} [Proc. Moscow
Math. Soc.], 1967, {\bf 17}, 323--350 (in Russian).

\bibitem{eng_Fad2}
L.~D.~Faddeev, ``Obratnaya zadacha kvantovoy teorii
rasseyaniya.~2'' [Inverse problem of quantum scattering
theory.~2], {\it Itogi nauki i tekhn. Sovrem. probl. mat.} [Totals
Sci. Tech. Contemp. Probl. Math.], 1974,  {\bf 3}, 93--180 (in
Russian).

\bibitem{eng_FeIv}
V.~E.~Fedorov, ``Nelokal'naya na poluosi zadacha dlya
vyrozhdennykh evolyutsionnykh uravneniy'' [Nonlocal problem for
degenerated evolution equations on semiaxis], {\it Mat. zametki
SVFU} [Math. Notes North-East Fed. Univ.], 2015, {\bf  22}, No.~1,
35--43 (in Russian).

\bibitem{eng_FGP}
V.~E.~Fedorov, D.~M.~Gordievskikh, and M.~V.~Plekhanova,
``Uravneniya v banakhovykh prostranstvakh s vyrozhdennym
operatorom pod znakom drobnoy proizvodnoy'' [Equations in Banach
spaces with a degenerated operator under fractional
differentiation symbol], {\it Diff. uravn.} [Differ. Equ.], 2015,
{\bf 51}, No.~10, 1367--1375 (in Russian).

\bibitem{eng_Fet}   V.~G.~Fetisov, ``Operatory i uravneniya v
lokal'no ogranichennykh prostranstvakh'' [Operators and equa\-tions
in locally bounded domains], In: {\it Issledovaniya po
funktsional'nomu analizu i ego pri\-lo\-zhe\-ni\-yam} [Investigations on
Functional Analysis and Its Applications], Nauka, Moscow, 2005,
pp.~249--292 (in Russian).

\bibitem{eng_Fish}     M.~K.~Fishman, ``Ob ekvivalentnosti nekotorykh lineynykh operatorov
v analiticheskom prostranstve'' [On equivalence of some linear
operators in analytic space], {\it Mat. sb.} [Math. Digest], 1965,
{\bf 68}, No.~1, 63--74 (in Russian).

\bibitem{eng_Fri}   K.~Friedrichs, {\it Vozmushchenie spektra operatorov v gil'bertovom
prostranstve} [Perturbation of Spectra in Hilbert Space], Mir,
Moscow, 1968 (Russian translation).

\bibitem{eng_Han}   A.~Kh.~Khanmamedov, ``Operatory
preobrazovaniya dlya vozmushchennogo raznostnogo uravneniya Khilla
i ikh odno prilozhenie'' [Transmutations for perturbed Hill
difference equation and one their application], {\it Sib.
mat.~zh.} [Siberian Math.~J.], 2003, {\bf 44}, No.~4, 926--937 (in
Russian).

\bibitem{eng_HLP}
G.~H.~Hardy, J.~E.~Littlewood, and G.~P\'olya, {\it Neravenstva}
[Inequalities], Inostr. lit., Moscow, 1948 (Russian translation).

\bibitem{eng_Hach1}   I.~G.~Khachatryan, ``Ob operatorakh preobrazovaniya dlya
differentsial'nykh uravneniy vysshikh po\-ryad\-kov'' [On
transmutations for higher-order differential equations], {\it Izv.
AN Armen. SSR. Ser. Mat.} [Bull. Acad. Sci. Armen. SSR. Ser.
Math.], 1978, {\bf 13}, No.~3, 215--236 (in Russian).

\bibitem{eng_Hach2}   I.~G.~Khachatryan, ``Ob
operatorakh preobrazovaniya dlya differentsial'nykh uravneniy
vysshikh poryadkov, sokhranyayushchikh asimptotiku resheniy'' [On
transmutations for higher-order differential equations preserving
asymptotics of solutions], {\it Izv. AN Armen. SSR. Ser. Mat.}
[Bull. Acad. Sci. Armen. SSR. Ser. Math.], 1979, {\bf 14}, No.~6,
424--445 (in Russian).

\bibitem{eng_Kan}  Kan Cher Khe, ``Singulyarnye kraevye zadachi dlya uravneniy
matematicheskoy fiziki s operatorami Besselya'' [Singular
boundary-value problems for equations of mathematical physics with
Bessel operators], {\it Doctoral Thesis}, Khabarovsk, 1991 (in
Russian).

\bibitem{eng_Hale}   J.~Hale, {\it Teoriya
funktsional'no-differentsial'nykh uravneniy} [Theory of Functional
Differential Equa\-tions], Mir, Moscow, 1984 (Russian translation).

\bibitem{eng_82}   S.~Helgason, {\it Preobrazovanie Radona} [The
Radon Transform], Mir, Moscow, 1983 (Russian translation).

\bibitem{eng_Hel1}  S.~Helgason, {\it Gruppy i geometricheskiy
analiz} [Groups and Geometric Analysis], Mir, Moscow, 1987
(Russian translation).

\bibitem{eng_83}   L.~H\"ormander, {\it Lineynye differentsial'nye operatory s chastnymi
proizvodnymi} [Linear Partial Differential Operators], Mir,
Moscow, 1965 (Russian translation).

\bibitem{eng_Hor}   R.~Horn, {\it Matrichnyy analiz} [Matrix Analysis], Mir, Moscow, 1989 (Russian translation).

\bibitem{eng_Hrom1}   A.~P.~Khromov, ``Konechnomernye vozmushcheniya
vol'terrovykh operatorov'' [Finite-dimensional per\-tur\-ba\-tions of
Volterra operators], {\it Sovrem. mat. Fundam. napravl.} [Contemp.
Math. Fundam. Directions], 2004, No.~10, 3--163 (in Russian).

\bibitem{eng_Cher}   V.~A.~Chernyatin, {\it Obosnovanie metoda Fur'e v smeshannoy
zadache dlya uravneniy v chastnykh proizvodnykh} [Substantiation
of the Fourier Method in Mixed Problem for Partial Differential
Equations], MGU, Moscow, 1991 (in Russian).

\bibitem{eng_ShSa}
K.~Chadan and P.~Sabatier, {\it Obratnye zadachi v kvantovoy
teorii rasseyaniya} [Inverse Problems in Quantum Scattering
Theory], Mir, Moscow, 1980 (Russian translation).

\bibitem{eng_Sha1}   V.~P.~Shatskiy, ``O giperbolicheskikh sistemakh s dvumya ploskostyami
osobennostey'' [On hyperbolic systems with two planes of
singularities], {\it Dokl. AN SSSR} [Rep. Acad. Sci. USSR], 1978,
{\bf 242}, No.~1, 56--59 (in Russian).

\bibitem{eng_Sha2}   V.~P.~Shatskiy, ``Ob odnoy kraevoy zadache dlya singulyarnykh
simmetricheskikh sistem nechetnogo poryadka'' [On one
boundary-value problem for singular symmetric systems of odd
order], {\it Dokl. AN SSSR} [Rep. Acad. Sci. USSR], 1979, {\bf
243}, No.~4, 806--809 (in Russian).

\bibitem{eng_Sha3}   V.~P.~Shatskiy, ``O nekotorykh vyrozhdayushchikhsya
sistemakh pervogo poryadka v oblastyakh s kharakteristicheskoy
chast'yu granitsy'' [On some degenerating first-order systems in
domains with characteristic part of the boundary], {\it Dokl. AN
SSSR} [Rep. Acad. Sci. USSR], 1982, {\bf 262}, No.~6, 1332--1335
(in Russian).

\bibitem{eng_ShiR1}
E.~L.~Shishkina, ``Obobshchennaya vesovaya funktsiya $r^\gamma$''
[Generalized weight function $r^\gamma$], {\it Vestn. VGU. Ser.
Fiz. Mat.} [Bull. Voronezh State Univ. Ser. Phys. Math.], 2006,
No.~1, 215--221 (in Russian).

\bibitem{eng_ShiR3}   E.~L.~Shishkina, ``Ravenstvo
dlya interirovannykh vesovykh sfericheskikh srednikh,
porozhdennykh obobshchennym sdvigom'' [An identity for iterated
weighted spherical means generated by generalized shift], {\it
Materialy nauch. konf. ``Gertsenovskie chteniya-2013.''} [Proc.
Sci. Conf. ``Gertsenovskie chteniya-2013,'']  RGPU im.
A.~I.~Gertsena, Saint-Petersburg, 2013, {\bf 66}, 143--145 (in
Russian).

\bibitem{eng_ShiR2}   E.~L.~Shishkina, ``O svoystvakh
odnogo usrednyayushchego yadra v vesovom klasse Lebega'' [On
properties of one averaging kernel in a weighted Lebesgue class],
{\it Nauch. vedom. Belgorod. gos. un-ta. Ser. Mat. Fiz.} [Sci.
Bull. Belgorod Univ. Ser. Math. Phys.], 2016, {\bf 42}, No.~6,
12--19 (in Russian).

\bibitem{eng_ShiR4}
E.~L.~Shishkina, ``Vesovye obobshchennye funktsii, otvechayushchie
kvadratichnoy forme s kompleksnymi koeffitsientami'' [Weighted
generalized functions corresponding to quadratic form with complex
coef\-fi\-ci\-ents], {\it Chelyabinsk. fiz.-mat.~zh.} [Chelyabinsk
Phys.-Math.~J.], 2017, {\bf 2}, No.~1, 88--98 (in Russian).

\bibitem{eng_85}   I.~A.~Shishmarev, {\it Vvedenie v teoriyu
ellipticheskikh uravneniy} [Introduction to Theory of Elliptic
Equations], MGU, Moscow, 1979 (in Russian).

\bibitem{eng_Shos}   R.~Ya.~Shostak, ``Aleksey Vasil'evich
Letnikov'' [Alexey Vasilievich Letnikov], In: {\it
Isto\-riko-mate\-ma\-ti\-ches\-kie issledovaniya. Trudy seminara MGU po
istorii matematiki} [Historical-Mathematical Investigations.
Proceedings of the MSU Seminar on History of Mathematics], GITTL,
Moscow--Leningrad, 1952, {\bf 5}, 167--238 (in Russian).

\bibitem{eng_Els}   L.~E.~El'sgol'ts, {\it Vvedenie v teoriyu differentsial'nykh uravneniy s otklonyayushchimsya
argumentom} [Introduction to the Theory of Differential Equations
with Deviating Agrument], Nauka, Moscow, 1964 (in Russian).

\bibitem{eng_ElNor}   L.~E.~El'sgol'ts and
S.~B.~Norkin, {\it Vvedenie v teoriyu differentsial'nykh uravneniy
s otklonyayushchimsya argumentom} [Introduction to the Theory of
Differential Equations with Deviating Agrument], Nauka, Moscow,
1971 (in Russian).

\bibitem{eng_Erg1}   T.~G.~Ergashev, ``Chetvertyy potentsial
dvoynogo sloya dlya obobshchennogo dvuosesimmetricheskogo
uravneniya Gel'mgol'tsa'' [Fourth potential of double layer for
generalized two-axial-symmetric Helmholtz equation], {\it Vestn.
Tomsk. gos. un-ta. Mat. i mekh.} [Bull. Tomsk State Univ. Math.
Mech.], 2017, {\bf 50}, 45--56 (in Russian).

\bibitem{eng_Yurko}   V.~A.~Yurko, {\it Vvedenie v teoriyu obratnykh spektral'nykh zadach}
[Introduction to Theory of Inverse Spectral Problems], Nauka,
Moscow, 2007 (in Russian).

\bibitem{eng_87}   G.~N.~Yakovlev, ``Neogranichennye resheniya vyrozhdayushchikhsya ellipticheskikh
uravneniy'' [Unboun\-ded solutions of degenerating elliptic
equations], {\it Tr. MIAN} [Proc. Math. Inst. Russ. Acad. Sci.],
1978, {\bf 117}, 312--320 (in Russian).

\bibitem{eng_Yanu}   A.~I.~Yanushauskas, {\it Analiticheskaya teoriya ellipticheskikh uravneniy}
[Analytical Theory of Elliptic Equations], Nauka, Novosibirsk,
1979 (in Russian).

\bibitem{eng_Yarem}   O.~E.~Yaremko, {\it Metod operatorov preobrazovaniya v
zadachakh matematicheskogo modelirovaniya} [The Transmutations
Method in Problems of Mathematical Modelling], Penz. Gos. un-t,
Penza, 2012; Lambert Academic Publishing, 2012 (in Russian).

\bibitem{eng_Yar1}   V.~Ya.~Yaroslavtseva, ``Ob odnom klasse operatorov
preobrazovaniya i ikh prilozhenii k differentsial'nym
uravneniyam'' [On one class of transmutations and its application
to differential equations], {\it Dokl. AN SSSR} [Rep. Acad. Sci.
USSR], 1976, {\bf 227}, No.~4, 816--819 (in Russian).

\bibitem{eng_Yar2}   V.~Ya.~Yaroslavtseva, ``Neodnorodnaya granichnaya zadacha v poluprostranstve dlya odnogo klassa
singulyarnykh uravneniy'' [Nonhomogeneous boundary-value problem
in a half-space for one class of singular equations], Submitted to
{\it Diff. uravn.} [Differ. Equ.], 1989 (in Russian).


\bibitem{eng_AKK}   I.~Ali, V.~Kiryakova, and S.~L.~Kalla, ``Solutions of fractional
multi-order integral and differential equations using a
Poisson-type transform,'' {\it J.~Math. Anal. Appl.}, 2002,  {\bf
269}, No.~1, 172--199.

\bibitem{eng_AlKi1} F.~Almalki and V.~Kisil, ``Geometric dynamics of a harmonic
oscillator, non-admissible mother wavelets and squeezed states,''
{\it arXiv:1805.01399v1}, 2018.

\bibitem{eng_AAR}
G.~E.~Andrews,  R.~Askey, and R.~Roy, {\it Special functions},
Cambridge University Press, Cambridge, 1999.

\bibitem{eng_AK}   M.~Ya.~Antimirov, A.~A.~Kolyshkin, and
R.~Vaillancourt, {\it Applied integral transforms}, Am. Math.
Soc., Providence, 1993.

\bibitem{eng_AKdF}   P.~Appell and J.~Kampe de Feriet, {\it Fonctions hypergeometriques et hyperspheriques;
polynomes d'Hermite}, Gauthier-Villars, Paris, 1926.

\bibitem{eng_AppHyp}
W.~C.~Connett, M.-O.~Gebuhrer, and F.~L.~Schwartz $($ed.$)$, {\it
Applications of hypergroups and related measure algebras}, Am.
Math. Soc., Providence, 1995.

\bibitem{eng_Bac}   C.~Baccar, N.~B.~Hamadi, and L.~T.~Achdi, ``Inversion formulas for
Riemann--Liouville transform and its dual associated with singular
partial differential operators,'' {\it Int. J. Math. Math. Sci.},
2006, 1--26.

\bibitem{eng_Bai}
W.~N.~Bailey, {\it Generalized hypergeometric series},
Stechert-Hafner Service Agency, New York--London, 1964.

\bibitem{eng_Bajlekova1}   E.~G.~Bajlekova, ``Subordination
principle for fractional evolution equations,'' {\it Fract. Calc.
Appl. Anal.}, 2000, {\bf 3}, No.~3, 213--230.

\bibitem{eng_Bajlekova0}   E.~G.~Bajlekova, ``Fractional evolution equations in Banach spaces,'' {\it Thesis}, Technische Universiteit Eindhoven, 2001.

\bibitem{eng_Bers1}  L.~Bers, ``On a class of differential equations
in mechanics of continua,'' {\it Quart. Appl. Math.}, 1943, {\bf
5}, No.~1, 168--188.

\bibitem{eng_Bers3}   L.~Bers, ``A remark on an applications of
pseudo-analytic functions,'' {\it Amer. J.~Math.}, 1956,  {\bf
78}, No.~3, 486--496.

\bibitem{eng_Bers2}   L.~Bers and A.~Gelbart, ``On a class of functions
defined by partial differential equations,'' {\it Trans. Am. Math.
Soc.}, 1944,  {\bf 56}, 67--93.

\bibitem{eng_BuKe1}  J.~Bourgain and C.~E.~Kenig, ``On localization
in the continuous Anderson--Bernoulli model in higher dimension,''
{\it Invent. Math.}, 2005,  {\bf 161}, No.~2, 389--426.

\bibitem{eng_Bozh}   N.~Bozhinov, {\it Convolution representations of commutants and multipliers}, Publishing House
Bulgarian Acad. Sci., Sofia, 1988.

\bibitem{eng_Bra}  B.~L.~Braaksma, {\it Asymptotic
expansions and analytic continuation for a class of
Barnes-integrals}, Noordhoff, Groningen, 1963.

\bibitem{eng_BD2}  Z.~R.~Bragg and J.~W.~Dettman, ``An operator
calculus for related partial differential equations,'' {\it
J.~Math. Anal. Appl.}, 1968, {\bf 22}, No.~2, 261--271.

\bibitem{eng_BD1}  Z.~R.~Bragg and J.~W.~Dettman, ``Related problems in partial differential equations,'' {\it Bull. Am. Math. Soc.}, 1968,
{\bf 74}, No.~2, 375--378.

\bibitem{eng_BD3}  Z.~R.~Bragg and J.~W.~Dettman, ``A
class of related Dirichlet and initial value problems,'' {\it
Proc. Am. Math. Soc.}, 1969,  {\bf 21}, No.~1, 50--56.

\bibitem{eng_Bresters2}  D.~W.~Bresters, ``On a
generalized Euler--Poisson--Darboux equation,'' {\it SIAM J.~Math.
Anal.}, 1978,  {\bf 9}, No.~5, 924--934.

\bibitem{eng_Bul2}   P.~S.~Bullen, {\it Handbook of means
and their inequalities}, Kluwer Academic Publishers,
Dordrecht--London, 2003.

\bibitem{eng_BMV}   P.~S.~Bullen, D.~S.~Mitrinovi\'c, P.~M.~Vasi\'c, {\it Means and
their inequalities}, D.~Reidel, Dordrecht, 1988.

\bibitem{eng_Bus2}
R.~G.~Buschman, ``An inversion integral for a Legendre
transformation,'' {\it Amer. Math. Monthly}, 1962,  {\bf  69},
No.~4, 288--289.

\bibitem{eng_Bus1}
R.~G.~Buschman, ``An inversion integral for a general Legendre
transformation,'' {\it SIAM Review}, 1963,  {\bf  5}, No.~3,
232--233.

\bibitem{eng_CKT1}
H.~Campos, V.~V.~Kravchenko, and S.~M.~Torba, ``Transmutations,
$L$-bases and complete families of solutions of the stationary
Schr\"odinger equation in the plane,'' {\it J.~Math. Anal. Appl.},
2012,  {\bf  389}, No.~2, 1222--1238.

\bibitem{eng_Car1}
R.~W.~Carroll, {\it Transmutation and operator differential
equations}, North Holland, Amsterdam--New York--Oxford, 1979.

\bibitem{eng_Car2}   R.~W.~Carroll, {\it Transmutation, scattering theory and special functions}, North Holland, Amsterdam--New
York--Oxford, 1982.

\bibitem{eng_Car3}   R.~W.~Carroll, {\it Transmutation
theory and applications}, North Holland, Amsterdam--New York,
1986.

\bibitem{eng_Car4}   R.~W.~Carroll, {\it Topics in soliton theory}, North Holland, 1991.

\bibitem{eng_Car10}   R.~W.~Carroll, {\it Calculus revisited}, Springer, Dordrecht--Boston--London, 2002.

\bibitem{eng_CB}  R.~W.~Carroll and A.~Boumenir, ``Toward a general theory of
transmutation,'' {\it arXiv:~funct-an/9501006}, 1995.

\bibitem{eng_CSh}   R.~W.~Carroll and R.~E.~Showalter, {\it Singular and degenerate
Cauchy problems}, Academic Press, N.Y., 1976.

\bibitem{eng_CKT2}
С.~R.~Castillo-P\'erez, V.~V.~Kravchenko, and  S.~M.~Torba,
``Spectral parameter power series for perturbed Bessel
equations,'' {\it Appl. Math. Comput.}, 2013,  {\bf 220},
676--694.

\bibitem{eng_ChCPR}   K.~Chadan, D.~Cotton, L.~Paivarinta, and W.~Rundell, {\it An
introduction to inverse scattering and inverse spectral problems},
SIAM, 1997.

\bibitem{eng_Che1}  H.~Chebli, ``Op\'erateurs de translation g\'en\'eralises et
semigroupes de convolution,'' {\it Springer Lect. Notes}, 1974,
{\bf 404}, 35--59.

\bibitem{eng_Che3}  H.~Chebli, ``Positivit\'e des op\'erateurs de
<<translation g\'en\'eralises>> associ\'e \`a un op\'erateur de
Sturm--Liouville et quelques applications a l'analyse
harmonique,'' {\it Thesis}, Strasbourg, 1974.

\bibitem{eng_Che2}  H.~Chebli, ``Sur un th\`eor\'eme de Paley--Winer associ\'e \`a la
d\'ecomposition spectrale d'un op\'erateur de Sturm--Liouville sur
$(0,\infty)$,'' {\it J.~Funct. Anal.}, 1974,  {\bf 17}, 447--461.

\bibitem{eng_Che4}
H.~Chebli, ``Th\`eor\'eme de Paley--Winer associ\'e \`a un
op\'erateur diff\'erentiel singulier sur $(0,\infty)$,'' {\it
J.~Math. Pures Appl.}, 1979,  {\bf  58}, 1--19.

\bibitem{eng_CFH}  H.~Chebli, A.~Fitouhi, and M.~M.~Hamza, ``Expansion in series
of Bessel functions and transmutations for perturbed Bessel
operators,'' {\it J.~Math. Anal. Appl.}, 1994,  {\bf  181}, No.~3,
789--802.

\bibitem{eng_Che}
B.~Cheikh, ``Relations between harmonic analysis associated with
two differential operators of different orders,'' {\it J.~Comput.
Appl. Math.}, 2003,  {\bf  153}, No.~1, 61--71.

\bibitem{eng_Col1}   D.~Colton, {\it Solution of boundary value problems by the
method of integral operators}, Pitman Press, London, 1976.

\bibitem{eng_Col2}
D.~Colton, {\it Analytic theory of partial differential
equations}, Pitman Press, London, 1980.

\bibitem{eng_Cop2}   E.~T.~Copson, ``On a singular boundary value problem for an
equation of hyperbolic type,'' {\it Arch. Ration. Mech. Anal.},
1957,  {\bf 1}, 349--356.

\bibitem{eng_Cop1}   E.~T.~Copson, ``On the Riemann--Green function,'' {\it Arch. Ration. Mech. Anal.}, 1957/58,  {\bf  1}, 324--348.

\bibitem{eng_Cop4}   E.~T.~Copson, {\it Partial differential equations}, Cambridge University Press, London--New
York--Melbourne, 1975.

\bibitem{eng_Cop3}
E.~T.~Copson and A.~Erd\'elyi, ``On a partial differential
equation with two singular lines,'' {\it Arch. Ration. Mech.
Anal.}, 1958, {\bf  2}, No.~1, 76--86.

\bibitem{eng_Darboux}   G.~Darboux, {\it Le\c cons sur la th\'eorie g\'en\'erale des
surfaces et les applications g\'eom\'etriques du calcul
infinit\'esimal. Vol.~2}, Gauthier-Villars, Paris, 1915.

\bibitem{eng_DKW1}  B.~Davey,  C.~Kenig, and J.-N.~Wang, ``The Landis conjecture for variable coeffcient second-order elliptic PDEs,'' {\it Trans. Am.
Math. Soc.}, 2017,  {\bf 369}, No.~11, 8209--8237.

\bibitem{eng_Deans}
S.~R.~Deans, {\it The Radon transform and some of its
applications}, Dover, Mineola--New York, 1990.

\bibitem{eng_Del1}   J.~Delsarte, ``Sur certaines transformation
fonctionnelles relative aux \'equations lin\'eares aux
d\`eriv\'ees partielles du seconde ordre,'' {\it C.~R. Acad.
Sci.}, 1938,  {\bf  206}, 1780--1782.

\bibitem{eng_Del2}   J.~Delsarte, ``Sur une extension de la formule de Taylor,'' {\it J.
Math. Pures Appl.}, 1938,  {\bf  17}, 217--230.

\bibitem{eng_Del3}
J.~Delsarte, ``Une extension nouvelle de la th\'eory de fonction
presque p\'eriodiques de Bohr,'' {\it Acta Math.}, 1939,  {\bf
69}, 259--317.

\bibitem{eng_Del4}
J.~Delsarte, ``Hypergroupes et operateurs de permutation et de
transmutation,'' {\it Colloq. Internat. Centre Nat. Rech. Sci.},
1956,  {\bf 71}, 29--44.

\bibitem{eng_Del7}   J.~Delsarte, {\it Lectures on topics in mean periodic functions and
the two-radius theorem}, Tata Inst. Fundam. Research, Bombay,
1961.

\bibitem{eng_Del5}   J.~Delsarte and J.~L.~Lions, ``Transmutations d'op\'erateurs diff\'erentiels dans le
domaine complexe,'' {\it Comment. Math. Helv.}, 1957,  {\bf 32},
113--128.

\bibitem{eng_Del6}   J.~Delsarte and J.~L.~Lions, ``Moyennes g\'en\'eralis\'ees,'' {\it Comment. Math. Helv.}, 1959,  No.~34, 59--69.

\bibitem{eng_Diaz}
J.~B.~Diaz and  H.~F.~Weinberger, ``A solution of the singular
initial value problem for the Euler--Poisson equation,'' {\it
Proc. Am. Math. Soc.}, 1953,  {\bf 4}, 703--715.

\bibitem{eng_Dim}   I.~Dimovski, {\it Convolutional calculus}, Kluwer Acad.
Publ., Dordrecht, 1990.

\bibitem{eng_DHS}   I.~Dimovski, V.~Hristov, and M.~Sifi, ``Commutants of the
Dunkl operators in $C(\R)$,'' {\it Fract. Calc. Appl. Anal.},
2006,  {\bf  9}, No.~3, 195--213.

\bibitem{eng_DK2}  I.~Dimovski and V.~S.~Kiryakova, ``Transmutations,
convolutions and fractional powers of Bessel-type operators via
Meijer $G$-functions,'' {\it Proc. conf. Complex Anal. and Appl.,
1983, Sofia}, Varna, 1985, pp.~45--66.

\bibitem{eng_DK1}   I.~Dimovski and V.~S.~Kiryakova, ``Generalized Poisson transmutations
and corresponding representations of hyper-Bessel functions,''
{\it C.~R. Acad. Bulgare Sci.}, 1986,  {\bf  39}, No.~10, 29--32.

\bibitem{eng_Djr2}
M.~M.~Djrbashian, {\it Harmonic analysis and boundary value
problems in the complex domain}, Birkh\"auser,
Basel--Boston--Berlin, 1993.

\bibitem{eng_DKN}   P.~Drabek,
A.~Kufner and F.~Nicolosi, {\it Quasilinear elliptic equations
with degenerations and singularities}, Walter de Gruyter,
Berlin--New York, 1997.

\bibitem{eng_Dun1}
Ch.~Dunkl, ``Differential-difference operators associated to
reflection groups,'' {\it Trans. Am. Math. Soc.}, 1989,  {\bf
311}, 167--183.

\bibitem{eng_Dun2}
Ch.~Dunkl, ``Intertwining operators associated to the group S3,''
{\it Trans. Am. Math. Soc.}, 1995,  {\bf  347}, 3347--3374.

\bibitem{eng_Dun3}  Ch.~Dunkl, ``An Intertwining operator for the group B2,'' {\it
arXiv:~math.~CA/~0607823}, 2006.

\bibitem{eng_Dwo1}   B.~Dwork, {\it Generalized hypergeometric functions}, Oxford, 1990.

\bibitem{eng_Dwo2}   B.~Dwork, S.~Gerotto, and F.~J.~Sullivan, {\it Introduction to
$G$-functions}, Princeton, 1994.

\bibitem{eng_Dzr}   M.~M.~Dzrbashian, {\it Harmonic
analysis and boundary value problems in the complex domain},
Birkh\"auser, Basel--Boston--Berlin, 1993.

\bibitem{eng_EiIvKoch}   S.~D.~Eidelman,  S.~D.~Ivasyshen, and
A.~N.~Kochubei, {\it Analytic methods in the theory of
differential and pseudo-differential equations of parabolic type},
Springer, Basel, 2004.

\bibitem{eng_86}
A.~Erd\'elyi, ``On fractional integration and its application to
the Hankel transforms,'' {\it Quart. J.~Math. $($Oxford$)$}, 1940,
{\bf  11}, 293--303.

\bibitem{eng_Erd1}   A.~Erd\'elyi, ``Some applications of
fractional integration,'' {\it Boeing Sci. Res. Labor.}, 1963,
Math. Note {\bf 316}, D1-82-0286.

\bibitem{eng_Erd2}
A.~Erd\'elyi, ``An integral equation involving Legendre
functions,'' {\it SIAM Rev.}, 1964,  {\bf  12}, No.~1, 15--30.

\bibitem{eng_Erd3}   A.~Erd\'elyi, ``An
application of fractional integrals,'' {\it J.~Analyse Math.},
1965,  {\bf  14}, 113--126.

\bibitem{eng_Erd4}   A.~Erd\'elyi, ``Some integral equations involving
finite parts of divergent integrals,'' {\it Glasgow Math.~J.},
1967,  {\bf  8}, No.~1, 50--54.

\bibitem{eng_Erd5}   A.~Erd\'elyi, ``On the Euler--Poisson--Darboux
equation,'' {\it J.~Analyse Math.}, 1970,  {\bf  23}, 89--102.

\bibitem{eng_Rub3}   R.~Estrada and B.~Rubin, ``Null spaces Of Radon transforms,'' {\it arXiv:1504.03766v1}, 2015.

\bibitem{eng_Euler}   L.~Euler, ``Institutiones calculi
integralis,'' {\it Opera Omnia}, 1914,  {\bf  1}, No.~13,
212--230.

\bibitem{eng_Ext}   H.~Exton, {\it Multiple hypergeometric functions and applications}, John Wiley and Sons, New York, 1976.

\bibitem{eng_FH}   A.~Fitouhi and M.~M.~Hamza, ``A
uniform expansion for the eigenfunction of a singular second-order
differential operator,'' {\it SIAM J.~Math. Anal.}, 1990,  {\bf
21}, No.~6, 1619--1632.

\bibitem{eng_FJSS}   A.~Fitouhi, I.~Jebabli, E.~Shishkina, and S.~M.~Sitnik, ``Applications
of integral transforms composition  method to wave-type singular
differential equations and index shift transmutations,'' {\it
Electron. J. Differ. Equ.}, 2018,  {\bf 2018}, No.~130, 1--27.

\bibitem{eng_Fox}  D.~N.~Fox, ``The solution and Huygens'
principle for a singular Cauchy problem,'' {\it J.~Math. Mech.},
1959,  {\bf 8}, 197--219.

\bibitem{eng_Gad1}  A.~D.~Gadjiev, {\it Selected works}, ELM, Baku, 2003.

\bibitem{eng_Gal1}   L.~Gallardo and Kh.~Trim\'eche, ``Un analogue d'un theoreme de
Hardy pour la transformation de Dunkl,'' {\it C.~R. Math. Acad.
Sci. Paris}, 2002, {\bf  334}, 849--854.

\bibitem{eng_Gal2}   L.~Gallardo and Kh.~Trim\'eche, ``A
version of Hardy's theorem for the Dunkl transform,'' {\it J.
Aust. Math. Soc.}, 2004,  {\bf  77}, 371--385.

\bibitem{eng_Rad1}   M.~Ghergu and V.~Radulescu, {\it Singular
elliptic problems. Bifurcation and asymptotic analysis}, Oxford
University Press, New York, 2008.

\bibitem{eng_Gil2}   R.~Gilbert, {\it Function theoretic methods in partial
differential equations}, Academic Press, N.Y., 1969.

\bibitem{eng_Gil1}
R.~Gilbert, {\it Constructive methods for elliptic equations},
Springer, Berlin--Heidelberg, 1974.

\bibitem{eng_GB}   R.~Gilbert and H.~Begehr, {\it Transformations,
transmutations and kernel functions. V.~1-2}, Longman, Harlow,
1992.

\bibitem{eng_GSS}   H.~J.~Glaeske, A.~P.~Prudnikov, and K.~A.~Skornik, {\it Operational calculus and related
topics}, Chapman \& Hall/CRC, Boca Raton, 2006.

\bibitem{eng_GSPP}  J.~Golenia,
A.~M.~Samoilenko, Ya.~A.~Prykarpatsky, and A.~K.~Prykarpatsky,
``The general differential-geometric structure of multidimensional
Delsarte transmutation operators in parametric functional spaces
and their applications in soliton theory,'' {\it
arXiv:~math-ph/0404016}, 2004.

\bibitem{eng_GKMR}   R.~Gorenflo, A.~A.~Kilbas,
F.~Mainardi, and S.~V.~Rogosin, {\it Mittag-Leffler functions,
related topics and applications}, Springer, Heidelberg--New
York--Dordrecht--London, 2014.

\bibitem{eng_GoMa}
R.~Gorenflo and  F.~Mainardi, ``Fractional calculus: integral and
differential equations of fractional order,'' In: {\it Fractals
and fractional calculus in continuum mechanics}, Springer,
Wien--New York, 1997, pp.~223--278.

\bibitem{eng_HaKa1}
A.~Hasanov and E.~T.~Karimov, ``Fundamental solutions for a class
of three-dimensional elliptic equations with singular
coefficients,'' {\it Appl. Math. Lett.}, 2009,   {\bf 22},
1828--1832.

\bibitem{eng_Hig1}   T.~P.~Higgins, ``A hypergeometric
function transform,'' {\it J.~SIAM}, 1964,  {\bf  12}, No.~3,
601--612.

\bibitem{eng_HT}   E.~Hille and J.~D.~Tamarkin, ``On the theory of linear integral
equations,'' {\it Ann. Math.}, 1930,  {\bf  31}, 479--528.

\bibitem{eng_Hol}   M.~Holzleitner, ``Transformation operators for spherical
Schr\"odinger operators,'' {\it arXiv:1805.10526v1}, 2018.

\bibitem{eng_Jaf}   K.~Jafford, ``Inversion of the Lions
transmutation operators using generalized wavelets,'' {\it Appl.
Comput. Harmon. Anal.}, 1997,  {\bf  4}, No.~1, 97--112.

\bibitem{eng_Kam}   L.~Kamoun and M.~Sifi, ``Bessel--Struve
intertwining operator and generalized Taylor series on the real
line,'' {\it Integral Transforms Spec. Funct.}, 2005,  {\bf  16},
No.~1, 39--55.

\bibitem{eng_KarST}  S.~T.~Karimov, ``Multidimensional generalized Erd\'elyi--Kober
operator and its application to solving Cauchy problems for
differential equations with singular coefficients,'' {\it Fract.
Calc. Appl. Anal.}, 2015,  {\bf  18}, No.~4, 845--861.

\bibitem{eng_KaSr}   P.~W.~Karlsson and H.~M.~Srivastava, {\it Multiple
Gaussian hypergeometric series}, Ellis Horwood, New York, 1985.

\bibitem{eng_KaLo1}  D.~Karp and J.~L.~L\'opez, ``Representations of hypergeometric
functions for arbitrary values of the parameters and their use,''
{\it J.~Approx. Theory}, 2017,  {\bf 218}, 42--70.

\bibitem{eng_KaLo2}  D.~Karp and J.~L.~L\'opez, ``A
class of Meijer’s $G$ functions and further representations of the
generalized hypergeometric functions,'' {\it arXiv:1801.08670v1},
2018.

\bibitem{eng_KaPr4}
D.~Karp and E.~Prilepkina, ``Hypergeometric functions as
generalized Stieltjes transforms,'' {\it J.~Math. Anal. Appl.},
2012,  {\bf 393}, No.~2, 348--359.

\bibitem{eng_KaPr2}   D.~Karp and E.~Prilepkina, ``Completely monotonic gamma ratio and
infinitely divisible $H$-function of Fox,'' {\it Comput. Methods
Funct. Theory}, 2016, {\bf 16}, No.~1, 135--153.

\bibitem{eng_KaPr3}   D.~Karp and E.~Prilepkina, ``Hypergeometric differential equation and new identities for the
coefficients of N\o rlund and B\"uhring,'' {\it SIGMA Symmetry
Integrability Geom. Methods Appl.}, 2016,  {\bf 052}.

\bibitem{eng_KaPr1}   D.~Karp and E.~Prilepkina, ``Some new facts concerning the
delta neutral case of Fox's $H$-function,'' {\it Comput. Methods
Funct. Theory}, 2017, {\bf 17}, No.~2, 343--367.

\bibitem{eng_S13}   D.~Karp, A.~Savenkova, and
S.~M.~Sitnik, ``Series expansions for the third incomplete
elliptic integral via partial fraction decompositions,'' {\it
J.~Comput. Appl. Math.}, 2007,  {\bf  207},  No.~2, 331--337.

\bibitem{eng_S12}   D.~Karp and S.~M.~Sitnik, ``Asymptotic approximations
for the first incomplete elliptic integral near logarithmic
singularity,'' {\it J.~Comput. Appl. Math.}, 2007,  {\bf  205},
No.~1, 186--206.

\bibitem{eng_S15}
D.~Karp and S.~M.~Sitnik, ``Inequalities and monotonicity of
ratios for generalized hypergeometric function,'' {\it J.~Approx.
Theory}, 2009,  {\bf  161}, No.~1, 337--352.

\bibitem{eng_S16}   D.~Karp and S.~M.~Sitnik, ``Log-convexity and
log-concavity of hypergeometric-like functions,'' {\it J.~Math.
Anal. Appl.}, 2010, {\bf  364}, No.~2, 384--394.

\bibitem{eng_Ken1}   C.~E.~Kenig, ``Some recent
quantitative unique continuation theorems,'' {\it S\'emin. \'Equ.
D\'eriv. Partielles. \'Ec. Polytech. Cent. Math. Palaiseau}, 2006,
{\bf 2005-2006}, XX1--XX10.

\bibitem{eng_Ken2}   C.~E.~Kenig, L.~Silvestre, and J.~N.~Wang, ``On Landis' conjecture in
the plane,'' {\it Commun. Part. Differ. Equ.}, 2015,  {\bf 40},
No.~4, 766--789.

\bibitem{eng_KiSa}   A.~A.~Kilbas and M.~Saigo, {\it $H$-transforms. Theory and applications}, Chapman \& Hall/CRC, Boca Raton, 2004.

\bibitem{eng_KiSk1}   A.~A.~Kilbas and
O.~V.~Skoromnik, ``Integral transforms with the Legendre function
of the first kind in the kernels on
$\mathcal{L}_{\protect\nu,r}$-spaces,'' {\it Integral Transforms
Spec. Funct.}, 2009,  {\bf  20}, No.~9, 653--672.

\bibitem{eng_KST}   A.~A.~Kilbas,
H.~M.~Srivastava, and J.~J.~Truhillo, {\it Theory and applications
of fractional differential equations}, Elsevier, Amsterdam, etc.,
2006.

\bibitem{eng_KT1}   A.~A.~Kilbas and
J.~J.~Truhillo, ``Differential equations of fractional order:
methods, results and problems. Part~I,'' {\it Appl. Anal.}, 2001,
{\bf  78}, No.~1-2, 153--192.

\bibitem{eng_KT2}   A.~A.~Kilbas and J.~J.~Truhillo, ``Differential equations of
fractional order: methods, results and problems. Part~II,'' {\it
Appl. Anal.}, 2001, {\bf  81}, No.~2, 435--493.

\bibitem{eng_KiZhu}   A.~A.~Kilbas and
N.~V.~Zhukovskaya, ``Euler-type non-homogeneous differential
equations with three Liouville fractional derivatives,'' {\it
Fract. Calc. Appl. Anal.}, 2009,  {\bf 12}, No.~2, 205--234.

\bibitem{eng_Kir2}   V.~Kiryakova, ``An explanation of Stokes phenomenon
by the method of transmutations,'' {\it Proc. conf. Diff.
Equations and Appl., Rousse}, 1982, 349--353.

\bibitem{eng_Kir1}   V.~Kiryakova, {\it Generalized fractional
calculus and applications}, Longman, Harlow, 1994.

\bibitem{eng_Kir3}
V.~Kiryakova, ``Applications of the generalized Poisson
transformation for solving hyper-Bessel differential equations,''
{\it Godishnik VUZ. Appl. Math.}, 1986,  {\bf 22}, No.~4, 129--140
(in Bulgarian).

\bibitem{eng_Kir4}   V.~Kiryakova, ``All
the special functions are fractional differintegrals of elementary
functions,'' {\it J.~Phys. A. Math. Gen.}, 1997,  {\bf  30},
No.~14, 5085--5103.

\bibitem{eng_Kir6}
V.~Kiryakova, ``Multiple (multiindex) Mittag-Leffler functions and
relations to generalized fractional calculus,'' {\it J.~Comput.
Appl. Math.}, 2000,  {\bf  118}, No.~1-2, 241-­259.

\bibitem{eng_Kir5}   V.~Kiryakova, ``The multi-index Mittag-Leffler
functions as an important class of special functions of fractional
calculus,'' {\it Computers Math. Appl.}, 2010,  {\bf  59}, No.~5,
1885--1895.

\bibitem{eng_47}
H.~Kober, ``On fractional integrals and derivatives,'' {\it Quart.
J. Math. $($Oxford$)$}, 1940,  {\bf  11}, 193--211.

\bibitem{eng_Kob1}   H.~Kober, ``On a theorem of Schur
and on fractional integrals of purely imaginary order,'' {\it
Trans. Am. Math. Soc.}, 1941,  {\bf  50}, No.~1, 160--174.

\bibitem{eng_Koepf}   W.~Koepf, {\it Hypergeometric summation. An algorithmic approach to summation and special function
identities}, Vieweg, Wiesbaden, 1998.

\bibitem{eng_KMRS}   V.~Kokilashvili,
A.~Meskhi, H.~Rafeiro, and S.~Samko, {\it Integral operators in
non-standard function spaces. Vol.~1-2}, Birkh\"auser/Springer,
Basel, 2016.

\bibitem{eng_Koo1}
T.~H.~Koornwinder, ``Fractional integral and generalized Stieltjes
transforms for hypergeometric functions as transmutation
operators,'' {\it SIGMA Symmetry Integrability Geom. Methods
Appl.}, 2015,  {\bf 11}, No.~074.

\bibitem{eng_Kravch}
V.~V.~Kravchenko, {\it Applied pseudoanalytic function theory},
Birkh\"auser, Basel, 2009.

\bibitem{eng_Krav8}   V.~V.~Kravchenko, ``Construction of a transmutation for the
one-dimensional Schr\"odinger operator and a representation for
solutions,'' {\it Appl. Math. Comput.}, 2018,  {\bf  328}, 75--81.

\bibitem{eng_Krav4}
V.~V.~Kravchenko, L.~J.~Navarro, and S.~M.~Torba, ``Representation
of solutions to the one-dimensional Schr\"odinger equation in
terms of Neumann series of Bessel functions,'' {\it Appl. Math.
Comput.}, 2017,  {\bf 314}, No.~1, 173--192.

\bibitem{eng_Krav6}   V.~V.~Kravchenko, J.~A.~Otero, and S.~M.~Torba, ``Analytic
approximation of solutions of parabolic partial differential
equations with variable coefficients,'' {\it Adv. Math. Phys.},
2017,  {\bf  2017}, 2947275.

\bibitem{eng_Krav1}   V.~V.~Kravchenko and   S.~M.~Torba, ``Transmutations for Darboux
transformed operators with applications,'' {\it J.~Phys. A. Math.
Theor.}, 2012,  {\bf 45}, No.~7, 075201.

\bibitem{eng_Krav3}   V.~V.~Kravchenko and   S.~M.~Torba, ``Analytic approximation of transmutation operators and applications to highly accurate
solution of spectral problems,'' {\it J.~Comput. Appl. Math.},
2015,  {\bf  275}, 1--26.

\bibitem{eng_Krav2}   V.~V.~Kravchenko and S.~M.~Torba, ``Construction of
transmutation operators and hyperbolic pseudoanalytic functions,''
{\it Complex Anal. Oper. Theory}, 2015,  {\bf  9}, No.~2,
379--429.

\bibitem{eng_Krav7}
V.~V.~Kravchenko and S.~M.~Torba, ``Asymptotics with respect to
the spectral parameter and Neumann series of Bessel functions for
solutions of the one-dimensional Schr\"odinger equation,'' {\it
J.~Math. Phys.}, 2017,  {\bf  58}, No.~12, 122107.

\bibitem{eng_Krav9}   V.~V.~Kravchenko and S.~M.~Torba, ``A Neumann series of Bessel
functions representation for solutions of Sturm--Liouville
equations,'' {\it Calcolo}, 2018,  {\bf  55}, No.~1, Paper No.~11.

\bibitem{eng_Krav10}
V.~V.~Kravchenko, S.~M.~Torba, and R.~Castillo-P\'erez, ``A
Neumann series of Bessel functions representation for solutions of
perturbed Bessel equations,'' {\it Appl. Anal.}, 2018,  {\bf  97},
No.~5, 677--704.

\bibitem{eng_Krav5}   V.~V.~Kravchenko,
S.~M.~Torba, and K.~V.~Khmelnytskaya, ``Transmutation operators:
construction and applications,'' {\it Proc. of the 17th Int. Conf.
on Comput. and Math. Methods in Sci. and Eng. CMMSE-2017}, Cadiz,
Andalucia, Espa\~na, 2017, pp.~1198--1206.

\bibitem{eng_KMP}
A.~Kufner, L.~Maligranda, and L.-E.~Persson, {\it The Hardy
inequality. About its history and some related results},
Vydavatelsk\'y Servis, Pilsen, 2007.

\bibitem{eng_Lev1995}
B.~M.~Levitan, ``Transmutation operators and the inverse spectral
problem,'' {\it Contemp. Math.}, 1995,  {\bf 183}, 237--244.

\bibitem{eng_Lind}   P.~Lindqvist, ``Notes
on the $p$-Laplace equation,'' {\it Report: Univ. of
Jyv\"askyl\"a. Dep. of Math. and Stat.}, 2006,  {\bf 102}.

\bibitem{eng_Lio2}   J.~L.~Lions, ``Op\'erateurs de Delsarte et
probl\`eme mixte,'' {\it Bull. Soc. Math. France}, 1956,  {\bf
84}, 9--95.

\bibitem{eng_Lio3}   J.~L.~Lions, ``Quelques applications  d'op\'erateurs de
transmutations,'' {\it Colloques Internat. Nancy}, 1956, 125--142.

\bibitem{eng_Lio1}   J.~L.~Lions, {\it Equations differentielles operationnelles et
probl\'emes aux limites}, Springer,
Berlin--G\"ottingen--Heidelberg, 1961.

\bibitem{eng_Love1}   E.~R.~Love, ``Some integral equations involving hypergeometric
functions,'' {\it Proc. Edinb. Math. Soc.}, 1967,  {\bf  15},
No.~3, 169--198.

\bibitem{eng_Love2}   E.~R.~Love, ``Two more hypergeometric integral equations,'' {\it Proc. Cambridge Phil. Soc.}, 1967,  {\bf  63}, No.~4, 1055--1076.

\bibitem{eng_Low1}   J.~S.~Lowndes, ``An application of some fractional integrals,'' {\it Glasg. Math. J.}, 1979,  {\bf  20}, No.~1, 35--41.

\bibitem{eng_Low2}
J.~S.~Lowndes, ``On some generalizations of Riemann--Liouville and
Weil fractional integrals and their applications,'' {\it Glasg.
Math. J.}, 1981,  {\bf  22}, No.~2, 73--80.

\bibitem{eng_Low3}   J.~S.~Lowndes, ``Cauchy problems for second order
hyperbolic differential equations with constant coefficients,''
{\it Proc. Edinb. Math. Soc.}, 1983,  {\bf  26}, No.~3, 97--105.

\bibitem{eng_Lud}   D.~Ludwig, ``The
Radon transform on Euclidean space,'' {\it Math. Methods Appl.
Sci.}, 1966,  {\bf 19}, 49--81.

\bibitem{eng_Luke2}   Y.~L.~Luke, {\it The special functions and their
approximations. V.~1}, Academic Press, New York--London, 1969.

\bibitem{eng_Luke1}   Y.~L.~Luke, {\it Mathematical functions and their approximations}, Academic Press, New York--San
Francisco--London, 1975.

\bibitem{eng_Matv}   V.~B.~Matveev, ``Intertwining relations between the Fourier transform and discrete Fourier transform, the
related functional identities and beyond,'' {\it Inverse
Problems}, 2001,  {\bf  17}, 633--657.

\bibitem{eng_MaSa}   V.~B.~Matveev and V.~B.~Salle, {\it Darboux--Backlund transformations and applications}, Springer, N.Y., 1991.

\bibitem{eng_McB}   A.~C.~McBride, {\it Fractional calculus and integral transforms of
generalized functions}, Pitman, San Francisco--London--Melbourne,
1979.

\bibitem{eng_MS1}   Kh.~Mehrez and S.~M.~Sitnik, ``On monotonicity of ratios of some
$q$-hypergeometric functions,'' {\it Mat. Vesn.}, 2016,  {\bf 68},
No.~3, 225--231.

\bibitem{eng_MS3}   Kh.~Mehrez and S.~M.~Sitnik, ``Proofs of some
conjectures on monotonicity of ratios of Kummer, Gauss and
generalized hypergeometric functions,'' {\it Analysis
$($Munich$)$}, 2016,  {\bf  36}, No.~4, 263--268.

\bibitem{eng_MS2}   Kh.~Mehrez and S.~M.~Sitnik, ``Functional Inequalities for the
Mittag-Leffler Functions,'' {\it Results Math.}, 2017,  {\bf  72},
No.~1-2, 703--714.

\bibitem{eng_Rad3}   E.~Mitidieri, J.~Serrin, and V.~Radulescu {\rm
(, {\it eds.{\rm )}} Recent trends in nonlinear partial
differential equations~I: Evolution problems}, Am. Math. Soc.,
Providence, 2013.

\bibitem{eng_Rad4}   E.~Mitidieri,
J.~Serrin, and V.~Radulescu {\rm (, {\it eds.{\rm )}} Recent
trends in nonlinear partial differential equations~II: Stationary
Problems}, Am. Math. Soc., Providence, 2013.

\bibitem{eng_MPF}   D.~S.~Mitrinovi\'c,  J.~E.~Pe\v cari\'c, and A.~M.~Fink, {\it Classical
and new inequalities in analysis}, Kluwer, Dordrecht, 1993.

\bibitem{eng_Mur7}
A.~B.~Muravnik, ``On weighted norm estimates for the mixed
Fourier--Bessel transforms on non-negative functions,'' In: {\it
Integral methods in science and engineering. Vol.~1: Analytic
methods}, Longman, Harlow, 1997, pp.~119--123.

\bibitem{eng_Mur6}
A.~B.~Muravnik, ``Fourier--Bessel transformation of measures and
singular differential operators,'' In: {\it Paul Erd\H{o}s and his
mathematics}, J'anos Bolyai Math. Soc., Budapest, 1999,
pp.~182--184.

\bibitem{eng_Mur5}   A.~B.~Muravnik, ``Fourier--Bessel transformation of measures with several special variables and properties
of singular differential equations,'' {\it J.~Korean Math. Soc.},
2000,  {\bf 37}, No.~6, 1043--1057.

\bibitem{eng_Mur4}   A.~B.~Muravnik, ``Fourier--Bessel
transformation of compactly supported non-negative functions and
estimates of solutions of singular differential equations,'' {\it
Funct. Differ. Equ.}, 2001,  {\bf  8}, No.~3-4, 353--363.

\bibitem{eng_Mur3}   A.~B.~Muravnik, ``Fourier--Bessel
transformation of measures and singular differential equations,''
In: {\it Recent progress in functional analysis}, Elsevier,
Amsterdam, 2001, pp.~335--345.

\bibitem{eng_Mur2}   A.~B.~Muravnik, ``Nonclassical Cauchy problem for singular
parabolic integro-differential equations,'' {\it Russ. J. Math.
Phys.}, 2002,  {\bf  9}, No.~3, 300--314.

\bibitem{eng_Mur1}   A.~B.~Muravnik, ``On stabilization of
solutions of elliptic equations containing Bessel operators,'' In:
{\it Integral methods in science and engineering. Analytic and
numerical techniques}, Birkh\"auser, Boston, 2004, pp.~157--162.

\bibitem{eng_NIST}    F.~W.~J.~Olver, D.~W.~Lozier,
R.~F.~Boisvert, and C.~W.~Clark, {\it NIST handbook  of
mathematical functions}, Cambridge University Press, Cambridge,
2010.

\bibitem{eng_OpKu}   B.~Opic and A.~Kufner, {\it Hardy-type
inequalities}, Longman, Harlow, 1990.

\bibitem{eng_OZK}   H.~Ozaktas, Z.~Zalevsky, and
M.~Kutay, {\it The fractional Fourier transform: with applications
in optics and signal processing}, John Wiley \& Sons, Chishester,
ets., 2001.

\bibitem{eng_Jor}
J.~Paneva-Konovska, {\it From Bessel to multi-index Mittag-Leffler
functions}, World Scientific, London, 2016.

\bibitem{eng_PWZ}   M.~Petkov\v sek,  H.~S.~Wilf, and D.~Zeilberger, {\it $A=B$}, A.K. Peters, Wellesley, 1996.

\bibitem{eng_PiSa}
S.~Pike and P.~Sabatier, {\it Scattering. Scattering and inverse
scattering in pure and applied science. Vol.~1-2}, Academic Press,
San Diego, 2002.

\bibitem{eng_Poisson}
S.~D.~Poisson, ``M\'emoire sur l'int\'egration des \'equations
lin\'eaires aux diff\'erences partielles,'' {\it J.~\'Ec. Roy.
Polytech.}, 1823,  {\bf 19}, No.~12, 215--248.

\bibitem{eng_Pou}   A.~D.~Poularicas $($ed.$)$, {\it The
transforms and applications handbook}, CRC Press, Boca Raton,
2010.

\bibitem{eng_Jan}  J.~Pr\"uss, {\it Evolutionary integral equations and applications}, Birkh\"auser, Basel, 2012.

\bibitem{eng_Pyat}   S.~G.~Pyatkov, {\it Operator theory.
Nonclassical problems}, VSP, Utrecht, 2002.

\bibitem{eng_Rad2}   V.~Radulescu and
D.~Repovs, {\it Partial differential equations with variable
exponents: Variational methods and qualitative analysis}, CRC
Press, Boca Raton, 2015.

\bibitem{eng_Riesz}
M.~Riesz, ``L'int\'egrale de Riemann--Liouville et le probl\`eme
de Cauchy,'' {\it Acta Math.}, 1949,  {\bf 81}, 1--223.

\bibitem{eng_Riman}   B.~Riemann, ``On the
propagation of flat waves of finite amplitude,'' In: {\it Ouvres},
OGIZ, Moscow--Leningrad, 1948, pp.~376--395.

\bibitem{eng_Rod}   J.~Rodrigues, ``Operational calculus for the generalized Bessel operator,'' {\it Serdica Math.~J.}, 1989,  {\bf  15}, 179--186.

\bibitem{eng_Ros1}   M.~R\"osler, ``Positivity of
Dunkl’s intertwining operator,'' {\it Duke Math.~J.}, 1999,  {\bf
98}, 445--463.

\bibitem{eng_Ros2}   M.~R\"osler, ``Dunkl operators: theory and
applications,'' In: {\it Orthogonal polynomials and special
functions}, Springer, Berlin, 2003, pp.~93--135.

\bibitem{eng_Rossi}  L.~Rossi, ``The Landis conjecture
with sharp rate of decay,'' {\it arXiv:1807.00341v1}, 2018.

\bibitem{eng_Rub4}
B.~Rubin, ``On the Funk--Radon--Helgason inversion method in
integral geometry,'' {\it Contemp. Math.}, 2013,  {\bf  599},
175--198.

\bibitem{eng_Rub2}   B.~Rubin, ``Gegenbauer--Chebyshev integrals and Radon transforms,'' {\it arXiv:1410.4112v2}, 2015.

\bibitem{eng_Rub1}   B.~Rubin, ``Radon transforms and Gegenbauer--Chebyshev
integrals.~I,'' {\it Anal. Math. Phys.}, 2017,  {\bf 7}, No.~2,
117--150.

\bibitem{eng_SaHa1}   M.~S.~Salakhitdinov and A.~Hasanov, ``The Dirichlet problem for the
generalized bi-axially symmetric Helmholtz equation,'' {\it
Eurasian Math.~J.}, 2012, {\bf  3}, No.~4, 99--110.

\bibitem{eng_SaKiMar}   S.~G.~Samko, A.~A.~Kilbas, and
O.~I.~Marichev, {\it Fractional integrals and derivatives: theory
and applications}, Gordon \& Breach, N.Y., 1993.

\bibitem{eng_SPP}   A.~M.~Samoilenko, Ya.~A.~Prykarpatsky, and
A.~K.~Prykarpatsky, ``The generalized de Rham--Hodge theory
aspects of Delsarte--Darboux type transformations in
multidimension,'' {\it Cent. Eur. J. Math.}, 2005,  {\bf  3},
No.~3, 529--557.

\bibitem{eng_74}   R.~T.~Seely, ``Extensions of $C^{\infty}$
functions defined in a half space,'' {\it Proc. Am. Math. Soc.},
1964,  {\bf  15}, 625--626.

\bibitem{eng_75}   R.~T.~Seely, ``Complex powers of an elliptic
operator,'' {\it Proc. Sympos. Pure Math.}, 1967,  {\bf  10},
288--307.

\bibitem{eng_ShiE1}  E.~L.~Shishkina, ``Inversion of integral of $B$-potential type
with density from $\Phi_\gamma$,'' {\it J.~Math. Sci.}, 2009, {\bf
160}, No.~1, 95--102.

\bibitem{eng_ShiE2}  E.~L.~Shishkina, ``On the boundedness of hyperbolic
Riesz $B$-potential,'' {\it Lith. Math.~J.}, 2016,  {\bf  56},
No.~4, 540--551.

\bibitem{eng_ShiE5}  E.~L.~Shishkina, ``Generalized Euler--Poisson--Darboux
equation and singular Klein--Gordon equation,'' {\it J.~Phys.
Conf. Ser.}, 2018,  {\bf 973}, 1--21.

\bibitem{eng_ShiE4}  E.~L.~Shishkina, ``Properties of mixed
hyperbolic $B$-potential,'' {\it Progr. Fract. Differ. Appl.},
2018,  {\bf  4}, No.~2, 1--16.

\bibitem{eng_ShiE6}  E.~L.~Shishkina, ``Singular Cauchy problem for
the general Euler--Poisson--Darboux equation,'' {\it Open Math.},
2018,  {\bf  16}, 23--31.

\bibitem{eng_ShiE3}  E.~L.~Shishkina and S.~M.~Sitnik, ``General form
of the Euler--Poisson--Darboux equation and application of the
transmutation method,'' {\it Electron. J.~Differ. Equ.}, 2017,
{\bf  177}, 1--20.

\bibitem{eng_SS}
E.~L.~Shishkina and S.~M.~Sitnik, ``On fractional powers of Bessel
operators,'' {\it J.~Inequal. Spec. Funct.}, 2017,  {\bf 8},
No.~1, 49--67.

\bibitem{eng_Sie}
J.~Siersma, {\it Thesis}, Groningen, 1979.

\bibitem{eng_S41}   S.~M.~Sitnik, ``Generalized Young and Cauchy--Bunyakowsky inequalities with
applications: a survey,'' {\it [math.CA] arXiv:1012.3864}, 2010.

\bibitem{eng_S42}   S.~M.~Sitnik, ``Transmutations and applications: a survey,'' {\it [math.CA]
arXiv:1012.3741}, 2010.

\bibitem{eng_S94}  S.~M.~Sitnik, ``Some problems in the modern theory of
transmutations,'' {\it Spectral theory and differential equations
(STDE-2012). Int. Conf. in honor of V.~A.~Marchenko's 90th
birthday}, Kharkiv, 2012, pp.~101-102.

\bibitem{eng_S38}  S.~M.~Sitnik, ``Buschman--Erd\'elyi transmutations,
classification and applications,'' In: {\it Analytic methods of
analysis   and differential equations}, Cambridge Scientific
Publishers, Cambridge, 2013, pp.~171--201.

\bibitem{eng_S401}   S.~M.~Sitnik, ``A short survey of recent results on
Buschman--Erd\'elyi transmutations,'' {\it J. Inequal. Spec.
Funct.}, 2017,  {\bf  8}, No.~1, 140--157.

\bibitem{eng_S402}   S.~M.~Sitnik, ``Buschman--Erd\'elyi
transmutations and  applications,'' {\it Abstr. of the 8th Int.
Conf. Transform Methods and Special Functions, Bulgaria, Sofia,
27--31 Aug. 2017}, Inst. Math. Inf. Bulg. Acad. Sci., 2017, p.~59.

\bibitem{eng_Skub}   A.~L.~Skubachevskii, {\it Elliptic functional
differential equations and applications}, Birkh\"auser, Basel,
1997.

\bibitem{eng_Sla}   L.~J.~Slater, {\it Generalized hypergeometric functions}, Cambridge
University Press, Cambridge, 1966.

\bibitem{eng_Spr}   I.~G.~Sprinkhuizen-Kuyper, ``A fractional
integral operator corresponding to negative powers of a certain
second-order differential operator,'' {\it J.~Math. Anal. Appl.},
1979,  {\bf 72}, 674--702.

\bibitem{eng_SvFe}   G.~A.~Sviridyuk and V.~E.~Fedorov, {\it Linear Sobolev type equations
and degenerate semigroups of operators}, VSP, Utrecht, 2003.

\bibitem{eng_Ta1}
Ta Li, ``A new class of integral transform,'' {\it Proc. Am. Math.
Soc.}, 1960,  {\bf  11}, No.~2, 290--298.

\bibitem{eng_Ta2}   Ta Li, ``A note on integral transform,'' {\it Proc. Am. Math. Soc.}, 1961,  {\bf  12}, No.~6, 556.

\bibitem{eng_Tri2}   Kh.~Trim\'eche, ``Transformation int\'egrale de Riemann--Liouville g\'en\'eralises et
convergence des series de Taylor g\'en\'eralis\'ees au sens de
Delsarte,'' {\it Rev. Fac. Sci. Tunis}, 1981,  No.~1, 7--14.

\bibitem{eng_Tri1}   Kh.~Trim\'eche, ``Transformation int\'egrale de Weil et th\`eor\'eme de Paley--Winer associ\'es \`a un
op\'erateur diff\'erentiel singulier sur $(0,\infty)$,'' {\it
J.~Math. Pures Appl.}, 1981,  {\bf 60}, 51--98.

\bibitem{eng_Tri3}   Kh.~Trim\'eche, ``Transmutation
operators and mean-periodic functions associated with differential
operators,'' {\it Math. Rep. Ser.~4}, 1988,  No.~1, i-xiv.

\bibitem{eng_Tri4}   Kh.~Trim\'eche, ``Inversion
of the Lions transmutation operators using generalized wavelets,''
{\it Appl. Comput. Harmon. Anal.}, 1997,  {\bf  4}, No.~1,
97--112.

\bibitem{eng_Tri5}
Kh.~Trim\'eche, {\it Generalized harmonic analysis and wavelet
packets}, Gordon and Breach, Amsterdam, 2001.

\bibitem{eng_Tri6}   Kh.~Trim\'eche, ``Inversion of the Dunkl
intertwining operator and its dual using Dunkl wavelets,'' {\it
Rocky Mountain J.~Math.}, 2002,  {\bf  32}, No.~2, 889--895.

\bibitem{eng_Vir2}   N.~Virchenko, ``On
some generalized symmetric integral operators of
Buschman--Erd\'elyi’s type,'' {\it J.~Nonlinear Math. Phys.},
1996,  {\bf  3}, No.~3-4, 421--425.

\bibitem{eng_Vir1}   N.~Virchenko and I.~Fedotova, {\it Generalized associated Legendre
functions and their applications}, World Scientific, Singapore,
2001.

\bibitem{eng_Volch}   V.~V.~Volchkov, {\it Integral geometry and convolution equations}, Kluwer, Dordrecht, 2003.

\bibitem{eng_14}   N.~Wiener, ``The Dirichlet problem,'' {\it J.~Math. Phys. Mass. Inst. Techn.}, 1924,  {\bf  3}, 127--146.

\bibitem{eng_Wei1}   A.~Weinstein, ``Discontinuous integrals and generalized theory of
potential,'' {\it Trans. Am. Math. Soc.}, 1948,  {\bf  63}, No.~2,
342--354.

\bibitem{eng_Wei2}   A.~Weinstein, ``Generalized axially symmetric potential theory,'' {\it Bull. Am. Math. Soc.}, 1953,  {\bf  59}, 20--38.

\bibitem{eng_Wei3}
A.~Weinstein, {\it Selecta}, Pitman, London--San Francisco, 1978.

}


\begin{thebibliography}{999}

\bibitem{AbOs}  {\it Абжандадзе~З.\,Л., Осипов~В.\,Ф.} Преобразование Фурье---Френеля и некоторые его приложения. "---  СПб.: Изд-во С.-Петерб. ун-та, 2000.

\bibitem{AbSi} {\it Абловиц~М., Сигур~Х.} Солитоны и метод обратной задачи. "--- М.: Мир, 1979.

\bibitem{AS} {\it Абрамовиц~М., Стиган~И.} Справочник по специальным функциям. "--- М.: Наука, 1979.

\bibitem{1} {\it Агмон~С., Дуглис~А., Ниренберг~Л.} Оценка решений эллиптических уравнений вблизи границы. "--- М.: Иностр. лит., 1962.

\bibitem{AM} {\it Агранович~З.\,С., Марченко~В.\,А.} Обратная задача теории рассеяния. "--- Харьков: Изд. ХГУ, 1960.

\bibitem{Agr} {\it Агранович~М.\,С.} Эллиптические псевдодифференциальные операторы. Ч.~1,~2. "--- М.: 2003, 2004.

\bibitem{AMR} {\it Азбелев~Н.\,В., Максимов~В.\,П., Рахматуллина~Л.\,Ф.}
Элементы современной теории функционально-дифференциальных
уравнений. Методы и приложения. "--- М.: Ин-т комп. иссл., 2002.

\bibitem{2} {\it  Алимов~Ш.\,А.}  Дробные степени эллиптических операторов и изоморфизм классов дифференцируемых функций//
Дифф. уравн. "--- 1972. "--- {\sl 8}, \No~9. "--- С.~1609--1626.

\bibitem{Arsh} {\it  Аршава~Е.\,А.} Обращение интегральных операторов методом
операторных тождеств// Науч. ведом. Белгород. гос. ун-та. Сер.
Мат. Физ. "--- 2009. "--- {\sl 17/2}, \No~13 (68). "--- С.~18--29.

\bibitem{Ahi1} {\it Ахиезер~Н.\,И.} К теории спаренных интегральных уравнений// Уч. зап. Харьков. гос. ун-та. "--- 1957. "--- {\sl 80}. "--- С.~5--21.

\bibitem{Ahi2} {\it Ахиезер~Н.\,И.} Лекции об интегральных преобразованиях. "--- Харьков: Вища школа. Изд-во при Харьк. ун-те, 1984.

\bibitem{3} {\it Бабиков~В.\,В.} Метод фазовых функций в квантовой механике. "--- М.: Наука, 1976.

\bibitem{BMYa} {\it Баврин~И.\,И., Матросов~В.\,Л.,  Яремко~О.\,Э.} Операторы преобразования для краевых задач,
интегральных представлений и восстановления зависимостей. "--- М.: Прометей, 2015.

\bibitem{BaSa} {\it Багров~В.\,Г.,  Самсонов~Б.\,A.} Преобразование Дарбу уравнения Шрёдингера//
Физ. элем. частиц и атом. ядра. "--- 1997. "--- {\sl 28}, \No~4.
"--- С.~951--1012.

\bibitem{Baid} {\it Байдаков~А.\,Н.} Априорные оценки гёльдеровых норм решений квазилинейных $B$-эллиптических уравнений//
Дифф. уравн. "--- 1987. "--- {\sl 23},  \No~11. "---
C.~1923--1930.

\bibitem{Bas} {\it Баскаков~А.\,Г.} Гармонический анализ линейных операторов. "--- Воронеж: ВГУ, 1987.

\bibitem{BE1} {\it Бейтмен~Г.,  Эрдейи~А.} Высшие трансцендентные
функции. Т.~1. "--- М.: Наука, 1966.

\bibitem{BE2}  {\it Бейтмен~Г.,  Эрдейи~А.} Высшие трансцендентные
функции. Т.~2. "--- М.: Наука, 1966.

\bibitem{BE3} {\it Бейтмен~Г.,  Эрдейи~А.} Высшие трансцендентные
функции. Т.~3. "--- М.: Наука, 1967.

\bibitem{BB} {\it Беккенбах~Э.,  Беллман~Р.} Неравенства. "--- М.: Мир, 1965.

\bibitem{Berg} {\it Бергман~С.} Интегральные операторы в теории уравнений с частными производными. "--- М.: Мир, 1964.

\bibitem{6} {\it Березанский~Ю.\,M.} Разложения по собственным функциям самосопряжённых операторов. "--- Киев: Наукова думка, 1965.

\bibitem{Berk} {\it Беркович~Л.\,М.} Факторизация и преобразования дифференциальных уравнений. Методы и приложения. "--- М.: РХД, 2002.

\bibitem{7} {\it  Берс~Л.} Математические вопросы дозвуковой и околозвуковой газовой динамики. "--- М.: Иностр. лит., 1961.

\bibitem{8} {\it Бесов~О.\,В., Ильин~В.\,П., Никольский~С.\,М.} Интегральные представления функций и теоремы вложения. "--- М.: Наука, 1975.

\bibitem{Bitz2} {\it Бицадзе~А.\,В.} Уравнения смешанного типа. "--- М.: Изд. АН СССР, 1959.

\bibitem{9} {\it Бицадзе~А.\,В.} Некоторые классы уравнений в частных производных. "--- М.: Наука, 1981.

\bibitem{Bitz1} {\it Бицадзе~А.\,В., Пашковский~В.\,И.}   К теории уравнений Максвелла---Эйнштейна// Докл. АН СССР. "--- 1974. "--- {\sl 216}, \No~2. "--- С.~9-10.

\bibitem{Bitz12} {\it Бицадзе~А.\,В., Пашковский~В.\,И.}   О некоторых классах решений уравнения Максвелла---Эйнштейна// Тр. МИАН. "--- 1975. "--- {\sl 134}. "--- C.~26--30.

\bibitem{Bloh} {\it Блох~А.\,Ш.} Об определении дифференциального оператора по его спектральной матрице-функции// Докл. АН СССР. "--- 1953. "--- {\sl 92}, \No~2. "--- С.~209--212.

\bibitem{Bor1} {\it Боровских~А.\,В.} Формула распространяющихся волн для одномерной неоднородной среды// Дифф. уравн. "--- 2002. "--- {\sl 38}, \No~6. "--- С.~758--767.

\bibitem{Bor2} {\it Боровских~А.\,В.} Метод распространяющихся волн//  Тр. сем. им. И.\,Г.~Петровского. "--- 2004. "--- {\sl 24}. "--- С.~3--43.

\bibitem{Boyar} {\it  Боярский~Б.} обобщённые решения системы дифференциальных
уравнений первого порядка эллиптического типа с разрывными
коэффициентами//  Мат. сб. "--- 1957. "--- {\sl 43}, \No~4. "---
С.~451--503.

\bibitem{11} {\it Брычков~Ю.\,А., Прудников~А.\,П.} Интегральные преобразования обобщённых функций. "--- М.: Наука, 1977.

\bibitem{Burb} {\it Бурбаки~Н.} Функции действительного переменного. "--- М.: Наука, 1965.

\bibitem{Bur1} {\it Буренков~В.\,И.} Функциональные пространства. "--- М.: РУДН, 1989.

\bibitem{Bur2} {\it Буренков~В.\,И.,  Гольдман~М.\,Л.}  Методические рекомендации к изучению курса <<Функциональные пространства>>. "--- М.: РУДН, 1989.

\bibitem{But} {\it Бутерин~С.\,А.} О восстановлении свёрточного возмущения оператора Штурма---Лиувилля по спектру// Дифф. уравн. "--- 2010. "--- {\sl 46}, \No~1. "--- С.~146--149.

\bibitem{Val}  {\it  Валицкий~Ю.\,Н.} Об операторе преобразования для интегро-дифференциальных операторов типа Вольтерра//
В сб.: <<Математическая физика>>. "--- Киев: Наукова Думка, 1965. "--- С.~23--36.

\bibitem{VaRo} {\it Варфоломеев~Е.\,М., Россовский~Л.\,Е.} Функционально-дифференциальные уравнения и их приложения к исследованию нейронных сетей и передаче информации нелинейными
лазерными системами с обратной связью. Учеб. пособие. "--- М.: РУДН, 2008.

\bibitem{Wat} {\it  Ватсон~Г.\,Н.} Теория бесселевых функций. Т.~1. "--- М.: Иностр. лит., 1949.

\bibitem{12}  {\it  Вашарин~А.\,А., Лизоркин~П.\,И.}  Некоторые краевые задачи для эллиптических уравнений с сильным вырождением на границе//
Докл. АН СССР. "--- 1961. "--- {\sl 137}, \No~5. "--- С.~1015--1018.

\bibitem{Vek1} {\it  Векуа~И.\,Н.}  О решениях уравнения $\Delta u +
\lambda^2 u$// Сообщ. АН Груз. ССР. "--- 1942. "---
{\sl 3}, \No~4. "--- С.~307--314.

\bibitem{Vek2} {\it  Векуа~И.\,Н.} Обращение одного интегрального
преобразования и его некоторые применения// Сообщ. АН Груз. ССР. "--- 1945. "--- {\sl 6}, \No~3. "--- C.~177--183.

\bibitem{13} {\it  Векуа~И.\,Н.} Об одном обобщении интеграла Пуассона для плоскости// Докл. АН СССР. "--- 1947. "--- {\sl 56}, \No~2. "--- С.~229--231.

\bibitem{Vek3} {\it  Векуа~И.\,Н.} Новые методы решения эллиптических уравнений. "--- М.---Л.: ГИТТЛ, 1948.

\bibitem{Vek4} {\it  Векуа~И.\,Н.} Обобщённые аналитические функции. "--- М.: Наука, 1988.

\bibitem{ViGa} {\it  Вирченко~Н.\,А., Гайдей~В.}   Классические и обобщённые многопараметрические функции. "--- Киев, 2008 (на украинском языке).

\bibitem{ViRy} {\it  Вирченко~Н.\,А., Рыбак~В.\,Я.} Основы дробного интегродифференцирования. "--- Киев, 2007 (на украинском языке).

\bibitem{15} {\it  Вишик~М.\,И., Грушин~В.\,В.} Краевые задачи для эллиптических уравнений, вырождающихся на границе области// Мат. сб. "--- 1969. "--- {\sl 80}. "--- С.~455--491.

\bibitem{16} {\it  Владимиров~B.\,C.} обобщённые функции в математической физике. "--- М.: Наука, 1979.

\bibitem{Volk} {\it  Волк~В.\,Я.} О формулах обращения для дифференциального уравнения с особенностью при $x=0$// Усп. мат. наук. "--- 1953. "--- {\sl 111}, \No~4. "--- С.~141--151.

\bibitem{VoKa} {\it  Волков~И.\,К., Канатников~А.\,Н.} Интегральные преобразования и операционное исчисление. "--- М.: Изд-во МГТУ им. Н.\,Э.~Баумана, 2002.

\bibitem{VZ} {\it Волкодавов~В.\,Ф., Захаров~В.\,Н.} Таблицы функций Римана и Римана---Адамара для некоторых дифференциальных уравнений в $n$-мерных евклидовых пространствах. "---
Самара, 1994.

\bibitem{VNN} {\it  Волкодавов~В.\,Ф., Лернер~М.\,Е., Николаев~Н.\,Я., Носов~В.\,А.} Таблицы некоторых функций Римана, интегралов и рядов. "---
Куйбышев: Изд. Куйбышев. гос. пед. ин-та, 1982.

\bibitem{VoNi} {\it  Волкодавов~В.\,Ф., Николаев~Н.\,Я.}  Краевые задачи для уравнения Эйлера---Пауссона---Дарбу. "--- Куйбышев: Изд. Куйбышев. гос. пед. ин-та, 1984.

\bibitem{VN} {\it  Волкодавов~В.\,Ф., Николаев~Н.\,Я.} Интегральные уравнения Вольтерра первого рода с некоторыми специальными функциями в ядрах и их приложения. "---
Самара: Изд.-во <<Самарский ун-т>>, 1992.

\bibitem{Vrag} {\it Врагов~В.\,Н.}  Краевые задачи для неклассических уравнений математической физики. "--- Новосибирск: НГУ, 1983.

\bibitem{Gan} {\it  Гантмахер~Ф.\,Р.} Теория матриц. "--- М.: Наука, 1988.

\bibitem{17} {\it  Гельфанд~И.\,М., Граев~М.\,И., Виленкин~Н.\,Я.} Интегральная геометрия и связанные с ней вопросы теории представлений. "--- М.: Гос. изд-во физ.-мат. лит., 1962.

\bibitem{Glu2} {\it  Глушак~А.\,В.} О возмущении абстрактного уравнения Эйлера---Пуассона---Дарбу// Мат. заметки. "--- 1996. "--- {\sl 60}, \No~3. "--- С.~363--369.

\bibitem{Glu4} {\it  Глушак~А.\,В.} О стабилизации решения задачи Дирихле для одного эллиптического уравнения в банаховом пространстве// Дифф. уравн. "--- 1997. "--- {\sl 33}, \No~4. "--- С.~433--437.

\bibitem{Glu3} {\it  Глушак~А.\,В.} Операторная функция Бесселя// Докл. РАН. "--- 1997. "--- {\sl 352}, \No~5. "--- С.~587--589.

\bibitem{Glu5} {\it  Глушак~А.\,В.} Операторная функция Бесселя и связанные с нею полугруппы и модифицированное преобразование Гильберта//
Дифф. уравн. "--- 1999. "--- {\sl 35}, \No~1. "--- С.~128--130.

\bibitem{Glu55} {\it  Глушак~А.\,В.} Регулярное и сингулярное возмущения абстрактного уравнения Эйлера---Пуассона---Дарбу// Мат. заметки. "--- 1999. "--- {\sl 66}, \No~3. "--- С.~364--371.

\bibitem{Glu6} {\it  Глушак~А.\,В.} Операторная функция Лежандра// Изв. РАН. Сер. Мат. "--- 2001. "--- {\sl 65}, \No~6. "--- С.~3--14.

\bibitem{Glu7} {\it  Глушак~А.\,В.} Задача типа Коши для абстрактного дифференциального уравнения с дробными производными// Мат. заметки. "--- 2005. "--- {\sl 77}, \No~1. "--- С.~28--41.

\bibitem{Glu8} {\it  Глушак~А.\,В.} О связи проинтегрированной косинус-оператор-функции с операторной функцией Бесселя// Дифф. уравн. "--- 2006. "--- {\sl 42}, \No~5. "--- С.~583--589.

\bibitem{Glu10} {\it  Глушак~А.\,В.} Начальная задача для слабо нагруженного уравнения Эйлера---Пуассона---Дарбу// Мат. межд. науч. конф.
<<Актуальные проблемы теории уравнений в частных производных>>. "--- М.:  МГУ,  2016. "--- С.~101.

\bibitem{Glu11} {\it  Глушак~А.\,В.} Нелокальная задача для абстрактного уравнения Эйлера---Пуассона---Дарбу// Изв. вузов. Сер. Мат. "--- 2016. "--- \No~6. "--- С.~1--9.

\bibitem{Glu12} {\it  Глушак~А.\,В.} Абстрактная задача Коши для уравнения Бесселя---Струве// Дифф. уравн. "--- 2017. "--- {\sl 53}, \No~7. "--- С.~891--905.

\bibitem{Glu1} {\it  Глушак~А.\,В., Кононенко~В.\,И., Шмулевич~С.\,Д.} Об одной сингулярной абстрактной задаче Коши// Изв. вузов. Сер. Мат. "--- 1986. "--- \No~6. "--- С.~55--56.

\bibitem{Glu9} {\it  Глушак~А.\,В., Покручин~О.\,А.} Критерий разрешимости задачи Коши для абстрактного уравнения Эйлера---Пуассона---Дарбу//
Дифф. уравн. "--- 2016. "--- {\sl 52}, \No~1. "--- С.~41--59.

\bibitem{Glu13} {\it  Глушак~А.\,В., Романченко~Т.\,Г.} Формулы связи между решениями абстрактных
сингулярных дифференциальных уравнений// Науч. ведом. Белгород.
гос. ун-та. Сер. Мат. Физ. "--- 2016. "--- {\sl 42}, \No~6. "---
С.~36--39.

\bibitem{Glushko} {\it  Глушко~В.\,П.} Линейные вырождающиеся дифференциальные
уравнения. "--- Воронеж: Изд-во ВГУ, 1972.

\bibitem{GKE} {\it Гноенский~Л.\,С., Каменский~Г.\,А., Эльсгольц~Л.\,Э.} Математические основы теории управляемых систем. "--- М.: Наука, 1969.

\bibitem{Grin} {\it  Гринберг~Г.\,А.} Избранные вопросы математической теории электрических и магнитных явлений. "--- М.: АН СССР, 1948.

\bibitem{Gul2} {\it Гулиев~В.\,С.} Интегральные операторы, функциональные пространства и вопросы аппроксимации на группе Гейзенберга. "--- Баку: <<ЭЛМ>>, 1996.

\bibitem{Gul1} {\it  Гулиев~В.\,С.} Функциональные пространства, интегральные операторы и двухвесовые оценки на однородных группах. Некоторые приложения. "--- Баку: Чашыоглы,  1999.

\bibitem{Gur} {\it  Гуревич~М.\,И.} Теория струй идеальной жидкости. "--- М.: Наука, 1979.

\bibitem{Gus1} {\it  Гусейнов~И.\,М.} Об одном операторе преобразования// Мат. заметки. "--- 1997. "--- {\sl 62}, \No~2. "--- С.~206--215.

\bibitem{Gus2} {\it  Гусейнов~И.\,М., Набиев~А.\,А., Пашаев~Р.\,Т.} Операторы преобразования и асимптотические формулы для собственных значений полиноминального пучка операторов
Штурма---Лиувилля// Сиб. мат.~ж. "--- 2000. "--- {\sl 41}, \No~3. "--- С.~554--566.

\bibitem{Dza2} {\it  Джаяни~Г.\,В.} Решение некоторых задач для одного вырождающегося эллиптического уравнения и их приложения к призматическим оболочкам. "---
Тбилиси: Изд-во Тбилис. ун-та, 1982.

\bibitem{Dza1} {\it  Джаяни~Г.\,В.} Уравнение Эйлера---Пауссона---Дарбу. "--- Тбилиси: Изд-во Тбилис. ун-та, 1984.

\bibitem{Gini} {\it  Джини~К.} Средние величины. "--- М.: Статистика, 1970.

\bibitem{Dzh1} {\it  Джрбашян~М.\,М.} Интегральные преобразования и представления функций в комплексной области. "--- М.: Наука, 1966.

\bibitem{Din} {\it  Динь~Х.\,А.} Интегральные уравнения с функцией Лежандра в ядрах в особых случаях// Докл. АН Белорус. ССР. "--- 1989. "--- {\sl 33}, \No~7. "--- С.~591--594.

\bibitem{EPP} {\it  Егоров~И.\,Е., Пятков~С.\,Г., Попов~С.\,В.} Неклассические дифференциально-операторные уравнения. "--- Новосибирск: Наука, 2000.

\bibitem{EgFed} {\it  Егоров~И.\,Е., Федоров~В.~Евс.} Неклассические уравнения математической физики высокого порядка. "--- Новосибирск: Изд-во ВЦ СО РАН, 1995.

\bibitem{ZhM1} {\it  Жегалов~В.\,И., Миронов~А.\,Н.} Дифференциальные уравнения со старшими частными производными. "--- Казань: Казан. мат. об-во, 2001.

\bibitem{ZhM2} {\it  Жегалов~В.\,И., Миронов~А.\,Н., Уткина~Е.\,А.} Уравнения с доминирующей частной производной. "--- Казань: Казан. ун-т, 2014.

\bibitem{Zhit} {\it  Житомирский~Я.\,И.} Задача Коши для систем линейных уравнений в частных производных с дифференциальными операторами типа Бесселя//
Мат. сб. "--- 1955. "--- {\sl 36}, \No~2. "--- С.~299--310.

\bibitem{ZhuSi1} {\it  Жуковская~Н.\,В., Ситник~С.\,М.} Дифференциальные уравнения типа Эйлера дробного порядка// Мат. заметки СВФУ. "--- 2018. "--- {\sl 25}, \No~2. "--- С.~27--39.

\bibitem{Zhu} {\it  Журавлёв~В.\,М.} Нелинейные волны. Точно решаемые задачи. "--- Ульяновск, 2001.

\bibitem{Zai1} {\it Зайцев~В.\,А.} О принципе Гюйгенса для некоторых уравнений с особенностями// Докл. АН СССР. "--- 1978. "--- {\sl 242}, \No~1. "--- C.~28--31.

\bibitem{Zai2} {\it Зайцев~В.\,А.} Слабые лакуны для одномерных строго гиперболических уравнений с постоянными коэффициентами//
Сиб. мат.~ж. "--- 1984. "--- {\sl 25}, \No~4. "--- C.~54--62.

\bibitem{ZMNP} {\it Захаров~В.\,Е.,  Манаков~С.\,В., Новиков~С.\,П., Питаевский~Л.\,П.} Теория солитонов: метод обратной задачи. "--- М.: Наука, 1980.

\bibitem{ZKNE} {\it Зверкин~А.\,М., Каменский~Г.\,А., Норкин~С.\,Б., Эльсгольц~Л.\,Э.}
Дифференциальные уравнения с отклоняющимся аргументом// Усп. мат. наук. "---
1962. "--- {\sl 17}, \No~2. "--- C.~77--164.

\bibitem{18} {\it  Ибрагимов~А.\,И.}  О поведении в окрестности граничных точек и теоремы об устранимых множествах для эллиптических уравнений второго порядка
с непрерывными коэффициентами// Докл. АН СССР. "--- 1980. "--- {\sl 250}, \No~1. "--- С.~25--28.

\bibitem{19} {\it Ивакин~В.\,М.} Видоизмененная задача Дирихле для вырождающихся на границе уравнений и систем// Дифф. уравн. "--- 1982. "--- {\sl 18}, \No~2. "--- С.~319--324.

\bibitem{Iva1} {\it Иванов~Л.\,А.}  О задаче Коши для операторов, распадающихся на множители Эйлера---Пуассона---Дарбу// Дифф. уравн. "--- 1978. "--- {\sl 14}, \No~4. "--- C.~736--739.

\bibitem{Iva2} {\it Иванов~Л.\,А.} Задача Коши для некоторых операторов с особенностями// Дифф. уравн. "--- 1982. "--- {\sl 18}, \No~6. "--- C.~1020--1028.

\bibitem{20} {\it Ильин~В.\,А.} Ядра дробного порядка// Мат. сб. "--- 1957. "--- {\sl 41}, \No~4. "--- С.~459--480.

\bibitem{21} {\it Казарян~К.\,С.} О задаче Дирихле в весовой метрике// В сб.: <<Применение методов теории функций и функционального анализа к задачам математической физики>>. "---
Ереван: Изд-во ЕГУ, 1982. "--- С.~134--136.

\bibitem{KaSk} {\it Каменский~Г.\,А., Скубачевский~А.\,Л.} Линейные краевые задачи для дифференциально-разностных уравнений. "--- М.: Изд-во МАИ, 1992.

\bibitem{Kap1} {\it Капцов~О.\,В.} Методы интегрирования уравнений с частными производными. "--- М.: Физматлит, 2009.

\bibitem{KaSa} {\it Карапетянц~Н.\,К., Самко~С.\,Г.} Уравнения с инволютивными операторами и их приложения. "--- Ростов-на-Дону: Изд-во Ростов. ун-та, 1988.

\bibitem{KarST3} {\it Каримов~Ш.\,Т.} Многомерный оператор Эрдейи---Кобера и его
приложение к решению задачи Коши для трехмерного гиперболического
уравнения с сингулярными коэффициентами// Узб. мат.~ж. "--- 2013. "--- \No~1. "---  C.~70--80.

\bibitem{KarST4} {\it Каримов~Ш.\,Т.} Об одном методе решения задачи Коши для
обобщённого уравнения Эйлера---Пуассона---Дарбу// Узб. мат.~ж.
"--- 2013. "--- \No~3. "---  C.~57--69.

\bibitem{KarST2} {\it Каримов~Ш.\,Т.} Решение задачи Коши для трехмерного
гиперболического уравнения с сингулярными коэффициентами и со
спектральным параметром// Узб. мат.~ж. "---
2014. "--- \No~2. "---  C.~55--65.

\bibitem{KarST1} {\it Каримов~Ш.\,Т.} Об одном методе решения задачи Коши
для одномерного поливолнового уравнения с сингулярным оператором
Бесселя// Изв. вузов. Сер. Мат. "---  2017. "--- \No~8. "---
C.~27--41.

\bibitem{Kar1} {\it  Карп~Д.\,Б.} Пространства с гипергеометрическими воспроизводящими ядрами
и дробные преобразования типа Фурье// Дисс. к.ф.-м.н. "---
Владивосток, 2000.

\bibitem{S142} {\it  Карп~Д.\,Б., Ситник~С.\,М.} Дробное преобразование Ханкеля и его приложения//
Тез. докладов. Воронеж. весен. мат.
школы (17--23  апреля 1996~г.). Соврем. методы в теор.
краевых задач. <<Понтрягинские чтения-VII>>. "--- Воронеж: ВГУ, 1996. "--- С.~92.

\bibitem{22} {\it  Катрахов~В.\,В.} О задаче на собственные значения для сингулярных эллиптических операторов// Докл. АН СССР. "--- 1972. "--- {\sl 207}, \No~2. "--- С.~284--287.

\bibitem{23} {\it  Катрахов~В.\,В.} К теории уравнений с частными производными с сингулярными коэффициентами// Докл. АН СССР. "--- 1974. "--- {\sl 218}, \No~1. "--- С.~17--20.

\bibitem{24} {\it  Катрахов~В.\,В.} Спектральная функция некоторых сингулярных дифференциальных операторов// Дифф. уравн. "--- 1976. "--- {\sl 12}, \No~7. "--- С.~1256--1266.

\bibitem{25} {\it  Катрахов~В.\,В.} Операторы преобразования в теории одномерных псевдодифференциальных операторов//
В сб.: <<Применение методов теории функций и функционального анализа к задачам математической физики>>. "--- Новосибирск: ИМ СО АН СССР, 1979. "--- С.~72--75.

\bibitem{26} {\it  Катрахов~В.\,В.} Операторы преобразования и псевдодифференциальные операторы// Сиб. мат.~ж. "--- 1980. "--- {\sl 21}, \No~1. "--- С.~86--97.

\bibitem{Kat1} {\it  Катрахов~В.\,В.}  Изометрические операторы преобразования и спектральная функция для одного класса одномерных сингулярных псевдодифференциальных операторов//
Докл. АН СССР. "--- 1980. "--- {\sl  251}, \No~5. "--- С.~1048--1051.

\bibitem{28} {\it  Катрахов~В.\,В.}  Общие краевые задачи для одного класса сингулярных и вырождающихся уравнений// Докл. АН СССР. "--- 1980. "--- {\sl 251}, \No~6. "--- С.~1296--1300.

\bibitem{29} {\it  Катрахов~В.\,В.} Общие краевые задачи для одного класса сингулярных и вырождающихся эллиптических уравнений//
Мат. сб. "--- 1980. "--- {\sl 112}, \No~3. "--- С.~354--379.

\bibitem{Kat2} {\it  Катрахов~В.\,В.}  Об одной краевой задаче для уравнения Пуассона// Докл. АН СССР. "--- 1981. "--- {\sl 259}, \No~5. "--- С.~1041--1045.

\bibitem{30} {\it  Катрахов~В.\,В.} Сингулярные краевые задачи и операторы преобразования//
В сб.: <<Корректные краевые задачи для неклассических уравнений математической физики>>. "--- Новосибирск: ИМ СО АН СССР, 1981. "---  С.~87--91.

\bibitem{32} {\it  Катрахов~В.\,В.} Метод операторов преобразования в теории общих весовых краевых задач для сингулярных и вырождающихся эллиптических уравнений с параметром//
Докл. АН СССР. "--- 1982. "--- {\sl  266}, \No~5. "--- С.~1037--1040.

\bibitem{KatDis} {\it  Катрахов~В.\,В.} Сингулярные эллиптические краевые задачи. Метод операторов преобразования// Дисс. д.ф.-м.н. "--- Новосибирск, 1989.

\bibitem{Kat3} {\it  Катрахов~В.\,В.}  Об одной сингулярной краевой задаче для уравнения Пуассона// Мат. сб. "--- 1991. "--- {\sl 182}, \No~6. "--- С.~849--876.

\bibitem{Kat4} {\it  Катрахов~В.\,В.} Сингулярные краевые задачи для некоторых эллиптических уравнений в областях с угловыми точками//
Докл. АН СССР. "--- 1991. "--- {\sl 316},  \No~5. "--- С.~1047--1050.

\bibitem{KaKa} {\it  Катрахов~В.\,В., Катрахова~А.\,А.} Формула Тэйлора с оператором
Бесселя для функций одной и двух переменных// Деп. в
ВИНИТИ. "--- Воронеж, 1982.

\bibitem{33} {\it  Катрахов~В.\,В., Киприянов~И.\,А.} Степени сингулярного эллиптического оператора//
В сб.: <<Теория кубатурных формул и приложения функционального анализа к задачам математической физики>>. "--- Новосибирск, 1980. "--- С.~60--80.

\bibitem{S1} {\it  Катрахов~В.\,В., Ситник~С.\,М.}  Краевая задача для  стационарного  уравнения  Шрёдингера  с сингулярным потенциалом//
Докл. АН СССР. "--- 1984. "--- {\sl 278}, \No~4. "--- С.~797--799.

\bibitem{S5}  {\it  Катрахов~В.\,В., Ситник~С.\,М.}
Метод факторизации в теории операторов преобразования//
В сб.: <<Мемориальный сборник памяти Бориса Алексеевича Бубнова:
неклассические уравнения и уравнения смешанного типа>>.
"--- Новосибирск, 1990.
"--- С.~104--122.

\bibitem{S7} {\it  Катрахов~В.\,В., Ситник~С.\,М.} Композиционный метод построения $B$-эллиптических,
$B$-па\-ра\-бо\-ли\-чес\-ких и $B$-гиперболических операторов
преобразования// Докл. РАН. "--- 1994. "---{\sl 337}, \No~3.
"--- C.~307--311.

\bibitem{S8} {\it  Катрахов~В.\,В., Ситник~С.\,М.}
Оценки решений Йоста  для одномерного  уравнения  Шрёдингера с сингулярным потенциалом//
 Докл. РАН. "--- 1995. "--- {\sl 340}, \No~1. "--- С.~18--20.

\bibitem{Kach1} {\it  Качалов~А.\,П., Курылёв~Я.\,В.} Метод операторов преобразования в обратной задаче рассеяния, одномерный Штарк-эффект//
Зап. науч. сем. ЛОМИ. "--- 1989. "--- {\sl 179}. "--- C.~73--87.

\bibitem{35} {\it  Квядарас~Б.\,В.} Решение задачи Дирихле для вырожденного эллиптического уравнения//
В сб.: <<Дифференциальные уравнения с частными производными: труды конференции по дифференциальным уравнениям и вычислительной математике>>. "---
Новосибирск: Наука, СО, 1980. "--- С.~35--36.

\bibitem{37} {\it  Келдыш~М.\,В.} О разрешимости и устойчивости задачи Дирихле//  Усп. мат. наук. "--- 1941. "--- {\sl  8}. "--- С.~171--292.

\bibitem{Kel} {\it  Келдыш~М.\,В.} О некоторых случаях вырождения уравнений эллиптического типа на границе области//  Докл. АН СССР. "--- 1951. "--- {\sl  77}, \No~1. "--- С.~181--183.

\bibitem{KSZ}  {\it  Килбас~А.\,А., Сайго~М., Жук~В.\,А.} О композиции операторов обобщённого дробного интегрирования с дифференциальным оператором осесимметрической теории потенциала//
Дифф. уравн. "--- 1991. "--- {\sl  27}, \No~9. "--- С.~1640--1642.

\bibitem{KiSk2} {\it  Килбас~А.\,А., Скоромник~О.\,В.} Решение многомерного интегрального уравнения первого рода с функцией Лежандра по пирамидальной области//
Докл. РАН. "--- 2009. "--- {\sl  429}, \No~4. "--- С.~442--446.

\bibitem{Kip2} {\it  Киприянов~И.\,А.} Преобразования Фурье---Бесселя и теоремы вложения для весовых классов//  Тр. МИАН. "--- 1967. "--- {\sl  89}. "--- С.~130--213.

\bibitem{39} {\it  Киприянов~И.\,А.} Краевые задачи для сингулярных эллиптических операторов в частных производных//  Докл. АН СССР. "--- 1970. "--- {\sl  195}, \No~1. "--- С.~32--35.

\bibitem{40} {\it  Киприянов~И.\,А.} Об одном классе сингулярных эллиптических операторов// Дифф. уравн. "--- 1971. "--- {\sl  7}, \No~11. "--- С.~2065--2077.

\bibitem{41} {\it  Киприянов~И.\,А.} Об одном классе сингулярных эллиптических уравнений//  Сиб. мат.~ж. "--- 1973. "--- {\sl  14}, \No~3. "--- С.~560--568.

\bibitem{Kip1} {\it  Киприянов~И.\,А.} Сингулярные эллиптические краевые задачи. "--- М.: Наука---Физматлит, 1997.

\bibitem{KipIv1} {\it  Киприянов~И.\,А., Иванов~Л.\,А.} О лакунах для некоторых классов уравнений с особенностями// Мат. сб. "--- 1979. "--- {\sl  110}, \No~2. "--- С.~235--250.

\bibitem{42} {\it  Киприянов~И.\,А., Иванов~Л.\,А.} Уравнение Эйлера---Пуассона---Дарбу в римановом пространстве// Докл. АН СССР. "--- 1981. "--- {\sl  260}, \No~4. "--- С.~790--794.

\bibitem{KipIv3} {\it  Киприянов~И.\,А., Иванов~Л.\,А.} Задача Коши для уравнения Эйлера---Пуассона---Дарбу в однородном симметрическом римановом пространстве.~I//
Тр. МИАН. "--- 1984. "--- {\sl  170}. "--- С.~139--147.

\bibitem{KipIv2} {\it  Киприянов~И.\,А., Иванов~Л.\,А.} Задача Коши для уравнения Эйлера---Пуассона---Дарбу в симметрическом пространстве//
Мат. сб. "--- 1984. "--- {\sl  124}, \No~1. "--- С.~45--55.

\bibitem{KipIv4} {\it  Киприянов~И.\,А., Иванов~Л.\,А.} Представление Даламбера и равнораспределение энергии// Дифф. уравн. "--- 1990. "--- {\sl  26}, \No~3. "--- С.~458--464.

\bibitem{KiKa1} {\it  Киприянов~И.\,А., Катрахов~В.\,В.} Об одном классе многомерных сингулярных псевдодифференциальных операторов//  Мат. сб. "--- 1977. "--- {\sl  104}, \No~1. "--- С.~49--68.

\bibitem{44} {\it  Киприянов~И.\,А., Катрахов~В.\,В.} Краевая задача для эллиптических уравнений второго порядка при наличии особенностей в изолированных граничных точках//
Докл. АН СССР. "--- 1984. "--- {\sl 276}, \No~2. "--- С.~274--276.

\bibitem{KiKa2} {\it  Киприянов~И.\,А., Катрахов~В.\,В.} Об одной сингулярной эллиптической краевой задаче в областях на сфере// Препринт ИПМ ДВО РАН. "--- 1989.

\bibitem{KiKa3} {\it  Киприянов~И.\,А., Катрахов~В.\,В.} Сингулярные краевые задачи для некоторых эллиптических уравнений высших порядков//  Препринт ИПМ ДВО РАН. "--- 1989.

\bibitem{KiKa4} {\it  Киприянов~И.\,А., Катрахов~В.\,В.} Об одной краевой задаче для эллиптических уравнений второго порядка в областях на сфере//
Докл. АН СССР. "--- 1990. "--- {\sl  313}, \No~3. "--- С.~545--548.

\bibitem{45} {\it  Киприянов~И.\,А., Ключанцев~М.\,И.} О ядрах Пуассона для краевых задач с дифференциальным оператором Бесселя//
В сб.: <<Дифференциальные уравнения с частными производными>>. "--- М., 1970. "---  C.~119--134.

\bibitem{KipKo1} {\it  Киприянов~И.\,А., Кононенко~В.\,И.} О фундаментальных решениях уравнений в частных производных с дифференциальным оператором Бесселя//
Докл. АН СССР. "--- 1966. "--- {\sl 170}, \No~2. "--- С.~261--264.

\bibitem{KipKo2} {\it  Киприянов~И.\,А., Кононенко~В.\,И.} Фундаментальные решения $B$-эллиптических уравнений// Дифф. уравн. "--- 1967. "--- {\sl 3}, \No~1. "--- C.~114--129.

\bibitem{KipKo3} {\it  Киприянов~И.\,А., Кононенко~В.\,И.} О фундаментальных решениях некоторых сингулярных уравнений в частных производных//
Дифф. уравн. "--- 1969. "--- {\sl 5}, \No~8. "--- C.~1470--1483.

\bibitem{KipKu1} {\it  Киприянов~И.\,А., Куликов~А.\,А.} Фундаментальные решения $B$-гипоэллиптических уравнений// Дифф. уравн. "--- 1991. "--- {\sl 27}, \No~8. "--- C.~1387--1395.

\bibitem{Kli2} {\it  Климентов~С.\,Б.} Классы Харди обобщённых аналитических функций//  Изв. вузов. Сев.-Кавказ. рег. Сер. Естеств. науки. "--- 2003. "--- \No~3. "--- С.~6--10.

\bibitem{Kli3} {\it  Климентов~С.\,Б.} Классы Смирнова обобщённых аналитических функций//  Изв. вузов. Сев.-Кавказ. рег. Сер. Естеств. науки. "--- 2005. "---\No~1. "--- С.~13--17.

\bibitem{Kli5} {\it  Климентов~С.\,Б.} Классы ВМО обобщённых аналитических функций// Владикавказ. мат.~ж. "--- 2006. "--- {\sl  8}, \No~1. "--- С.~27--39.

\bibitem{Kli1} {\it  Климентов~С.\,Б.} Граничные свойства обобщённых аналитических функций. "--- Владикавказ: Изд. Южного мат. ин-та ВНЦ РАН и РСО-А, 2014.

\bibitem{Kly1} {\it  Ключанцев~М.\,И.} О построении $r$-чётных решений сингулярных дифференциальных уравнений//  Докл. АН СССР. "--- 1975. "--- {\sl  224}, \No~5. "--- С.~1004--1007.

\bibitem{Kly2} {\it  Ключанцев~М.\,И.} Интегралы дробного порядка и сингулярные краевые задачи//   Дифф. уравн. "--- 1976. "--- {\sl  12}, \No~6. "--- С.~983--990.

\bibitem{KF} {\it  Колмогоров~А.~Н, Фомин~С.\,В.} Элементы теории функций и функционального анализа. "--- М.: Наука, 1981.

\bibitem{Kol1} {\it  Колтон~Д., Кресс~Р.} Методы интегральных уравнений в теории рассеяния. "--- М.: Мир, 1987.

\bibitem{Kor1} {\it  Коробейник~Ю.\,Ф.} Операторы сдвига на числовых семействах. "--- Ростов-на-Дону: Изд. Ростов. ун-та, 1983.

\bibitem{Kor2} {\it  Коробейник~Ю.\,Ф.}  О разрешимости в комплексной области некоторых общих классов линейных интегральных уравнений. "--- Ростов-на-Дону: Изд. Ростов. ун-та, 2005.

\bibitem{Koch1} {\it  Кочубей~А.\,Н.} Задача Коши для эволюционных уравнений дробного порядка//  Дифф. уравн. "--- 1989. "--- {\sl  25}, \No~8. "--- С.~1359--1369.

\bibitem{Koch2} {\it  Кочубей~А.\,Н.} Диффузия дробного порядка//  Дифф. уравн. "--- 1990. "--- {\sl  26}, \No~4. "--- С.~660--770.

\bibitem{KGS} {\it  Кошляков~Н.\,С., Глинер~Э.\,Б., Смирнов~М.\,М.} Уравнения в частных производных математической физики. "--- М.: Высшая школа, 1962.

\bibitem{Kra2} {\it  Кравченко~В.\,Ф. $($ред.$)$} Цифровая обработка сигналов и изображений в радиофизических приложениях. "--- М.: Физматлит, 2007.

\bibitem{Kra1} {\it  Кравченко~В.\,Ф., Рвачёв~В.\,Л.}  Алгебра логики, атомарные функции и вейвлеты в физических приложениях. "--- М.: Физматлит, 2006.

\bibitem{KPS} {\it  Крейн~С.\,Г., Петунин~Ю.\,И., Семёнов~Е.\,М.} Интерполяция линейных операторов. "--- М.: Наука, 1978.

\bibitem{48} {\it  Кудрявцев~Л.\,Д.} Прямые и обратные теоремы вложения. Приложения к решению вариационным методом эллиптических уравнений. "--- М.: Наука, 1959.

\bibitem{KN} {\it  Кудрявцев~Л.\,Д., Никольский~С.\,М.} Пространства дифференцируемых функций многих переменных и теоремы вложения//
Итоги науки и техн. Соврем. пробл. мат. Фундам. направл. "--- М.: ВИНИТИ, 1988. "--- {\sl 26}. "--- С.~5--157.

\bibitem{Kuz3} {\it  Кузнецов~Н.\,В.} О собственных функциях одного интегрального уравнения//  Зап. науч. сем. ЛОМИ. "--- 1970. "--- {\sl  17}, \No~3. "--- С.~66--149.

\bibitem{Kuz2} {\it  Кузнецов~Н.\,В.} Гипотеза Петерсона для параболических форм веса нуль и гипотеза Линника. Суммы сумм Клоостермана//
Мат. сб. "--- 1980. "--- {\sl  111}, \No~3. "--- С.~334--383.

\bibitem{Kuz1} {\it  Кузнецов~Н.\,В.} Формулы следа и некоторые их приложения в теории чисел. "--- Владивосток: Дальнаука, 2003.

\bibitem{49} {\it  Курант~Р.} Уравнения с частными производными. "--- М.: Мир, 1979.

\bibitem{Kus} {\it  Кусраев~А.\,Г.} Мажорируемые операторы. "--- М.: Наука, 2003.

\bibitem{Lav1} {\it Лаврентьев~М.\,М.} Одномерные обратные задачи математической физики. "--- Новосибирск: Наука, 1982.

\bibitem{51} {\it  Ладыженская~O.\,A., Уральцева~Н.\,Н.} Линейные и квазилинейные уравнения эллиптического типа. "--- М.: Наука, 1973.

\bibitem{Laks1} {\it  Лакс~П., Филлипс~Р.} Теория рассеяния. "--- М.: Мир, 1971.

\bibitem{Laks2} {\it  Лакс~П.} Теория рассеяния для автоморфных функций. "--- М.: Мир, 1979.

\bibitem{53} {\it  Ландис~Е.\,М.} Уравнение второго порядка эллиптического и параболического типов. "--- М.: Наука, 1971.

\bibitem{Lan} {\it  Ландис~Е.\,М.} Задачи Е.\,М.~Ландиса//  Усп. мат. наук. "--- 1982. "--- {\sl  37}, \No~6. "--- С.~278--281.

\bibitem{54} {\it  Ландкоф~Н.\,С.} Основы современной теории потенциала. "--- М.: Наука, 1966.

\bibitem{Lar6} {\it  Ларин~А.\,A.} О спектральных разложениях, отвечающих самосопряжённым расширениям некоторых сингулярных эллиптических операторов//
Докл. АН СССР. "--- 1987. "--- {\sl  293}, \No~2. "--- С.~309--312.

\bibitem{Lar5} {\it  Ларин~А.\,A.} О свойствах собственных функций некоторых сингулярных эллиптических операторов//
Дифф. уравн. "--- 1991. "--- {\sl  27}, \No~5. "--- С.~849--856.

\bibitem{Lar4} {\it  Ларин~А.\,A.} Об ограниченности степеней самосопряжённых расширений сингулярных эллиптических операторов, действующих в весовых классах//
Дифф. уравн. "--- 1992. "--- {\sl  28}, \No~3. "--- С.~528--529.

\bibitem{Lar3} {\it  Ларин~А.\,A.} О представлении решений одного сингулярного эллиптического уравнения второго порядка в окрестности угловой точки//
Дифф. уравн. "--- 2000. "--- {\sl  36}, \No~4. "--- С.~566--568.

\bibitem{Lar2} {\it  Ларин~А.\,A.} Об одной краевой задаче в плоском угле для сингулярного эллиптического уравнения второго порядка//
Дифф. уравн. "--- 2000. "--- {\sl  36}, \No~12. "--- С.~1687--1694.

\bibitem{Lar1} {\it  Ларин~А.\,A.} О теореме сужения на сферическую поверхность для
преобразований Фурье---Бесселя//  Докл. Адыгской (Черкесской) межд. акад. наук. "--- 2014. "--- {\sl  16}, \No~3.
"--- С.~22--29.

\bibitem{Levin2} {\it  Левин~Б.\,Я.} Преобразования типа Фурье и Лапласа при помощи решений дифференциального уравнения второго порядка//
Докл. АН СССР. "--- 1956. "--- {\sl  106}, \No~2. "--- С.~187--190.

\bibitem{Levin1} {\it  Левин~Б.\,Я.} Распределение корней целых функций. "--- М.: ГИТТЛ, 1956.

\bibitem{Lev1} {\it  Левитан~Б.\,М.} Разложение по собственным функциям дифференциальных уравнений второго порядка. "--- М.: Гостехиздат, 1950.

\bibitem{Lev7} {\it  Левитан~Б.\,М.}  Разложения по функциям Бесселя в ряды и интегралы Фурье//  Усп. мат. наук. "--- 1951. "--- {\sl 6}, \No~2. "--- С.~102--143.

\bibitem{Lev4} {\it  Левитан~Б.\,М.}  Почти-периодические функции. "--- М.: ГИТТЛ, 1953.

\bibitem{Lev2} {\it  Левитан~Б.\,М.}  Операторы обобщённого сдвига и некоторые их применения. "--- М.: ГИФМЛ, 1962.

\bibitem{Lev3} {\it  Левитан~Б.\,М.}  Теория операторов обобщённого сдвига. "--- М.: Наука, 1973.

\bibitem{Lev5}  {\it  Левитан~Б.\,М.}  Обратные задачи Штурма---Лиувилля. "--- М.: Наука, 1984.

\bibitem{Lev8} {\it  Левитан~Б.\,М., Повзнер~А.\,Я.} Дифференциальные уравнения  Штурма---Лиувилля на полуоси и теорема Планшереля//
Докл. АН СССР. "--- 1946. "--- {\sl  52}, \No~6. "--- С.~483--486.

\bibitem{Lev6} {\it  Левитан~Б.\,М., Саргсян~И.\,С.}  Операторы Штурма---Лиувилля и Дирака. "--- М.: Наука, 1988.

\bibitem{Lei1} {\it  Лейзин~М.\,А.} К теоремам вложения для одного класса сингулярных дифференциальных операторов в полупространстве//
Дифф. уравн. "--- 1976. "--- {\sl  12}, \No~6. "--- С.~1073--1083.

\bibitem{Lei2} {\it  Лейзин~М.\,А.} О вложении некоторых весовых классов//  В сб.: <<Методы решений операторных уравнений>>. "--- Воронеж: Изд-во ВГУ, 1978. "--- С.~96--103.

\bibitem{Leo} {\it  Леонтьев~А.\,Ф.} Оценка роста решения одного дифференциального уравнения при больших по модулю значениях параметра и её применения к некоторым вопросам теории функций//
Сиб. мат.~ж. "--- 1960. "--- {\sl  1}, \No~3. "--- С.~456--487.

\bibitem{Ler} {\it Лернер~М.\,Е.} Принципы максимума для уравнений гиперболического типа и новые свойства функции Римана. "--- Самара: Самар. гос. тех. ун-т, 2001.

\bibitem{56} {\it  Лизоркин~П.\,И.} обобщённое лиувиллевское дифференцирование и метод мультипликаторов в теории вложений классов дифференцируемых функций//
Тр. МИАН. "--- 1969. "--- {\sl  105}. "--- С.~89--167.

\bibitem{Lis} {\it  Лизоркин~П.\,И.} Классы функций, построенные на основе усреднений по сферам. Случай пространств Соболева//  Тр. МИАН. "--- 1990. "--- \No~192. "--- С.~122--139.

\bibitem{57} {\it  Лизоркин~П.\,И., Никольский~С.\,М.}  Эллиптическое уравнение с вырождением. Вариационный метод//
Докл. АН СССР. "--- 1981. "--- {\sl  257}, \No~1. "--- С.~42--45.

\bibitem{58}  {\it  Лизоркин~П.\,И., Никольский~С.\,М.}  Эллиптические уравнения с вырождением. Дифференциальные свойства решений//
Докл. АН СССР. "--- 1981. "--- {\sl  257}, \No~2. "--- С.~278--282.

\bibitem{59}  {\it  Лизоркин~П.\,И., Никольский~С.\,М.}  Коэрцитивные свойства эллиптического уравнения с сильным вырождением (случай обобщённых решений)//
Докл. АН СССР. "--- 1981. "--- {\sl  259}, \No~1. "--- С.~28--30.

\bibitem{62}{\it Лионс~Ж.-Л.,  Мадженес~Э.} Неоднородные граничные задачи и их приложения. "--- М.: Мир, 1971.

\bibitem{Lit1}{\it Литвинчук~Г.\,С.} Краевые задачи и сингулярные интегральные уравнения со сдвигом. "--- М.: Наука, 1977.

\bibitem{Luke3}{\it Люк~Ю.} Специальные математические функции и их аппроксимации. "--- М.: Мир, 1980.

\bibitem{Lyah3} {\it Ляхов~Л.\,Н.} Обращение $B$-потенциалов// Докл. АН СССР. "--- 1991. "--- {\sl 321}, \No~3. "---
C.~466--469.

\bibitem{Lyah1}{\it Ляхов~Л.\,Н.} Весовые сферические функции и потенциалы Рисса, порожденные обобщённым сдвигом. "--- Воронеж: ВГТА, 1997.

\bibitem{Lyah2} {\it Ляхов~Л.\,Н.} $B$-гиперсингулярные интегралы и их приложения к описанию функциональных классов Киприянова и к интегральным уравнениям с $B$-потенциальными ядрами.
 "--- Липецк: Изд-во ЛГПУ, 2007.

\bibitem{LPSh1} {\it  Ляхов~Л.\,Н., Половинкин~И.\,П., Шишкина~Э.\,Л.} Об одной задаче И.\,А.~Киприянова для сингулярного ультрагиперболического уравнения//
Дифф. уравн. "--- 2014. "--- {\sl  50}, \No~4. "--- С.~516--528.

\bibitem{LPSh2} {\it  Ляхов~Л.\,Н., Половинкин~И.\,П., Шишкина~Э.\,Л.} Формулы решения задачи Коши для сингулярного волнового уравнения с оператором Бесселя по времени//
Докл. РАН. "--- 2014. "--- {\sl  459}, \No~5. "--- С.~533--538.

\bibitem{LShFrac} {\it  Ляхов~Л.\,Н., Шишкина~Э.\,Л.} Дробные производные и интегралы и их приложения. "--- Воронеж: Изд-во ВГУ, 2011.

\bibitem{Mal1} {\it  Маламуд~М.\,М.} Об операторах преобразования для обыкновенных дифференциальных уравнений высших порядков//
В сб.: <<Математический анализ и теория вероятностей>>. "--- Киев: Наукова думка, 1978. "--- С.~108--111.

\bibitem{Mal2} {\it  Маламуд~М.\,М.}  Необходимые условия существования оператора преобразования для уравнений высших порядков//
Функц. анализ и его прилож. "--- 1982. "--- {\sl  16}, \No~3. "--- С.~74--75.

\bibitem{Mal3} {\it  Маламуд~М.\,М.} К вопросу об операторах преобразования//  Препринт ИМ АН УССР. "---  Киев, 1984.

\bibitem{Mal4} {\it  Маламуд~М.\,М.} Операторы преобразования для  уравнений высших порядков// Мат. физ. и нелин. мех. "--- 1986. "--- \No~6. "--- С.~108--111.

\bibitem{Mal5} {\it  Маламуд~М.\,М.} К вопросу об операторах преобразования для обыкновенных дифференциальных уравнений//
Тр. Моск. мат. об-ва. "--- 1990. "--- {\sl  53}. "--- С.~68--97.

\bibitem{Marich1}{\it Маричев~О.\,И.} Метод вычисления интегралов от специальных функций. "--- Минск: Наука и техника, 1978.

\bibitem{MaKiRe}{\it Маричев~О.\,И., Килбас~А.\,А., Репин~О.\,А.} Краевые задачи для уравнений в частных производных с разрывными коэффициентами. "---
Самара: Изд-во Самар. гос. эконом. ун-та, 2008.

\bibitem{Mar3} {\it Марченко~В.\,А.} Некоторые вопросы теории дифференциального оператора второго порядка//
Докл. АН СССР. "--- 1950. "--- {\sl 72}, \No~3. "--- С.~457--460.

\bibitem{Mar4} {\it Марченко~В.\,А.} Операторы преобразования//
Докл. АН СССР. "--- 1950. "--- {\sl  74}, \No~2. "--- С.~185--188.

\bibitem{Mar5} {\it Марченко~В.\,А.} О формулах обращения, порождаемых линейным дифференциальным оператором второго порядка//
Докл. АН СССР. "--- 1950. "--- {\sl  74}, \No~4. "--- С.~657--660.

\bibitem{Mar6} {\it Марченко~В.\,А.} Некоторые вопросы теории одномерных дифференциальных операторов второго порядка.~I//
Тр. Моск. мат. об-ва. "--- 1952. "--- {\sl 1}. "--- С.~327--420.

\bibitem{Mar7} {\it Марченко~В.\,А.} Некоторые вопросы теории одномерных дифференциальных операторов второго порядка.~II//
Тр. Моск. мат. об-ва. "--- 1953. "--- {\sl 2}. "--- С.~3--82.

\bibitem{Mar1}{\it Марченко~В.\,А.} Спектральная теория операторов
Штурма---Лиувилля. "--- Киев: Наукова Думка, 1972.

\bibitem{Mar2}{\it Марченко~В.\,А.} Операторы
Штурма---Лиувилля и их приложения. "--- Киев: Наукова Думка, 1977.

\bibitem{Mar8}{\it Марченко~В.\,А.} Нелинейные уравнения и операторные алгебры. "--- Киев: Наукова Думка, 1986.

\bibitem{Mar9}{\it Марченко~В.\,А.} Обобщённый сдвиг, операторы преобразования и обратные задачи// В сб.: <<Математические события ХХ века>>. "--- М.: Фазис, 2003.

\bibitem{Mat1}  {\it  Матiйчук~М.\,I.} Параболiчнi сингулярнi крайовi задачi. "---  Ки\"\iв: Iн-т математики НАН Укра\"\iни, 1999.

\bibitem{Mat2}  {\it  Матiйчук~М.\,I.} Параболічні та еліптичні крайові задачі з особливостями. "--- Чернівці: Прут, 2003.

\bibitem{Mats} {\it  Мацаев~В.\,И.} О существовании оператора преобразования для дифференциальных уравнений высших порядков//
Докл. АН СССР. "--- 1960. "--- {\sl  130}, \No~3. "--- С.~499--502.

\bibitem{S24}{\it  Мехрез~Х., Ситник~С.\,М.} Монотонность отношений некоторых гипергеометрических функций// Сиб. электрон. мат. изв. "--- 2016. "--- {\sl  13}. "---  С.~260--268.

\bibitem{Mesh1} {\it  Мешков~В.\,З.} Весовые дифференциальные неравенства и их применение для оценок скорости убывания на бесконечности решений эллиптических уравнений второго порядка//
Тр. МИАН. "--- 1989. "--- {\sl 190}. "--- С.~139--158.

\bibitem{Mesh2} {\it  Мешков~В.\,З.}  О возможной скорости убывания на бесконечности решений уравнений в частных производных второго порядка//
Мат. сб. "--- 1991. "--- {\sl  182}, \No~3. "--- С.~364--383.

\bibitem{Me1} {\it Мещеряков~В.\,В.} Дифференциально-разностные операторы, ассоциированные с системами корней коксетеровского типа// Дисс. к.ф.-м.н. "--- Коломна, 2008.

\bibitem{Miz} {\it Мизохата~С.} Теория уравнений с частными производными. "--- М.: Мир, 1977.

\bibitem{Miran} {\it Миранда~К.} Уравнения с частными производными эллиптического типа. "--- М.: Мир, 1957.

\bibitem{Moi} {\it Моисеев~Е.\,И.} Уравнения смешанного типа со спектральным параметром. "--- М.: Изд-во МГУ, 1988.

\bibitem{Mur9} {\it  Муравник~А.\,Б.}  О стабилизации решений некоторых сингулярных квазилинейных параболических задач//  Мат. заметки. "--- 2003. "--- {\sl  74}, \No~6. "--- С.~858--865.

\bibitem{Mur8} {\it  Муравник~А.\,Б.}  О стабилизации решений сингулярных эллиптических уравнений//  Фундам. и прикл. мат. "--- 2006. "--- {\sl  12}, \No~4. "--- С.~169--186.

\bibitem{Mur} {\it  Муравник~А.\,Б.} Функционально-дифференциальные параболические уравнения: интегральные представления и качественные свойства решений задачи Коши//
Соврем. мат. Фундам. направл. "--- 2014. "--- {\sl  52}. "--- С.~3--141.

\bibitem{Mys} {\it  Мышкис~А.\,Д.} Линейные дифференциальные уравнения с запаздывающим аргументом. "--- М.---Л.: Гостехиздат, 1951.

\bibitem{Nai} {\it Наймарк~М.\,А.} Линейные дифференциальные операторы. "--- М.: Наука, 1969.

\bibitem{65} {\it Нарасимхан~Р.} Анализ на действительных и комплексных многообразиях. "--- М.: Мир, 1971.

\bibitem{Nat1} {\it Наттерер~Ф.} Математические аспекты компьютерной томографии. "--- М.: Мир, 1990.

\bibitem{Nah1} {\it Нахушев~А.\,М.} Уравнения математической биологии. "--- М.: Высшая Школа, 1995.

\bibitem{Nah2}  {\it Нахушев~А.\,М.}  Элементы дробного исчисления и их применение. "--- Нальчик, 2000.

\bibitem{Nah3}  {\it Нахушев~А.\,М.} Дробное исчисление и его применение. "--- М.: Физматлит, 2003.

\bibitem{Nah4}  {\it Нахушев~А.\,М.} Нагруженные уравнения и их применения. "--- М.: Наука, 2012.

\bibitem{Nizh1}   {\it Нижник~Л.\,П.} Обратная нестационарная задача теории рассеяния. "--- Киев: Наукова Думка, 1973.

\bibitem{Nizh2}  {\it Нижник~Л.\,П.}  Обратные задачи рассеяния для гиперболических уравнений. "--- Киев:  Наукова думка, 1990.

\bibitem{66} {\it  Никольский~С.\,М.} Приближение функций многих переменных и теоремы вложения. "--- М.: Наука, 1977.

\bibitem{67} {\it  Никольский~С.\,М.} Вариационная проблема для уравнения эллиптического типа с вырождением на границе//  Тр. МИАН. "--- 1979. "--- {\sl  150}. "--- С.~212--238.

\bibitem{68} {\it  Никольский~С.\,М., Лизоркин~П.\,И.} О некоторых неравенствах для функций из весовых классов и краевых задачах с сильным вырождением на границе//
Докл. АН СССР. "--- 1964. "--- {\sl  159}, \No~3. "--- С.~512--515.

\bibitem{Nov} {\it  Новоженова~О.\,А.}  Биография и научные труды Алексея Никифоровича Герасимова. О линейных операторах, упруго-вязкости, элевтерозе и дробных производных. "---
М.: Перо, 2018.

\bibitem{69} {\it  Новрузов~А.\,А.} О задачах Дирихле для эллиптических уравнений второго порядка//  Докл. АН СССР. "--- 1979. "--- {\sl  246}, \No~1. "--- С.~11--14.

\bibitem{Nogin} {\it Ногин~В.\,А., Сухинин~Е.\,В.} Обращение и описание гиперболических потенциалов с $L_p$-плотностями// Докл. РАН. "--- 1993. "---  {\sl  329}, \No~5. "--- С.~550--552.

\bibitem{New}  {\it  Ньюэлл~А.} Солитоны в математике и физике. "--- М.: Мир, 1989.

\bibitem{OlRa} {\it  Олейник~О.\,А., Радкевич~Е.\,В.} Уравнения с неотрицательной характеристической
формой. "--- М.: МГУ, 2010.

\bibitem{Ome} {\it  Омельченко~А.\,В.} Методы интегральных преобразований в задачах математической физики. "--- М.: МЦНМО, 2010.

\bibitem{Os}  {\it  Осипов~В.\,Ф.} Почти периодические функции Бора---Френеля. "--- СПб.: Изд-во С.-Петерб. ун-та, 1992.

\bibitem{Pas1} {\it  Пасенчук~А.\,Э.} Абстрактные сингулярные операторы. "--- Новочеркасск, 1993.

\bibitem{Plat2} {\it  Платонов~С.\,С.}   Обобщённые сдвиги Бесселя и некоторые задачи теории приближения функций в метрике $L_2$.~1//
Тр. ПетрГУ. Сер. Мат. "--- 2000. "--- {\sl 7}. "--- С.~70--82.

\bibitem{Plat3} {\it  Платонов~С.\,С.} Обобщённые сдвиги Бесселя и некоторые задачи теории приближения функций в метрике $L_2$.~2//
Тр. ПетрГУ. Сер. Мат. "--- 2001. "--- {\sl 8}. "---  С.~20--36.

\bibitem{Plat1} {\it  Платонов~С.\,С.} Гармонический анализ Бесселя и приближение функций на полупрямой//  Изв. РАН. Сер. Мат. "--- 2007. "--- {\sl  71}, \No~5. "--- С.~149--196.

\bibitem{Povz} {\it  Повзнер~А.\,Я.} О дифференциальных уравнениях типа Штурма---Лиувилля на полуоси//  Мат. сб. "--- 1948. "--- {\sl  23}, \No~1. "--- С.~3--52.

\bibitem{Pol1} {\it  Положий~Г.\,Н.} Уравнения математической физики. "--- М.: Высшая школа, 1964.

\bibitem{Pol2} {\it  Положий~Г.\,Н.} Обобщение теории аналитических функций комплексного переменного. $P$-ана\-ли\-ти\-чес\-кие и $(P,Q)$-аналитические функции и некоторые их применения.
"--- Киев: Изд-во КГУ, 1965.

\bibitem{Pol3} {\it  Положий~Г.\,Н.}  Теория и применение $p$-аналитических функций. "--- Киев: Наукова думка, 1973.

\bibitem{ZaiPol} {\it  Полянин~А.\,Д., Зайцев~В.\,Ф.} Cправочник по нелинейным уравнениям
математической физики: точные решения. "--- М.: Физматлит, 2002.

\bibitem{Pol} {\it  Поляцкий~В.\,Т.} О свойствах решений некоторого уравнения//  Усп. мат. наук. "--- 1965. "--- {\sl  17}, \No~4. "--- С.~119--124.

\bibitem{71} {\it  Пресдорф~З.} Некоторые классы сингулярных уравнений. "--- М.: Мир, 1979.

\bibitem{PBM123} {\it  Прудников~А.\,П., Брычков~Ю.\,А., Маричев~О.\,И.} Интегралы и ряды. Т.~1,\,2,\,3. "--- М.: Наука, 1981, 1983, 1986.

\bibitem{PBM} {\it  Прудников~А.\,П., Брычков~Ю.\,А., Маричев~О.\,И.} Вычисление интегралов и преобразование Меллина//
Итоги науки и техн. Мат. анализ. "--- 1989. "--- {\sl  27}. "--- С.~3--146.

\bibitem{Pshu2} {\it  Псху~А.\,В.} Интегральные преобразования с функцией Райта в ядре// Докл. Адыгской (Черкесской) межд. акад. наук. "--- 2002. "--- {\sl  6}, \No~1. "--- С.~35--47.

\bibitem{Pshu1} {\it  Псху~А.\,В.} Краевые задачи для дифференциальных уравнений с частными производными дробного и континуального порядка. "--- Нальчик, 2005.

\bibitem{Pul1} {\it  Пулькин~С.\,П.}  Некоторые краевые задачи для уравнения $u_{xx}\pm u_{yy}+\dfrac{p}{x}u_x$//  Уч. зап. Куйбышев.
пед. ин-та. "--- 1958. "--- {\sl 21}. "--- С.~3--54.

\bibitem{Pul2} {\it  Пулькин~С.\,П.} Избранные труды. "--- Самара: Универс групп, 2007.

\bibitem{Pulkina}  {\it  Пулькина~Л.\,С.} Об одной неклассической задаче для вырождающегося гиперболического уравнения//
Изв. вузов. Сер. Мат. "--- 1991. "--- \No~11. "--- С.~48--51.

\bibitem{Ram1}{\it  Рамм~А.\,Г.} Многомерные обратные задачи теории рассеяния. "--- М.: Мир, 1994.

\bibitem{Rv3}{\it  Рвачёв~В.\,А.} Финитные решения функционально-дифференциальных
уравнений и их применения// Усп. мат. наук. "--- 1990. "--- {\sl  45},
\No~1. "--- С.~77--103.

\bibitem{Rv1}{\it  Рвачёв~В.\,Л., Рвачёв~В.\,А.} Теория приближений и атомарные функции. "--- М.: Знание, 1978.

\bibitem{Rv2}{\it  Рвачёв~В.\,Л., Рвачёв~В.\,А.} Неклассические методы теории приближений в краевых задачах. "--- Киев: Наукова думка, 1979.

\bibitem{Rep} {\it  Репин~О.\,А.} Краевые задачи со смещением для уравнений
гиперболического и смешанного типов. "--- Самара, 1992.

\bibitem{RiNa} {\it  Рисс~Ф., Сёкефальви-Надь~Б.} Лекции по функциональному анализу. "--- М.: Мир, 1979.

\bibitem{72} {\it  Ройтберг~Я.\,А., Шефтель~З.\,Г.} Об общих эллиптических задачах с сильным вырождением// Докл. АН СССР. "--- 1980. "--- {\sl  254}, \No~6. "--- С.~1336--1341.

\bibitem{Ros3} {\it Россовский~Л.\,Е.} Эллиптические функционально-дифференциальные уравнения со сжатием и растяжением аргументов неизвестной функции//
Соврем. мат. Фундам. направл. "--- 2014. "--- {\sl  54}. "--- С.~3--138.

\bibitem{Ross1} {\it Россовский~Л.\,Е., Скубачевский~Л.\,Е.} Разрешимость и регулярность решений некоторых классов эллиптических функционально-дифференциальных уравнений//
Итоги науки и техн. Сер. Соврем. мат. и её прил. "--- 1999. "---
{\sl  66}. "--- С.~114--192.

\bibitem{Rut2}{\it  Руткаускас~С.}  Задачи Дирихле с асимптотическими условиями для вырождающейся в точке эллиптической системы.~I//
Дифф. уравн. "--- 2002. "--- {\sl  38}, \No~3. "--- С.~385--392.

\bibitem{Rut3}{\it  Руткаускас~С.} Задачи Дирихле с асимптотическими условиями для вырождающейся в точке эллиптической системы.~II//
Дифф. уравн. "--- 2002. "--- {\sl  38}, \No~5. "--- С.~681--686.

\bibitem{Rut1}{\it  Руткаускас~С.}    О задаче типа Дирихле для эллиптических систем с вырождением на прямой//
Мат. заметки. "--- 2016. "--- {\sl  100}, \No~2. "--- С.~270--278.

\bibitem{Ryko1}{\it  Рыко~В.\,С.} Композиционная структура рядов Фурье, дискретных преобразований
Фурье и Меллина и вычисление их сумм// Деп. в ВИНИТИ АН
СССР. "--- Минск, 1987. "--- \No~2826--В~87.

\bibitem{Ryko2}{\it  Рыко~В.\,С.} Метод суммирования и улучшения сходимости функциональных рядов//
Деп. в ВИНИТИ АН СССР. "--- Вологда, 1988. "--- \No~3542--В~88.

\bibitem{SaIl}  {\it  Сабитов~К.\,Б., Ильясов~Р.\,Р.} Решение задачи Трикоми для уравнения смешанного типа с сингулярным коэффициентом спектральным методом//
 Изв. вузов. Сер. Мат. "--- 2004. "--- \No~2. "--- С.~64--71.

\bibitem{SKM} {\it  Самко~С.\,Г., Килбас~А.\,А., Маричев~О.\,И.} Интегралы и производные дробного порядка и некоторые их приложения. "--- Минск: Наука и техника, 1987.

\bibitem{Sah1} {\it  Сахнович~Л.\,А.} Спектральный анализ вольтерровских операторов и обратные задачи//  Докл. АН СССР. "--- 1957. "--- {\sl  115}, \No~4. "--- С.~666--669.

\bibitem{Sah2} {\it  Сахнович~Л.\,А.}   Обратная задача для дифференциальных операторов порядка $n>2$ с аналитическими коэффициентами//
Мат. сб. "--- 1958. "--- {\sl  46}, \No~1. "--- С.~61--76.

\bibitem{Sah3} {\it  Сахнович~Л.\,А.}   Необходимые условия наличия операторов преобразования для уравнения четвёртого порядка//
Усп. мат. наук. "--- 1961. "--- {\sl  16}, \No~5. "--- С.~199--205.

\bibitem{S75} {\it  Ситник~С.\,М.} О скорости убывания решений стационарного уравнения
Шрёдингера с потенциалом, зависящим от одной переменной// В сб.: <<Краевые задачи для неклассических уравнений
математической физики>>. "--- Новосибирск, 1985. "---   С.~139--147.

\bibitem{S71} {\it  Ситник~С.\,М.}  О скорости убывания решений некоторых эллиптических и  ультраэллиптических уравнений//
Деп. в ВИНИТИ. "---  Воронеж: ВГУ, 1986. "--- 13.11.1986, \No~7771--В86.

\bibitem{S72} {\it  Ситник~С.\,М.}  Об унитарных операторах преобразования//
Деп. в ВИНИТИ. "--- Воронеж: ВГУ, 1986. "---  13.11.1986, \No~7770--В86.

\bibitem{S70} {\it  Ситник~С.\,М.}  Операторы преобразования для дифференциального выражения
Бесселя// Деп. в ВИНИТИ. "--- Воронеж: ВГУ, 1986. "--- 23.01.1987, \No~535--В87.

\bibitem{S2} {\it  Ситник~С.\,М.}  Об одной паре операторов преобразования//
В сб.: <<Краевые задачи для неклассических уравнений математической физики>>. "--- Новосибирск, 1987. "--- С.~168--173.

\bibitem{S73} {\it  Ситник~С.\,М.}  Метод операторов преобразования для стационарного уравнения
Шрёдингера// Дисс. к.ф.-м.н. "--- Воронеж, 1987.

\bibitem{S3}     {\it  Ситник~С.\,М.}  О скорости убывания решений некоторых эллиптических  и  ультраэллиптических уравнений//
Дифф. уравн. "--- 1988. "--- {\sl  24}, \No~3. "--- С.~538--539.

\bibitem{S4} {\it  Ситник~С.\,М.}  Операторы преобразования для сингулярных дифференциальных уравнений
с оператором Бесселя// В сб.: <<Краевые задачи для
неклассических уравнений математической физики>>. "--- Новосибирск, 1989. "---
С.~179--185.

\bibitem{S66} {\it  Ситник~С.\,М.}  Унитарность и  ограниченность операторов Бушмана---Эрдейи
нулевого порядка гладкости// Препринт Ин-та автоматики и
процессов управл. ДВО РАН. "---  Владивосток, 1990.

\bibitem{S6} {\it  Ситник~С.\,М.} Факторизация и оценки норм в весовых  лебеговых  пространствах операторов Бушмана---Эрдейи//
Докл. АН СССР. "--- 1991. "--- {\sl  320}, \No~6. "--- С.~1326--1330.

\bibitem{S63} {\it  Ситник~С.\,М.}  Оператор преобразования  и  представление Йоста для уравнения с
сингулярным потенциалом// Препринт Ин-та автоматики и процессов управл. ДВО РАН. "---  Владивосток,  1993.

\bibitem{S61}  {\it  Ситник~С.\,М.}  Неравенства для полных эллиптических интегралов Лежандра//
Препринт Ин-та автоматики и процессов управл. ДВО РАН. "---
Владивосток, 1994.

\bibitem{S9} {\it  Ситник~С.\,М.}  Неравенства для функций Бесселя//  Докл. РАН. "--- 1995. "--- {\sl  340}, \No~1. "--- С.~29--32.

\bibitem{S140} {\it  Ситник~С.\,М.}   Метод получения  последовательных уточнений неравенства
Коши---Буняковского и его применения к оценкам специальных
функций// Тез. докладов. Воронеж. весен.
мат. школы <<Современные методы в теории краевых задач.
Понтрягинские чтения-VII>>. "---  Воронеж: ВГУ, 1996. "---
С.~164.

\bibitem{S140p} {\it  Ситник~С.\,М.}  Формула Тэйлора для операторов типа Бесселя//
Тез. докладов. Воронеж. весен.
мат. школы <<Современные методы в теории краевых задач.
Понтрягинские чтения-VII>>. "---  Воронеж: ВГУ, 1996. "---  С.~102.

\bibitem{S135} {\it  Ситник~С.\,М.}   О некоторых обобщениях дробного интегро-дифферен\-цирования// Материалы межд. Российско-Узбекского симп.
<<Уравнения смешанного типа и родственные проблемы анализа и
информатики>>. "--- Нальчик, 2003. "--- C.~86--88.

\bibitem{S133} {\it  Ситник~С.\,М.}  Дробное интегродифференцирование для дифференциального оператора
Бесселя// Материалы межд. Российско-Казахского симп.
<<Уравнения смешанного типа и родственные проблемы анализа и
информатики>>. "--- Нальчик, 2004. "---
С.~163--167.

\bibitem{S53}  {\it  Ситник~С.\,М.}  Обобщения неравенств Коши---Буняковского методом средних значений
и их приложения// Чернозёмный альманах науч. иссл. Сер.
Фундам. мат. "--- 2005. "---  \No~1 (1). "--- C.~3--42.

\bibitem{S127} {\it  Ситник~С.\,М.}  Об обобщении формулы Хилле---Тамаркина для резольвенты на случай операторов дробного интегрирования Бесселя//
III Межд. конф. <<Нелокальные краевые задачи и
родственные проблемы математической биологии, информатики и
физики>>. "--- Нальчик, 2006. "--- С.~269--270.

\bibitem{S123} {\it  Ситник~С.\,М.}   Операторы дробного интегро-дифференцирования для дифференциального
оператора Бесселя// Тр. Четвёртой всерос. науч.
конф. с межд. участием <<Математическое
моделирование и краевые задачи>>, Ч.~3. "--- Самара, 2007.
"--- С.~158--160.

\bibitem{S125} {\it  Ситник~С.\,М.}  Построение операторов преобразования Векуа---Эрдейи---Лаундеса//
Межд. конф., посвящ. 100-летию со дня рождения
акад. И.\,Н.~Векуа, <<Дифференциальные уравнения, теория
функций и приложения>>. Тезисы докладов. "---
Новосибирск, 2007. "--- С.~469--470.

\bibitem{S46} {\it  Ситник~С.\,М.}  Операторы преобразования и их приложения// В сб.: <<Исследования по современному анализу и математическому моделированию>>. "---
Владикавказ: Владикавказ. науч. центр РАН и РСО-А, 2008. "---
C.~226--293.

\bibitem{S14} {\it  Ситник~С.\,М.} Метод факторизации операторов преобразования в теории дифференциальных уравнений//
Вестн. Самар. гос. ун-та. Естественнонауч. сер. "--- 2008. "---
\No~8/1 (67). "--- С.~237--248.


\bibitem{S45} {\it  Ситник~С.\,М.}  Уточнения и обобщения классических неравенств//
В сб.:  <<Итоги науки. Южный федеральный округ>>. Сер.
<<Математический форум>>. Т.~3. <<Исследования по математическому
анализу>>. "---
Владикавказ: Южный мат. ин-т ВНЦ  РАН и РСО Алания, 2009.
"--- С.~221--266.

\bibitem{S18}  {\it  Ситник~С.\,М.} О явных реализациях дробных степеней дифференциального оператора Бесселя и их приложениях к дифференциальным уравнениям//
Докл. Адыгской (Черкесской) межд. акад. наук. "--- 2010. "--- {\sl  12}, \No~2. "--- С.~69--75.

\bibitem{S19} {\it  Ситник~С.\,М.} О представлении в интегральном виде решений одного дифференциального уравнения с особенностями в коэффициентах//
Владикавказ. мат.~ж. "--- 2010. "--- {\sl  12}, \No~4. "--- С.~73--78.

\bibitem{S43}  {\it  Ситник~С.\,М.}   Оператор преобразования специального вида для дифференциального оператора с сингулярным в нуле потенциалом//
 В сб.: <<Неклассические уравнения математической физики>>. "--- Новосибирск: Изд-во
ин-та мат. им. С.\,Л.~Соболева СО РАН, 2010. "--- С.~264--278.

\bibitem{S103} {\it  Ситник~С.\,М.}  Различные классы операторов преобразования Бушмана---Эрдейи//
Межд. конф. <<Дифференциальные уравнения и смежные
вопросы>>. Сборник тезисов. "---  Москва: МГУ, 2011. "---  С.~344--345.

\bibitem{S92} {\it  Ситник~С.\,М.}  Новые краевые задачи с $K$-следом для решений с существенными особенностями в работах В.\,В.~Катрахова//
Межд. конф. <<Обратные и некорректные задачи математической физики>>. Тезисы докладов. "---
Новосибирск: Сиб. науч. изд-во, 2012. "--- С.~439.

\bibitem{S95} {\it  Ситник~С.\,М.}  Работы В.\,В.~Катрахова по теории операторов преобразования//
Материалы второго Российско-Узбекского симп. <<Уравнения
смешанного типа и родственные проблемы анализа и информатики>>.
"---  Нальчик: КБР, Институт прикладной математики и
автоматизации КБНЦ РАН, 2012. "--- С.~241--243.

\bibitem{S400} {\it  Ситник~С.\,М.}   Обзор основных свойств операторов преобразования Бушмана---Эрдейи//  Челябинск. физ.-мат.~ж. "--- 2016. "--- {\sl  1}, \No~4. "--- С.~63--93.

\bibitem{SitDis} {\it Ситник~С.\,М.}  Применение операторов преобразования Бушмана---Эрдейи и их обобщений в теории дифференциальных уравнений с особенностями в коэффициентах//
Дисс. д.ф.-м.н. "--- Воронеж, 2016.

\bibitem{S62} {\it  Ситник~С.\,М., Карп~Д.\,Б.}  Формулы композиций для интегральных преобразований
с  функциями Бесселя в ядрах// Препринт Ин-та автоматики и
процессов управл. ДВО РАН. "---  Владивосток, 1993.

\bibitem{S60} {\it  Ситник~С.\,М., Карп~Д.\,Б.}  Дробное преобразование Ханкеля и его приложения
в математической физике// Препринт Ин-та автоматики и процессов управл. ДВО РАН. "---  Владивосток,  1994.

\bibitem{S65} {\it  Ситник~С.\,М., Ляховецкий~Г.\,В.}  Формулы композиций для операторов Бушмана---Эрдейи//
Препринт Ин-та автоматики и процессов управл. ДВО РАН. "---
Владивосток, 1991.

\bibitem{S59} {\it  Ситник~С.\,М., Ляховецкий~Г.\,В.}  Операторы преобразования Векуа---Эрдейи---Лаундеса//
Препринт Ин-та автоматики и процессов управл. ДВО РАН. "---  Владивосток,  1994.

\bibitem{SitShishSemi} {\it  Ситник~С.\,М., Шишкина~Э.\,Л.}   Об одном тождестве для итерированного весового сферического среднего и его приложениях//
Сиб. электрон. мат. изв. "--- 2016. "--- {\sl  13}. "--- С.~849--860.

\bibitem{SSfiz} {\it Ситник~С.\,М., Шишкина~Э.\,Л.} Метод операторов преобразования для дифференциальных уравнений с операторами Бесселя. "--- М.: Физматлит, 2018.

\bibitem{S700} {\it Ситник~С.\,М., Шишкина~Э.\,Л.}  О дробных степенях оператора Бесселя на полуоси//  Сиб. электрон. мат. изв. "--- 2018. "--- {\sl  15}. "--- С.~1--10.

\bibitem{Sku1}  {\it  Скубачевский~А.\,Л.} Неклассические краевые задачи.~I//  Соврем. мат. Фундам. направл. "--- 2007. "--- {\sl  26}. "--- С.~3--132.

\bibitem{Sku4}  {\it  Скубачевский~А.\,Л.} Нелокальные краевые задачи и их приложения к исследованию
многомерных диффузионных процессов и процессов терморегуляции
живых клеток. Учеб. пособие. "--- М.: РУДН, 2008.

\bibitem{Sku2}  {\it  Скубачевский~А.\,Л.} Неклассические краевые задачи.~II// Соврем. мат. Фундам. направл. "--- 2009. "--- {\sl 33}. "--- С.~3--179.

\bibitem{Sku3}  {\it  Скубачевский~А.\,Л., Гуревич~П.\,Л.} Применение методов нелинейного функционального анализа к нелокальным проблемам процессов распределения тепла. Учеб. пособие.
"--- М.: РУДН, 2008.

\bibitem{76} {\it  Смирнов~М.\,М.} Вырождающиеся эллиптические и гиперболические уравнения. "--- М.: Наука, 1966.

\bibitem{Smi} {\it  Смирнов~М.\,М.} Уравнения смешанного типа. "--- М.: Наука, 1970.

\bibitem{77} {\it  Соболев~С.\,Л.} Введение в теорию кубатурных формул. "--- М.: Наука, 1974.

\bibitem{Sob}{\it  Соболев~С.\,Л.} Уравнения математической физики. "--- М.: Наука, 1992.

\bibitem{Sol} {\it  Солдатов~А.\,П.} Одномерные сингулярные операторы и краевые задачи теории функций. "--- М.: Высшая школа, 1991.

\bibitem{Son1} {\it  Сонин~Н.\,Я.} Исследования о цилиндрических  функциях и специальных полиномах. "--- М.: Гостехтеоретиздат, 1954.

\bibitem{Soh1} {\it  Сохин~А.\,С.} Об одном классе операторов преобразования//  Тр. физ.-тех. ин-та низк. темп. АН УССР. "--- 1969. "--- {\sl 1}. "--- С.~117--125.

\bibitem{Soh2} {\it  Сохин~А.\,С.} Обратные задачи рассеяния для уравнений с особенностью//  Тр. физ.-тех. ин-та низк. темп. АН УССР. "--- 1971. "--- {\sl 2}. "--- С.~182--233.

\bibitem{Soh3} {\it  Сохин~А.\,С.} Обратные задачи рассеяния для уравнений с особенностями специального вида//
Теор. функций, функц. анализ и их прилож. "--- 1973.
"--- {\sl 17}. "--- С.~36--64.

\bibitem{Soh4} {\it  Сохин~А.\,С.} О преобразовании операторов  для уравнений с особенностью специального вида//  Вестн. Харьков. ун-та. "--- 1974. "--- {\sl 113}. "--- C.~36--42.

\bibitem{Sta1} {\it Сташевская~В.\,В.} Метод операторов преобразования//  Докл. АН СССР. "--- 1953. "--- {\sl  113}, \No~3. "--- С.~409--412.

\bibitem{Sta2} {\it Сташевская~В.\,В.} Об обратной задаче спектрального анализа для дифференциального оператора с особенностью в нуле//
Уч. зап. Харьков. мат. об-ва. "--- 1957. "--- \No~5. "--- С.~49--86.

\bibitem{78} {\it Стейн~И., Вейс~Г.} Введение в гармонический анализ на евклидовых пространствах. "--- М.: Мир, 1974.

\bibitem{Ter1} {\it  Терсенов~С.\,А.} Введение в теорию уравнений, вырождающихся на границе. "--- Новосибирск: НГУ, 1973.

\bibitem{Tit1} {\it  Титчмарш~Е.} Введение в теорию интегралов Фурье. "--- М.---Л.: ГИТТЛ, 1948.

\bibitem{79} {\it  Тихонов~А.\,Н., Самарский~А.\,А.} Уравнения математической физики. "--- М.: Наука, 1972.

\bibitem{Trib1} {\it  Трибель~Х.} Теория интерполяции. Функциональные пространства. Дифференциальные операторы. "--- М.: Мир, 1980.

\bibitem{Tricomi1} {\it  Трикоми~Ф.} Лекции по уравнениям в частных производных. "--- М.: Иностр. лит., 1957.

\bibitem{Tyr} {\it  Тыртышников~Е.\,Е.} Матричный анализ и линейная алгебра. "--- М.: Физматлит, 2007.

\bibitem{Wit1} {\it  Уиттекер~Э., Ватсон~Дж.} Курс современного анализа. Ч.~2.  Трансцендентные функции. "--- М.: ГИФМЛ, 1963.

\bibitem{Usp} {\it  Успенский~С.\,В.} О теоремах вложения для весовых классов//  Тр. МИАН. "--- 1961. "--- {\sl 60}. "--- С.~282--303.

\bibitem{Fage3} {\it  Фаге~М.\,К.} Построение операторов преобразования и решение одной проблемы моментов для обыкновенных линейных дифференциальных уравнений произвольного порядка//
Усп. мат. наук. "--- 1957. "--- {\sl  12}, \No~1. "--- С.~240--245.

\bibitem{Fage1} {\it  Фаге~М.\,К.} Операторно-аналитические функции одной независимой переменной//
Докл. АН СССР. "--- 1957. "--- {\sl 112}, \No~6. "--- С.~1008--1011.

\bibitem{Fage2} {\it  Фаге~М.\,К.} Интегральные представления операторно-аналитических функций одной независимой переменной//
Докл. АН СССР. "--- 1957. "--- {\sl  115}, \No~5. "--- С.~874--877.

\bibitem{Fage4} {\it  Фаге~М.\,К.} Операторно-аналитические функции одной независимой переменной//  Тр. Моск. мат. об-ва. "--- 1958. "--- {\sl  7}. "--- С.~227--268.

\bibitem{Fage5} {\it  Фаге~М.\,К.} Интегральные представления операторно-аналитических функций одной независимой переменной//  Тр. Моск. мат. об-ва. "--- 1958. "--- {\sl  8}. "--- С.~3--48.

\bibitem{Fage6} {\it  Фаге М.\,К.} Операторно-аналітичні функції однієї незалежної змінної. "--- Львов: Изд-во Львовск. ун-та, 1959.

\bibitem{FN} {\it  Фаге~Д.\,К., Нагнибида~Н.\,И.} Проблема эквивалентности обыкновенных
дифференциальных операторов. "--- Новосибирск: Наука, 1977.

\bibitem{Fad1} {\it  Фаддеев~Л.\,Д.} Обратная задача квантовой теории рассеяния.~1//  Усп. мат. наук. "--- 1959. "--- {\sl  14}, \No~4. "--- С.~57--119.

\bibitem{81} {\it  Фаддеев~Л.\,Д.} Разложение по собственным функциям оператора Лапласа на фундаментальной области дискретной группы на плоскости Лобачевского//
Тр. Моск. мат. об-ва. "--- 1967. "--- {\sl  17}. "--- С.~323--350.

\bibitem{Fad2} {\it  Фаддеев~Л.\,Д.}  Обратная задача квантовой теории рассеяния.~2// Итоги науки и техн. Соврем. пробл. мат. "--- 1974. "--- {\sl 3}. "--- С.~93--180.

\bibitem{FeIv} {\it  Фёдоров~В.\,Е.}   Нелокальная на полуоси задача для вырожденных эволюционных уравнений//
Мат. заметки СВФУ. "--- 2015. "--- {\sl  22}, \No~1. "---
С.~35--43.

\bibitem{FGP} {\it  Фёдоров~В.\,Е., Гордиевских~Д.\,М., Плеханова~М.\,В.} Уравнения в банаховых пространствах с вырожденным оператором под знаком дробной производной//
Дифф. уравн. "--- 2015. "--- {\sl  51}, \No~10. "--- С.~1367--1375.

\bibitem{Fet} {\it  Фетисов~В.\,Г.} Операторы и уравнения в локально ограниченных пространствах//
В сб.: <<Исследования по функциональному анализу и его
приложениям>>. "--- М.: Наука, 2005. "--- С.~249--292.

\bibitem{Fish}   {\it  Фишман~М.\,К.} Об эквивалентности некоторых линейных операторов в аналитическом пространстве//  Мат. сб. "--- 1965. "--- {\sl  68}, \No~1. "--- С.~63--74.

\bibitem{Fri} {\it  Фридрихс~К.} Возмущение спектра операторов в гильбертовом пространстве. "--- М.: Мир, 1968.

\bibitem{Han} {\it  Ханмамедов~А.\,Х.} Операторы преобразования для возмущённого разностного уравнения Хилла и их одно приложение//
Сиб. мат.~ж. "--- 2003. "--- {\sl  44}, \No~4. "--- С.~926--937.

\bibitem{HLP} {\it  Харди~Г.\,Г., Литтлвуд~Дж.~Е., Полиа~Г.} Неравенства. "--- М.: Иностр. лит., 1948.

\bibitem{Hach1} {\it  Хачатрян~И.\,Г.} Об операторах преобразования для дифференциальных уравнений высших порядков//
Изв. АН Армен. ССР. Сер. Мат. "--- 1978. "--- {\sl  13}, \No~3. "--- С.~215--236.

\bibitem{Hach2} {\it  Хачатрян~И.\,Г.} Об операторах преобразования для дифференциальных уравнений высших порядков, сохраняющих асимптотику решений//
Изв. АН Армен. ССР. Сер. Мат. "--- 1979. "--- {\sl  14}, \No~6. "---  С.~424--445.

\bibitem{Kan} {\it Хе Кан Чер.} Сингулярные краевые задачи для уравнений математической физики с операторами Бесселя// Дисс. д.ф.-м.н. "--- Хабаровск, 1991.

\bibitem{Hale} {\it  Хейл~Дж.} Теория функционально-дифференциальных уравнений. "---
 М.: Мир, 1984.

\bibitem{82} {\it  Хелгасон~С.}  Преобразование Радона. "--- М.: Мир, 1983.

\bibitem{Hel1} {\it  Хелгасон~С.} Группы и геометрический анализ. "--- М.: Мир, 1987.

\bibitem{83} {\it  Хермандер~Л.} Линейные дифференциальные операторы с частными производными. "--- М.: Мир, 1965.

\bibitem{Hor} {\it  Хорн~Р.} Матричный анализ. "--- М.: Мир, 1989.

\bibitem{Hrom1} {\it  Хромов~А.\,П.} Конечномерные возмущения вольтерровых операторов//  Соврем. мат. Фундам. направл. "--- 2004. "--- \No~10. "--- С.~3--163.

\bibitem{Cher} {\it  Чернятин~В.\,А.} Обоснование метода Фурье в смешанной задаче для уравнений в частных производных. "--- М.: МГУ, 1991.

\bibitem{ShSa} {\it  Шадан~К., Сабатье~П.} Обратные задачи в квантовой теории рассеяния. "--- М.: Мир, 1980.

\bibitem{Sha1} {\it  Шацкий~В.\,П.} О гиперболических системах с двумя плоскостями особенностей//    Докл. АН СССР. "---  1978. "--- {\sl 242}, \No~1. "--- С.~56--59.

\bibitem{Sha2} {\it  Шацкий~В.\,П.} Об одной краевой задаче для сингулярных симметрических систем нечётного порядка// Докл. АН СССР. "--- 1979. "--- {\sl 243}, \No~4. "--- С.~806--809.

\bibitem{Sha3} {\it  Шацкий~В.\,П.} О некоторых вырождающихся системах первого порядка в областях с характеристической частью границы//
Докл. АН СССР. "--- 1982. "--- {\sl 262}, \No~6. "--- С.~1332--1335.

\bibitem{ShiR1} {\it  Шишкина~Э.\,Л.} Обобщённая весовая функция $r^\gamma$//  Вестн. ВГУ. Cер. Физ. Мат. "--- 2006. "--- \No~1. "---  С.~215--221.

\bibitem{ShiR3} {\it  Шишкина~Э.\,Л.}  Равенство для интерированных весовых сферических средних, порожденных обобщённым сдвигом//
Материалы науч. конф. <<Герценовские чтения-2013>>. "--- СПб.:
Изд-во РГПУ им. А.\,И.~Герцена, 2013. "--- {\sl 66}. "---
C.~143--145.

\bibitem{ShiR2} {\it  Шишкина~Э.\,Л.}  О свойствах одного усредняющего ядра в весовом классе Лебега//
Науч. ведом. Белгород. гос. ун-та. Сер. Мат. Физ. "--- 2016. "---
{\sl  42}, \No~6. "--- С.~12--19.

\bibitem{ShiR4}{\it Шишкина~Э.\,Л.}   Весовые обобщённые функции, отвечающие квадратичной форме с комплексными коэффициентами//
Челябинск. физ.-мат.~ж. "--- 2017. "--- {\sl  2}, \No~1. "---
С.~88--98.

\bibitem{85} {\it  Шишмарев~И.\,А.} Введение в теорию эллиптических уравнений. "--- М.: Изд-во МГУ, 1979.

\bibitem{Shos} {\it  Шостак~Р.\,Я.} Алексей Васильевич Летников// В сб.: <<Историко-математические исследования.
Труды семинара МГУ по истории математики>>. "--- М.-Л.: ГИТТЛ,
1952. "--- {\sl 5}. "--- С.~167--238.

\bibitem{Els} {\it  Эльсгольц~Л.\,Э.} Введение в теорию дифференциальных уравнений с отклоняющимся аргументом. "--- М.: Наука, 1964.

\bibitem{ElNor} {\it  Эльсгольц~Л.\,Э., Норкин~С.\,Б.} Введение в теорию дифференциальных уравнений с отклоняющимся аргументом. "--- М.: Наука, 1971.


\bibitem{Erg1} {\it  Эргашев~Т.\,Г.} Четвертый потенциал двойного слоя для обобщённого двуосесимметрического уравнения Гельмгольца//
Вестн. Томск. гос. ун-та. Мат. и мех. "---  2017. "--- {\sl 50}.
"--- C.~45--56.

\bibitem{Yurko} {\it  Юрко~В.\,А.} Введение в теорию обратных спектральных задач. "--- М.: Наука, 2007.

\bibitem{87} {\it  Яковлев~Г.\,Н.} Неограниченные решения вырождающихся эллиптических уравнений//  Тр. МИАН. "--- 1978. "--- {\sl  117}. "--- С.~312--320.

\bibitem{Yanu} {\it  Янушаускас~А.\,И.} Аналитическая теория эллиптических уравнений. "--- Новосибирск: Наука, 1979.

\bibitem{Yarem} {\it  Яремко~О.\,Э.} Метод операторов преобразования в задачах математического
моделирования. "--- Пенза: Пенз. гос. ун-т, 2012; Lambert Academic Publishing, 2012.

\bibitem{Yar1} {\it  Ярославцева~В.\,Я.} Об одном классе операторов преобразования и их приложении к дифференциальным уравнениям//
Докл. АН СССР. "--- 1976. "--- {\sl  227}, \No~4. "--- С.~816--819.

\bibitem{Yar2} {\it  Ярославцева~В.\,Я.}  Неоднородная граничная задача в полупространстве для одного класса сингулярных уравнений//
Деп. ред. Дифф. уравн. "--- 1989.






\bibitem{AKK} {\it  Ali~I., Kiryakova~V., Kalla~S.\,L.} Solutions of fractional multi-order integral and differential equations using a Poisson-type transform//
J.~Math. Anal. Appl. "--- 2002. "--- {\sl  269}, \No~1. "---
С.~172--199.

\bibitem{AlKi1} {\it Almalki~F., Kisil~V.} Geometric dynamics of a harmonic oscillator,
non-admissible mother wavelets and squeezed states//
arXiv:1805.01399v1. "--- 2018.

\bibitem{AAR} {\it  Andrews~G.\,E.,  Askey~R., Roy~R.} Special functions. "--- Cambridge: Cambridge University Press, 1999.

\bibitem{AK} {\it  Antimirov~M.\,Ya., Kolyshkin~A.\,A., Vaillancourt~R.} Applied integral transforms. "--- Providence: Am. Math. Soc., 1993.

\bibitem{AKdF} {\it  Appell~P., Kampe de Feriet~J.} Fonctions hypergeometriques et hyperspheriques; polynomes d'Hermite. "--- Paris: Gauthier-Villars, 1926.

\bibitem{AppHyp} {\it Connett~W.\,C., Gebuhrer~M.-O., Schwartz~F.\,L. $($ред.$)$} Applications of hypergroups and related measure algebras. "--- Providence: Am. Math. Soc., 1995.

\bibitem{Bac} {\it  Baccar~C., Hamadi~N.\,B., Achdi~L.\,T.}  Inversion formulas for Riemann---Liouville transform and its dual associated with singular partial differential operators//
Int. J. Math. Math. Sci. "--- 2006. "--- C.~1--26.

\bibitem{Bai}  {\it  Bailey~W.\,N.} Generalized hypergeometric series. "--- New York---London: Stechert-Hafner Service Agency, 1964.

\bibitem{Bajlekova1} {\it  Bajlekova~E.\,G.} Subordination principle for fractional evolution equations//  Fract. Calc. Appl. Anal. "---  2000. "---  {\sl 3}, \No~3. "---   C.~213--230.

\bibitem{Bajlekova0}  {\it  Bajlekova~E.\,G.} Fractional evolution equations in Banach spaces// Thesis. "---  Technische Universiteit Eindhoven, 2001.

\bibitem{Bers1} {\it  Bers~L.} On a class of differential equations in mechanics of continua//  Quart. Appl. Math. "--- 1943. "--- {\sl  5}, \No~1. "--- С.~168--188.

\bibitem{Bers3} {\it  Bers~L.}   A remark on an applications of pseudo-analytic functions//  Amer. J.~Math. "--- 1956. "--- {\sl  78}, \No~3. "--- С.~486--496.

\bibitem{Bers2} {\it  Bers~L., Gelbart~A.}    On a class of functions defined by partial differential equations//  Trans. Am. Math. Soc. "--- 1944. "--- {\sl  56}. "--- С.~67--93.

\bibitem{BuKe1} {\it Bourgain~J., Kenig~C.\,E.} On localization in the continuous Anderson---Bernoulli model in higher dimension//
Invent. Math. "---  2005. "--- {\sl 161}, \No~2. "--- C.~389--426.

\bibitem{Bozh}  {\it  Bozhinov~N.} Convolution representations of commutants and multipliers. "--- Sofia: Publishing House Bulgarian Acad. Sci., 1988.

\bibitem{Bra}  {\it Braaksma B.\,L.\,J.} Asymptotic expansions and analytic continuation for a class of Barnes-integrals. "--- Groningen: Noordhoff, 1963.

\bibitem{BD2}  {\it Bragg~Z.\,R., Dettman~J.\,W.} An operator calculus for related partial differential equations// J.~Math. Anal. Appl. "--- 1968. "--- {\sl 22}, \No~2. "--- C.~261--271.

\bibitem{BD1}  {\it Bragg~Z.\,R., Dettman~J.\,W.} Related problems in partial differential equations// Bull. Am. Math. Soc. "--- 1968. "--- {\sl 74}, \No~2. "--- C.~375--378.

\bibitem{BD3}  {\it Bragg~Z.\,R., Dettman~J.\,W.} A class of related Dirichlet and initial value problems// Proc. Am. Math. Soc. "--- 1969. "--- {\sl 21}, \No~1. "--- C.~50--56.

\bibitem{Bresters2}  {\it Bresters~D.\,W.} On a generalized Euler---Poisson---Darboux equation// SIAM J. Math. Anal. "--- 1978. "--- {\sl 9}, \No~5. "---  C.~924--934.

\bibitem{Bul2} {\it  Bullen~P.\,S.} Handbook of means and their inequalities. "--- Dordrecht---London: Kluwer Academic Publishers, 2003.

\bibitem{BMV}  {\it  Bullen~P.\,S., Mitrinovi\'c~D.\,S., Vasi\'c~P.\,M.} Means and their inequalities. "--- Dordrecht: D.~Reidel, 1988.

\bibitem{Bus2} {\it  Buschman~R.\,G.}  An inversion integral for a Legendre transformation//  Amer. Math. Monthly. "--- 1962. "--- {\sl  69}, \No~4. "--- С.~288--289.

\bibitem{Bus1} {\it  Buschman~R.\,G.} An inversion integral for a general Legendre transformation//  SIAM Review. "--- 1963. "--- {\sl  5}, \No~3. "--- С.~232--233.

\bibitem{CKT1}  {\it  Campos~H., Kravchenko~V.\,V., Torba~S.\,M.}
 Transmutations, $L$-bases and complete families of solutions of the stationary Schr\"odinger equation in the plane//
 J.~Math. Anal. Appl. "--- 2012. "--- {\sl  389}, \No~2. "--- С.~1222--1238.

\bibitem{Car1} {\it  Carroll~R.\,W.} Transmutation and operator differential
equations. "--- Amsterdam---New York---Oxford: North Holland,
1979.

\bibitem{Car2} {\it  Carroll~R.\,W.} Transmutation, scattering theory and special functions. "--- Amsterdam---New York---Oxford: North Holland, 1982.

\bibitem{Car3} {\it  Carroll~R.\,W.} Transmutation theory and
applications. "--- Amsterdam---New York: North Holland, 1986.

\bibitem{Car4} {\it  Carroll~R.\,W.} Topics in soliton theory. "--- North Holland, 1991.

\bibitem{Car10} {\it  Carroll~R.\,W.} Calculus revisited. "--- Dordrecht---Boston---London: Springer, 2002.

\bibitem{CB} {\it  Carroll~R.\,W., Boumenir~A.} Toward a general theory of
transmutation//  arXiv:~funct-an/9501006. "--- 1995.

\bibitem{CSh} {\it  Carroll~R.\,W., Showalter~R.\,E.} Singular and degenerate Cauchy problems. "--- N.Y.: Academic Press, 1976.

\bibitem{CKT2} {\it  Castillo-P\'erez~С.\,R., Kravchenko~V.\,V.,  Torba~S.\,M.} Spectral parameter power series for perturbed Bessel equations//
Appl. Math. Comput. "--- 2013. "--- {\sl 220}. "--- C.~676--694.

\bibitem{ChCPR} {\it  Chadan~K., Cotton~D., Paivarinta~L., Rundell~W.} An introduction to inverse scattering and inverse spectral problems. "--- SIAM, 1997.

\bibitem{Che1} {\it Chebli~H.} Op\'erateurs de translation g\'en\'eralises et semigroupes de convolution// Springer Lect. Notes. "--- 1974. "--- {\sl 404}. "--- C.~35--59.

\bibitem{Che3} {\it Chebli~H.}  Positivit\'e des op\'erateurs de <<translation g\'en\'eralises>> associ\'e \`a un op\'erateur de Sturm---Liouville et quelques applications
a l'analyse harmonique// Thesis. "--- Strasbourg, 1974.

\bibitem{Che2} {\it Chebli~H.}  Sur un th\`eor\'eme de Paley---Winer associ\'e \`a la d\'ecomposition spectrale d'un op\'erateur de Sturm---Liouville sur $(0,\infty)$//
J.~Funct. Anal. "--- 1974. "--- {\sl 17}. "--- С.~447--461.

\bibitem{Che4} {\it Chebli~H.}   Th\`eor\'eme de Paley---Winer associ\'e \`a un op\'erateur diff\'erentiel singulier sur $(0,\infty)$//
J.~Math. Pures Appl. "--- 1979. "--- {\sl  58}. "--- С.~1--19.

\bibitem{CFH} {\it Chebli~H., Fitouhi~A., Hamza~M.\,M.}  Expansion in series of Bessel functions and transmutations for perturbed Bessel operators//
J.~Math. Anal. Appl. "--- 1994. "--- {\sl  181}, \No~3. "--- С.~789--802.

\bibitem{Che} {\it  Cheikh~B.} Relations between harmonic analysis associated with two differential operators of different orders//
J.~Comput. Appl. Math. "--- 2003. "--- {\sl  153}, \No~1. "--- С.~61--71.

\bibitem{Col1} {\it  Colton~D.} Solution of boundary value problems by the method of integral operators. "--- London: Pitman Press, 1976.

\bibitem{Col2} {\it  Colton~D.} Analytic theory of partial differential equations. "--- London: Pitman Press, 1980.

\bibitem{Cop2}  {\it  Copson~E.\,T.} On a singular boundary value problem
for an equation of hyperbolic type//  Arch. Ration. Mech. Anal. "--- 1957. "--- {\sl 1}. "--- С.~349--356.

\bibitem{Cop1} {\it  Copson~E.\,T.} On the Riemann---Green function// Arch. Ration. Mech. Anal. "--- 1957/58. "--- {\sl  1}. "--- С.~324--348.

\bibitem{Cop4} {\it  Copson~E.\,T.} Partial differential equations. "--- London---New York---Melbourne: Cambridge University Press, 1975.

\bibitem{Cop3} {\it  Copson~E.\,T., Erd\'elyi~A.} On a partial differential equation
with two singular lines// Arch. Ration. Mech. Anal. "---
1958. "--- {\sl  2}, \No~1. "--- С.~76--86.

\bibitem{Darboux} {\it  Darboux~G.}  Le\c cons sur la th\'eorie g\'en\'erale des surfaces et les applications g\'eom\'etriques du calcul
        infinit\'esimal. Vol.~2. "--- Paris: Gauthier-Villars, 1915.

\bibitem{DKW1} {\it Davey~B.,  Kenig~C., Wang~J.-N.} The Landis conjecture for variable coeffcient second-order elliptic PDEs//
Trans. Am. Math. Soc. "--- 2017. "--- {\sl 369}, \No~11. "--- C.~8209--8237.

\bibitem{Deans} {\it  Deans~S.\,R.} The Radon transform and some of its applications. "--- Mineola---New York: Dover, 1990.

\bibitem{Del1} {\it  Delsarte~J.} Sur certaines transformation fonctionnelles relative aux \'equations lin\'eares aux d\`eriv\'ees partielles du seconde ordre//
C.\,R. Acad. Sci. "--- 1938. "--- {\sl  206}. "--- С.~1780--1782.

\bibitem{Del2} {\it  Delsarte~J.}  Sur une extension de la formule de Taylor//  J. Math. Pures Appl. "--- 1938. "--- {\sl  17}. "--- С.~217--230.

\bibitem{Del3} {\it  Delsarte~J.}  Une extension nouvelle de la th\'eory de fonction presque p\'eriodiques de Bohr//  Acta Math. "--- 1939. "--- {\sl  69}. "--- С.~259--317.

\bibitem{Del4} {\it  Delsarte~J.} Hypergroupes et operateurs de permutation et de transmutation// Colloq. Internat. Centre Nat. Rech. Sci. "--- 1956. "--- {\sl 71}. "--- С.~29--44.

\bibitem{Del7} {\it  Delsarte~J.} Lectures on topics in mean periodic functions
and the two-radius theorem. "--- Bombay: Tata Inst.
Fundam. Research, 1961.

\bibitem{Del5} {\it  Delsarte~J., Lions~J.\,L.}  Transmutations d'op\'erateurs diff\'erentiels dans le domaine complexe//  Comment. Math. Helv. "--- 1957. "--- {\sl 32}. "--- С.~113--128.

\bibitem{Del6}  {\it  Delsarte~J., Lions~J.\,L.} Moyennes g\'en\'eralis\'ees// Comment. Math. Helv. "--- 1959. "--- \No~34. "--- C.~59--69.

\bibitem{Diaz}  {\it  Diaz~J.\,B.,  Weinberger~H.\,F.} A solution of the singular initial value problem for the Euler---Poisson equation//
Proc. Am. Math. Soc. "--- 1953. "--- {\sl 4}. "--- C.~703--715.

\bibitem{Dim} {\it  Dimovski~I.} Convolutional calculus. "--- Dordrecht: Kluwer Acad. Publ., 1990.

\bibitem{DHS} {\it  Dimovski~I., Hristov~V., Sifi~M.}  Commutants of the Dunkl operators in $C(\R)$// Fract. Calc. Appl. Anal. "--- 2006. "--- {\sl  9}, \No~3. "--- С.~195--213.

\bibitem{DK2}  {\it  Dimovski~I., Kiryakova~V.\,S.} Transmutations, convolutions and fractional powers of Bessel-type operators via Meijer $G$-functions//
Proc. conf. <<Complex Anal. and Appl.>>, 1983, Sofia. "--- Varna, 1985. "--- C.~45--66.

\bibitem{DK1} {\it  Dimovski~I., Kiryakova~V.\,S.}  Generalized Poisson transmutations and corresponding representations of hyper-Bessel functions//
C.\,R. Acad. Bulgare Sci. "--- 1986. "--- {\sl  39}, \No~10. "--- С.~29--32.

\bibitem{Djr2} {\it  Djrbashian~M.\,M.} Harmonic analysis and boundary value problems in the complex domain. "--- Basel---Boston---Berlin: Birkh\"auser, 1993.

\bibitem{DKN} {\it  Drabek~P., Kufner~A., Nicolosi~F.} Quasilinear elliptic equations
with degenerations and singularities. "--- Berlin---New York: Walter de Gruyter, 1997.

\bibitem{Dun1} {\it  Dunkl~Ch.} Differential-difference operators associated to reflection groups// Trans. Am. Math. Soc. "--- 1989. "--- {\sl  311}. "--- С.~167--183.

\bibitem{Dun2} {\it  Dunkl~Ch.}   Intertwining operators associated to the group S3//  Trans. Am. Math. Soc. "--- 1995. "--- {\sl  347}. "--- С.~3347--3374.

\bibitem{Dun3} {\it  Dunkl~Ch.}   An Intertwining operator for the group B2//  arXiv:~math.~CA/~0607823. "--- 2006.

\bibitem{Dwo1} {\it  Dwork~B.} Generalized hypergeometric functions. "--- Oxford, 1990.

\bibitem{Dwo2}  {\it  Dwork~B., Gerotto~S., Sullivan~F.\,J.} Introduction to $G$-functions. "--- Princeton, 1994.

\bibitem{Dzr} {\it  Dzrbashian~M.\,M.} Harmonic analysis and boundary value problems in the complex domain. "--- Basel---Boston---Berlin: Birkh\"auser, 1993.

\bibitem{EiIvKoch} {\it  Eidelman~S.\,D.,  Ivasyshen~S.\,D.,  Kochubei~A.\,N.} Analytic methods in the theory of differential and pseudo-differential equations of parabolic type.
"--- Basel: Springer, 2004.

\bibitem{86} {\it  Erd\'elyi~A.} On fractional integration and its application to the Hankel transforms//  Quart. J.~Math. (Oxford). "--- 1940. "--- {\sl  11}. "--- С.~293--303.

\bibitem{Erd1} {\it  Erd\'elyi~A.} Some applications of fractional integration//  Boeing Sci. Res. Labor. "--- 1963. "---  Math. Note {\sl 316}. "--- D1-82-0286.

\bibitem{Erd2} {\it  Erd\'elyi~A.}  An integral equation involving Legendre functions//  SIAM Rev. "--- 1964. "--- {\sl  12}, \No~1. "--- С.~15--30.

\bibitem{Erd3} {\it  Erd\'elyi~A.}  An application of fractional integrals//  J.~Analyse Math. "--- 1965. "--- {\sl  14}. "--- С.~113--126.

\bibitem{Erd4} {\it  Erd\'elyi~A.}  Some integral equations involving finite parts of divergent integrals//  Glasgow Math.~J. "--- 1967. "--- {\sl  8}, \No~1. "--- С.~50--54.

\bibitem{Erd5} {\it  Erd\'elyi~A.} On the Euler---Poisson---Darboux equation//  J.~Analyse Math. "--- 1970. "--- {\sl  23}. "--- С.~89--102.

\bibitem{Rub3}  {\it  Estrada~R., Rubin~B.} Null spaces Of Radon transforms//  arXiv:1504.03766v1. "--- 2015.

\bibitem{Euler} {\it  Euler~L.} Institutiones calculi integralis// Opera Omnia. "--- 1914. "--- {\sl  1}, \No~13. "--- C.~212--230.

\bibitem{Ext} {\it  Exton~H.} Multiple hypergeometric functions and applications. "--- New York: John Wiley and Sons, 1976.

\bibitem{FH}  {\it  Fitouhi~A., Hamza~M.\,M.} A uniform expansion for the eigenfunction of a singular second-order differential operator//
SIAM J. Math. Anal. "--- 1990. "--- {\sl  21}, \No~6. "--- С.~1619--1632.

\bibitem{FJSS} {\it  Fitouhi~A., Jebabli~I., Shishkina~E., Sitnik~S.\,M.}
Applications of integral transforms composition  method to
wave-type singular differential equations and index shift
transmutations// Electron. J. Differ. Equ.
"--- 2018. "--- {\sl  2018}, \No~130. "--- C.~1--27.

\bibitem{Fox} {\it Fox~D.\,N.} The solution and Huygens' principle for a singular Cauchy problem//
J.~Math. Mech. "--- 1959. "--- {\sl 8}. "--- C.~197--219.

\bibitem{Gad1}  {\it Gadjiev~A.\,D.} Selected works. "--- Baku: <<ELM>>, 2003.

\bibitem{Gal1} {\it  Gallardo~L., Trim\'eche~Kh.} Un analogue d'un theoreme de Hardy pour la transformation de Dunkl//
C.\,R. Math. Acad. Sci. Paris. "--- 2002. "--- {\sl  334}. "--- С.~849--854.

\bibitem{Gal2} {\it  Gallardo~L., Trim\'eche~Kh.} A version of Hardy's theorem for the Dunkl transform//  J. Aust. Math. Soc. "--- 2004. "--- {\sl  77}. "--- С.~371--385.

\bibitem{Rad1} {\it  Ghergu~M., Radulescu~V.} Singular elliptic problems. Bifurcation and asymptotic analysis. "---  New York: Oxford University Press, 2008.

\bibitem{Gil2} {\it  Gilbert~R.}  Function theoretic methods in partial differential equations. "--- N.Y.: Academic Press, 1969.

\bibitem{Gil1} {\it  Gilbert~R.}  Constructive methods for elliptic equations. "---  Berlin---Heidelberg: Springer, 1974.

\bibitem{GB} {\it  Gilbert~R., Begehr~H.} Transformations, transmutations and kernel functions. V.~1-2. "--- Harlow: Longman, 1992.

\bibitem{GSS}   {\it  Glaeske~H.\,J., Prudnikov~A.\,P., Skornik~K.\,A.} Operational calculus and related topics. "--- Boca Raton: Chapman \& Hall/CRC, 2006.

\bibitem{GSPP} {\it Golenia~J., Samoilenko~A.\,M., Prykarpatsky~Ya.~A., Prykarpatsky~A.\,K.} The general differential-geometric structure of
multidimensional Delsarte transmutation operators in parametric
functional spaces and their applications in soliton theory//
arXiv:~math-ph/0404016. "--- 2004.

\bibitem{GKMR} {\it  Gorenflo~R., Kilbas~A.\,A., Mainardi~F., Rogosin~S.\,V.} Mittag-Leffler functions,
related topics and applications. "--- Heidelberg---New York---Dordrecht---London: Springer, 2014.

\bibitem{GoMa} {\it  Gorenflo~R.,  Mainardi~F.}   Fractional calculus: integral and differential equations of fractional order//
В сб.: <<Fractals and fractional calculus in continuum mechanics>>. "--- Wien---New York: Springer, 1997. "--- C.~223--278.

\bibitem{HaKa1} {\it Hasanov~A., Karimov~E.\,T.} Fundamental solutions for a class of three-dimensional elliptic equations with singular coefficients//
Appl. Math. Lett. "---  2009. "---  {\sl  22}. "--- С.~1828--1832.

\bibitem{Hig1} {\it  Higgins~T.\,P.} A hypergeometric function transform//  J.~SIAM. "---
1964. "--- {\sl  12}, \No~3. "--- С.~601--612.

\bibitem{HT} {\it  Hille~E., Tamarkin~J.\,D.} On the theory of linear integral equations//  Ann. Math. "--- 1930. "--- {\sl  31}. "--- C.~479--528.

\bibitem{Hol} {\it  Holzleitner~M.} Transformation operators for spherical
Schr\"odinger operators// arXiv:1805.10526v1. "--- 2018.

\bibitem{Jaf} {\it  Jafford~K.} Inversion of the Lions transmutation operators using generalized wavelets//  Appl. Comput. Harmon. Anal. "--- 1997. "--- {\sl  4}, \No~1. "--- С.~97--112.

\bibitem{Kam} {\it  Kamoun~L., Sifi~M.} Bessel---Struve intertwining operator and generalized Taylor series on the real line//  Integral Transforms Spec. Funct. "--- 2005. "--- {\sl  16}, \No~1. "--- С.~39--55.

\bibitem{KarST} {\it Karimov~S.\,T.} Multidimensional generalized Erd\'elyi---Kober
operator and its application to solving Cauchy problems for
differential equations with singular coefficients// Fract. Calc. Appl. Anal. "--- 2015. "--- {\sl  18}, \No~4. "--- C.~845--861.

\bibitem{KaSr} {\it  Karlsson~P.\,W., Srivastava H.\,M.} Multiple Gaussian hypergeometric series. "--- New York: Ellis Horwood,  1985.

\bibitem{KaLo1}{\it Karp~D., L\'opez~J.\,L.} Representations of hypergeometric
functions for arbitrary values of the parameters and their use//
J.~Approx. Theory. "---  2017. "--- {\sl 218}. "---
C.~42--70.

\bibitem{KaLo2}{\it Karp~D., L\'opez~J.\,L.}  A class of Meijer’s $G$ functions and further representations of the generalized hypergeometric functions// arXiv:1801.08670v1. "--- 2018.

\bibitem{KaPr4}{\it  Karp~D., Prilepkina~E.}  Hypergeometric functions as generalized
Stieltjes transforms// J.~Math. Anal. Appl. "--- 2012. "--- {\sl
393}, \No~2. "--- C.~348--359.

\bibitem{KaPr2} {\it  Karp~D., Prilepkina~E.}  Completely monotonic gamma ratio and
infinitely divisible $H$-function of Fox// Comput. Methods Funct. Theory. "--- 2016. "--- {\sl 16}, \No~1. "---
C.~135--153.

\bibitem{KaPr3}{\it  Karp~D., Prilepkina~E.}  Hypergeometric differential equation and
new identities for the coefficients of N\o rlund and
B\"uhring// SIGMA Symmetry Integrability Geom. Methods Appl. "--- 2016. "--- {\sl 052}.

\bibitem{KaPr1} {\it  Karp~D., Prilepkina~E.}  Some new facts concerning the delta neutral case of Fox's $H$-function//
Comput. Methods Funct. Theory. "--- 2017. "--- {\sl 17}, \No~2. "---  C.~343--367.

\bibitem{S13} {\it  Karp~D., Savenkova~A., Sitnik~S.\,M.}  Series expansions for the third incomplete elliptic integral via partial fraction decompositions//
J.~Comput. Appl. Math. "--- 2007. "--- {\sl  207},  \No~2. "--- С.~331--337.

\bibitem{S12} {\it  Karp~D., Sitnik~S.\,M.}  Asymptotic approximations for the first incomplete elliptic integral near logarithmic singularity//
J.~Comput. Appl. Math. "--- 2007. "--- {\sl  205}, \No~1. "--- С.~186--206.

 \bibitem{S15} {\it  Karp~D., Sitnik~S.\,M.}  Inequalities and monotonicity of ratios for generalized  hypergeometric function//
 J.~Approx. Theory. "--- 2009. "--- {\sl  161}, \No~1. "--- С.~337--352.

\bibitem{S16} {\it  Karp~D., Sitnik~S.\,M.}  Log-convexity and log-concavity of hypergeometric-like functions//  J.~Math. Anal. Appl. "--- 2010. "--- {\sl  364}, \No~2. "--- С.~384--394.

\bibitem{Ken1}{\it  Kenig~C.\,E.} Some recent quantitative unique continuation theorems//
S\'emin. \'Equ. D\'eriv. Partielles. \'Ec. Polytech. Cent. Math. Palaiseau "---   2006. "---  {\sl 2005-2006}, XX1--XX10.

\bibitem{Ken2}{\it  Kenig~C.\,E., Silvestre~L., Wang~J.\,N.} On Landis' conjecture in the plane// Commun. Part. Differ. Equ. "--- 2015. "--- {\sl 40}, \No~4. "--- C.~766--789.

\bibitem{KiSa} {\it  Kilbas~A.\,A., Saigo~M.} $H$-transforms. Theory and applications. "--- Boca Raton: Chapman \& Hall, CRC, 2004.

\bibitem{KiSk1} {\it  Kilbas~A.\,A., Skoromnik~O.\,V.} Integral transforms
with the Legendre function of the first kind in the kernels on
$\mathcal{L}_{\protect\nu,r}$-spaces//  Integral Transforms Spec. Funct. "--- 2009. "--- {\sl  20}, \No~9. "---
С.~653--672.

\bibitem{KST} {\it  Kilbas~A.\,A., Srivastava~H.\,M., Truhillo~J.\,J.} Theory and applications of fractional differential equations. "--- Amsterdam, etc.: Elsevier, 2006.

\bibitem{KT1} {\it  Kilbas~A.\,A., Truhillo~J.\,J.}   Differential equations of fractional order: methods, results and problems. Part~I//
Appl. Anal. "--- 2001. "--- {\sl  78}, \No~1-2. "--- С.~153--192.

\bibitem{KT2} {\it  Kilbas~A.\,A., Truhillo~J.\,J.}  Differential equations of fractional order: methods, results and problems. Part~II//
Appl. Anal. "--- 2001. "--- {\sl  81}, \No~2. "--- С.~435--493.

\bibitem{KiZhu} {\it  Kilbas~A.\,A., Zhukovskaya~N.\,V.} Euler-type non-homogeneous differential equations with three Liouville fractional derivatives//
Fract. Calc. Appl. Anal. "--- 2009. "--- {\sl 12}, \No~2. "--- C.~205--234.

\bibitem{Kir2} {\it  Kiryakova~V.}  An explanation of Stokes phenomenon by the method of transmutations// Proc. conf. Diff. Equations and Appl., Rousse. "--- 1982. "--- C.~349--353.

\bibitem{Kir1} {\it  Kiryakova~V.} Generalized fractional calculus and applications. "--- Harlow: Longman, 1994.

\bibitem{Kir3} {\it  Kiryakova~V.}  Applications of the generalized Poisson transformation for solving hyper-Bessel differential equations//
Godishnik VUZ. Appl. Math. "--- 1986. "--- {\sl  22}, \No~4. "--- С.~129--140 (in Bulgarian).

\bibitem{Kir4} {\it  Kiryakova~V.}  All the special functions are fractional differintegrals of elementary functions//
J.~Phys. A. Math. Gen. "--- 1997. "--- {\sl  30}, \No~14. "--- С.~5085--5103.

\bibitem{Kir6} {\it  Kiryakova~V.} Multiple (multiindex) Mittag-Leffler functions and relations to generalized fractional calculus//
J.~Comput. Appl. Math. "--- 2000. "--- {\sl  118}, \No~1-2. "--- С.~241-­259.

\bibitem{Kir5} {\it  Kiryakova~V.} The multi-index Mittag-Leffler functions as an important class of special functions of fractional calculus//
Computers Math. Appl. "--- 2010. "--- {\sl  59}, \No~5. "--- С.~1885--1895.

\bibitem{47}{\it  Kober~H.} On fractional integrals and derivatives//  Quart. J. Math. (Oxford). "--- 1940. "--- {\sl  11}. "--- С.~193--211.

\bibitem{Kob1} {\it  Kober~H.} On a theorem of Schur and on fractional integrals of purely imaginary order//
Trans. Am. Math. Soc. "--- 1941. "--- {\sl  50}, \No~1. "--- С.~160--174.

\bibitem{Koepf} {\it  Koepf~W.} Hypergeometric summation. An algorithmic approach to summation
and special function identities.
"--- Wiesbaden: Vieweg, 1998.

\bibitem{KMRS} {\it  Kokilashvili~V., Meskhi~A., Rafeiro~H., Samko~S.} Integral operators in non-standard function spaces. Vol.~1-2. "--- Basel: Birkh\"auser/Springer, 2016.

\bibitem{Koo1} {\it  Koornwinder~T.\,H.} Fractional integral and generalized Stieltjes transforms for
hypergeometric functions as transmutation operators// SIGMA Symmetry Integrability Geom. Methods Appl. "---
2015. "--- {\sl 11}, \No~074.

\bibitem{Kravch} {\it  Kravchenko~V.\,V.} Applied pseudoanalytic function theory. "--- Basel: Birkh\"auser, 2009.

\bibitem{Krav8}{\it  Kravchenko~V.\,V.} Construction of a transmutation for the one-dimensional Schr\"odinger operator and a representation for solutions//
Appl. Math. Comput. "--- 2018. "--- {\sl  328}. "--- C.~75--81.

\bibitem{Krav4}  {\it  Kravchenko~V.\,V., Navarro~L.\,J., Torba~S.\,M.} Representation of solutions to the one-dimensional Schr\"odinger equation in terms of Neumann series
of Bessel functions//  Appl. Math. Comput. "--- 2017. "--- {\sl  314}, \No~1. "--- С.~173--192.

\bibitem{Krav6}{\it  Kravchenko~V.\,V., Otero~J.\,A., Torba~S.\,M.} Analytic approximation of solutions of parabolic partial differential equations with variable coefficients//
Adv. Math. Phys. "--- 2017. "--- {\sl  2017}. "--- 2947275.

\bibitem{Krav1}  {\it  Kravchenko~V.\,V.,   Torba~S.\,M.}  Transmutations for Darboux transformed operators with applications//
J.~Phys. A. Math. Theor. "--- 2012. "--- {\sl  45}, \No~7. "--- 075201.

\bibitem{Krav3}  {\it  Kravchenko~V.\,V.,   Torba~S.\,M.}  Analytic approximation of transmutation operators and applications to highly accurate solution of spectral problems//
J.~Comput. Appl. Math. "--- 2015. "--- {\sl  275}. "--- С.~1--26.

\bibitem{Krav2}  {\it  Kravchenko~V.\,V., Torba~S.\,M.}  Construction of transmutation operators and hyperbolic pseudoanalytic functions//
Complex Anal. Oper. Theory. "--- 2015. "--- {\sl  9}, \No~2. "--- С.~379--429.

\bibitem{Krav7}{\it  Kravchenko~V.\,V., Torba~S.\,M.} Asymptotics with respect to the spectral parameter and Neumann series of Bessel functions for solutions
of the one-dimensional Schr\"odinger equation// J.~Math. Phys. "--- 2017. "--- {\sl  58}, \No~12. "--- 122107.

\bibitem{Krav9}{\it  Kravchenko~V.\,V., Torba~S.\,M.} A Neumann series of Bessel functions representation for solutions of Sturm---Liouville equations//
Calcolo. "--- 2018. "--- {\sl  55}, \No~1. "--- Paper No.~11.

\bibitem{Krav10}{\it  Kravchenko~V.\,V., Torba~S.\,M., Castillo-P\'erez~R.} A Neumann series of Bessel functions representation for solutions of perturbed Bessel equations//
Appl. Anal. "--- 2018. "--- {\sl  97}, \No~5. "--- С.~677--704.

\bibitem{Krav5}{\it  Kravchenko~V.\,V., Torba~S.\,M., Khmelnytskaya~K.\,V.} Transmutation operators: construction and applications//
Proc. of the 17th Int. Conf. on Comput. and Math. Methods in Sci. and Eng. CMMSE-2017. "--- Cadiz, Andalucia, Espa\~na, 2017. "---  C.~1198--1206.

\bibitem{KMP} {\it  Kufner~A., Maligranda~L., Persson~L.-E.}
The Hardy inequality. About its history and some related results. "--- Pilsen: Vydavatelsk\'y Servis, 2007.

\bibitem{Lev1995} {\it Levitan~B.\,M.} Transmutation operators and the inverse spectral problem// Contemp. Math. "--- 1995. "--- {\sl 183}. "--- C.~237--244.

\bibitem{Lind} {\it  Lindqvist~P.} Notes on the $p$-Laplace equation// Report: Univ. of Jyv\"askyl\"a. Dep. of Math. and Stat. "--- 2006. "--- {\sl 102}.

\bibitem{Lio2} {\it  Lions~J.\,L.}  Op\'erateurs de Delsarte et probl\`eme mixte//  Bull. Soc. Math. France. "--- 1956. "--- {\sl 84}. "--- С.~9--95.

\bibitem{Lio3} {\it  Lions~J.\,L.}   Quelques applications  d'op\'erateurs de transmutations// Colloques Internat. Nancy. "--- 1956. "--- C.~125--142.

\bibitem{Lio1}  {\it  Lions~J.\,L.} Equations differentielles operationnelles et probl\'emes aux limites. "---  Berlin---G\"ottingen---Heidelberg: Springer, 1961.

\bibitem{Love1} {\it  Love~E.\,R.} Some integral equations involving hypergeometric functions//  Proc. Edinb. Math. Soc. "--- 1967. "--- {\sl  15}, \No~3. "--- С.~169--198.

\bibitem{Love2} {\it  Love~E.\,R.}  Two more hypergeometric integral equations//  Proc. Cambridge Phil. Soc. "--- 1967. "--- {\sl  63}, \No~4. "--- С.~1055--1076.

\bibitem{Low1} {\it  Lowndes~J.\,S.} An
application of some fractional integrals//  Glasg. Math. J. "---
1979. "--- {\sl  20}, \No~1. "--- С.~35--41.

\bibitem{Low2} {\it  Lowndes~J.\,S.} On some generalizations of Riemann---Liouville
and Weil fractional integrals and their applications//  Glasg. Math. J. "--- 1981. "--- {\sl  22}, \No~2. "--- С.~73--80.

\bibitem{Low3} {\it  Lowndes~J.\,S.} Cauchy problems for second order hyperbolic differential
equations with constant coefficients//  Proc. Edinb. Math. Soc.
"--- 1983. "--- {\sl  26}, \No~3. "--- С.~97--105.

\bibitem{Lud} {\it  Ludwig~D.}  The Radon transform on Euclidean space//  Math. Methods Appl. Sci. "--- 1966. "--- {\sl  19}. "--- C.~49--81.

\bibitem{Luke2}  {\it  Luke~Y.\,L.} The special functions and their approximations. V.~1. "--- New York---London: Academic Press, 1969.

\bibitem{Luke1} {\it  Luke~Y.\,L.} Mathematical functions and their approximations. "--- New York---San Francisco---London: Academic Press, 1975.

\bibitem{Matv} {\it  Matveev~V.\,B.} Intertwining relations between the Fourier transform
and discrete Fourier transform, the related functional identities
and beyond//  Inverse Problems. "--- 2001. "--- {\sl  17}. "---
С.~633--657.

\bibitem{MaSa}  {\it  Matveev~V.\,B., Salle~V.\,B.} Darboux---Backlund transformations and applications. "--- N.Y.: Springer, 1991.

\bibitem{McB} {\it  McBride~A.\,C.} Fractional calculus and integral transforms of generalized functions. "--- San Francisco---London---Melbourne: Pitman, 1979.

\bibitem{MS1} {\it  Mehrez~Kh., Sitnik~S.\,M.}   On monotonicity of ratios of some $q$-hypergeometric functions// Mat. Vesn. "--- 2016. "--- {\sl  68}, \No~3. "--- С.~225--231.

\bibitem{MS3} {\it  Mehrez~Kh., Sitnik~S.\,M.}   Proofs of some conjectures on monotonicity of ratios of Kummer, Gauss and generalized hypergeometric functions//
Analysis (Munich). "--- 2016. "--- {\sl  36}, \No~4. "--- С.~263--268.

\bibitem{MS2} {\it  Mehrez~Kh., Sitnik~S.\,M.}    Functional Inequalities for the Mittag-Leffler Functions//  Results Math. "--- 2017. "--- {\sl  72}, \No~1-2. "--- С.~703--714.

\bibitem{Rad3} {\it  Mitidieri~E., Serrin~J., Radulescu~V. {\rm (}eds.{\rm )}} Recent trends in nonlinear partial differential equations~I: Evolution problems.
 "--- Providence: Am. Math. Soc., 2013.

\bibitem{Rad4} {\it  Mitidieri~E., Serrin~J., Radulescu~V. {\rm (}eds.{\rm )}} Recent trends in nonlinear partial differential equations~II: Stationary Problems.
"--- Providence: Am. Math. Soc., 2013.

\bibitem{MPF} {\it  Mitrinovi\'c~D.\,S., Pe\v cari\'c~J.\,E., Fink~A.\,M.} Classical and new inequalities in analysis. "--- Dordrecht: Kluwer, 1993.

\bibitem{Mur7} {\it  Muravnik~A.\,B.} On weighted norm estimates for the mixed Fourier---Bessel transforms on non-negative functions//
В сб.: <<Integral methods in science and engineering. Vol.~1: Analytic methods>>. "---  Harlow: Longman, 1997. "--- С.~119--123.

\bibitem{Mur6} {\it  Muravnik~A.\,B.}  Fourier---Bessel transformation of measures and singular differential operators//
В сб.: <<Paul Erd\H{o}s and his mathematics>>. "---
Budapest: J'anos Bolyai Math. Soc., 1999. "--- С.~182--184.

\bibitem{Mur5} {\it  Muravnik~A.\,B.}  Fourier---Bessel transformation of
measures with several special variables and properties of singular
differential equations// J.~Korean Math. Soc. "--- 2000. "--- {\sl
37}, \No~6. "--- С.~1043--1057.

\bibitem{Mur4} {\it  Muravnik~A.\,B.}  Fourier---Bessel transformation of compactly
supported non-negative functions and estimates of solutions of
singular differential equations//  Funct. Differ. Equ. "--- 2001. "--- {\sl  8}, \No~3-4. "--- С.~353--363.

\bibitem{Mur3} {\it  Muravnik~A.\,B.}  Fourier---Bessel transformation of measures and singular differential equations//
В сб.: <<Recent progress in functional analysis>>. "--- Amsterdam: Elsevier, 2001. "--- С.~335--345.

\bibitem{Mur2} {\it  Muravnik~A.\,B.}  Nonclassical Cauchy problem for singular parabolic integro-differential equations//
Russ. J. Math. Phys. "--- 2002. "--- {\sl  9}, \No~3. "--- С.~300--314.

\bibitem{Mur1} {\it  Muravnik~A.\,B.}  On stabilization of solutions of elliptic equations containing Bessel operators//
В сб.: <<Integral methods in science and engineering. Analytic and numerical techniques>>. "---
Boston: Birkh\"auser, 2004. "--- С.~157--162.

\bibitem{NIST}  {\it  Olver F.\,W.\,J., Lozier~D.\,W., Boisvert~R.\,F., Clark~C.\,W.} NIST handbook  of  mathematical  functions. "--- Cambridge: Cambridge University Press,  2010.

\bibitem{OpKu} {\it  Opic~B., Kufner~A.} Hardy-type inequalities. "--- Harlow: Longman, 1990.

\bibitem{OZK} {\it  Ozaktas~H., Zalevsky~Z., Kutay~M.} The fractional Fourier transform: with applications in optics and signal processing. "--- Chishester, ets.: John Wiley \& Sons, 2001.

\bibitem{Jor}  {\it  Paneva-Konovska~J.} From Bessel to multi-index Mittag-Leffler functions. "--- London: World Scientific, 2016.

\bibitem{PWZ}  {\it  Petkov\v sek~M.,  Wilf~H.\,S., Zeilberger~D.} $A=B$. "--- Wellesley: A.\,K. Peters, 1996.

\bibitem{PiSa} {\it  Pike~S., Sabatier~P.} Scattering. Scattering and inverse scattering in pure and applied science. Vol.~1-2. "--- San Diego: Academic Press, 2002.

\bibitem{Poisson} {\it  Poisson~S.\,D.}  M\'emoire sur l'int\'egration des \'equations lin\'eaires
            aux diff\'erences partielles//   J.~\'Ec. Roy. Polytech. "--- 1823. "--- {\sl 19}, \No~12. "--- C.~215--248.

\bibitem{Pou} {\it  Poularicas~A.\,D. $($ed.$)$} The transforms and applications handbook. "--- Boca Raton: CRC Press, 2010.

\bibitem{Jan} {\it Pr\"uss~J.} Evolutionary integral equations and applications. "--- Basel: Birkh\"auser, 2012.

\bibitem{Pyat} {\it  Pyatkov~S.\,G.}  Operator theory. Nonclassical problems. "--- Utrecht: VSP, 2002.

\bibitem{Rad2} {\it  Radulescu~V., Repovs~D.} Partial differential equations with variable exponents: Variational methods and qualitative analysis. "---
Boca Raton: CRC Press, 2015.

\bibitem{Riesz} {\it  Riesz~M.} L'int\'egrale de Riemann---Liouville et le probl\`eme de Cauchy//  Acta Math. "--- 1949. "--- {\sl 81}. "--- C.~1--223.

\bibitem{Riman} {\it  Riemann~B.}  On the propagation of flat waves of finite amplitude// В сб.: <<Ouvres>>. "--- Moscow---Leningrad: OGIZ, 1948. "--- C.~376--395.

\bibitem{Rod} {\it  Rodrigues~J.} Operational calculus for the generalized Bessel operator// Serdica Math.~J. "--- 1989. "--- {\sl  15}. "--- С.~179--186.

\bibitem{Ros1} {\it  R\"osler~M.} Positivity of Dunkl’s intertwining operator//  Duke Math.~J. "--- 1999. "--- {\sl  98}. "--- С.~445--463.

\bibitem{Ros2} {\it  R\"osler~M.}  Dunkl operators: theory and applications//
В сб.: <<Orthogonal polynomials and special functions>>. "--- Berlin: Springer, 2003. "--- C.~93--135.

\bibitem{Rossi} {\it Rossi~L.}  The Landis conjecture with sharp rate of decay// arXiv:1807.00341v1. "--- 2018.

\bibitem{Rub4} {\it  Rubin~B.} On the Funk---Radon---Helgason inversion method in integral geometry// Contemp. Math. "--- 2013. "--- {\sl  599}. "--- С.~175--198.

\bibitem{Rub2}  {\it  Rubin~B.} Gegenbauer---Chebyshev integrals and Radon transforms// arXiv:1410.4112v2. "--- 2015.

\bibitem{Rub1}  {\it  Rubin~B.} Radon transforms and Gegenbauer---Chebyshev integrals.~I//  Anal. Math. Phys. "--- 2017. "--- {\sl 7}, \No~2. "--- С.~117--150.

\bibitem{SaHa1} {\it  Salakhitdinov~M.\,S., Hasanov~A.} The Dirichlet problem for the generalized bi-axially symmetric Helmholtz equation//
Eurasian Math.~J. "--- 2012. "--- {\sl  3}, \No~4. "--- С.~99--110.

\bibitem{SaKiMar} {\it  Samko~S.\,G., Kilbas~A.\,A., Marichev~O.\,I.} Fractional integrals and derivatives: theory and applications. "--- N.Y.: Gordon \& Breach, 1993.

\bibitem{SPP} {\it  Samoilenko~A.\,M., Prykarpatsky~Ya.\,A., Prykarpatsky~A.\,K.} The generalized de Rham---Hodge theory aspects of Delsarte---Darboux type transformations
in multidimension//  Cent. Eur. J. Math. "--- 2005. "--- {\sl  3}, \No~3. "--- С.~529--557.

\bibitem{74} {\it  Seely~R.\,T.} Extensions of $C^{\infty}$ functions defined in a half space//  Proc. Am. Math. Soc. "--- 1964. "--- {\sl  15}. "--- С.~625--626.

\bibitem{75} {\it  Seely~R.\,T.}  Complex powers of an elliptic operator//  Proc. Sympos. Pure Math. "--- 1967. "--- {\sl  10}. "--- C.~288--307.

\bibitem{ShiE1}  {\it Shishkina~E.\,L.}   Inversion of integral of $B$-potential type with density from $\Phi_\gamma$//  J.~Math. Sci. "--- 2009. "--- {\sl  160}, \No~1. "--- С.~95--102.

\bibitem{ShiE2} {\it Shishkina~E.\,L.}   On the boundedness of hyperbolic Riesz $B$-potential//  Lith. Math.~J. "--- 2016. "--- {\sl  56}, \No~4. "--- С.~540--551.

\bibitem{ShiE5} {\it Shishkina~E.\,L.}    Generalized Euler---Poisson---Darboux equation and
 singular Klein---Gordon equation//  J.~Phys. Conf. Ser. "--- 2018. "--- {\sl  973}. "--- С.~1--21.

 \bibitem{ShiE4} {\it Shishkina~E.\,L.}  Properties of mixed hyperbolic $B$-potential//  Progr. Fract. Differ. Appl. "--- 2018. "--- {\sl  4}, \No~2. "--- С.~1--16.

\bibitem{ShiE6}  {\it Shishkina~E.\,L.}   Singular Cauchy problem for the general Euler---Poisson---Darboux equation//  Open Math. "--- 2018. "--- {\sl  16}. "--- С.~23--31.

\bibitem{ShiE3}  {\it Shishkina~E.\,L., Sitnik~S.\,M.}   General form of the Euler---Poisson---Darboux equation and application of the transmutation method//
Electron. J.~Differ. Equ. "--- 2017. "--- {\sl  177}. "--- С.~1--20.

\bibitem{SS} {\it Shishkina~E.\,L., Sitnik~S.\,M.}
On fractional powers of Bessel operators//  J.~Inequal. Spec. Funct. "--- 2017. "--- {\sl  8}, \No~1.
"--- С.~49--67.

\bibitem{Sie} {\it  Siersma~J.} Thesis. "--- Groningen, 1979.

\bibitem{S41} {\it  Sitnik~S.\,M.} Generalized Young and Cauchy---Bunyakowsky inequalities with applications: a survey// [math.CA] arXiv:1012.3864. "--- 2010.

\bibitem{S42}  {\it  Sitnik~S.\,M.} Transmutations and applications: a survey// [math.CA] arXiv:1012.3741. "---  2010.

\bibitem{S94} {\it  Sitnik~S.\,M.} Some problems in the modern theory of transmutations//
Spectral theory and differential equations (STDE - 2012).
Int. Conf. in honor of V.\,A.~Marchenko's 90~th
birthday. "--- Kharkiv, 2012. "---
C.~101-102.

\bibitem{S38} {\it  Sitnik~S.\,M.}  Buschman---Erd\'elyi transmutations, classification and applications//
В сб.: <<Analytic methods of analysis   and differential
equations>>. "--- Cambridge: Cambridge Scientific Publishers,
2013. "--- C.~171--201.

\bibitem{S401} {\it  Sitnik~S.\,M.}  A short survey of recent results on Buschman---Erd\'elyi transmutations//
J. Inequal. Spec. Funct. "--- 2017. "--- {\sl  8}, \No~1. "--- С.~140--157.

\bibitem{S402} {\it  Sitnik~S.\,M.}   Buschman---Erd\'elyi transmutations and  applications//
Abstr. of the 8th Int. Conf. <<Transform Methods
and Special Functions>>, Bulgaria, Sofia, 27--31 Aug. 2017.
"--- Inst. Math. Inf. Bulg. Acad. Sci., 2017. "--- C.~59.

\bibitem{Skub} {\it  Skubachevskii~A.\,L.} Elliptic functional differential equations and applications. "--- Basel: Birkh\"auser, 1997.

\bibitem{Sla} {\it  Slater~L.\,J.} Generalized hypergeometric functions. "--- Cambridge: Cambridge University Press, 1966.

\bibitem{Spr} {\it  Sprinkhuizen-Kuyper~I.\,G.} A fractional integral operator corresponding to negative powers of a certain second-order differential operator//
J.~Math. Anal. Appl. "--- 1979. "--- {\sl 72}. "--- С.~674--702.

\bibitem{SvFe} {\it  Sviridyuk~G.\,A., Fedorov~V.\,E.} Linear Sobolev type equations and degenerate semigroups of operators. "--- Utrecht: VSP, 2003.

\bibitem{Ta1} {\it  Ta Li} A new class of integral transform// Proc. Am. Math. Soc. "--- 1960. "--- {\sl  11}, \No~2. "--- С.~290--298.

\bibitem{Ta2} {\it  Ta Li}  A note on integral transform//  Proc. Am. Math. Soc. "--- 1961. "--- {\sl  12}, \No~6. "--- С.~556.

\bibitem{Tri2} {\it  Trim\'eche~Kh.} Transformation int\'egrale de Riemann---Liouville g\'en\'eralises et convergence des series de Taylor g\'en\'eralis\'ees au sens de Delsarte//
Rev. Fac. Sci. Tunis. "--- 1981. "--- \No~1. "--- С.~7--14.

\bibitem{Tri1} {\it  Trim\'eche~Kh.} Transformation int\'egrale de Weil et th\`eor\'eme de Paley---Winer associ\'es \`a un op\'erateur diff\'erentiel singulier sur $(0,\infty)$//
J.~Math. Pures Appl. "--- 1981. "--- {\sl 60}. "--- С.~51--98.

\bibitem{Tri3} {\it  Trim\'eche~Kh.}  Transmutation operators and mean-periodic functions associated with differential ope\-ra\-tors//
Math. Rep. Ser.~4. "--- 1988. "--- \No~1, i-xiv.

\bibitem{Tri4} {\it  Trim\'eche~Kh.} Inversion of the Lions transmutation operators using generalized wavelets//  Appl. Comput. Harmon. Anal. "---
1997. "--- {\sl  4}, \No~1. "--- С.~97--112.

\bibitem{Tri5} {\it  Trim\'eche~Kh.} Generalized harmonic analysis and wavelet packets. "--- Amsterdam: Gordon and Breach, 2001.

\bibitem{Tri6} {\it  Trim\'eche~Kh.}  Inversion of the Dunkl intertwining
operator and its dual using Dunkl wavelets//  Rocky Mountain J.~Math. "--- 2002. "--- {\sl  32}, \No~2. "---
С.~889--895.

\bibitem{Vir2} {\it  Virchenko~N.} On some generalized symmetric integral operators of Buschman---Erd\'elyi’s type//
J.~Nonlinear Math. Phys. "--- 1996. "--- {\sl  3}, \No~3-4. "--- С.~421--425.

\bibitem{Vir1} {\it  Virchenko~N., Fedotova~I.} Generalized associated Legendre functions and their applications. "--- Singapore: World Scientific, 2001.

\bibitem{Volch} {\it  Volchkov~V.\,V.} Integral geometry and convolution equations. "---  Dordrecht: Kluwer, 2003.

\bibitem{14} {\it  Wiener~N.} The Dirichlet problem//  J.~Math. Phys. Mass. Inst. Techn. "--- 1924. "--- {\sl  3}. "--- С.~127--146.

\bibitem{Wei1} {\it  Weinstein~A.} Discontinuous integrals and generalized theory of potential// Trans. Am. Math. Soc. "--- 1948. "--- {\sl  63}, \No~2. "--- С.~342--354.

\bibitem{Wei2} {\it  Weinstein~A.}  Generalized axially symmetric potential theory//  Bull. Am. Math. Soc. "--- 1953. "--- {\sl  59}. "--- С.~20--38.

\bibitem{Wei3} {\it  Weinstein~A.} Selecta. "--- London---San Francisco: Pitman, 1978.

\label{BibEnd}
\end{thebibliography}
\end{document}